\documentclass[11pt]{scrartcl}

\usepackage[utf8]{inputenc}
\usepackage[T1]{fontenc}
\usepackage{lmodern}
\usepackage{alphabeta}

\usepackage{titling}
\usepackage{xspace}
\usepackage{soul} 

\usepackage[letterpaper, top=25.4mm, bottom=25.4mm, left=25.4mm, right=25.4mm, includefoot]{geometry}

\DeclareSectionCommand[%
level=4,
indent=0pt,
beforeskip=1ex plus 1ex minus .2ex,
afterskip=-1em,
font={},
tocindent=7em,
tocnumwidth=4.1em,
counterwithin=subsubsection
]{paragraph}

\usepackage[inline]{enumitem}

\usepackage{dsfont}
\usepackage{setspace}

\usepackage{bm}
\usepackage{microtype}
\usepackage{amsmath}
\usepackage{amssymb}
\usepackage{amsfonts}
\usepackage{mathtools}
\usepackage[mathscr]{euscript}

\usepackage[hyphens]{url}
\usepackage{amsthm}

\usepackage{tikz}
\usepackage{xcolor}
\usepackage{subcaption}
\usepackage{graphicx}
\usepackage{wrapfig}

\usepackage{hyperref}
\usepackage{cleveref}

\usepackage{nicefrac}

\colorlet{myGreen}{green!50!black}
\colorlet{myLightgreen}{green}
\colorlet{myRed}{red!90!black}
\definecolor{myBlue}{rgb}{0.25, 0.0, 1.0}
\definecolor{myLightBlue}{rgb}{0.39, 0.58, 0.93}
\colorlet{myViolet}{myBlue!55!myRed}
\definecolor{myOrange}{rgb}{1.0, 0.66, 0.07}

\definecolor{CornflowerBlue}{rgb}{0.39, 0.58, 0.93}
\definecolor{DarkGoldenrod}{rgb}{0.72, 0.53, 0.04}
\definecolor{BritishRacingGreen}{rgb}{0.0, 0.26, 0.15}
\definecolor{DarkMagenta}{rgb}{0.55, 0.0, 0.55}
\definecolor{AO}{rgb}{0.0, 0.5, 0.0}
\definecolor{BostonUniversityRed}{rgb}{0.8, 0.0, 0.0}
\definecolor{myRed}{rgb}{0.8, 0.0, 0.0}
\definecolor{DarkMidnightBlue}{rgb}{0.0, 0.2, 0.4}
\definecolor{DarkTangerine}{rgb}{1.0, 0.66, 0.07}
\definecolor{AppleGreen}{rgb}{0.55, 0.71, 0.0}
\definecolor{BrightUbe}{rgb}{0.82, 0.62, 0.91}
\definecolor{Amethyst}{rgb}{0.6, 0.4, 0.8}
\definecolor{DarkGray}{rgb}{0.52, 0.52, 0.51}
\definecolor{Gray}{rgb}{0.66, 0.66, 0.66}
\definecolor{BananaYellow}{rgb}{1.0, 0.88, 0.21}
\definecolor{Amber}{rgb}{1.0, 0.75, 0.0}
\definecolor{LightGray}{rgb}{0.83, 0.83, 0.83}
\definecolor{PrincetonOrange}{rgb}{1.0, 0.56, 0.0}
\definecolor{DeepCarrotOrange}{rgb}{0.91, 0.41, 0.17}
\definecolor{CarrotOrange}{rgb}{0.93, 0.57, 0.13}
\definecolor{MidnightBlue}{rgb}{0.1, 0.1, 0.44}
\definecolor{Magenta}{rgb}{0.50, 0.0, 0.50}
\definecolor{BrightPink}{rgb}{1.0, 0.0, 0.5}
\definecolor{BrilliantRose}{rgb}{1.0, 0.33, 0.64}
\definecolor{ChromeYellow}{rgb}{1.0, 0.65, 0.0}
\definecolor{HotMagenta}{rgb}{1.0, 0.11, 0.81}
\definecolor{Amethyst}{rgb}{0.6, 0.4, 0.8}

\vfuzz \hfuzz
\setlist[itemize]{topsep=0pt,partopsep=0pt,itemsep=0pt,parsep=0pt}
\setlist[itemize,1]{label={\small\textbullet}}
\setlist[itemize,2]{label={\tiny\textbullet}}
\setlist[itemize,3]{label=$\cdot$}
\setlist[enumerate]{topsep=0pt,partopsep=0pt,itemsep=0pt,parsep=0pt}
\setlist[enumerate,1]{label=\roman*)}
\setlist[enumerate,2]{label=\alph*)}
\setlist[enumerate,3]{label=\arabic*)}

\hypersetup{
colorlinks=true,
linkcolor=AO!65!black,
citecolor=AO!65!black,
urlcolor=AppleGreen!65!black,
bookmarksopen=true,
bookmarksnumbered,
bookmarksopenlevel=2,
bookmarksdepth=3
}

\theoremstyle{definition}

\newtheorem{environment}{Environment}[section]

\newtheorem{lemma}[environment]{Lemma}
\crefname{lemma}{lemma}{lemmata}
\crefformat{lemma}{#2Lemma~#1#3}
\Crefformat{lemma}{#2Lemma~#1#3}

\newtheorem*{lemma*}{Lemma}
\crefname{lemma*}{lemma}{lemmata}
\crefformat{lemma*}{#2Lemma~#1#3}
\Crefformat{lemma*}{#2Lemma~#1#3}

\newtheorem{proposition}[environment]{Proposition}
\crefname{proposition}{proposition}{propositions}
\crefformat{proposition}{#2Proposition~#1#3}
\Crefformat{proposition}{#2Proposition~#1#3}

\newtheorem{corollary}[environment]{Corollary}
\crefname{corollary}{corollary}{corollaries}
\crefformat{corollary}{#2Corollary~#1#3}
\Crefformat{corollary}{#2Corollary~#1#3}

\newtheorem{theorem}[environment]{Theorem}
\crefname{theorem}{theorem}{Theorems}
\crefformat{theorem}{#2Theorem~#1#3}
\Crefformat{theorem}{#2Theorem~#1#3}

\newtheorem*{theorem*}{Theorem}
\crefname{theorem*}{theorem}{Theorems}
\crefformat{theorem*}{#2Theorem~#1#3}
\Crefformat{theorem*}{#2Theorem~#1#3}

\newtheorem{conjecture}[environment]{Conjecture}
\crefname{conjecture}{conjecture}{Conjectures}
\crefformat{conjecture}{#2Conjecture~#1#3}
\Crefformat{conjecture}{#2Conjecture~#1#3}

\newtheorem*{hypothesis*}{Hypothesis}
\crefname{hypothesis*}{conjecture}{Conjectures}
\crefformat{hypothesis*}{#2Conjecture~#1#3}
\Crefformat{hypothesis*}{#2Conjecture~#1#3}

\newtheorem{observation}[environment]{Observation}
\crefname{observation}{observation}{Observations}
\crefformat{observation}{#2Observation~#1#3}
\Crefformat{observation}{#2Observation~#1#3}

\crefname{example}{example}{examples}
\crefformat{example}{#2Example~#1#3}
\Crefformat{example}{#2Example~#1#3}

\newtheorem{remark}[environment]{Remark}
\crefname{remark}{remark}{remarks}
\crefformat{remark}{#2Remark~#1#3}
\Crefformat{remark}{#2Remark~#1#3}

\crefname{figure}{figure}{figures}
\crefformat{figure}{#2Figure~#1#3}
\Crefformat{figure}{#2Figure~#1#3}

\crefname{equation}{equation}{Equations}
\crefformat{equation}{#2#1#3}
\Crefformat{equation}{#2#1#3}

\crefname{chapter}{chapter}{chapters}
\crefformat{chapter}{#2Chapter~#1#3}
\Crefformat{chapter}{#2Chapter~#1#3}

\crefname{section}{section}{sections}
\crefformat{section}{#2Section~#1#3}
\Crefformat{section}{#2Section~#1#3}

\crefname{algorithm}{algorithm}{algorithms}
\crefformat{algorithm}{#2Algorithm~#1#3}
\Crefformat{algorithm}{#2Algorithm~#1#3}

\crefname{notation}{notation}{notations}
\crefformat{notation}{#2Notation~#1#3}
\Crefformat{notation}{#2Notation~#1#3}

\crefname{question}{question}{questions}
\crefformat{question}{#2Question~#1#3}
\Crefformat{question}{#2Question~#1#3}

\newtheorem{problem}[environment]{Problem}
\crefname{problem}{problem}{problem}
\crefformat{problem}{#2problem~#1#3}
\Crefformat{problem}{#2problem~#1#3}

\newtheorem{claim}{Claim}[environment]
\crefname{claim}{claim}{claims}
\crefformat{claim}{#2Claim~#1#3}
\Crefformat{claim}{#2Claim~#1#3}

\newtheorem{definition}[environment]{Definition}
\crefname{definition}{definition}{definitions}
\crefformat{definition}{#2Definition~#1#3}
\Crefformat{definition}{#2Definition~#1#3}

\usetikzlibrary{calc}
\usetikzlibrary{fit}
\usetikzlibrary{decorations}
\usetikzlibrary{decorations.pathmorphing}
\usetikzlibrary{decorations.text}
\usetikzlibrary{shapes,hobby}

\tikzset{
	position/.style args={#1:#2 from #3}{
		at=($(#3)+(#1:#2)$)
	}
}

\tikzset{
  v:main/.style = {draw, circle, scale=0.8, thick,fill=black,inner sep=0.7mm},
  v:ghost/.style = {inner sep=0pt,scale=1},
  >={latex},
  e:marker/.style = {line width=8.5pt,line cap=round,opacity=0.35,color=DarkGoldenrod},
  e:main/.style = {line width=1pt},
}

\title{Polynomial bounds for the Graph Minor Structure Theorem}
\predate{}
\date{}
\postdate{}

\preauthor{}
\DeclareRobustCommand{\authorthing}{
	\begin{center}
		Maximilian Gorsky\thanks{Supported by the Institute for Basic Science (IBS-R029-C1).}~~\!\footnote{\href{mailto:m.gorsky@pm.me}{m.gorsky@pm.me}} \\
		{\small Discrete Mathematics Group, Institute for Basic Science (IBS), Daejeon, South Korea} \\
        \medskip
		Micha\l{} T.\ Seweryn\thanks{Supported by ERC-CZ project LL2328 of the Ministry of Education of Czech Republic.}~~\!\footnote{\href{mailto:michal.seweryn@matfyz.cuni.cz}{michal.seweryn@matfyz.cuni.cz}} \\
		{\small Charles University in Prague, Czech Republic} \\
		  \medskip
		Sebastian Wiederrecht\footnote{\href{mailto:sebastian.wiederrecht@gmail.com}{wiederrecht@kaist.ac.kr}} \\
		{\small KAIST, South Korea} \\
\end{center}}
\author{\authorthing}
\postauthor{}

\begin{document}
\maketitle

\begin{abstract}
The Graph Minor Structure Theorem, originally proven by Robertson and Seymour [JCTB, 2003], asserts that there exist functions $f_1, f_2 \colon \mathbb{N} \to \mathbb{N}$ such that for every non-planar graph $H$ with $t\coloneq |V(H)|$, every $H$-minor-free graph can be obtained via the clique-sum operation from graphs which embed into surfaces where $H$ does not embed after deleting at most $f_1(t)$ many vertices with up to at most $t^2-1$ many ``vortices'' which are of ``depth'' at most $f_2(t)$.
In the proof presented by Robertson and Seymour the functions $f_1$ and $f_2$ are non-constructive.
Kawarabayashi, Thomas, and Wollan [arXiv, 2020] found a new proof showing that $f_1(t), f_2(t) \in 2^{\mathbf{poly}(t)}$.
While believing that this bound was the best their methods could achieve, Kawarabayashi, Thomas, and Wollan conjectured that $f_1$ and $f_2$ can be improved to be polynomials.

In this paper we confirm their conjecture and prove that $f_1(t), f_2(t) \in \mathbf{O}(t^{2300})$.
Our proofs are fully constructive and yield a polynomial-time algorithm that either finds $H$ as a minor in a graph $G$ or produces a clique-sum decomposition for $G$ as above.
\end{abstract}
\let\sc\itshape
\thispagestyle{empty}

\newpage

\setcounter{page}{1}


\newpage
\thispagestyle{empty}
\tableofcontents
\thispagestyle{empty}
\newpage

\section{Introduction}\label{sec:intro}

A graph $H$ is a \emph{minor} of a graph $G$ if it can be obtained from a subgraph of $G$ by \textsl{contracting} edges.
There are many instances of graph classes where graphs excluding particular minors share a specific, restricted structure.
For example, $K_3$-minor-free graphs are forests, and by Wagner's Theorem \cite{Wagner1937Ueber}, graphs excluding $K_5$ and $K_{3,3}$ as minors are planar.
These two examples can be seen as special cases of two theorems that sit at the heart of the monumental Graph Minors Series by Robertson and Seymour (R\&S), namely the \emph{Grid Theorem} \cite{RobertsonS1986Grapha} and the \emph{Graph Minor Structure Theorem} (\textsl{GMST}) \cite{RobertsonS2003Grapha} (see \cite{RobertsonS1985Graph} for a summary of both results written by R\&S).
The former describes the structure of graphs excluding planar graphs as minors (such as \(K_3\)), and the latter explains the structure of graphs excluding non-planar graphs as minors.

These two theorems are widely regarded to be among the deepest results in graph theory and their importance for structural graph theory, graph algorithms, and parametrised complexity can hardly be overstated \cite{RobertsonS1984Grapha,DowneyF1995Parameterized,Mohar2001Graph,Lovasz2005Graph,KawarabayashiM2007Recent,Norin2015New,LokshtanovSZ2020Efficient,Wollan2022Explicit}.
There has been considerable effort in the community to simplify the proofs and to find explicit and close to optimal bounds for the involved functions.
These goals were posed as important challenges by Lov{\'a}sz \cite{Lovasz2005Graph}.
Both of these efforts have been highly successful for the Grid Theorem due to the works of many authors \cite{RobertsonST1994Quickly,LeafS2015Treewidth,ChekuriC2016Polynomial,KawarabayashiK2020Linear,ChuzhoyT2021Tighter}.
In the case of the GMST, Lov{\'a}sz writes in his survey \cite{Lovasz2005Graph} on the Graph Minors Series: ``\textsl{It would be quite important to have simpler proofs with more explicit bounds. Warning: many of us have tried, but only few successes can be reported.}''

So far this statement has remained true as only a single set of authors, namely Kawarabayashi, Thomas, and Wollan (KTW) \cite{KawarabayashiTW2021Quickly} has managed, through immense effort,\footnote{The project took the trio roughly a decade \cite{KawarabayashiW2011Simpler,KawarabayashiTW2018New,KawarabayashiTW2021Quickly,Wollan2024Personal} and the central paper \cite{KawarabayashiTW2021Quickly} is still in the process of being published at the time of writing this article.}
to obtain explicit bounds, which are of the form $2^{\mathbf{poly}(t)}$.
In their paper, KTW conjecture that polynomial bounds should be possible.
\smallskip

\textbf{Our result.}
We show that the functions involved in the GMST are \emph{polynomials} in the size of the excluded minor, confirming the conjecture of KTW.
Our proofs are fully constructive and lead to a polynomial-time algorithm.

In addition, we also confirm two conjectures of Wollan \cite{Wollan2022Explicit}:
The first conjecture is a further strengthening of the result of KTW \cite{KawarabayashiTW2021Quickly}, confirming polynomial bounds and tight bounds on the Euler-genus of the surfaces involved.
The second asks for polynomial bounds on a theorem of R\&S on minors of graphs with highly representative embeddings (Theorem 4.3 in \cite{RobertsonS1995Graphb}) for which previously no constructive bounds were known.
\smallskip

Both the Grid Theorem and the GMST act as rough descriptions of graphs excluding a fixed graph $H$ as a minor and a simple way to reformulate those theorems is through the use of so-called ``clique-sums''.
A graph $G$ is a \emph{clique-sum} of two graphs $G_1$ and $G_2$ if there exist cliques\footnote{A \emph{clique} is a set of pairwise adjacent vertices.} $X_i\subseteq V(G_i)$, $i\in\{ 1,2\}$, where $|X_1|=|X_2|$ and $G$ can be obtained from $G_1$ and $G_2$ by identifying the vertices of $X_1$ and $X_2$ and then possibly deleting some of the edges between the vertices of $X_1=X_2$.

Roughly speaking, the Grid Theorem says that there exists a function $\mathsf{grid}$ such that for every planar graph $H$, any graph excluding $H$ as a minor can be obtained via clique-sums from graphs of size\footnote{The \emph{size} of a graph $G$ is the number of its vertices which we denote by $|G|$.} at most $\mathsf{grid}(|H|)$, i.e.\@ the \textsl{treewidth} of $G$ is at most $\mathsf{grid}(|H|)$.
The GMST on the other side says that there exists a function $\mathsf{clique}$ such that excluding a \textsl{non-planar} graph $H$ results in a class of graphs which can be obtained by means of clique-sums from graphs that are \textsl{$\mathsf{clique}(|H|)$-near embeddable} into surfaces\footnote{A \emph{surface} is a $2$-dimensional manifold, with or without boundary.} of small Euler-genus.\footnote{The \emph{Euler-genus} of a surface $\Sigma$ is the number $g=2h+c$ such that $\Sigma$ is homeomorphic to the surface obtained by attaching $h$ handles and $c$ crosscaps to the sphere \cite{Dyck1888Beitraege,FrancisW1999Conways}.}
Here a graph $G$ is said to be \emph{$k$-near embeddable}\footnote{See \cref{subsec:localstructure} for additional explanations, \cref{fig:intro-example} for an illustration, and \cref{sec:localtoglobal} for a rigorous formalisation of near embeddability.} into a surface $\Sigma$ if it can be obtained from \textsl{(i)} a graph $G_0$ embedded in $\Sigma$ by \textsl{(ii)} attaching \emph{vortices} \(G_1, \ldots, G_{\ell}\), where $\ell\leq k$, to distinct faces, each obtained by clique-summing graphs of size at most \(k\) in a particular way, and \textsl{(iii)} adding a set $A$ of at most $k$ \emph{apex vertices} and any set of edges with at least one endpoint in $A$.

Our main result and thus the GMST reads as follows (see \cref{thm:GMST} for a precise statement).

\begin{theorem}\label{thm:mainthm_simplest}
There exists a function $\mathsf{clique}\colon\mathbb{N}\to\mathbb{N}$, with $\mathsf{clique}(t)\in\mathbf{O}(t^{2300})$, such that for every non-planar graph $H$, every $H$-minor-free graph $G$ can be obtained by means of clique-sums from graphs that are $\mathsf{clique}(|H|)$-near embeddable into surfaces where $H$ does not embed.
\end{theorem}

In order to prove \cref{thm:mainthm_simplest} we begin with a simplified proof for the so-called \textsl{Flat Wall Theorem}, also known as the ``Weak Structure Theorem'', which first appeared in \cite{RobertsonS1995Graph} as Theorem 9.8.
In addition to this, we obtain a new proof for the so-called \textsl{Local Structure Theorem} (LST), also known as the ``structure with respect to a tangle'', which appeared in its original form as Theorem 3.1 in \cite{RobertsonS2003Graph}.
Our version of the LST enjoys the same qualities, including the polynomial bounds, as our main theorem.
Indeed one may regard the new bounds for the LST as our true main result.

Both theorems, the Flat Wall and the LST, are much more technical in nature than the statement of \Cref{thm:mainthm_simplest}.
Nonetheless they are key tools for the design of many algorithms on $H$-minor-free graphs (see for example \cite{RobertsonS1995Graph,LokshtanovSZ2020Efficient,GolovachST2023ModelChecking}) as they give rise to the so-called \emph{Irrelevant Vertex Technique} \cite{RobertsonS2012Graph}.
Beyond algorithmic applications, even the proof of well-quasi-ordering for graph minors builds on a different kind of ``globalisation'' of the LST \cite{RobertsonS2003Graph} rather than the GMST itself.
R\&S point this out themselves in \cite{RobertsonS2003Grapha} where they refer to the GMST as a ``red herring'' and explain that all of their own future applications are based on the LST.

Similar to the situation regarding the Grid Theorem, the Flat Wall Theorem has enjoyed quite a lot of improvements over the time \cite{GiannopoulouT2013Optimizing,Chuzhoy2015Improved,KawarabayashiTW2018New,SauST2024More}.
The LST on the other hand is so closely related to the full GMST that the sole significant progress is again the one of KTW \cite{KawarabayashiTW2021Quickly}.
Our new proof of the LST establishes polynomial bounds on \textsl{all} involved functions.
This has immediate consequences on the parametric dependencies in the running times of most (if not all) algorithmic applications of this theorem.

\begin{figure}
	\begin{center}
        \input{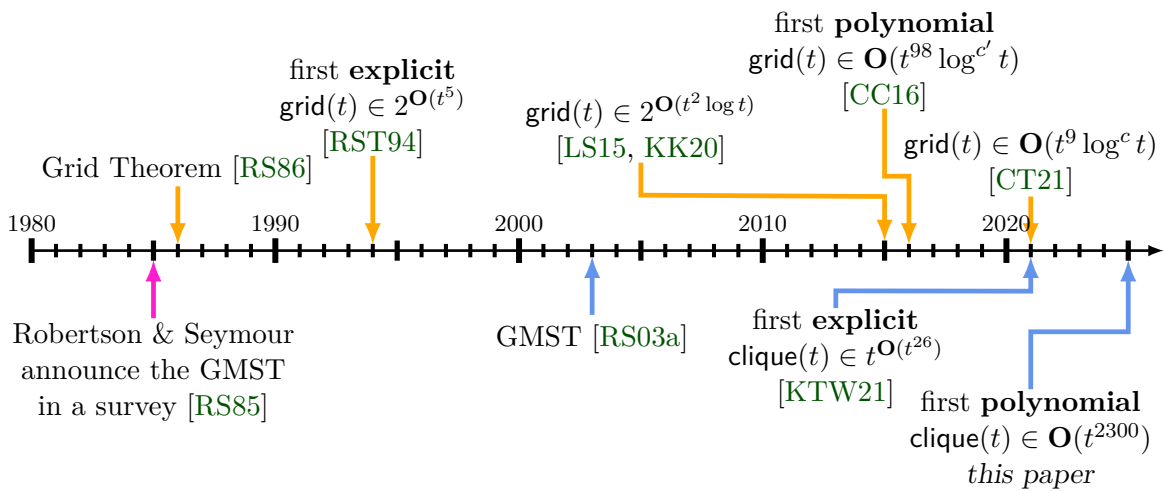}
	\end{center}
    \caption{A timeline on explicit bounds for R\&S' theory of graph minors.\vspace{-3mm}}
    \label{fig:timeline}
\end{figure}

\paragraph{A short history of explicit bounds for R\&S' theory of graph minors.}
The original proof for the Grid Theorem can be found in \cite{RobertsonS1986Grapha}, while the GMST was established in \cite{RobertsonS2003Grapha}, but requires the theory and many results proven in the preceding 15 papers of the series.
Both proofs only establish \textsl{existence} for the functions $\mathsf{grid}$ and $\mathsf{clique}$.
\Cref{fig:timeline} presents a timeline of the milestones regarding explicit bounds for the Grid Theorem and the GMST.
The first explicit bound for $\mathsf{grid}$ was proven by Robertson, Seymour, and Thomas \cite{RobertsonST1994Quickly}, establishing $\mathsf{grid}(t)\in 2^{\mathbf{poly}(t)}$.
Some follow-up work was concerned with improving the polynomial in the exponent \cite{LeafS2015Treewidth,KawarabayashiK2020Linear}.
In their breakthrough result, Chekuri and Chuzhoy \cite{ChekuriC2016Polynomial} achieved the next leap by proving $\mathsf{grid}(t)\in\mathbf{poly}(t)$.
This has since been improved by Chuzhoy and Tan \cite{ChuzhoyT2021Tighter} to $\mathbf{O}(t^9\operatorname{polylog} t)$.

The original statement of the GMST stems from R\&S' 16th entry in the Graph Minors Series \cite{RobertsonS2003Graph} and its proof contains a single non-constructive argument towards the existence of $\mathsf{clique}$.
This argument was made constructive by Geelen, Huynh, and Richter \cite{GeelenHR2018Explicit} without making the bound for $\mathsf{clique}$ explicit.
In their milestone papers, KTW \cite{KawarabayashiTW2018New,KawarabayashiTW2021Quickly} provide\footnote{We count \cite{KawarabayashiTW2018New} where KTW give a simple proof for the Flat Wall Theorem as part of their proof of the GMST. In \cref{sec:flattening} we prove our own variant of the Flat Wall Theorem.} a proof of the GMST giving a first explicit bound with $\mathsf{clique}(t)\in 2^{\mathbf{poly}(t)}$.
With their two papers, KTW also achieved the feat of greatly simplifying the proof of the GMST.
However, the bound on the Euler-genus of the surfaces into which their graphs may near embed is slightly worse than the optimal bound on the Euler-genus of the surfaces ``where $H$ does not embed''.

In \cite{KawarabayashiTW2021Quickly} the authors conjecture that the bounds for the GMST should be polynomials.
Moreover, Wollan conjectured \cite{Wollan2022Explicit} that in addition to polynomial bounds for the GMST, one can simultaneously achieve the optimal bound on the Euler-genus.
This would then further imply a third conjecture,
which asserts that there exists a polynomial $\mathsf{p}$ such that for every graph $H$, any graph $G$ that can be embedded with representativity.\footnote{Any non-contractible closed curve in $\Sigma$ intersecting $G$'s embedding only in vertices contains at least $\mathsf{p}(|H|)$ vertices.} at least $\mathsf{p}(|H|)$ into some surface $\Sigma$ where $H$ also embeds must contain $H$ as a minor \cite{Wollan2022Explicit}
In this paper, we resolve all three conjectures.

\paragraph{Applications and consequences of the GMST.}

Given the vast amount of applications of R\&S' Graph Minor Series we can only provide a glimpse at the mountains
of work that have since been built on this foundation.
For this reason, we present a few selected examples for applications of the GMST and immediate consequences of \cref{thm:mainthm_simplest} for them.

Without a doubt, the main algorithmic consequence is the Graph Minors Algorithm by R\&S \cite{RobertsonS1995Graph}.
R\&S used the LST to prove a bound on the so called
\emph{linkage function} $\mathsf{linkage}$ \cite{RobertsonS2009Graph,RobertsonS2012Graph} that can then be used to show that, given a graph $G$ together with vertices $s_1,\dots,s_k$, $t_1,\dots,t_k$, if $G$ has treewidth larger than $\mathsf{linkage}(k)$, then there exists a vertex $v\in V(G)$ such that there exist disjoint paths $P_1,\dots P_k$ in $G-v$ where $P_i$ has endpoints $s_i$ and $t_i$, if and only if such paths exist in $G$ \cite{RobertsonS1995Graph,RobertsonS2012Graph}.
This result implies a polynomial-time algorithm to test if a graph $G$ contains a graph $H$ as a minor \cite{RobertsonS1995Graph,KawarabayashiW2010Shorter,KorhonenPS2024Minora}.

We next highlight a structural contribution to the modern theory of graph minors in form of
a recent development regarding the \emph{Erd\H{o}s-P{\'o}sa property}.\footnote{We say that a graph $H$ has the \emph{Erd\H{o}s-P{\'o}sa property} if there exists a function $f_H\colon\mathbb{N}\to\mathbb{N}$ such that for every $k\in\mathbb{N}$ and every graph $G$, either (i) $G$ has $k$ subgraphs which are vortex-disjoint, each with an $H$-minor, or (ii) there is a set $S\subseteq V(G)$, $|S|\leq f_H(k)$, such that $G-S$ is $H$-minor-free. The \emph{half-integral Erd\H{o}s-P\'osa property} is defined by replacing ``vertex-disjoint''-condition in (i) by ``any vertex of $G$ appears in at most two of these subgraphs''.}
As a consequence of the Grid Theorem, R\&S proved that a graph $H$ has the Erd\H{o}s-P{\'o}sa property if and only if $H$ is planar \cite{RobertsonS1986Grapha}.
More generally, Thomas conjectured that every graph has the half-integral Erd\H{o}s-P{\'o}sa property \cite{Kawarabayashi2007Half,KawarabayashiM2007Recent,Norin2015New}.
This conjecture was confirmed without explicit bounds by Liu \cite{Liu2022Packing}.
Paul, Protopapas, Thilikos, and Wiederrecht \cite{PaulPTW2024Obstructionsa} gave an independent proof of Thomas' Conjecture with explicit bounds and showed that for every non-planar graph $H$ there are only finitely many families of counterexamples to the Erd\H{o}s-P{\'o}sa property of $H$.

It is important to point out that the functions $f_H$ from \cite{Liu2022Packing,PaulPTW2024Obstructionsa}, heavily depend on the functions of the GMST.

Some applications use generalisations or specialisations of the GMST as a departure point.
One such generalisation of the GMST is the one of Grohe and Marx \cite{GroheM2012Structure} to so-called \emph{topological minors}.\footnote{A graph $H$ is said to be a \emph{topological minor} of a graph $G$ if $G$ contains a subgraph isomorphic to a graph that can be obtained from $H$ by replacing its edges with paths or arbitrary length.}
They prove that there exists a function $\mathsf{top}(t)$ such that, for every integer $t\geq 1$, every graph that excludes $K_t$ as a topological minor can be obtained via the clique-sum operation from $\mathsf{top}(t)$-near embeddable graphs and graphs that have at most $\mathsf{top}(t)$ vertices of degree more than $\mathsf{top}(t)$. 
They apply this theorem to develop an algorithm for \textsc{Graph Isomorphism} for graphs excluding a graph $H$ as a topological minor. Our result implies $\mathsf{top}(t)\in\mathbf{poly}(t)$.

An intriguing specialisation of the GMST can be found in the work of Aprile, Fiorini, Joret, Kober, Seweryn, Weltge, and Yuditsky \cite{FioriniJWY2022Integer,AprileFJKSWY2025Integer}.
A deep open question in the field of integer programming asks whether there exists a polynomial time algorithm to solve integer programs whose coefficient matrices are \emph{totally $\Delta$-modular}.\footnote{A matrix $A$ is \emph{totally $\Delta$-modular}, $\Delta\in\mathbb{N}$, if for every square submatrix $A'$ of $A$ it holds that $-\Delta\leq \mathsf{Det}(A')\leq \Delta$ and $\mathsf{Det}(A)\in\mathbb{Z}$.}
The interplay with the GMST is hidden within a deeper connection to matroid minors \cite{GeelenGW2015Highly,GeelenNW2021Excluding}.
In \cite{FioriniJWY2022Integer,AprileFJKSWY2025Integer}, Aprile et al.\@ made use of this connection to provide polynomial-time algorithms for integer programming on important classes of totally $\Delta$-modular matrices.

Thilikos and Wiederrecht \cite{ThilikosW2024Killing} proved a dichotomy result for the $\#\mathsf{P}$-hard problem \cite{Valiant1979Complexity} of \textsc{Counting Perfect Matchings} on minor-closed graph classes.
They achieved this through a characterisation of all graphs $H$ where vortices may be omitted from the structure promised by the GMST for $H$-minor-free graphs.
Combining this structural result with a hardness result of Curticapean and Xia \cite{CurticapeanX2022Parameterizing} yields the dichotomy.

To contextualise this with our work, we stress that the only reason why $\mathsf{top}$ is not polynomial is its reliance on $\mathsf{clique}$.
Hence, \cref{thm:mainthm_simplest} implies that all bounds on the structure theorem of Grohe and Marx for excluding topological minors are polynomials.
Moreover, the algorithms from \cite{GroheM2012Structure,FioriniJWY2022Integer,AprileFJKSWY2025Integer,ThilikosW2024Killing} have running times of the form $n^{f(|H|)}$ where $H$ is the excluded (topological) minor.
In all cases, the function $f$ heavily depends on $\mathsf{clique}$ and thus, \cref{thm:mainthm_simplest}, together with the accompanying algorithm, implies drastic improvements on the running times of these algorithms.

Lastly, a broad field of applications can be found in alternative ways to decompose $H$-minor-free graphs.
Such decompositions usually aim to represent the $H$-minor-free input graph $G$ as a (family of) graph(s) of small treewidth by means of vertex-partitioning \cite{DeVosDOSRSV2004Excluding,Grohe2003Local,DemaineHK2005Algorithmic,DemaineHK2010Decomposition}, edge-contraction \cite{DemaineHK2011Contraction} or embedding \cite{Cohen-AddadLPP2023Planar}.
These types of approaches usually lead to constant-factor approximations and (quasi\footnote{In the case of \cite{Cohen-AddadLPP2023Planar}.}-)polynomial approximation schemes (PTAS) for various optimisation problems.
The range of applications for such techniques encompasses classical problems like \textsc{Chromatic Number} \cite{DemaineHK2005Algorithmic,DemaineHK2010Decomposition} for constant-factor approximation as well as PTASs for problems such as \textsc{Maximum Independent Set} and \textsc{Minimum Dominating Set} \cite{DemaineHK2010Decomposition}.
Moreover, routing problems such as \textsc{TSP} \cite{DemaineHK2011Contraction,Cohen-AddadLPP2023Planar}, as well as clustering, and network design problems \cite{Cohen-AddadLPP2023Planar} can be approached from this angle.
The role of the GMST in these applications lies in proving the existence of the required properties.
Often, the algorithms for finding these decompositions are based on the structure provided by the GMST.
Thus, \cref{thm:mainthm_simplest} can be readily applied to reduce the immense dependencies on $H$ in the running times of these algorithms.

\paragraph{Computing the GMST.}
Almost all algorithmic applications of the GMST require direct access to the clique-sum structure.
We  call the arising decomposition the \emph{structural decomposition} in the following.
R\&S' original proof for the GMST is mostly constructive but also contains a non-constructive argument.
For this reason it took considerable effort to make it algorithmic, a feat that was first achieved by Demaine, Hajiaghayi, and Kawarabayashi \cite{DemaineHK2005Algorithmic}.
The running time of their algorithm is estimated to be of the form $n^{f(|H|)}$ for some function $f$ where $H$ is the excluded minor.
Dawar, Grohe, and Kreutzer \cite{DawarGK2007Locally} gave a first algorithm with a running time of the form $f(|H|)n^{\mathbf{O}(1)}$ but for a much weaker variant of the structural decomposition.
In their paper, KTW \cite{KawarabayashiTW2021Quickly} claim a running time of the form $2^{t^{\mathbf{poly}(t)}}n^3$ to find the structural decomposition.
The best running time known at the time of writing this article is $f(|H|)n^2$ which was obtained by Grohe, Kawarabayashi, and Reed \cite{GroheKR2013Simple} through the use of irrelevant vertices and logic.
In light of recent progress in the algorithmic theory of graph minors \cite{GolovachST2023ModelChecking,KorhonenPS2024Minora}, Korhonen, Pilipczuk, and Stamoulis \cite{KorhonenPS2024Minora} hinted at the possibility for a running time of the form $f(|H|)n^{1+\mathbf{o}(1)}$ to compute the structural decomposition of $H$-minor-free graphs.

The logic-based approaches (\cite{GroheKR2013Simple} and the one hinted at by Korhonen et al.\@ in \cite{KorhonenPS2024Minora}) employ variants of Courcelle's Theorem \cite{Courcelle1990Monadic} which implies non-elementary dependencies hidden in the function $f$.
Moreover, the algorithm of KTW \cite{KawarabayashiTW2021Quickly} can only guarantee exponential bounds for the GMST.
It is also unclear if logic-based approaches can do any better than exponential since they rely on the Irrelevant Vertex Technique and therefore on the \emph{linkage function} \cite{RobertsonS2009Graph,KawarabayashiW2010Shorter} which is known to be at least exponential on planar graphs \cite{AdlerK2010Lower,AdlerKKLST2011Tight} and the best known general upper bounds are estimated to be at least quadruple exponential \cite{KawarabayashiW2010Shorter,Wollan2024Personal}.

\textbf{Our contribution.}
Our proof of \cref{thm:mainthm_simplest} is fully constructive and we obtain an algorithm that, given an $n$-vertex graph $G$ as input, either finds a witness that $G$ contains the graph $H$ as a minor or computes the structural decomposition for $H$-minor-free graphs within the bounds from \cref{thm:mainthm_simplest} for $G$ in time $2^{\mathbf{poly}(|H|)}n^{3}m\log n$.
Moreover, we provide a randomised algorithm that takes as input graphs $H$ and $G$ and finds either $H$ as a minor of $G$ or the structural decomposition for $H$-minor-free graphs for $G$ with high probability in time $(|H|+n)^{\mathbf{O}(1)}$.

\subsection{Sufficiency and necessity}\label{subsec:sufficiency}

Both the Grid Theorem and the GMST act as approximate \textsl{characterisations} for $H$-minor-free graphs.
However, so far we have only emphasised the \textsl{necessity} for the arising clique-sum structure of $H$-minor-free graphs.
A key aspect of both theorems is that these rough descriptions are also \textsl{sufficient} to capture the structure of $H$-minor-free graphs.
To make this precise, the Grid Theorem says that any graph with treewidth at least $\mathsf{grid}(t)$ must contain the \emph{$(t\times t)$-grid}\footnote{The $(t\times t)$-grid is the graph with vertex set $\{1,2,\dots,t\}^2$, $t\geq 3$, and where $(h,i)$ and $(j,k)$ are adjacent if and only if ($h=j$ and $|j-k|=1$) or ($j=k$ and $|h-j|=1$).} as a minor.
Now firstly, the $(t\times t)$-grid is a planar graph for all $t$ and for every planar graph $H$ there is $t_H\in\mathbf{O}(|H|)$ \cite{RobertsonST1994Quickly} such that $H$ is a minor of the $(t_H\times t_H)$-grid.
Secondly, and this is the sufficiency part, it is well known \cite{RobertsonS1986Grapha} that \textbf{no} graph that contains the $(t\times t)$-grid as a minor has treewidth less than $t$.

A similar statement holds for the GMST.
However, instead of the $(t\times t)$-grid as the universal graph we now have to consider the (non-)existence of a $K_t$-minor.
Joret and Wood \cite{JoretW2013Complete} proved that there exists a constant $c$ such that \textbf{no} graph with a $K_t$-minor can be obtained via clique-sums from graphs that are $k$-near embeddable into surfaces of Euler-genus $g$ for any $k$ with $t > c(k+1)\sqrt{g+k}$.

\paragraph{Tree structure and the three ingredients of near embeddings.}

\begin{figure}
    \centering
    \begin{tikzpicture}[scale=0.9]

        \pgfdeclarelayer{background}
		\pgfdeclarelayer{foreground}
			
		\pgfsetlayers{background,main,foreground}

        \begin{pgfonlayer}{background}
            \pgftext{\includegraphics[width=17cm]{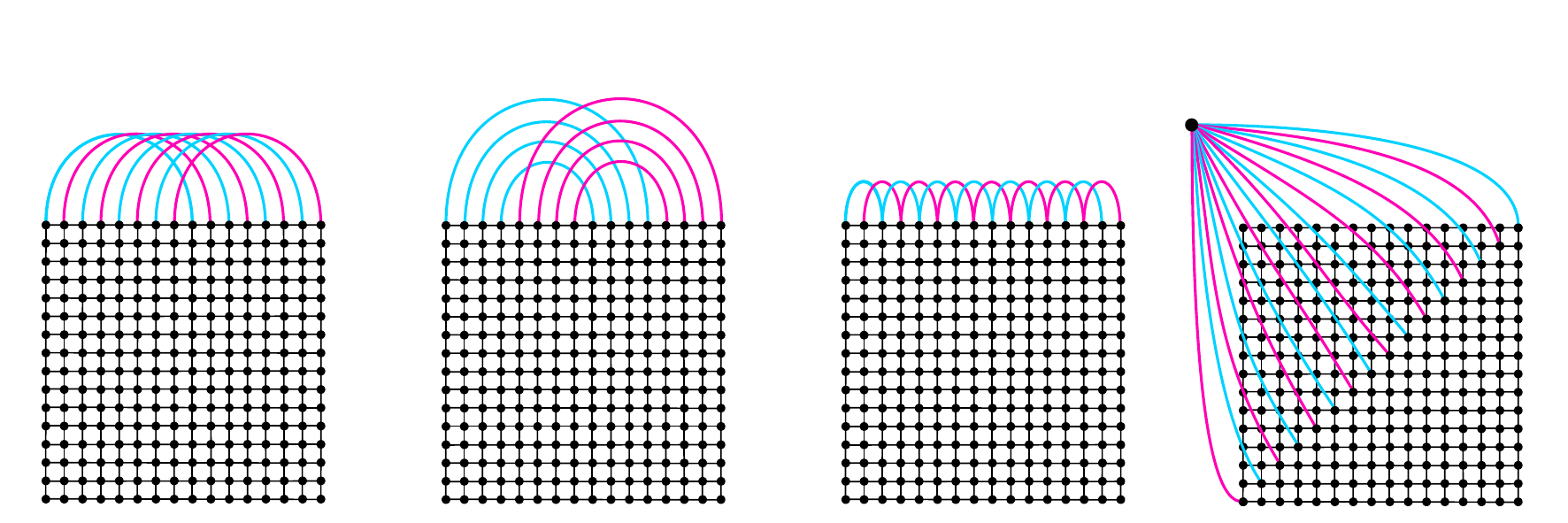}} at (C.center);
        \end{pgfonlayer}{background}
			
        \begin{pgfonlayer}{main}
        \node (C) [v:ghost] {};

        \node (M) [v:ghost,position=270:31.5mm from C] {};
        \node (i) [v:ghost,position=180:65mm from M] {\textsl{(i)}};
        \node (ii) [v:ghost,position=180:22.5mm from M] {\textsl{(ii)}};
        \node (iii) [v:ghost,position=0:22mm from M] {\textsl{(iii)}};
        \node (iv) [v:ghost,position=0:65mm from M] {\textsl{(iv)}};

        \end{pgfonlayer}{main}
        
        \begin{pgfonlayer}{foreground}
        \end{pgfonlayer}{foreground}

    \end{tikzpicture}
    \caption{Four examples illustrating that it is not possible to omit any of the three main ingredients of the GMST, namely (i) and (ii) Euler-genus, (iii) vortices, and (iv) apices.}
    \label{fig:three_ingredients}
\end{figure}
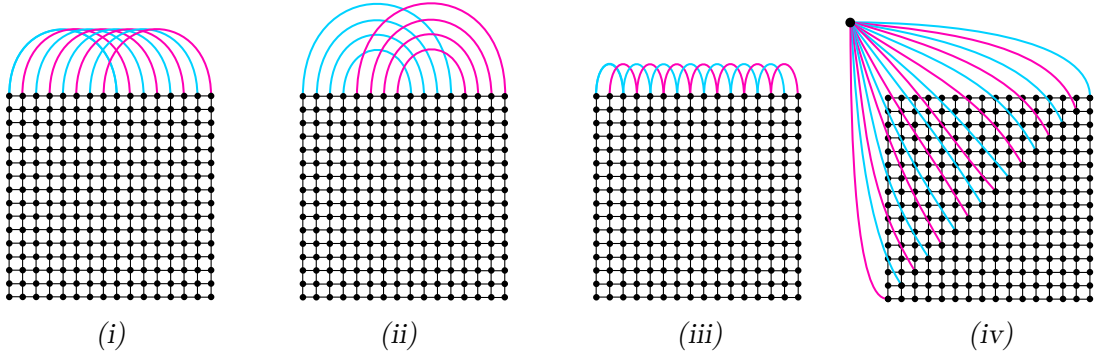

In the original paper for the GMST \cite{RobertsonS2003Graph} R\&S spent some time explaining why their three main ingredients are unavoidable for the notion of near embeddings to capture the structure of all $H$-minor-free graphs when $H$ is non-planar.
We present a short summary of their arguments.
The need for a tree-like structure of somewhat structured pieces was discovered early on.
Wagner \cite{Wagner1937Ueber} already proved that the $K_5$-minor-free graphs are exactly those graphs that can be obtained via clique-sums on cliques of size at most three from planar graphs and the Wagner graph.\footnote{The \emph{Wagner graph} is obtained from $C_8$ by adding edges between vertices at distance exactly four.}

In \cref{fig:three_ingredients} we present four graphs that can be scaled in straightforward ways to obtain examples of arbitrary sizes.
Thus, one can think of these graphs as representatives of infinite families with the following properties:
The family represented by \textsl{(i)} consists of grid-like graphs that embed in the projective plane while the family represented by \textsl{(ii)} consists of grid-like graphs that embed on the torus.
Such graphs are known to be $K_8$-minor-free.\footnote{For the projective plane this follows from Euler's Formula and for toroidal graphs, this is implied by the fact that they have chromatic number at most $7$ \cite{Heawood1890MapColour,Kauffman2011Seven}.}
The family of graphs in \textsl{(iii)} has been dubbed the ``shallow vortex grids'' and they are known to be $K_8$-minor-free \cite{SeeseW1989Grids}.
Finally, the graphs of the family represented by \textsl{(iv)} are known as \emph{apex grids}.
Since these graphs are one vertex-deletion away from being planar, they cannot contain a $K_6$-minor.
So each of these families excludes some $K_t$-minor for small $t$.
At the same time, all graphs from \textsl{(i)} and \textsl{(ii)} embed into surfaces of bounded Euler genus.
In contrast, the graphs from \textsl{(iii)} and \textsl{(iv)} can be seen to contain arbitrarily large $K_{3,k}$-minors which have arbitrarily large Euler-genus.
However, if we allow for a single vortex of depth $4$, the graphs from \textsl{(iii)} have an near embedding on the sphere and if we allow for an apex vertex, the graphs from \textsl{(iv)} also near embed on the sphere.

\subsection{The Local Structure Theorem}\label{subsec:localstructure}

Considering the discussion in \cref{subsec:sufficiency} and the Grid Theorem itself, it becomes clear that the only graphs where the notion of near embeddability becomes relevant are those of large treewidth, which is not visible in the statement of \cref{thm:mainthm_simplest}.
As discussed above, the centrepiece to achieving the GMST is the LST \cite{RobertsonS1995Graph,RobertsonS2003Grapha,RobertsonS2009Graph,RobertsonS2012Graph}.
Indeed, we believe that the LST is the heart of the entire graph minors theory of R\&S and its role is pivotal to most applications including the majority of those mentioned above.
Moreover, the LST explicitly captures the interesting parts of the GMST as those that exhibit large treewidth.

\begin{figure}
    \centering
    \begin{tikzpicture}[scale=0.9]

        \pgfdeclarelayer{background}
		\pgfdeclarelayer{foreground}
			
		\pgfsetlayers{background,main,foreground}

        \begin{pgfonlayer}{background}
            \pgftext{\includegraphics[width=12.5cm]{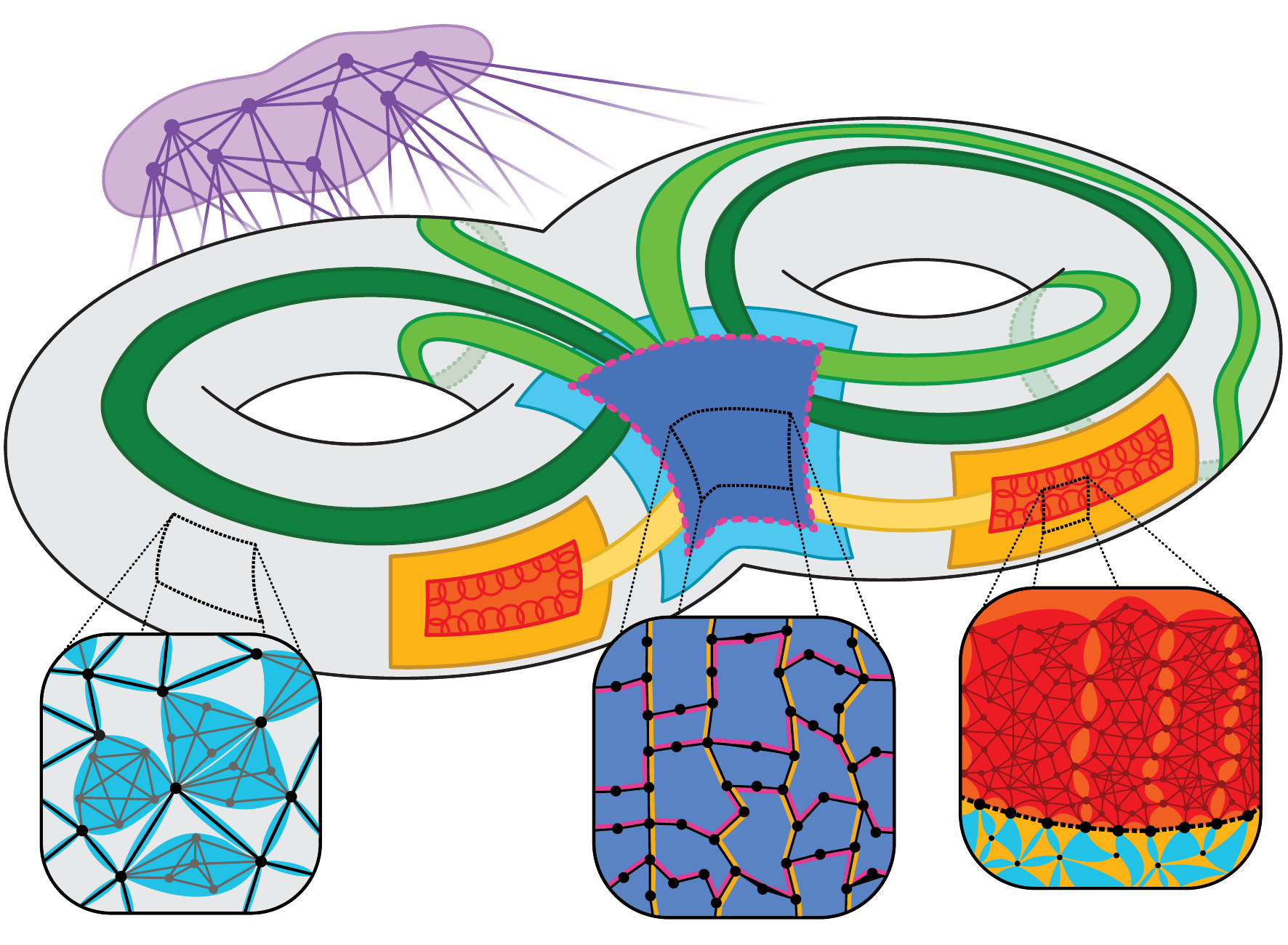}} at (C.center);
        \end{pgfonlayer}{background}
			
        \begin{pgfonlayer}{main}
        \node (C) [v:ghost] {};

        \node (apex) [v:ghost,position=135:57mm from C] {$A$};
        \node (wall) [v:ghost,position=45:13mm from C] {$M$};

        \node (ii) [v:ghost,position=247:24mm from C] {\textsl{(ii)}};
        \node (iii) [v:ghost,position=10.5:37.3mm from ii] {\textsl{(iii)}};
        \node (i) [v:ghost,position=178:52mm from ii] {\textsl{(i)}};

        \node (H1) [v:ghost,position=152:29.5mm from C] {$\mathcal{P}_1$};
        \node (H2) [v:ghost,position=60:25.5mm from C] {$\mathcal{P}_2$};

        \node (V1) [v:ghost,position=204:24.7mm from C] {$v_1$};
        \node (V2) [v:ghost,position=359:31.7mm from C] {$v_2$};

        \end{pgfonlayer}{main}
        
        \begin{pgfonlayer}{foreground}
        \end{pgfonlayer}{foreground}

    \end{tikzpicture}
    \caption{A $\Sigma$-decomposition $\delta$ of a graph $G-A$ with two vortices, namely $v_1$ and $v_2$, and an apex set $A$. Here $\Sigma$ is the double torus, $\delta$ is centred at the flat mesh $M$, (i) illustrates the cells of $\delta$, (ii) shows the two intertwining linkages that form the mesh $M$, and (iii) illustrates the inner, path-like, structure of a vortex.}
    \label{fig:intro-example}
\end{figure}
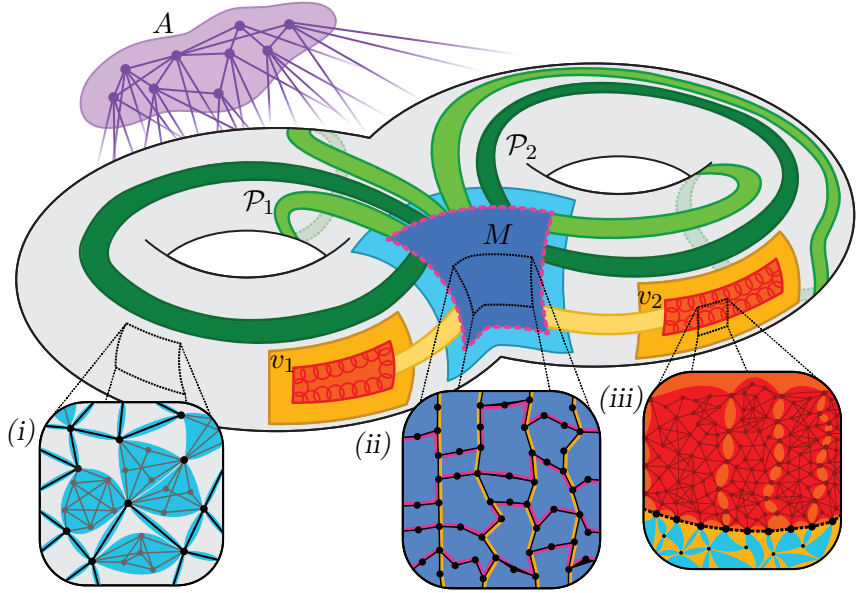

The LST describes the structure of graphs excluding a minor ``locally''. Informally, the theorem decomposes an \(H\)-minor-free ``highly linked part'' of a graph. In the original formulation of the LST, a ``highly-linked part'' is formalised by the concept of a ``tangle''.\footnote{See \cref{sec:tangles} for a definition of tangles.}
Tangles found their use not only in graphs, but also matroids, and more generally, artificial intelligence and empirical sciences (see the recent book by Diestel \cite{Diestel2024Tangles}).
The statement of the LST due to KTW is expressed in terms of walls instead of tangles. The benefit of using walls is that they are more concrete, and in particular they straightforwardly allow for efficient representation in computer memory.  
However, we found it to be more natural to work with an even more concrete object called a ``mesh''.
Hence, our local theorem describe the structure of a graph relative to a large mesh.\footnote{Every wall is a mesh and it is relatively straight forward to find a big wall in a large mesh, so the two notions are not to far apart and we prefer meshes purely for convenience. One can restate all of our results by replacing ``mesh'' with ``wall'' without changing the order of the functions.}

For an integer $k \geq 2$, a \emph{$k$-mesh} in a graph $G$ is the union of two families of pairwise disjoint paths $P_1,\dots, P_k$ and $Q_1,\dots,Q_k$, all in $G$, where the intersection of every $P_i$ with any $Q_j$ is a non-empty path, each $P_i$ intersects the paths $Q_j$ in order, and each $Q_i$ intersects the paths $P_j$ in order.

It is well-known that a graph $H$ is a minor of a graph $G$ if and only if there exists a function $\mu$ that maps the vertices of $H$ to connected subsets of $V(G)$ such that \textsl{(i)} if $u\neq v\in V(H)$ then $\mu(u)\cap\mu(v) = \emptyset$ and \textsl{(ii)} if $uv\in E(H)$ then there exists an edge $ab\in E(G)$ with $a\in \mu(u)$ and $b\in \mu(v)$.
We call $\mu$ a \emph{model} of $H$ in $G$ and the sets $\mu(u)$ the \emph{branch sets} of $\mu$.

Let $M$ be a $k$-mesh in a graph $G$, and let $P_1, \ldots, P_k, Q_1, \ldots, Q_k$ be as in the definition of a $k$-mesh.
Let $H$ be a graph with $|H| \leq k$, and let $\mu$ be a model of $H$ in $G$.
We say that $M$ \emph{controls} $\mu$ if for every set of vertices $X\subseteq V(G)$ with $|X| < |H|$, there exists a component of $G - X$ containing all branch sets of $\mu$ that are disjoint from $X$ and all
paths $P_i$ and $Q_j$ that are disjoint from $X$.

To fully understand the LST, we have to understand the notions of near embeddings and clique-sums from a different angle.
Let $\Sigma$ be a surface.
A \emph{$\Sigma$-decomposition}\footnote{We present here a simplified variant of $\Sigma$-decompositions for the sake of brevity. For a full definition see \cref{sec:decompositions}.} $\delta$ of \emph{breadth $b$} and \emph{depth $d$} of a graph $G$ is a collection of graphs $G_0,J_1,\dots,J_{\ell},G_1,\dots,G_b$ such that
\textsl{(i)} $G_0$ has an embedding in $\Sigma$ with $b$ \emph{(vortex) faces} bounded by the cycles $F_1,\dots,F_b$, \textsl{(ii)} each $G_i$ shares exactly the vertices of $F_i$ with $G_0$, \textsl{(iii)} for \textbf{no} $F_i$ there exist vertex-disjoint subpaths $P_1$ and $P_2$ of $F_i$ such that there are $d+1$ vertex-disjoint paths between $V(P_1)$ and $V(P_2)$ in $G_i$, and \textsl{(iv)} there exists a graph $G'$ that can be obtained via clique-sums of order at most three from $G_0$ and $J_1,\dots,J_{\ell}$ such that\footnote{We allow $J_i$ to be isomorphic to $K_2$. Hence, we may assume that every edge of $G_0$ appears in at least one of the $J_i$.} $G$ can be obtained from $G'\cup G_1 \cup \dots \cup G_b$ by deleting some of the edges of the cycles $F_1,\dots,F_b$.
We call the graphs $J_1,\dots,J_{\ell}$ the \emph{cells} of $\delta$ and the graphs $G_1,\dots,G_b$ the \emph{vortices} of $\delta$.
A $k$-mesh $M\subseteq G$ is said to be \emph{grounded} in $\delta$ if every cycle of $M$ uses edges from at least two different cells of $\delta$ and none of its edges belong to a vortex.
See \cref{fig:intro-example} for an illustration of a $\Sigma$-decomposition, with several additional features.
Most importantly, \cref{fig:intro-example} depicts a $\Sigma$-decomposition of a graph $G-A$ where $A$ is an apex set that has been deleted.
Moreover, there exist two pairs of linkages, $\mathcal{P}_1$ and $\mathcal{P}_2$ which force $G_0$ to be embedded in the double torus (or a surface of even higher Euler genus). 

With these definitions, our version of the LST reads as follows.

\begin{theorem}\label{thm:intro_local_structure}
There exist functions $\mathsf{local}\colon\mathbb{N}\to\mathbb{N}$ and $\mathsf{mesh}\colon\mathbb{N}^2\to\mathbb{N}$, with $\mathsf{local}(t)\in\mathbf{poly}(t)$ and $\mathsf{mesh(t,r)}\in\mathbf{poly}(t)\cdot r$, such that for all graphs $H$ with $|H|=t$ and $G$ with a $\mathsf{mesh}(t,r)$-mesh $M\subseteq G$ one of the following holds.
\begin{enumerate}
    \item there exists a model of $H$ in $G$ that is controlled by $M$, or
    \item there exists a surface $\Sigma$ into which $H$ does not embed, a set $A\subseteq V(G)$, $|A|\leq\mathsf{local}(t)$, and an $r$-mesh $M'\subseteq M$ such that $G-A$ has a $\Sigma$-decomposition $\delta$ of breadth at most $\nicefrac{1}{2}(t-3)(t-4)$, depth at most $\mathsf{local}(t)$, and where $M'$ is grounded.
\end{enumerate}
Moreover, there exists an algorithm that, given $G$ and $M$ as input, finds one of these outcomes in time $\mathbf{poly}(t)mn^2$. 
\end{theorem}

We believe \cref{thm:intro_local_structure} to be the key towards major improvements of the known bounds for $\mathsf{linkage}$, and note that the running time for \cref{thm:intro_local_structure} is polynomial in \textsl{both} $t$ and $n$.

\subsection{A short overview of our proof}\label{subsec:overview}

In the following, we first present a short high-level overview of our proof.
Then we provide some discussion on the key differences between our proof and the one of KTW.
We elaborate further on the key technique towards polynomial bounds in the next subsection.

We call a pair $(G,\Omega)$ where $G$ is a graph and $\Omega$ is a cyclic ordering of a subset of $V(G)$ a \emph{society}.
A \emph{transaction}\footnote{See the corresponding subsections of \cref{sec:decompositions} for a complete introduction to societies and their transactions.} in a society is a family of vertex-disjoint paths connecting disjoint segments of $\Omega$.
The proof for the GMST by R\&S was based on the study of societies and their transactions. 
KTW then proved, as their key lemma, a classification of societies.

\begin{wrapfigure}{l}{0.3\textwidth}
\centering
    \vspace{2mm}
    \scalebox{0.96}{
    \begin{tikzpicture}[scale=0.95]

        \pgfdeclarelayer{background}
		\pgfdeclarelayer{foreground}
			
		\pgfsetlayers{background,main,foreground}

        \begin{pgfonlayer}{background}
        \end{pgfonlayer}{background}
			
        \begin{pgfonlayer}{main}
        \node (x_0) [v:ghost] {};
        \node (main_1) [v:ghost,position=0:5mm from x_0] {\footnotesize graph $G$};
        \node (x_1) [v:ghost,position=270:16mm from x_0] {};
        
        \node (main_2) [v:ghost,position=270:3mm from main_1.center] {\footnotesize \& large mesh $M$};
        \node (x_2) [v:ghost,position=270:2mm from x_1.center] {};
        \node (main_4a) [v:ghost,position=0:5mm from x_2.center, align=center] {\footnotesize few apices};
        \node (main_4b) [v:ghost,position=270:3mm from main_4a.center, align=center] {\footnotesize \& flat mesh};
        \node (x_3) [v:ghost,position=270:12mm from x_2.center] {};
        \node (main_5a) [v:ghost,position=0:5mm from x_3.center, align=center] {\footnotesize society \& nest};
        \node (main_5b) [v:ghost,position=270:3mm from main_5a.center, align=center] {\footnotesize around a vortex};
        \node (main_5c) [v:ghost,position=270:7.5mm from main_5b.center, align=center] {\pgftext{\includegraphics[width=12mm]{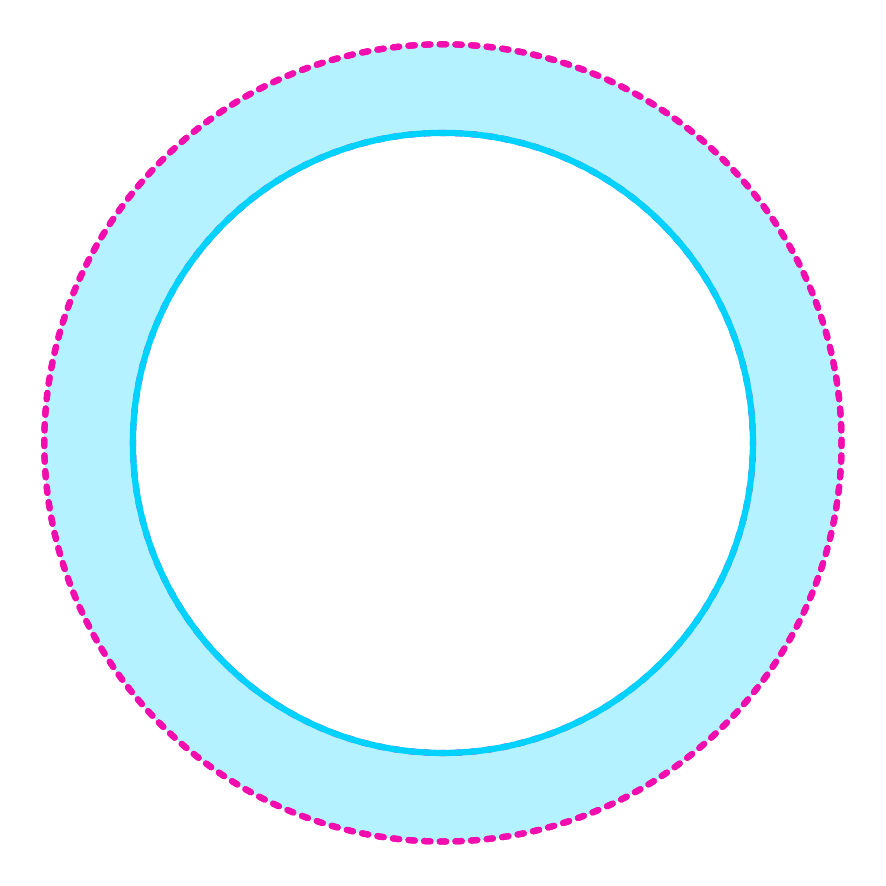}}};
        \node (x_4a) [v:ghost,position=270:16.5mm from x_3.center] {};
        \node (x_4) [v:ghost,position=270:9mm from x_4a.center] {};
        \node(main_6a) [v:ghost, position=180:5mm from x_4] {\footnotesize few apices};
        \node(main_6b) [v:ghost, position=270:3mm from main_6a] {\footnotesize \& crosscap};
        \node(x_5) [v:ghost, position=0:17mm from x_4] {};
        \node(main_7a) [v:ghost, position=270:3mm from x_5] {\footnotesize few apices};
        \node(main_7b) [v:ghost, position=270:3mm from main_7a] {\footnotesize \& flat transaction};
        \node(main_7c) [v:ghost, position=270:3mm from main_7b] {\footnotesize mesh of bounded};
        \node(main_7d) [v:ghost, position=270:3mm from main_7c] {\footnotesize depth};
        \node(main_7e) [v:ghost, position=270:7.5mm from main_7d] {\pgftext{\includegraphics[width=12mm]{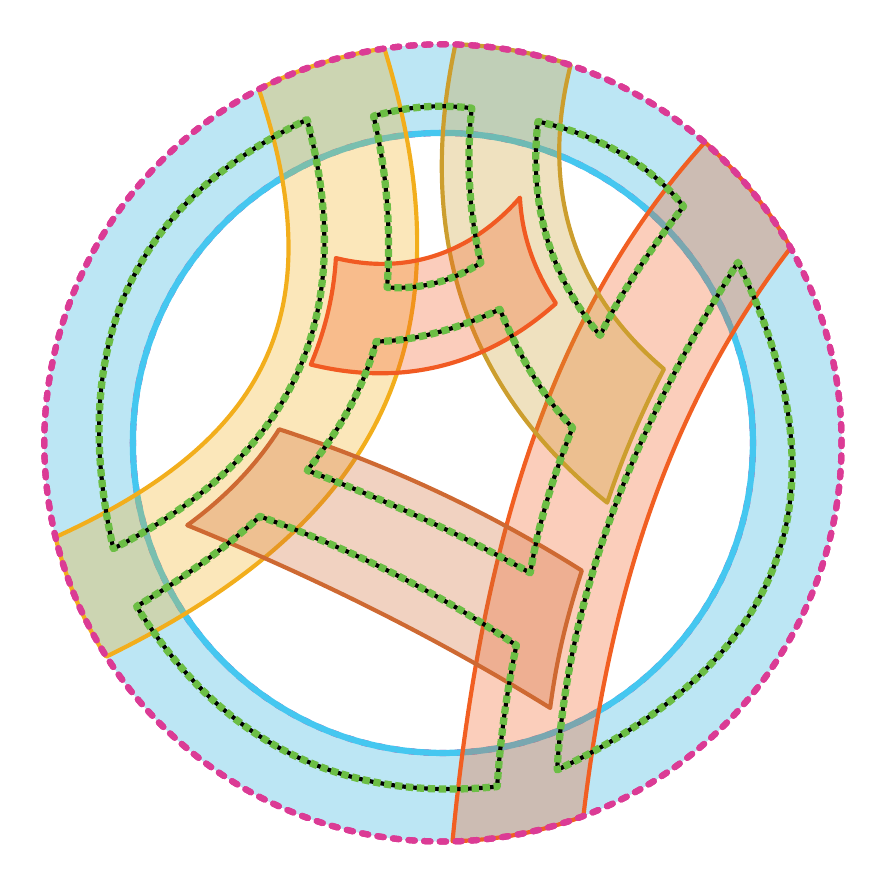}}};
        \node(x_6) [v:ghost,position=270:13mm from x_4.center] {};
        \node(main_8a) [v:ghost, position=180:5mm from x_6.center] {\footnotesize increase};
        \node(main_8b) [v:ghost, position=270:3mm from main_8a.center] {\footnotesize the genus};
        \node(main_8c) [v:ghost, position=270:7.5mm from main_8b.center] {\pgftext{\includegraphics[width=12mm]{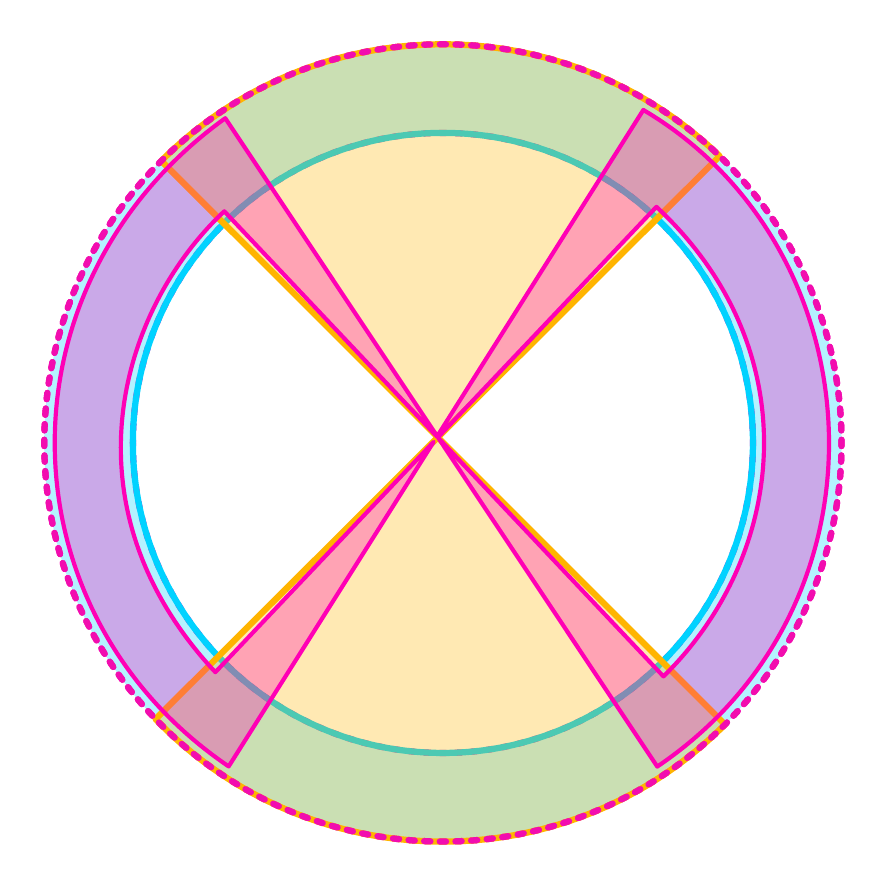}}};
        \node (x_7a) [v:ghost,position=270:12.5mm from x_6.center] {};
        \node (x_7) [v:ghost,position=270:13mm from x_7a.center] {};
        \node (main_9a) [v:ghost,position=0:5mm from x_7.center] {\footnotesize few apices \&};
        \node (main_9b) [v:ghost,position=270:3mm from main_9a.center] {\footnotesize all but a few disctrics};
        \node (main_9c) [v:ghost,position=270:3mm from main_9b.center] {\footnotesize are separated};
        \node (x_8a) [v:ghost,position=270:7.5mm from x_7] {};
        \node (x_8) [v:ghost,position=270:10mm from x_8a.center] {};
        \node (main_10a) [v:ghost,position=0:17mm from x_8.center] {\footnotesize few apices \&};
        \node (main_10b) [v:ghost,position=270:3mm from main_10a.center] {\footnotesize crosscap or};
        \node (main_10c) [v:ghost,position=270:2.5mm from main_10b.center] {\footnotesize handle};
        \node (x_9a) [v:ghost,position=270:8mm from x_8.center] {};
        \node (x_9) [v:ghost,position=270:8mm from x_9a.center] {};
        \node(main_11a) [v:ghost, position=0:17mm from x_9.center] {\footnotesize increase};
        \node(main_11b) [v:ghost, position=270:3mm from main_11a.center] {\footnotesize the genus};
        \node(main_11c) [v:ghost, position=270:7.5mm from main_11b.center] {\pgftext{\includegraphics[width=12mm]{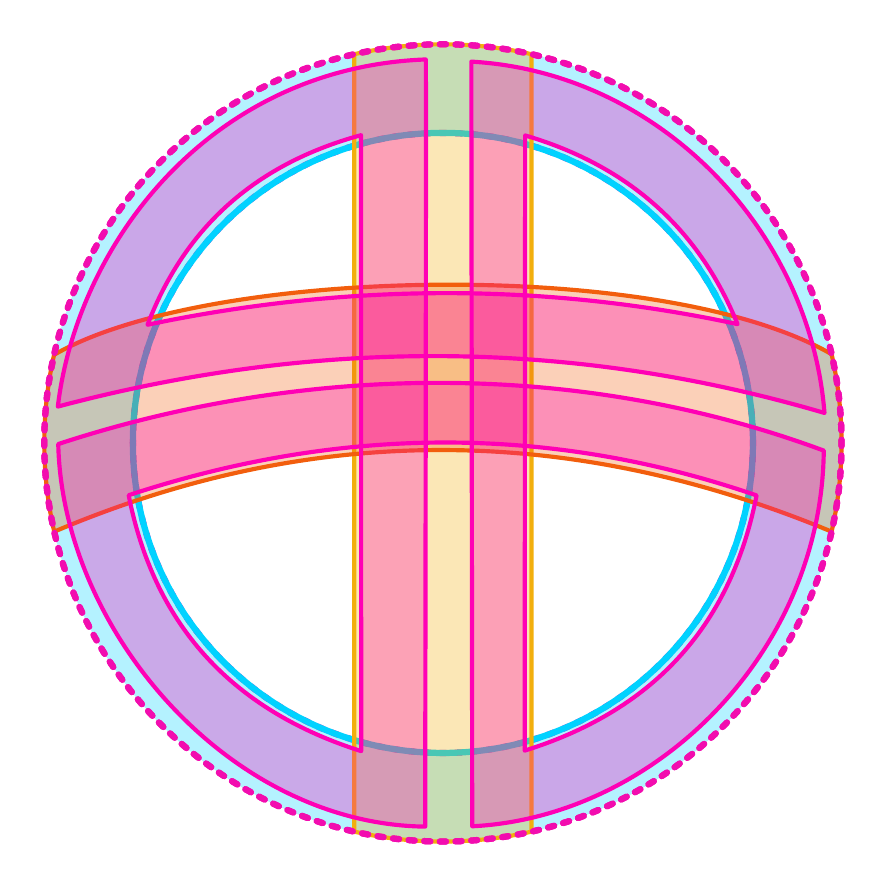}}};
        \node (x_10) [v:ghost,position=270:18mm from x_8a] {};
        \node (main_12a) [v:ghost,position=180:3mm from x_10] {\footnotesize few apices};
        \node (main_12b) [v:ghost,position=270:3mm from main_12a.center] {\footnotesize \& \textbf{all} districts};
        \node (main_12c) [v:ghost,position=270:3mm from main_12b.center] {\footnotesize are separated};
        \node (x_11a) [v:ghost,position=270:8mm from x_10.center] {};
        \node (x_11) [v:ghost,position=270:28mm from x_11a.center] {};
        \node (main_13) [v:ghost,position=180:12mm from x_11.center] {\footnotesize $K_t$-minor};
        \node (x_12) [v:ghost,position=270:10mm from x_11] {};
        \node (main_14) [v:ghost,position=0:6mm from x_12.center] {\footnotesize \cref{thm:intro_local_structure}};

        \node (y_2) [v:ghost,position=270:2mm from main_4b] {};

        \node (z_3a) [v:ghost,position=0:0mm from y_2] {};
        \node (z_3b) [v:ghost,position=90:1.5mm from main_5a] {};
        \draw [line width=1.2pt,->,color=DarkGray] (z_3a.center) to (z_3b.center);

        \node (y_3) [v:ghost,position=270:2mm from main_2] {};

        \node (z_4a) [v:ghost,position=0:0.5mm from y_3] {};
        \node (z_4b) [v:ghost,position=270:5mm from z_4a.center] {};
        \node (z_4c) [v:ghost,position=180:0.5mm from z_4b.center] {};
        \node (z_4d) [v:ghost,position=90:1.5mm from main_4a] {};
        \draw [line width=1.2pt,->,color=DarkGray] (z_4a.center) to (z_4b.center) to (z_4c.center) to (z_4d.center);

        \node (y_4) [v:ghost,position=270:7mm from main_5c] {};

        \node (z_5a) [v:ghost,position=0:0mm from y_4] {};
        \node (z_5b) [v:ghost,position=270:3.5mm from z_5a] {};
        \node (z_5c) [v:ghost,position=180:10mm from z_5b] {};
        \node (z_5d) [v:ghost,position=90:1.5mm from main_6a] {};
        \draw [line width=1.2pt,->,color=DarkGray] (z_5a.center) to (z_5b.center) to (z_5c.center) to (z_5d.center);

        \node (z_6a) [v:ghost,position=0:1mm from y_4] {};
        \node (z_6b) [v:ghost,position=270:3.5mm from z_6a] {};
        \node (z_6c) [v:ghost,position=0:11mm from z_6b] {};
        \node (z_6d) [v:ghost,position=90:1.5mm from main_7a] {};
        \draw [line width=1.2pt,->,color=DarkGray] (z_6a.center) to (z_6b.center) to (z_6c.center) to (z_6d.center);

        \node (y_5) [v:ghost,position=270:2mm from main_6b] {};

        \node (z_7a) [v:ghost,position=0:0mm from y_5] {};
        \node (z_7b) [v:ghost,position=90:1.5mm from main_8a] {};
        \draw [line width=1.2pt,->,color=DarkGray] (z_7a.center) to (z_7b.center);

        \node (y_6) [v:ghost,position=270:7mm from main_7e] {};

        \node (z_8a) [v:ghost,position=0:0mm from y_6] {};
        \node (z_8b) [v:ghost,position=270:2.5mm from z_8a] {};
        \node (z_8c) [v:ghost,position=180:12mm from z_8b] {};
        \node (z_8d) [v:ghost,position=90:1.5mm from main_9a] {};
        \draw [line width=1.2pt,->,color=DarkGray] (z_8a.center) to (z_8b.center) to (z_8c.center) to (z_8d.center);

        \node (y_7) [v:ghost,position=270:2mm from main_9c] {};

        \node (z_9a) [v:ghost,position=0:1mm from y_7] {};
        \node (z_9b) [v:ghost,position=270:4mm from z_9a] {};
        \node (z_9c) [v:ghost,position=0:11mm from z_9b] {};
        \node (z_9d) [v:ghost,position=90:1.5mm from main_10a] {};
        \draw [line width=1.2pt,->,color=DarkGray] (z_9a.center) to (z_9b.center) to (z_9c.center) to (z_9d.center);

        \node (z_10a) [v:ghost,position=0:0mm from y_7] {};
        \node (z_10b) [v:ghost,position=270:6mm from z_10a] {};
        \node (z_10c) [v:ghost,position=180:8mm from z_10b] {};
        \node (z_10d) [v:ghost,position=90:1.5mm from main_12a] {};
        \draw [line width=1.2pt,->,color=DarkGray] (z_10a.center) to (z_10b.center) to (z_10c.center) to (z_10d.center);

        \node (y_8) [v:ghost,position=270:2mm from main_10c] {};

        \node (z_11a) [v:ghost,position=0:0mm from y_8] {};
        \node (z_11b) [v:ghost,position=90:1.5mm from main_11a] {};
        \draw [line width=1.2pt,->,color=DarkGray] (z_11a.center) to (z_11b.center);

        \node (y_9) [v:ghost,position=270:2mm from main_12c] {};

        \node (z_12a) [v:ghost,position=0:0mm from y_9] {};
        \node (z_12b) [v:ghost,position=270:10mm from z_12a] {};
        \node (z_12c) [v:ghost,position=0:9mm from z_12b] {};

        \node (y_10) [v:ghost,position=270:7mm from main_8c] {};
        \node (y_11) [v:ghost,position=270:7mm from main_11c] {};

        \node (m_0) [v:ghost,position=90:2mm from main_13] {};
        \node (m_1) [v:ghost,position=90:6mm from m_0.center] {};
        \node (m_2) [v:ghost,position=0:28.5mm from m_1.center] {};
        \node (m_3) [v:ghost,position=180:0.5mm from y_11] {};

        \node (m_4) [v:ghost,position=90:6mm from main_14] {};
        \node (m_5) [v:ghost,position=180:18mm from m_4] {};
        \node (m_6) [v:ghost,position=270:2mm from main_13] {};

        \node (m_7) [v:ghost,position=180:6mm from m_1.center] {};
        \node (m_8) [v:ghost,position=90:41.5mm from m_7.center] {};
        \node (m_9) [v:ghost,position=0:22mm from m_8.center] {};
        \node (m_10) [v:ghost,position=180:1mm from y_7] {};

        \node (m_11) [v:ghost,position=90:17mm from m_8.center] {};
        \node (m_12) [v:ghost,position=0:13.5mm from m_11.center] {};
        \node (m_13) [v:ghost,position=0:0.5mm from y_10.center] {};

        \node (m_14) [v:ghost,position=90:39mm from m_11.center] {};
        \node (m_15) [v:ghost,position=0:22mm from m_14.center] {};
        \node (m_16) [v:ghost,position=180:1mm from y_4.center] {};

        \node (m_17) [v:ghost,position=90:41mm from m_14.center] {};
        \node (m_18) [v:ghost,position=0:22.5mm from m_17.center] {};
        \node (m_19) [v:ghost,position=180:0.5mm from y_3.center] {};

        \node (b_1) [v:ghost,position=0:0.5mm from y_11] {};
        \node (b_2) [v:ghost,position=270:2.5mm from b_1.center] {};
        \node (b_3) [v:ghost,position=0:13.5mm from b_2.center] {};
        \node (b_4) [v:ghost,position=90:117.5mm from b_3.center] {};
        \node (b_5) [v:ghost,position=0:11.5mm from main_5a.center] {};

        \draw [line width=1.2pt,->,color=CornflowerBlue] (b_1.center) to (b_2.center) to (b_3.center) to (b_4.center) to (b_5.center);

        \node (a_1) [v:ghost,position=180:0.5mm from y_10] {};
        \node (a_2) [v:ghost,position=270:2mm from a_1.center] {};
        \node (a_3) [v:ghost,position=180:10mm from a_2.center] {};
        \node (a_4) [v:ghost,position=90:58mm from a_3.center] {};
        \node (a_5) [v:ghost,position=180:11.5mm from main_5a.center] {};

        \draw [line width=1.2pt,->,color=CornflowerBlue] (a_1.center) to (a_2.center) to (a_3.center) to (a_4.center) to (a_5.center);
        
        \draw [line width=5pt,color=white] (m_14.center) to (m_15.center);
        
        \draw [line width=1.2pt,->,color=DarkGray] (m_3.center) to (m_2.center) to (m_1.center) to (m_0.center);
        \draw [line width=1.2pt,color=DarkGray] (m_1.center) to (m_7.center) to (m_8.center) to (m_9.center) to (m_10.center);
        \draw [line width=1.2pt,color=DarkGray] (m_8.center) to (m_11.center) to (m_12.center) to (m_13.center);
        \draw [line width=1.2pt,color=DarkGray] (m_11.center) to (m_14.center) to (m_15.center) to (m_16.center);
        \draw [line width=1.2pt,color=DarkGray] (m_14.center) to (m_17.center) to (m_18.center) to (m_19.center);
        
        \node (z_12d) [v:ghost,position=90:1.5mm from main_14] {};
        \draw [line width=5pt,color=white] (z_12a.center) to (z_12b.center) to (z_12c.center) to (z_12d.center);

        \draw [line width=1.2pt,color=DarkGray] (m_6.center) to (m_5.center) to (m_4.center);
        
        \draw [line width=1.2pt,->,color=DarkGray] (z_12a.center) to (z_12b.center) to (z_12c.center) to (z_12d.center);
          
        \end{pgfonlayer}{main}
        
        \begin{pgfonlayer}{foreground}
        \end{pgfonlayer}{foreground}

    \end{tikzpicture}}
\caption{The proof of \cref{thm:intro_local_structure} (simplified).}
\label{fig:intro_flowchart}
\end{wrapfigure}

In a generalisation of the so-called \textsl{Two Paths Theorem}\footnote{See \cref{prop:TwoPaths}.} \cite{Jung1970Verallgemeinerung,Seymour1980Disjoint,Shiloach1980Polynomial,Thomassen19802Linked} KTW proved (see also \cite{RobertsonS1990Graph} for a weaker but more fundamental variant) that any society equipped with a $\Delta$-decomposition\footnote{Here $\Delta$ is a disk, that is a special case of a surface with boundary.} with a single vortex and a set of concentric cycles around it, called the \emph{nest}, either contains \textsl{(i)} a $K_t$-minor, \textsl{(ii)} a ``crosscap-'' or ``handle-transaction'',\footnote{See for example the coloured linkages in (i) and (ii) of \cref{fig:three_ingredients} or the linkages $\mathcal{P}_1$ and $\mathcal{P}_1$ in \cref{fig:intro-example}.} or \textsl{(iii)} a small set $A$ such that $(G-A,\Omega)$ has a $\Delta$-decomposition of bounded depth and breadth.

A main source for the exponential growth of their bounds is that whenever KTW find a transaction that ``splits'' the vortex, they either find a handle or crosscap, or they separate both vortices by adding more vertices to the apex set.
Please note that this is only one of three sources of exponential growth in the proof of KTW.
We discuss the other two sources further below.

Separating the two newly born vortices after each ``split'' requires the sacrifice of roughly half of the initial nest in each iteration.
To avoid this problem, we do not ``separate'' our vortices from each other.
Instead, we build a meta-structure, called a \emph{transaction mesh},\footnote{See \cref{sec:transactionmeshes}.} by iteratively adding as many transactions as possible without increasing the genus.
Some of these transactions may carry local witnesses of non-planarity hidden within a private apex set.
However, if there are too many such transaction we are able to construct a $K_t$-minor and thus most of these transactions can be added without contributing to the apex set at all.
While this approach generates an unbounded number of ``vortices'' (or better \textsl{districts}), by extending arguments used in the modern proofs for the Flat Wall Theorem \cite{GiannopoulouT2013Optimizing,Chuzhoy2015Improved,KawarabayashiTW2018New}, namely variations of Gallai's \textsl{$A$-Paths Theorem} \cite{Gallai1964MaximumMinimum,BruhnHJ2018Frames}, we show that only a small number of these districts are potential vortices as otherwise $G$ contains a $K_t$-minor.
It then suffices to separate these potential vortices, or find a handle or crosscap instead, in the very last step of our procedure.
As a result, this technique alone is able to avoid most of the telescoping effects that are responsible for the exponentiality in the proof of \cite{KawarabayashiTW2021Quickly}.
See \cref{fig:intro_flowchart} for a flow chart depicting the general strategy of our proof.
We show that the genus-increasing step can only occur a bounded number of times by iteratively constructing a minor-universal graph for graphs embedded in surfaces derived from a recent discovery by Gavoille and Hilaire \cite{GavoilleH2023MinorUniversal}.
This last part allows us to improve over the suboptimal bound on the Euler-genus proven by KTW and recover the original statement of R\&S while still guaranteeing overall polynomial bounds.

\paragraph{Crooked transactions, leap patterns and Gallai's $A$-paths theorem.}

We pause for a second to give a warning to the reader.
The outline of our proof given above is already very rough and only vaguely captures the countless technicalities of our proof.
In what follows we try to shed some additional light on our approach, the challenges we encounter, and how we overcome them.
However, this is to be treated cautiously as we make use of, sometimes gross, simplifications and are only able to convey so much information in only a few pages.
The following discussion is meant to instil a small amount of intuition and the reader is advised not to dwell to long over a seemingly mysterious paragraph.
All mysteries are fully unravelled in the main body of the paper.
In \cref{subsec:organisation} we provide a comprehensive overview on where to find which part of the proof, so anyone who wants to see a certain point explained in full detail can quickly find their way.

Before we discuss the last, and probably deepest reason for the exponentiality of the bounds obtained by KTW, let us unravel the idea of transaction meshes a bit more.
As described above, the proof of KTW either finds a way to increase the Euler-genus of the current surface by adding a crosscap or handle transaction, or it splits the current vortex into two and then separates the two sides.
 
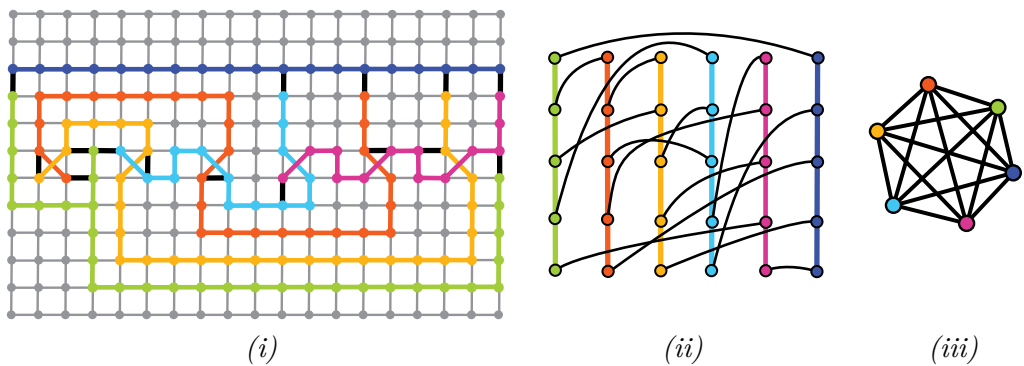
\begin{figure}[ht]
    \centering
    \begin{tikzpicture}

        \pgfdeclarelayer{background}
		\pgfdeclarelayer{foreground}
			
		\pgfsetlayers{background,main,foreground}

        \begin{pgfonlayer}{background}
        \node (C) [v:ghost] {{\includegraphics[width=15cm]{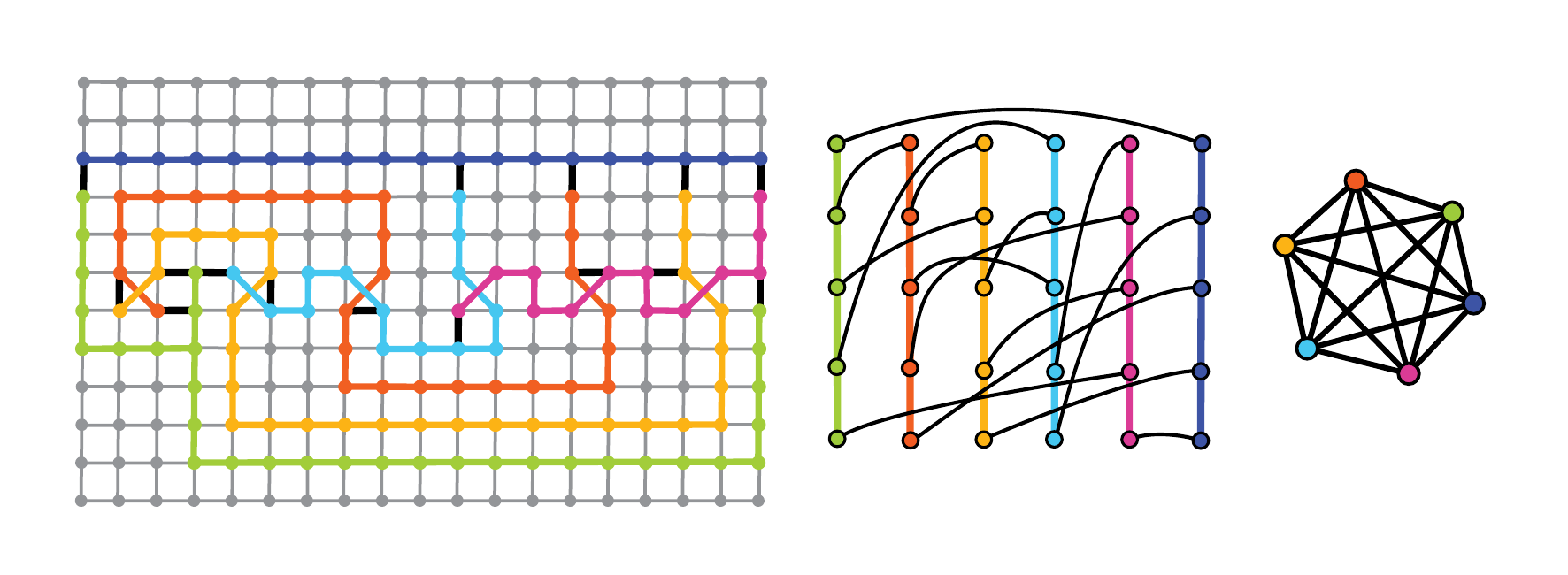}}};
        \end{pgfonlayer}{background}
			
        \begin{pgfonlayer}{main}

        \node (X) [v:ghost,position=270:25mm from C.center] {};
        
        \node (C_label) [v:ghost,position=0:22mm from X.center] {\textit{(ii)}};
        \node (L_label) [v:ghost,position=180:34mm from X.center] {\textit{(i)}};
        \node (R_label) [v:ghost,position=0:58mm from X.center] {\textit{(iii)}};
        
        \end{pgfonlayer}{main}

        \begin{pgfonlayer}{foreground}
        \end{pgfonlayer}{foreground}

    \end{tikzpicture}
    \caption{\textsl{(i)} A mesh with some crosses spread out over its middle row together with branchsets of a $K_6$-minor, \textsl{(ii)} a diagram of the branchsets of the $K_6$-minor, and \textsl{(iii)} the final $K_6$ obtained by contracting each of the branchsets into a single vertex.}
    \label{fig:K6inGrid}
\end{figure}

A major part of the analysis boils down to the concept of ``crosses'' on a society.
Given a society $(G,\Omega)$, a \emph{cross} on $(G,\Omega)$ is a pair of disjoint paths $P_1$ and $P_2$ each with endpoints $s_i$, $t_i$, $i\in[2]$, such that the four vertices $s_1, s_2, t_2, t_2$ occur in $\Omega$ in the order listed.
A nice preliminary observation, also providing some intuition for the proof of the Flat Wall Theorem, is that many crosses spread out over a large mesh give rise to a clique minor as depicted in \cref{fig:K6inGrid} (see \cite{KawarabayashiTW2018New} or \cref{sec:forcing}).

By following the strategy of KTW as described above, it might happen that after the separation step, one of the two sides ceases to be a vortex, i.e.\@ one of the two sides could be \textsl{essentially} planar.
This poses a problem because a priori the only ways to construct a clique-minor are either finding one through an application of the Flat Wall Theorem, by increasing the Euler-genus of the surface too much, or by having too many vortices, each with a private cross in its interior.
Notice that the third outcome would again yield a situation akin to the one in \cref{fig:K6inGrid}.
If however, in each step we simply find a flat transaction that does not allow us to increase the Euler-genus and then delete a few vertices to separate both sides, only to end up with one side being what is called \emph{flat}, i.e.\@ without a private cross, in the process, we do not work towards any of the three goals above in any way.
This is a fundamental problem that R\&S were already faced with.

To overcome this issue, R\&S introduced the concept of \emph{crooked transactions}.
Up to some small linear dependency on the order of the transaction, one may classify crooked transactions to be exactly those transactions that fall into one of the two following types (see \cite{RobertsonS1990Graph} and \Cref{subsec:normalisingcrooked}):\\
\noindent \textsl{(i)} \emph{leaps}; transactions containing a path that jumps all the way from one side of the transaction to the other, thereby forming a cross with all other paths involved in the transaction, and\\
\noindent \textsl{(ii)} \emph{doublecrosses}; transactions in which the two outermost paths on either side form a cross.

\begin{figure}[ht]
    \centering
    \begin{tikzpicture}

        \pgfdeclarelayer{background}
		\pgfdeclarelayer{foreground}
			
		\pgfsetlayers{background,main,foreground}

        \begin{pgfonlayer}{background}
        \node (C) [v:ghost] {{\includegraphics[width=15cm]{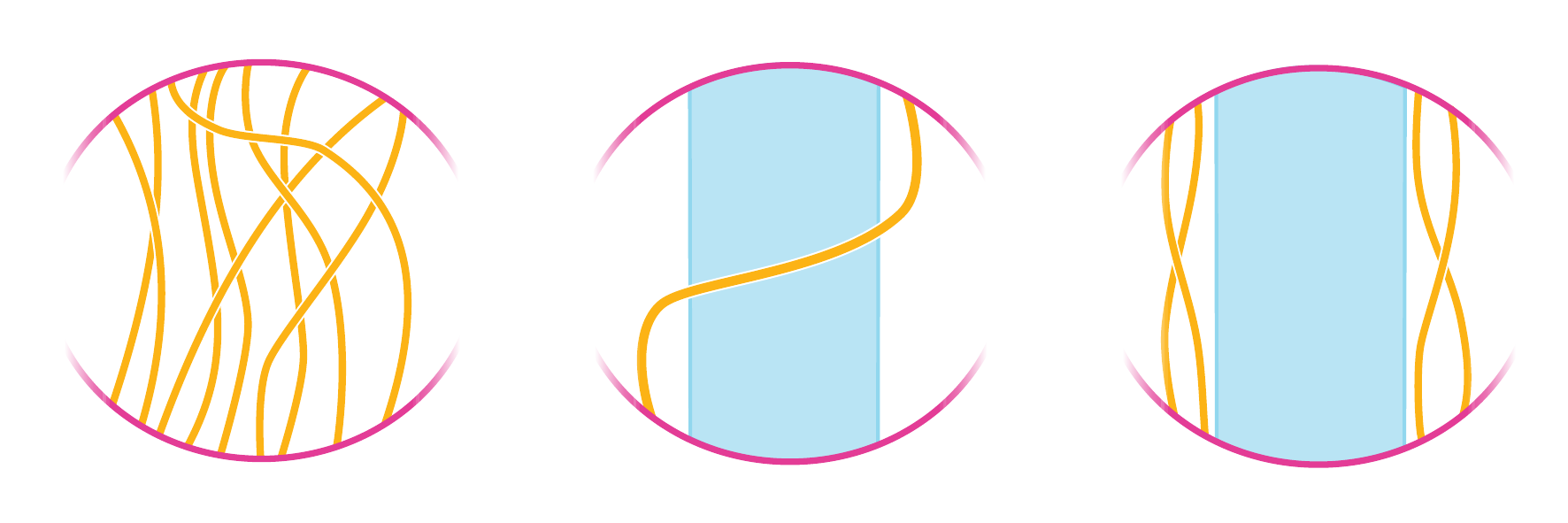}}};
        \end{pgfonlayer}{background}
			
        \begin{pgfonlayer}{main}

        \node (X) [v:ghost,position=270:25mm from C.center] {};
        
        \node (C_label) [v:ghost,position=0:0mm from X.center] {\textit{(ii)}};
        \node (L_label) [v:ghost,position=180:50mm from X.center] {\textit{(i)}};
        \node (R_label) [v:ghost,position=0:50mm from X.center] {\textit{(iii)}};
        
        \end{pgfonlayer}{main}

        \begin{pgfonlayer}{foreground}
        \end{pgfonlayer}{foreground}

    \end{tikzpicture}
    \caption{\textsl{(i)} A crooked transaction in a society (see \cref{sec:findcrooked} for a definition), \textsl{(ii)} diagram of a crooked transaction of type \textsl{leap}, and \textsl{(iii)} a diagram of a crooked transaction of type \textsl{doublecross}.}
    \label{fig:crookedIntro}
\end{figure}
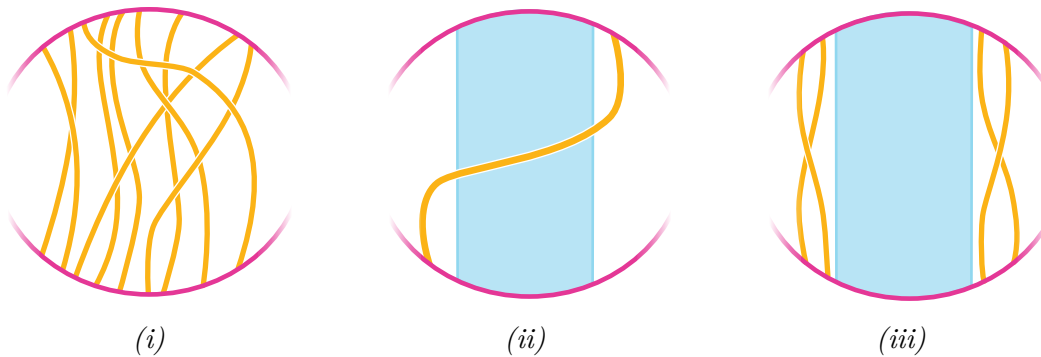

See \cref{fig:crookedIntro} for an illustration.

\begin{figure}[ht]
    \centering
    \begin{tikzpicture}

        \pgfdeclarelayer{background}
		\pgfdeclarelayer{foreground}
			
		\pgfsetlayers{background,main,foreground}

        \begin{pgfonlayer}{background}
        \node (C) [v:ghost] {{\includegraphics[width=16cm]{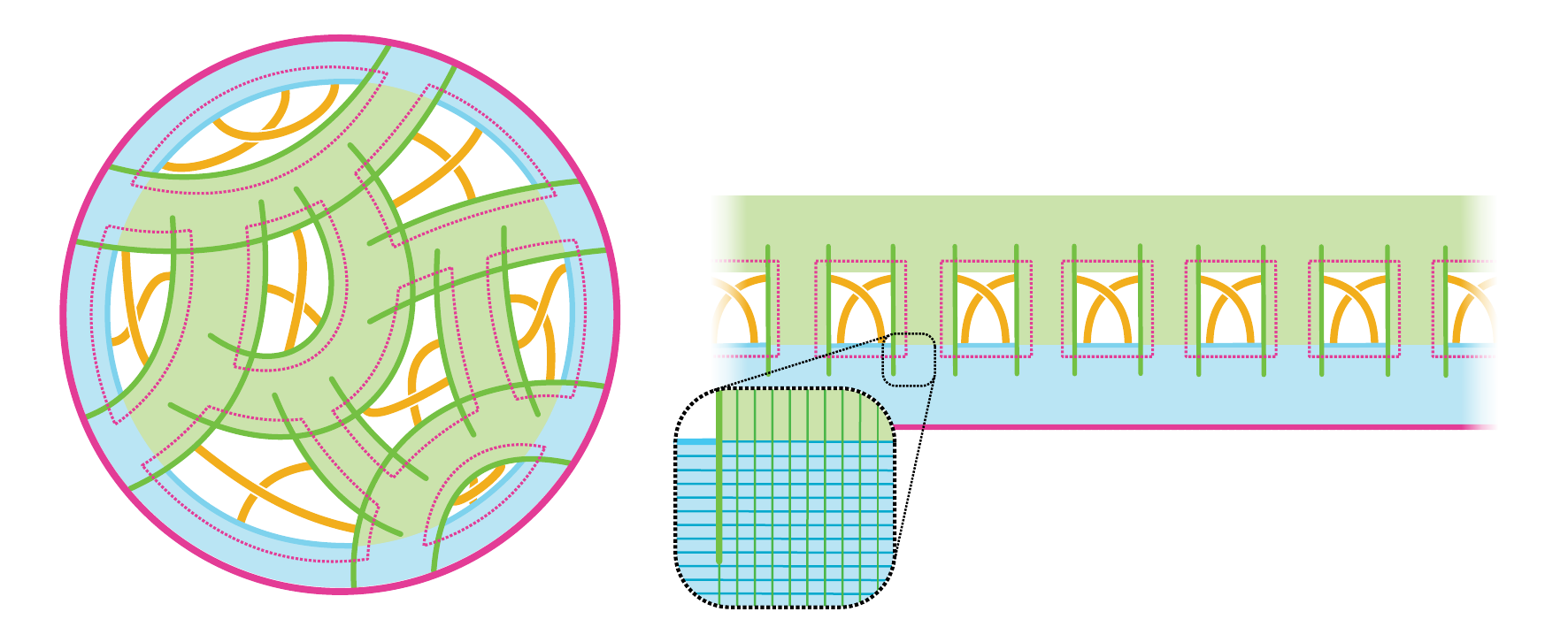}}};
        \end{pgfonlayer}{background}
			
        \begin{pgfonlayer}{main}

        \node (X) [v:ghost,position=270:25mm from C.center] {};
        
        \end{pgfonlayer}{main}

        \begin{pgfonlayer}{foreground}
        \end{pgfonlayer}{foreground}

    \end{tikzpicture}
    \caption{A diagram of a transaction mesh with several districts, each with a private cross (on the left).
    Each district is indicated by a dashed \textcolor{HotMagenta}{magenta} curve.
    The diagram on the right shows a way to normalise the districts using the infrastructure provided by the large transactions in our transaction mesh.
    Notice the similarity to the structure from \cref{fig:K6inGrid} used to create a clique minor.}
    \label{fig:crossesInTransactionmesh}
\end{figure}

Now assume that we are in the situation where we have found a large flat transaction $\mathcal{P}$ in a society $(G,\Omega)$ that is equipped with a big nest.
Moreover, let us assume that there is an even larger crooked transaction $\mathcal{Q}$ such that $\mathcal{P}\subseteq \mathcal{Q}$.
We discuss how to ensure this situation in \cref{subsec:crooks!}.
In this case, we know one of two things:\\
\noindent \textsl{(a)} Either there is at least one path connecting the two newly formed districts of $(G,\Omega)$ and therefore we would need to delete at least one vertex to separate the two sides, or\\
\noindent \textsl{(b)} both sides are already separated.
However, in the second case, because $\mathcal{P}\subseteq \mathcal{Q}$ and $\mathcal{Q}$ is crooked, each of the two newly formed districts must have a private cross.
As discussed above, this situation cannot occur many times without forcing the existence of a $K_t$-minor.
See \cref{fig:crossesInTransactionmesh} for a diagram of this situation.
Hence, most of the time we must be faced with the first outcome.

\begin{figure}[ht]
    \centering
    \begin{tikzpicture}

        \pgfdeclarelayer{background}
		\pgfdeclarelayer{foreground}
			
		\pgfsetlayers{background,main,foreground}

        \begin{pgfonlayer}{background}
        \node (C) [v:ghost] {{\includegraphics[width=10cm]{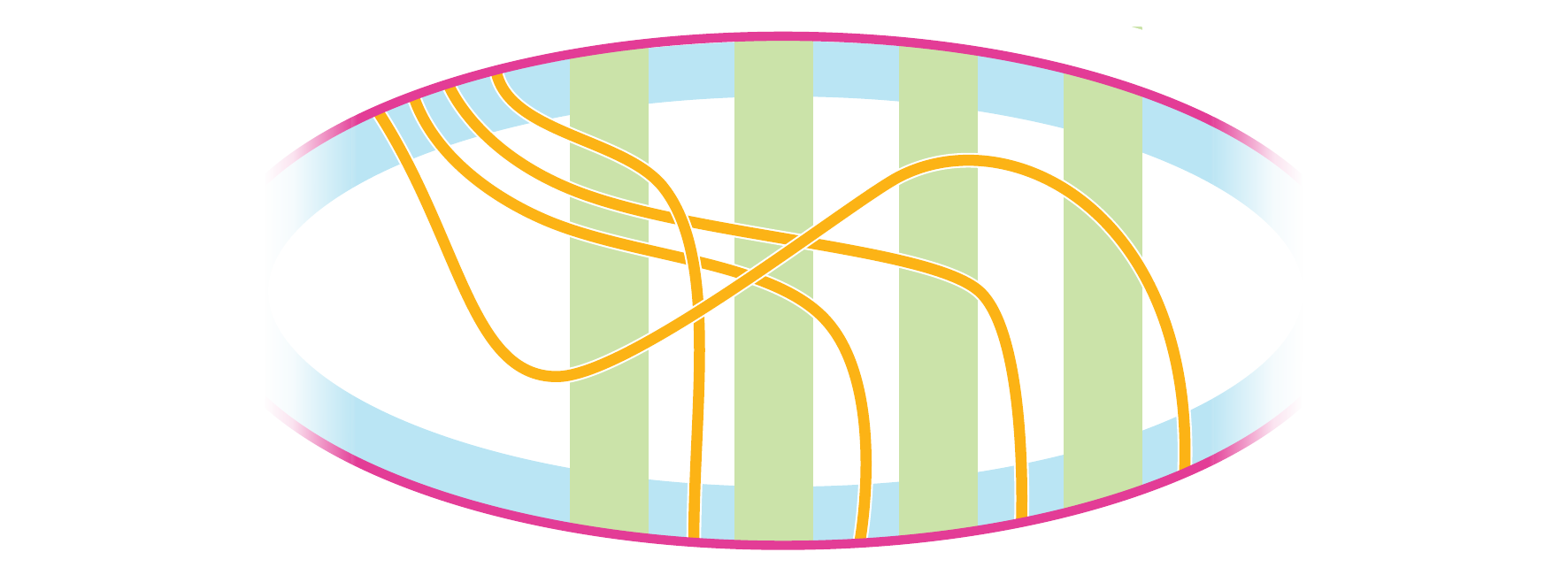}}};
        \end{pgfonlayer}{background}
			
        \begin{pgfonlayer}{main}

        \node (X) [v:ghost,position=270:25mm from C.center] {};
        
        \end{pgfonlayer}{main}

        \begin{pgfonlayer}{foreground}
        \end{pgfonlayer}{foreground}

    \end{tikzpicture}
    \caption{A diagram of a \textsl{leap pattern}, a structure that was used in the proof by KTW \cite{KawarabayashiTW2021Quickly} to force clique minors.}
    \label{fig:leapPattern}
\end{figure}

The way KTW deal with the first outcome is as follows:
Each time two newly formed districts are separated by actually deleting at least one vertex, they keep track of the connection between the two sides, each of them forming a \textsl{leap}.
They then show that, in case they find many such leaps, they can route those paths to form a structure called \emph{leap pattern} (see \cref{fig:leapPattern}).

It is then shown in \cite{KawarabayashiTW2021Quickly} that leap patterns of order polynomial in $t$ contain a $K_t$-minor.
This fact may make the reader expect that we will use leap patterns in our approach.
However, the issue with this is that the strategy of KTW to extract a leap pattern from a large collection of leaps is an inductive one.
In particular, they need to show that whenever they find one additional leap, they can integrate it into an already existing leap pattern.
To ensure the feasibility of this procedure, in each step where they find a leap, they sacrifice roughly half of their current infrastructure to ensure that there is enough space that allows for the extraction of the leap pattern.
Inductively dividing by $2$, however, is a clear cause of exponential blow-up.

This is the second cause of exponentiality in the proof of KTW and a second reason for us relying on the transaction mesh meta structure.
Recall that, unless we are able to increase the Euler-genus of our inductively constructed surface, we do not separate newly formed districts at all.
Instead we keep track of all flat transactions found since the last time we increased the Euler-genus.
Similar to the strategy of R\&S and KTW we make sure to only use flat transactions arising from crooked transactions to construct our transaction mesh.
As a consequence, each district in the transaction mesh arose from a situation of type \textsl{(a)} or \textsl{(b)} as above.
Moreover, similar to before it is relatively easy to argue that the situation of type \textsl{(b)} only occurs few times (recall \cref{fig:K6inGrid}).
But now we also have to show that situations of type \textsl{(a)} force a $K_t$-minor if they occur more than $\mathsf{poly}(t)$ times.

Here we employ a generalisation of a classic result due to Gallai \cite{Gallai1964MaximumMinimum}.
Let $G$ be a graph and $A\subseteq V(G)$ be a set of vertices.
An \emph{$A$-path} in $G$ is a path which has both endpoints in $A$ but is otherwise disjoint from $A$.
Gallai proved that for every integer $k$, there always either exist $k$ pairwise disjoint $A$-paths of a set $S\subseteq V(G)$ with $|S|\in\mathbf{O}(k)$ such that there is no $A$-path in $G-S$.
We prove the following very useful generalisation.

\begin{theorem}\label{GallaiGeneralization_simplified}
    Let $G$ be a graph, let $\mathcal{A}$ be a family of (possibly intersecting) subsets of $V(G)$, and let $k$ be a non-negative integer.
    Then there either exists a family $\mathcal{P}$ of $k$ pairwise disjoint $\bigcup\mathcal{A}$-paths such that each $A\in\mathcal{A}$ contains at most one endpoint of at most one path \(P \in \mathcal{P}\), or there exists a set $S\subseteq V(G)$ with $|S|\in\mathbf{O}(k)$ such that each $\bigcup\mathcal{A}$-path in $G-S$ has an endpoint in some $A\in\mathcal{A}$ with $A\cap S\neq\emptyset$.
\end{theorem}

\begin{figure}[ht]
    \centering
    \begin{tikzpicture}

        \pgfdeclarelayer{background}
		\pgfdeclarelayer{foreground}
			
		\pgfsetlayers{background,main,foreground}

        \begin{pgfonlayer}{background}
        \node (C) [v:ghost] {{\includegraphics[width=16cm]{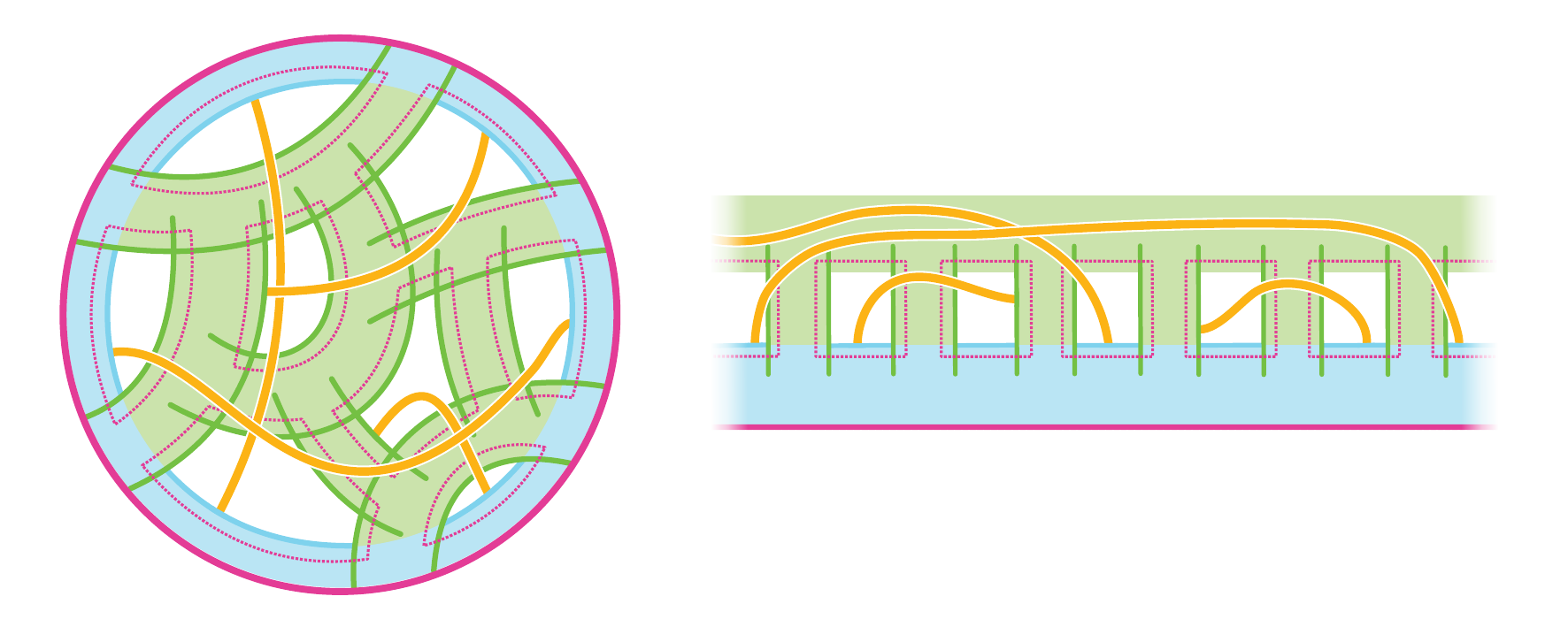}}};
        \end{pgfonlayer}{background}
			
        \begin{pgfonlayer}{main}

        \node (X) [v:ghost,position=270:25mm from C.center] {};
        
        \end{pgfonlayer}{main}

        \begin{pgfonlayer}{foreground}
        \end{pgfonlayer}{foreground}

    \end{tikzpicture}
    \caption{A diagram of a transaction mesh with several districts and a family of pairwise disjoint paths \textsl{jumping} between distinct districts.
    The diagram on the right indicates a noramlised structure that can be extracted from the transaction mesh and, after some further normalisation, yields a clique-minor.}
    \label{fig:jumpsInTransactionmesh}
\end{figure}

Further analysis allows us to strengthen \cref{GallaiGeneralization_simplified} for our setting (see \cref{CorollaryGG} and \cref{PPjumps}).
Equipped with these powerful tools we are now able to exploit the fact that most of our iterations are situations of type \textsl{(a)} as follows.
Consider a transaction mesh with many districts, all arising from situations of type \textsl{(a)}.
Then each such district is connected to at least one other district through a leap.
If we now apply our strengthened version of \cref{GallaiGeneralization_simplified} where each set $A\in\mathcal{A}$ is the boundary of one of our districts, we either can find many disjoint jumps between distinct districts, directly leading to a clique-minor (see \cref{fig:jumpsInTransactionmesh}) or we know that, after deleting few vertices, all remaining jumps must be between a small number of districts.

By building on a construction from KTW \cite{KawarabayashiTW2021Quickly} (see \cref{sec:Mtree,sec:signatures}) we ensure that the number of times we encounter situation \textsl{(a)} is roughly proportional to the number of paths an application of \cref{GallaiGeneralization_simplified} would find between the districts.
Hence, if the number of times we encounter situation \textsl{(a)} grows above a certain threshold in $\mathsf{poly}(t)$, we are sure to find a $K_t$-minor.
Consequently, the total number of iterations we have to do before our induction stops must be bounded by some polynomial in $t$.
See \cref{sec:signatures} for the rigorous analysis of the procedure described above.
The situation we find ourselves in after the construction of the transaction mesh comes to a halt, if it did not find a $K_t$-minor, is the one we described above as a ``transaction mesh of bounded depth''.
As explained above, from here we either find a way to increase the Euler-genus of our embedding in a meaningful way, or we are able to separate all districts from one another.
The latter case would mean that the proof is complete as depicted in \cref{fig:intro_flowchart}.

\subsection{Orthogonalisation of crooked transactions: The key to polynomial bounds}\label{subsec:crooks!}

As described in \cref{subsec:overview}, two out of three sources for exponential growth for the functions presented by KTW \cite{KawarabayashiTW2021Quickly} are their strategy to separate newly formed vortices in every step, and their way to deal with leaps by inductively constructing \textsl{leap patterns}.
Both of these can be fully avoided by our approach through transaction meshes.
However, as in the original strategy of R\&S and the proof of KTW, our own approach hinges on the use of \textsl{crooked transactions}.

A fundamental step in each iteration is the one of \textsl{orthogonalisation}.
Let us first describe what we mean by that and then explain why this step is crucial.

\begin{figure}[ht]
    \centering
    \begin{tikzpicture}

        \pgfdeclarelayer{background}
		\pgfdeclarelayer{foreground}
			
		\pgfsetlayers{background,main,foreground}

        \begin{pgfonlayer}{background}
        \node (C) [v:ghost] {{\includegraphics[width=17.5cm]{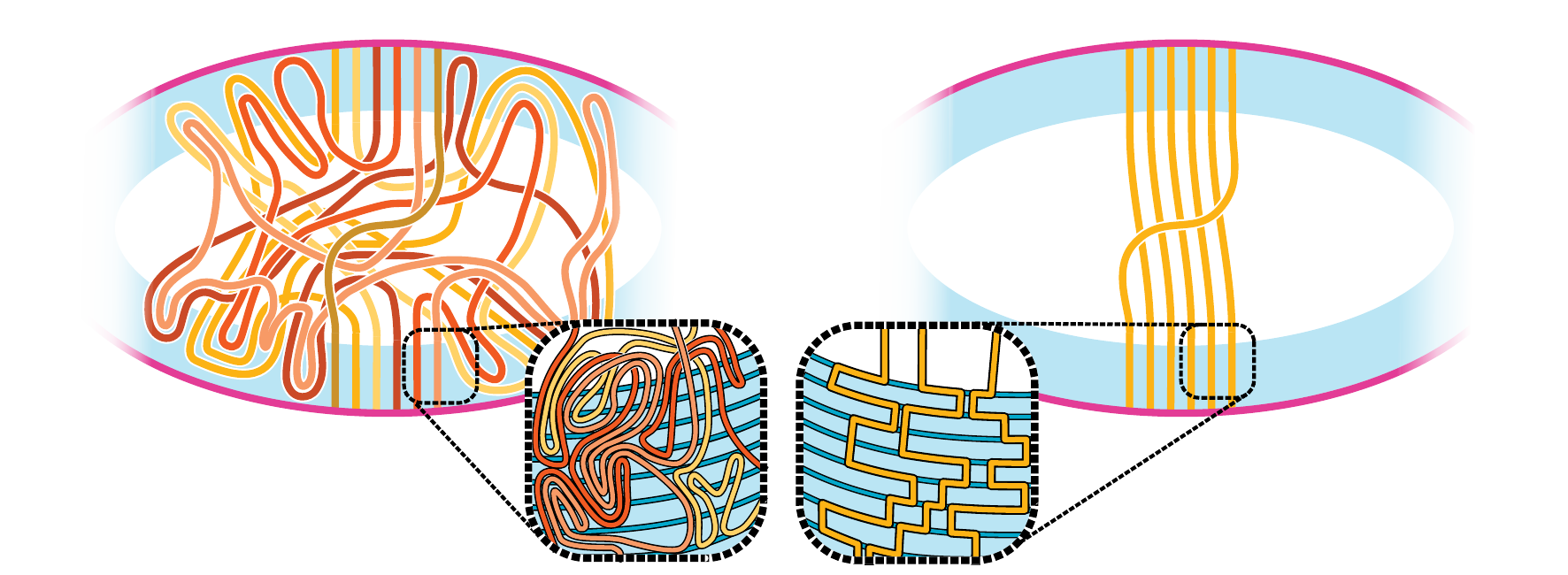}}};
        \end{pgfonlayer}{background}
			
        \begin{pgfonlayer}{main}

        \node (X) [v:ghost,position=270:25mm from C.center] {};
        
        \end{pgfonlayer}{main}

        \begin{pgfonlayer}{foreground}
        \end{pgfonlayer}{foreground}

    \end{tikzpicture}
    \caption{A diagram of two crooked transactions of type 
    \textsl{leap}. The transaction on the right is orthogonal to the nest (indicated by the \textcolor{CornflowerBlue}{blue} area) and the transaction on the left is not.}
    \label{fig:chaoticTransactions}
\end{figure}

As described above, both proofs, the one of KTW and ours, rely fundamentally on iteratively finding large crooked transactions.
Once such a transaction is found, it has to go through several ``clean-up'' steps to ensure that part of it can be used to augment the current embedding.
This is required, among other things, for an application of a variant of the Flat Wall Theorem.
However, the keen observer will notice that a transaction on its own neither is a mesh, nor does it offer enough ``horizontal'' infrastructure to construct something resembling a mesh.
This is where the need for maintaining a large \textsl{nest} comes from.
Recall the diagram on the right of \cref{fig:crossesInTransactionmesh}.
Here we depict a grid-like infrastructure created by vertical paths stemming from a transaction, and horizontal paths contributed by the original nest (indicated by the colours \textcolor{AppleGreen}{green} for the transaction and \textcolor{CornflowerBlue}{blue} for the nest).
This ``grid'' formed by the two families of paths (and cycles) is the result of the transaction behaving in an ``orthogonal'' way with respect to the nest.

Formally, let $(G,\Omega)$ be a society with a nest $\mathcal{C}$ and a transaction $\mathcal{P}$.
We say that $\mathcal{P}$ is \emph{orthogonal} to $\mathcal{C}$ if each path $P\in\mathcal{P}$ intersects every cycle $C\in\mathcal{C}$ exactly twice (once on the way from the nest into the vortex and once on the way back).
Note that this definition only makes sense if a transaction properly enters the vortex.
This is another issue that the property of being crooked solves for us, as this forces our transaction to behave this way, under the assumption that the nest is already embedded in a ``flat'' fashion.
However, an arbitrary crooked transaction in a society does not \textsl{have} to be orthogonal to a given nest.
See \cref{fig:chaoticTransactions} for a depiction of two societies, one with an arbitrary crooked transaction, and one with an orthogonal crooked transaction.

The orthogonalisation serves two goals.
First, if we have a transaction orthogonal to a nest, then the two parts where they intersect form meshes.
This facilitates an application of the Flat Wall Theorem to either find a $K_t$-minor, or -- after deleting a small apex set -- obtain a still large subtransaction $\mathcal{Q}$ whose \textsl{strip}, i.e. the part of the graph that attaches to the ``inner'' paths of the transaction, is flat and does not communicate with the part of $\Omega$ that is not involved in the transaction.
Second, it allows us to take this newly obtained embedding of $\mathcal{Q}$ and \textsl{re-conciliate} it with the original flat embedding of (most of) the nest.
This last step is crucial because this is what finally allows us to consistently grow the embedding.

Let us quickly discuss \textsl{why} it is possible to talk about crooked transactions in the first place.
In the 9th entry of the Graph Minors Series \cite{RobertsonS1990Graph}, R\&S proved the following key theorem.

\begin{proposition}[Robertson and Seymour \cite{RobertsonS1990Graph} (see also KTW, Lemma 3.6 \cite{KawarabayashiTW2021Quickly})]\label{prop:3.6}
Let $(G,\Omega)$ be a society and $k\geq 2$ be an integer.
Then either
\begin{enumerate}
    \item $G$ can be drawn, up to $\leq 3$-sums, on a disk except for one vortex of depth in $\mathbf{O}(k)$, and with the vertices $\Omega$ on the boundary and in the correct order, or
    \item there exists a crooked transaction of order $k$ in $(G,\Omega)$.
\end{enumerate}
\end{proposition}

Interestingly, while \cref{prop:3.6} is the main theorem of \cite{RobertsonS1990Graph}, this paper appears to be mostly cited for the proof of the Two Paths Theorem it contains.
Notice, that \cref{prop:3.6} can be seen as a generalisation of the Two Paths Theorem itself, since a cross is the unique smallest crooked transaction (of type \emph{leap}).

Given a society with a nest and a large transaction, it is relatively simple to find a still large crooked transaction which is orthogonal to all cycles of the next except for the inner-most one (see \cref{lemma:orthogonal_transaction}).
However, we are not asking for \textsl{any} large orthogonal transaction, it is crucial for us that the transaction \textsl{remains} crooked after the process of orthogonalising.

\paragraph{Orthogonal crooked transactions {\`a} la KTW.}
The way KTW maintain crookedness is due to a mix of techniques from the proof of \cref{prop:3.6}, an augmenting path argument from the proof of Menger's Theorem, and Menger's Theorem itself.
In a nutshell, KTW show that given a crooked transaction of order $k$ in a society with a nest of order $\mathbf{O}(k)$, then there also exists a crooked transaction of order $k$ that is orthogonal to the outermost $k$ cycles of the nest.
To guarantee this, they make use of the fact that in the innermost $\mathbf{O}(k)$ cycles of the nest, while the transaction is not necessarily orthogonal here, it at least establishes large connectivity between the vertices of $\Omega$ and the parts of the transaction living inside the vortex.
This large amount of connectivity guarantees that a clever application of Menger's Theorem yields a large linkage that is then further refined to yield the sought-after transaction.

\begin{figure}[ht]
    \centering
    \begin{tikzpicture}

        \pgfdeclarelayer{background}
		\pgfdeclarelayer{foreground}
			
		\pgfsetlayers{background,main,foreground}

        \begin{pgfonlayer}{background}
        \node (C) [v:ghost] {{\includegraphics[width=17.5cm]{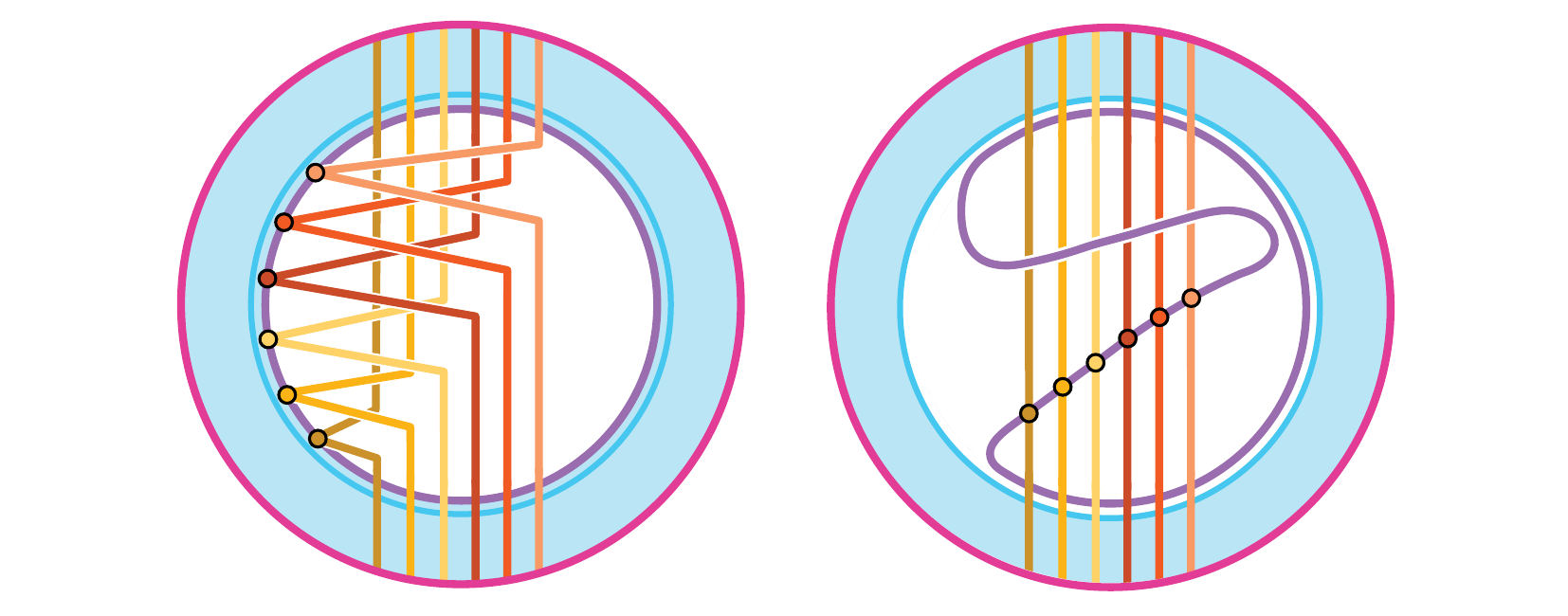}}};
        \end{pgfonlayer}{background}
			
        \begin{pgfonlayer}{main}

        \node (X) [v:ghost,position=270:25mm from C.center] {};

        \end{pgfonlayer}{main}

        \begin{pgfonlayer}{foreground}
        \end{pgfonlayer}{foreground}

    \end{tikzpicture}
    \caption{A transaction $\mathcal{P}$ in a society that is not orthogonal to the nest $\mathcal{C}$ on the left. On the right, the same transaction $\mathcal{P}$ is drawn in a way that highlights that it is orthogonal to all but the innermost cycle of $\mathcal{C}$ (depicted in \textcolor{Amethyst}{purple}).
    Notice, that there does not exist any transaction of large order in the graph depicted that is orthogonal to the entire nest.
    It follows, that it is not possible to prove a statement that finds, given the existence of a large transaction, a large transaction that is orthogonal to the \textsl{entire} nest. 
    As a side note, observe that in the drawing on the right, $\mathcal{P}$ \textsl{and} all of the nest except for the innermost cycle are embedded in a planar way.
    This is impossible if we also want to embed the innermost cycle.}
    \label{fig:sacrificeIsNecessary}
\end{figure}

While their lemma always guarantees the existence of a crooked transaction orthogonal to a large portion of the nest, and even of the same size as the original, it has one fatal flaw.
Each application of the KTW orthogonalisation lemma requires the sacrifice of $\Omega(k)$ many cycle from the nest.
Since a part of the crooked transaction is then used to build the nests for the next iteration, this forces a telescoping effect on the required size of the starting nest, hence making exponential bounds unavoidable.
Consequently, not only the lemma of KTW, but this entire method of orthogonalisation is unsuitable for the purpose of staying within polynomial bounds.

\paragraph{Orthogonal crooked transaction: The hard way.}

It is worth pointing out that sacrificing \textsl{some} part of the nest is completely unavoidable as there might exist large (crooked) transactions that interacts with the innermost cycle of the nest in a way that fundamentally obstructs the existence of a flat embedding of the entire nest and part of the transaction.
See \cref{fig:sacrificeIsNecessary} for an illustration of such a situation.
However, it would be enough if we managed to only sacrifice a \textsl{constant} number of cycles from the nest while orthogonalising an \textsl{arbitrarily} large crooked transaction.

\begin{figure}[ht]
    \centering
    \begin{tikzpicture}

        \pgfdeclarelayer{background}
		\pgfdeclarelayer{foreground}
			
		\pgfsetlayers{background,main,foreground}

        \begin{pgfonlayer}{background}
        \node (C) [v:ghost] {{\includegraphics[width=17.5cm]{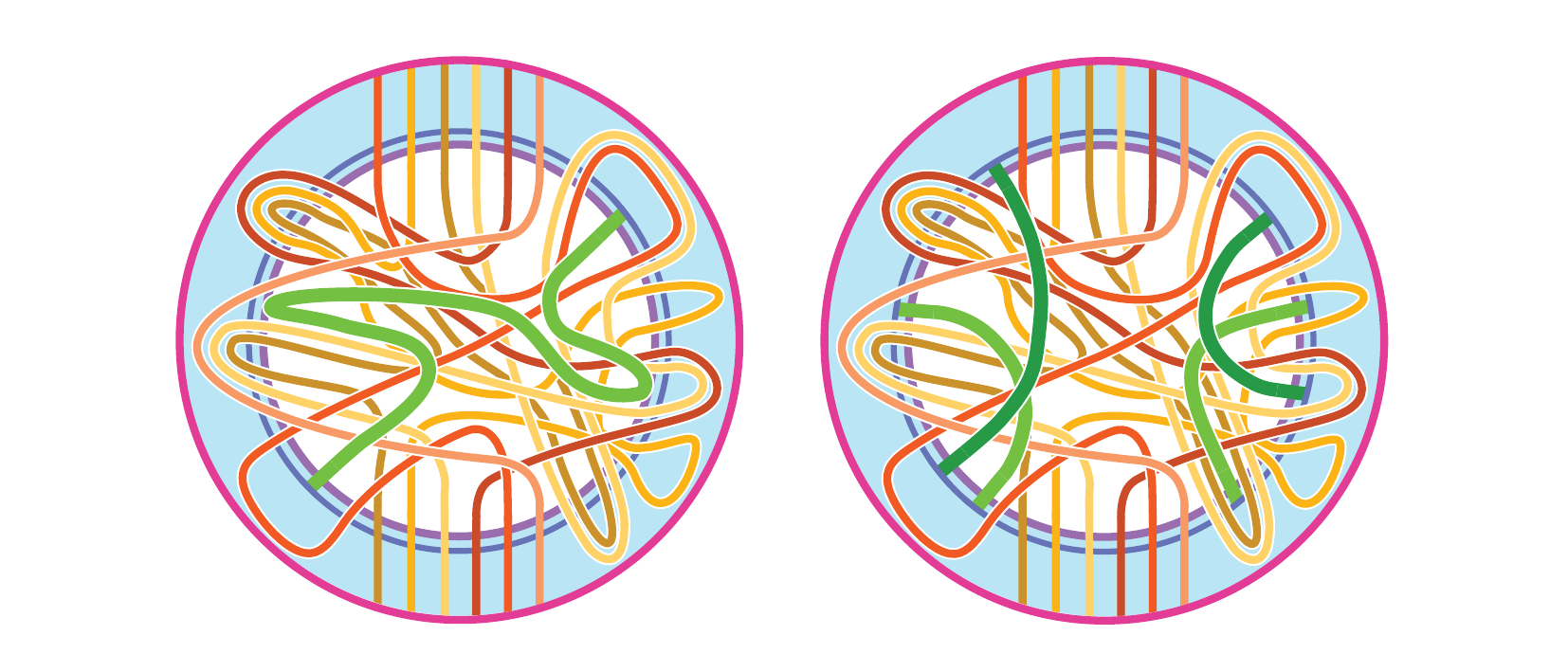}}};
        \end{pgfonlayer}{background}
			
        \begin{pgfonlayer}{main}

        \node (X) [v:ghost,position=270:25mm from C.center] {};

        \end{pgfonlayer}{main}

        \begin{pgfonlayer}{foreground}
        \end{pgfonlayer}{foreground}

    \end{tikzpicture}
    \caption{Two transactions (left and right), each in their own society with a large nest.
    On the left, there is a single leap (coloured \textcolor{AppleGreen}{green}) ``over'' the transaction with both endpoints on the second innermost cycle of the nest.
    Similarly, on the right there are two crosses, one ``left'' and one ``right'' of the transaction.
    Both of these crosses have their endpoint on the second innermost cycle of the nest and are otherwise disjoint from this cycle.
    Notice that in neither of these situations, this ``witness'' is actually enough to confirm that the transaction is indeed crooked.
    However, they provide some evidence of non-planar behaviour that turns out to be strong enough to facilitate our arguments.}
    \label{fig:introWitnesses}
\end{figure}

For this, we follow the route as laid out below.\\
\textsl{(i)} We find a large subtransaction $\mathcal{P}_1$ of the initial crooked transaction $\mathcal{Q}$ together with a \textsl{witness} $\mathcal{W}$ of the crookedness of $\mathcal{P}$ in the form of a few paths, all of which only interact with the innermost two cycles of the nest.
See \cref{fig:introWitnesses} for an illustration for this initial situation.

\begin{figure}[ht]
    \centering
    \begin{tikzpicture}

        \pgfdeclarelayer{background}
		\pgfdeclarelayer{foreground}
			
		\pgfsetlayers{background,main,foreground}

        \begin{pgfonlayer}{background}
        \node (C) [v:ghost] {{\includegraphics[width=15cm]{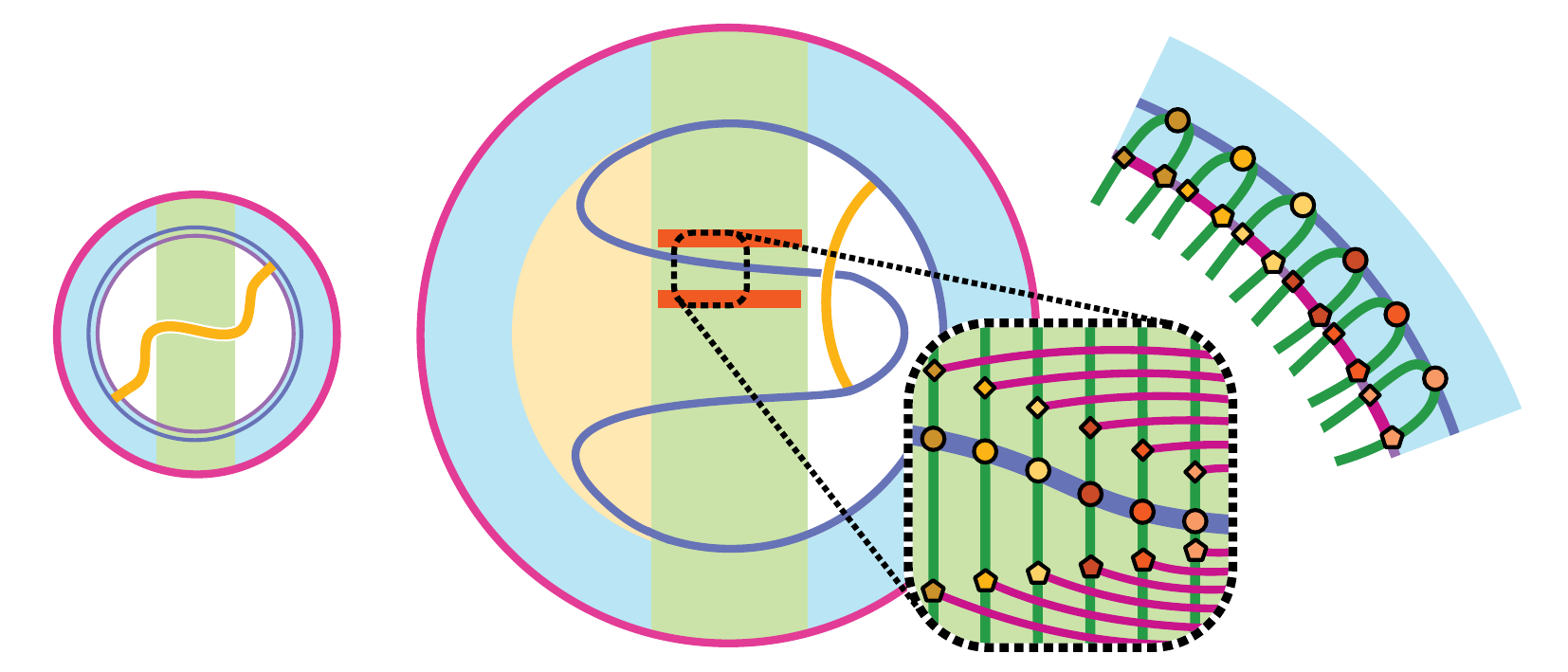}}};
        \end{pgfonlayer}{background}
			
        \begin{pgfonlayer}{main}

        \node (X) [v:ghost,position=270:25mm from C.center] {};

        \node (A) [v:ghost,position=117:30mm from X] {$A$};
        \node (B) [v:ghost,position=90:9.5mm from A] {$B$};

        \node (C2mid) [v:ghost,position=270:23mm from A] {$C_2$};

        \node (C1right) [v:ghost,position=10:60mm from X] {$C_1$};
        \node (C2right) [v:ghost,position=22:7mm from C1right] {$C_2$};

        \end{pgfonlayer}{main}

        \begin{pgfonlayer}{foreground}
        \end{pgfonlayer}{foreground}

    \end{tikzpicture}
    \caption{The situation before (on the left) and after (in the middle) orthogonalising and then flattening a transaction under the presence of a witness. If this process does not lead to an orthogonal crooked transaction, the cycle $C_2$ of the original nest where the endpoints of the witness are rooted must be ``pulled'' through the entire transaction (this is similar to the situation in \cref{fig:sacrificeIsNecessary}) in order to facilitate a flat embedding. In the situation depicted, the left side of the transaction (marked in \textcolor{BananaYellow}{yellow}) is flat. Finally, the combination of the flatness of the transaction itself and $C_2$ being forced to intersect every path of the transaction implies that there is a large linkage $\mathcal{L}$, fully contained in the cycle $C_1$, connecting the two ``sides'' ($A$ and $B$) of $C_2$ (as indicated on the right).
    With $\mathcal{L}$ being contained in $C_1$, it is fully disjoint from $C_2$ in particular.}
    \label{fig:flatteningWitnesses}
\end{figure}
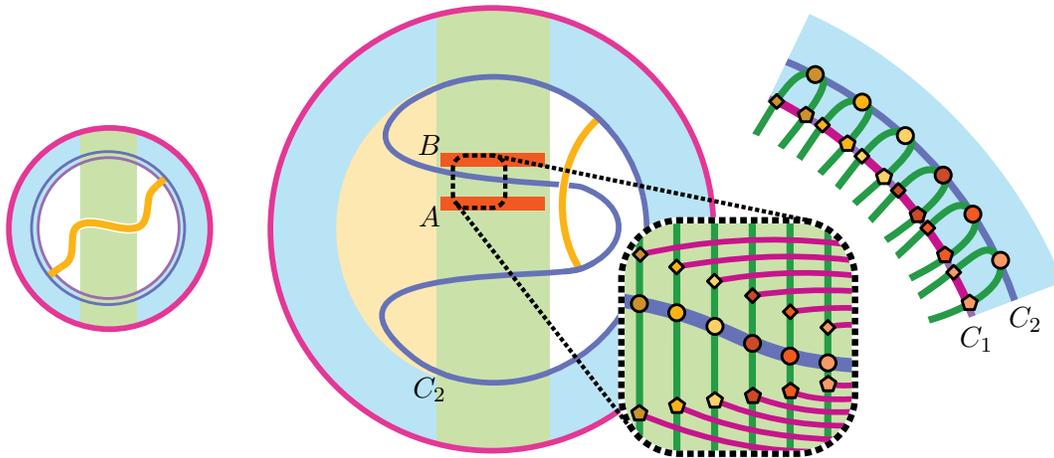

\textsl{(ii)} Then, we apply our simple orthogonalisation lemma (i.e.\@ \cref{lemma:orthogonal_transaction}) to the graph where $\mathcal{W}$ is removed to obtain a large orthogonal transaction $\mathcal{P}_2$ that starts and ends on the same vertices as the transaction $\mathcal{P}_1$.\\
\textsl{(iii)} In case $\mathcal{P}_2$ is already crooked, we immediately win.
Otherwise, we now begin to apply the techniques from the proof of \cref{prop:3.6} to $\mathcal{P}_2$ after we returned $\mathcal{W}$ to the graph.
This either yields a big crooked transaction -- and thus we win -- or it finds a large subtransaction $\mathcal{P}_3$ of $\mathcal{P}_2$ such that $\mathcal{P}_3$ itself and one of the two sides of $\mathcal{P}_3$ are entirely flat, let us call this one the \emph{flat side}.
Notice that it is impossible for both sides to be flat since $\mathcal{Q}$ was a huge crooked transaction.\\
\textsl{(iv)} If we now examine our witness $\mathcal{W}$, we can observe that the flat side of $\mathcal{P}_3$ must now either contain the entire witness, or be disjoint from it.
Indeed, if the witness is a doublecross, it cannot be contained in the flat side at all.
Since $\mathcal{W}$ has at least one vertex, say $w$, of the cycle $C_2$, i.e.\@ the second innermost cycle of the nest, such that $w$ \textsl{should} belong to the flat side, this observation indicates that $C_2$ must intersect every path in $\mathcal{P}_3$ several times.
See the left and middle part of \cref{fig:flatteningWitnesses} for an illustration of this situation.

We now focus on the the innermost cycle of the original nest, namely $C_1$.
The cycles $C_1$ and $C_2$ were both originally part of the flat area of the embedding.
While this might have changed now, we can still access this previous information to understand how $\mathcal{P}_3$ can intersect $C_2$.
In order to do so, each path $P\in\mathcal{P}_3$ must always first pass through $C_1$ before it can make contact with $C_2$, then it must return to $C_1$.
With both the strip of $\mathcal{P}_3$ and the area between the two cycles $C_1$ and $C_2$ being flat, the only way the entire transaction can touch $C_2$ is if all paths of $\mathcal{P}_3$ -- in order -- realise this kind of pattern.
This situation is indicated in the right part of \cref{fig:flatteningWitnesses}.

Notice now, that for each $P\in\mathcal{P}_3$ there is a subpath $L_P$ of $C_1$ connecting the ``entry'' and ``exit'' point that $P$ must pass through on $C_1$ in order to touch $C_2$.
The collection $\mathcal{L}$ of these paths is a huge linkage, entirely disjoint from $C_2$, in itself.
We may now partition the set of endpoints of $\mathcal{L}$ into sets $A$ and $B$ such that the linkage $\mathcal{L}$ has to ``cross over'' parts of $C_2$ in order to reach from $A$ to $B$.
See \cref{fig:flatteningWitnesses} for a sketch of the situation.
It turns out, that this ``crossing'' of $\mathcal{L}$ over $C_2$ is enough to guarantee that this situation will always give rise to a large crooked transaction orthogonal to all but a constant number of cycles of the original nest.

We strongly believe that overcoming the issues in finding an orthogonal crooked transaction by sacrificing a portion of the nest that does not depend on the size of the transaction was the most challenging part of our proof.
With this, we have given a (very) rough overview on all major steps in the proof of the LST.
Please notice that, due to the complexity of the topic, this overview is far from perfect and relies on countless convenient simplifications.
We advise the cautious reader to use this overview as a guide and a (sometimes a bit faulty) source of intuition for approaching the main body of this paper.

\subsection{The structure of this paper}\label{subsec:organisation}
In the following we provide a structured overview on the organisation of this paper to allows for better orientation.

\paragraph{Preliminaries.}
There are two bulky sets of definitions we need to import.
The first are the general notions of graph theory together with definitions from the theory of graph minors that are known to a wider audience such as \textsl{treewidth}, \textsl{walls}, and \textsl{tangles}.
These definitions can be found in \cref{sec:prelim}.

We then provide an introduction to the definitions and tools used by KTW in their proof of the GMST in \cref{sec:decompositions}.
Here, a statement of the Two Paths Theorem can be found alongside formalisations of drawings and a rigorous definition of \textsl{$\Sigma$-embeddings}.

\paragraph{Finding a $K_t$-minor.}
As described above, there are, in essence, three ways how we can find a clique minor: Through many crosses, many disjoint jumps, and through large Euler-genus.
The corresponding lemmas are proven in \cref{sec:forcing} for crosses and jumps, and \cref{sec:lowerbound} for Euler-genus.

\paragraph{Flat meshes and transactions.}
In this paper we provide a stand-alone proof of the Flat Wall Theorem (or better the Flat \textsl{Mesh} Theorem).
This is done to streamline some of our arguments, but also to provide constructive and non-randomised bounds, strictly improving over those found in \cite{KawarabayashiTW2018New}.
The Flat Mesh Theorem is proven in \cref{sec:flattening} and the ``Flat Transaction Theorem'' can be found in \cref{sec:flatteTransaction}.

\paragraph{A cozy nest.}
Most of our arguments about orthogonalising transactions are based on certain extremal assumptions on the structure of our nests called \emph{cozy}.
In \cref{sec:cozy} we provide a definition for this property and show it it can be achieved and maintained throughout our construction.

\paragraph{Orthogonal crooked transactions.}
The next major point on the agenda is finally the theorem about orthogonal crooked transactions.
In \cref{sec:findcrooked} we open up the black box of \cref{prop:3.6} and provide a constructive way to either find a large crooked transaction in a society, or an embedding with a single vortex of bounded depth.
We then employ many of these tools in \cref{sec:crookedmagic} to show that we can always find a large crooked transaction orthogonal to all but a constant number of cycles of the nest.
Finally, in \cref{sec:multisocietycrooked} we extend this result to also work on transaction meshes.
This last step is necessary because our approach necessitates that we keep the connectivity between our vortices whilst we are trying to extend the $\Sigma$-decomposition we have built up until that point.
As a consequence, any transaction we find might end up accidentally embedding parts of a foreign district.

\paragraph{Society classification.}
We are now approaching the proof of one of the major theorems of this paper: A society classification theorem with polynomial bounds.
First, we provide a formal definition for transaction meshes in \cref{sec:transactionmeshes}.
Then, in \cref{sec:ApathsArgument}, we prove several routing lemmas necessary to find a $K_t$-minor in a transaction mesh with many districts that either have private crosses, or many jumps between them.
Finally, \cref{sec:societyclassification} brings everything together and proves that in any society with large nest we can either find a $K_t$-minor, a way to increase the Euler-genus by means of adding a handle or a crosscap, or we find a small number of vortices of bounded depth after deleting few vertices.

\paragraph{The Local Structure Theorem.}
In \cref{sec:localstructure} we finally derive the LST.
Here, we first prove a slightly weaker variant with polynomial bounds but without tight bounds in the Euler-genus.
We then apply this theorem to solve two conjectures of Wollan, the first of which finally implies \cref{thm:intro_local_structure}.

\paragraph{The Graph Minor Structure Theorem.}
The final proof of the GMST, and thus \cref{thm:mainthm_simplest}, can be found in \cref{sec:localtoglobal}.

\paragraph{Concluding remarks.}
Our conclusion is two-fold.
In \cref{sec:polytime} we describe a \textsl{fully} polynomial-time randomised algorithm that finds either a forbidden graph $H$ or the structural decomposition with high probability.
We conclude this paper with a collection of open problems in \cref{sec:conclusion}.

\section{Preliminaries}\label{sec:prelim}
We effectively split our preliminaries into two major sections.
First, in this section, we introduce most of the basic notions associated with graphs and particularly minors.
Additionally, we introduce what is needed to state the grid theorem and some notions related to the Graph Minor Structure Theorem that may be fairly familiar to most readers, like meshes and tangles.
For the more involved definitions, especially those regarding embeddings in surfaces, we reserve \Cref{sec:decompositions}.

While the original version of the Graph Minor Structure Theorem (GMST) was introduced by Robertson and Seymour \cite{RobertsonS2003Graph}, we adapt most definitions related to the GMST from the work of Kawarabayashi et al., who have pioneered modernised versions of the major tools from Robertson and Seymour's series of papers (see \cite{KawarabayashiTW2018New} and \cite{KawarabayashiTW2021Quickly}).
Definitions unique to our work will be introduced throughout later sections whenever they become necessary.

\medskip

By $\mathbb{Z}$ we denote the set of integers and by $\mathbb{R}$ the set of reals.
Given any two integers $a,b\in\mathbb{Z}$, we write $[a,b]$ for the set $\{z\in\mathbb{Z} ~\!\colon\!~ a\leq z\leq b\}$.
Notice that the set $[a,b]$ is empty whenever $a>b$.
For any positive integer $c$ we set $[c]\coloneqq [1,c]$.
In order to avoid ambiguity for this notation, we denote closed intervals over the reals by $[x,y]_{\mathbb{R}}$.
We however rarely need these.

On many occasions we will deal with a set $\mathcal{G}$ of graphs and want to form the union of all of them.
We will often write $\bigcup \mathcal{G}$ to denote $\bigcup_{G \in \mathcal{G}} G$ for the sake of simplicity.
In particular, if $\mathfrak{G} = \{ \mathcal{G}_1, \ldots , \mathcal{G}_k \}$ is a set of sets of graphs, we further let $\bigcup \mathfrak{G}$ denote $\bigcup_{\mathcal{G} \in \mathfrak{G}} \bigcup \mathcal{G}$.

\paragraph{Separations.}
A \emph{separation} in a graph $G$ is a pair $(A,B)$ of vertex sets such that $A \cup B=V(G)$ and there is no edge in $G$ with one endpoint in $A\setminus B$ and the other in $B\setminus A$.
The \emph{order} of $(A,B)$ is $|A\cap B|$.
Given sets $X,Y\subseteq V(G)$ an \emph{$X$-$Y$-separation} in $G$ is a separation $(A,B)$ such that $X\subseteq A$ and $Y\subseteq B$.
\medskip

\paragraph{Linkages and paths.}
A \emph{linkage} $\mathcal{L}$ in a graph $G$ is a set of pairwise vertex-disjoint paths and, in a slight abuse of notation, we use $V(\mathcal{L})$ and $E(\mathcal{L})$ to denote $\bigcup_{L \in \mathcal{L}} V(L)$ and $\bigcup_{L \in \mathcal{L}} E(L)$ respectively.

We say that a path $P$ is \emph{internally disjoint} from a set $X$, or a subgraph $H$, if $V(P) \cap X$, respectively $V(P) \cap V(H)$, does not contain any non-endpoint vertex of $P$.
Let $A, B \subseteq V(G)$, then an $A$-path is a path of length at least one with both endpoints in $A$ that is internally disjoint from $A$.
If $H$ is a subgraph of $G$, we also denote $V(H)$-paths simply as $H$-paths.
Furthermore, an \emph{$A$-$B$-path} is a path \(P = x_0 \cdots x_\ell\) with $V(P) \cap A = \{x_0\}$ and $V(P) \cap B = \{x_\ell\}$.
If $A$ contains only a single element $a \in A$, we also denote an $A$-$B$-path as an $a$-$B$-path.
Both $a$-$b$-paths and $A$-$b$-paths similarly denote paths between sets which may be singletons.
An \emph{$A$-$B$-linkage} for two vertex sets $A,B \subseteq V(G)$ is a linkage consisting of $A$-$B$-paths.

\begin{proposition}[Menger's Theorem \cite{Menger1927Zur}]\label{prop:mengersthm}
Let $G$ be a graph and $X,Y\subseteq V(G)$ be two sets of vertices.
Then the minimum order of a $X$-$Y$-separation in $G$ equals the maximum number of paths in an $X$-$Y$-linkage in $G$.
\end{proposition}

A well-known algorithm of Ford and Fulkerson \cite{FordF1956Maximal} takes time $\mathbf{O}(|E(G)|)$ to find, given an $X$-$Y$-linkage of size $k$, an $X$-$Y$-linkage of size $k+1$ if one exists.
It follows that we can find an $X$-$Y$-linkage of order $k$ in time $\mathbf{O}(k|E(G)|)$ or determine that no such linkage exists and find an $X$-$Y$-separation of minimum order instead.

If $P$ is a path and $x,y \in V(P)$ are vertices on $P$, we denote by $xPy$ the subpath of $P$ with the endpoints $x$ and $y$.
Let $P$ be a path from $s$ to $t$ and $Q$ be a path from $q$ to $p$.
If $x$ is a vertex in $V(P) \cap V(Q)$ and $sPx$ and $xQp$ only intersect in $x$, then we let $sPxQp$ be the path obtained from the union of $sPx$ and $xQp$.

\subsection{Treewidth, minors, and walls}\label{sec:walls}
On the most basic level, we are interested in decomposing a graph into simple parts that are easy to separate from each other.
A pair of notions that embody this intuition well, are central to the theory of graph minors, and have by now established themself as classic concepts in graph theory are tree-decompositions and the associated treewidth of a graph.

\paragraph{Tree-decompositions and treewidth.}
A \emph{tree-decomposition} of a graph $G$ is a tuple $(T, \beta)$, where $T$ is a tree and $\beta \colon V(T) \rightarrow 2^{V(G)}$ is a function such that
    \begin{itemize}
        \item for each $e \in E(G)$ there exists a vertex $t \in V(G)$ such that $e \subseteq \beta(t)$, and

        \item for each $v \in V(G)$ the set $\beta^{-1}(v)$ induces a connected, non-empty subtree in $T$.
    \end{itemize}
The \emph{width} of $(T,\beta)$ is $\max(\{ |\beta(t)| - 1 ~\!\colon\!~ t \in V(T)\})$ and the \emph{treewidth} $\mathrm{tw}(G)$ of $G$ is the minimum width of any tree-decomposition of $G$.
We call the sets $\beta(t)$ for $t \in V(T)$ the \emph{bags} of $(T, \beta)$ and the value $\max_{tt' \in E(T)}( | \beta(t) \cap \beta(t') | ) $ is called the \emph{adhesion} of the tree-decomposition.\footnote{For tree-decompositions with a single vertex we define the adhesion to be 0.}\medskip

This notion first appeared in a different guise as \emph{dimension} in the work of Bertelè and Brioschi \cite{BerteleB1972Nonserial,BerteleB1973Nonserial}.
Later Halin introduced the equivalent concept of \emph{$S$-functions} \cite{Halin1976$S$functions}, which was followed by Arnborg defining \emph{partial $k$-trees} \cite{Arnborg1985Efficient}, which are equivalent to all of the notions we just mentioned.
Finally, the definition we gave above is what Robertson and Seymour introduced as treewidth and tree-decompositions in \cite{RobertsonS1983Graph,RobertsonS1986Grapha}.

As we will see in a moment, the treewidth of $G$ and the adhesion of a tree-decomposition of minimum width for $G$ can differ substantially.
The reason as to why this value is called adhesion is given by the following standard lemma.

\begin{observation}\label{lem:treedecompositionseparator}
    Let $G$ be a graph and let $(T, \beta)$ be a tree-decomposition of $G$ of adhesion $k \in \mathbb{N}$.
    Let \(x_1 x_2 \in E(T)\), and let \(T_1\) and \(T_2\) be the two components of \(T - x_1 x_2\) where each \(T_i\) is the component containing \(x_i\).
    The tuple $(\bigcup_{u \in V(T_t)} \beta(u) , \bigcup_{v \in V(T_{t'})} \beta(v))$ is a separation of order at most $k$.
\end{observation}

\paragraph{Minors.}
Given a graph $G$ with an edge $e = uv$, the \emph{contraction} of $e$ is the operation of identifying $u$ and $v$ into a single vertex and subsequently deleting all loops and parallel edges.
Note that, if we find a subset of vertices $U \subseteq V(G)$ that induces a connected subgraph in $G$, we can contract this set into a single vertex via repeatedly contracting edges.
We say that $H$ is a \emph{minor} of $G$ and conversely we that $G$ has an \emph{$H$-minor} if a graph isomorphic to $H$ can be obtained from a subgraph of $G$ by repeatedly contracting edges.
\medskip

It will often be useful to have an exact description of how we find a minor within a graph and the following notion makes this possible.
Let $G$ be a graph, let $uv\in E(G)$, and let $P$ be a path with endpoints $u$ and $v$ such that no internal vertex of $P$ belongs to $G$.
We say that the graph $G'\coloneqq (G-uv)+P$ is obtained from $G$ by \emph{subdividing} the edge $uv$.
A graph $G''$ is called a \emph{subdivision} of $G$ if it can be obtained from $G$ by subdividing aby subset of edges of $G$.
Let $v\in V(G)$ be a vertex of degree at least four.

We will rarely want to talk about minors in our proofs -- and if we do, it is mainly because some lemma we use guarantees the existence of a helpful minor in some situation.
The reason why we can largely avoid this is found in the following lemma, which guarantees that as long as we are interested in subcubic graphs,\footnote{A graph is called \emph{subcubic} if its maximum degree is at most 3.} we can find the graph in question as a subdivision if we can find it as a minor.

\begin{lemma}[folklore]\label{lem:subcubicminorassubdivision}
    Let $G$ and $H$ be graphs where $\Delta(H)\leq 3$.
    Then $G$ contains $H$ as a minor if and only if $G$ contains a subdivision of $H$ as a subgraph.
\end{lemma}

\begin{figure}[ht]
	\begin{center}
            		\begin{tikzpicture}[scale=0.75]
			
			\foreach\i in {0,1,2,3,4,5}
			{
                    \node (A\i) at (\i,0) [draw, circle, scale=0.6, fill] {};
                    \node (B\i) at (\i,1) [draw, circle, scale=0.6, fill] {};
                    \node (C\i) at (\i,2) [draw, circle, scale=0.6, fill] {};
                    \node (D\i) at (\i,3) [draw, circle, scale=0.6, fill] {};
                    \node (E\i) at (\i,4) [draw, circle, scale=0.6, fill] {};
            }

            \foreach\i in {0,1,2,3,4}
			{
		    	\path
                    (A\i) edge[very thick] (B\i)
                    (B\i) edge[very thick] (C\i)
                    (C\i) edge[very thick] (D\i)
                    (D\i) edge[very thick] (E\i)
			    ;
                \pgfmathtruncatemacro\iplus{\i+1}
                \path
                    (A\i) edge[very thick] (A\iplus)
                    (B\i) edge[very thick] (B\iplus)
                    (C\i) edge[very thick] (C\iplus)
                    (D\i) edge[very thick] (D\iplus)
                    (E\i) edge[very thick] (E\iplus)
			    ;
			}
                \path
                    (A5) edge[very thick] (B5)
                    (B5) edge[very thick] (C5)
                    (C5) edge[very thick] (D5)
                    (D5) edge[very thick] (E5)
			    ;
			\end{tikzpicture} \qquad \qquad
					\begin{tikzpicture}[scale=0.75]
			
			\foreach\i in {1,2,3,4}
			{
                    \node (A\i) at (\i,0) [draw, circle, scale=0.6, fill] {};
                    \node (B\i) at (\i,1) [draw, circle, scale=0.6, fill] {};
                    \node (C\i) at (\i,2) [draw, circle, scale=0.6, fill] {};
                    \node (D\i) at (\i,3) [draw, circle, scale=0.6, fill] {};
                    \node (E\i) at (\i,4) [draw, circle, scale=0.6, fill] {};
            }
            \node (A0) at (0,0) [draw, circle, scale=0.6, fill, opacity=0.1] {};
            \node (B0) at (0,1) [draw, circle, scale=0.6, fill] {};
            \node (C0) at (0,2) [draw, circle, scale=0.6, fill] {};
            \node (D0) at (0,3) [draw, circle, scale=0.6, fill] {};
            \node (E0) at (0,4) [draw, circle, scale=0.6, fill] {};
            \node (A5) at (5,0) [draw, circle, scale=0.6, fill] {};
            \node (B5) at (5,1) [draw, circle, scale=0.6, fill] {};
            \node (C5) at (5,2) [draw, circle, scale=0.6, fill] {};
            \node (D5) at (5,3) [draw, circle, scale=0.6, fill] {};
            \node (E5) at (5,4) [draw, circle, scale=0.6, fill, opacity=0.1] {};

            \foreach\i in {1,3}
			{
		    	\path
                    (A\i) edge[very thick] (B\i)
                    (B\i) edge[very thick, opacity=0.1] (C\i)
                    (C\i) edge[very thick] (D\i)
                    (D\i) edge[very thick, opacity=0.1] (E\i)
			    ;
			}
            \foreach\i in {0,2,4}
			{
		    	\path
                    (A\i) edge[very thick, opacity=0.1] (B\i)
                    (B\i) edge[very thick] (C\i)
                    (C\i) edge[very thick, opacity=0.1] (D\i)
                    (D\i) edge[very thick] (E\i)
			    ;
			}
            \foreach\i in {1,2,3}
			{
                \pgfmathtruncatemacro\iplus{\i+1}
                \path
                    (A\i) edge[very thick] (A\iplus)
                    (B\i) edge[very thick] (B\iplus)
                    (C\i) edge[very thick] (C\iplus)
                    (D\i) edge[very thick] (D\iplus)
                    (E\i) edge[very thick] (E\iplus)
			    ;
			}
                \path
                    (A0) edge[very thick, opacity=0.1] (A1)
                    (B0) edge[very thick] (B1)
                    (C0) edge[very thick] (C1)
                    (D0) edge[very thick] (D1)
                    (E0) edge[very thick] (E1)
                    (A4) edge[very thick] (A5)
                    (B4) edge[very thick] (B5)
                    (C4) edge[very thick] (C5)
                    (D4) edge[very thick] (D5)
                    (E4) edge[very thick, opacity=0.1] (E5)
			    ;
                \path
                    (A5) edge[very thick] (B5)
                    (B5) edge[very thick, opacity=0.1] (C5)
                    (C5) edge[very thick] (D5)
                    (D5) edge[very thick, opacity=0.1] (E5)
			    ;
			\end{tikzpicture}
	\end{center}
    \caption{A $(5 \times 6)$-grid to the left and a $(5 \times 3)$-wall to the right, in which the deleted edges and vertices are still faintly visible.}
    \label{fig:gridandwallexmp}
\end{figure}
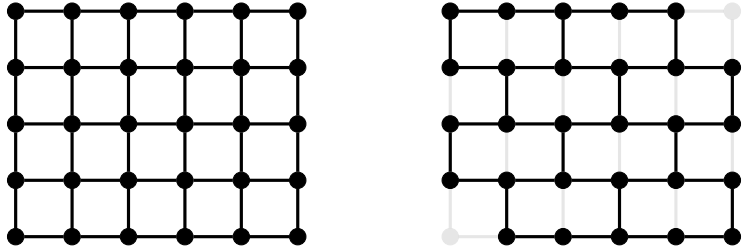

\paragraph{Walls.}
There is one more notable way in which we can connect treewidth and minors, which popped up early on in Robertson and Seymour's series of papers leading up to the GMST.
Let $n,m \in \mathbb{N}$ be two positive integers.
The \emph{$(n \times m)$-grid} is the graph $G$ with the vertex set $V(G) = [n] \times [m]$ and the edges
\begin{align*}
E(G) = \big\{ \{ (i, j) , (\ell , k) \} ~\!\colon\!~    & i, \ell \in [n], \ j,k \in [m], \text{ and } \ \\
                                                        & ( |i - \ell| = 1  \text{ and } j = k ) \text{ or } ( |j - k| = 1 \text{ and } i = \ell ) \big\} .
\end{align*}
We give a drawing of the elementary $(5 \times 6)$-grid in \Cref{fig:gridandwallexmp}.
The \emph{elementary $(n \times m)$-wall} is in turn derived from the $(n \times 2m)$-grid by deleting all edges in the following set
\begin{align*}
\big\{ \{ (i, j) , (i+1 , j) \} ~\!\colon\!~ i \in [n - 1], \ j \in [m], \text{ and } {i \not\equiv j \mod{2} } \big\}
\end{align*}
and removing all vertices of degree at most 1 in the resulting graph.
In \Cref{fig:gridandwallexmp} we show the elementary $(5 \times 3)$-wall and indicate how it is constructed from the $(5 \times 6)$-grid.
An \emph{$(n \times m)$-wall} is a subdivision of the elementary $(n \times m)$-wall.
An \emph{$n$-wall} is an $(n \times n)$-wall.

As we mentioned in the introduction, the Grid Theorem is one of the most well-known results from the Graph Minors series, and asserts that for any planar graph \(H\), every graph of a sufficiently large treewidth has an \(H\)-minor. The name of the theorem arises from the fact that the theorem usually focuses on the case when \(H\) is a grid since any \(k\)-vertex graph is a minor of a \((\mathcal{O}(k) \times \mathcal{O}(k))\)-grid~\cite{RobertsonST1994Quickly}.
For us, it is important that the Grid Theorem guarantees the presence of a large grid -- or a large wall, should we want it -- in a graph with high treewidth.

\begin{proposition}[Robertson and Seymour \cite{RobertsonS1986Grapha}]\label{thm:undirgridtheorem}
    There exists a function $g_{\ref{thm:undirgridtheorem}} \colon \mathbb{N} \rightarrow \mathbb{N}$ such that for any positive integer $k \in \mathbb{N}$, any graph $G$ with treewidth at least $g_{\ref{thm:undirgridtheorem}}(k)$ contains a $(k \times k)$-grid as a minor.
    Furthermore, if $G$ has treewidth at least $2g_{\ref{thm:undirgridtheorem}}(k)$ then $G$ contains a $(k \times k)$-wall as a subgraph.
\end{proposition}

Currently, the best known bound for the function $g_{\ref{thm:undirgridtheorem}}(k)$ is $\mathbf{O}(k^9\operatorname{poly log} k)$.
This bound is due to Chuzhoy and Tan \cite{ChuzhoyT2021Tighter}, who remark that they believe it can be further improved.

In general it is easier to work with walls since they can be found as subgraphs thanks to \Cref{lem:subcubicminorassubdivision}.
If we instead wanted to use the slightly better function for the grid, we would have to work with the subgraph that produces this grid minor.
In several cases we use a sort of hybrid between a grid and a wall which can still be found as a subgraphs; a ``mesh''.

\subsection{Meshes}\label{sec:mesh}
For several of our more technical statements we will need the notion of a \emph{mesh} that generalises both grids and walls.
Meshes allow us to use several theorems on families of intersecting paths and cycles without having to first build a wall or grid from them.

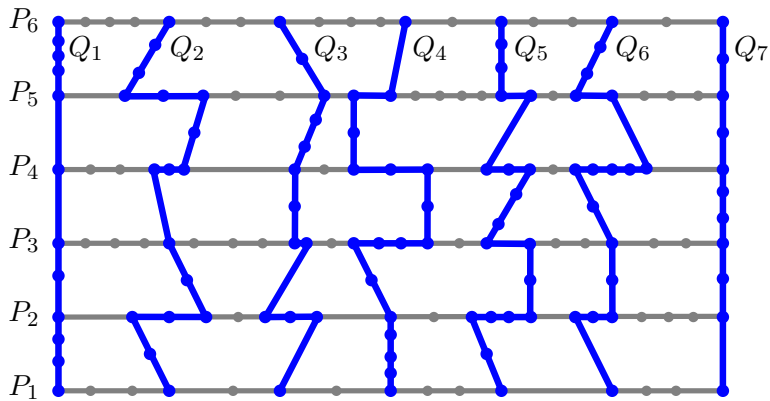
\begin{figure}[ht]
    \centering
    \begin{tikzpicture}[scale=1.5]

        \pgfdeclarelayer{background}
		\pgfdeclarelayer{foreground}
			
		\pgfsetlayers{background,main,foreground}

        \begin{pgfonlayer}{background}
            \pgftext{\includegraphics[width=6cm]{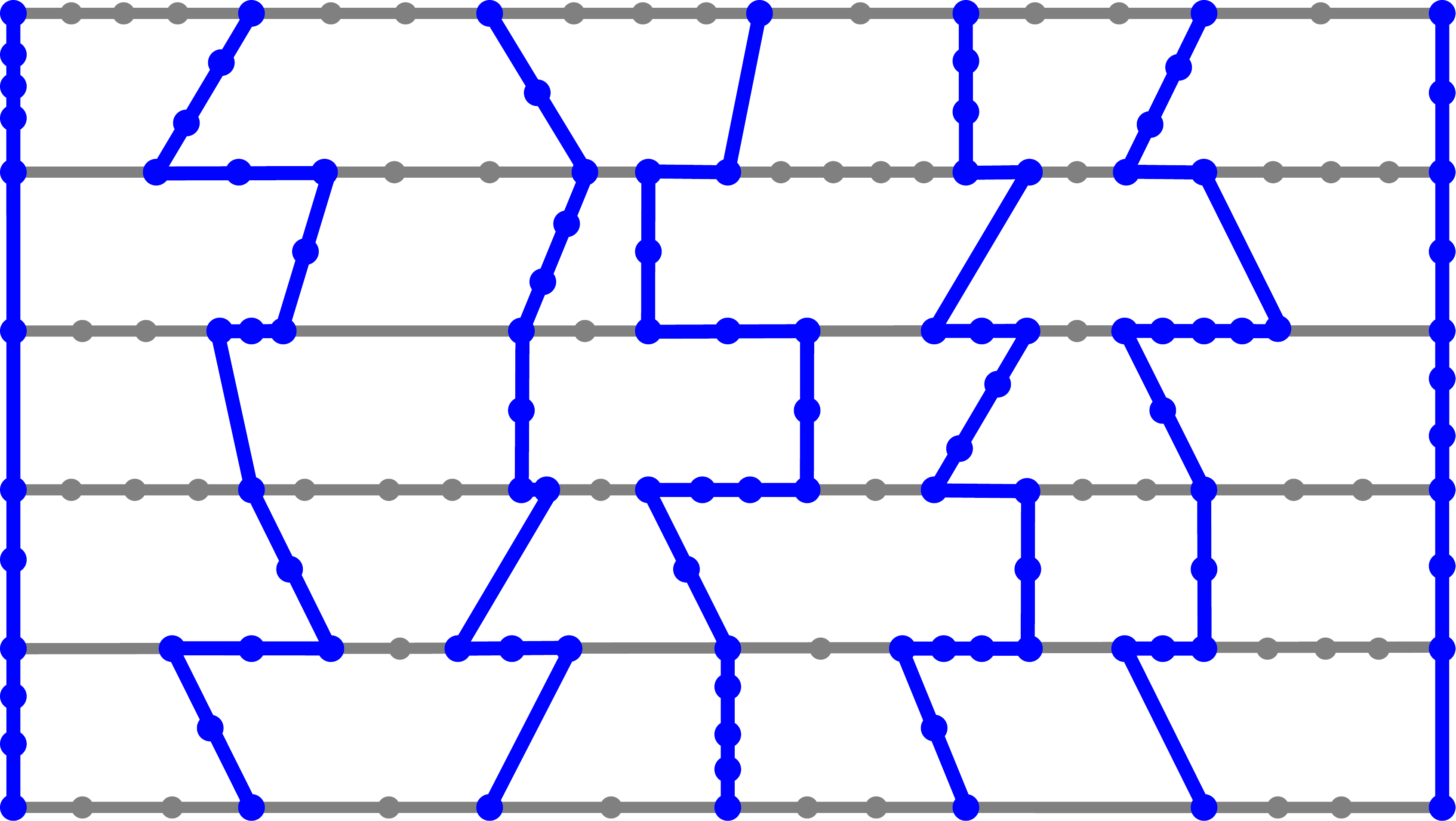}} at (C.center);
        \end{pgfonlayer}{background}
			
        \begin{pgfonlayer}{main}
        \node (C) [v:ghost] {};
            
        \end{pgfonlayer}{main}
        
        \begin{pgfonlayer}{foreground}

            \node (NAME1) at (-3.25,-1.6) [draw=none] {$P_1$};
            \node (NAME2) at (-3.25,-0.95) [draw=none] {$P_2$};
            \node (NAME3) at (-3.25,-0.3) [draw=none] {$P_3$};
            \node (NAME4) at (-3.25,0.35) [draw=none] {$P_4$};
            \node (NAME5) at (-3.25,1) [draw=none] {$P_5$};
            \node (NAME6) at (-3.25,1.65) [draw=none] {$P_6$};
            \node (NAME7) at (-2.7,1.4) [draw=none] {$Q_1$};
            \node (NAME8) at (-1.8,1.4) [draw=none] {$Q_2$};
            \node (NAME9) at (-0.525,1.4) [draw=none] {$Q_3$};
            \node (NAME10) at (0.35,1.4) [draw=none] {$Q_4$};
            \node (NAME11) at (1.25,1.4) [draw=none] {$Q_5$};
            \node (NAME12) at (2.15,1.4) [draw=none] {$Q_6$};
            \node (NAME13) at (3.2,1.4) [draw=none] {$Q_7$};
        \end{pgfonlayer}{foreground}

    \end{tikzpicture}
    \caption{An illustration of a $(6 \times 7)$-mesh with the horizontal paths drawn in gray and the vertical paths drawn in blue.}
    \label{fig:mesh}
\end{figure}

\paragraph{Meshes.}
Let $w, h$ be integers with \(w, h \ge 2\). A \emph{\((w \times h)\)-mesh} is a graph \(M\) which is the union of paths \(M = P_1 \cup \cdots \cup P_w \cup Q_1 \cup \cdots \cup Q_h\) where
\begin{itemize}
    \item $P_1, \ldots , P_w$ are pairwise vertex-disjoint, and $Q_1, \ldots, Q_h$ are pairwise vertex-disjoint,

    \item for any \(i \in [w]\) and \(j \in [h]\), the intersection $P_i \cap Q_j$ is a path,

    \item each $P_i$ is a \(V(Q_1)\)-\(V(Q_h)\)-path intersecting the paths $Q_1, \ldots , Q_h$ in the order listed, and each $Q_j$ is a \(V(P_1)\)-\(V(P_w)\)-path intersecting the paths $P_1, \ldots , P_w$ in the order listed.
\end{itemize}
(See \Cref{fig:mesh} for an example.)
The paths $P_1, \ldots, P_w$ are called the \emph{vertical paths} and the paths $Q_1, \ldots , Q_h$ are called the \emph{horizontal paths} of \(M\).
Note that the vertical and horizontal paths of a mesh \(M\) in general are not uniquely determined, but we will always assume that a mesh \(M\) comes with a fixed choice of horizontal and vertical paths. 
The union $P_1\cup P_w \cup Q_1\cup Q_h$ is a cycle called the \emph{perimeter} of $M$.
An \emph{$n$-mesh} is an $(n\times n)$-mesh. 
We observe that according to this definition both grids and walls are meshes.

A \(w' \times h'\) mesh \((\mathcal{P}', \mathcal{Q}')\) is a \emph{submesh} of \((\mathcal{P}, \mathcal{Q})\) if every path in \(\mathcal{P}'\) is a subpath of a path in \(\mathcal{P}\), and every path in \(\mathcal{Q}'\) is a subpath of a path in \(\mathcal{Q}\).

\paragraph{Modelling a minor.}
In some situations we are interested in how exactly a minor is found in the graph.
For example when we find a $K_t$-minor linked to a wall.
For this purpose we also introduce a notion of minor models.

A function $\varphi \colon V(H) \rightarrow 2^{V(G)}$ is called a \emph{minor model (of $H$ in $G$)} (or simply a \emph{model}) if $\varphi(V(H))$ is a set of pairwise disjoint vertex sets, each inducing a connected subgraph of $G$, and for each edge $uv \in E(H)$ there exists an edge $ab \in E(G)$ with $a \in \varphi(u)$ and $b \in \varphi(v)$.
Given this definition, we observe that a graph $H$ is a minor of a graph $G$ if and only if there exists a minor model of $H$ in $G$.

\paragraph{Minor models controlled by meshes.}
Let \(G\) be a graph, and let \(M\) be a \((w \times h)\)-mesh in \(G\) with vertical paths \(P_1, \ldots, P_w\) and horizontal paths \(Q_1, \ldots, Q_h\). In this setting, for any set \(Z \subseteq V(G)\) with \(|Z| < \min\{w, h\}\),
at least one vertical path is disjoint from \(Z\), and at least one horizontal path is disjoint from \(Z\).
Since every vertical path intersects every horizontal path, there exists a component of \(G - Z\) containing all vertical and horizontal paths of the mesh that are disjoint from \(Z\).
Now, let \(H\) be a \(t\)-vertex graph with \(t \le \min\{w, h\}\), and let \(\varphi\) be a minor model of \(H\) in the graph \(G\).
We say that \(\varphi\) is \emph{controlled} by \(M\) if for every set \(Z \subseteq V(G)\) with \(|Z| < t\), and every \(x \in V(H)\),
if the branch set \(\varphi(x)\) is disjoint from \(Z\), then it is contained in the same component of \(G - Z\) as the vertical and horizontal paths of \(M\) which are disjoint from \(Z\).
By Menger's Theorem this is equivalent to saying that for each \(x \in V(H)\) and each vertical or horizontal path \(P\) of the mesh, there exist \(t\) internally disjoint \(\varphi(x)\)--\(V(P)\) paths \(R_1, \ldots, R_t\) in \(G\).
(We allow \(R_1, \ldots, R_t\) to be the same path of length \(\le 1\), if the sets \(\varphi(x)\) and \(V(P)\) intersect or there exists an edge between them.)

\subsection{Highly linked sets and tangles}\label{sec:tangles}
While walls represent amazing witnesses for large treewidth which impose a lot of additional structure onto a graph, they are also relatively precise in nature and for this reason not always good candidates for working with treewidth in a more abstract way.
Two more abstract alternatives are given by so-called ``highly linked sets'' which are sets that cannot be spread out over several bags of any tree-decomposition of small width, and ``tangles'' which constitute possibly the most abstract and versatile obstruction for treewidth.

\paragraph{The set of separations.}
Let $G$ be a graph and $k$ be a positive integer.
We denote by $\mathcal{S}_k(G)$ the collection of all separations $(A,B)$ of order $< k$ in $G$.

\paragraph{Tangles.}
An \emph{orientation} of $\mathcal{S}_k(G)$ is a set $\mathcal{O}$ such that for all $(A,B)\in\mathcal{S}_k(G)$ exactly one of $(A,B)$ and $(B,A)$ belongs to $\mathcal{O}$. 
A \emph{tangle} of order $k$ in $G$ is an orientation $\mathcal{T}$ of $\mathcal{S}_k(G)$ such that for all $(A_1,B_1),(A_2,B_2),(A_3,B_3)\in\mathcal{T}$, it holds that $G[A_1]\cup G[A_2]\cup G[A_3]\neq G$.
If $\mathcal{T}$ is a tangle and $(A,B)\in\mathcal{T}$ we call $A$ the \emph{small side} and $B$ the \emph{big side} of $(A,B)$.

Let $G$ be a graph and $\mathcal{T}$ and $\mathcal{D}$ be tangles of $G$.
We say that $\mathcal{D}$ is a \emph{truncation} of $\mathcal{T}$ if $\mathcal{D}\subseteq\mathcal{T}$.
\medskip

Let $r \in \mathbb{N}$ with $r\geq 3$, let $G$ be a graph, and $M$ be an $r$-mesh in $G$.
Let $\mathcal{T}_M$ be the orientation of $\mathcal{S}_r$ such that for every $(A,B)\in\mathcal{T}_M$, the set $B\setminus A$ contains the vertex set of both a horizontal and a vertical path of $M$, we call $B$ the \emph{$M$-majority side} of $(A,B)$.
Then $\mathcal{T}_M$ is the tangle \emph{induced} by $M$.
If $\mathcal{T}$ is a tangle in $G$, we say that $\mathcal{T}$ \emph{controls} the wall $M$ if $\mathcal{T}_M$ is a truncation of $\mathcal{T}$.

Let $G$ and $H$ be graphs as well as $\mathcal{T}$ be a tangle in $G$.
We say that a minor model $\mu$ of $H$ in $G$ is \emph{controlled} by $\mathcal{T}$ if there do not exist a separation $(A,B)\in\mathcal{T}$ of order $<|V(H)|$ and $x\in V(H)$ such that $\mu(x)\subseteq A\setminus B$.

Notice that, if a minor model $\mu$ of $H$ in $G$ is controlled by some mesh $M$ in $G$, then $\mu$ is also controlled by $\mathcal{T}_M$.

\paragraph{Highly linked sets.}
Let $\alpha \in [2/3, 1)$.
Moreover, let $G$ be a graph and $X \subseteq V(G)$ be a vertex set. 
A set $S \subseteq V(G)$ is said to be an \emph{$\alpha$-balanced separator} for $X$ if for every component $C$ of $G - S$ it holds that $|V(C) \cap X| \leq \alpha|X|$. 
Let $k$ be a non-negative integer.
We say that $X$ is a \emph{$(k, \alpha)$-linked set} of $G$ if there is no $\alpha$-balanced separator of size at most $k$ for $X$ in $G$.
\medskip

Given a $(3k, \alpha)$-linked set $X$ of $G$ we define $$\mathcal{T}_{X} \coloneqq \{ (A, B) \in \mathcal{S}_{k+1}(G) ~\!\colon\!~ |X \cap B| > \alpha|X| \}.$$ 
It is not hard to see that $\mathcal{T}_{S}$ is a tangle of order $k+1$ in $G$.

We are also interested in the computational aspects of the GMST, which includes an algorithm to compute the structural decomposition as a whole.
Usually the corresponding proofs are formulated in terms of tangles.
Tangles, however, constitute objects which are somewhat hard to handle computationally since, a priori, one would need to have knowledge of the entire set $\mathcal{S}_k(G)$.
Since this set can be of order $|V(G)|^k$ this is not desirable.
Instead, it is more convenient to work with highly linked sets as described above.
Indeed, since highly linked sets induce tangles it is possible to phrase most of the theory of graph minors (as long one does not worry about a linear blowup for the numbers) in terms of highly linked sets.

Finally, we need an algorithmic way to find, given a highly linked set, a large wall whose tangle is a truncation of the tangle induced by the highly linked set.
This is done in \cite{ThilikosW2024Excluding} by algorithmatising a proof of Kawarabayashi et al.\ from \cite{KawarabayashiTW2021Quickly}. 

\begin{proposition}[Thilikos and Wiederrecht \cite{ThilikosW2024Excluding} (see Theorem 4.2.)]\label{thm:algogrid}
Let $k\geq 3$ be an integer and $\alpha\in [2/3,1)$.
There exist universal constants $c_1, c_2\in\mathbb{N}\setminus\{ 0\}$, and an algorithm that, given a graph $G$ and a $(c_1k^{20},\alpha)$-linked set $X\subseteq V(G)$ computes in time $2^{\mathbf{O}(k^{c_2})}|V(G)|^2|E(G)|\log(|V(G)|)$ a $k$-wall $W\subseteq G$ such that $\mathcal{T}_W$ is a truncation of $\mathcal{T}_X$.
\end{proposition}

\section{Drawings, embeddings, and \texorpdfstring{$\Sigma$}{Sigma}-decompositions}\label{sec:decompositions}
In the following we provide an extensive introduction to the definitions found in the work by Kawarabayashi et al.
Most of these definitions are very technical.
While it is certainly possible to present the GMST with simpler terms, these concepts facilitate the proof in many ways by allowing us to localise a lot of the arguments.
Moreover, these definitions provide a precise encoding of ``near embeddings'' which allows us to work with our graph mostly as if it was actually embedded in some surface.

\paragraph{Surfaces.}
By a \emph{surface} we mean a compact $2$-dimensional manifold with or without boundary.

Given a pair $(\mathsf{h}, \mathsf{c}) \in \mathbb{N} \times [0,2]$ we define $\Sigma^{(\mathsf{h}, \mathsf{c})}$ to be the two-dimensional surface without boundary created from the sphere by adding $\mathsf{h}$ handles and $\mathsf{c}$ crosscaps (see \cite{MoharT2001Graphs} for more details).
If $\mathsf{c} = 0$ the surface $\Sigma^{(\mathsf{h}, \mathsf{c})}$ is an \emph{orientable} surface, otherwise it is called \emph{non-orientable}.
By Dyck's theorem \cite{Dyck1888Beitraege,FrancisW1999Conways}, two crosscaps are equivalent to a handle in the presence of a third crosscap.
Thus the notation $\Sigma^{(\mathsf{h}, \mathsf{c})}$ is sufficient to capture all two-dimensional surfaces without boundary.\medskip

In a first step, we introduce a framework to talk about ``drawings'' of graphs in a surface.
The final notion of $\Sigma$-decompositions, as introduced by Kawarabayashi et al., will make use of drawings as a basis but impose additional structure and information on top of them.

\paragraph{Drawing in a surface.}
A \emph{drawing} (with crossings) \emph{in a surface $\Sigma$} is a triple $\Gamma=(U,V,E)$ such that
\begin{itemize} 
\item $V$ and $E$ are finite, 
\item $V\subseteq U\subseteq\Sigma$, 
\item $V\cup\bigcup_{e\in E}e=U$ and $V\cap (\bigcup_{e\in E}e)=\emptyset$, 
\item for every $e\in E$, it holds that $e=h([0,1]_{\mathbb{R}})\setminus\{h(0),h(1)\}$, where $h\colon[0,1]_{\mathbb{R}}\rightarrow U$ is a homeomorphism onto its image with $h(0),h(1)\in V$, and 
\item if $e,e'\in E$ are distinct, then $|e\cap e'|$ is finite.
\end{itemize}
We call the set $V$, sometimes referred to by $V(\Gamma)$, the \emph{vertices of $\Gamma$} and the set $E$, referred to by $E(\Gamma)$, the \emph{edges of $\Gamma$}.
If $G$ is graph and $\Gamma=(U,V,E)$ is a drawing with crossings in a surface $\Sigma$ such that $V$ and $E$ naturally correspond to $V(G)$ and $E(G)$ respectively, we say that $\Gamma$ is a \emph{drawing of $G$ in $\Sigma$ (possibly with crossings)}.
In case no two edges of $E(\Gamma)$ have a common point, we say that $\Gamma$ is a \emph{drawing of $G$ in $\Sigma$ without crossings}, this is also sometimes referred to as an \emph{embedding of $G$ in $\Sigma$}.
In this last case, the connected components of $\Sigma\setminus U$ are the \emph{faces} of $\Gamma$.

Let $G$ be a graph and $\Gamma$ be a drawing of $G$ in some surface $\Sigma$.
Now let $H$ be a subgraph of $G$.
Then there exists a drawing $\Gamma'$ of $H$ in $\Sigma$ where $V(\Gamma')\subseteq V(\Gamma)$ and $E(\Gamma')\subseteq E(\Gamma)$.
We call $\Gamma'$ the \emph{subdrawing of $\Gamma$ induced by $H$} or, if $\Gamma$ is understood from the context, the \emph{induced drawing of $H$}.
\medskip

The next step is to properly introduce the concept of ``$\Sigma$-decompositions''.
The main difference between a drawing (without crossings) of a graph and a ``$\Sigma$-decomposition'' where $\Sigma$ is the surface where our graph ``near embeds'' is, that we have to account for the possibility that it is not possible to draw the entire graph in the surface without crossings.
Similar to Wagner's Theorem on $K_5$-minor-free graphs \cite{Wagner1937Ueber} mentioned earlier, graphs without large clique-minors can be created from a small collection of ``base classes'' by using the notion of clique sums.
This means that, behind the boundary of a clique sum of order at most three one may find an arbitrarily large additional part of the graph, while the boundary itself lives on a face in $\Sigma$.
Moreover, we need to account for the presence of ``vortices'' which can be seen as areas, or faces, of the ``near embedding'' where the rest of the graphs attaches via a sequence of clique-sums in a more non-trivial manner compared to the $3$-clique-sums mentioned above.
The purpose of $\Sigma$-decompositions is to not only provide a data structure to deal with these objects while maintaining some sense of the embedding, but to also provide a way to precisely define the notions mentioned above.

\begin{figure}[ht]
    \centering
    \scalebox{1}{
    \begin{tikzpicture}[scale=1]

        \pgfdeclarelayer{background}
		\pgfdeclarelayer{foreground}
			
		\pgfsetlayers{background,main,foreground}
			
        \begin{pgfonlayer}{main}
        \node (C) [v:ghost] {};

            \pgftext{\includegraphics[width=8cm]{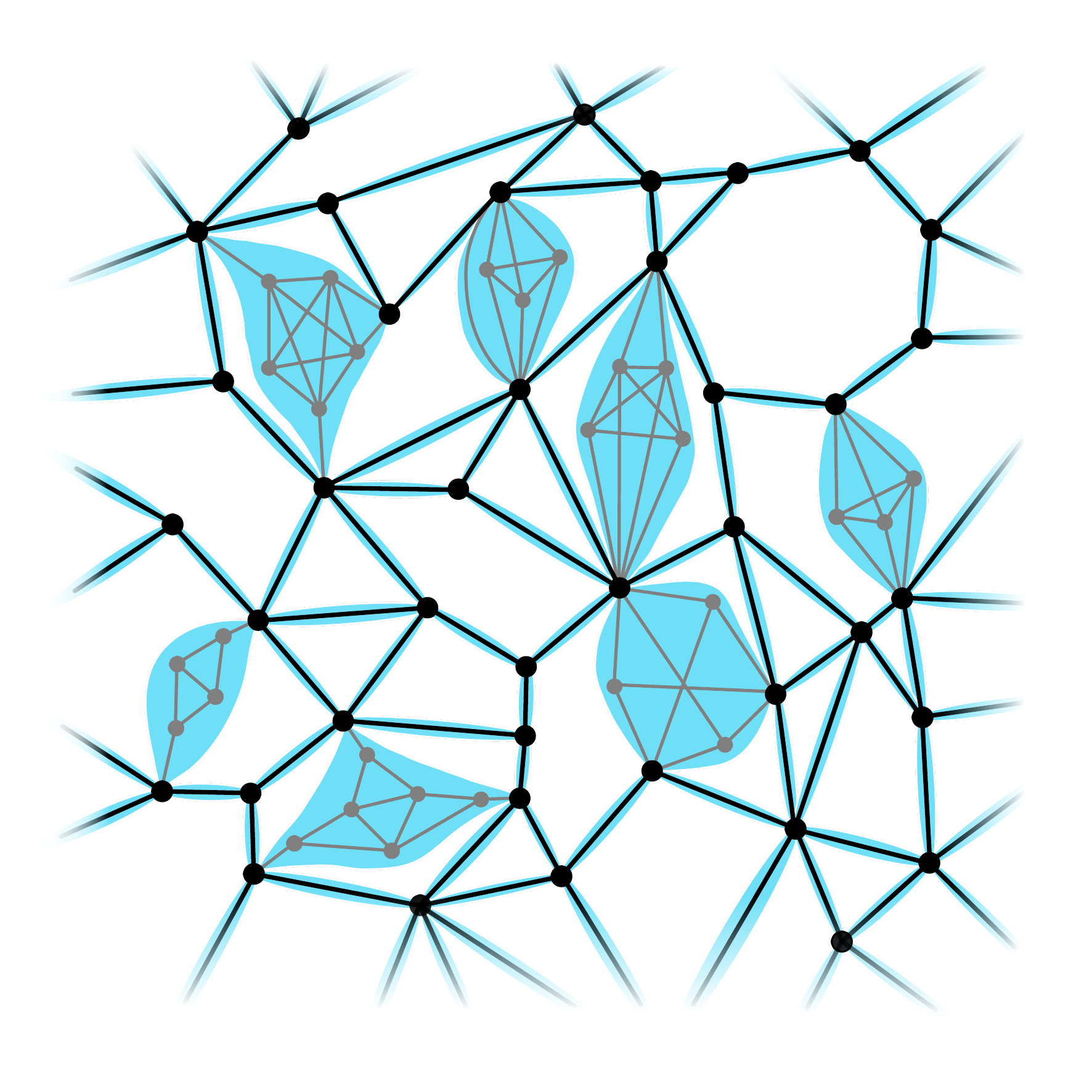}} at (C.center);
            
        \end{pgfonlayer}{main}
        
        \begin{pgfonlayer}{foreground}
        \end{pgfonlayer}{foreground}

        \begin{pgfonlayer}{background}
        \end{pgfonlayer}{background}
        
    \end{tikzpicture}
    }
    \caption{Part of a $\Sigma$-decomposition $\delta$ of some graph $G$. The blue shapes are the cells of $\delta$, the black vertices are its nodes, and the gray vertices are those fully contained within the cells of $\delta$. Notice that in this example, the boundary of each cell has size at most three.}
    \label{fig:sigma_decomposition}
\end{figure}

\paragraph{$\Sigma$-decompositions.}
Let $\Sigma$ be a surface. 
A \emph{$\Sigma$-decomposition} of a graph $G$ is a pair $ \delta=(\Gamma,\mathcal{D})$, where $\Gamma$ is a drawing of $G$ in $\Sigma$ with crossings, and $\mathcal{D}$ is a collection of closed disks, each a subset of $\Sigma$ such that 
\begin{enumerate}[label=\textit{\roman*})]
\item the disks in $\mathcal{D}$ have pairwise disjoint interiors, 
\item the boundary of each disk in $\mathcal{D}$ intersects $\Gamma$ in vertices only, 
\item if $ \Delta_1, \Delta_2\in\mathcal{D}$ are distinct, then $ \Delta_1\cap \Delta_2\subseteq V(\Gamma)$, and
\item every edge of $\Gamma$ belongs to the interior of one of the disks in $\mathcal{D}$. 
\end{enumerate} 
Let $N$ be the set of all vertices of $\Gamma$ that do not belong to the interior of any of the disks in $\mathcal{D}$. 
We refer to the elements of $N$ as the \emph{nodes} of $\delta$. 
For every $\Delta\in\mathcal{D}$, we refer to the set $\Delta-N$ as a \emph{cell} of $\Delta$. 
We denote the set of cells of $\delta$ by $C(\delta)$. 
For a cell $c\in C(\delta)$ the set of nodes that belong to the closure of $c$ is denoted by $\widetilde{c}$. 
For a cell $c\in C(\delta)$ we define $\sigma_{\delta}(c)$, or $ \sigma(c)$ in case $\delta$ is clear from the context, to be the subgraph of $G$ consisting of those vertices and edges which are drawn in the closure of $c$ by $\Gamma$.
See \cref{fig:sigma_decomposition} for an illustration of a (part of a) $\Sigma$-decomposition.
 
We define $\pi_{\delta}\colon N(\delta)\rightarrow V(G)$ to be the mapping which assigns to every node in $N(\delta)$ the corresponding vertex of $G$.

Isomorphisms between two $\Sigma$-decompositions are defined in the natural way. That is given $\Sigma$-decompositions $\delta=(\Gamma,\mathcal{D})$ and $\delta'=(\Gamma',\mathcal{D}')$, the drawing $\Gamma=(U,V,E)$ is mapped to a drawing $\Gamma'=(U',V',E')$ where the elements of $V$ are in bijection with the elements of $V'$, similarly for $E$ and $E'$ such that the map between the corresponding graphs is a graph isomorphism, and the elements of $\mathcal{D}$ are mapped to the disks in $\mathcal{D}'$ while agreeing with the map between $\Gamma$ and $\Gamma'$.

Let $G$ be a graph, let $\Sigma$ be a surface, and let $\delta=(\Gamma,\mathcal{D})$ be a $\Sigma$-decomposition for $G$.
Moreover, let $H$ be a subgraph of $G$ and let $\Gamma'$ be the subdrawing of $\Gamma$ induced by $H$.
Then there exists a $\Sigma$-decomposition $\delta'=(\Gamma',\mathcal{D}')$ of $H$ where $\mathcal{D}'\subseteq \mathcal{D}$ is obtained as follows.
For every cell $c\in C(\delta)$, if there is no edge or vertex of $H$ drawn within $c$, then the closure $\Delta_c\in\mathcal{D}$ of $c$ is discarded and does not belong to $\mathcal{D}'$.
Otherwise, the closure $\Delta_c\in\mathcal{D}$ of $c$ does belong to $\mathcal{D}$.
We call $\delta'$ a \emph{$\Sigma$-subdecomposition} of $\delta$ or the \emph{$\Sigma$-decomposition of $H$ inherited from $\delta$}.

Let $\delta$ and $\delta'$ be $\Sigma$-decompositions of $G$ in a surface $\Sigma$ and let $H$ be $G - \bigcup \{ \sigma(c) ~\!\colon\!~ c \in C(\delta) \text{ is a vortex} \}$.
We say that a $\delta'$ is a \emph{refinement of $\delta$} if the $\Sigma$-decomposition of $H$ inherited from $\delta$ is also the $\Sigma$-decomposition of $H$ inherited from $\delta'$.

\begin{figure}[ht]
    \centering
    \scalebox{1}{
    \begin{tikzpicture}[scale=1]

        \pgfdeclarelayer{background}
		\pgfdeclarelayer{foreground}
			
		\pgfsetlayers{background,main,foreground}
			
        \begin{pgfonlayer}{main}
        \node (C) [v:ghost] {};

            \pgftext{\includegraphics[width=7cm]{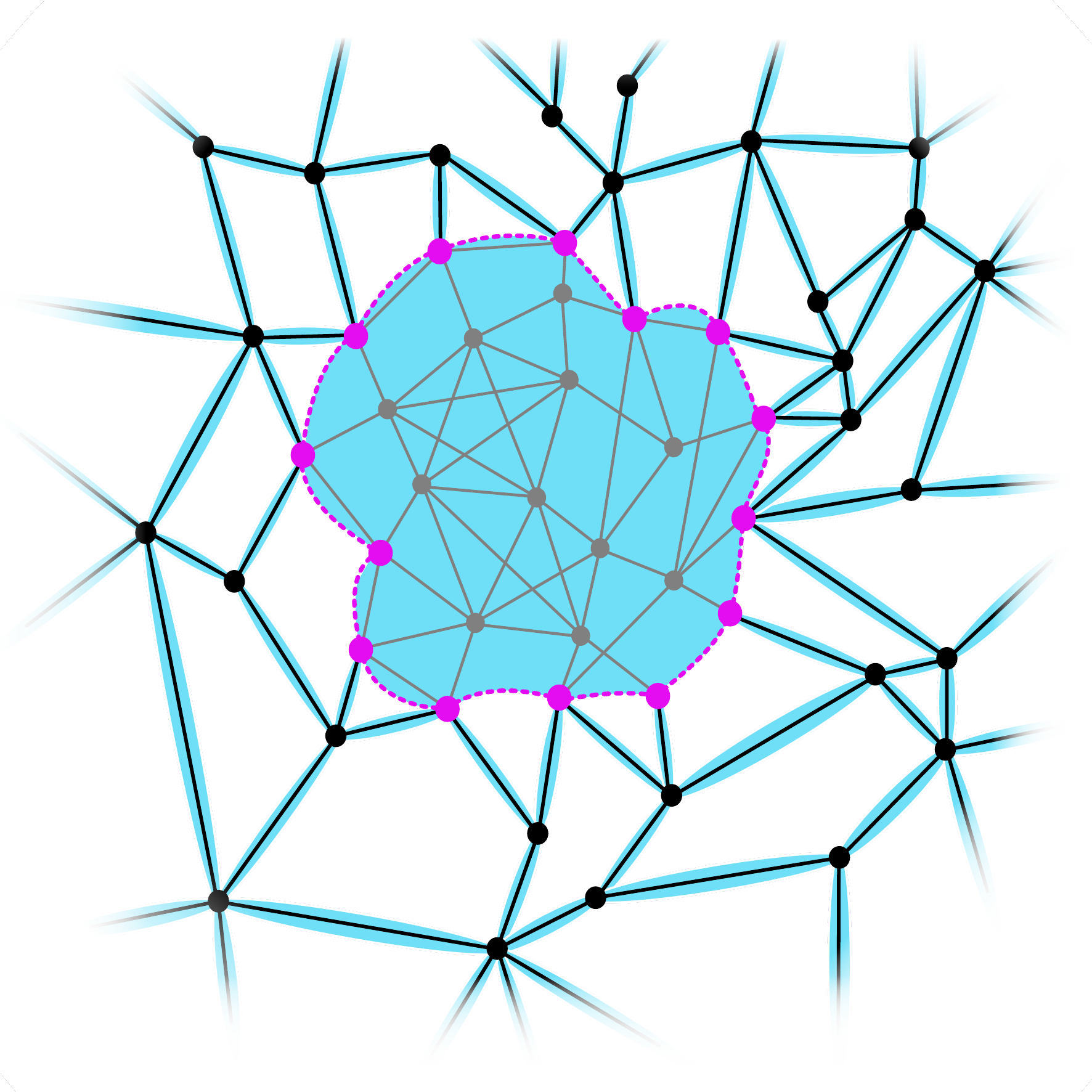}} at (C.center);
            
        \end{pgfonlayer}{main}
        
        \begin{pgfonlayer}{foreground}
        \end{pgfonlayer}{foreground}

        \begin{pgfonlayer}{background}
        \end{pgfonlayer}{background}
        
    \end{tikzpicture}
    }
    \caption{A vortex cell and its immediate vicinity in a $\Sigma$-decomposition.}
    \label{fig:vortex}
\end{figure}

Another way to interpret the definition of $\Sigma$-decompositions is to see the cells of a $\Sigma$-decomposition $\delta$ as hyperedges of a hypergraph with vertex set $N(\delta)$ where $\widetilde{c}$ is the set of vertices incident with $c$ for each $c\in C(\delta)$.
In this interpretation, the hypergraph can be embedded in $\Sigma$ such that hyperedges only intersect in common vertices and all vertices are drawn in the boundaries of their hyperedges.
Recall \cref{fig:sigma_decomposition} and also consider \cref{fig:vortex} for indicative illustrations.
In both figures, the \textcolor{CornflowerBlue}{blue} areas mark the cells which become the hyperedges of the hypergraph described above.
Please note that, for the sake of a simpler presentation, the graphs in the interiors of the cells are chosen to be somewhat small.
In general no such restriction exists, the graphs within cells may be arbitrarily complex.

\paragraph{Vortices}
Let $G$ be a graph, $\Sigma$ be a surface and $ \delta=(\Gamma,\mathcal{D})$ be a $\Sigma$-decomposition of $G$.
A cell $c \in C(\delta)$ is called a \emph{vortex} if $|\widetilde{c}|\geq 4$.
Moreover, we call $\delta$ \emph{vortex-free} if none of its cells is a vortex.\medskip

See \cref{fig:vortex} for an illustration of a part of some $\Sigma$-decomposition which includes a vortex.

The origin of the threshold for the boundary size of a vortex to be four can be found in the so-called ``\textit{Two Paths Theorem}'' as presented in \cref{sec:societies}.
In simple terms, one need at least four vertices in the boundary of a cell to route a ``cross'', i.e.\ two vertex-disjoint paths with antipodal endpoints (see \cref{fig:society_cross}) through the interior of the cell.
This means that any ``local'' obstruction to embeddability in a surface can be fully separated from the embedded part by deleting at most three vertices while these three boundary vertices of the separation may be completed to a triangle which then, in a reduced graph, bounds a face in the embedding.
This intuition is the basis of the notion of ``flatness'' and is elaborated on further in the next subsection.

\subsection{Societies and flatness}\label{sec:societies}
Central to understanding the differentiation between ``normal'' cells and vortices in a $\Sigma$-decomposition is the notion of flatness which is linked to the so called ``Two Paths Theorem'' as mentioned above.
This theorem can be understood as a bridge between the topological arguments involved in the theory of graph minors, and its combinatorial aspects.

Let $G$ be a graph and let $s_1,s_2,t_1,t_2\in V(G).$
The \textsc{Two Disjoint Paths Problem} (\textsc{TDPP}) with \emph{terminals} $s_1,s_2,t_1,t_2$ is the decision problem asking for the existence of two vertex-disjoint paths $P_1$ and $P_2$ in $G$ such that $P_i$ has endpoints $s_i$ and $t_i$ for both $i\in[2]$.

\textsc{TDPP} famously admits a polynomial time algorithm.
This algorithm is based on a full characterisation of the \textsc{yes}-instances of the problem and it is exactly this characterisation that gives us a notion of ``flatness'' and introduces topology to the theory of graph minors.
The statement of the Two Paths Theorem we present here is based on the notion of ``societies'' which have already played a crucial role in the proof of Kawarabayashi et al.\ for the GMST \cite{KawarabayashiTW2021Quickly}, though they also already appear in the work of Robertson and Seymour (see for example \cite{RobertsonS1990Graph}).

To provide a rough intuition, a society is a graph $G$ together with a cyclic ordering $\Omega$ of a vertex set $V(\Omega)\subseteq V(G)$.
This cyclic ordering is meant to indicate that we are interested in some particular topological properties of the set $V(\Omega)$.
In most cases, for a society, we are interested in embedding the vertices in $V(\Omega)$ on the boundary of a disk $\Delta$, respecting the ordering $\Omega$, and drawing the rest of $G$ in the interior of $\Delta$.
This is particularly useful as it allows us to maintain a $\Sigma$-decomposition $\delta$ in the background, carve out a subgraph $G'$ of the graph which is currently drawn in a disk $\Delta$, and then further refine how $\delta$ decomposes $G'$ in $\Delta$.
In particular, this further refinement might add new handles or crosscap and thereby also edit $\Sigma$ itself.
Relating to the language of $\Sigma$-decompositions, each cell $c$, including vortex cells, defines as so-called ``vortex society'' which is given by the cyclic ordering of the vertices drawn in the boundary of $c$ together with the graph $\sigma(c)$.
From this perspective, we will use the notion of societies to iteratively refine the vortices of an initial $\Sigma$-decomposition.

\begin{figure}[ht]
    \centering
    \begin{tikzpicture}

        \pgfdeclarelayer{background}
		\pgfdeclarelayer{foreground}
			
		\pgfsetlayers{background,main,foreground}

        \begin{pgfonlayer}{background}
            \pgftext{\includegraphics[width=10cm]{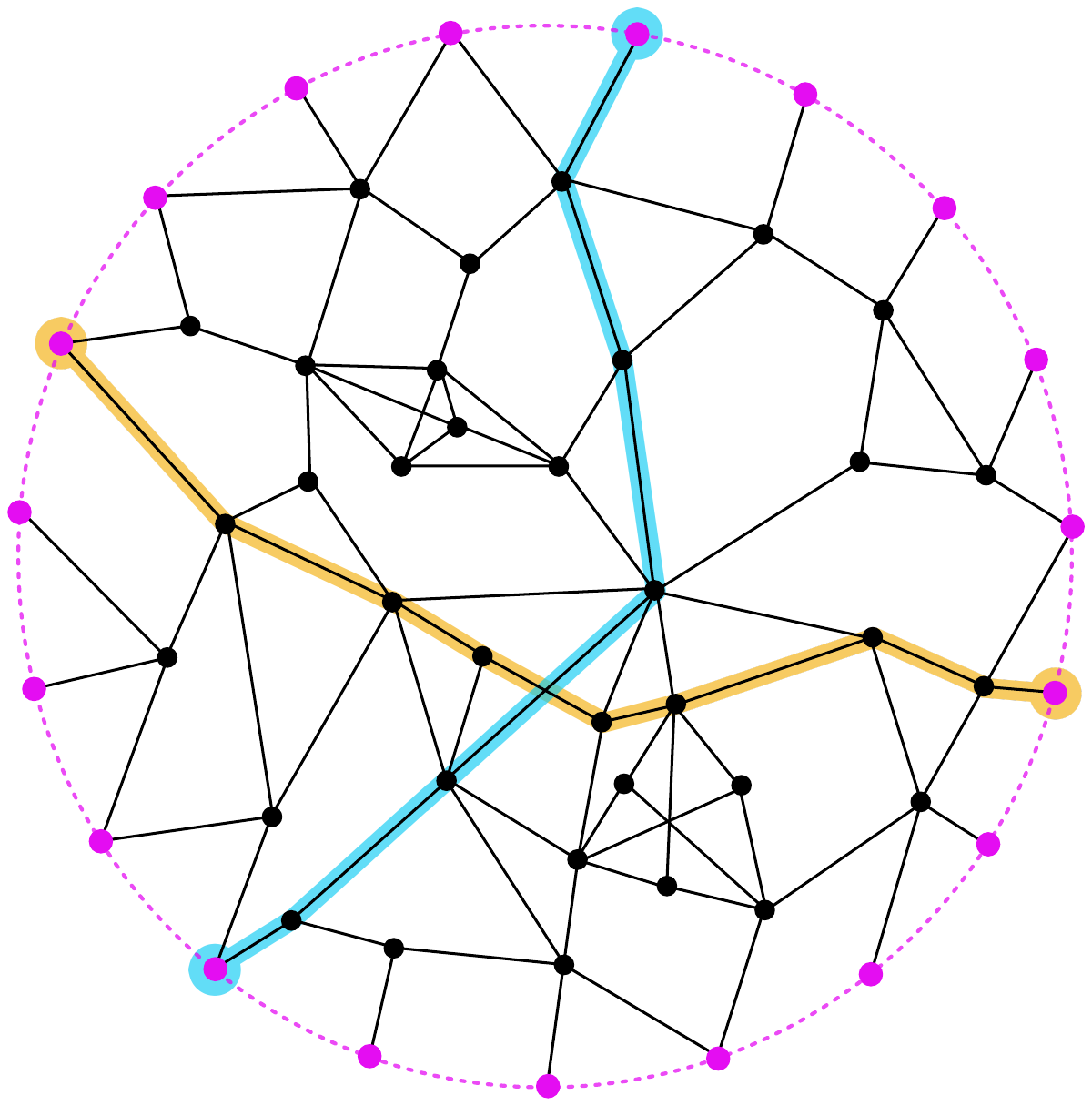}} at (C.center);
        \end{pgfonlayer}{background}
			
        \begin{pgfonlayer}{main}
        \node (C) [v:ghost] {};

            \node (s_1) [v:ghost,position=155:35.5mm from C] {$s_1$};
            \node (t_1) [v:ghost,position=345:39mm from C] {$t_1$};
            \node (s_2) [v:ghost,position=75:37.5mm from C] {$s_2$};
            \node (t_2) [v:ghost,position=235:35.6mm from C] {$t_2$};

        \end{pgfonlayer}{main}
        
        \begin{pgfonlayer}{foreground}
        \end{pgfonlayer}{foreground}

    \end{tikzpicture}
    \caption{A society $(G,\Omega)$ together with a cross.
    The set $V(\Omega)$ is depicted in magenta.}
    \label{fig:society_cross}
\end{figure}

\paragraph{Societies.}
Let $\Omega$ be a cyclic ordering of the elements of some set which we denote by $V(\Omega)$.
A \emph{society} is a pair $(G,\Omega),$ where $G$ is a graph and $\Omega$ is a cyclic ordering with $V(\Omega)\subseteq V(G)$.
A \emph{cross} in a society $(G,\Omega)$ is a pair $(P_1,P_2)$ of vertex-disjoint paths in $G$ such that $P_i$ has endpoints $s_i,t_i\in V(\Omega)$ and is otherwise disjoint from $V(\Omega)$, and the vertices $s_1,s_2,t_1,t_2$ occur in $\Omega$ in the order listed.

Let $(G, \Omega)$ be a society.
For a given set $S \subseteq V(\Omega)$ a vertex $s \in S$ is an \emph{endpoint} of $S$ if there exists a vertex $t \in V(\Omega) \setminus S$ that immediately precedes or succeeds $s$ in $\Omega$.
We call $S$ a \emph{segment} of $\Omega$ if $S$ has two or less endpoints.
For $s,t \in V(\Omega)$ we denote by $s\Omega t$ the uniquely determined segment with the first vertex $s$ and the last vertex $t$ according to $\Omega$.
\medskip

Notice, that with the terminology above, we may translate an instance $(G,s_1,s_2,t_1,t_2)$ of \textsc{TDPP} to a society $(G,\Omega)$ where $V(\Omega)=\{ s_1,s_2,t_1,t_2\}$ and these four vertices occur in $\Omega$ in the order listed.
Then $(G,s_1,s_2,t_1,t_2)$ is a \textsc{yes}-instance of \textsc{TDPP} if and only if the society $(G,\Omega)$ has a cross.
See \cref{fig:society_cross} for an illustration of a society with a cross.

\paragraph{Renditions of societies.}
Let $(G,\Omega)$ be a society, and let $\Sigma$ be a surface with one boundary component $B$ homeomorphic to the unit circle.
A \emph{rendition} of $G$ in $\Sigma$ is a $\Sigma$-decomposition $\rho$ of $G$ such that the image under $\pi_{\rho}$ of $N(\rho) \cap B$ is $V(\Omega)$ and $\Omega$ is one of the two cyclic orderings of $V(\Omega)$ defined by the way the points of $\pi_{\delta}(V(\Omega))$ are arranged in the boundary $B$.
\medskip

\begin{figure}[ht]
    \centering
    \begin{tikzpicture}

        \pgfdeclarelayer{background}
		\pgfdeclarelayer{foreground}
			
		\pgfsetlayers{background,main,foreground}

        \begin{pgfonlayer}{background}
            \pgftext{\includegraphics[width=6cm]{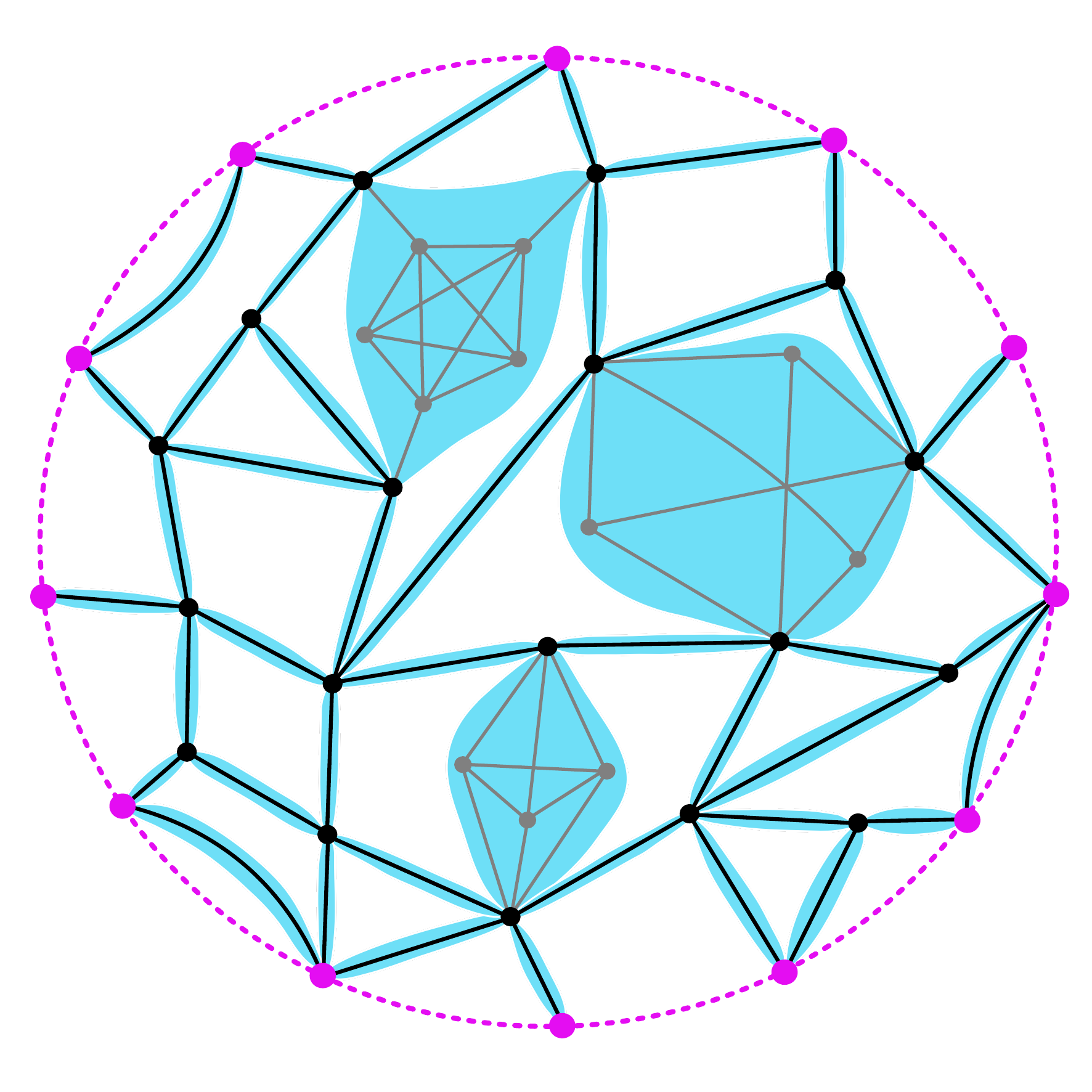}} at (C.center);
        \end{pgfonlayer}{background}
			
        \begin{pgfonlayer}{main}
        \node (C) [v:ghost] {};
            
        \end{pgfonlayer}{main}
        
        \begin{pgfonlayer}{foreground}
        \end{pgfonlayer}{foreground}

    \end{tikzpicture}
    \caption{A vortex-free rendition of a society.}
    \label{fig:vortex_free_rendition}
\end{figure}

With these definitions at hand we are able to state the Two Paths Theorem in the general context of the GMST as follows.
For an illustration of a vortex-free rendition and the absence of a cross in a society $(G,\Omega),$ even when $G$ is non-planar, see \cref{fig:vortex_free_rendition}.

\begin{proposition}[Two Paths Theorem, \cite{Jung1970Verallgemeinerung,Seymour1980Disjoint,Shiloach1980Polynomial,Thomassen19802Linked,RobertsonS1990Graph}]\label{prop:TwoPaths}
A society $(G,\Omega)$ has no cross if and only if it has a vortex-free rendition in a disk.
\end{proposition}

In particular these outcomes can be found algorithmically.
This is implicit in the way Robertson and Seymour prove this theorem \cite{RobertsonS1990Graph}, with explicit algorithms with an almost linear runtime given by Tholey \cite{Tholey2006Solving,Tholey2009Improved}, and an algorithm with linear runtime is presented by Kawarabayashi, Li, and Reed \cite{KawarabayashiLR2015Connectivity} (see also \cite[Theorem 7.5]{KawarabayashiTW2018New}).

\subsection{Transactions}\label{sec:transaction}
Though vortices represent non-planar parts of a $\Sigma$-decomposition that cannot easily be separated from the remaining embedded graph, the vortices in the GMST are still fairly simple parts of the graph.
This is because the vortices appearing in the statement of the GMST no longer allow many disjoint paths to move through them.
Of course the vortices we encounter throughout our proof are not guaranteed to have this property and thus we must introduce tools to discuss large sets of paths interacting with societies, particularly around vortices.

\paragraph{Transactions and their types.}
Given a society $(G,\Omega)$, a \emph{transaction} $\mathcal{P}$ in $(G,\Omega)$ is an $A$-$B$-linkage in $G$ such that $A,B$ are disjoint segments of $\Omega$.
The inclusion-wise minimal segments $X$ and $Y$ of $\Omega$ for which $\mathcal{P}$ is an $X$-$Y$-linkage are called the \emph{end segments of $\mathcal{P}$ in $(G,\Omega)$}.
\medskip

\begin{figure}[ht]
    \centering
    \begin{tikzpicture}

        \pgfdeclarelayer{background}
		\pgfdeclarelayer{foreground}
			
		\pgfsetlayers{background,main,foreground}

        \begin{pgfonlayer}{background}
            \pgftext{\includegraphics[width=6cm]{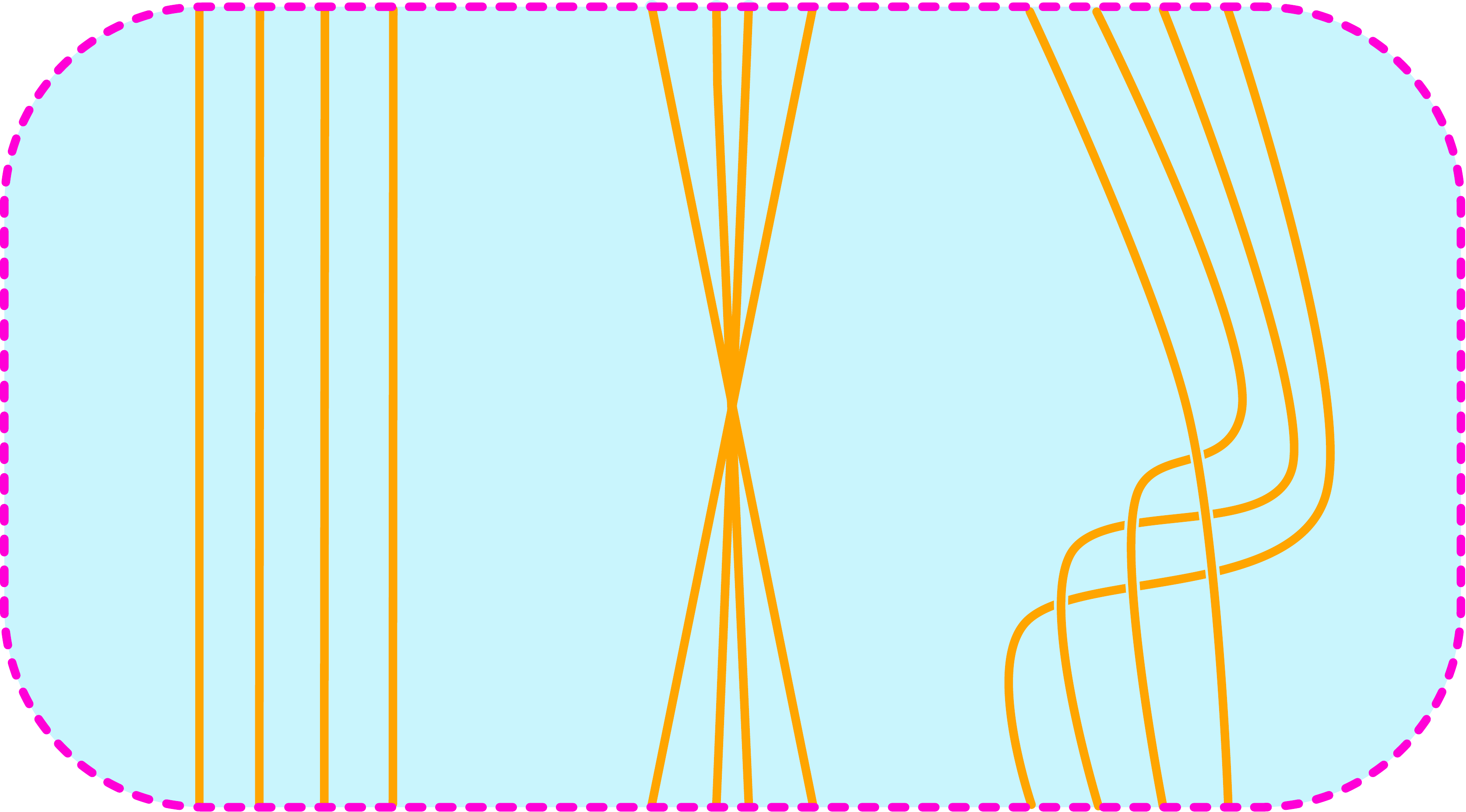}} at (C.center);
        \end{pgfonlayer}{background}
			
        \begin{pgfonlayer}{main}
        \node (C) [v:ghost] {};
            
        \end{pgfonlayer}{main}
        
        \begin{pgfonlayer}{foreground}
        \end{pgfonlayer}{foreground}

    \end{tikzpicture}
    \caption{Three monotone transactions each of order 4. To the left a planar transaction is represented, the transaction in the middle is a crosscap transaction in a common representation, and the right transaction is an example of a crosscap transaction drawn to help distinguish the paths.}
    \label{fig:monotone_transactions}
\end{figure}

Let \(\mathcal{P}\) be a transaction in a society \((G, \Omega)\) with the end segments \(X\) and \(Y\).
Suppose that the members of \(\mathcal{P}\) can be enumerated as \(P_1, \ldots, P_n\) so that if \(x_i \in X\) and \(y_i \in Y\) denote the endpoints of \(P_i\), then the vertices \(x_1, \ldots, x_n\) appear in the segment \(X\) in the listed order or the reverse one, and the vertices \(y_1, \ldots, y_n\) appear in \(Y\) in the listed order or the reverse one. Then we say that \(\mathcal{P}\) is \emph{monotone} (see \Cref{fig:monotone_transactions}), and if \(P_1, \ldots, P_n\) are ordered as above, they are \emph{indexed naturally}.

Should the vertices \(x_1, \ldots, x_n, y_n, \ldots, y_1\) appear in \(\Omega\) in the listed cyclic ordering or its reverse, we call \(\mathcal{P}\) a \emph{planar transaction}, and if the vertices \(x_1, \ldots, x_n, y_1, \ldots, y_n\) appear in \(\Omega\) in the listed cyclic ordering or its reverse, we call \(\mathcal{P}\) a \emph{crosscap transaction}.
The paths \(P_1\) and \(P_n\) are called the \emph{boundary paths} of \(\mathcal{P}\).

Monotone transactions are easy to find using a classic theorem by Erd\H{o}s and Szekeres.

\begin{proposition}[Erd\H{o}s and Szekeres \cite{ErdosS1935Combinatorial}]\label{prop:ErdosSzekeres}
Let $r,s$ be positive integers.
Every sequence of length at least $(r-1)(t-1)+1$ contains a monotonic decreasing subsequence of length $r$ or a monotonic increasing subsequence of length $s$.
\end{proposition}

\begin{lemma}\label{lem:monotonetransaction}
    Let $p,q$ be a positive integers, let $(G,\Omega)$ be a society, and let $\mathcal{P}$ be a transaction of order $(p-1)(q-1) + 1$.
    Then there exists a planar transaction $\mathcal{P}' \subseteq \mathcal{P}$ of order $p$ or a crosscap transaction $\mathcal{Q} \subseteq \mathcal{P}$ of order $q$, which can be found in $\mathbf{O}(pq|V(G)|)$-time.
\end{lemma}
\begin{proof}
    Let $r = (p-1)(q-1)^2 + 1$ and let $X_1$ and $X_2$ be the two end segments of $\mathcal{P}$.
    Furthermore, let $\mathcal{P} = \{ P_1, \ldots , P_r \}$, we label the endpoints of $P_i$ as $a_i \in X_1$ and $b_i \in X_2$, and choose the indices for $\mathcal{P}$ such that $a_1, a_2 \ldots , a_r$ appear in $X_1$ in the given order according to $\Omega$.
    Thus $a_1, \ldots , a_r$ induce an order on the paths in $\mathcal{P}$ and we can use \Cref{prop:ErdosSzekeres} on $b_1, \ldots , b_r$ to guarantee the existence of the desired monotone transaction, which we can find easily in linear time.
\end{proof}

\paragraph{Strips and their societies.}
Let $H$ be a subgraph of a graph $G$.
An \emph{$H$-bridge} in $G$ is a connected subgraph $B$ of $G$ such that $E(B) \cap E(H) = \emptyset$ and either $E(B)$ consists of a unique edge with both ends in $H$, or 
$B$ is constructed from a component $C$ of $G - V(H)$ and the non-empty set of edges $F \subseteq E(G)$ with one end in $V(C)$ and the other in $V(H)$, by taking the union of $C$, the endpoints of the edges in $F$, and $F$ itself.
The vertices in $V(B) \cap V(H)$ are called the \emph{attachments} of $B$.

\begin{figure}[ht]
    \centering
    \begin{tikzpicture}[scale=1.25]

        \pgfdeclarelayer{background}
		\pgfdeclarelayer{foreground}
			
		\pgfsetlayers{background,main,foreground}

        \begin{pgfonlayer}{background}
            \pgftext{\includegraphics[width=6cm]{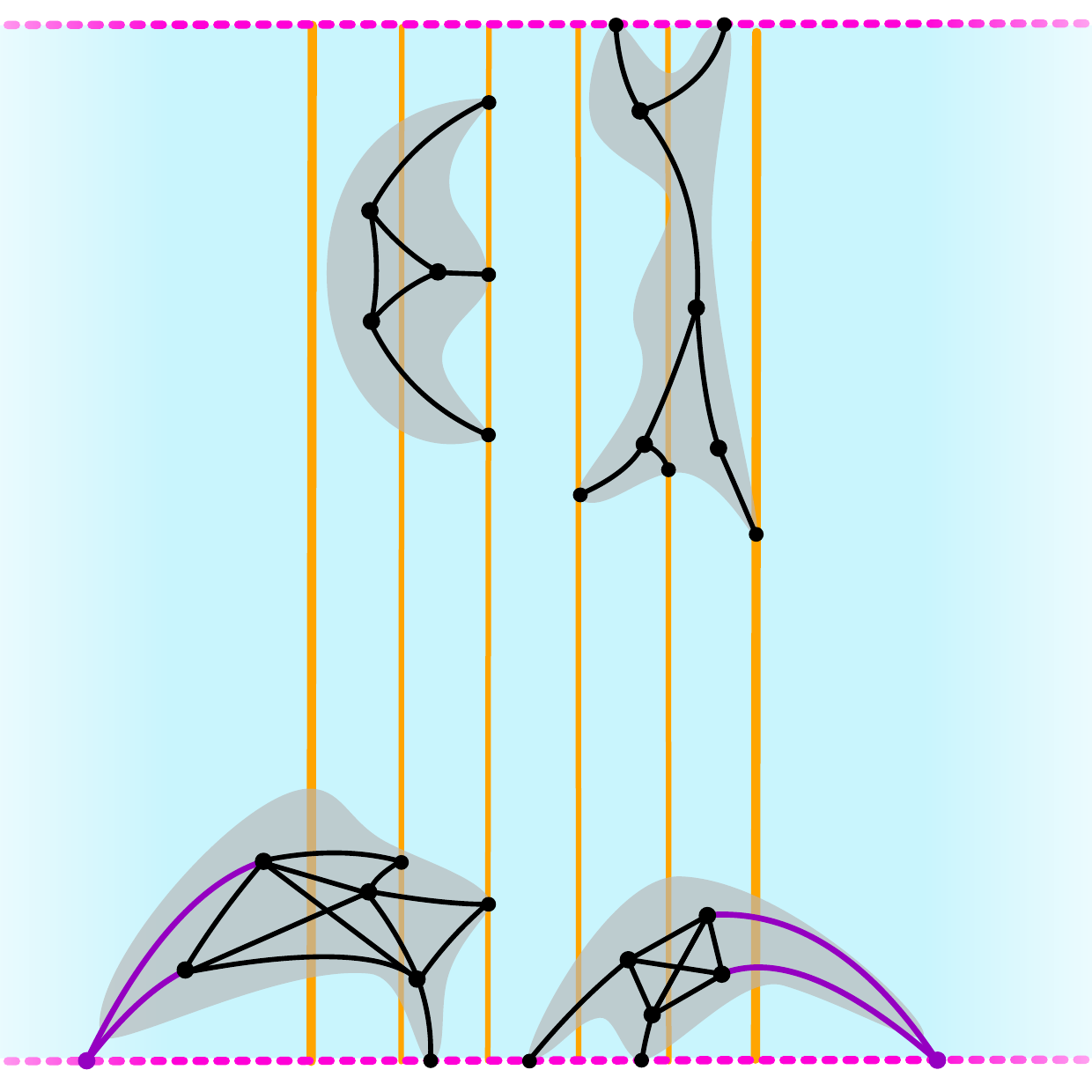}} at (C.center);
        \end{pgfonlayer}{background}
			
        \begin{pgfonlayer}{main}
        \node (C) [v:ghost] {};
            
        \end{pgfonlayer}{main}
        
        \begin{pgfonlayer}{foreground}
        \end{pgfonlayer}{foreground}

    \end{tikzpicture}
    \caption{A representation of which kind of bridges (indicated in black) contribute to the definition of a strip based on a transaction (marked in orange). For the two bridges with edges to attachments that need to be removed, the respective attachments and edges are coloured in purple.}
    \label{fig:strip}
\end{figure}

We let $H$ denote the subgraph of $G$ obtained from the union of elements of $\mathcal{P}$ by adding the elements of $V(\Omega)$ as isolated vertices.
Further, we define $H'$ as the subgraph of $H$ consisting of $\mathcal{P}$ and all vertices of $X \cup Y$.
Consider all $H$-bridges of $G$ with at least one attachment in $V(H') \setminus V(P_1 \cup P_n)$, and for each such $H$-bridge $B$ let $B'$ denote the graph obtained from $B$ by deleting all attachments that do not belong to $V(H')$.
We let $G_1$ denote the union of $H'$ and all graphs $B'$ as above (see \Cref{fig:strip} for an example of which bridges are included).

The \emph{\(\mathcal{P}\)-strip society in \((G, \Omega)\)} is defined as the society \((G_1, \Omega_1)\), where \(\Omega_1\) is the concatenation of the segment \(X\) ordered from \(x_1\) to \(x_n\), and the segment \(Y\) ordered from \(y_n\) to \(y_1\).
If the \(\mathcal{P}\)-strip society admits a vortex-free rendition in a disk, we call \(\mathcal{P}\) a \emph{flat transaction}.
Further, if no edge of $G$ has an endpoint in $V(G_1) \setminus V(P_1 \cup P_n)$ and the other point $V(G) \setminus V(G_1)$, then we call $\mathcal{P}$ \emph{isolated}.
Let $X',Y'$ be the two distinct segments of $\Omega$ that have one endpoint in $X$ and the other in $Y$.
Note that $V(\Omega) = X \cup Y \cup X' \cup Y'$.
We say that \(\mathcal{P}\) is \emph{separating} if it is isolated and there exists no $X'$-$Y'$-path in $G - V(G_1)$.

\begin{figure}[ht]
    \centering
    \begin{tikzpicture}

        \pgfdeclarelayer{background}
		\pgfdeclarelayer{foreground}
			
		\pgfsetlayers{background,main,foreground}

        \begin{pgfonlayer}{background}
            \pgftext{\includegraphics[width=6cm]{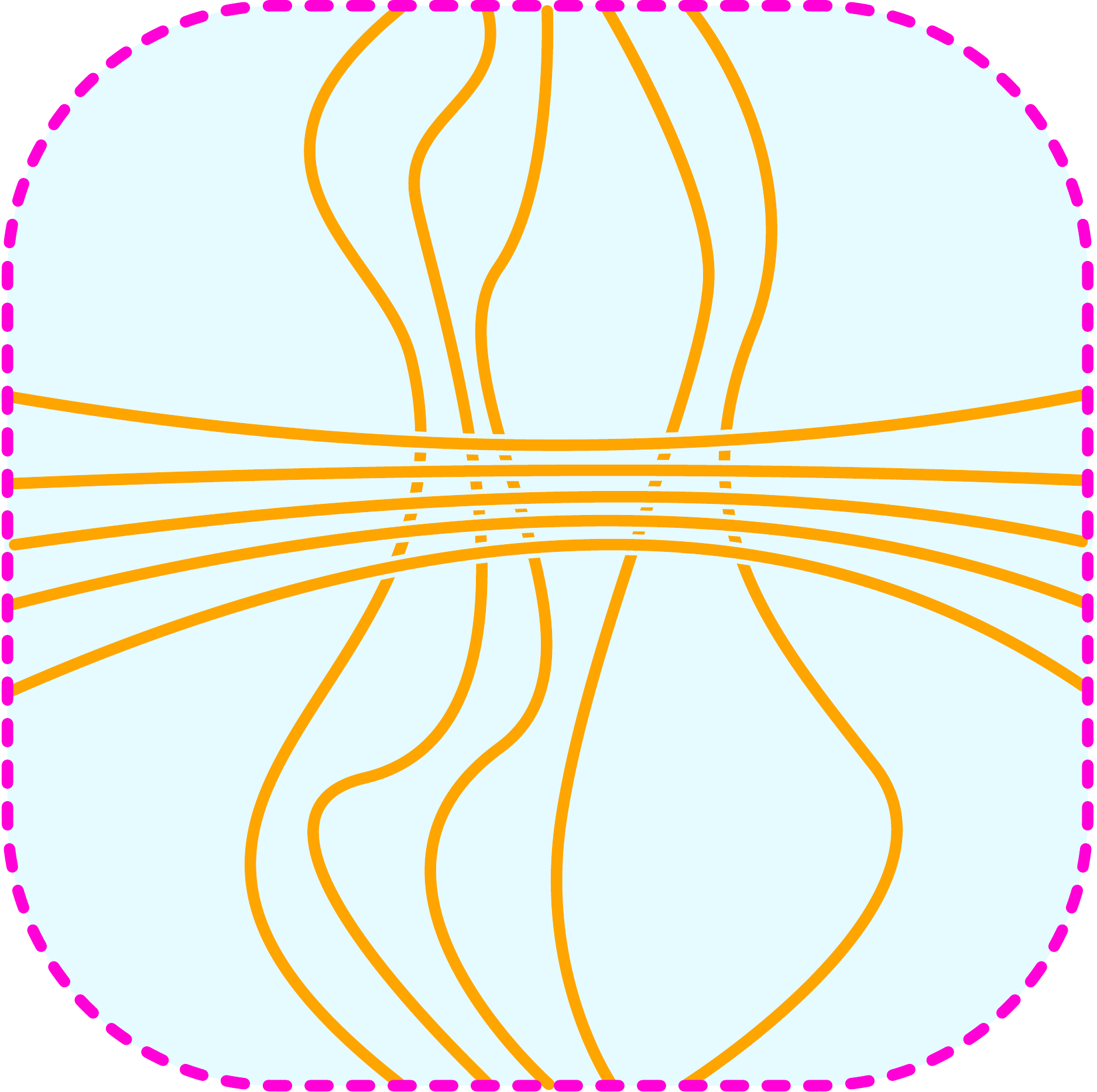}} at (C.center);
        \end{pgfonlayer}{background}
			
        \begin{pgfonlayer}{main}
        \node (C) [v:ghost] {};
            
        \end{pgfonlayer}{main}
        
        \begin{pgfonlayer}{foreground}
        \end{pgfonlayer}{foreground}

    \end{tikzpicture}
    \caption{A handle transaction of order 10.}
    \label{fig:handle_transaction}
\end{figure}

\paragraph{Handle transactions.}
A transaction $\mathcal{P}$ of order $2n$ in a society $(G, \Omega)$, for a positive integer $n \in \mathbb{N}$, is called a \emph{handle transaction} if $\mathcal{P}$ can be partitioned into two transactions $\mathcal{R}, \mathcal{Q}$ each of order $n$, such that both are planar, and $S^\mathcal{R}_1, S^\mathcal{Q}_1, S^\mathcal{R}_2, S^\mathcal{Q}_2$ are segments partitioning $\Omega$, with $\mathcal{X}$ being a $S^\mathcal{X}_1$-$S^\mathcal{X}_2$-linkage for both $\mathcal{X} \in \{ \mathcal{R} , \mathcal{Q} \}$, and the segments are found on $\Omega$ in the order they were listed above (see \Cref{fig:handle_transaction} for an example).
We call a handle transaction $\mathcal{P} = \mathcal{R} \cup \mathcal{Q}$ \emph{isolated in $G$}, respectively \emph{flat in $G$}, if both the $\mathcal{R}$-strip society and the $\mathcal{Q}$-strip society of $(G, \Omega)$ are isolated in $G$, respectively flat in $G$.
The transactions $\mathcal{R}$ and $\mathcal{Q}$ are called the two \emph{constituent} transactions of $\mathcal{P}$.

\subsection{Cylindrical renditions, nests, and more on societies}\label{sec:nest}

The proof of the GMST as imagined by Kawarabayashi et al.\ \cite{KawarabayashiTW2021Quickly} is based on the systematic study of societies.
Our new proof does not change this perspective but rather modifies the method of how societies are studied.
This allows us to further adopt more involved definitions from \cite{KawarabayashiTW2021Quickly} which facilitate our own approach.

Let $\rho=(\Gamma,\mathcal{D})$ be a rendition of a society $(G,\Omega)$ in a surface $\Sigma$.
For every cell $c\in C(\rho)$ with $|\widetilde{c}|=2$, we select one of the components of $\mathsf{boundary}(c)-\widetilde{c}$.
This selection will be called a \emph{tie-breaker in $\rho$}, and we assume that every rendition comes equipped with a tie-breaker.

This assumption is mostly for convenience to have the following notion of ``traces'' be well-defined.
While $\Sigma$-decompositions bear close resemblance to actual embeddings of graphs in surfaces, they are better understood as an imperfect variant of such embeddings.
Still, some of the topological qualities of embeddings are also witnessed by $\Sigma$-decompositions.
One particularly important feature is, that paths and cycles encode (closed) curves in the surface which can, in some cases, be used to define territories or even vertex separators.
While the graph $G$ is not necessarily drawn by $\Gamma$ in a way that reveals these curves immediately, the structure of the cells of of the $\Sigma$-decomposition $\delta$ allows to ``trace'' along any path or cycle that contains nodes of $\delta$ and thereby define such a curve.

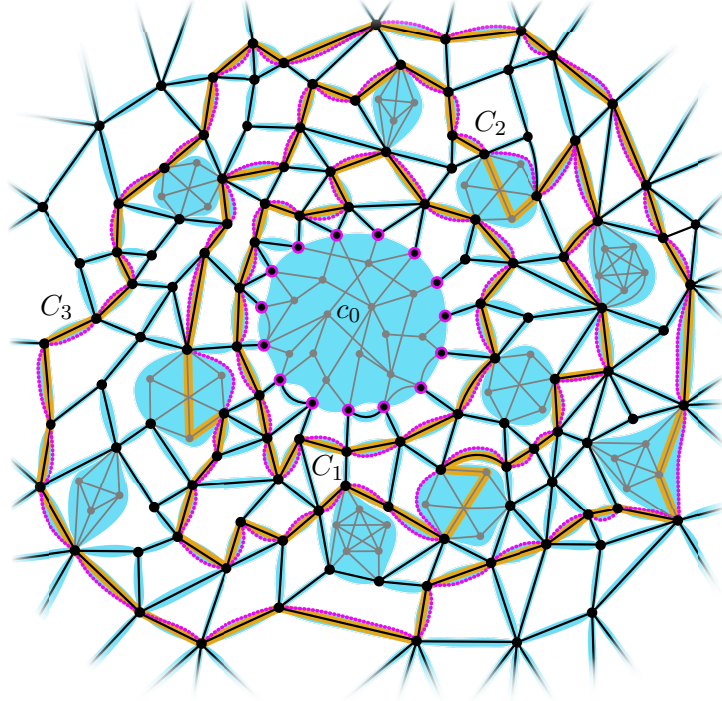
\begin{figure}[ht]
    \centering
    \begin{tikzpicture}

        \pgfdeclarelayer{background}
		\pgfdeclarelayer{foreground}
			
		\pgfsetlayers{background,main,foreground}

        \begin{pgfonlayer}{background}
            \pgftext{\includegraphics[width=10cm]{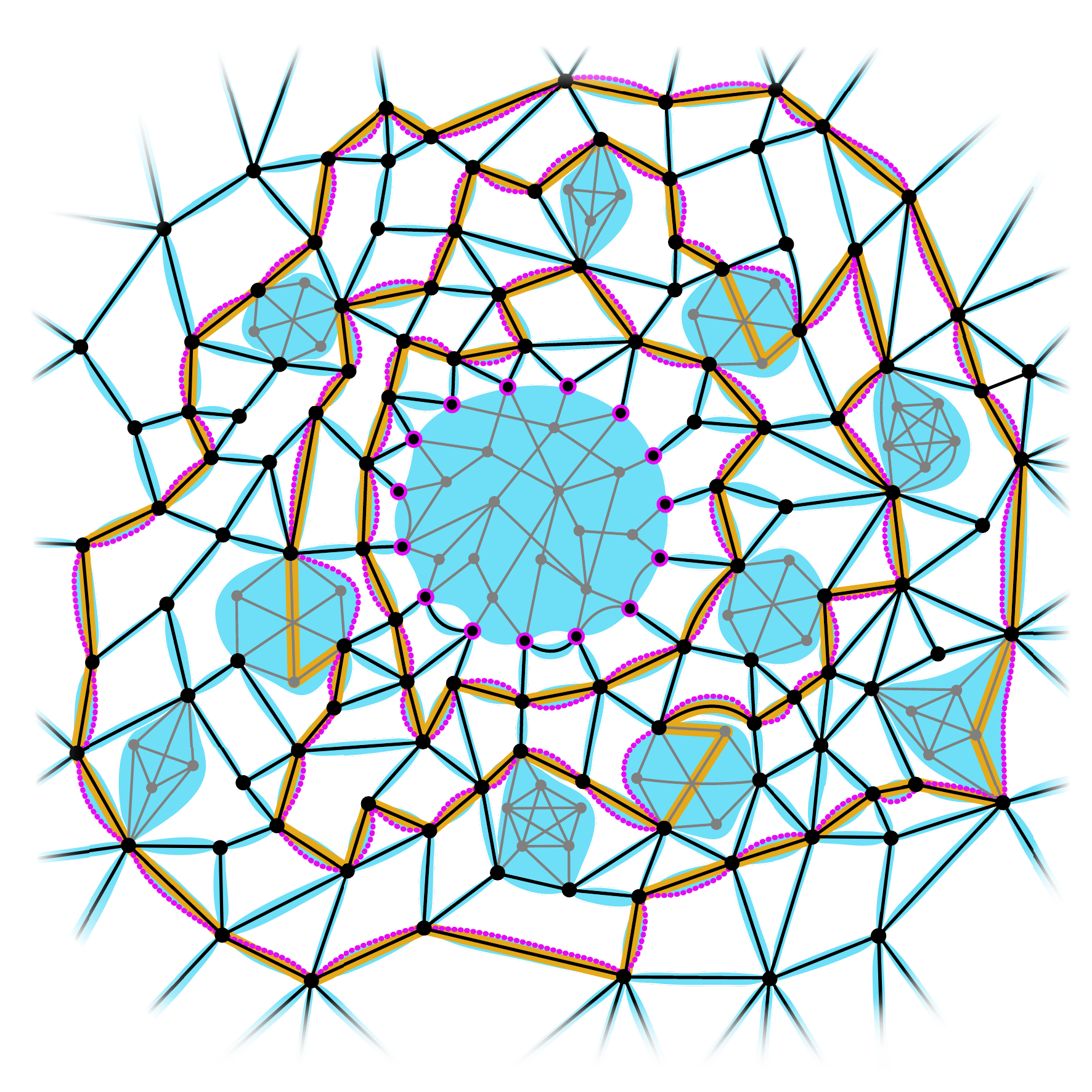}} at (C.center);
        \end{pgfonlayer}{background}
			
        \begin{pgfonlayer}{main}
        \node (C) [v:ghost] {};

        \node (vortex) [v:ghost,position=115:4.5mm from C] {$c_0$};

        \node (C1) [v:ghost,position=254:17mm from C] {$C_1$};
        \node (C2) [v:ghost,position=60:34mm from C] {$C_2$};
        \node (C3) [v:ghost,position=173:41mm from C] {$C_3$};
            
        \end{pgfonlayer}{main}
        
        \begin{pgfonlayer}{foreground}
        \end{pgfonlayer}{foreground}

    \end{tikzpicture}
    \caption{A vortex cell $c_0$ in a $\Sigma$-decomposition surrounded by three cycles of the nest $\{ C_1, C_2, C_3 \}$. The dotted line indicates the traces of the cycles $C_1$, $C_2$, and $C_3$.}
    \label{fig:vortex_with_nest}
\end{figure}

\paragraph{Traces of paths and cycles.}
Let $G$ be a graph and $\delta$ be a $\Sigma$-decomposition for $G$.
Let $Q$ be either a cycle or a path in $G$ that uses no edge of $\sigma(c)$ for every vortex $c\in C(\rho)$.
We say that $Q$ is \emph{grounded} in $\delta$ if either $Q$ is a non-zero length path with both endpoints in $\pi_{\delta}(N(\delta))$, or $Q$ is a cycle and it uses edges of $\sigma(c_1)$ and $\sigma(c_2)$ for two distinct cells $c_1, c_2 \in C(\delta)$.
If $Q$ is grounded we define the \emph{trace}\footnote{This is called a \emph{track} in \cite{KawarabayashiTW2021Quickly}. We rename this concept as we believe that the word trace more accurately captures the idea that we are following a line in a drawing to define this object.} of $Q$ as follows.
Let $P_1,\dots,P_k$ be distinct maximal subpaths of $Q$ such that $P_i$ is a subgraph of $ \sigma(c)$ for some cell $c$.
Fix an index $i \in [k]$.
The maximality of $P_i$ implies that its endpoints are $\pi(n_1)$ and $\pi(n_2)$ for distinct nodes $n_1,n_2\in N(\rho)$.
If $|\widetilde{c}|=2$, let $L_i$ be the component of $\mathsf{boundary}(c)-\{n_1,n_2\}$ selected by the tie-breaker, and if $|\widetilde{c}|=3$, let $L_i$ be the component of $\mathsf{boundary}(c)-\{n_1,n_2\}$ that is disjoint from $\widetilde{c}$.
Finally, we define $L_i'$ by pushing $L_i$ slightly so that it is disjoint from all cells in $C(\delta)$.
We define such a curve $L_i'$ for all $i \in [k]$, while maintaining that these curves intersect only at a common endpoint.
The \emph{trace} of $Q$ is defined to be $\bigcup_{i\in[k]}L_i'$.
So the trace of a cycle is the homeomorphic image of the unit circle, and the trace of a path is an arc in $\Sigma$ with both endpoints in $N(\delta)$.
In \cref{fig:vortex_with_nest} the closed dotted curves illustrate the traces of three highlighted cycles.

\paragraph{Cylindrical renditions.} 
Let $(G,\Omega)$ be a society, let $\rho=(\Gamma,\mathcal{D})$ be a rendition of $(G,\Omega)$ in a disk, and let $c_0\in C(\rho)$ be the unique vortex in $\rho$.
In those circumstances we say that the triple $(\Gamma,\mathcal{D},c_0)$ is a \emph{cylindrical rendition} of $(G,\Omega)$ around $c_0$.
\medskip

If $(G, \Omega)$ is a society with a cylindrical rendition $(\Gamma,\mathcal{D},c_0)$ and $\mathcal{P}$ is a transaction in $(G,\Omega)$, we call $\mathcal{P}$ \emph{exposed} if for every path $P \in \mathcal{P}$ there exists an edge $e \in E(P) \cap \sigma(c_0)$ and if for every path $P \in \mathcal{P}$ we have $E(P) \cap E(\sigma(c_0)) = \emptyset$, we call $\mathcal{P}$ \emph{unexposed}.\footnote{This flips the meaning of these two terms in comparison to the definition Kawarabayashi et al.\ chose in \cite{KawarabayashiTW2021Quickly}. We make this change since we believe it to be more intuitive to call a transaction that touches the undecomposed, possibly very chaotic part of our rendition exposed and consider a transaction that is already embedded in an almost planar fashion unexposed.}
A transaction may be neither exposed nor unexposed.
But at least half of any transaction will fall into one of these two categories.

\paragraph{Nests.}
Let $\rho = (\Gamma, \mathcal{D}, c_0)$ be a cylindrical rendition of a society $(G, \Omega)$ in a disk $\Delta$, then a \emph{nest (in $\rho$) around $c_0$} is a set of disjoint cycles $\mathcal{C} = \{ C_1, \ldots , C_s \}$ in $G$ such that each of them is grounded in $\rho$ and the trace of $C_i$ bounds a closed disk $\Delta_i$ in such a way that $c_0 \subseteq \Delta_1 \subseteq \cdots \subseteq \Delta_s \subseteq \Delta$.
We will sometimes use the letter $\mathcal{C}$ to denote the graph $\bigcup_{i=1}^s C_i$.
If $c_0$ is clear from the context, we call it a nest in $\rho$ and drop the reference to $c_0$.
See \cref{fig:vortex_with_nest} for an illustration of a nest around a vortex.
\medskip

\paragraph{Grounded subgraphs and aligned disks.}
Let $G$ be a graph and let $\delta=(\Gamma,\mathcal{D})$ be a $\Sigma$-decomposition for $G$, with $H$ being a 2-connected subgraph of $G$.

We say that $H$ is \emph{grounded} in $\delta$ if every cycle in $H$ is grounded and no vertex of $H$ belongs to a vortex of $\delta$.
A planar and $2$-connected graph $H\subseteq G$ is said to be \emph{flat} in $\delta$ if it is grounded and there exists a disk $\Delta\subseteq\Sigma$ whose boundary intersects $\Gamma$ only in nodes of $\delta$, such that $H$ is drawn in $\Delta$ and no vortex of $\delta$ is contained in $\Delta$.

Given a closed curve $\gamma$ in $\Sigma$, we say that $\gamma$ is \emph{$\delta$-aligned} if it only intersects $\delta$ in nodes.
Furthermore, a disk in $\Sigma$ is called \emph{$\delta$-aligned} if its boundary is $\delta$-aligned.
For any $\delta$-aligned disk $\Delta$, we call the subgraph of $G$ that is drawn by $\Gamma$ onto $\Delta$ the \emph{crop of $H$ by $\Delta$ (in $\delta$)}.
Furthermore, the \emph{restriction $\delta'$ of $\delta$ by $\Delta$} is defined as the $\Delta$-decomposition that consists of the restriction of both $\Gamma$ and $\mathcal{D}$ to $\Delta$.

This then further defines a society as follows.
Let $V(\Omega_{\Delta})$ be the set of all vertices whose corresponding nodes are drawn in the boundary of $\Delta$ and let $\Omega_{\Delta}$ be the cyclic ordering of $V(\Omega_{\Delta})$ obtained by traversing along $\mathsf{boundary}(\Delta)$ in the anticlockwise direction.
Now, let $G_{\Delta}$ be the crop of $G$ by $\Delta$.
We now call the society $(G_{\Delta}, \Omega_{\Delta})$ the \emph{$\Delta$-society (in $\rho$)}.
If $\rho$ is clear from the context, we do not mention it.

Let $C$ be a grounded cycle in $G$.
We say that $C$ is \emph{contractible} if its trace defines a contractible curve in $\Sigma$ and note that the trace of $C$ bounds at least one disk $\Delta_C$ that allows us to define a society $(G_{\Delta_C}, \Omega_{\Delta_C})$.

\medskip
Importantly, not all parts of $C$ have to be contained in $G_{\Delta_C}$, as the subpaths of $C$ that are drawn within cells of $\delta$ could be ``pushed out'' of $\Delta_C$ by the tie-breakers associated with the trace of $C$.

\paragraph{Inner and outer graphs of a cycle.}
Let $(G, \Omega)$ be a society and $\rho$ a rendition of $(G, \Omega)$ in a disk.
Let $C$ be a cycle in $G$ that is grounded in $\rho$ such that the trace $T$ of $C$ bounds a closed disk $\Delta$.
Let $L$ be the subgraph of $G$ with vertex set $\pi(N(\rho) \cap T)$ and no edges.
We define the \emph{outer graph of $C$ in $\rho$} as
\[ L \cup \bigcup \{ \sigma (c) ~\!\colon\!~ c \in C(\rho) \text{ and } c \not\subseteq \Delta \} ,\]
and we also call the society $(G', \Omega')$ associated with $\Delta$ the \emph{$C$-society in $\rho$}, with $G'$ being the \emph{inner graph of $C$ in $\rho$}.
For the outer graph $G''$ of $C$ in $\rho$ we note that $G'$ and $G''$ are edge-disjoint, $G' \cup G'' = G$, and $G' \cap G'' = L$.

\paragraph{Linkages in relation to a nest.}
Let $(G, \Omega)$ be a society, let $\mathcal{L}$ and $\mathcal{P}$ be two linkages in $G$, and let $C$ be a grounded cycle in a rendition $\rho$ of $(G,\Omega)$ in a disk.
We say that $\mathcal{L}$ is \emph{coterminal with $\mathcal{P}$ up to level $C$} if there exists a subset $\mathcal{P}'$ of $\mathcal{P}$ such that $H \cap \mathcal{L} = H \cap \mathcal{P}'$, for the outer graph $H$ of $C$.
Whenever $\mathcal{C} = \{ C_1, \ldots , C_s \}$ is a nest in $\rho$ and it is clear from the context which nest we refer to, we will abbreviate ``coterminal with $\mathcal{P}$ up to level $C_i$'' with ``coterminal with $\mathcal{P}$ up to level $i$''.
\medskip

Of particular importance to our approach are linkages that connect the innermost cycle $C_1$ of a nest $\mathcal{C} = \{ C_1, \ldots , C_s \}$ with the society associated with this nest.
These linkages appear in \cite{KawarabayashiTW2021Quickly} as a tool that allows them to reroute the crosses, crosscap transactions, or handle transactions that they find throughout their proof into parts of the graph that permit them to build a $K_t$-minor, should too many of these objects accumulate.
This is an essential part of their proof strategy that we do not dispense with.
However, these linkages are quite problematic for the bounds that their proof yields, since at several points they must ensure the existence of sizeable instances of these linkages without knowing whether they need all of the paths of the linkage that they find.
An important contribution of ours is that we defer the determination of the needed parts of these linkages until the end of our proof, at the expense of having to introduce an additional data structure that we must maintain.
We will introduce the required tools for this later, but first, we define the linkages we just discussed.

If $(G, \Omega)$ is a society with a nest $\mathcal{C} = \{ C_1, \ldots , C_s \}$ in a rendition $\rho$ of $(G, \Omega)$ in a disk, we call a $V(\Omega)$-$V(C_1)$-linkage $\mathcal{R}$ a \emph{radial linkage (in $\rho$) for $\mathcal{C}$}.
We say that $\mathcal{R}$ is \emph{orthogonal to $\mathcal{C}$} if for all $C \in \mathcal{C}$ and all $R \in \mathcal{R}$ the graph $C \cap R$ is a path.
Similarly, we say that a transaction $\mathcal{P}$ in $(G,\Omega)$ is \emph{orthogonal to $\mathcal{C}$} if for all $C\in\mathcal{C}$ and all $P\in\mathcal{P}$ the graph $C\cap P$ consists of exactly two paths.

Please note that the definitions above must work, in particular, in subcubic graphs.
This means that it is not possible to replace the requirement of the intersection of a path with a cycle to be a path by ``any path should intersect any cycle in exactly one vertex''.

\subsection{Vortices and their depth}\label{subsec:vortex}
Finally, we need a notion of ``depth'' for our vortices.
So far, the definitions allow us to hide almost anything within a vortex, as long as its interface with the rest of the graphs agrees on the cylindrical ordering of some set of vertices.
In the process of the proof, such vortices will be refined further until the graphs drawn in their interiors do not provide enough infrastructure to continue this refinement process.

\paragraph{Depth of vortices.}
Let $G$ be a graph and $\delta$ be a $\Sigma$-decomposition of $G$ with a vortex cell $c_0$.
Notice that $c_0$ defines a society $(\sigma(c_0),\Omega_{c_0})$ where $V(\Omega_{c_0})$ is the set of vertices of $G$ corresponding to the nodes of $\delta$ drawn in the boundary of the closure of $c_0$.
The ordering $\Omega$ is obtained by traversing along the boundary of the closure of $c_0$ in anti-clockwise direction.
We call $(\sigma(c_0),\Omega_{c_0})$ as obtained above the \emph{vortex society} of $c_0$.

We define the \emph{depth} of a society $(G,\Omega)$ as the maximum cardinality of a transaction in $(G,\Omega)$.
The \emph{depth} of the vortex $c_0$ is thereby defined as the depth of its vortex society.
\medskip

All vortices that remain in the structural decomposition that the GMST offers are of bounded depth, where the bound is some function in the size of the $K_t$-minor we exclude.
While the definition above is very much aligned with the way we further refine a vortex in our proof, that is by iteratively finding large transactions passing through it, it is not yet good enough to provide the information necessary to extend a vortex of bounded depth further into a tree-decomposition with small adhesion.
For this we need one additional definition from the original statement of the GMST.

\paragraph{Linear decompositions of vortices.}
Let $(G,\Omega)$ be a society.
A \emph{linear decomposition} of $(G,\Omega)$ is a labeling $v_1,v_2,\dots,v_n$ of $V(\Omega)$ such that $v_1,v_2,\dots,v_n$ appear in $\Omega$ in the order listed, together with sets $(X_1,X_2,\dots,X_n)$ such that
\begin{enumerate}
    \item $X_i\subseteq V(G)$ and $v_i\in X_i$ for all $i\in[n]$,
    \item $\bigcup_{i\in[n]}X_i=V(G)$ and for every $uv\in E(G)$ there exists $i\in[n]$ such that $u,v\in X_i$, and
    \item for every $x\in V(\Omega)$ the set $\{ i\in[n] ~\!\colon\!~ x\in X_i \}$ forms an interval in $[n]$.
\end{enumerate}
The \emph{adhesion} of a linear decomposition is $\max \{ |X_i\cap X_{i+1}| ~\!\colon\!~ i\in[n-1] \}$.
The \emph{width} of a linear decomposition is $\max \{ |X_i| ~\!\colon\!~ i\in[n] \}$.
\medskip

What is called ``adhesion'' above was originally called the ``depth'' in \cite{RobertsonS1990Graph}.
Since we follow the path of Kawarabayashi et al., we adopt their terminology and stick to the term ``adhesion'' \cite{KawarabayashiTW2021Quickly}, which also aligns with the terminology we use for tree-decompositions.

It is easy to see that every society with a linear decomposition of adhesion at most $k$ has depth at most $2k$, where depth here is used in our sense of the definition.
The (partial) reverse of this observation was shown by Robertson and Seymour in \cite{RobertsonS1990Graph}.

\begin{proposition}[Robertson and Seymour \cite{RobertsonS1990Graph}]\label{prop:depth_to_lin_decomp}
Let $k$ be a non-negative integer and $(G,\Omega)$ be a society of depth at most $k$.
Then $(G,\Omega)$ has a linear decomposition of adhesion at most $k$.
\end{proposition}

We give a new algorithmic proof of this fact in \Cref{lem:findtransactionorlineardecomp}.

\subsection{Extracting a cross from a vortex}\label{subsec:crossExtraction}

In the following, we extend the \hyperref[prop:TwoPaths]{Two Paths Theorem} by adding a small amount of infrastructure to the assumptions on our society $(G,\Omega)$.
In return, this strengthened version allows us to force the endpoints of the cross that must occur if $(G,\Omega)$ does not have a vortex-free rendition in the disk onto four specified vertices of $\Omega$.

\begin{lemma}\label{lemma:reroute_crosses}
Let $(G,\Omega)$ be a society with a cylindrical rendition $\rho$, a nest $\mathcal{C}$ of order $16$ around the vortex $c_0$, a radial linkage $\mathcal{R}$ of order $10$ orthogonal to $\mathcal{C}$, another\footnote{We allow $\mathcal{Q}$ to share paths with $\mathcal{R}$ or even be contained in $\mathcal{R}$.} radial linkage $\mathcal{Q}$ of order $4$, also orthogonal to $\mathcal{C}$, and let $x_1,x_2,x_3,x_4$ be the endpoints of the paths $\mathcal{Q}$ in $V(\Omega)$ numbered in order of their appearance in $\Omega$.
Then exactly one of the following is true.
\begin{enumerate}
    \item $(G,\Omega)$ has a vortex-free rendition in the disk.
    \item There exists two vertex-disjoint paths $P_1,P_2$ in $G$ such that $P_i$ has endpoints $x_i$ and $x_{i+3}$ for both $i\in[2]$, and $P_1$ and $P_2$ are otherwise disjoint from $V(\Omega)$.
\end{enumerate}
Moreover, there exists an algorithm that finds one of these outcomes in time $\mathcal{O}(|V(G)|)$.
\end{lemma}

Please note that the constants $16$ and $10$ here are suboptimal.
We chose to take this route purely for convenience as \cref{lemma:reroute_crosses} follows almost immediately from the following special case of Lemma 4.5 from \cite{KawarabayashiTW2021Quickly}.
The original statement is much more powerful than needed for our purpose, the version we state here is the one where we set $r=2$.

\begin{proposition}[Kawarabayashi, Thomas, and Wollan \cite{KawarabayashiTW2021Quickly}]\label{prop:oldOrthogonality}
Let $s\geq 11$ be a positive integer.
Let $(G,\Omega)$ be a society with a cylindrical rendition and a nest $\mathcal{C}=\{ C_1,\dots,C_s\}$ of order $s$ around the vortex $c_0$.
Further, let $\mathcal{R}$ be a radial linkage of order $10$.
If $(G,\Omega)$ has a cross, then there exists a transaction $\mathcal{P}=\{ P_1,P_2\}$ in $(G,\Omega)$ that is coterminal with $\mathcal{R}$ up to level $11$ such that $P_1$ and $P_2$ form a cross in $(G,\Omega)$.
\end{proposition}

Kawarabayashi et al.\ explain in \cite{KawarabayashiTW2021Quickly} that \cref{prop:oldOrthogonality} yields an algorithm with linear running time since their proof is only a combination of certain exchange properties of paths and applications of Menger's Theorem for a fixed number of paths.

With this we are ready for the proof of \cref{lemma:reroute_crosses}.

\begin{proof}[Proof of \cref{lemma:reroute_crosses}]
Let $\mathcal{C}=\{ C_1,\dots,C_s\}$.
Let $(G',\Omega')$ be the $C_{12}$-society in $\rho$.
We begin by applying \cref{prop:TwoPaths} to $(G',\Omega')$.
If the outcome is that $(G',\Omega')$ has a vortex-free rendition in the disk, we can combine this rendition with $\rho$ to obtain a vortex-free rendition of $(G,\Omega)$ in the disk as desired.
Hence, we may assume that $(G',\Omega')$ has a cross.

We now apply \cref{prop:oldOrthogonality} to $(G',\Omega')$ with nest $\{C_1,\dots,C_{11} \}$ and obtain a transaction $\mathcal{P}''=\{ P_1'',P_2''\}$ in $(G',\Omega')$ which is coterminal up to level $11$ with the radial linkage $\mathcal{R}'$ in $(G',\Omega')$ obtained from $\mathcal{R}$.
Moreover, the paths in $\mathcal{P}''$ form a cross in $(G',\Omega')$.
We may extend $\mathcal{P}''$ along the paths in $\mathcal{R}$ to a transaction $\mathcal{P}$ in $(G,\Omega)$ which is still coterminal with $\mathcal{R}$ up to level $11$ and which still forms a cross in $(G,\Omega)$.
Notice that $\mathcal{P}'$ is orthogonal to $\{ C_{12},C_{13},C_{14},C_{15},C_{16}\}$.

Let $H'$ be the outer graph of $C_{11}$ and let $H''$ be the inner graph of $C_{16}$.
Moreover, let $X$ be the set of all vertices $x$ such that $x$ is the first vertex of some $P\in\mathcal{P}'$ on $C_{11}$ when traversing along $P$ starting from one of its endpoints.
Similarly, let $Y$ be the set of all vertices $y$ such that $y$ is the first vertex of some path $Q\in\mathcal{Q}$ on $C_{15}$ we encounter when traversing along $Q$ starting from a vertex on $C_1$.
In $H\coloneqq H'\cap H''$ there must exist an $X$-$Y$-linkage $\mathcal{L}$ of order $4$.
Combining $\mathcal{P}''$ with $\mathcal{L}$ and what remains of $\mathcal{Q}$ in the outer graph of $C_{16}$ now yields a transaction $\mathcal{P}=\{ P_1,P_2\}$ which is coterminal with $\mathcal{Q}$ up to level $16$ and whose paths form a cross in $(G,\Omega)$.

To see that our claim on the running time of the implied algorithm holds notice the following.
Let $x_1,x_2,x_3,x_4$ be the endpoints of the paths in $\mathcal{Q}$ in $V(\Omega)$ numbered in order of appearance in $\Omega$.
Let $\Psi\coloneqq (x_1,x_2,x_3,x_4)$ be the derived cyclic order of these four vertices.
What we have proven above is that the existence of a cross in the society $(G',\Omega')$ defined by the inner graph of $C_{12}$ implies the existence of a cross in the society $(G-(V(\Omega)\setminus V(\Psi)),\Psi)$.
This means that we can find the desired paths $P_1$ and $P_2$ or the vortex-free-rendition of $(G,\Omega)$ by two applications of \cref{prop:TwoPaths}.
\end{proof}

\subsection{Surface walls}\label{sec:surfacwalls}

As a final piece of infrastructure, we extend the notion of walls to include surfaces and additionally, vortices.
This part is inspired by the work of Thilikos and Wiederrecht on excluding graphs of bounded genus \cite{ThilikosW2024Excluding}.

\paragraph{Annulus walls.}
Let $m,n$ be positive integers.
The \emph{$(n\times m)$-annulus grid} is the graph obtained from the $(n\times m)$-grid by adding the edges $\{ \{(i,1),(i,n)\} ~\!\colon\!~\ i\in[n] \}$.
The \emph{elementary $(n\times m)$-annulus wall} is the graph obtained from the $(n\times 2m)$-annulus grid by deleting all edges in the following set
\begin{align*}
    \big\{  \{(i,j),(i+1,j) \}  ~\!\colon\!~ i\in[n-1],\text{ }j\in[2m]\text{, and }i\not\equiv j\mod 2 \big\}.
\end{align*}
An \emph{$n$-annulus wall} is a subdivision of the elementary $n$-annulus wall.
\medskip

One can also see an annulus $(n\times m)$-wall as the graph obtained by completing the horizontal paths of a wall to cycles instead of discarding the vertices of degree one.
This viewpoint will be very helpful in the following.
An $n$-annulus wall contains $n$ cycles, $C_1,\dots,C_n$, such that $C_i$ consists exactly of the vertices of the $i$th row of the original wall.
We refer to these cycles as the \emph{base cycles} of the $n$-annulus wall.

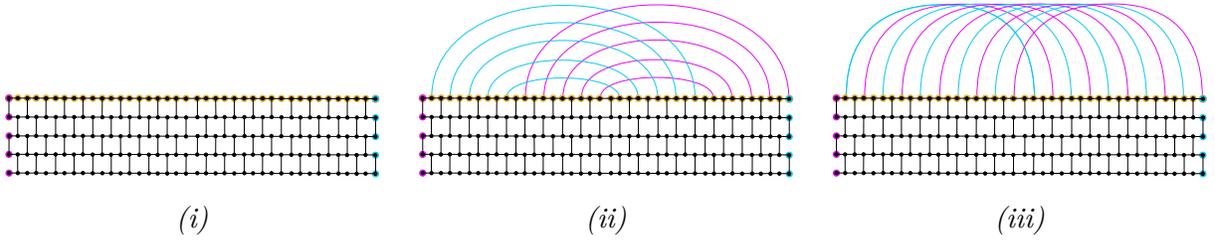
\begin{figure}[ht]
    \centering
    \begin{tikzpicture}

        \pgfdeclarelayer{background}
		\pgfdeclarelayer{foreground}
			
		\pgfsetlayers{background,main,foreground}

        \begin{pgfonlayer}{background}
        \node (C) [v:ghost] {{\includegraphics[width=5cm]{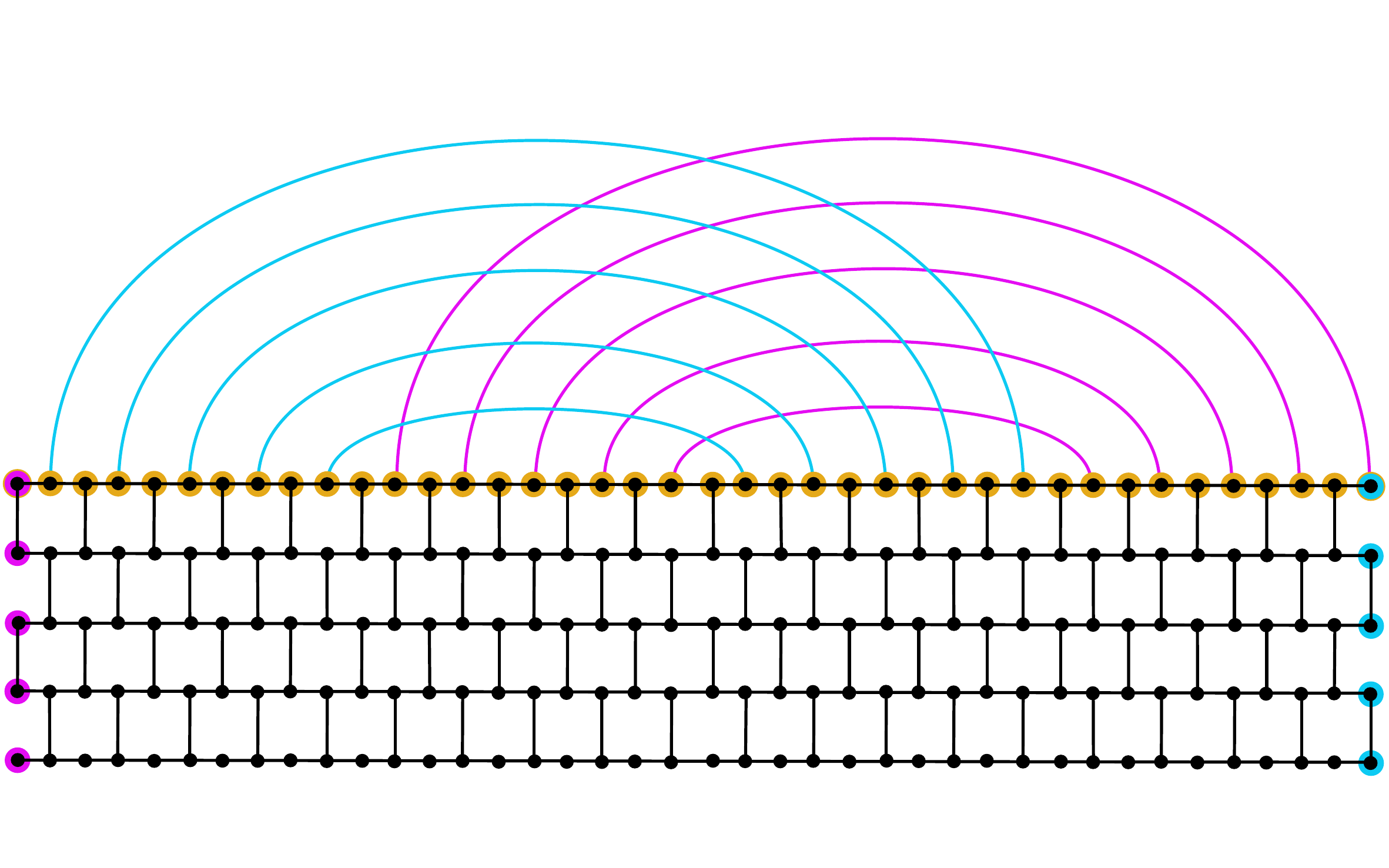}}};
        \node (L) [v:ghost,position=180:55mm from C] {\includegraphics[width=5cm]{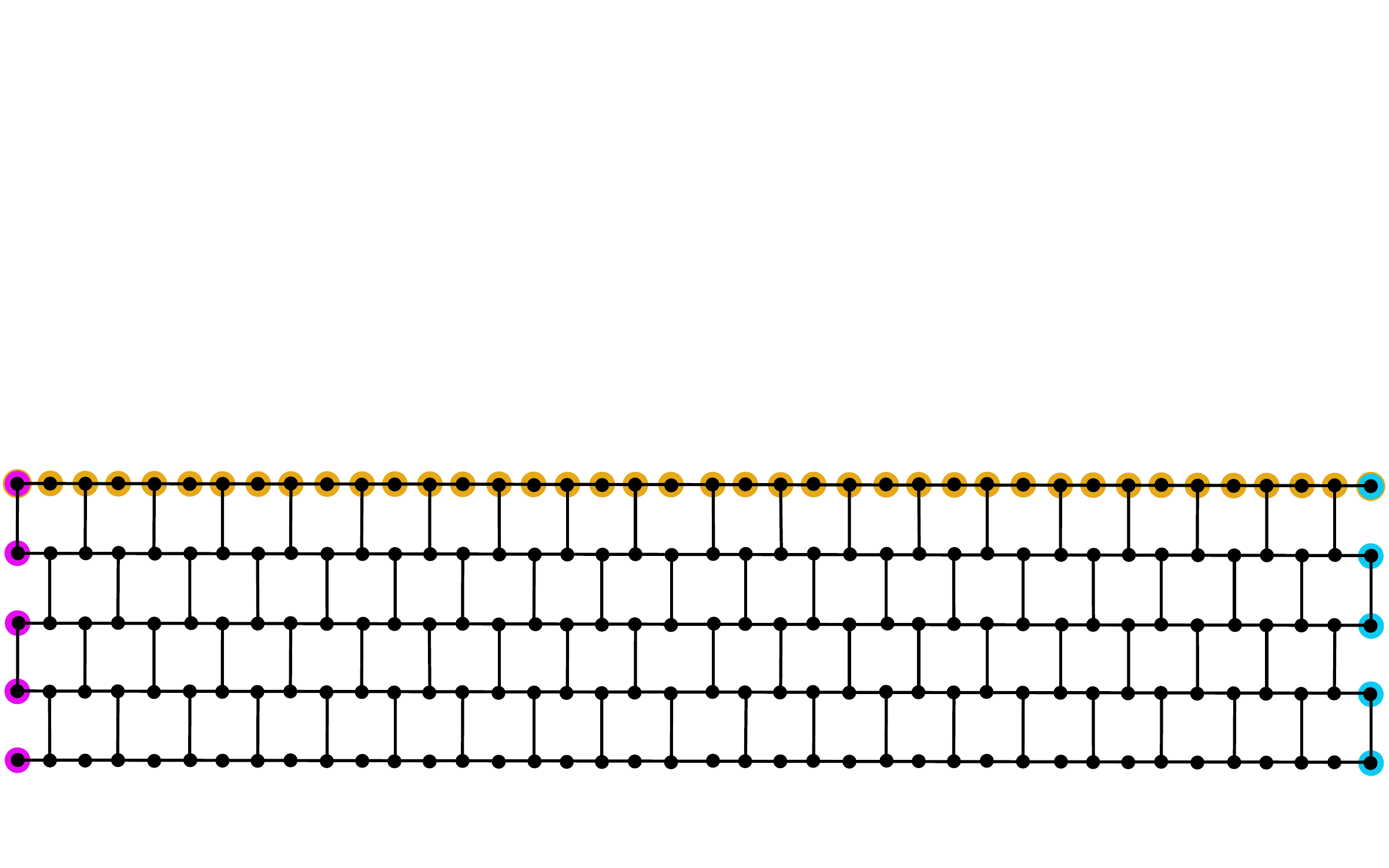}};
        \node (R) [v:ghost,position=0:55mm from C] {{\includegraphics[width=5cm]{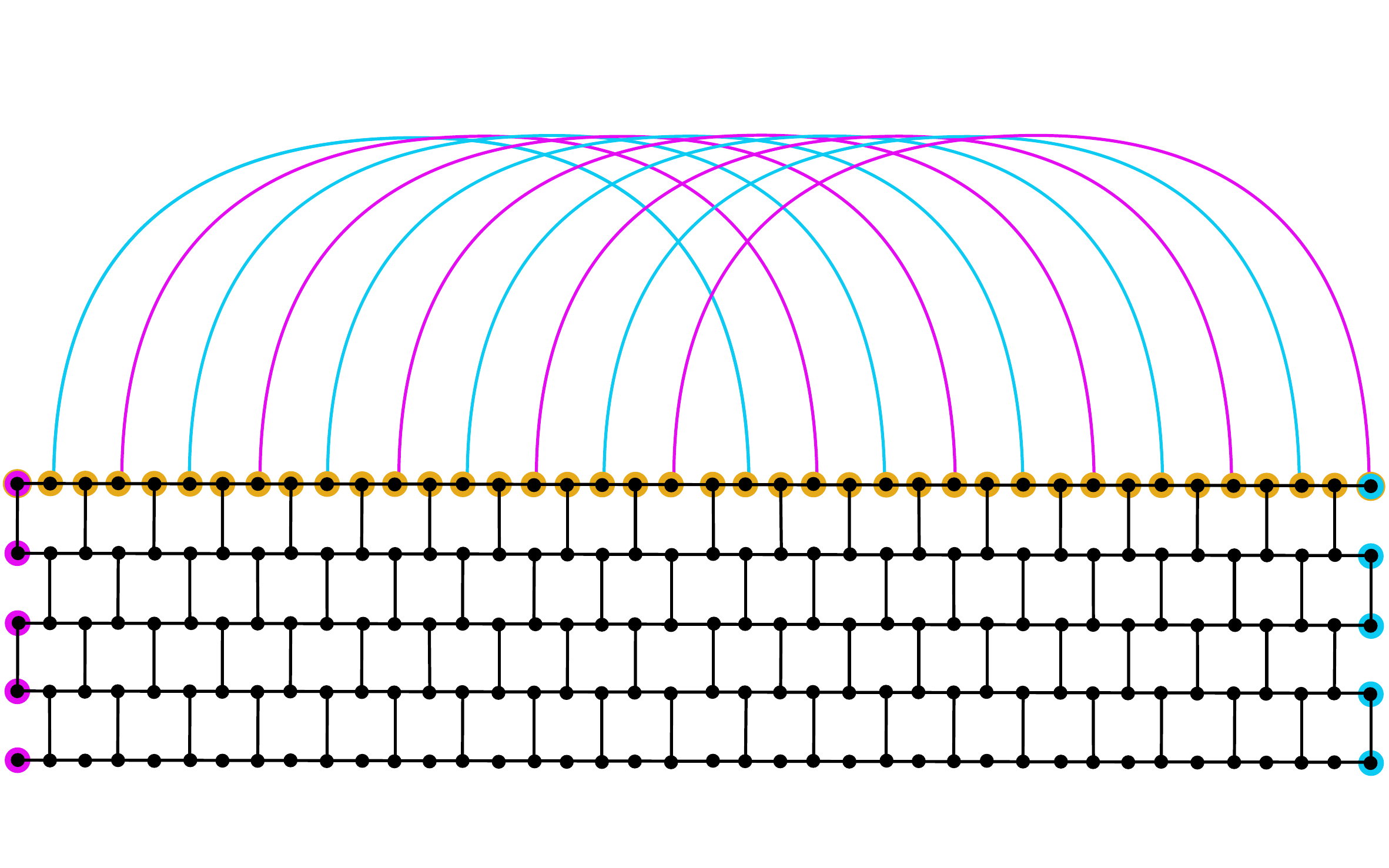}}};
        \end{pgfonlayer}{background}
			
        \begin{pgfonlayer}{main}
        
        \node (C_label) [v:ghost,position=270:18mm from C.center] {\textit{(ii)}};
        \node (L_label) [v:ghost,position=270:18mm from L.center] {\textit{(i)}};
        \node (R_label) [v:ghost,position=270:18mm from R.center] {\textit{(iii)}};
        
        \end{pgfonlayer}{main}

        \begin{pgfonlayer}{foreground}
        \end{pgfonlayer}{foreground}

    \end{tikzpicture}
    \caption{The elementary $5$-wall-segment \textit{(i)}, the elementary $5$-handle-segment \textit{(ii)}, and the elementary $5$-crosscap-segment \textit{(iii)}.}
    \label{fig:basic_segments}
\end{figure}

\paragraph{Wall segments.}
Let $n$ be a positive integer.
An \emph{elementary $n$-wall-segment} is the graph $W_0$ obtained from the $(n\times 8n)$-grid by deleting all edges in the following set
\begin{align*}
    \big\{  \{(i,j),(i+1,j) \}  ~\!\colon\!~ i\in[n-1],\text{ }j\in[8n]\text{, and }i\not\equiv j\mod 2 \big\}.
\end{align*}
The vertices in $\{ (i,1)  ~\!\colon\!~ i\in[n] \}$ are said to be the \emph{left boundary} of the segment, while the vertices in $\{ (i,n)  ~\!\colon\!~ i\in[n] \}$ form the \emph{right boundary of the segment}.
Finally, we refer to the vertices in $\{ (1,i) ~\!\colon\!~ i=2j\text{, }j\in[1,8n] \}$ as the \emph{top boundary}.
See \cref{fig:basic_segments}\textit{(i)} for an illustration.
\smallskip

An \emph{elementary $n$-handle-segment} is the graph obtained from the elementary $n$-wall-segment $W_0$ by adding the edges of the set
\begin{align*}
    &\big\{ \{(1,2i),(1,6n+2-2i)\}  ~\!\colon\!~  i\in[1,n]  \big\}\\
    \cup \ &\big\{ \{(1,2i),(1,8n+2-2i)\}  ~\!\colon\!~ i\in[n+1,2n]  \big\}.
\end{align*}
In \cref{fig:basic_segments}, \textit{(ii)} we present an illustration of the elementary $5$-handle-segment.
\smallskip

An \emph{elementary $n$-crosscap-segment} is the graph obtained from the elementary $n$-wall-segment $W_0$ by adding the edges of the set
\begin{align*}
    \big\{ \{(1,2i),(1,4n+2i)\}  ~\!\colon\!~  i\in[1,2n]  \big\}.
\end{align*}
The elementary $5$-crosscap-segment is depicted in \cref{fig:basic_segments}, \textit{(iii)}.

\begin{figure}[ht]
    \centering
    \begin{tikzpicture}

        \pgfdeclarelayer{background}
		\pgfdeclarelayer{foreground}
			
		\pgfsetlayers{background,main,foreground}

        \begin{pgfonlayer}{background}
        \node (C) [v:ghost] {{\includegraphics[width=12cm]{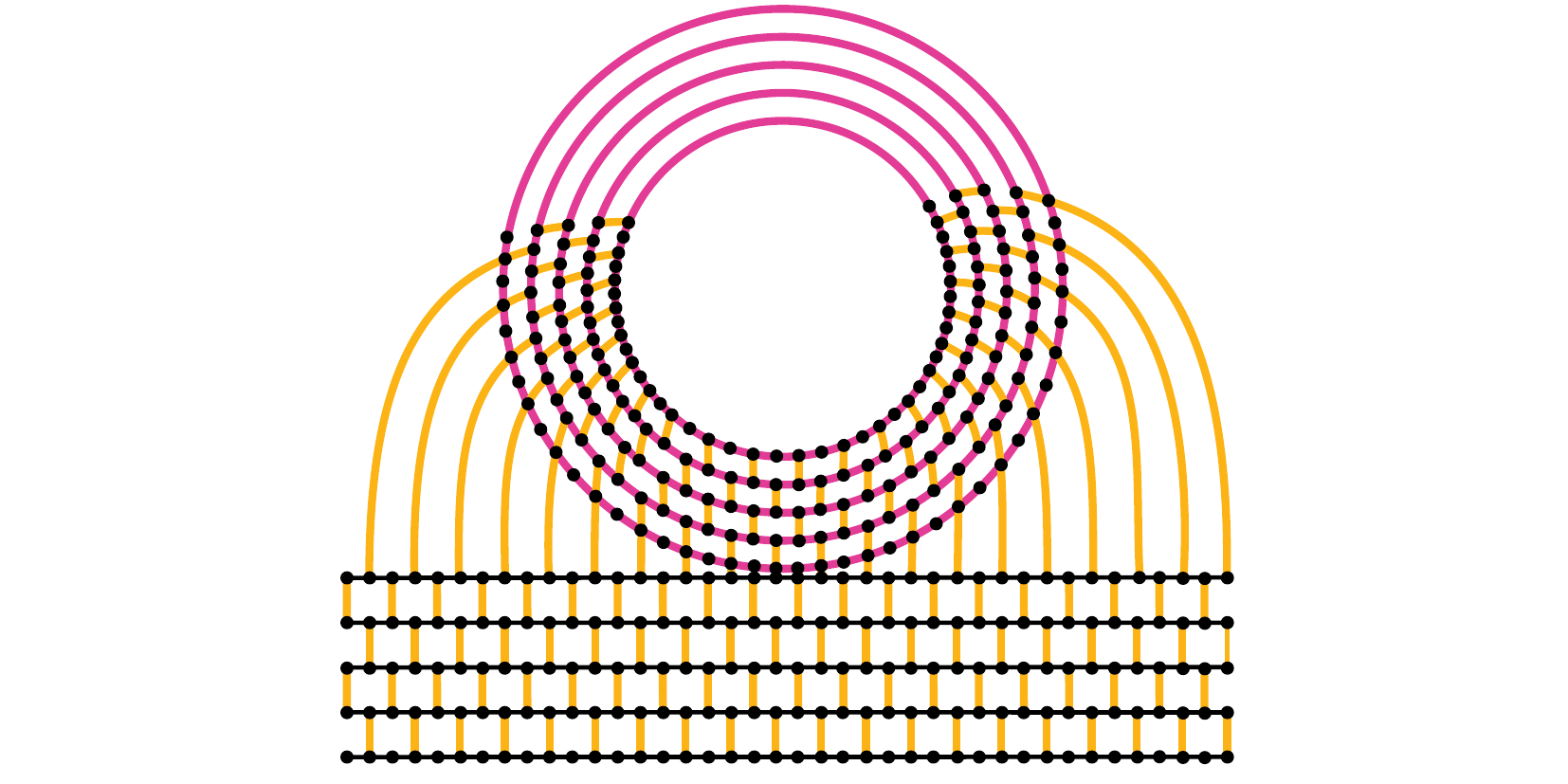}}};
        \end{pgfonlayer}{background}
			
        \begin{pgfonlayer}{main}
        \end{pgfonlayer}{main}

        \begin{pgfonlayer}{foreground}
        \end{pgfonlayer}{foreground}

    \end{tikzpicture}
    \caption{The elementary $5$-vortex-segment. The rails are depicted in \textcolor{BananaYellow}{yellow}, while the nest is highlighted in \textcolor{HotMagenta}{magenta}.}
    \label{fig:vortex_segment}
\end{figure}

An \emph{elementary $n$-vortex-segment} is the graph obtained from the disjoint union of two elementary $n$-wall segments $W_0$ and $W_1$ by making the $i$th top boundary vertex of $W_0$ adjacent to the $i$th top boundary vertex of $W_1$ for each $i \in [4n]$, and by making the $j$th left boundary vertex of $W_1$ adjacent to the $j$th right boundary vertex of $W_1$ for each $j \in [n]$.
See \cref{fig:vortex_segment} for an illustration.

We denote the $(n\times 4n)$-annulus wall defined on the vertex set of $W_1$ as above by $W$.
The base cycle $C_1,\dots C_n$ of $W$ are assumed to be ordered such that $C_n$ contains all vertices adjacent to $W_0$ and we call $\{ C_1,\dots, C_s\}$ the \emph{nest} of the elementary $n$-vortex segment.
We refer to $C_n$ as the \emph{outer cycle} and to $C_1$ as the \emph{inner cycle} of the elementary $n$-vortex segment.
Finally notice that there exist $4n$ pairwise disjoint ``vertical'' paths which are orthogonal to both the horizontal paths of $W_0$ and the cycles $C_1,\dots, C_n$.
The family of these paths is called the \emph{rails} of the elementary $n$-vortex-segment.
\smallskip

In all four types of segments we refer to the elementary wall segment $W_0$ as the \emph{base}.
If we do not want to specify the \textit{type} of an elementary wall-, handle-, crosscap-, or vortex-segment, we simply refer to the graph as a \emph{elementary $n$-segment}, or \emph{elementary segment} if $n$ is not specified.

Let $n$ and $\ell$ be positive integers and let $S_1,\dots,S_\ell$ be elementary $n$-segments.
The \emph{cylindrical closure} of $S_1,\dots,S_{\ell}$ is the graph obtained by introducing, for every $i\in[\ell-1]$ and every $j\in[n]$ and edge between the $j$th vertex of the right boundary of $S_i$ and the $j$th vertex of the left boundary of $S_{i+1}$ together with edges between the $j$th vertex of the right boundary of $S_{\ell}$ and the $j$th vertex of the left boundary of $S_1$.

\paragraph{Surface-walls.}
Let $n$, $h$, and $c$ be non-negative integers.
A \emph{elementary $n$-surface-wall} with \emph{$h$ handles and $c$ crosscaps} is the cylindrical closure of elementary $n$-segments $S_1,\dots,S_{h+c+1}$ such that there is exactly one elementary $n$-wall-segment, $h$ elementary $n$-handle-segments, and $c$ elementary $n$-crosscap-segments among the $S_i$.
An \emph{$n$-surface-wall} is a subdivision of an elementary $n$-surface wall.
We refer to the tuple $(h,c)$ as the \emph{signature} of the $n$-surface wall with $h$ handles and $c$ crosscaps.
\smallskip

Notice that every $n$-surface-wall with signature $(h,c)$ contains an $(n\times (h+c+1))$-annulus wall consisting of $n$ cycles.
We refer to this wall as the \emph{base wall} of the $n$-surface-wall.
Let $C_1,\dots,C_n$ be the base cycles of the base wall.
We will usually assume that $C_1$ is the cycle that contains all top boundary vertices of all elementary segments involved, while $C_n$ can be seen as the ``outermost'' cycle.
We refer to $C_n$ as the \emph{simple cycle} of an $n$-surface-wall.
\medskip

The Euler-genus of an $n$-surface-wall with signature $(h,c)$ is $2h+c$.
By the surface classification theorem \cite{Dyck1888Beitraege,FrancisW1999Conways}, any surface (without boundary) can be generated by adding at most two cross-caps and handles to the sphere.
A similar theorem giving a characterisation of surfaces via the embeddability of surface walls was proven by Thilikos and Wiederrecht in \cite{ThilikosW2024Excluding}.
Their bound, however, is exponential in the Euler-genus which would not work for us as our goal is a polynomial bound.
For this reason, we present a new proof for their result later on (see \cref{subsec:extendeddyckwall}).
\medskip

The \emph{elementary extended $n$-surface wall} with \emph{$h$ handles, $c$ crosscaps, and $b$ vortices} is the graph obtained from the cylindrical closure of $n$-segments $S_1,\dots,S_{h+c+b+1}$ such that there are exactly one elementary $n$-wall-segment, $h$ elementary $n$-handle-segments, $c$ elementary $n$-crosscap-segments, and $b$ elementary $n$-vortex-segments among the $S_i$.
An \emph{extended $n$-surface-wall} is a subdivision of an elementary $n$-surface-wall.
We refer to the tuple $(h,c,b)$ as the \emph{signature} of the extended $n$-surface-wall with $h$ handles, $c$ crosscaps, and $b$ vortices.
Moreover, we adapt the notions of the \emph{base wall} and the \emph{simple cycle} from $n$-surface-walls to extended $n$-surface-walls in the natural way.

\paragraph{Dyck-walls.}
Let $n$, $g$ be non-negative integers.
The \emph{elementary $n$-Dyck-wall of orientable Euler-genus $g$}, where $g=2h$, is the elementary $n$-surface wall with signature $(h,0)$.
The \emph{elementary $n$-Dyck-wall of non-orientable Euler-genus $g$} is the elementary $n$-surface wall with signature $(0,g)$.

Notice that each Dyck-wall corresponds to a fixed surface $\Sigma$.
So we will sometimes write \emph{$n$-$\Sigma$-wall} for an $n$-Dyck-wall of (orientable/non-orientable) Euler-genus $g$ where $\Sigma$ is the surface of (orientable/non-orientable) Euler-genus $g$.

\begin{figure}[ht]
    \centering
    \begin{tikzpicture}

        \pgfdeclarelayer{background}
		\pgfdeclarelayer{foreground}
			
		\pgfsetlayers{background,main,foreground}

        \begin{pgfonlayer}{background}
        \node (C) [v:ghost] {{\includegraphics[width=6cm]{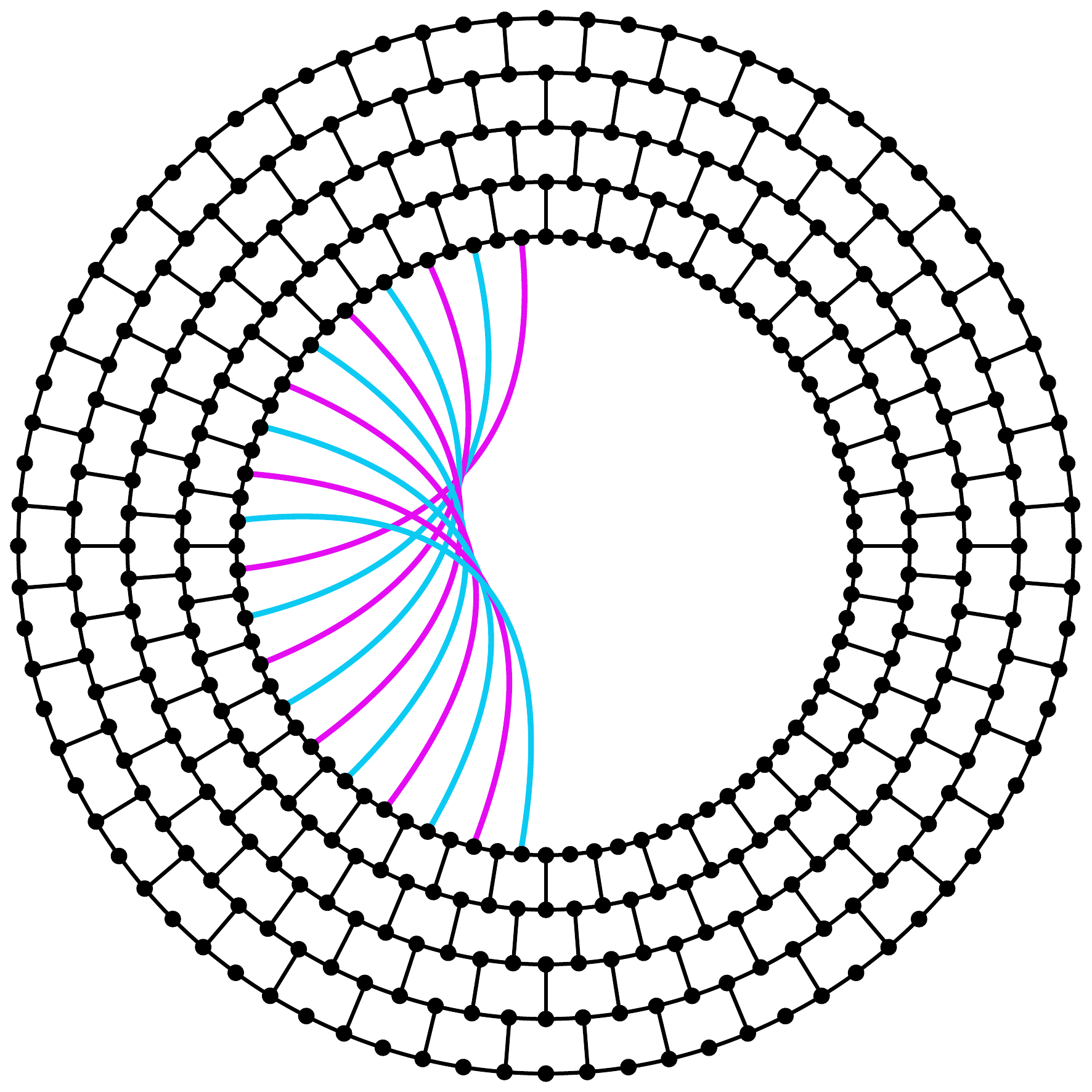}}};
        \node (L) [v:ghost,position=180:65mm from C] {\includegraphics[width=6cm]{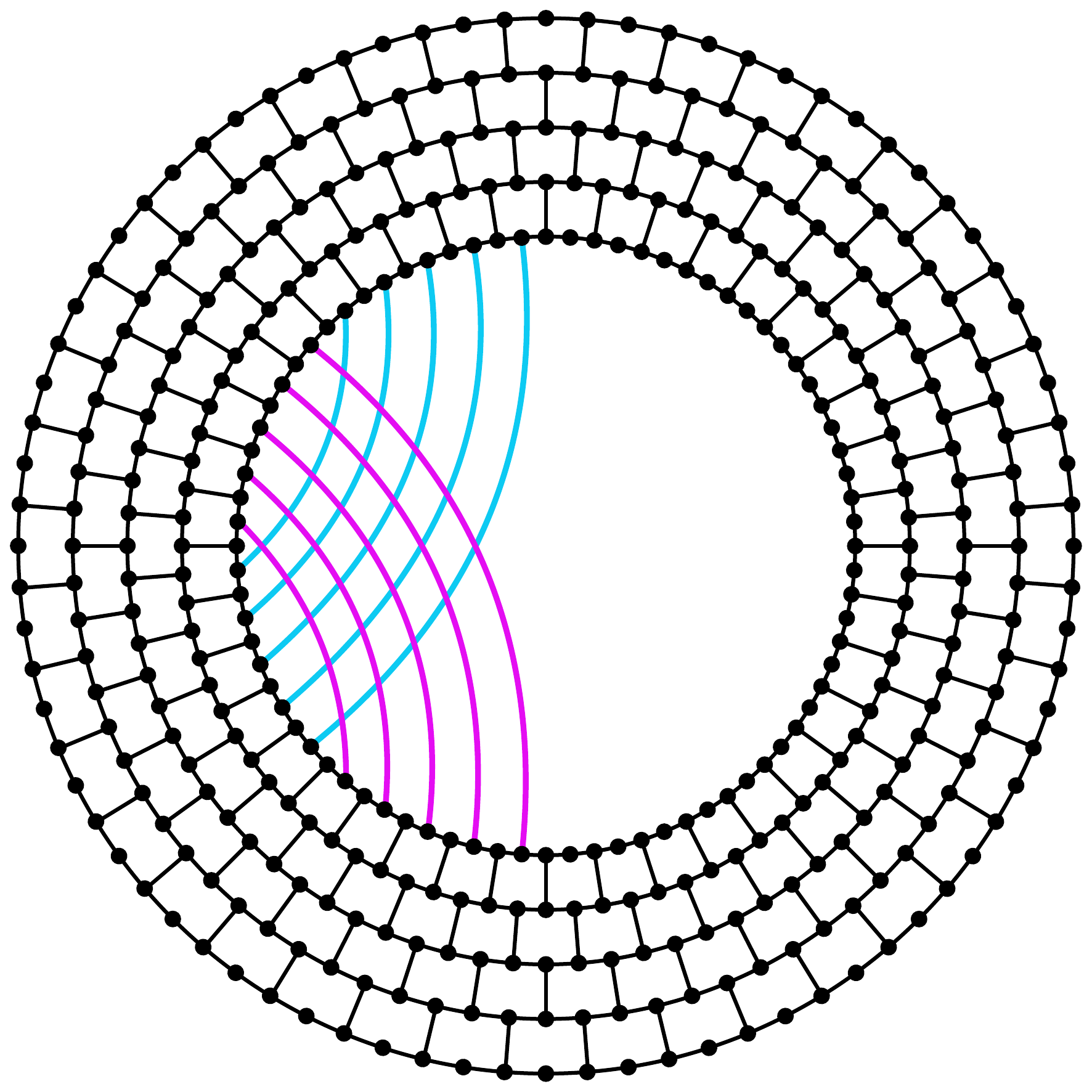}};
        \end{pgfonlayer}{background}
			
        \begin{pgfonlayer}{main}
        
        \node (C_label) [v:ghost,position=270:35mm from C.center] {\textit{(ii)}};
        \node (L_label) [v:ghost,position=270:35mm from L.center] {\textit{(i)}};
        
        \end{pgfonlayer}{main}

        \begin{pgfonlayer}{foreground}
        \end{pgfonlayer}{foreground}

    \end{tikzpicture}
    \caption{The $5$-Dyck-walls corresponding to the torus \textit{(i)} and the projective plane \textit{(ii)}.}
    \label{fig:dyck_walls}
\end{figure}

We extend the notions of $n$-surface-walls of, signatures, and the simple cycle to (orientable/non-orientable) $n$-Dyck-walls in the natural way.
See \cref{fig:dyck_walls} for an illustration.

Notice that there are exactly two non-isomorphic elementary $n$-Dyck-walls of Euler-genus $2g$ for every positive integer $g$, there is exactly one elementary $n$-Dyck-wall of odd Euler-genus $g$ for any $g$, and the elementary $n$-Dyck-wall of Euler-genus $0$ is the elementary $(n\times 4n)$-annulus wall.

It follows from a recent result of Gavoille and Hilare \cite{GavoilleH2023MinorUniversal} that a $\mathbf{O}_g(n^2)$-Dyck-wall with signature $(h,c)$, $2h+c=g$, contains every graph on at most $n$ vertices which embeds in the surface with $h$ handles and $c$ crosscaps as a minor.

\begin{proposition}[Gavoille and Hilare \cite{GavoilleH2023MinorUniversal}]\label{prop:dyckwall_minors}
There exists a universal constant $c_{\ref{prop:dyckwall_minors}}>0$ such that for every positive integer $n$ and every non-negative integer $g=2h+c$, $h,c\in\mathbb{N}$, every graph on at most $n$ vertices which embeds in the surface $\Sigma$ with $h$ handles and $c$ crosscaps is a minor of the $(c_{\ref{prop:dyckwall_minors}}\cdot g^2(n+g)^2)$-$\Sigma$-wall.
\end{proposition}

\subsection{Landscapes and layouts}\label{subsec:landscapes}

We are finally ready to define the full data structure that describes our local theorem.

\paragraph{Landscapes.}
Let $k,w\geq 4$ be integers, let $G$ be a graph, and let $\Sigma$ be a surface of Euler-genus $g$.
Let $h$, $c$, and $b$ be a non-negative integers where $g=2h+c$ and $c\neq 0$ if and only if $\Sigma$ is non-orientable.
Moreover, let $D\subseteq G$ be an extended $k$-surface-wall with signature $(h,c,b)$, and let $W\subseteq G-A$ be a $w$-mesh in $G$.
Finally, let $A\subseteq V(G)\setminus V(D)$.
The tuple $\Lambda=(A,W,D,\delta)$ is called a \emph{$\Sigma$-landscape} of \emph{detail $k$}
\begin{description}
    \item[L1~~] $\delta$ is a $\Sigma$-decomposition for $G-A$,
    \item[L2~~] $D$ and $W$ are grounded in $\delta$,
    \item[L3~~] the trace of the perimeter of $W$ bounds a disk that contains all grounded vertices of $W$ while being disjoint from $D$,
    \item[L4~~] the trace of the simple cycle of $D$ bounds a disk that avoids the traces of the other base cycles of $D$,
    \item[L5~~] the tangle $\mathcal{T}_D$ is a truncation of the tangle $\mathcal{T}_W$,
    \item[L6~~] if $C$ is a cycle from the nest of some vortex-segment of $D$, then the trace of $C$ is a contractible closed curve in $\Sigma$,
    \item[L7~~] $\delta$ has exactly $b$ vortices and there exists a bijection between the vortices $v$ of $\delta$ and the vortex segments $S_v$ of $D$ such that $v$ is the unique vortex of $\delta$ which is contained in the disk $\Delta_{C_1}$ defined by the trace of the inner cycle of $S_v$ where $\Delta_{C_1}$ avoids the trace of the simple cycle of $D$, and
    \item[L8~~] for every vortex $v$ of $\delta$, the society $(G_{\Delta^v_k},\Omega_{\Delta^v_k})$, induced by the outer cycle from the nest of the corresponding vortex segment has a cross.
\end{description}
We refer to the vortex segments as the \emph{vortices} of $\Lambda$.
The integer $b$ is called the \emph{breadth} of $\Lambda$ and we refer to the maximum depth over the depths of the vortex societies of the vortices of $\delta$ as the \emph{depth} of $\Lambda$.
The set $A$ is called \emph{apex set}.
We also say that $\Lambda$ is \emph{centred} at the mesh $W$.

\paragraph{Layouts.}
Let $k\geq 4$, $d$, $b$, $r$, and $a$ be non-negative integers and $\Sigma$ be a surface.
We say that a graph $G$ with a mesh $M$ has a \emph{$k$-$(a,b,d,r)$-$\Sigma$-layout} \emph{centred at $M$} if there exists a set $A\subseteq V(G)$ of size at most $a$ and a submesh $M'\subseteq M$ such that there exists a $\Sigma$-landscape $(A,M',D,\delta)$ of detail $k$ for $G$ where every vortex of $\delta$ has a linear decomposition of adhesion at most $d$, and $M'$ is a $w$-mesh with $w\geq a+b(2d+1)+6+r$.
\medskip

The order of the mesh $W$ will ensure that the tangle of $W$ ``respects'' the surface.
That is, any separation created by one of the cells of $\delta$ or by the linear decomposition of one of its vertices must be oriented \textsl{towards} $W$.
This will come into play much later when we derive the original GMST of Robertson and Seymour.


\section{Forcing a model of \texorpdfstring{\(K_t\)}{Kₜ}}\label{sec:forcing}

In the proof of the Graph Minor Structure Theorem, a model of \(K_t\) is not built directly. Instead, we find either (1) a model of an intermediate graph which contains \(K_t\) as a minor, or (2) a large mesh and a certain configuration of paths such that the union of the mesh and the paths contains \(K_t\) as a minor.
In this section, we define these two structures, and show that they indeed get a model of \(K_t\), and the model of \(K_t\) is controlled by a certain mesh.

\subsection{Models of \texorpdfstring{\(K_t\)}{Kₜ} forced by crosses}

For integers \(w, h \ge 1\) with \(h\) even, let \(H_{w,h}\) denote the graph obtained from the \((w \times h)\)-grid by adding two crossing edges \((i, h/2)(i+1, h/2+1)\) and \((i, h/2+1)(i+1, h/2)\) for each \(i \in [w-1]\). In order to bound the number of vortices in the graph minor structure theorem by \(t(t-1)\) Kawarabayashi et al.\ proved the following.

\begin{figure}
    \centering
    \includegraphics{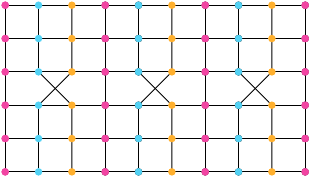}
    \caption{The graph \(H'_{3, 6}\)}
    \label{fig:crossed-grid}
\end{figure}

\begin{proposition}[Kawarabayashi, Thomas, and Wollan~\cite{KawarabayashiTW2018New}]\label{prop:cross_row_grid}
    Let $t\geq 2$ be an integer. The graph $H_{t(t-1),t(t-1)}$ has a $K_t$-model controlled by the \((t(t-1) \times t(t-1))\)-grid contained in it.
\end{proposition}
In this subsection we prove a variant of this result, which will allow us to reduce the final number of vortices from \(t(t-1)\) to \((t-3)(t-4)/2\).
For integers \(c, h \ge 1\) with \(h\) even, let \(H'_{c, h}\) denote the graph obtained from the \((3c+1 \times h)\)-grid by deleting all edges \((i, h/2)(i, h/2+1)\) with \(i \not \equiv 1 \pmod 3\), and adding, for each \(i \in [3c]\) with \(i \equiv 2 \pmod 3\), two ``crossing'' edges \((i, h/2)(i+1, h/2+1)\) and \((i, h/2+1)(i+1, h/2)\) (see \cref{fig:crossed-grid}).
The two \((w \times h/2)\)-subgrids of the \((w \times h)\)-grid induced by \([w] \times [t]\) and \([w] \times [t+1,2t]\) are called the \emph{underlying grids of \(H'_{w, 2t}\)}. We show the following.

\begin{lemma}\label{lem:Kt_from_crosses}
    Let \(t \ge 5\) be an integer, let \(c=(t-3)(t-4)/2\), and let \(h=2t\). The graph \(H'_{c, h}\) has a model of \(K_t\) controlled by each of the underlying grids, and such a model can be found in \(\mathbf{O}(t^2)\) time.
\end{lemma}
\begin{proof}
    Let \(w'=3t-11\), let \(h'=2t-6\), and let \(H'\) denote the graph obtained from the \((w' \times h')\)-grid by adding all possible edges \((i, j)(i', j')\) with \(|i-i'|=1\) and \(|j - j'| = 1\).
    We start from defining a model \(\mu\) of \(K_t\) in \(H'\), see \cref{fig:kt_in_strong_product}.
    \begin{figure}
        \centering
        \includegraphics{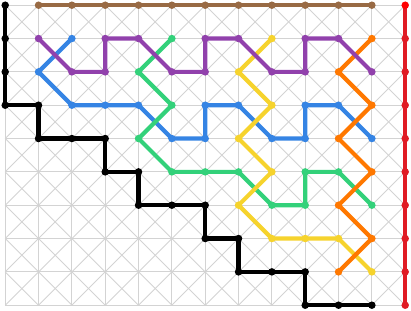}
        \caption{A construction of a \(K_8\)-model in the \((13 \times 10)\)-grid with crosses.}
        \label{fig:kt_in_strong_product}
    \end{figure}
    
    Each branch set of \(\mu\) induces a path in the graph.
    Only \(t-3\) branch sets make an actual use of the crossing edges , and we call these \emph{crossing} branch sets. The remaining three branch sets are called \emph{outer} as their union induces a cycle ``surrounding'' the remaining branch sets.
    
    For each \(a \in [t-3]\), the \(a\)th crossing branch set is contained in an associated \emph{region}. For \(a = 1\), the region is \([2,w'-1 ] \times [2, 3]\),  for \(a \in [2,t-4]\), the region is \([3a-4, 3a-3] \times [2, 2a+1] \cup [3a-4, w'-1] \times [2a, 2a+1]\), and for \(a=t-3\), the region is \([w'-2, w'-1] \times [2, h'-1]\). The intersection of the regions of two distinct crossing branch sets is \([i,i+1] \times [j,j+1]\) for some \((i, j) \in [w'-1] \times [h'-1]\)
    with \(i \equiv 2 \pmod 3\) and \(j \equiv 2 \pmod 2\),
    and each of the two branch sets contains one of the crossing edges \((i, j)(i+1,j+1)\) and \((i, j+1)(i+1, j)\).

    The three outer branch sets are \([2, w'-1] \times [1]\), \(\{w'\} \times [h']\),
    and the third branch set which is adjacent to all remaining branch sets (the black path in \cref{fig:kt_in_strong_product}).

    Observe that the model \(\mu\) of \(K_t\) is controlled by the \((w' \times h')\)-grid  --  this is due to the fact that every branch set of \(\mu\) intersects at least \(t\) vertical paths or at least \(t\) horizontal paths of the grid.

    We observe that the model \(\mu\) constructed as above can be embedded in \(H_{c, h}'\). This model makes an essential use only of \((t-3)(t-4)/2\) pairs of crossing edges \((i, j)(i+1,j+1), (i, j+1)(i+1, j)\), and we can ``stretch'' the model so that it uses the
    \(c = (t-3)(t-4)/2\) crosses in \(H'_{c, h}\) (see \cref{fig:kt_in_separated_crosses}).
    Observe that every branch set of the model of \(K_t\) intersects more than \(t\) vertical paths of \(H_{c, 2t}'\), so it is controlled by the
    two underlying \((3c+1 \times t)\)-grids.\qedhere
    
    \begin{figure}
        \centering
        \includegraphics{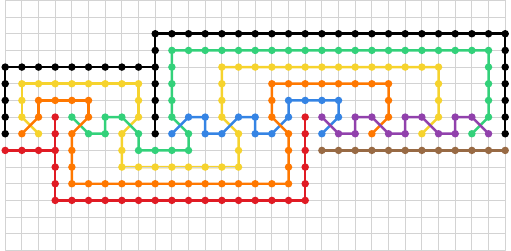}
        \caption{A model of \(K_8\) in \(H_{10,16}'\)}
        \label{fig:kt_in_separated_crosses}
    \end{figure}
\end{proof}

\subsection{Models of \texorpdfstring{\(K_t\)}{Kₜ} forced by jumps}
Let \(M\) we be a \((w \times 2t)\)-mesh with vertical paths \(P_1, \ldots, P_w\)
and horizontal paths \(Q_1, \ldots, Q_{2t}\). A \emph{middle vertex} in \(M\) is an inner vertex of the \(Q_t\)--\(Q_{t+1}\) subpath of a vertical path \(P_i\).
Two middle vertices \(x\) and \(y\) with \(x \in V(P_i)\) and \(y \in V(P_j)\)
are \emph{independent} in \(M\) if \(|i-j| > 1\).
An \emph{\(M\)-middle jump} is a \(V(M)\)-path \(R\) between independent middle vertices.
In this subsection we prove the following. 
\begin{lemma}\label{KtFromJumps}
    Let \(w \ge 1, t \ge 5\) be an integer, and let \(M\) be a \((w \times 2t)\)-mesh in \(G\).
    If there exist \(t^3 - 2\) pairwise vertex-disjoint \(M\)-middle jumps \(R_1, \ldots, R_{t^3-2}\) in \(G\) with pairwise independent endpoints, then \(M\) controls a model of \(K_t\) in \(G\).

    Moreover, such a model can be found in \(\mathbf{O}(|E(G)|)\) time.

\end{lemma}
Our proof actually gives a \(t^3 - O(t^2)\) bound which is slightly better than \(t^3 - 2\).
We choose the stated bound for the sake of convenience in the context of our applications of it in some of the later lemmas.
To prove our lemma, we need two lemmas about sets of intervals, and a lemma about finding linkages in a grid.

In what follows, intervals have integer endpoints.
For an integer \(d \ge 2\), we say that a set of intervals \(\mathcal{I}\) is \emph{\(d\)-independent} if for any distinct \(I, J \in \mathcal{I}\), any endpoint \(i\) of \(I\), and any endpoint \(j\) of \(J\) we have \(|i-j| \ge d\).

\begin{lemma}\label{lem:disjoint-short-intervals}
    Let \(k \ge 1, \ell \ge 4\) be integers.
    Let \(\mathcal{I}\) be a \(2\)-independent set of intervals with \(|\mathcal{I}| \ge \ell(k-1)/2 + 1\) such that for each \([i, j] \in \mathcal{I}\), we have \(2 \le j - i \le \ell\). Then \(\mathcal{I}\) contains \(k\) pairwise disjoint intervals.

    Moreover, such intervals can be found in \(\mathbf{O}(\ell k^2)\) time.
\end{lemma}
\begin{proof}
    After discarding some of the intervals, we assume that \(|\mathcal{I}| = \ell(k-1)/2 + 1\).

    We prove the lemma by induction. In the base case \(k=1\), the lemma holds trivially since \(|\mathcal{I}| \ge \ell(k-1)/2+1>0\).

    For the induction step, suppose that \(k \ge 2\). Let \([i, j] \in \mathcal{I}\) be the interval with smallest \(i\), and partition \(\mathcal{I}\) into sets \(\mathcal{I}_0\) and \(\mathcal{I}_1\) by assigning each \([i', j'] \in \mathcal{I}\) to \(\mathcal{I}_0\) if \([i', j'] \cap [i, j] = \emptyset\),
    and to \(\mathcal{I}_1\) otherwise. Since \(\mathcal{I}\) is \(2\)-independent, every interval \([i', j'] \in \mathcal{I}_1 \setminus \{[i, j]\}\), has an endpoint in \([i+2,j-2]\). Hence,
    there are at most \(((j-2) - (i+2))/2 + 1\) such intervals, so
    \[|\mathcal{I}_1| = |\mathcal{I}_1 \setminus \{[i, j]\}|+1 \le (((j-2) - (i+2))/2+1)+1 \le (\ell-4)/2+1+1 = \ell/2.\] Thus, \[|\mathcal{I}_0| \ge |\mathcal{I}| - |\mathcal{I}_1| \ge \ell(k-1)/2 + 1 - \ell/2 = \ell(k - 2)+1,\] and by the induction hypothesis, the set \(\mathcal{I}_0\) contains \(k-1\) pairwise disjoint intervals. Together with \([i, j]\), they make \(k\) pairwise disjoint intervals in \(\mathcal{I}\), completing the inductive proof.
\end{proof}
\begin{lemma}\label{lem:d-independent-intervals}
    Let \(k \ge 1, d \ge 1\) be integers, and let \(\mathcal{I}\) be a \(2\)-independent set of intervals with \(|\mathcal{I}| \ge (3d-1)(k-1)/2+1\). Then there exists a \(d\)-independent subset \(\mathcal{I}' \subseteq \mathcal{I}\) with \(|\mathcal{I}'| = k\).

    Moreover, such a subset can be found in \(\mathbf{O}(dk^2)\) time.
\end{lemma}
\begin{proof}
    After discarding some of the intervals, we may assume that \(|\mathcal{I}| = \lceil(3d-1)(k-1)/2+1\rceil\).

    We prove the lemma by induction, similarly to the previous lemma. In the base case \(k=1\), the lemma holds trivially since \(|\mathcal{I}| \ge \ell(k-1)/2+1>0\).
    
    For the induction step, suppose that \(k \ge 2\). Let \([i, j] \in \mathcal{I}\) be the interval with the smallest \(i\), and partition \(\mathcal{I}\) into sets \(\mathcal{I}_0\) and \(\mathcal{I}_1\) by assigning each \([i', j'] \in \mathcal{I}\) to \(\mathcal{I}_0\) if
    \([i, j] \neq [i', j']\) and
    the set \(\{[i, j], [i', j']\}\) is \(d\)-independent, and otherwise assigning it to \(\mathcal{I}_1\).
    Since \(\mathcal{I}\) is \(2\)-independent, every interval \([i', j'] \in \mathcal{I}_1 \setminus \{[i, j]\}\) has an endpoint in one of the intervals \([i+2, i+d-1]\), \([j-(d-1), j-2]\) and \([j + 2, j+d-1]\). Each of these intervals contains endpoints of at most \((d-3)/2+1 = (d-1)/2\) of intervals from \(\mathcal{I}_1 \setminus \{[i, j]\}\), so
    \[
      |\mathcal{I}_1| = |\mathcal{I}_1 \setminus \{[i, j]\}|+1\le 3(d-1)/2+1 = (3d-1)/2,
    \]
    and thus
    \[
      |\mathcal{I}_0| = |\mathcal{I}| - |\mathcal{I}_1| \ge (3d-1)(k-1)/2+1 - (3d-1)/2 = (3d-1)(k-2)/2+1.
    \]
    By the induction hypothesis, the set \(\mathcal{I}_0\) contains \(k - 1\) intervals forming a \(d\)-independent set.
    These intervals together with \([i, j]\) satisfy the lemma.
\end{proof}

\begin{lemma}\label{lem:ltorlinkage}
    Let \(w, t \ge 1\) be integers with \(w \ge t+2\), and let \(h = 2t\).
    Let \(M\) be a \((w \times h)\)-mesh with vertical paths \(P_1, \ldots, P_w\) and horizontal paths \(Q_1, \ldots, Q_h\). Then, for any integers
    \(j', j'' \in [0, t]\), there exist pairwise disjoint \(P_1\)--\(P_w\) paths \(S_1, \ldots, S_t\).
    in \(M\) such that each \(S_a\) has endpoints in \(P_1 \cap Q_{j'+a}\) and \(P_w \cap Q_{j''+a}\).

    Moreover, such paths can be found in \(\mathbf{O}(wt)\) time.
\end{lemma}
\begin{proof}
    After possibly reversing the order of the vertical paths \(Q_1, \ldots, Q_h\), we may assume that \(j' \ge j''\). For each \(a \in [t]\), let \(S_a\)
    be the path contained in \(Q_{j'+a} \cup P_{a+1} \cup Q_{j''+a}\)
    between \(P_1 \cap Q_{j'+a}\) and \(P_w \cap Q_{j''+a}\).
    Such paths \(S_1, \ldots, S_t\) satisfy the lemma (see \cref{fig:ltorlinkage}).
    \begin{figure}
        \centering
        \includegraphics{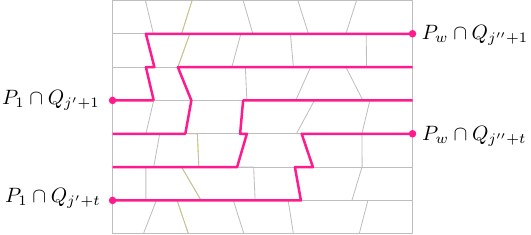}
        \caption{The construction of the paths \(S_1, \ldots, S_t\) in \cref{lem:ltorlinkage}}
        \label{fig:ltorlinkage}
    \end{figure}
\end{proof}

\begin{proof}[Proof of \cref{KtFromJumps}]
    We need to show that if a graph \(G\) contains a \((w \times 2t)\)-mesh
    \(M\) and \(t^3 - 2\) pairwise disjoint \(M\)-middle jumps with pairwise independent endpoints, then \(M\)
    controls a model of \(K_t\) in \(G\).
    Let \(P_1, \ldots, P_w\) be the vertical paths of \(M\). For an \(M\)-middle jump \(R\) with endpoints on paths \(P_i\) and \(P_j\) with \(1 \le i < j \le w\), we define an interval \(I(R) \coloneqq  [i, j]\). By definition of \(M\)-middle jump, we have \(j - i \ge 2\), and we call \(R\) a \emph{short jump} if \(j-i \le t\), and a \emph{long jump} if
    \(j-i \ge t+1\). Since there are at least \(t^3 - 2\) pairwise disjoint \(M\)-middle jumps in \(G\), there are at least \(t^3/4-1\) pairwise disjoint short jumps or at least \(3t^3/4-1\) pairwise disjoint long paths in \(G\).

    Suppose that there are pairwise disjoint short jumps \(R_1, \ldots, R_m\) in \(G\) where \(m \ge t^3/4-1\).
    At most two of these
    paths have an endpoint on \(P_1\) or \(P_w\), so we may assume that
    \(R_1, \ldots, R_{m-2}\) are disjoint from \(P_1\) and \(P_w\).
    the
    Let \(c = (t-3)(t-4)/2\), let \(\ell = t\), and note that \(m - 2 \ge t^3/4-3 \ge \ell(c-1)/2+1\).
    By \cref{lem:disjoint-short-intervals} applied to the intervals \(I(R_1), \ldots, I(R_m)\), there exists a subset \(\{R_1', \ldots, R_{c}'\} \subseteq \{R_1, \ldots, R_m\}\) such that the intervals \(I(R_1'), \ldots, I(R_{c}')\) are pairwise disjoint.
    For each \(\alpha \in [c]\), let \([i_\alpha, j_\alpha] = I(R_{\alpha}')\). Without loss of generality, the intervals are sorted,
    that is, \(2 \le i_1 < j_1 < \cdots < i_c < j_c \le w - 1\).
    Note that for each \(\alpha \in [c]\) we have \(j_\alpha - i_\alpha \ge 2\), and for each \(\alpha \in [c-1]\) we have \(i_{\alpha+1} - j_\alpha \ge 2\).
    
    We will now use the mesh together with the jumps to build a model \(\eta\) of \(H'_{c, 2t}\) (see \cref{fig:crossedgridmodel}).
    \begin{figure}
        \centering
        \includegraphics{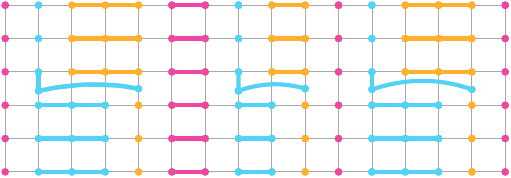}
        \caption{A model of \(H'_{3, 6}\) in a mesh \(M\) with three \(M\)-middle jumps. The mesh \(M\) is the \((w \times 6)\)-grid with some of the edges \((i,3)(i,4)\) subdivided once.}
        \label{fig:crossedgridmodel}
    \end{figure}    
    Let \(\eta_0\) be a model of the \((w \times 2t)\)-grid in \(M\) such that for any \(i \in [w]\), \(j \in [2t]\) we have \(P_i \cap Q_j \subseteq \eta_0((i, j))\). We assume that for each \(i \in \{i_1, j_1, \ldots, i_c, j_c\}\), if \(x\) is the endpoint of a jump \(R_\alpha\) on \(P_i\), then all vertices of the \(Q_t\)--\(x\) subpath of \(P_i\) belong to \(\eta_0((i, t))\), and the remaining vertices of the \(Q_t\)--\(Q_{t+1}\) subpath of \(P_i\) belong to \(\eta_0((i, t+1))\).

    We construct the model \(\eta\) of \(H'_{c, 2t}\) as follows.
    For any \(\alpha \in [c+1]\) and \(j \in [2t]\), we define \(\eta((3\alpha-2, j)) = \bigcup_{i=j_{\alpha-1}+1}^{i_\alpha-1}\eta_0((i, j))\) where \(j_0 = 0\) and \(i_{c+1} = w+1\). The definitions of \(\eta((3\alpha-1, j))\) and \(\eta((3\alpha, j))\) depend on the value of \(j\).
    If \(j \le t - 1\), then \(\eta((3\alpha-1, j)) = \eta_0((i_\alpha, j))\) and \(\eta((3\alpha-1, j)) = \bigcup_{i=i_{\alpha}+1}^{j_\alpha}\eta_0((i, j))\).
    If \(j = t\), then \(\eta((3\alpha-1, j)) = \eta_0(i_\alpha, j) \cup V(R_\alpha)\), and \(\eta((3\alpha, j)) = \bigcup_{i=i_\alpha+1}^{j_\alpha}\eta_0(i, j) \setminus V(R_\alpha)\).
    If \(j \ge t + 1\), then \(\eta((3\alpha-1, j)) = \bigcup_{i=i_\alpha}^{j_\alpha - 1}\eta_0((i, j))\) and \(\eta((3\alpha, j)) = \eta_0((j_\alpha, j))\). It is straightforward to verify that
    \(\eta\) is indeed a model of \(H'_{c, 2t}\). By \cref{lem:Kt_from_crosses}, there is a model \(\mu_0\) of \(K_t\) in
    \(H'_{c, 2t}\) which is controlled by the underlying grids of \(H'_{c, 2t}\). Thus, for any \(u \in V(K_t)\) and \(j \in [2t]\) there exist \(t\) pairwise internally disjoint paths between \(\mu_0(u)\) and \(\bigcup_{i=1}^w \eta_0((i, j))\). Consequently, in the graph \(G\) there exist \(t\) pairwise internally disjoint paths, each between some \(\eta_0(x)\) with \(x \in \mu_0(u)\) and \(P_i\). Therefore we can
    lift \(\mu_0\) to a model \(\mu\) of \(K_t\) controlled by \(M\) in \(G\) by defining for each \(u \in V(K_t)\), \(\mu(u) = \bigcup_{x \in \mu_0(u)} \eta(x)\).

    Now, suppose that there exist pairwise disjoint long jumps \(R_1, \ldots, R_m\) with pairwise independent endpoints, where \(m \ge 3t^3/4-1\). Let \(k = t(t-1)/2\) and \(d = t+1\), and observe that \(m \ge (3d-1)(k-1)/2+1\).
    By \cref{lem:d-independent-intervals} applied to the intervals \(I(R_1), \ldots, I(R_m)\), there exists a subset \(\{R_1', \ldots, R_k'\}\) such that for any distinct \(i, j \in [w]\) if each of \(P_i\)
    and \(P_j\) contains an endpoint of one of the paths \(R_1', \ldots, R_k'\), then \(|i-j|\ge t\).
    We will next use \cref{lem:ltorlinkage} to find \(P_1\)--\(P_w\) paths \(S_1, \ldots, S_t\) such that they can be pairwise connected with \(\binom{t}{2}\) disjoint paths, each obtained by extending a distinct jump \(R_\alpha'\).
    
    Let \(i_1, \ldots, i_{2k}\) be the indices of the vertical paths \(P_i\) containing the endpoints of \(R_1', \ldots, R_k'\), listed in the increasing order. For each \(\alpha \in [2k-1]\), let \(w_\alpha = i_{\alpha+1} - i_\alpha + 1\), and let \(M^\alpha\) denote the \((w_\alpha \times 2t)\)-mesh with vertical paths \(P_{i_\alpha}, \ldots, P_{i_{\alpha+1}}\) and horizontal paths
    \(Q_1^\alpha, \ldots, Q_{2t}^\alpha\) where each \(Q_j^\alpha\) is the
    \(P_{i_\alpha}\)--\(P_{i_{\alpha}+1}\) subpath of \(Q_j\).
    Let \(j_1', \ldots, j_{2k}'\) be a sequence of integers from \([t]\) such that each \(j' \in [t]\) appears in the sequence exactly \(t-1\) times, and for any distinct \(j', j'' \in [t]\), there exists exactly
    one jump \(R_{\alpha}'\) such that it has endpoints on vertical paths \(i_{\alpha'}\) and \(i_{\alpha''}\) with \(\{j'_{\alpha'}, j'_{\alpha''}\} = \{j', j''\}\). Apply \cref{lem:ltorlinkage} to each \(M^\alpha\), and join the resulting paths to obtain disjoint paths \(S_1, \ldots, S_t\) such that
    for each \(\alpha \in [2k]\), the paths \(S_{t+1-j_\alpha'}\) contains
    \(Q_{t+1} \cap P_{i_\alpha}\), and the \(Q_t\)--\(Q_{t+1}\) subpath of \(P_{i_\alpha}\) is internally disjoint from the paths \(S_1, \ldots, S_t\). Hence, by our choice of the sequence \(j_1', \ldots, j_{2k}'\),
    we can extend the jumps \(R_1', \ldots, R_k'\) to paths
    \(R_1'', \ldots, R_k''\), each between a distinct pair of paths among \(S_1, \ldots, S_t\), so that the paths \(R_1'', \ldots, R_k''\) are pairwise disjoint and internally disjoint from \(S_1, \ldots, S_t\).
    Thus, \(G\) contains a model of \(K_t\) such that each branch set contains the vertex set of a distinct path \(S_\alpha\).
    This model is controlled by \(M\) since each \(S_\alpha\) intersects more than \(t\) vertical paths of \(M\).
\end{proof}

\section{Clique-minors in landscapes}\label{sec:lowerbound}

The purpose of this section is to establish bounds on the Euler-genus and the breadth of landscapes of large enough detail.
We do this in three steps.
To bound the total number of crosscaps and handles on the Dyck-wall of our surface we make use of a result of Gavoille and Hilaire \cite{GavoilleH2023MinorUniversal}.
To bound the breadth, that is the total number of vortices in the landscape we show directly how to construct the graph $H'_{\nicefrac{1}{2}(t-3)(t-4),2t}$ as defined in \cref{sec:forcing} from the extended surface grid.
With \cref{lem:Kt_from_crosses} this yields the second bound.
Finally, returning to surfaces, we show show that any large enough surface wall contains the Dyck-wall representing the same surface as a minor.
This allows us to extend the result of Gavoille and Hilaire \cite{GavoilleH2023MinorUniversal} to arbitrary surface-walls by paying a small polynomial factor for the order.

To fully implement these bounds, we will show in \cref{sec:localstructure} that, given a landscape of unbounded depth and applying the Society Classification Theorem from \cref{sec:societyclassification} to one of its vortices of unbounded depth results in a landscape with either higher Euler-genus, or with larger breadth.
Assuming that we can never find a $K_t$-minor, this implies that the process must terminate after a bounded number of applications of the Society Classification Theorem.

Let us start with quickly discussing the case of bounded Euler-genus.
Let $k\geq 1$ and $b$ be non-negative integers and $\Sigma$ be a surface.
Notice that every extended $k$-$\Sigma$-wall with vortices contains the elementary $k$-$\Sigma$-wall as a minor.
Hence, it suffices to show that that $\Sigma$-walls of large, in $n$, order contain every graph on $n$ vertices that embeds in $\Sigma$ as a minor.
This is directly implied by \cref{prop:dyckwall_minors}.
Moreover, for complete graphs this is a straight forward consequence of the naive upper bound that says that $K_t$ can be embedded in any surface of Euler-genus at least $t^2$ we obtain the following corollary.
It is important to note that the construction given by Gavoille and Hilaire in \cite{GavoilleH2023MinorUniversal} finds a model $\mu$ of the $H$-minor where, for each $v\in V(H)$, $\mu(v)$ contains the intersection of $C_k$ -- that is the inner most base cycle of the $k$-$\Sigma$-wall $D$ -- with some path from one of the handle and crosscap transactions.
Indeed, they show that the part of the graph embedded in the cylindrical mesh made from the base cycles of $D$ can be considered to be outerplanar.
It follows that the minor model found by Gavoille and Hilaire is indeed \textsl{controlled} by the cylindrical mesh of $D$.

\begin{corollary}\label{cor:embedKt}
Let $t\geq 1$ be a non-negative integer, $\Sigma$ be a surface, and let $c_{\ref{prop:dyckwall_minors}}$ be the constant from \cref{prop:dyckwall_minors}.
Then every graph $H$ is a minor of the $c_{\ref{prop:dyckwall_minors}}\cdot g^2(|V(H)|+g)^2$-$\Sigma$-wall $D$ and, if $H$ is a complete graph, there exists a minor model of $H$ controlled by $D$.

In particular, $K_t$ is a minor of the $c_{\ref{prop:dyckwall_minors}}\cdot g^2(t+g)^2$-$\Sigma$-wall.
\end{corollary}

So, for now, we may focus on showing that an extended surface-wall of large enough order, relative to $t$, with $\nicefrac{1}{2}(t-3)(t-4)+1$ vortex segments (each with a cross) contains a $K_t$-minor.

\begin{lemma}\label{lemma:cliqes_in_extended_Dyckwalls}
Let $t$ be a positive integer, $\Sigma$ be the sphere, and let $W$ be an extended $(2t+8)$-surface-wall with  $\nicefrac{1}{2}(t-3)(t-4)$ vortices.
Moreover, let $G$ be the union of $W$ and $(t-3)(t-4)$ vertex-disjoint $V(W)$-paths $L_i,R_i$, $i\in[\nicefrac{1}{2}(t-3)(t-4)]$ such that for every $i\in[\nicefrac{1}{2}(t-3)(t-4)]$, $L_i$ and $R_i$ form a cross on the innermost cycle of the nest of the $i$th vortex-segment of $W$ and the endpoints of $L_i$ and $R_i$ are endpoints of rails of the $i$th vortex-segment of $W$.
Then $G$ contains a $K_t$-minor controlled by $W$.

Moreover, a model of such a $K_t$-minor can be found in time $\mathbf{O}(|V(G)|)$.
\end{lemma}

\begin{proof}
Instead of constructing the $K_t$-minor directly, we show that $G$ contains $H'_{\nicefrac{1}{2}(t-3)(t-4)+1,2t}$ as a minor.
The assertion then follows from \cref{lem:Kt_from_crosses}.

Let $S_1,\dots,S_{\nicefrac{1}{2}(t-3)(t-4)+1}$ be the $\nicefrac{1}{2}(t-3)(t-4)+1$ segments such that $W$ is their cylindrical closure.
It follows that $S_1$ is a $(2t+8)$-wall segment while $S_2,\dots,S_{\nicefrac{1}{2}(t-3)(t-4)+1}$ are $(2t+8)$-vortex segments.

Let $C_1,\dots C_{2t+8}$ be the base cycles of $D$ in order such that $C_{2t+8}$ is the simple cycle of $D$.

Notice that the wall-segment $S_1$ has order at least two, so it has a ``left most'' vertical path, let us call this path $P_{\nicefrac{3}{2}(t-3)(t-4)+1}$, and a ``right most'' vertical path, say $P_1$.
Let $G'$ and $W'$ be obtained from $G$ and $W$ respectively by deleting all vertices of $S_1$ that do not belong to any of these two paths.

\begin{figure}[ht]
    \centering
    \begin{tikzpicture}

        \pgfdeclarelayer{background}
		\pgfdeclarelayer{foreground}
			
		\pgfsetlayers{background,main,foreground}

        \begin{pgfonlayer}{background}
            \pgftext{\includegraphics[width=15cm]{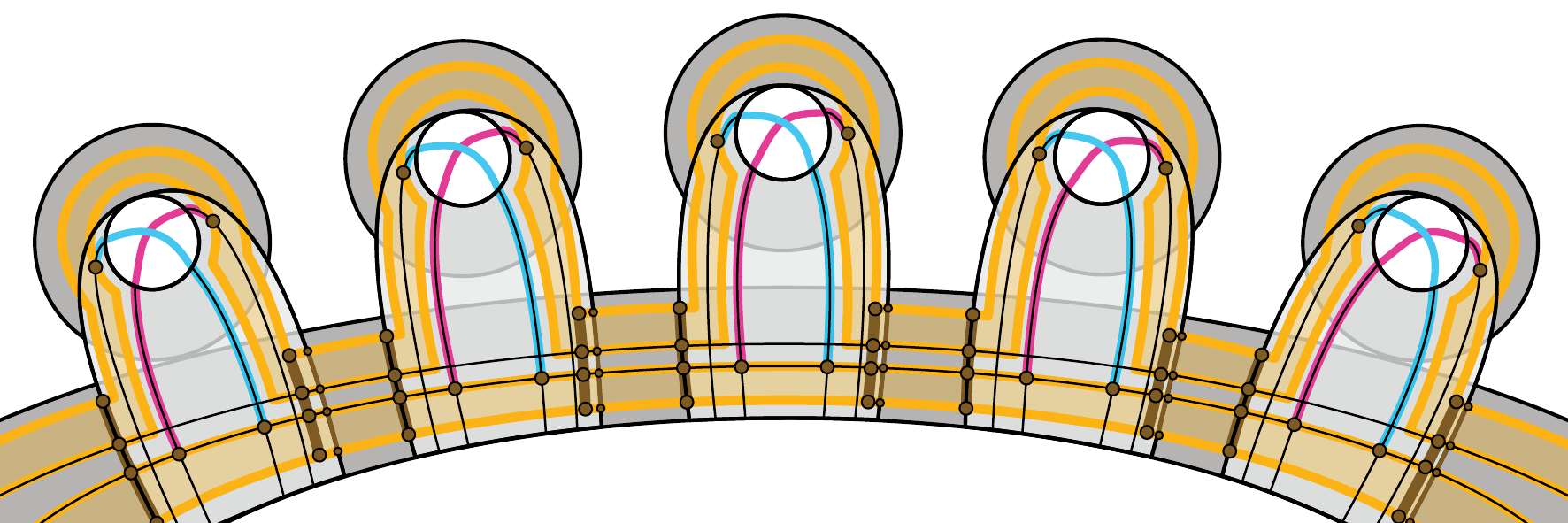}};
        \end{pgfonlayer}{background}
			
        \begin{pgfonlayer}{main}
        \node (C) [v:ghost] {};
   
        \end{pgfonlayer}{main}
        
        \begin{pgfonlayer}{foreground}
        \end{pgfonlayer}{foreground}

    \end{tikzpicture}
    \caption{A sketch of how to find the graph $H'_{w,2t}$ as a minor in an extended surface-wall by routing through its vortex segments.}
    \label{fig:clique_through_vortex}
\end{figure}

These two paths will be the first and last of the vertical paths of $H'_{\nicefrac{1}{2}(t-3)(t-4)+1,2t}$.
Before we explain how the remaining $P_i$ are found let us quickly discuss the $2t$ horizontal paths.
These paths, we will name them $Q_1,\dots,Q_{2t}$, are divided into two groups, namely $Q_1\dots Q_t$ which will be the \emph{top paths} and $Q_{t+1},\dots Q_{2t}$ which will be the \emph{bottom paths} of $H'_{\nicefrac{1}{2}(t-3)(t-4)+1,2t}$.
The paths $Q_{t+1},\dots,Q_{2t}$ are straightforward to find as these are simply the subpaths of $C_{t+9},\dots,C_{2t+8}$ that form a $V(P_1)$-$V(P_{\nicefrac{3}{2}(t-3)(t-4)+1})$-linkage in $W'$.

For the top paths we start at the last vertex of $C_i$ on $P_1$ -- when traversing along $C_i$ within $W'$ -- for each $i\in[t]$.
We then traverse along the cycles $C_1,\dots,C_t$ and pass through each of the vortex segments as follows.
Say we enter $S_j$ with $j\in[2,\nicefrac{1}{2}(t-3)(t-4)]$ while traversing along the cycles $C_1,\dots,C_t$.
For each of the $Q_i$, with $i\in[t]$, we follow along $C_i$ until we have met the $(i+3)$th path of the rails of $S_j$.
We then follow along these rails until, for each $i\in[t]$, we have met the outermost $i$th cycle of the nest of $S_j$.
From here there are two possible ways to travel ``around'' the vortex in the middle of the nest of $S_j$.
We choose the direction of traversal that lets us encounter the last (or ``right most'') rail before we \textsl{would} encounter the $(t+4)$th rail.
However, before that happens each $Q_i$ stops when it first meets the $(i+3)$th right most path of the rails of $S_j$.
From here our $Q_i$ can safely pass down back onto the cycles $C_1,\dots,C_t$ and leave $S_j$ to the right.
See \cref{fig:clique_through_vortex} for an illustration.

For the remaining vertical paths, we collect $3$ vertical paths from each vortex but the very last one, here we only collect two.
In total this will amount to
\begin{align*}
    2 + 3\cdot (\frac{1}{2}(t-3)(t-4)-1) +2 = \frac{3}{2}(t-3)(t-4)+1
\end{align*}
vertical paths as required.

Recall that we require a cross between the $i$th and $(i+1)$th vertical path for each $i\in[\nicefrac{1}{2}(t-3)(t-4)]$ with $i\equiv 2 \pmod 3$.
So for each $j\in[2,\nicefrac{1}{2}(t-3)(t-4)]$ we select $P_{3(j-1)+2}$ to be the $C_1$-$C_{2t+8}$-subpath of the first (or ``left most'') path from the rails.
We also select $P_{3(j-1)+3}$ and $P_{3(j-1)+4}$ to be the $C_1$-$C_{2t+8}$-subpaths of the second to last and last (or the two ``right most'') rails of $S_j$.
Finally, $P_{\nicefrac{3}{2}(t-3)(t-4)-1}$ and $P_{\nicefrac{3}{2}(t-3)(t-4)}$ are the $C_1$-$C_{2t+8}$-subpath of the first and last rail from $S_{\nicefrac{1}{2}(t-1)(t-3)}$.

Notice that each vortex segment has at least four cycles in its nest that separate the innermost cycle of the nest from those cycles selected for the construction of the $Q_i$.
Moreover, there are at least four rails that have not been used for our construction.
It follows that the two paths forming the cross at the centre of each vortex can be routed such that these crosses together with $P_i$ and $Q_j$ form a graph that must contain $H'_{\nicefrac{1}{2}(t-3)(t-4),2t}$ as a minor.
\end{proof}

An immediate consequence of \cref{lemma:cliqes_in_extended_Dyckwalls} is the following.
Here, the additional ``$+16$'' on the detail is used to guarantee that the crosses provided by the vortices can be routed into the rails of their respective vortex-segments by using \cref{lemma:reroute_crosses}.

\begin{corollary}\label{cor:cliquesinlandscapes}
Let $t\geq 4$, $b\geq \nicefrac{1}{2}(t-3)(t-4)$, and $d\geq 2t+24$ be positive integers, $\Sigma$ be a surface, and $G$ be a graph.
If there exists a $\Sigma$-landscape $\Lambda=(A,W,D,\delta)$ of detail $d$ and breadth $b$ for $G$, then $G$ contains a $K_t$-minor controlled by $D$.

Moreover, a model of such a $K_t$-minor can be found in time $\mathbf{O}(|V(G)|)$.
\end{corollary}

\subsection{Bounding the Euler-genus of surfaces with landscapes}\label{subsec:extendeddyckwall}

With \cref{cor:embedKt} we have a powerful tool to find a graph as a minor whenever we consider $\Sigma$-walls of moderately large order. 
However, repeated applications of \cref{thm:societyclassification} will not necessarily give us access to the correct surface-wall right away.

To give more context to why this might happen, we give a bit more context.
Gavoille and Hilaire do in fact prove that for every ``\emph{signature}'' of a surface $\Sigma$, a graph that can be embedded into $\Sigma$ with such a signature is indeed a minor of the corresponding surface-wall.
Here ``signature'' (roughly) refers to the order in which we the handle- and crosscap-segments are arranged to make up the surface wall.
They then make use of a result by Lazarus, Pocchiola, Vegter, and Verroust \cite{LazarusPVV2001Computing} stating that every graph embeddable on an \textsl{orientable} surface has an embedding that fits into their framework.
Notice that here the signature does not play a role.
For \textsl{non-orientable} surfaces, Fuladi, Hubard, and de Mesmay \cite{FuladiHD2024Short} proved that every graph embeddable on a non-orientable surface embeds, in a way compatible with the work of Gavoille and Hilaire, in that surface with a signature corresponding to only using crosscap-segments.
Fouladi \cite{Fuladi2023Embedded} extended the result from \cite{FuladiHD2024Short} in her PhD thesis to also allow for signatures that consist of handle- and crosscap-segments, but here the handles need to form a segment in the cylindrical concatenation.

In order to obtain a $\Sigma$-wall for some surface $\Sigma$ with $h$ handles and $c\leq 2$ crosscaps, we have to prove that any surface-wall with $h'$ handles and $'c\geq c$ crosscaps such that $2h'+c' = 2h+c$ is also minor-universal for the graphs that embed into $\Sigma$.
To establish this, and to ultimately bound the Euler-genus of our surface $\Sigma$ for which we find a landscape, we prove that, for any choice of non-negative integers $h,h',c,'$ as above, the $k$-$\Sigma$-wall is a minor of the $k'$-surface-wall with $h'$ handles and $c'$ where $k'\in\mathbf{poly}(k)$.
This means that the family $\{H^{(h,c)}_k\}_{k\in\mathbb{N}}$ is minor-universal for all graphs that embed into a surface that is homeomorphic to the surface that can be obtained by adding $h$ handles and $c$ crosscaps to the sphere.
Notice that this is a strengthening of the result of Gavoille and Hilaire \cite{GavoilleH2023MinorUniversal} with the downside that our universal graph must be slightly larger than theirs.

Of course, we cannot proceed without introducing a dedicated data structure for this situation.
The arguments below, in their purest form, have appeared in the work of Thilikos and Wiederrecht \cite{ThilikosW2024Excluding}.
However, their analysis yields bounds exponential in the Euler-genus.
The purpose of the data structure below is to allow for a finer analysis of the argument to ensure that all dependencies remain polynomial.

Let $(G,\Omega)$ be a society with a cylindrical rendition and a nest $\mathcal{C}$ around the vortex $c_0$.
Further let $\mathcal{P}_1,\dots,\mathcal{P}_{\ell}$ be a set of transactions on $(G,\Omega)$ as well as $\mathcal{R}$ be a radial linkage such that
\begin{itemize}
    \item $V(\mathcal{P}_i)\cap V(\mathcal{P}_j)=\emptyset$ for all $i\neq j\in[\ell]$ as well as $V(\mathcal{R})\cap V(\mathcal{P}_i)=\emptyset$ for all $i\in[\ell]$,
    \item for each $i\in[\ell]$, $\mathcal{P}_i$ is orthogonal to $\mathcal{C}$ and $\mathcal{R}$ is orthogonal to $\mathcal{C}$, and
    \item there exist pairwise disjoint segments $I_1,J_1,I_2,J_2,\dots,I_{\ell},J_{\ell},R$ of $\Omega$ such that these segments appear on $\Omega$ in the order listed, for each $i\in[\ell]$, $\mathcal{P}_i$ is a $V(I_i)$-$V(J_i)$-linkage, and each path from $\mathcal{R}$ has one endpoint in $V(R)$.
\end{itemize}
We call the tuple $(G,\Omega,\mathcal{C},\mathcal{R},\mathfrak{P})=\{ \mathcal{P}_1,\dots,\mathcal{P}_{\ell}\})$ a \emph{surface configuration} of $(G,\Omega)$ and the sequence $(I_1,J_1,I_2,J_2,\dots,I_{\ell},J_{\ell},R)$ the \emph{signature} of $(G,\Omega,\mathcal{C},\mathcal{R},\mathfrak{P})$.
Finally we say that $(s,r,p_1,\dots,p_{\ell})$ is the \emph{strength} of $(G,\Omega,\mathcal{C},\mathcal{R},\mathfrak{P})$ if $|\mathcal{C}|=s$, $|\mathcal{R}|=r$ and $|\mathcal{P}_i|=p_i$ for all $i\in[\ell]$.

\begin{observation}\label{obs:surface-configs-to-walls}
Let $k\geq 3$, $\ell$, $h$, and $c$ be non-negative integers with $\ell=h+c$.
Moreover, let $s\geq k$, $r\geq 4k$, $p_i\geq 4k$ for all $i\in[\ell]$.

Then, for every surface configuration $(G,\Omega,\mathcal{C},\mathcal{R},\mathfrak{P})$ of strength $(s,r,p_1,\dots,p_{\ell})$ with $h$ handle-transactions and $c$ crosscap-transactions, $G$ contains a $k$-surface-wall $W$ with $h$ handles and $k$ crosscaps as a subgraph such that the base cycles of $W$ are cycles from $\mathcal{C}$ and the vertical paths of the $k$-wall-segment of $W$ are subpaths of the paths from $\mathcal{R}$.
\end{observation}

\begin{observation}\label{obs:suface-walls-to-configs}
Let $h$, $c$, and $k\geq 4$ be non-negative integers.
Also, let $W$ be a $k$-surface wall with $h$ handles and $c$ crosscaps.
Finally, let $\Omega$ be a cyclic ordering of the vertices of the simple cycle of $W$ obtained by traversing along this cycle.
Then the society $(W,\Omega)$ has a surface configuration of strength $(k,4k,2k,\dots,2k)$ with $h$ handle and $c$ crosscaps transactions where the nest consists of the base cycles of $W$.
\end{observation}

\begin{lemma}\label{lemma:dissolve_crosscap}
Let $a_0,b_0,c_0$ be non-negative integers and $(G,\Omega)$ be a society.
If $(G,\Omega,\mathcal{C}=\{ C_1,\dots,C_s\},\mathcal{R},\mathfrak{P}=\{ \mathcal{P}_1,\dots,\mathcal{P}_{\ell}\})$ is a surface configuration of strength $(s,r,p_1,\dots,p_{\ell})$ such that $\mathcal{P}_i$, $i\in[2,\ell]$ ($i\in[\ell-1]$), is a crosscap transaction, $\mathcal{P}_{i-1}$ ($\mathcal{P}_{i+1}$) is a handle transaction, $s\geq 2a_0+2b_0+c_0$, $p_i\geq a_0+b_0+c_0$, and $p_{i-1} \text{($p_{i+1}$)} \geq 6a_0+6b_0+4c_0$.

Then there exist crosscap transactions $\mathcal{Q}_1,\mathcal{Q}_2,\mathcal{Q}_3$ such that $$(G,\Omega,\{ C_{2a_0+2b_0+c_0+1},\dots, C_s\},\mathcal{R},\{ \mathcal{P}'_1,\dots,\mathcal{P}_{\ell+1}'\})$$ with
\begin{itemize}
    \item $\mathcal{P}'_j=\mathcal{P}_j$ for all $j\in[1,i-2]$ ($j\in[1,i-1]$),
    \item $\mathcal{P}'_{i-1} = \mathcal{Q}_1$, $\mathcal{P}'_{i} = \mathcal{Q}_2$, and $\mathcal{P}'_{i+1} = \mathcal{Q}_3$ ($\mathcal{P}'_{i} = \mathcal{Q}_3$, $\mathcal{P}'_{i+1} = \mathcal{Q}_2$, and $\mathcal{P}'_{i+2} = \mathcal{Q}_1$), and
    \item $\mathcal{P}'_j=\mathcal{P}_{j+1}$ for all $j\in[i+1,\ell]$ ($j\in[i+2,\ell]$)
\end{itemize}
is a surface configuration of strength $(s-2a_0-2b_0-c_0,r,p_1,\dots,p_{i-2},a_0,b_0,c_0,p_{i+1},\dots,p_{\ell})$ ($(s-2a_0-2b_0-c_0,r,p_1,\dots,p_{i-1},c_0,b_0,a_0,p_{i+2},\dots,p_{\ell})$).
\end{lemma}

\begin{proof}
Notice that it suffices to prove the case where $\mathcal{P}_{i-1}$ is a handle transaction and $\mathcal{P}_i$ is a crosscap transaction.
The situation where $\mathcal{P}_{i+1}$ is the handle transaction that should be turned into two crosscaps follows along symmetric arguments.

Since $\mathcal{P}_{i-1}$ is a handle transaction, it consists of two planar transactions $\mathcal{L}_1$ and $\mathcal{L}_2$, each of order $\frac{1}{2}p_{i-1}\geq 3a_0+3b_0+2c_0$.
We will assume to be working on a linearisation $\lambda$ of $\Omega$ such that the segments $I_{i-1},J_{i-1},I_i,J_i$ as in the definition of surface configurations appear from left to right in this order.
For each $j\in[2]$ let $L_j^1$ and $L_j^2$ be the two end-segments of $\mathcal{L}_j$ such that $L_j^1$ appears before (or ``left'' of) $L_j^2$ in $\lambda$.
We also assume that the paths in $\mathcal{L}_1$, $\mathcal{L}_2$, $\mathcal{P}_i$ are numbered according the the appearance of their endpoints in $L_1^1$, $L_2^1$, and $I_i$ respectively.

The proof is dived into two steps.

First we move one of the two end-segments of $\mathcal{P}_i$ in a way that $L_2^2$ separates $I_j$ and $J_j$ on $\lambda$.
See \cref{fig:crosscapification1} for an illustration of this first step.

\paragraph{Step 1.}
We begin by constructing the planar transaction $\mathcal{A}_1$ of order $a_0$ as follows.
Simply define $\mathcal{A}_1\subseteq \mathcal{L}_1$ to be the left most $a_0$ paths of $\mathcal{L}_1$.
Similarly, we define $\mathcal{B}_1$ to be the planar transaction of order $a_0+b_0$ consisting of the first $a_0+b_0$ paths of $\mathcal{L}_2$ starting from the $a_0+b_0+c_0+1$-left most path.

\begin{figure}[ht]
    \centering
    \begin{tikzpicture}

        \pgfdeclarelayer{background}
		\pgfdeclarelayer{foreground}
			
		\pgfsetlayers{background,main,foreground}

        \begin{pgfonlayer}{background}
            \pgftext{\includegraphics[width=15cm]{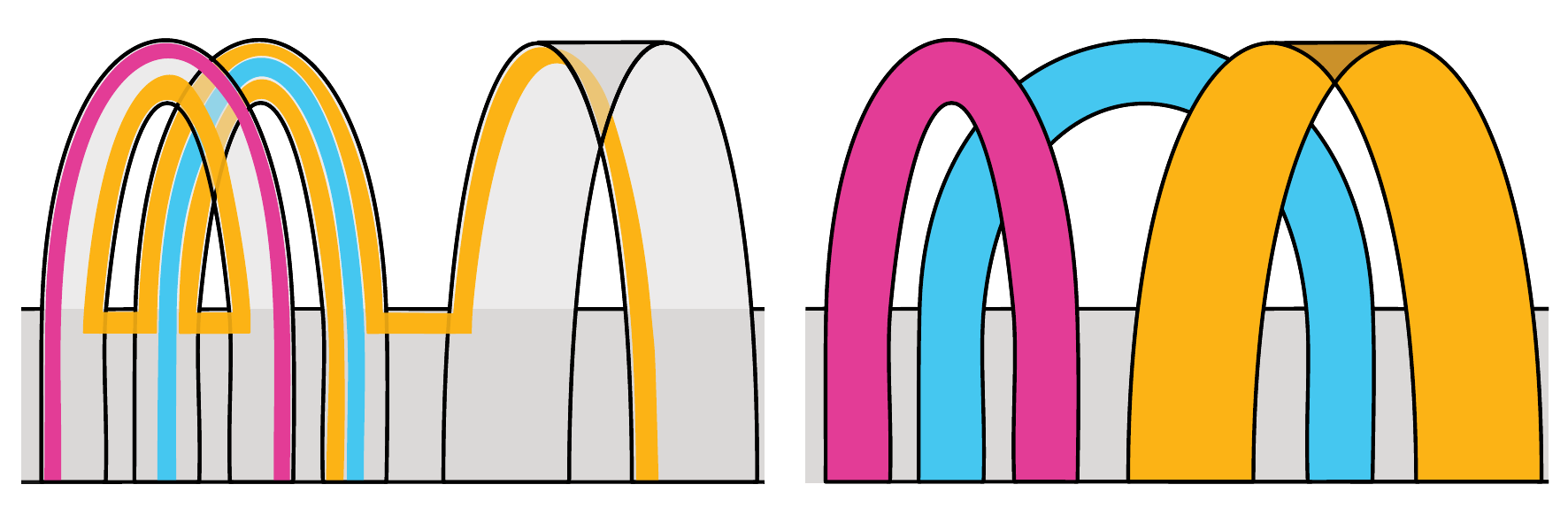}};
        \end{pgfonlayer}{background}
			
        \begin{pgfonlayer}{main}
        \node (C) [v:ghost] {};

        \node (x) [v:ghost,position=270:24.5mm from C] {};

        \node (a1) [v:ghost,position=180:70mm from x] {$\mathcal{A}_1$};
        \node (a2) [v:ghost,position=180:48mm from x] {$\mathcal{A}_1$};

        \node (b1) [v:ghost,position=180:58mm from x] {$\mathcal{B}_1$};
        \node (b2) [v:ghost,position=180:39.5mm from x] {$\mathcal{B}_1$};

        \node (c1) [v:ghost,position=180:43.5mm from x] {$\mathcal{C}_1$};
        \node (c2) [v:ghost,position=180:13mm from x] {$\overline{\mathcal{C}_1}$};

        \node (a1) [v:ghost,position=0:7mm from x] {$\mathcal{A}_1$};
        \node (a2) [v:ghost,position=0:26mm from x] {$\mathcal{A}_1$};

        \node (b1) [v:ghost,position=0:16mm from x] {$\mathcal{B}_1$};
        \node (b2) [v:ghost,position=0:54mm from x] {$\mathcal{B}_1$};

        \node (c1) [v:ghost,position=0:39mm from x] {$\mathcal{C}_1$};
        \node (c2) [v:ghost,position=0:67mm from x] {$\overline{\mathcal{C}_1}$};
        \end{pgfonlayer}{main}
        
        \begin{pgfonlayer}{foreground}
        \end{pgfonlayer}{foreground}

    \end{tikzpicture}
    \caption{The first step in the proof of \cref{lemma:dissolve_crosscap}, moving one end of the crosscap transaction in-between the handle transaction.
    On the left the routing of the three transactions $\mathcal{A}_1$, $\mathcal{B}_1$, $\mathcal{C}_1$ through the transactions $\mathcal{P}_{i-1}$ and $\mathcal{P}_i$ is depicted.
    The result is shown on the right.
    We write $\overline{\mathcal{X}}$ to indicate that the order of the endpoints of $\mathcal{X}$ is reversed between the two end-segments of the transaction $\mathcal{X}$.}
    \label{fig:crosscapification1}
\end{figure}
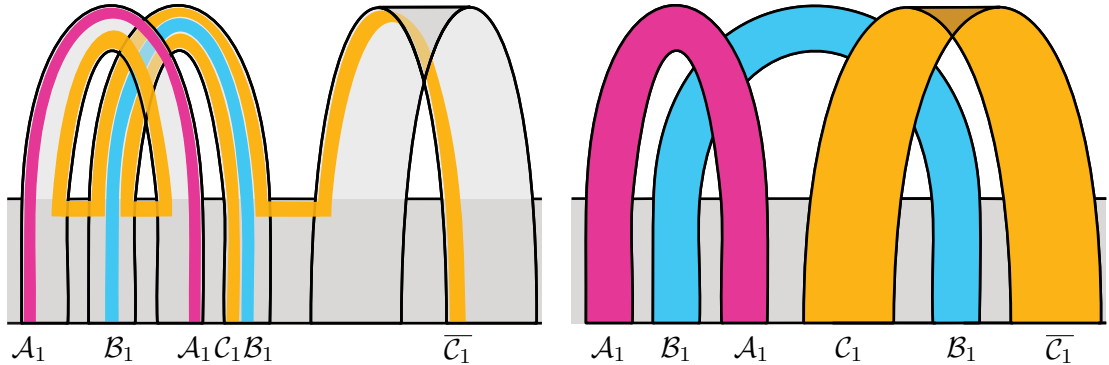

All that is left is the crosscap transaction $\mathcal{C}_1$ of order $a_0+b_0+c_0$.
In the definition of $\mathcal{B}_1$ we have reserved the first $a_0+b_0+c_0$ paths of $\mathcal{L}_2$.
Moreover, as $|\mathcal{L}_2|\geq 3a_0+3b_0+2c_0$, also the last $a_0+b_0+c_0$ paths of $\mathcal{L}_2$ are not used yet.
Hence, we may safely start from the left most $a_0+b_0+c_0$ endpoints from the right end-segment of $\mathcal{L}_2$, these belong to the $a_0+b_0+c_0$ right most paths of $\mathcal{L}_2$.
We follow along these paths until the $j$th path has met the $j$th cycle of $\mathcal{C}$ for the second time for each $j\in[a_0+b_0+c_0]$.
From here we follow along the cycles of $\mathcal{C}$ towards the right until the $j$th path of our routing has met the path $L_1^{\frac{1}{2}p_{i-1}-j}$, that is the $j$th path of $\mathcal{L}_1$ from the right, for each $j\in[a_0+b_0+c_0]$.
Notice that the interlinking between the cycles of the nest $\mathcal{C}$ and the paths from our transactions is effortlessly possible because the transactions are orthogonal to $\mathcal{C}$.
Now follow along $\mathcal{L}_1$ ``backwards'' until the $j$th paths of our linkage has met the $j$th cycle of $\mathcal{C}$ again for each $j\in[a_0+b_0+c_0]$.
Then traverse along the cycles of $\mathcal{C}$ from left to right until the $j$th path of $\mathcal{C}_1$ has met the $j$th path of $\mathcal{L}_2$.
As mentioned before, these paths were reserved when we defined $\mathcal{B}_1$ and thus they are free to use now.
Follows $\mathcal{L}_2$, this time in forward direction, until, again, the $j$th path of our linkage has met the $j$th cycle of $\mathcal{C}$ for the third time.
Follow along these cycles until the $j$th path has met the $j$th path of $\mathcal{P}_i$.
Then, finally, follow along $\mathcal{P}_i$ until the segment $J_i$ is reached.

Notice that this first step of the construction has produced three transactions, namely $\mathcal{A}_1$ of order $a_0$, $\mathcal{B}_1$ of order $a_0+b_0$, and $\mathcal{C}_1$, so far the only crosscap, of order $a_0+b_0+c_0$.
Moreover, all of these transactions are orthogonal to the cycles $\{ C_{a_0+b_0+c_0+1},\dots C_s\}$
\medskip

In the second step we will now ``unravel'' the newly created transactions $\mathcal{A}_1$, $\mathcal{B}_2$, and $\mathcal{C}_1$ into three individual crosscaps that occur in $\Omega$ in order.
See \cref{fig:crosscapification2} for an illustration.

\paragraph{Step 2.}
Let us select $\mathcal{A}_0$ to be a transaction of order $a_0$ in $(G,\Omega)$ as follows.
We begin with the $a_0$ endpoints of $\mathcal{A}_1$.
We then follows along $\mathcal{A}_1$ until the $j$th path of the linkage we are constructing meets $C_{j+2a_0+b_0+c_0}$, $j\in[a_0]$ for the second time.
Our routing then follows along the first $a_0$ cycles of $\mathcal{C}$ after $C_{2a_0+b_0+c_0}$, until the $j$th path of $\mathcal{A}_1$ meets the $j$th path of $\mathcal{C}_1$ for the first time.
Now, follow along $\mathcal{C}_1$ until the $j$th path has met, yet again, the $j$th cycle of $\mathcal{C}$, after $C_{2a_0+b_0+c_0}$, for each $j\in[a_0]$.
We then route our paths to the left, until the $j$th path has met the $j$th path, for each $j\in[a_0]$, of $\mathcal{B}_1$.
Finally, we may follow along $\mathcal{B}_1$ ``backwards'' until we meet its left end-segment on $\Omega$.

Notice that $\mathcal{A}_0$, as constructed above, has passed through the crosscap $\mathcal{C}_1$ exactly once.
Thus, $\mathcal{A}_0$ must be a crosscap transaction itself.

\begin{figure}[ht]
    \centering
    \begin{tikzpicture}

        \pgfdeclarelayer{background}
		\pgfdeclarelayer{foreground}
			
		\pgfsetlayers{background,main,foreground}

        \begin{pgfonlayer}{background}
            \pgftext{\includegraphics[width=15cm]{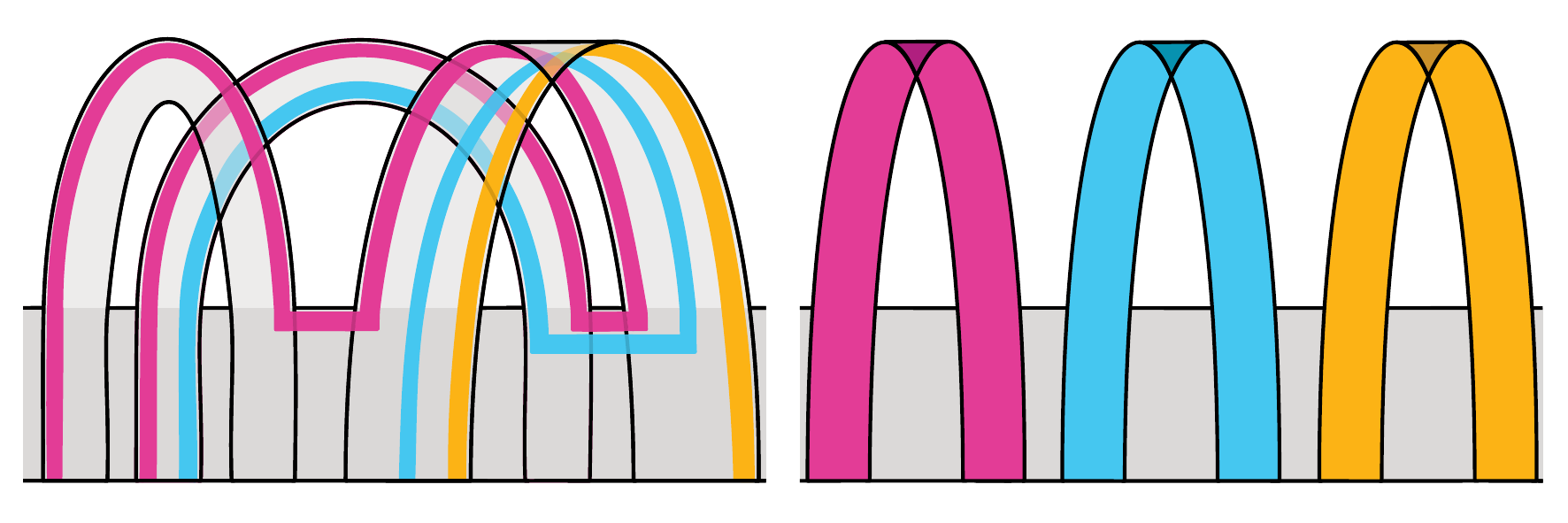}};
        \end{pgfonlayer}{background}
			
        \begin{pgfonlayer}{main}
        \node (C) [v:ghost] {};

        \node (a1) [v:ghost,position=180:70mm from x] {$\mathcal{A}_0$};
        \node (a2) [v:ghost,position=180:61mm from x] {$\overline{\mathcal{A}_0}$};

        \node (b1) [v:ghost,position=180:55.5mm from x] {$\mathcal{B}_0$};
        \node (b2) [v:ghost,position=180:36mm from x] {$\overline{\mathcal{B}_0}$};

        \node (c1) [v:ghost,position=180:30mm from x] {$\mathcal{C}_0$};
        \node (c2) [v:ghost,position=180:3mm from x] {$\overline{\mathcal{C}_0}$};

        \node (a1) [v:ghost,position=0:6mm from x] {$\mathcal{A}_0$};
        \node (a2) [v:ghost,position=0:20mm from x] {$\overline{\mathcal{A}_0}$};

        \node (b1) [v:ghost,position=0:30mm from x] {$\mathcal{B}_0$};
        \node (b2) [v:ghost,position=0:44.5mm from x] {$\overline{\mathcal{B}_0}$};

        \node (c1) [v:ghost,position=0:55mm from x] {$\mathcal{C}_0$};
        \node (c2) [v:ghost,position=0:69mm from x] {$\overline{\mathcal{C}_0}$};

        \end{pgfonlayer}{main}
        
        \begin{pgfonlayer}{foreground}
        \end{pgfonlayer}{foreground}

    \end{tikzpicture}
    \caption{The second step in the proof of \cref{lemma:dissolve_crosscap}, untangling the handle and the crosscap transaction into three individual crosscap transactions.
    On the left the routing of the three transactions $\mathcal{A}_0$, $\mathcal{B}_0$, $\mathcal{C}_0$ through the transactions $\mathcal{A}_1$, $\mathcal{B}_1$, and $\mathcal{C}_1$ is depicted.
    The resulting three crosscap transactions are shown on the right.}
    \label{fig:crosscapification2}
\end{figure}
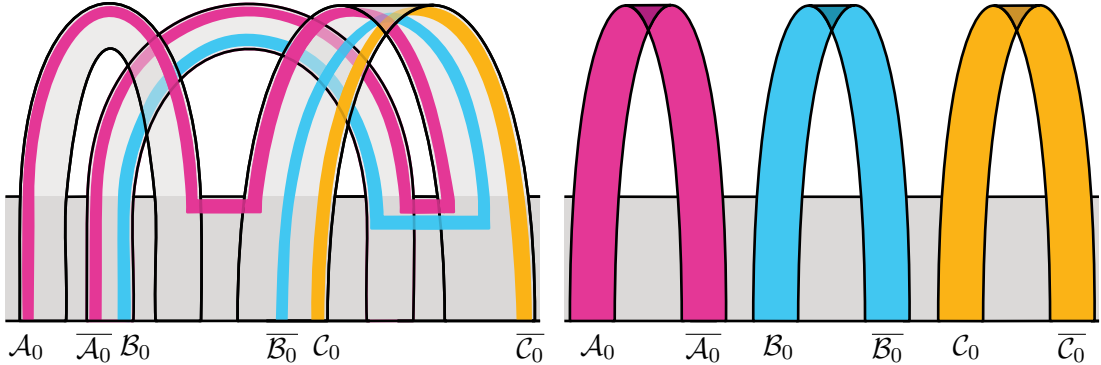

To construct $\mathcal{B}_0$ we start with the right most $b_0$ endpoints from the left end-segment of $\mathcal{B}_1$.
Notice that, with $|\mathcal{B}_1|=a_0+b_0$, none of the paths involved intersect $\mathcal{A}_0$.
We follow along the right most $b_0$ paths of $\mathcal{B}_1$ until the $j$th path, $j\in[b_0]$, has met the cycle $C_{2a_0+b_0+c_0+j}$ for the second time.
Notice that $\mathcal{A}_1$ makes use of a similar section of the nest, but only uses the cycles $C_{a_0+b_0+c_0+1}$ to $C_{2a_0+b_0+c_0}$.
From here we follow the cycles towards the right until, for the last time, the $j$th path, $j\in[b_0]$, has met the $a_0+j$th path of $\mathcal{C}_1$ for the first time.
Since $\mathcal{A}_0$ intersects $\mathcal{C}_1$ exactly in its first $a_0$ paths and $|\mathcal{C}_1|=a_0+b_0+c_0$, there is still enough space to now follow along $\mathcal{C}_1$ ``backwards'' until we have reached the left end-segment of $\mathcal{C}_1$.

As before for $\mathcal{A}_0$, $\mathcal{B}_0$ has passed through $\mathcal{C}_1$ exactly one and is therefore a crosscap transaction itself.

To complete our construction, we define $\mathcal{C}_0$ to be the crosscap transaction consisting of the right most $c_0$ paths of $\mathcal{C}_1$.
\medskip

What is left is to observe that, by construction, the transactions $\mathcal{A}_0$, $\mathcal{B}_0$, and $\mathcal{C}_0$ are orthogonal to the cycles $C_{2a_0+2b_0+c_0+1},\dots,C_s$, they are pairwise disjoint and disjoint from all $\mathcal{P}_j$, $j\in[1,i-2]\cup[i+1,\ell]$.
Moreover, their end segments are contained in $I_{i-1}\cup J_{i-1}\cup I_i\cup J_i$.
Thus, this proves the assertion.
\end{proof}

\begin{lemma}\label{lemma:transforming-surface-configs}
Let $h$, $c\geq 1$, $s$, $s_0$, $r$, and $k$ be non-negative integers and $K=(G,\Omega,\mathcal{C},\mathcal{R},\mathfrak{P})$ be a surface configuration with $h$ handles and $c$ crosscaps.

If the strength of $K$ is $(s,r,p_1,\dots,p_{h+c})$ such that $s\geq (6h+7)hk+s_0$ and for each $i\in[h+c]$ 
\begin{itemize}
    \item $p_i\geq (20h+28)k$ if $\mathcal{P}_i$ is a handle transaction,
    \item $p_i\geq k$ and for all $j\in\{ i-1,i+1\}\cap[h+c]$, $\mathcal{P}_{j}$ is a crosscap transaction,
    \item $p_i\geq (4h+5)k$ if $i\in[2,h+c-1]$, $\mathcal{P}_{i}$ is a crosscap transaction, and $\mathcal{P}_{i-1}$ and $\mathcal{P}_{i+1}$ are handle transactions, and
    \item $p_i\geq (2h+3)k$ $\mathcal{P}_{i}$ is a crosscap transaction, and $\mathcal{P}_{j}$ is a handle transaction for exactly one $j\in\{ i-1,i+1\}$,
\end{itemize}
then there exists crosscap transactions $\mathfrak{Q}=\{\mathcal{Q}_1,\dots,\mathcal{Q}_{2h+c}\}$ and a nest $\mathcal{N}=\{ C_{s-s_0+1},\dots,C_{s_0}\}\subseteq \mathcal{C}$ such that $(G,\Omega, \mathcal{N},\mathcal{R},\mathfrak{Q})$ is a surface configuration of strength $(s_0,r,k,\dots,k)$.
\end{lemma}

\begin{proof}
We prove the claim by induction on the number $h$ of handles.
In case $h=0$ there is nothing to show, so we may assume that $h\geq 1$ there is at least one handle transaction in $\mathfrak{P}$.

Let $i\in[h+c]$ be chosen such that $\mathcal{P}_i$ is a crosscap transaction and for some $j\in\{ i-1,i+1\}$, $\mathcal{P}_j$ is a handle transaction.
In total, there are three cases to consider, namely the cases where \textsl{(i)} only $\mathcal{P}_{i-1}$ is a handle transaction, \textsl{(ii)} only $\mathcal{P}_{i+1}$ is a handle transaction, and \textsl{(iii)} both $\mathcal{P}_{i-1}$ and $\mathcal{P}_{i+1}$ are handle transactions.

Cases \textsl{(i)} and \textsl{(ii)} are symmetric and can be treated in one.
Moreover, we will see in our treatment of \textsl{(i)} that case \textsl{(iii)} can be reduced to case \textsl{(i)} by requiring that the right most crosscap transaction produced by our argument has order $(2h+3k)$ instead of order $k$.
As we will see, the argument works, with minor adjustments, in both cases.
Hence, in the following we will assume the following
\begin{itemize}
    \item $i\in[2,h+c]$ and $\mathcal{P}_{i-1}$ is a handle transaction of order at least $(20h+24)k$, and
    \item $|\mathcal{P}_i|\geq x\cdot k$ where $x\in\{4h+5,2h+3\}$.
\end{itemize}
And deduce both cases \textsl{(i)} and \textsl{(iii)} later.
\smallskip

Now let
\begin{align*}
    a_0 &\coloneqq  (2h+1)k,\\
    b_0 &\coloneqq  k\text{, and}\\
    c_0 &\coloneqq  (x-2h-2)k.
\end{align*}
Observe that
\begin{align*}
    p_{i-1} \geq (20h+24)k &\geq 6(2h+1)k + 6k + 4(2h+3)k\\
    & \geq 6(2h+1)k + 6k + 4(x-2h-2)k\\
    & = 6a_0 + 6b_0 + 4c_0.
\end{align*}
By \cref{lemma:dissolve_crosscap} there exist crosscap transactions $\mathcal{Q}_1$ of order $a_0$, $\mathcal{Q}_2$ of order $b_0$, and $\mathcal{Q}_3$ of order $c_0$ such that removing the innermost
\begin{align*}
    2a_0+2b_0+c_0 &= (4h+2)k + 2k + (x-2h-2)k\\
    &\geq (6h+7)k
\end{align*}
cycles from the nest and replacing $\mathcal{P}_{i-1}$ and $\mathcal{P}_i$ with $\mathcal{Q}_1$, $\mathcal{Q}_2$, and $\mathcal{Q}_3$ in this order yields a surface configuration $K'$ of strength (at least)
\begin{align*}
(s-(6h+7)k,r,q_1,\dots,q_{h+c+1}).
\end{align*}
Here $q_j=p_j$ for all $j\in[i-2]$, $q_j=p_{j-1}$ for all $j\in[i+2,h+c+1]$, $q_{i-1}=a_0$, $q_i=b_0$, and $q_{i+1}=c_0$.

Notice that the number of handle transactions of $K'$ is exactly $h-1$.
Moreover, $\mathcal{Q}_2$ borders on both sides to crosscap transactions, $\mathcal{Q}_1$ borders on at least one side to a crosscap transaction, and $\mathcal{Q}_3$ borders to a crosscap transaction either side if and only if $\mathcal{P}_i$ did and otherwise it boarders to exactly one handle transaction.
By our choice of $x$ and the facts that $a_0=(2h+1)k=(2(h-1)+3)k$ and $s-(6h+7)k \geq (6(h-1)+7)(h-1)k+s_0$, our claim now follows by applying the induction hypothesis to $K'$.
\end{proof}

Notice that the routing through the nest in the proof of \cref{lemma:dissolve_crosscap} occurs only locally.
This means that, instead of throwing away $(6h+7)k$ cycles in each step in the proof of \cref{lemma:transforming-surface-configs} we could use the same $(6h+7)k$ cycles in each step.
This would reduce the lower bound on $s$ in the statement of \cref{lemma:transforming-surface-configs} from $(6h+7)hk$ to $(6h+7)k$.

By combining \cref{cor:embedKt} and \cref{obs:suface-walls-to-configs,obs:surface-configs-to-walls} with \cref{lemma:transforming-surface-configs} we obtain the following.

\begin{corollary}\label{cor:universal-surface-walls}
There exists a universal constant $c_{\ref{cor:universal-surface-walls}}$ such that for all non-negative integers $c,h,t$, every graph $H$ that embeds in a surface homeomorphic to the surface obtained from the sphere by adding $h$ handles and $c$ crosscaps is a minor of the $\big( c_{\ref{cor:universal-surface-walls}}g^4(|V(H)|+g)^2\big)$-surface wall $W$ with $h$ handles and $c$ crosscaps, where $g=2h+c$.
Moreover, if $H$ is a complete graph, then there exists an $H$-minor model controlled by $W$.

In particular, $K_t$ is a minor for every $k$-surface-wall with $h$ handles and $c$ crosscaps where $2h+c=t^2$ and $k\in \mathbf{\Omega}(t^{12})$.
\end{corollary}

\section{Flattening a mesh}\label{sec:flattening}

\begin{definition}[Flat mesh]
    Let $n \ge 2$ be an integer.
    Let $G$ be a graph, let $M \subseteq G$ be an \((n \times n)\)-mesh.
    We say that $M$ is \emph{flat} in $G$ if there exists a \(\Sigma\)-decomposition \(\delta\) of \(G\) in a sphere \(\Sigma\) with a single vortex $c_0$, such that \(M\) is grounded in \(\delta\), and the trace of the perimeter of \(M\) bounds a closed disk \(\Delta \subseteq \Sigma\) which is disjoint from $c_0$, such that every vertex \(x \in V(M)\) that is a node in \(\delta\) is embedded in \(\Delta\).
\end{definition}

In this section we show the following.

\begin{theorem}[Flat Mesh Theorem]\label{thm:flatmesh}
    Let \(t, n'\) be integers with \(t \ge 5\), \(n' \ge 2\), let \(n = 100t^3(n'+2t+2)\), and let \(G\) be a graph with an \((n \times n)\)-mesh \(M\). Then either
    \begin{itemize}
        \item there exists a model \(\mu\) of \(K_t\) in \(G\)  which is controlled by \(M\), or
        \item there exist a set \(Z\) with \(|Z| < 16t^3\) and an \((n' \times n')\)-submesh \(M'\) of \(M\) which is disjoint from \(Z\) and flat in \(G - Z\).
    \end{itemize}
    Furthermore, there exists a \(\mathbf{O}(|E(G)|)\)-time algorithm which finds either the model \(\mu\) or the set \(Z\), the submesh \(M'\), and a \(\Sigma\)-decomposition witnessing that \(M'\) is flat in \(G - Z\).
\end{theorem}
Note that the \(\mathbf{O}(|E(G)|)\) running time above does not depend on \(t\).

The main tool used in the proof will be a generalisation of Gallai's Theorem on \(A\)-paths. 

\subsection{Generalising Gallai's theorem on \texorpdfstring{$A$}{A}-paths}\label{subsec:Apaths}

Gallai~\cite{Gallai1964MaximumMinimum} showed that for every graph \(G\) with a subset \(A \subseteq V(G)\) and every integer \(q \ge 1\),
either there are \(q\) pairwise vertex-disjoint \(A\)-paths in \(G\), or there exists a set \(Z \subseteq V(G)\) with \(|Z| \in \mathbf{O}(q)\) such that there is no \(A\)-path in \(G - Z\).
In this subsection we prove a generalisation of this theorem. The proof of this generalisation is identical to the proof of Gallai's theorem by Bruhn, Heinlein and Joos~\cite{BruhnHJ2018Frames}, except that the set \(A = \{x_1, \ldots, x_w\}\) is replaced by a family of sets \(\mathcal{A} = \{X_1, \ldots, X_w\}\).

For a family \(\mathcal{A}\) of sets of vertices, an \emph{\(\mathcal{A}\)-path} is a \((\bigcup\mathcal{A})\)-path \(R\) such that no set \(X \in \mathcal{A}\) contains both endpoints of \(R\). An \emph{\(\mathcal{A}\)-linkage} is a set \(\mathcal{R}\) of pairwise vertex-disjoint \(\mathcal{A}\)-paths such that no \(X \in \mathcal{A}\) contains an endpoint of more than one path from \(\mathcal{R}\). We denote by \(\|\mathcal{A}\|\) the greatest integer \(c \ge 0\) such that some \(x \in V(G)\)
belongs to \(c\) sets from \(\mathcal{A}\). We will mainly consider families \(\mathcal{A}\) with \(\|\mathcal{A}\| = 2\).
In this subsection we prove the following.

\begin{theorem}[\cref{GallaiGeneralization_simplified} reformulated]\label{GallaiGeneralization}
    Let \(G\) be a graph, let \(\mathcal{A}\) be a family of subsets of \(V(G)\), and let \(q \ge 0\) be an integer. Then
    \begin{enumerate}
        \item there exists an \(\mathcal{A}\)-linkage \(\mathcal{R}\) with \(|\mathcal{R}| = q\), or
        \item there exists a set \(Z \subseteq V(G)\) with \(|Z|< 4q\) such that every \(\mathcal{A}\)-path in \(G - Z\) has an endpoint in a set \(X \in \mathcal{A}\) with \(X \cap Z \neq \emptyset\).
    \end{enumerate}
    Moreover, there exists an algorithm with running time \(\mathbf{O}(\|\mathcal{A}\|\cdot|E(G)|)\) which finds either \(\mathcal{R}\) or \(Z\) as above.
\end{theorem}

We use the following simple observation from~\cite{BruhnHJ2018Frames}.

\begin{proposition}[Bruhn, Heinlein, and Joos \cite{BruhnHJ2018Frames}]\label{LeafToLeafPaths}
    A  subcubic tree with \(p \ge 2\) leaves contains at least \((p-1)/2\) pairwise vertex-disjoint leaf-to-leaf paths.
\end{proposition}

This enables us to give a fairly simple proof of \cref{GallaiGeneralization}.

\begin{proof}[Proof of \cref{GallaiGeneralization}]
    Let \(F \subseteq G\) be an inclusion-maximal subcubic forest without isolated vertices such that the set of leaves of \(F\) is \(V(F) \cap (\bigcup\mathcal{A})\), and no \(X \in \mathcal{A}\) contains two leaves of \(F\).
    We assume that the forest \(F\) is not empty because otherwise the theorem is satisfied by \(Z = \emptyset\).
    Let \(p\) be the number of leaves in \(F\), let \(T_1, \ldots, T_m\) be the components of \(F\), and for each \(i \in \{1, \ldots, m\}\), let \(p_i\) be the number of leaves in \(T_i\). By \cref{LeafToLeafPaths}, each \(T_i\) contains a set \(\mathcal{R}_i\) of at least \((p_i-1)/2\) pairwise vertex-disjoint leaf-to-leaf paths. Since \(p_1 + \cdots + p_m = p\), the union \(\mathcal{R} = \mathcal{R}_1 \cup \cdots \cup \mathcal{R}_m\) contains at least \((p-m)/2\) paths, and \(\mathcal{R}\) is an \(\mathcal{A}\)-linkage by our choice of \(F\).
    If \((p-m)/2 \ge q\), then the lemma is satisfied by a \(q\)-element subset of \(\mathcal{R}\). Hence we assume that \((p-m)/2 < q\).

    Let \(Z_1\) denote the set of leaves in \(F\), let \(Z_3\) denote the set of degree-3 vertices in \(F\), and let \(Z = Z_1 \cup Z_3\). We have \(|Z_1| = p\) and \(|Z_3| = \sum_{i=1}^m |Z_3 \cap V(T_i)| = \sum_{i=1}^m (p_i-2) = p - 2m\), so \(|Z| = |Z_1|+|Z_3| \le 2p-2m =4(p-m)/2< 4q\).
    
    Let \(R\) be an \(\mathcal{A}\)-path in \(G - Z\).
    By maximality of \(F\), the path \(R\) intersects \(F\) in at least one vertex.
    Since \(R\) is disjoint from \(Z\), every vertex in \(V(R\cap F)\) has degree two in \(F\).
    Let \(R'\) be an initial segment of \(R\) which is a \(\bigcup\mathcal{A}\)--\(V(F)\) path. By maximality of \(F\), the forest \(F \cup R'\) has two leaves in one set \(X \in \mathcal{A}\).
    One of these leaves is the common endpoint of \(R\) and \(R'\), and the other is a leaf of \(F\) witnessing that \(X \cap Z \neq \emptyset\).
    Thus \(Z\) has the desired properties and we are done.

    The only non-trivial step in turning the above proof into a linear-time algorithm is finding \(F\). We construct the forest greedily with a modified DFS algorithm. The algorithm will keep track of the set \(\mathcal{A}_{Z_1} \coloneqq  \{X \in \mathcal{A}: X \cap Z_1 \neq \emptyset\}\), where \(Z_1\) is the current leaf set of \(F\). 
    This costs us \(\mathbf{O}(\|A\|)\) of running time each time we add a new leaf to \(F\), and allows us to test in \(\mathbf{O}(\|A\|)\)-time whether a vertex of \(G\) belongs to \(\bigcup\mathcal{A}_{Z_1}\).
    The algorithm finding the forest \(F\) exhaustively calls a modified recursive DFS procedure on yet unvisited starting vertices \(z_0 \in\bigcup\mathcal{A}\setminus\bigcup\mathcal{A}_{Z_1}\).
    Each time the algorithm makes a recursive call to (``visits'') a vertex \(x \in V(G)\), the algorithm tests if \(x \in \bigcup\mathcal{A}_{Z_1}\) and proceeds as follows:
    \begin{itemize}
        \item if \(x \in \bigcup\mathcal{A}_{Z_1}\), then the recursive call terminates immediately,
        \item if \(x \in \bigcup\mathcal{A}\setminus\bigcup\mathcal{A}_{Z_1}\), then the forest \(F\) is extended by adding the path with which \(x\) was reached from the starting vertex \(z_0\), and the set \(\mathcal{A}_{Z_1}\) is updated, and
        \item if \(x \in V(G) \setminus \bigcup\mathcal{A}\), the algorithm makes recursive calls to the yet unvisited neighbours, with an early-termination after the second recursive subcall which results in extending \(F\) (if \(x\) is the starting vertex \(z_0\), the early-termination is after the first such recursive call).
    \end{itemize}
    This algorithm will find a forest \(F\) as in the proof with the claimed running time.
\end{proof}

\subsection{One dimensional jumps}

A flat mesh \(M\) in a graph \(G\) has the property that every \(M\)-bridge in \(G\) which does not have all its attachments on the perimeter of \(M\), has attachments on at most two consecutive vertical paths of \(M\).
In this subsection, we develop tools allowing us to have control over the jumps between non-consecutive vertical paths of a mesh which does not control a model of \(K_t\).

Let \(\mathcal{G} = (G_1, \ldots, G_w)\) be a sequence of pairwise disjoint connected subgraphs of a graph \(G\). We write \(V({\mathcal{G}}) \coloneqq  V(G_1 \cup \cdots \cup G_w)\). Two graphs \(G_i\) and \(G_j\) are \emph{\(\mathcal{G}\)-independent} if \(|i-j| > 1\).
A \emph{\(\mathcal{G}\)-jump} is a \(V(\mathcal{G})\)-path between \(\mathcal{G}\)-independent graphs \(G_i\) and \(G_j\).
Let \(\mathcal{A} = \{X_1, \ldots, X_{w-1}\}\) where \(X_i = V(G_i \cup G_{i+1})\). Then \(\|\mathcal{A}\|=2\), and \(\mathcal{G}\)-jumps are exactly \(\mathcal{A}\)-paths.
For each \(i \in [w]\) there are at most three indices \(j \in [w]\) such that \(G_i\) is not \(\mathcal{G}\)-independent with \(G_j\).
Therefore, from \cref{GallaiGeneralization} applied to \(G\) and \(\mathcal{A}\), we deduce the following.
\begin{corollary}\label{CorollaryGG}
    Let \(q \ge 1\) be an integer, let \(G\) be a graph, and let \(\mathcal{G} = (G_1, \ldots, G_w)\) be a sequence of pairwise disjoint connected subgraphs of \(G\). Then, either
    \begin{enumerate}
        \item there exist \(q\) pairwise disjoint \({\mathcal{G}}\)-jumps \(R_1, \ldots, R_q\) in \(G\) with endpoints in \(2q\) pairwise \(\mathcal{G}\)-independent graphs from \(\mathcal{G}\), or
        \item there exists a set \(Z \subseteq V(G)\) with \(|Z| < 4q\) and a subset \(\mathcal{Y} \subseteq \{G_1, \ldots, G_w\}\) with \(|\mathcal{Y}| < 12q\) and \(Z \cap V(\mathcal{G}) \subseteq V(\bigcup\mathcal{Y})\), such that every \(\mathcal{G}\)-jump in \(G - Z\) has an endpoint in \(\bigcup\mathcal{Y}\).
    \end{enumerate}
    Moreover there exists a \(\textbf{O}(|E(G)|)\)-time algorithm which finds either \((R_1, \ldots, R_q)\) or \((Z, \mathcal{Y})\) as above.
\end{corollary}

We derive the following variant of this corollary for the setting where the graphs from \(\mathcal{G}\) have bounded degree. Its strength lies in the fact that every \(\mathcal{G}\)-jump in \(G - Z\) has not one, but two endpoints in \(\bigcup\mathcal{Y}\).
\begin{lemma}\label{PPjumps}
    Let \(d, q \ge 1\) be integers, let \(G\) be a graph, and let \(\mathcal{G} = (G_1, \ldots, G_w)\) be a sequence of pairwise disjoint connected subgraphs of \(G\) each with maximum degree at most \(d\). Then either
    \begin{itemize}
        \item there exists a subsequence \(\mathcal{G}' = (G_1', \ldots, G_{w'}')\) of \(\mathcal{G}\) and \(q\) pairwise vertex-disjoint \({\mathcal{G}'}\)-jumps \(R_1, \ldots, R_q\) in \(G\)  with endpoints in \(2q\) pairwise \(\mathcal{G}'\)-independent graphs from \(\mathcal{G}'\), or
        \item there exist a set \(Z \subseteq V(G)\) with \(|Z| < 8q\) and a subset \(\mathcal{Y} \subseteq \{G_1, \ldots, G_w\}\) with \(|\mathcal{Y}| < 16(d+2)q\) and \(Z \cap V(\mathcal{G}) \subseteq V(\bigcup\mathcal{Y})\) such that every \(\mathcal{G}\)-jump in \(G - Z\) has both endpoints in \(\bigcup\mathcal{Y}\).
    \end{itemize}
    Moreover there exists a \(\textbf{O}(|E(G)|)\)-time algorithm which finds \((\mathcal{G'}, R_1, \ldots, R_q)\) or \((Z, \mathcal{Y})\) as above.
\end{lemma}
\begin{proof}
    Apply \cref{CorollaryGG} to \(q\), \(G\), and \(\mathcal{G}\).
    Lest the lemma be satisfied with \(\mathcal{G}' = \mathcal{G}\),
    we assume that there are a set \(Z_1 \subseteq V(G)\) with
  \(|Z_1| < 4q\) and a set \(\mathcal{Y}_1 \subseteq \{G_1, \ldots, G_w\}\)
    with \(|\mathcal{Y}_1| < 12q\) and \(Z_1 \cap V(\mathcal{G}) \subseteq V(\bigcup\mathcal{Y}_1)\) such that every \(\mathcal{G}\)-jump in \(G - Z_1\) has an endpoint in \(\bigcup\mathcal{Y}_1\).
    Let \(\mathcal{G}_1\) denote the subsequence of \(\mathcal{G}\) obtained by deleting the graphs belonging to \(\mathcal{Y}_1\).

    Next, apply \cref{CorollaryGG} to \(q\), \(G - Z_1\), and \(\mathcal{G}_1\).
    Lest the lemma be satisfied with \(\mathcal{G}' = \mathcal{G}_1\),
    we assume that there are a set \(Z_2 \subseteq V(G - Z_1)\) with \(|Z_2| < 4q\) and a set \(\mathcal{Y}_2 \subseteq \{G_1, \ldots, G_w\} \setminus \mathcal{Y}_1\)
    with \(|\mathcal{Y}_2| < 12q\) and \(Z_2 \cap V(\mathcal{G}_1) \subseteq V(\bigcup\mathcal{Y}_1)\) such that every \(\mathcal{G}_1\)-jump in \(G - (Z_1\cup Z_2)\) has an endpoint in \(\bigcup\mathcal{Y}_2\).
    Let \(\mathcal{G}_2\) denote the subsequence of \(\mathcal{G}_1\) obtained by deleting the graphs belonging to \(\mathcal{Y}_2\).

    Observe that deleting a vertex from a graph with maximum degree at most \(d\) increases the number of connected components by at most \(d-1\).
    Every graph in \(\mathcal{Y}_1\) has maximum degree at most \(d\), \(|\mathcal{Y}_1| < 12q\) and \(|Z_1 \cup Z_2| < 8q\), so the graph \(\bigcup\mathcal{Y}_1 - (Z_1 \cup Z_2)\) has less than \(12q + 8q\cdot(d-1) = (8d+4)q\) components.
    
    Ley us say that a component \(G'\) of \(\bigcup\mathcal{Y}_1 - (Z_1 \cup Z_2)\) \emph{sees}
    a graph \(G_i\) belonging to \(\mathcal{G}_2\) if there exists a \(V(G')\)--\(V(G_i)\) path in \(G - (Z_1 \cup Z_2)\) without inner vertices in \(V(\mathcal{G}_1)\).
    Observe that if \(G'\) sees two graphs \(G_i\) and \(G_j\) belonging to \(\mathcal{G}_2\), then there exists
    a \(V(G_i)\)--\(V(G_j)\) path in \(G - (Z_1 \cup Z_2)\) without inner vertices in \(V(\mathcal{G}_1)\). As every \(\mathcal{G}_1\)-jump in \(G - (Z_1 \cup Z_2)\) has an endpoint in \(\mathcal{Y}_2\), this implies \(G_i\) and \(G_j\) are not \(\mathcal{G}_1\)-independent. There do not exist three pairwise non-\(\mathcal{G}_1\)-independent graphs in \(\mathcal{G}_1\), so each component \(G'\) of \(\bigcup\mathcal{Y}_1-(Z_1 \cup Z_2)\) sees at most two graphs from \(\mathcal{G}_2\).
    Since \(\bigcup\mathcal{Y}_1-(Z_1 \cup Z_2)\) has at most \((8d+4)q\) components,
    we conclude that the components of \(\bigcup\mathcal{Y}_1-(Z_1 \cup Z_2)\) see in total at most \((16d+8)q\) graphs from \(\mathcal{G}_2\). Let \(\mathcal{Y}_3\) denote the set of all these graphs.

    Let \(\mathcal{Y} = \mathcal{Y}_1 \cup \mathcal{Y}_2 \cup \mathcal{Y}_3\), and let \(Z = Z_1 \cup Z_2\). We have \(|\mathcal{Y}| = |\mathcal{Y}_1| + |\mathcal{Y}_2| + |\mathcal{Y}_3| < 12q+12q+(16d+8)q=16(d+2)q\), \(|Z| = |Z_1| + |Z_2| < 4q+4q = 8q\), and \(Z \cap V(\mathcal{G}) \subseteq V(\bigcup\mathcal{Y})\).
    Suppose towards a contradiction that there exists a \(\mathcal{G}\)-jump \(R\) in \(G - Z\) with an endpoint \(x \in V(\mathcal{G}) \setminus V(\bigcup\mathcal{Y})\).
    Let \(R'\) be an initial segment of \(R\) between \(x\) and a vertex \(y \in V(\bigcup\mathcal{G})\) such that no inner vertex of \(R'\)
    belongs to \(\bigcup\mathcal{G}\). Then \(R'\) is an \(\mathcal{G}\)-jump in \(G - Z_1\). Hence some endpoint of \(R'\) lies in a graph from \(\mathcal{Y}_1\).
    Since \(x \not \in V(\bigcup\mathcal{Y})\) and \(R\) is disjoint from \(Z\), the endpoint \(y\) belongs to a component \(G'\) of \(\bigcup\mathcal{Y}_1 - Z\).
    Since \(x\) does not lie on a graph from \(\mathcal{Y}\), it lies on a graph from \(\mathcal{G}_2\). That path is seen by \(G'\), so \(x\) lies on a path from \(\mathcal{Y}_3\),
    This contradiction proves that \(Z\) and \(\mathcal{Y}\) satisfy the lemma.
\end{proof}

\subsection{Isolating a mesh}

We now generalise the results from the previous subsection from one to two dimensions, so that, for a mesh \(M\) which does not control a model of \(K_t\), we have a control over \(V(M)\)-paths with endpoints ``far apart'' in \(M\). 

Let \(w, h \ge 2\) be integers, let \(\eta\) be a model of the \((w \times h)\)-grid in a graph \(G\). We denote the union of all branch sets of \(\eta\) by \(\bigcup\eta\). The \emph{coordinates} of \(x \in \bigcup\eta\) in \(\eta\) is the pair \((i, j) \in [w]\times [h]\) with \(x \in \eta((i, j))\). Two vertices \(x, x' \in \bigcup \eta\) are \emph{independent} in \(\eta\) if their coordinates do not lie on one face of the \((w \times h)\)-grid. In other words, if \(x\) and \(x'\) have coordinates \((i, j)\) and \((i', j')\) respectively, then \(x\) and \(x'\) are independent in \(\eta\) if and only if at least one of \((i, j)\) and \((i', j')\) belongs to \({[2,w-1]}\times [2, h-1]\), and we have \(|i-i'| \ge 2\) or \(|j-j'| \ge 2\). An \emph{\(\eta\)-jump} is a \((\bigcup\eta)\)-path between two independent vertices in \(\eta\). We say that \(\eta\) is \emph{isolated} in \(G\) if there is no \(\eta\)-jump in \(G\).

\begin{figure}
    \centering
    \includegraphics{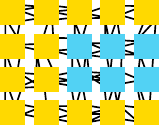}
    \caption{The \(([3, 5] \times [2, 3])\)-submodel of a model of the \((5 \times 4)\) grid is an \((3 \times 2)\)-submodel.}
    \label{fig:grid-submodel}
\end{figure}

For a model \(\eta\) of the \((w \times h)\)-grid in a graph \(G\), we say that a model \(\eta'\) of the \((w' \times h')\)-grid with \((w', h') \in [w] \times [h]\) is an \emph{\((w' \times h')\)-submodel} of \(\eta\) if there exist \(i_0 \in [0, w-w']\) and \(j_0 \in [0,h-h']\) such that \(\eta'((i, j)) = \eta((i_0+i, j_0+j))\) for each \((i, j) \in {[w']} \times [h']\). Such \(\eta'\) is also called an \emph{\(([i_0+1,i_0+w'] \times [j_0+1,j_0+h'])\)-submodel} of \(\eta\). See \cref{fig:grid-submodel}.
If \(M\) is a \((w \times h)\)-mesh with vertical paths \(P_1, \ldots, P_n\) and horizontal path \(Q_1, \ldots, Q_n\), then we say that a model \(\eta\) of the \((w \times h)\) grid in \(M\) is \emph{associated} with \(M\) if \(V(P_i \cap Q_j) \subseteq \eta((i, j))\) for \((i, j) \in {[w]\times[h]}\).

In this subsection, we show that for every graph \(G\) with an \((n \times n)\)-mesh \(M\) which does not control a model of \(K_t\), and with an associated model \(\eta\) of the \((n \times n)\)-grid,  we can delete a set \(Z\) of \(\mathbf{O}(t^3)\) vertices, so that every \(\eta\)-jump has an endpoint in a set of \(\mathbf{O}(t^3)\) exceptional vertical and horizontal paths of \(\eta\) (the first \(t+1\) and the last \(t+1\) vertical/horizontal paths are always implicitly considered exceptional). This implies existence of an isolated \((n/\mathbf{O}(t^3) \times n/\mathbf{O}(t^3))\)-submodel of \(\eta\) in \(G - Z\).

\begin{lemma}\label{IsolatedGridSubmodel}
    Let \(t, n \ge 1\) be integers with \(n \ge 2t+5\), let \(G\) be a graph, let \(M\) be an \(n\)-mesh in \(G\), and let \(\eta\) be a model of the \((n \times n)\)-grid in \(G\) associated with \(M\) such that \(\bigcup \eta = V(M)\). Then either
    \begin{itemize}
        \item there exists a model \(\mu\) of \(K_t\) in \(G\) which is controlled by \(M\), or
        \item there exist a set \(Z \subseteq V(G)\) with \(|Z|<16t^3\) and subsets \(I^*, J^* \subseteq [t+2, n-t-1]\) with \(|I^*|, |J^*| < 96t^3\) such that for each \((i, j) \in {[t+2, n-t-1]}^2\) with \(i \not\in I^*\) and \(j \not\in J^*\), the branch set \(\eta((i, j))\) contains no vertex from \(Z\) and no endpoint of an \(\eta\)-jump in \(G - Z\).
    \end{itemize}
    Moreover there exists \(\textbf{O}(|E(G)|)\)-time algorithm which finds \(\mu\) or \((Z, I^*, J^*)\) as above.
\end{lemma}

\begin{proof}
    Let us call an \(\eta\)-jump between branch set \(\eta((i, j))\)
    and \(\eta((i', j'))\) \emph{wide} if \(|i-i'| \ge 2\), and \emph{tall} if \(|j-j'| \ge 2\). The idea of the proof is to handle the wide and tall jumps separately, with two symmetric applications of \cref{PPjumps}. We will consider the \((n-2t \times n-2t)\)-submodel of \(\eta\) obtained by discarding \(t\) vertical/horizontal paths from each side. In order to handle, say, wide jumps, we will consider the sequence of subgraphs of \(M\) induced by the vertical paths of the submodel, and apply \cref{PPjumps} to that sequence. If as an outcome we get many jumps, then we extend them to middle jumps in a \((w' \times 2t)\)
    mesh contained in \(M\), and we get a model of \(K_t\). Otherwise,
    we can delete a small set of vertices so that the jumps have endpoints in a small set of vertical paths. Then we repeat the argument for horizontal paths.
    
    We want to ensure that if the outcome of either application is a set of \(\mathcal{G}\)-jumps, then we get a model of \(K_t\), so we reserve a small number of outer layers of \(\eta\) for that purpose.
    
    Let \(I = [n]\), and split \(I\) into three intervals \(I_{-1} = [1, t]\), \(I_0 = [t+1, n-t]\), \(I_1 = [n-t+1, n]\).
    By partitioning the vertical and horizontal paths of the \((n \times n)\)-grid according to the partition \(\{I_{-1}, I_0, I_1\}\),
    we obtain a partition of \(\bigcup\eta\) into nine parts, which we name naturally \emph{top-left}, \emph{top}, \emph{top-right}, \emph{left}, \emph{center}, \emph{right}, \emph{bottom-left}, \emph{bottom}, \emph{bottom-right}.

    We first focus on the vertical paths.
    Let \(S_{\updownarrow}\) denote the union of the top part and the bottom part of \(\eta\), i.e.\ the union of all branch sets \(\eta((i, j))\) with \(i \in I_0, j \in I_{-1} \cup I_1\). Further, we partition the ``center'' part into \(|I_0|\) graphs \((G_{t+1}, \ldots, G_{n-t})\), where for each \(i \in I_0\), \(G_i\) is the subgraph of \(M\) induced by the branch sets of the form \(\eta((i, j))\) where \(j \in I_0\).
    Furthermore, we let \(G_t\) denote the subgraph of \(M\) induced by the union of the top-left, left, and bottom-left parts, i.e.\ the union of all branch sets \(\eta((i, j))\) with \(i \in I_{-1}, j \in I\), and symmetrically, we let \(G_{n-t+1}\) denote the subgraph of \(M\) induced by the union of the top-right, right and bottom-right paths, i.e.\ the union of all branch sets \(\eta((i, j))\) with \(i \in I_{1}, j \in I\). See \cref{fig:mesh-partition}.
    Let \(\mathcal{G} = (G_t, \ldots, G_{n-t+1})\).
    Note that each \(G_i\) has maximum degree at most four.
    \begin{figure}
        \centering
        \includegraphics{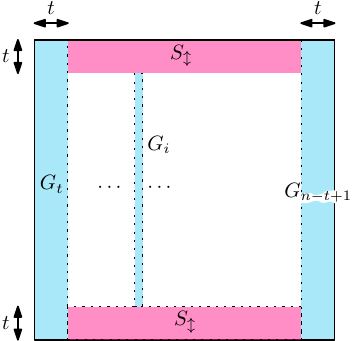}
        \caption{A schematic representation of the partition of the \(n\)-mesh into the set \(S_{\updownarrow}\) and the graphs \(G_t, \ldots, G_{n-t+1}\)}
        \label{fig:mesh-partition}
    \end{figure}
    Let \(d \coloneqq  4\) and \(q \coloneqq  t^3\).
    Apply \cref{PPjumps} to \(d\), \(q\), \(G - S_{\updownarrow}\) and \(\mathcal{G}\).
    
    Suppose first that there exists a subsequence
    \(\mathcal{G}' = (G_1', \ldots, G_{w'}')\) of \(\mathcal{G}\) and pairwise disjoint \(\mathcal{G}'\)-jumps \(R_1, \ldots, R_{q}\) in \(G - S_{\updownarrow}\) with endpoints in pairwise \(\mathcal{G}'\)-independent graphs from \(\mathcal{G}'\).
    Let \(M^*\) be the \((w' \times 2t)\)-submesh of \(M\) whose vertical paths are those vertical paths of \(M\) which intersect \(V(\mathcal{G}')\), and whose horizontal paths are subpaths of the horizontal paths \(Q_1, \ldots, Q_t, Q_{n-t+1}, \ldots, Q_n\) of \(M\) respectively. Then the paths \(R_1, \ldots, R_{q-2}\) can be extended to pairwise disjoint \(M^*\)-middle jumps in \(G\) with pairwise independent endpoints.
    By \cref{KtFromJumps}, \(M^*\) controls a model \(\mu\) of \(K_t\) in \(G\), and thus \(M\) controls \(\mu\) in \(G\), so the lemma follows.

    Hence we assume that there exist a set \(Z_1 \subseteq V(G)\) with \(|Z_1| < 8t^3\), and a subset \(\mathcal{Y}_1 \subseteq \{G_{t}, \ldots, G_{n-t+1}\}\) with \(|\mathcal{Y}_1| < 16(d+2)q=96q=96t^3\) 
    and \(Z_1 \cap V(\mathcal{G}) \subseteq V(\bigcup\mathcal{Y}_1)\) such that every \(\mathcal{G}\)-jump in \(G - S_{\updownarrow} - Z_1\) has both endpoints in \(V(\bigcup\mathcal{Y}_1)\).
    Hence every \(V(M)\)-path \(R\) in \(G - Z_1\) with an endpoint in a graph \(G_i\) which does not belong to \(\mathcal{Y}_1\) has the other endpoint either in \(S_{\updownarrow}\) or in a graph \(G_{i'}\) with \(|i-i'| \le 1\). In particular, if \(i \in [t+2, n-t-1]\), then \(i' \in [t+1, n-t]\).
    Let \({I}^*\) be the set of all \(i \in [t, n-t+1]\) with \(G_i \in \mathcal{Y}_1\), so that \(|{I}^*| < 96t^3\), and for each \((i, j) \in {[t+2, n-t-1]}^2\) with \(i \not \in I^*\), we have
    \begin{itemize}
        \item the branch set \(\eta((i, j))\) contains no vertex from \(Z_1\), and 
        \item every \(V(M)\)-path in \(G - Z_1\) with an endpoint in \(\eta((i, j))\) has the other endpoint in a branch set \(\eta((i', j'))\) such that \(i' \in [t+1, n-t]\) and if \(j' \in [t+1, n-1]\), then \(|i-i'| \le 1\).
    \end{itemize}

    Now for the vertical paths, we can repeat the above arguments with
    ``horizontal'' and ``vertical'' swapped. Hence we assume that there exist a subset \(Z_2 \subseteq V(G)\) and a subset \({J}^* \subseteq [t, n-t+1]\) with properties symmetric to the properties of \(Z_1\) and \(I^*\), that is, \(|Z_2| < 8t^3\), \(|{J}^*| < 96t^3\), and for each \((i, j) \in {[t+2, n-t-1]}^2\) with \(j \not \in J^*\), we have
    \begin{itemize}
        \item the branch set \(\eta((i, j))\) contains no vertex from \(Z_2\), and 
        \item every \(V(M)\)-path in \(G - Z_2\) with an endpoint in \(\eta((i, j))\) has the other endpoint in a branch set \(\eta((i', j'))\) such that \(j' \in [t+1, n-t]\), and if \(i' \in [t+1, n-1]\), then \(|j-j'| \le 1\).
    \end{itemize}

    Let \(Z = Z_1 \cup Z_2\), and note that \(|Z| \le |Z_1| + |Z_2| < 16q\). By combining the properties listed above, we conclude that for each \((i, j) \in {[t+2, n-t-1]}^2\) with \(i \not \in I^*\) and \(j \not \in J^*\), the branch set \(\eta((i, j))\) contains no vertex from \(Z\), and if \(\eta((i, j))\) contains an endpoint of an \(V(M)\)-path in \(G - Z\), then the other endpoint belongs to a branch set \(\eta((i', j'))\) with \((i', j') \in {[t+1, n-t]}^2\) and \(|i-i'| \le 1\) and \(|j - j'| \le 1\). Hence the lemma is satisfied by \(Z\), \(I^*\) and \(J^*\).
\end{proof}

\subsection{Flattening an isolated mesh}

We are now ready to prove the Flat Mesh Theorem. The proof proceeds as follows. Given a graph \(G\) with an \((n \times n)\)-mesh \(M\),
we fix a model \(\eta\) of the \((n \times n)\)-grid associated with \(M\), and we use the results from the previous subsection to find a small set of vertices \(Z\) such that \(\eta\) is isolated in \(G - Z\) except for a small set of vertical and horizontal paths. Next, we consider appropriately chosen \((t - 3)(t-4)/2\) isolated \((n' \times n')\)-submeshes, and for each of them we consider a society containing it. We apply the Two Paths Theorems to each of these societies. If each application yields a cross, we get a model of \(K_t\) controlled by \(M\), and if some application yields a vortex-free rendition in a disk, the \((n' \times n')\)-submesh contained in 
that society is flat in \(G - Z\).

\begin{proof}[Proof of \cref{thm:flatmesh}]
    Fix integers \(t, n'\) with \(t \ge 5\), \(n' \ge 2\), let \(n = 100t^3(n'+2t+2)\), and fix a graph \(G\) and an \((n \times n)\)-mesh \(M\) in \(G\). We need to show that either \(M\) controls a model of \(K_t\) or there exist a set \(Z \subseteq V(G)\) with \(|Z| < 16t^3\) and an \(n'\)-mesh \(M'\) disjoint from \(Z\) which is disjoint from \(Z\) and flat in \(G - Z\).

    Let \(\eta\) be any model of the \(n \times n\) grid associated with \(M\).
    Apply \cref{IsolatedGridSubmodel} to \(t, n\), \(G\), \(M\), and \(\eta\). If \(M\) controls a model of \(K_t\) in \(G\)
    then the proof follows, so we assume that there exist \(Z \subseteq V(G)\) and \(I^*,J^* \subseteq [t+2,n-t-1]\) such that \(|Z| < 16t^3\), \(|I^*|, |J^*|< 96t^3\), and no branch set \(\eta((i, j))\) with \(i \in [t+2, n-t-1] \setminus I^*\) and \(j \in [t+2, n-t-1] \setminus J^*\) contains a vertex from \(Z\) or an endpoint of an \(\eta\)-jump in \(G - Z\).

    Let \(\bar{\eta}\) be an \((n \times n'+2t+2)\)-submodel of \(\eta\) which does not contain any branch set \(\eta((i, j))\) with \(j \in [t]\cup J^* \cup [n-t+1, n]\).
    Let \(c = \nicefrac{1}{2}(t-3)(t-4)\), and choose pairwise disjoint \((n'+2 \times n'+2)\)-submodels \(\bar{\eta}_1\), \ldots, \(\bar{\eta}_c\)
    such that for each \(a \in [c]\), we have
    \begin{itemize}
        \item each \(\bar{\eta}_a\) has a ``centered vertical alignment'', i.e.\ there exists \(i_0\) such that for any \((i, j) \in [n'+2]^2\) we have \(\bar{\eta}_a((i, j)) = \bar{\eta}((i_0 + i, t+j))\) 
        \item each \(\bar{\eta}_a\) is disjoint from \(Z\) and there is no  \(\bar{\eta}_a\)-jump in \(G - Z\), i.e, \(\bar{\eta}_a\) is disjoint from the branch sets \(\bar{\eta}((i, j))\) with \(i \in [t+2] \cup I^* \cup [n-t-1, n]\).  
    \end{itemize}
    See \cref{fig:clean-subgrid}.

    \begin{figure}
        \centering
        \includegraphics{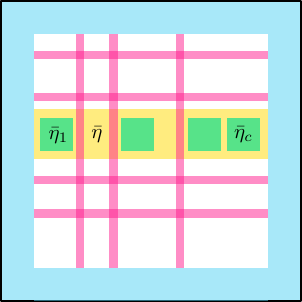}
        \caption{A schematic representation of the construction of the submodels \(\bar{\eta}_a\). The pink region represents the branch sets \(\eta((i, j))\) with \(i \in I^*\) or \(j \in J^*\). We will seek for crosses in the models \(\bar{\eta}_a\). If some \(\bar{\eta}_a\) has no cross, we will find a flat mesh in it, otherwise we get \((t-3)(t-4)/2\) crosses, which force a \(K_t\)-minor.}
        \label{fig:clean-subgrid}
    \end{figure}

    For each \(a \in [c]\), let \(J_a = M[\bigcup\bar{\eta}_a]\), and let \(G_a\) be the union of \(J_a\) and all \(M\)-bridges
    in \(G - Z\) with at least one attachment vertex in a branch set \(\bar{\eta}_a((i, j))\) with \((i, j) \in [2, n'+1]^2\).
    Note that by our choice of \(\bar{\eta}_a\), all remaining attachment vertices lie in \(\bigcup\bar{\eta}_a\).
    Next, for each \(a \in [c]\), let \(\Omega_a = (x_1, x_2, x_3, x_4)\) be a cylindrical ordering consisting of one vertex from each ``corner'' branch set \(\bar{\eta}_a((i, j))\) with \(i, j \in \{1, n'+2\}\), which is adjacent to a branch set \(\bar{\eta}((i, j))\) with \(j \in \{t, t+(n'+2)+1\}\).
    We assume that the vertices \(x_1, x_2, x_3, x_4\) are in the same cyclic order as the corresponding corners of the \((n'+2 \times n'+2)\)-grid.

    If every \((G_a, \Omega_a)\) has a cross, then we can build a model \(\nu\) of the graph \(H'_{c, 2t}\), as defined in \cref{sec:forcing}, such that there exists an integer \(j_0 \in [0, n-t]\)
    such that for all \(i \in [3c+1], j \in [t]\), the branch set \(\nu((i, j))\) intersects the \((j_0+j)\)th horizontal path of \(M\). Then by \cref{lem:Kt_from_crosses}, \(\nu\) controls a model of \(K_t\), and, by our choice of \(\nu\) so does \(M\).
    Therefore we assume that not all \((G_a, \Omega_a)\) have a cross.

    Fix \(a \in [c]\) such that there is no cross in \((G_a, \Omega_a)\).
    Let \(\bar{\eta}_a'\) denote the \(([2, n'+1] \times [2, n'+1])\)-submodel of \(\bar{\eta}_a\), and let \(M'\) be the \((n' \times n')\)-mesh contained in \(M[\bigcup\bar{\eta}_a'])\). Then \(M'\)
    is disjoint from \(Z\). We show that \(M'\) is flat in \(G - Z\).
    By the Two Paths Theorem, there is a vortex-free rendition \(\rho_0\) of \((G_a, \Omega_a)\) in a disk.
    Let \(C\) be the perimeter of \(M'\).
    Observe that \(M '\subseteq G_a\) and there exist four pairwise vertex-disjoint \(V(\Omega)\)--\(V(C)\) paths, and thus the cycle \(C\) is grounded in \(\rho_0\). Let \(\Delta\) denote the disk bounded by the trace of \(C\), and let \(\Delta'\) be a disk in \(\rho_0\) which contains \(\Delta\) in its interior, intersects the rendition only in nodes, and contains exactly those cells of \(\rho_0\) which intersect \(\Delta\). See \cref{fig:flat-mesh-subgraph}. Now, we can obtain a \(\Sigma\)-decomposition of \(G - Z\) witnessing that \(M'\) is flat by restricting \(\rho_0\) to \(\Delta'\), gluing a closed disk \(c_0\) to its boundary, and making \(c_0\) the only vortex of the rendition, which contains the part of \(G - Z\) not embedded in \(\Delta'\). 
    \begin{figure}
        \centering
        \includegraphics{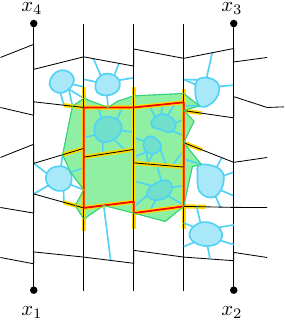}
        \caption{A society \((G_a, \Omega_a)\) which does not have a cross (\(\Omega_a = (x_1, x_2, x_3, x_4)\)). The graph \(G_a\) is obtained from the \((n'+2 \times n'+2)\)-mesh \(J_a\) by adding all \(M\)-bridges in \(G - Z\) with at least one attachment in \(\bigcup\bar{\eta}_a'\) (yellow). The perimeter of \(M'\) (red) is grounded in the vortex-free rendition of \((G_a, \Omega_a)\) in a disk. The green subgraph represents the subgraph induced by the non-vortex cells in the constructed \(\Sigma\)-decomposition of \(G\) witnessing that \(M'\) is flat.}
        \label{fig:flat-mesh-subgraph}
    \end{figure}

\end{proof}

\section{Flattening transactions}\label{sec:flatteTransaction}

In the previous section, we developed a tools for finding a flat submesh. In this section, we use these tools to find a flat subtransactions in societies with a cylindrical rendition and a large nest, and we state a result from \cite{KawarabayashiTW2018New} which allow us to extend a \(\Sigma\)-decomposition by a flat transaction.

Let \(\mathcal{C} = \{C_1, \ldots, C_s\}\) be a nest in some cylindrical rendition of a society \((G, \Omega)\). We say that a model \(\mu\) of \(K_t\) is \emph{controlled} by \(\mathcal{C}\) if \(s \ge t\) and for any \(i \in [s]\), \(x \in V(K_t)\) there exist \(t\) internally disjoint \(V(C_i)\)--\(\mu(x)\) paths in \(G\) (if \(V(C_i)\) intersects or is adjacent to \(\mu(x)\) we assume there are infinitely many internally disjoint \(V(C_i)\)--\(\mu(x)\) paths). Note that if
\(\mathcal{C}\) is a nest in a cylindrical rendition \(\rho\) which admits a radial linkage of order \(t\), and \(M\) is a mesh which has \(t\) horizontal paths contained in distinct cycles of \(\mathcal{C}\), then \(\mu\) is controlled by \(\mathcal{C}\) if and only if it is controlled by \(M\).

\begin{lemma}\label{lemma:findflattransaction}
    Let \(t \ge 5, w' \ge 2\) be integers, let \(w = 97t^3(w'+40) + 1\), let \((G, \Omega)\) be a society with a cylindrical rendition \(\rho\) around a vortex \(c_0\), and let
    \(\mathcal{C}\) be a nest in \(\rho\) with \(|\mathcal{C}| \ge t+16\). Let \(\mathcal{P}\) be a transaction in \((G, \Omega)\)  with \(|\mathcal{P}| = w\) which is monotone and orthogonal to \(\mathcal{C}\).
    Then, either there exists a model of \(K_t\) controlled by \(\mathcal{C}\), or there exist a set \(A \subseteq V(\sigma(c_0))\) with \(|A| \le 8t^3\), and a subtransaction \(\mathcal{P}' \subseteq \mathcal{P}\) with \(|\mathcal{P}'| = w'\) which consists of consecutive paths from \(\mathcal{P}\), is disjoint from \(A\) and isolated and flat in \((G - A, \Omega)\).
    Furthermore, there exists an \(\mathbf{O}(|E(G)|)\)-time algorithm, which finds \(\mu\) or \((A, \mathcal{P}')\) as above.
\end{lemma}
\begin{proof}
    By applying Menger's Theorem to \(\tilde{c_0}\) and \(V(C_1)\), we may assume (after possibly extending the vortex \(c_0\)) that there exist \(|\tilde{c_0}|\) pairwise vertex-disjoint \(\tilde{c_0}\)--\(V(C_1)\) paths in \(G\).
    Let \(\{S_x : x \in \tilde{c_0}\}\) be a set of such paths, with each \(S_x\) having an endpoint in \(x\).
    Let \(G_0 = \sigma(c_0)\), and let \(G_1\) denote the part of \(G\) embedded outside the interior of \(c_0\), i.e.\ let \(G_1 = \bigcup_{c \in C(\rho) \setminus \{c_0\}}\sigma(c)\).
    Note that \(V(G_0) \cap V(G_1) = \tilde{c_0}\).

    Let \(P_1, \ldots, P_w\) be the paths in \(\mathcal{P}\), naturally indexed.
    Let \(X\) and \(Y\) denote the end segments of \(\mathcal{P}\), and for each \(i \in [w]\), let \(P_i^X\)
    and \(P_i^Y\) denote the components of \(P_i \cap G_1\) intersecting \(X\) and \(Y\) respectively.
    Let \(\{Q_i^Z : i \in [w], Z \in \{X, Y\}\}\) be a family of pairwise vertex-disjoint paths in \(C_1\) covering all vertices of \(C_1\) with each \(Q_i^Z\) containing \(P_i^Z \cap C_1\). For each \(x \in \tilde{c_0}\), let \(S_x'\)
    denote the shortest subpath of \(S_x\) from \(x\) to \(V(\bigcup\mathcal{P}) \cup V(C_1)\). For each \(i \in [w]\), , let \(H_i^0 = P_i \cup Q_i^X \cup Q_i^Y\), and let \(H_i\) denote the union of the graph \(H_i^0\) and all paths \(S_x'\) intersecting it. Observe that the maximum degree of each \(H_i\) is at most four and \(\tilde{c_0} \subseteq V(G_0) \cap V(H_1 \cup \cdots \cup H_w)\). Let \(\mathcal{H} = (H_1, \ldots, H_w)\). Note that if \(\mathcal{P}\) is crosscap, then \(C_1\) contains two edges between \(H_1\) and \(H_w\).

    Let \(G_0^+ = G_0 \cup \bigcup\mathcal{H}\).
    Apply \cref{PPjumps} to the graph \(G_0^+\) and the sequence \(\mathcal{H}\) with \(d=4\) and \(q=t^3\). Suppose first that there are a subsequence \(\mathcal{H}' = (H_{i_1}, \ldots, H_{i_{w'}})\) of \(\mathcal{H}\) and pairwise vertex-disjoint \(\mathcal{H}'\)-jumps \(R_1, \ldots, R_{t^3}\) in \(G\) with endpoints in \(2t^3\) pairwise independent graphs from \(\mathcal{H}'\).
    Let \(M\) be the \((w' \times 2t)\)-mesh such that for each \(i \in [w']\) the \(i\)th vertical path of \(M\) is a subpath of \(P_{i_i}\), for each \(j \in [t]\), the \(j\)th horizontal path of \(M\) is a path in \(C_{t+2-j}\) intersecting the paths \(P_{i_1}^X, \ldots, P_{i_{w'}}^X\), and for each \(j \in [t+1,2t]\), the \(j\)th horizontal path
    of \(M\) is a path in \(C_{j-t+1}\) intersecting the paths \(P_{i_1}^Y, \ldots, P_{i_{w'}}^Y\).
    For each path \(R_a\) with endpoints in, say, \(V(H_i)\) and  \(V(H_j)\), let \(R_a'\) be the \(V(P_i)\)--\(V(P_j)\) path in \(H_i \cup H_j \cup R_a\). This way, the paths
    \(R_1', \ldots, R_{t^3}'\) are pairwise vertex-disjoint \(M\)-middle jumps in \(G\) with endpoints in pairwise non-consecutive vertical paths of \(M\).
    Then, by \cref{KtFromJumps}, \(M\) controls a model of
    \(K_t\), and since \(M\) has \(t\) distinct vertical paths
    contained in distinct cycles from \(\mathcal{C}\), that model is controlled by \(\mathcal{C}\) as well.

    Therefore we assume that the other outcome of \cref{PPjumps} holds: there exist \(A \subseteq V(G_0^+)\) with \(|A| < 8t^3\) and a subset \(\mathcal{Y} \subseteq \{H_1, \ldots, H_w\}\) with \(|\mathcal{Y}| < 16\cdot 6t^3 = 96t^3\) and with \(A \cap V(\mathcal{H})\subseteq V(\bigcup\mathcal{Y})\) such that every \(\mathcal{H}\)-jump in \(G-A\) has both endpoints in \(\bigcup\mathcal{Y}\). No vertex \(z \in V(\bigcup\mathcal{H}) \setminus V(G_0)\) is an endpoint of a \(\mathcal{H}\)-jump in \(G_0\cup\bigcup\mathcal{H}\), so we may assume that \(A \subseteq V(G_0) = V(\sigma(c_0))\).
    Consider the subsequence \((H_1^*, \ldots, H_{97t^3+1}^*)\) of \(\mathcal{H}\) where \(H_i^* = H_{(i-1)(w'+40)+1}\).
    
    The graphs \(H_i^*\) split \(\mathcal{H}\) into \(97t^3\) blocks of length \(w'+39\) each, and since \(|\mathcal{Y}| < 96t^3\), at least \(t^3\) of them
    are disjoint of \(\mathcal{Y}\).
    Let \((H_1', \ldots, H'_{w'+39})\) be such a block disjoint from \(\mathcal{Y}\), and thus from \(A\).
    Let \(\{P_1', \ldots, P_{w'+39}'\}\) be the subtransaction of \(\mathcal{P}\) where each \(P_i'\) is the path such that \(P_i' \subseteq H_i'\).
    Let \(X'\) and \(Y'\) denote the end segments of \(\{P_1', \ldots, P_{w'+39}'\}\) so that \(X' \subseteq X\) and \(Y' \subseteq Y\).
    For each \(Z' \in \{X', Y'\}\), let \(G_1^{Z'}\) denote the union of \(P_1' \cup P_{w'+39}'\) and all \(V(P_1') \cup V(P_{w'+39}') \cup V(Z')\)-bridges in \(G_1\) with an attachment in an inner vertex of \(Z'\). Next, let \(G'\) be the graph obtained as the union of \(G_1^{X'} \cup G_1^{Y'} \cup P_1' \cup \cdots \cup P_{w'+39}'\) and all \(V(\bigcup\mathcal{H} - A)\)-bridges in \(G_0^+ - A\) with an attachment in \(H_{20}' \cup \cdots \cup H_{w'+20}'\).
    Not that any such a bridge has all its attachments in \(H_{19}' \cup \cdots \cup H_{w'+21}'\).
    Apply the Two Paths Theorem to \((G', \Omega')\) where \(\Omega'\) is the concatenation of \(X'\) and \(Y'\).
    If the society admits a vortex-free rendition on a disk, then it contains the strip
    society of the transaction \(\mathcal{P}' = \{P_{20}', \ldots, P'_{w'+19}\}\) in \((G - A, \Omega)\), and thus \(\mathcal{P'}\) is flat in \((G - A, \Omega)\), so the lemma holds.
    Therefore we assume that there is a cross in the \((G', \Omega')\). Observe that \((G', \Omega')\) admits a cylindrical rendition with a nest \(\mathcal{C}' = (C_1', \ldots, C_{16}')\) where \(C_j' \subseteq C_j \cup P_{18-j}' \cup P_{w'+22+j}'\).
    By \cref{lemma:reroute_crosses}, there exist two disjoint paths \(R_1\) and \(R_2\) in the inner graph of \(C_1'\), with each \(R_a\) having endpoints in vertices \(x_i\) and \(x_{2+i}\) where \(x_1 \in V(G_1^{X'}) \cap V(C_{16}) \cap V(P_2')\), \(x_2 \in V(G_1^{X'}) \cap V(C_{16}) \cap V(P_{w'+38}')\),
    \(x_3 \in V(G_1^{Y'}) \cap V(C_{16}) \cap V(P_{w'+38}')\) and \(x_3 \in V(G_1^{Y'}) \cap V(C_{16}) \cap V(P_{w'+2}')\).
    We may assume that such paths \(R_1, R_2\) exists for each block
    \((H_1', \ldots, H'_{w'+39})\) disjoint from \(\mathcal{Y}\).
    Now the pairs of paths \((R_1, R_2)\) and the graphs \(H_i^*\) can be used to build a model \(\nu\) of the graph \(H_{(p-3)(p-4)/2, 2t}'\) from \cref{sec:forcing} such that for any \(i \in [3(p-3)(p-4)/2+1]\) and \(j \in [t]\) we have \(\nu((i, j)) \cap V(C_{16+t+1-j}) \neq \emptyset\).
    By \cref{lem:Kt_from_crosses}, \(\nu\) controls a model of \(K_t\),
    and hence so does \(\mathcal{C}\).
\end{proof}

\subsection{Reconciling a flat transaction with a cylindrical rendition}
When we find a large monotone isolated and flat transaction using \cref{lemma:findflattransaction}, we will need to extend the rendition with the strip society of the transaction. This can be achieved with Lemma 5.15 from \cite{KawarabayashiTW2021Quickly}.
\begin{proposition}[Kawarabayashi, Thomas, and Wollan \cite{KawarabayashiTW2021Quickly}]\label{lem:reconciliation}
    Let \((G, \Omega)\) be a society and \(\rho\) a cylindrical rendition around a vortex \(c_0\). Let \(\mathcal{C} = (C_1, \ldots, C_s)\) be a nest in \(\rho\) and let \(\mathcal{Q}\) be an exposed%
    \footnote{Recall that a transaction which is called unexposed in \cite{KawarabayashiTW2021Quickly} is called exposed in this paper.}
    monotone transaction in \((G, \Omega)\) which is orthogonal to \(\mathcal{C}\) of order at least three. Let \(Y_1\) and \(Y_2\)
    be the two segments obtained by deleting the end segments of \(\mathcal{Q}\) from \(\Omega\). Assume that there exists a linkage \(\mathcal{P} = \{P_1, P_2\}\) such that \(P_i\) links \(Y_i\) and \(V(\sigma(c_0))\) for \(i \in [2]\), \(\mathcal{P}\) is disjoint from \(\mathcal{Q}\), and \(\mathcal{P}\) is orthogonal to \(\mathcal{C}\).

    Let \(i \ge 7\), let \((G', \Omega')\) be the inner society of \(C_i\), and let \(\mathcal{Q}'\) be the restriction of \(\mathcal{Q}\) to \((G', \Omega')\). Let \(\rho' \) be the restriction of \(\rho\) to be a cylindrical rendition of \((G', \Omega')\). Then \(\mathcal{Q}'\) is exposed and monotone in \(\rho'\); moreover, \(\mathcal{Q}'\) is a crosscap transaction if and only if \(\mathcal{Q}\) is a crosscap transaction. If the \(\mathcal{Q}\) strip society is rural and isolated in \((G, \Omega)\), then the \(\mathcal{Q}'\) strip society is rural and isolated in \((G', \Omega')\).
\end{proposition}


\section{\texorpdfstring{$\heartsuit$}{<3} Making things cozy \texorpdfstring{$\heartsuit$}{<3}}\label{sec:cozy}
Throughout our proofs we will endeavour to ensure that the nest we work with is as close as possible to the boundary of the society we are considering.
Intuitively we are trying to be very conservative with our claims of what part of the graph is already embedded in the surface, since everything that is contained in the outer graph of the innermost cycle of a given nest in a cylindrical rendition of a society is already embedded outside of a vortex in a $\Sigma$-decomposition.
This property is surprisingly powerful given its relative simplicity.

In this section we will show that we can make any nest cozy in an efficient manner and then introduce analogous notions for transactions that are used in the more involved parts of our proof.

\subsection{Making a nest cozy}
Towards the goal of defining coziness formally, we first need to introduce a notion of direction inside a cylindrical rendition.
For nests this is only required with respect to a vortex $c_0$.
Thus the following notion defines disks that are described by traces in our cylindrical rendition that contain $c_0$.

\begin{definition}[$c_0$-disk]\label{def:cdisk}
    Let $(G,\Omega)$ be a society, $\rho = (\Gamma, \mathcal{D}, c_0)$ be a cylindrical rendition of $(G,\Omega)$ in a disk $\Delta$, and $T$ be a closed curve in $\Delta$, or a curve with both endpoints in $\mathsf{boundary}(\Delta)$.
    Then, if one of the two regions of $\Delta - T$ is a disk whose closure $d$ contains $c_0$, we call $d$ the \emph{$c_0$-disk of $T$}.
    
    If $T$ is the trace of a grounded cycle or a grounded $V(\Omega)$-path $C$ and $T$ has a $c_0$-disk $d$, then we say that $d$ is the \emph{$c_0$-disk of $C$}.
    For a grounded subgraph $H$ of $G$ containing a cycle that possesses a $c_0$-disk, the inclusion-wise minimal $c_0$-disk of any cycle in $H$ is the \emph{$c_0$-disk of $H$}.
\end{definition}

This allows us to define $C_i$-paths, for some $C_i$ in the nest, that stick out away from the vortex $c_0$.
Note that each cycle of a nest in a cylindrical rendition around a vortex $c_0$ has a $c_0$-disk.

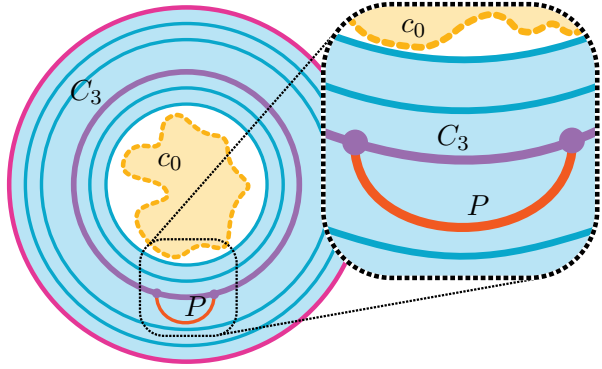
\begin{figure}[ht]
    \centering
    \begin{tikzpicture}

        \pgfdeclarelayer{background}
		\pgfdeclarelayer{foreground}
			
		\pgfsetlayers{background,main,foreground}

        \begin{pgfonlayer}{background}
        \node (C) [v:ghost] {{\includegraphics[width=12cm]{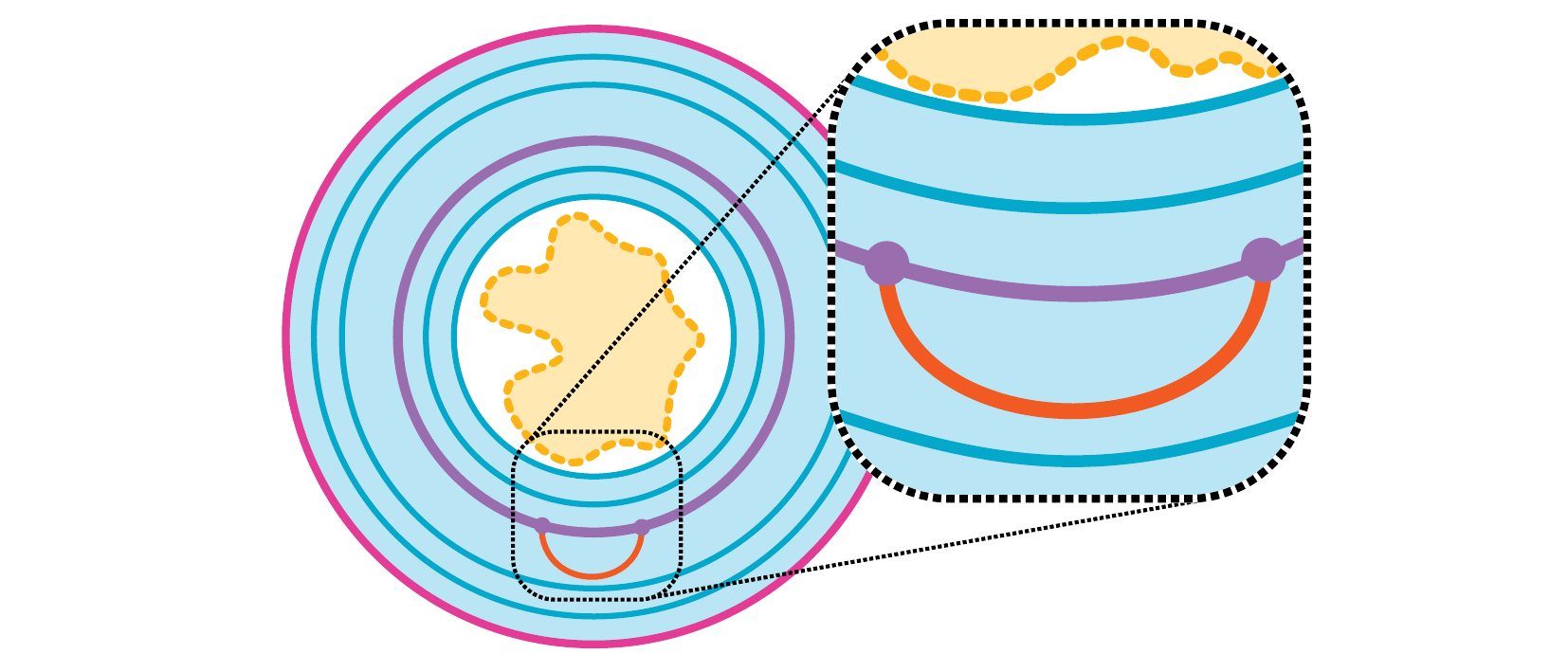}}};
        \end{pgfonlayer}{background}
			
        \begin{pgfonlayer}{main}

            \node (X) [v:ghost,position=0:0mm from C] {};

            \node (c0L) [v:ghost,position=170:17mm from X] {$c_0$};
            \node (c0R) [v:ghost,position=53:26mm from X] {$c_0$};

            \node (C3L) [v:ghost,position=140:14.5mm from c0L] {$C_3$};
            \node (C3R) [v:ghost,position=290:15.5mm from c0R] {$C_3$};

            \node (PL) [v:ghost,position=280:19.5mm from c0L] {$P$};
            \node (PR) [v:ghost,position=290:25mm from c0R] {$P$};
        
        \end{pgfonlayer}{main}

        \begin{pgfonlayer}{foreground}
        \end{pgfonlayer}{foreground}

    \end{tikzpicture}
    \caption{A path $P$ (drawn in \textcolor{CarrotOrange}{orange}) that sticks out away from the vortex $c_0$ (the \textcolor{BananaYellow}{yellow} area).}
    \label{fig:StickOut}
\end{figure}

\begin{definition}[Paths that stick out]\label{def:pathsstickingout}
    Let $(G,\Omega)$ be a society, let $\rho$ be a cylindrical rendition of $(G,\Omega)$ in a disk $\Delta$ with the vortex $c_0$, and let $C$ be a grounded cycle or $V(\Omega)$-path with a $c_0$-disk.

    Given a grounded $C$-path $P$ of length at least one, we say that $P$ \emph{sticks out towards $c_0$ (in $\rho$)} if the $c_0$-disk of $P \cup C$ is not the $c_0$-disk of $C$ and otherwise we say that $P$ \emph{sticks out away from $c_0$ (in $\rho$)}.
    See \cref{fig:StickOut} for an illustration.
\end{definition}

The definition of coziness is now simple to state.

\begin{definition}[Cozy nests]\label{def:cozynest}
Let $(G,\Omega)$ be a society, let $\rho$ be a cylindrical rendition of $(G,\Omega)$ in a disk $\Delta$ with a nest $\mathcal{C } = \{ C_1, \dots , C_s \}$ around the vortex $c_0$.
We say that $\mathcal{C}$ is \emph{cozy} if for every $i \in [s]$ and every grounded $C_i$-path $P$ that sticks out away from $c_0$ we have $V(P) \cap V(\Omega) \neq \emptyset$ or there exists a $j \in [s] \setminus \{ i \}$ such that $V(P) \cap V(C_j) \neq \emptyset$.
\end{definition}

One may imagine a cozy nest as one that has been stretched outwards, away from $c_0$ as much as possible to fully encompass, within its grasp, all paths that might leave the nest and return at any level without touching $\Omega$.
It will be convenient for us to always work with cozy nests as they allow for additional assumptions on the behaviour of radial linkages and transactions.
This is not difficult to guarantee once we have found a nest.

\begin{lemma}\label{lemma:makenestcozy}
Let $s \geq 1$ be an integer, $(G,\Omega)$ be a society, $\rho$ be a cylindrical rendition of $(G,\Omega)$ in a disk with a nest $\mathcal{C} = \{ C_1, \dots , C_s \}$ around the vortex $c_0$.
Then there exists a cozy nest $\mathcal{C}'$ of order $s$ in $(G,\Omega)$ around $c_0$, such that the union of $C_1$ and the outer graph of $C_1$ in $\rho$ contains $\bigcup \mathcal{C}'$, and an algorithm that finds $\mathcal{C}'$ in time $\mathbf{O}(s|E(G)|^2)$.
\end{lemma}
\begin{proof}
To prove the assertion we show by induction on $j$ that, given that $\{ C'_j, \dots , C'_s \}$ is a cozy nest such that $\{ C_1, \dots , C_{j-1}, C_j', \dots , C'_s \}$ is a nest around $c_0$ for some $j \in [2,s]$ we can find a cycle $C'_{j-1}$ such that $\{ C'_{j-1}, \dots , C'_s\}$ is cozy and $\{ C_1, \dots , C_{j-2}, C_{j-1}', \dots , C'_s \}$ is a nest around $c_0$.

For this, we first have to show that we can find a cycle $C'_s$ such that $\{ C'_s \}$ is a cozy nest and $\{ C_1, \dots , C_{s-1}, C'_s \}$ is a nest around $c_0$ if $s \geq 2$.
For this let $G^0_s$ be the graph built from outer graph of $C_s$ by deleting the vertices in $V(\Omega)$ and let $P$ be an arbitrary grounded $C_s$-path in $G^0_s$.
If no such path exists we set $C'_s \coloneqq C_s$ and we are done.
Thus $P$ exists and sticks out away from $c_0$.
Moreover, there exists a grounded cycle $C^1_s$ in $C_s \cup P$ whose outer graph $G^1_s$ is a proper subgraph of $G^0_s$.
Suppose the cycle $C^j_s$ and its outer graph $G^j_s$ have been constructed for $j \geq 1$.
Then, if there exists a grounded $C^j_s$-path $P$ in $G^j_s$ we define $C^{j+1}_s$, as before, to be the cycle in $C^j_s \cup P$ whose outer graph $G^{j+1}_s$ is a proper subgraph of $G^j_s$.
Otherwise we set $C'_s \coloneqq C^j_s$.
Notice that it takes only $\mathbf{O}(|E(G)|)$ steps to find a grounded $C^j_s$-path in $G^j_s$ and $j \leq |E(G)|$.
So $C'_s$ can be found in time $\mathbf{O}(|E(G)|^2)$.

Now suppose the cycles $C_j', \dots , C'_s$ have already been constructed.
Let $G^0_{j-1}$ be the graph obtained from the intersection of the outer graph of $C_{j-1}$ and the inner graph of $C'_j$ by deleting the vertices of $C'_j$.
Notice, that every grounded $C_{j-1}$-path that sticks out away from $c_0$ and from all cycles in $\{ C_1, \dots , C_{j-2} \} \cup \{ C'_j, \dots , C'_s\}$ must be a path in $G^0_{j-1}$.
In case no such path exists we set $C'_{j-1} \coloneqq C_{j-1}$ and otherwise we call this path $P$.
Then we define $C^1_{j-1}$ to be the unique cycle in $C_{j-1} \cup P$ whose outer graph is a proper subgraph of the outer graph of $C_{j-1}$ and we define $G^1_{j-1}$ to be the intersection of the outer graph of $C^1_{j-1}$ with $G^0_{j-1}$.
In particular, every $C^1_{j-1}$-path that sticks out away from $c_0$ and is disjoint from the cycles in $\{ C_1, \dots , C_{j-2} \} \cup \{ C'_j, \dots , C'_s \}$ must be a subgraph of $G^1_{j-1}$.

Similar to before, for $i \geq 1$ we can assume $C^i_{j-1}$ and $G^i_{j-1}$ to already be constructed.
We may further assume that every $C^i_{j-1}$-path that sticks out away from $c_0$ must be a subgraph of $G^i_{j-1}$.
If no such path exists we set $C'_{j-1} \coloneqq C^i_{j-1}$, otherwise let $P$ be such a path.
Finally, we set $C^{i+1}_{j-1}$ to be the unique cycle in $C^i_{j-1} \cup P$ whose outer graph is a proper subgraph of the outer graph of $C^i_{j-1}$ and $G^{i+1}_{j-1}$ to be the intersection of the outer graph of $C^{i+1}_{j-1}$ and $G^i_{j-1}$.
It follows that every $C^{i+1}_{j-1}$-path that sticks out away from $c_0$ and is disjoint from the cycles in $\{ C_1, \dots , C_{j-2} \} \cup \{ C'_j, \dots , C'_s\}$ must be a subgraph of $G^{i+1}_{j-1}$.

To complete the proof assume that the cycles $\{ C'_1, \dots , C'_s \}$ have been constructed by the mechanism above.
This took at most $\mathbf{O}(s|E(G)|^2)$ time.
For each $j \in [s]$ there exists some $i_j \in [0,|E(G)|]$ such that $C'_j = C^{i_j}_j$ and the graph $G^{i_j}_j$ has been defined.
From the construction, every $C'_j$-path that sticks out away from $c_0$ and which is disjoint from the cycles in $\{ C'_1, \dots , C'_{j-1} \} \cup \{ C'_{j+1}, \dots , C'_s \}$ must be a subgraph of $G^{i_j}_j$.
By choice of $i_j$ we have that $C'_j = C^{i_j}_j$.
This can only occur in the case where no such path exists.
Hence, $\{ C'_1, \dots , C'_s\}$ is a cozy nest around $c_0$.
\end{proof}

\subsection{Making transactions cozy (and tautening them)}
This notion can be translated to transactions without much trouble.
Unlike with nests, we also needs cozy transactions in vortex-free renditions.
Thus we designate some point in our disk as the reference standing in for the vortex a cylindrical rendition might have.

\begin{definition}[Cozy transaction]
    Let $(G,\Omega)$ be a society, $\rho$ be a rendition of $(G,\Omega)$ in a disk $\Delta$ with at most one vortex $c_0$, and let $\mathcal{P}$ be an transaction in $(G,\Omega)$ that is unexposed if $\rho$ has a vortex.
    If $\rho$ is vortex-free, let $c_0$ be a point in $\Delta$.
    
    We say that $\mathcal{P}$ is \emph{cozy (with respect to $c_0$)} if for every $P \in \mathcal{P}$ and every grounded $P$-path $Q$ that sticks out away from $c_0$ we have $V(P) \cap V(\Omega) \neq \emptyset$ or there exists a path $P' \in \mathcal{P} \setminus \{ P \}$ such that $V(Q) \cap V(P') \neq \emptyset$.
\end{definition}

We can efficiently transform unexposed transactions into a cozy transaction.
The proof for this fact is analogous to our proof for \Cref{lemma:makenestcozy}.
Thus we omit it and only provide the statement.

\begin{lemma}\label{lemma:maketransactioncozy}
Let $p \geq 1$ be an integer, $(G,\Omega)$ be a society, $\rho$ be a rendition of $(G,\Omega)$ in a disk $\Delta$ with at most one vortex $c_0$, and let $\mathcal{P}$ be a transaction of order $p$ in $(G,\Omega)$ that is unexposed if $\rho$ has a vortex.
Furthermore, if $\rho$ is vortex-free, let $c_0$ be a point in $\Delta$ such that the set of $c_0$-disks defined by the paths in $\mathcal{P}$ is totally ordered by the containment relation.

Then there exists a transaction $\mathcal{P}'$ of order $p$ in $(G,\Omega)$ with the same endpoints as $\mathcal{P}$ that is cozy (with respect to $c_0$) and unexposed, such that the set of $c_0$-disks defined by the paths in $\mathcal{P}'$ is totally ordered by the containment relation and there exists a path $P' \in \mathcal{P}'$ whose $c_0$-disk contains the traces of all paths in $\mathcal{P}$.
Moreover, there exists an algorithm that finds $\mathcal{P}'$ in time $\mathbf{O}(p|E(G)|^2)$.
\end{lemma}

Cozy transactions are an essential ingredient for the structure that we introduce in \Cref{def:transactionmesh}.
However, we will need a type of transaction that is cozy towards its centre in that definition.

\begin{definition}[Taut transaction]\label{def:tauttransaction}
    Let $p$ be an integer with $p \geq 3$, let $(G,\Omega)$ be a society with a vortex-free rendition $\rho$ in a disk $\Delta$ and let $\mathcal{P} = \{ P_1, \ldots , P_p \}$ be a planar transaction in $(G,\Omega)$ of order $p$ that is indexed naturally.
    Further, let $q_1 = \lfloor \nicefrac{p}{2} \rfloor$ and $q_p = \lceil \nicefrac{(p+1)}{2} \rceil$, and let $T_i$ be the trace of $P_{q_i}$ and $T_i'$ be the trace of $P_i$ in $\rho$ for both $i \in \{ 1, p \}$.

    For both $i \in \{ 1, p \}$, let $\Delta_i$ be the closure of the disk in $\Delta - T_i$ that contains the trace of $P_i$ and let $(G_i, \Omega_i)$ be society associated with the restriction of $\rho$ to $\Delta_i$.
    Further, let $v_i$ be a point in the disk of $\Delta - T_i'$ that does not contain any trace of a path in $\mathcal{P}$.

    We say that $\mathcal{P}$ is \emph{taut in $(G, \Omega)$} if $\{ P_1, \ldots , P_{q_1 - 1} \}$ is cozy in $(G_1, \Omega_1)$ with respect to $v_1$ and $\{ P_{q_2 +1}, \ldots , P_p \}$ is cozy in $(G_p, \Omega_p)$ with respect to $v_p$.
\end{definition}

Via an easy application of \Cref{lemma:maketransactioncozy}, we can also turn any transaction into a taut transaction.

\begin{corollary}\label{corollary:maketransactiontaut}
Let $p \geq 3$ be an integer, $(G,\Omega)$ be a society, $\rho$ be a vortex-free rendition of $(G,\Omega)$ in a disk $\Delta$, let $\mathcal{P} = \{ P_1, \ldots , P_p \}$ be a transaction of order $p$ in $(G,\Omega)$ that is indexed naturally, and let $q_1 = \lfloor \nicefrac{p}{2} \rfloor$ and $q_p = \lceil \nicefrac{(p+1)}{2} \rceil$.

Then there exists a taut transaction $\mathcal{P}' = \{ P_1', \ldots , P_p' \}$ of order $p$ with the same endpoints as $\mathcal{P}$ in $(G,\Omega)$ that is indexed naturally, with $P_{p_1} = P_{p_1}'$ and $P_{p_2} = P_{p_2}'$ and an algorithm that finds $\mathcal{P}'$ in time $\mathbf{O}(p|E(G)|^2)$.
\end{corollary}

\section{Finding crooked transactions}\label{sec:findcrooked}
We will need a special type of transaction introduced in \cite{RobertsonS1990Graph}.
Let $\mathcal{P}$ be a transaction in a society $(G,\Omega)$.
A path $P \in \mathcal{P}$ is called \emph{peripheral} if the endpoints of $P$ split $\Omega$ into two segments, one of which contains all of the endpoints of the other paths in $\mathcal{P}$.
We call a transaction \emph{crooked} if it does not contain a peripheral path.
A cross is itself an example of a crooked transaction.

There are essentially two ways in which crooked transactions are more helpful to us than regular transactions.
The first is neatly captured in the following observation that is a consequence of the fact that any transaction with an unexposed path must contain a peripheral path.

\begin{observation}\label{obs:crookedisexposed}
    Let $(G,\Omega)$ be a society with a cylindrical rendition $\rho$ around the vortex $c_0$ and let $\mathcal{P}$ be a crooked transaction in $(G,\Omega)$.
    Then $\mathcal{P}$ is an exposed transaction.
\end{observation}

This fundamental property of exposed transactions saves us some work, as we would otherwise have to recover an exposed transaction from an arbitrary transaction.
Though this would still be comparatively easy to accomplish if we did not insist on the exposed transaction we will find being crooked.
The second and more important reason we care about crooked transactions is buried deep in a single claim in \Cref{sec:societyclassification} (see \Cref{claim:auseforcrookedtransactions}).
Put simply, crooked transactions witness non-planarity in our cylindrical renditions in a very easy to grasp fashion without requiring us to carry around a complex structure.
This fact was also of use to Kawarabayashi et al., though their use of this second property is similarly buried deep within one of their longer proofs.

\subsection{The key lemma for finding crooked transactions}
Within any transaction one can easily find either a large crooked or a large planar transaction.
The following lemma is stated as Lemma 3.1 in \cite{KawarabayashiTW2021Quickly} without a running time.
Its proof consists of a simple counting argument, which easily witnesses the runtime we claim for it.

\begin{proposition}[Kawarabayashi, Thomas, and Wollan \cite{KawarabayashiTW2021Quickly}]\label{lem:crookedorplanar}
    Let $p \geq 1$ and $q \geq 3$ be integers, let $(G,\Omega)$ be a society, and let $\mathcal{P}$ be a transaction of order $p+q-2$ in $(G,\Omega)$.
    Then there exists $\mathcal{P}' \subseteq \mathcal{P}$ such that $\mathcal{P}'$ is either a planar transaction of order $p$ or a crooked transaction of order $q$.
    
    Furthermore, we can find one of these two transactions in time $\mathbf{O}((p+q)|V(G)|)$.
\end{proposition}

Our main interest is being able to find large crooked transactions in societies which do not have cylindrical renditions of small depth.
In \cite{RobertsonS1990Graph}, Robertson and Seymour go through significant efforts to prove that this is possible with a very modest bound on the required depth (see (6.1) and (11.11) in \cite{RobertsonS1990Graph}).
In \cite{KawarabayashiTW2021Quickly}, Kawarabayashi et al.\ provide a simpler proof of the same statement with slightly worse bounds (see Lemma 3.6 in \cite{KawarabayashiTW2021Quickly}).

\begin{proposition}[Kawarabayashi, Thomas, and Wollan \cite{KawarabayashiTW2021Quickly}]\label{lem:crookedexistence}
   Let $q \geq 4$ be an integer and let $(G,\Omega)$ be a society.
   Either there exists a crooked transaction of order $q$ in $(G,\Omega)$ or $(G,\Omega)$ has a cylindrical rendition of depth at most $6q$.
\end{proposition}

This lemma is a very interesting generalisation of the Two Paths Theorem that has so far gone somewhat underappreciated.\footnote{\cite{RobertsonS1990Graph} is generally only cited for the proof of a specific version of the Two Paths Theorem that it contains. However Robertson and Seymour mainly present this proof as a lead-up towards their versions of \Cref{lem:crookedexistence}, which itself has been largely overlooked in the literature.}
We want two things from this lemma.
First, we of course want to apply it.
However, secondly, we will also want to use all of the ideas that Kawarabayashi et al.\ present in their proof of the above lemma to prove a substantially stronger version of it in \Cref{sec:crookedmagic}.
Furthermore, we also want to convince the reader that there exists a version of \Cref{lem:crookedexistence} that is algorithmically efficient.
The proof presented in \cite{KawarabayashiTW2021Quickly} uses several arguments with unclear running times, such as several minimality assumptions.
We therefore devote this section to a thorough examination of the ideas of Robertson and Seymour, and Kawarabayashi, Thomas, and Wollan with the goal of adding a runtime-bound to \Cref{lem:crookedexistence}.

\subsection{Normalising a crooked transaction}\label{subsec:normalisingcrooked}
We start by showing that any given crooked transaction has to contain one of two prototypical types of crooked transactions.
Similar results are presented by Robertson and Seymour in section 7 of \cite{RobertsonS1990Graph}.
This subsection is mainly meant to give the reader a better intuition for what a crooked transaction is, though we will find a use for this intuition in \Cref{sec:multisocietycrooked}.

The types of crooked transactions Robertson and Seymour desire are much more specific than what we need, as they require the transaction to essentially be monotone except for at most four paths.
We will not need to be as restrictive, but we nonetheless adopt some of the terminology from \cite{RobertsonS1990Graph}, particularly the terms \emph{overpass}, \emph{leap}, and \emph{doublecross}.

\begin{figure}[ht]
    \centering
    \begin{tikzpicture}

        \pgfdeclarelayer{background}
		\pgfdeclarelayer{foreground}
			
		\pgfsetlayers{background,main,foreground}

        \begin{pgfonlayer}{background}
        \node (C) [v:ghost] {{\includegraphics[width=12cm]{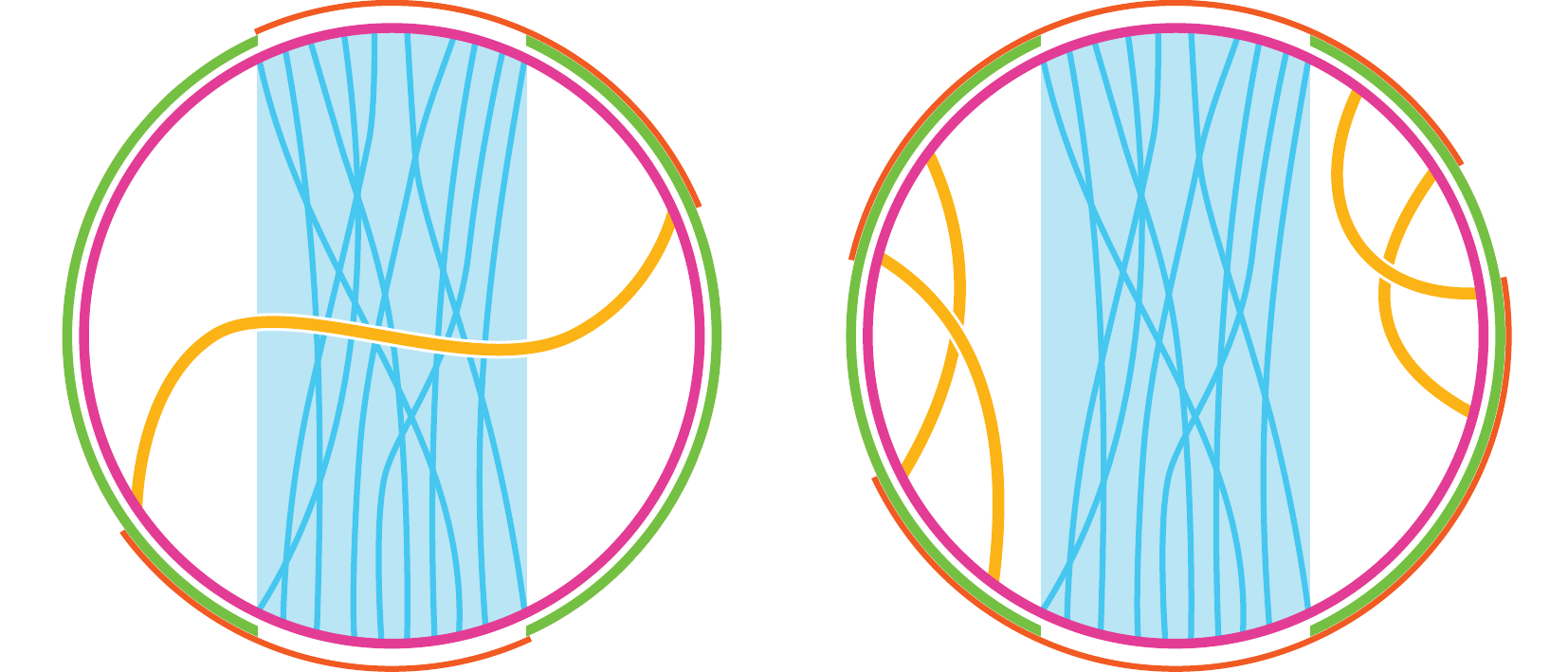}}};
        \end{pgfonlayer}{background}
			
        \begin{pgfonlayer}{main}

            \node (X) [v:ghost,position=0:0mm from C] {};

            \node (L) [v:ghost,position=180:30mm from X] {};
            \node (R) [v:ghost,position=0:30mm from X] {};

            \node (AL) [v:ghost,position=50:29mm from L] {$A'$};
            \node (BL) [v:ghost,position=235:29mm from L] {$B'$};

            \node (AR) [v:ghost,position=130:29mm from R] {$A'$};
            \node (BR) [v:ghost,position=305:29mm from R] {$B'$};

            \node (XL) [v:ghost,position=15:29mm from L] {$X$};
            \node (YL) [v:ghost,position=195:29mm from L] {$Y$};

            \node (XR) [v:ghost,position=15:29mm from R] {$X$};
            \node (YR) [v:ghost,position=195:29mm from R] {$Y$};

            \node (P1) [v:ghost,position=90:7mm from L] {$\mathcal{P}_1$};
            \node (P2) [v:ghost,position=90:7mm from R] {$\mathcal{P}_2$};

            \node (PL) [v:ghost,position=173:15mm from L] {$P$};

            \node (PR) [v:ghost,position=150:17mm from R] {$P$};
            \node (PPR) [v:ghost,position=220:14mm from R] {$P'$};

            \node (QR) [v:ghost,position=50:15mm from R] {$Q$};
            \node (QQR) [v:ghost,position=344:16mm from R] {$Q'$};

        \end{pgfonlayer}{main}

        \begin{pgfonlayer}{foreground}
        \end{pgfonlayer}{foreground}

    \end{tikzpicture}
    \caption{Two societies with crooked transactions (between the segments $A'$ and $B'$ marked in \textcolor{CarrotOrange}{orange}). The transaction $\mathcal{P}_1$ has an overpass $P$ and is a leap, the transaction $\mathcal{P}_2$ has two crosses, one in each of the two segments $X$ and $Y$, and is a doublecross.}
    \label{fig:LeapAndDoublecross}
\end{figure}
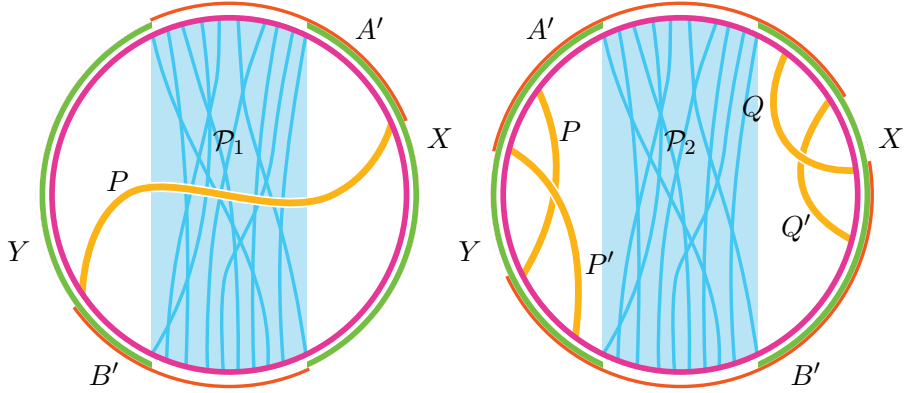

Let $(G,\Omega)$ be society, let $\mathcal{P}$ be a transaction of order at least 2 in $(G,\Omega)$, and let $P \in \mathcal{P}$ be some path.
We let $X,Y$ be the two distinct end segments of $\mathcal{P} \setminus \{ P \}$ and let $A,B$ be the two distinct segments of $\Omega$ that have one endpoint in $X$ and the other in $Y$.
If $P$ is an $A$-$B$-path, we call it an \emph{overpass} and if $\mathcal{P}$ contains an overpass, we call $\mathcal{P}$ a \emph{leap} (see \cref{fig:LeapAndDoublecross} for an illustration).
Note that the unique leap of order 2 is a cross.

Suppose instead that $\mathcal{P}$ has order at least 5 and let $P,P',Q,Q' \in \mathcal{P}$ be four distinct paths.
Let $X,Y$ be the two distinct end segments of $\mathcal{P} \setminus \{ P,P',Q,Q' \}$ and let $A,B$ be the two distinct segments of $\Omega$ that have one endpoint in $X$ and the other in $Y$.
If $P,P'$ and $Q,Q'$ respectively are a cross in $(G,\Omega)$, all endpoints of $P,P'$ are found in $A$, and all endpoints of $Q,Q'$ are found in $B$, we call $\mathcal{P}$ a \emph{doublecross} (see \cref{fig:LeapAndDoublecross} for an illustration).

Let us now demonstrate that any crooked transaction contains a large leap or a large doublecross.

\begin{lemma}\label{lem:leapordoublecross}
    Let $\ell,d$ be integers with $\ell \geq 2$ and $d \geq 1$, let $(G,\Omega)$ be a society, and let $\mathcal{P}$ be a crooked transaction of order $4(\ell - 2) + d$ in $(G,\Omega)$.
    
    Then there exists a leap of order at least $\ell$ or a doublecross of order at least $d+4$ in $(G,\Omega)$.
    
    Furthermore, finding the leap or the doublecross can be accomplished in time $\mathbf{O}(\ell|V(G)|)$.
\end{lemma}
\begin{proof}
    Let $X$ and $Y$ be the end segments of $\mathcal{P}$ in $\Omega$ and let $u$ be one of the two endpoints of $Y$, with $P$ being the path in $\mathcal{P}$ that also has the endpoint $u$.
    We let $v \in X$ be the other endpoint of $P$.
    Consider the segment of $\Omega$ with the same endpoints as $P$ that only intersects $Y$ in $u$ and let $X'$ be the intersection of that segment and $X$.

    We let $\mathcal{P}_{X'}$ be the set of paths in $\mathcal{P} \setminus \{ P \}$ that has an endpoint in $X'$.
    If $|\mathcal{P}_{X'}| \geq \ell - 1$, we have found a leap of order at least $\ell$.
    Thus we may suppose that $|\mathcal{P}_{X'}| \leq \ell - 2$.
    
    Let $w \in Y$ be the endpoint of a path $P'$ in $\mathcal{P}_{X'}$ that minimises the size of the subsegment $Y'$ of $Y$ with the endpoints $u$ and $w$.
    We let $\mathcal{P}_{Y'}$ be those paths in $\mathcal{P} \setminus \{ P' \}$ that have one endpoint in $Y'$.
    Note that $P \in \mathcal{P}_{Y'}$ and, due to the minimality of $Y'$, all endpoint of the paths in $\mathcal{P}_{Y'} \setminus \{ P \}$ in $X$ are found outside of $X'$.
    Thus, if $|\mathcal{P}_{Y'}| \geq \ell - 1$, we have found a leap of order $\ell$, with $P'$ being the overpass.
    We may therefore suppose that $|\mathcal{P}_{Y'}| \leq \ell - 2$.

    Note that $\mathcal{P}_{X'} \cup \mathcal{P}_{Y'}$ therefore contain at most $2\ell - 4$ paths, which already includes both $P$ and $P'$.
    We let $\mathcal{Q} = \mathcal{P} \setminus (\mathcal{P}_{X'} \cup \mathcal{P}_{Y'})$ and let $X'',Y''$ be the two end segments of $\mathcal{Q}$ in $\Omega$.
    Consider the two distinct segments $A,B$ of $\Omega$ that have one endpoint in $X''$ and the other in $Y''$ and note that all endpoints of $P,P'$ are contained in one of these two.
    W.l.o.g.\ we suppose this is true for $A$ and note that $P,P'$ is a cross in $(G,\Omega)$.

    We can now repeat the arguments we have just laid out, starting with choosing the common endpoint $u'$ of $B$ and $Y$.
    This allows us to either find a leap of order $\ell$ -- in which case we are done -- or find a new transaction $\mathcal{Q}' \subseteq \mathcal{Q}$ with $|\mathcal{Q}'| \geq |\mathcal{Q}| - (2\ell - 4) = d$ and a cross consisting of the two paths $Q,Q'$, such that $\mathcal{Q}' \cup \{ P,P',Q,Q' \}$ are a doublecross of order $d+4$.
    This completes our proof.
\end{proof}

\subsection{Finding a transaction in a society}\label{subsec:findingtransaction}
Before going further, we first address the fundamental task of finding a transaction in a society with depth above a certain bound.
The following is very straightforward, though somewhat inefficient.

\begin{lemma}\label{lem:findtransactioninefficient}
    Let $p$ be a positive integer and let $(G,\Omega)$ be a society with depth at least $p$.
    Then it is possible to find a transaction of order $p$ in $(G,\Omega)$ in time $\mathbf{O}(p|E(G)||V(G)|^2)$.
\end{lemma}
\begin{proof}
    For all vertices $u \in V(\Omega)$ and for all vertices $v \in V(\Omega) \setminus \{ u \}$ let $I$ and $I'$ be the two distinct segments of $\Omega$ whose endpoints are $u$ and $v$.
    We test via an algorithmic version of \Cref{prop:mengersthm}, whether there exists an $I$-$I'$-linkage in $G$.
    Since the depth of $(G,\Omega)$ is at least $p$, this will yield a transaction of order $p$ in $(G,\Omega)$.
\end{proof}

We can substantially reduce the runtime if we are willing to sacrifice the tight relation between the depth of the society and the size of the transaction we find.
In fact, we can even recover a witness for the non-existence of a large transaction, if we fail to find one, in the form of a linear decomposition of small adhesion.

\begin{lemma}\label{lem:findtransactionorlineardecomp}
    Let $p$ be a positive integer and let $(G,\Omega)$ be a society.
    Then it is possible to find either a transaction of order $p$ in $(G,\Omega)$ or a linear decomposition of \((G, \Omega)\)
    with adhesion less than \(p\) in time $\mathbf{O}(p|E(G)| \log |V(G)|)$.
\end{lemma}
\begin{proof}
    We apply the divide-and-conquer approach.
    Let \(v_1, \ldots, v_n\) denote the vertices of \(\Omega\) in order.
    We will first test whether there is an \(\{v_1, \ldots, v_{\lfloor n/2 \rfloor}\}\)-\(\{v_{\lfloor n/2 \rfloor+1}, \ldots, v_n\}\)-linkage of order \(p\) in \(G\). If yes, then the linkage is a transaction satisfying the lemma and the algorithm terminates. Otherwise, we get a separation \((G_1, G_2)\) of \(G\), and we recurse on the graphs \(G_1\) and \(G_2\). A more detailed description of the algorithm follows.

    During the run of the algorithm, we maintain \textsl{(i)} a partition of \(\Omega\) into pairwise disjoint segments \(I_1, \ldots, I_k\) such that for any \(v_{i'} \in I_{j'}\) and \(v_{i''} \in I_{j''}\), if \(j' < j''\), then \(i' < i''\), and \textsl{(ii)} a sequence \((G_1, \ldots, G_k)\) of pairwise edge-disjoint subgraphs of \(G\) with \(G_1 \cup \cdots \cup G_k = G\) such that for each \(j \in [k]\) we have \(I_j \subseteq V(G_j)\), and for any \(j_1, j_2, j_3 \in [k]\) with \(j_1 < j_2 < j_3\) we have \(V(G_{j_1}) \cap V(G_{j_3}) \subseteq V(G_{j_2})\) (hence for each \(j \in [k-1]\), the pair \((G_1 \cup \cdots \cup G_j, G_{j+1} \cup \cdots \cup G_k)\) is a separation of \(G\) of order \(|V(G_j) \cap V(G_{j+1})|\)).
    Additionally, we keep as an invariant that for each \(j \in [k-1]\) we have \(|V(G_j) \cap V(G_{j+1})| < p\) and there exist an \((I_1 \cup \cdots \cup I_j)\)--\((I_{j+1} \cup \cdots \cup I_n)\)-linkage \(\mathcal{P}_j\) in \(G\) of order \(|V(G_j) \cap V(G_{j+1})|\) (the linkage will not be stored explicitly, but will be easy to retrieve in \(\mathbf{O}(|E(G)|)\)-time). In each turn, we either find a transaction of order \(p\) and terminate the algorithm, or we partition a segment \(I_j\) with the corresponding graph \(G_j\) into two shorter segments \(I'\) and \(I''\) with corresponding graphs \(G'\) and \(G''\).
    Once each segment \(I_j\) contains only one vertex, we have \(k = n\), and
    \((V(G_1), \ldots, V(G_n))\) is a linear decomposition of \((G, \Omega)\) with adhesion less than \(p\), as desired.

    Next, we describe the procedure of splitting a segment \(I_j = v_{i_1} \cdots v_{i_2}\).
    Let \(i = \lfloor(i_1+i_2-1)/2\rfloor\), let \(I' = v_{i_1} \cdots v_{i}\), and let \(I'' = v_{i+1} \cdots v_{i_2}\).
    Furthermore, let \(X' = V(G_j) \cap V(G_{j-1})\), and let \(X'' = V(G_j) \cap V(G_{j+1})\) (assume \(V(G_0)=V(G_{k+1})=\emptyset\)).
    By our invariant, there is an \((I_1 \cup \cdots \cup I_{j-1})\)-\((I_j \cup \cdots \cup I_k)\)-linkage \(\mathcal{P}_{j-1}\) of order \(|X'|\) in \(G\). Each path in \(\mathcal{P}_{j-1}\) contains a distinct vertex from \(X'\), so there is an \((I_1 \cup \cdots \cup I_{j-1})\)-\(X'\) linkage \(\mathcal{P}'\) of order \(|X'|\) in \(G_1 \cup \cdots \cup G_{j-1}\). Symmetrically, there is an \(X''\)--\((I_{j+1} \cup \cdots \cup I_k)\)-linkage \(\mathcal{P}''\) of order \(|X''|\) in \(G_{j+1} \cup \cdots \cup G_k\).
    Apply Menger's Theorem to \(G_j\) and the sets \(X' \cup V(I')\) and \(X'' \cup V(I'')\). If we get a linkage of order \(p\), then we can use the paths from \(\mathcal{P}'\) and \(\mathcal{P}''\) to extend the paths with endpoints in \(X'\) or \(X''\) so that they have endpoints in \(v_1 \cdots v_{i_1-1}\) or \(v_{i_2+1} \cdots v_n\) respectively, and we obtain a desired transaction of order \(p\).
    Otherwise, we get an \((X' \cup V(I'))\)-\((X'' \cup V(I''))\)-separation \((G', G'')\) of \(G_j\) with order less than \(p\), and an \((X' \cup V(I'))\)-\((X'' \cup V(I''))\)-linkage of the same order. As before, we can extend it to a \(\{v_1, \ldots, v_i\}\)--\(\{v_{i+1}, \ldots, v_i\}\) linkage, so the invariant is satisfied, and we can safely replace \(I_j\) with \(I'\) and \(I''\), and \(G_j\) with \(G'\) and \(G''\).

    Let us analyse the time complexity of the described algorithm.
    Consider the sequence \((G_1, \ldots, G_k)\) at any moment of the execution.
    We refer to the graphs \(G_j\) as the \emph{parts} of \(G\).
    Observe that if the sequence of parts was \((G_1', \ldots, G_{k'}')\) at an earlier stage of the execution, then each new part \(G_j\) is a subgraph of some old part \(G_{j'}'\). After an execution of Menger's Theorem, the part \(G'_{j'}\) splits into two parts whose corresponding segments are roughly two times smaller length, and \(G_j\) is contained in one of these parts.
    Hence, during the whole execution of the algorithm,  Menger's Theorem is applied to a part containing \(G_j\) only \(\mathbf{O}(\log |V(G)|)\) times. Since the time complexity of running Menger's Theorem on a part \(G_{j'}'\) is \(\mathbf{O}(p|E(G_{j'}')|)\), the total time complexity is \(\mathbf{O}(p|E(G)|\log |V(G)|)\).
\end{proof}

The following is an easy consequence of this lemma and the fact that a society with a linear decomposition of adhesion at most $p$ has depth at most $2p$, which we already noted in \Cref{subsec:vortex}.

\begin{corollary}\label{lem:findtransactionefficient}
    Let $p$ be a positive integer and let $(G,\Omega)$ be a society with depth at least $2p$.
    Then we can find a transaction of order $p$ in $\mathbf{O}(p|E(G)| \log |V(G)|)$-time.
\end{corollary}

\subsection{Constructively finding crooked transactions}
The proof of \Cref{lem:crookedexistence} that Kawarabayashi et al.\ present in \cite{KawarabayashiTW2021Quickly} can be split into three major components.
To describe them, we first need a technical definition that describes a planar, flat, isolated, separating transaction in reference to which one ``side'' of the society is flat.

\begin{definition}[Transactions that are flat on one side]\label{def:flatononeside}
    Let $p$ be an integer with $p \geq 2$, let $(G,\Omega)$ be a society, and let $\mathcal{P} = \{ P_1, \ldots , P_p \}$ be a planar, flat, isolated, separating transaction in $(G,\Omega)$ that is indexed naturally, with $\rho$ being a vortex-free rendition of the $\mathcal{P}$-strip society of $(G,\Omega)$.
    We let $I$, respectively $I'$, be the unique segment of $\Omega$ that has the same endpoints as $P_{p-1}$, respectively $P_2$, and does not contain the endpoints of $P_p$, respectively $P_1$.
    
    Let $H_1$ be the graph consisting of $\mathcal{P} \setminus \{ P_p \}$ and the vertices in $I$, and let $G_1$ be the graph consisting of the union of $H_1$ and all $H_1$-bridges with at least one attachment in $V(H_1) \setminus V(P_{p-1})$.\footnote{Due to $\mathcal{P}$ being flat and separating, no such bridge can contain a vertex of $V(P_p)$ or $V(\Omega) \setminus I$.}
    The society $(G_1, \Omega_1)$ is defined by choosing $\Omega_1$ such that $V(\Omega_1) = I \cup (N(\rho) \cap V(P_{p-1}))$, with the order on $I$ being obtained naturally from $\Omega$ and the order on $N(\rho) \cap V(P_{p-1})$ being derived from a traversal of $P_{p-1}$.
    We define $H_p$, $G_p$, and $\Omega_p$ analogously based on the $I'$ and the paths $\mathcal{P} \setminus \{ P_1 \}$, with $P_2$ playing the role of $P_{p-1}$.

    We say that $\mathcal{P}$ is \emph{flat on one side} if $(G_1, \Omega_1)$ or $(G_p, \Omega_p)$ has a vortex-free rendition in a disk.
\end{definition}

The three components of the proof of \Cref{lem:crookedexistence} in \cite{KawarabayashiTW2021Quickly} can now be described as follows:
\begin{enumerate}
    \item Find a large flat, planar transaction that is flat on one side in the absence of a large crooked transaction.

    \item Iteratively build a well-structured cylindrical rendition whenever we fail to find a large crooked transaction.

    \item Reroute a crooked transaction onto the original society using the well-structured cylindrical rendition built in the second step.
\end{enumerate}
Both the second and third point are taken care of by using minimality arguments in \cite{KawarabayashiTW2021Quickly}.
We will demonstrate constructive procedures with explicit runtimes for both of these steps.

Despite the somewhat complicated definition the first component entails, the statements directly associated with it are simple to prove.
We will directly follow the arguments presented in \cite{KawarabayashiTW2021Quickly}.
The bounds for our version of \Cref{lem:crookedexistence} will be slightly worse, since splitting the proof of \Cref{lem:crookedexistence} into several statements prevents us from reusing certain renditions with convenient properties.
Our first lemma shows that any large planar transaction yields a large crooked transaction or a large planar, flat transaction that is flat on one side.
The specificity in the first point of the statement will be useful later on in \Cref{sec:multisocietycrooked}.

\begin{proposition}[Kawarabayashi, Thomas, and Wollan \cite{KawarabayashiTW2021Quickly}]\label{lem:planartransactionflatting}
    Let $p,q$ be positive integers with $p,q \geq 2$, let $(G,\Omega)$ be a society, and let $\mathcal{P}$ be a planar transaction in $(G,\Omega)$ of order $3q + 2p - 3$.
    Then in $(G,\Omega)$ there exists
    \begin{enumerate}
        \item a crooked transaction $\mathcal{Q}$ of order $q$ such that
        \begin{itemize}
             \item $\mathcal{Q}$ is a leap with the overpass $Q$ and $\mathcal{Q} \setminus \{ Q \} \subseteq \mathcal{P}$, or 

             \item $\mathcal{Q}$ is a doublecross containing two crosses $P,P'$ and $R,R'$, and $\mathcal{Q} \setminus \{ P,P',R,R' \} \subseteq \mathcal{P}$, or
        \end{itemize}

        \item a planar, flat, isolated, separating transaction $\mathcal{P}' \subseteq \mathcal{P}$ of order $p$ that is flat on one side.
    \end{enumerate}
    Furthermore, we can find one of these two transactions in time $\mathbf{O}(pq|E(G)|)$.
\end{proposition}
\begin{proof}
    Let $p' \coloneqq 3q+2p-3$ and let $\mathcal{P} = \{ P_1, \ldots , P_{p'} \}$ be indexed naturally.
    Let $A,B$ be the end segments of $\mathcal{P}$, let $a_i \in A$ be an endpoint of $P_i$, and let $b_i \in B$ be an endpoint of $P_i$ for each $i \in [p']$.

    We let $H$ be the subgraph of $G$ obtained from the union of the elements of $\mathcal{P}$ and $V(\Omega)$.
    Further, let $G_1'$ be the subgraph of $H$ consisting of $\bigcup_{i=1}^{q+p-1} P_i$ and all vertices in the segment $S_1$ of $\Omega$ that ends on $a_{q+p-1}$ and $b_{q+p-1}$ and contains the endpoints of $P_1$.
    Consider all $H$-bridges with at least one attachment in $V(G_1') \setminus V(P_{q+p-1})$, and for each such bridge $B$, let $B'$ denote the graph obtained from $B$ by deleting all attachments that are not contained in $V(G_1')$.
    Finally, let $G_1$ denote the union of $G_1'$ and all graphs $B'$ as above.
    The cyclic permutation $\Omega_1$ is then defined by saying that $V(\Omega_1) = S_1$, with the order of $\Omega_1$ being obtained naturally from $\Omega$.
    Thus $(G_1,\Omega_1)$ is a society and we define the society $(G_2, \Omega_2)$ analogously using the paths $P_{2q+p-1}, \ldots , P_{p'}$ instead, with $P_{2q+p-1}$ playing the role of $P_{q+p-1}$.

    If $G_1$ and $G_2$ are not disjoint then $G_1 \cup G_2$ includes a path $P$ from $S_1$ to the segment of $\Omega$ with the endpoints $a_{2q+p-1}, b_{2q+p-1}$ that contains the endpoints of $P_{p'}$.
    The path $P$ must be disjoint from $P_{q+p}, \ldots , P_{2q+p-2}$, which means that $P$ together with these paths forms a crooked transaction of order $q$, as we desire.
    Thus we can assume that $G_1$ and $G_2$ are disjoint.

    Should both societies $(G_1,\Omega_1)$ and $(G_2, \Omega_2)$ contain a cross, then the union of the two crosses and $P_{q+p}, \ldots , P_{2q+p-2}$ form a crooked transaction of order more than $q$, again satisfying the first item of our statement.
    We can therefore assume without loss of generality that $(G_1,\Omega_1)$ has a vortex-free rendition $\rho_1'$ in a disk $\Delta_1'$.

    For each $i \in [p]$ we let $T_i$ be the trace of $P_i$ in $\rho_1'$ and let $D_i$ be the disk in $\Delta_1' \setminus T_i$ that does not contain $\pi^{-1}(a_{p+1})$.
    We let $J_i$ be the crop of $G_1$ to $D_i$.
    Since $\mathcal{L}'$ is indexed naturally, we have $J_i \subseteq J_j$ whenever $i,j \in [p]$ and $i < j$.

    Suppose that there exists some $i \in [p]$ and an edge in $G$ with one endpoint $u$ in $J_i - \pi(N(\rho_1') \cap T_i)$ and the other endpoint $v$ in $G - V(J_i)$.
    Due to $T_i$ being a trace in $\rho_1'$, we conclude that $v \not\in V(G_1)$.
    Thus there exists a path $Q_1$ in $J_i$ connecting $u$ and a vertex in $V(\Omega)$ and a path $Q_2$ connecting $v$ and a vertex of $V(\Omega)$ in $G - V(G_1)$.
    Since $Q_1$ and $Q_2$ must be disjoint, $Q_1 \cup Q_2$ together with $uv$ make up a path that forms a crooked transaction of order $q$ together with the paths $P_{p+1}, \ldots , P_{q+p-1}$.
    Thus, we may suppose that this is not the case for any $i \in [p]$.

    Let $\mathcal{Q} = \{ P_1, \ldots , P_p \}$.
    Using $\rho_1'$, we can construct a vortex-free rendition $\rho''$ of the $\mathcal{Q}$-strip society of $(G,\Omega)$ that is a restriction of $\rho_1'$.
    In particular, our arguments imply that $\mathcal{Q}$ is a planar, flat, isolated, separating transaction of order $p$ in $(G,\Omega)$.
    Furthermore, $\mathcal{Q}$ is flat on one side, deriving from the fact that $\rho_1'$ already yields a vortex-free rendition of $(G_1,\Omega_1)$ and our above observations on the graphs $J_1, \ldots , J_p$.
    This concludes our proof.
\end{proof}

This allows us to show that in a society with large depth we can find a large crooked transaction or a large flat transaction that is flat on one side.
We could give a $6q+2$ bound on the depth of the society if we were to apply \Cref{lem:findtransactioninefficient}, which is slightly worse than the $6q$ bound given in \cite{KawarabayashiTW2021Quickly}.
However, we will be using \Cref{lem:findtransactionefficient} and thus incur another factor 2 in this bound.

\begin{proposition}[Kawarabayashi, Thomas, and Wollan \cite{KawarabayashiTW2021Quickly}]\label{lem:planarandflatononesideorcrooked}
    Let $q$ be an integer with $q \geq 2$ and let $(G,\Omega)$ be a society with depth at least $12q+4$.
    Then in $(G,\Omega)$ there exists
    \begin{enumerate}
        \item a crooked transaction of order $q$, or

        \item a planar, flat, isolated, separating transaction of order $q+4$ that is flat on one side.
    \end{enumerate}
    Furthermore, we can find one of these two transactions in time $\mathcal{O}(q|E(G)|(q+ \log |V(G)|))$.
\end{proposition}
\begin{proof}
    Since $(G,\Omega)$ has depth at least $12q+4$, we can find a transaction $\mathcal{P}$ in $(G,\Omega)$ of order $6q+2$ using \Cref{lem:findtransactionefficient}.
    Furthermore, using \Cref{lem:crookedorplanar}, we can either find a crooked transaction of order $q$ within $\mathcal{P}$, which would mean we are done, or there exists a planar transaction $\mathcal{P}' \subseteq \mathcal{P}$ of order $5q +5$.
    Applying \Cref{lem:planartransactionflatting} to $\mathcal{P}'$ now yields one of the desired outcomes.
\end{proof}

The second component of the proof of \Cref{lem:crookedexistence} in \cite{KawarabayashiTW2021Quickly} requires us to introduce a new data structure to track the iterative progress we make towards embedding the society.

\begin{figure}[ht]
    \centering
    \begin{tikzpicture}

        \pgfdeclarelayer{background}
		\pgfdeclarelayer{foreground}
			
		\pgfsetlayers{background,main,foreground}

        \begin{pgfonlayer}{background}
        \node (C) [v:ghost] {{\includegraphics[width=12cm]{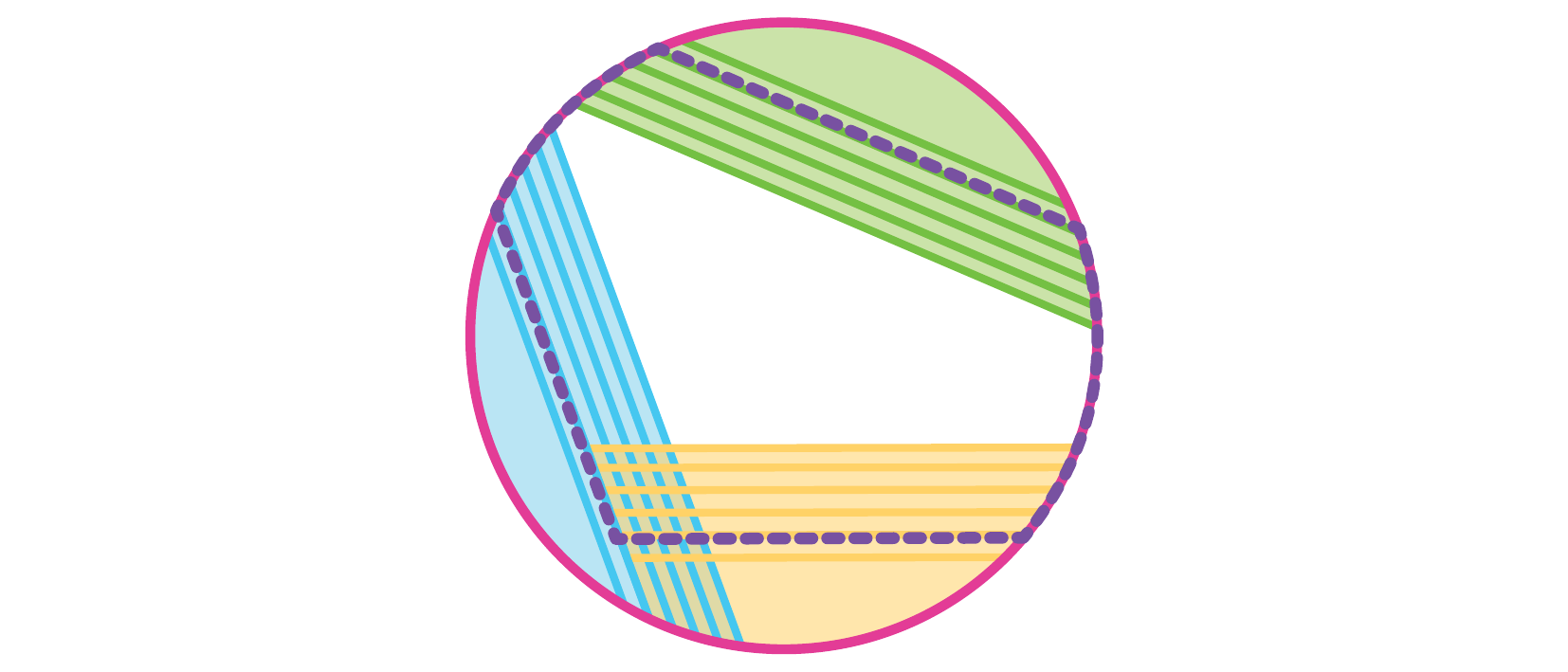}}};
        \end{pgfonlayer}{background}
			
        \begin{pgfonlayer}{main}

            \node (X) [v:ghost,position=0:0mm from C] {};

            \node (c0) [v:ghost,position=0:0mm from X] {$c_3$};

            \node (P11) [v:ghost,position=242:27.4mm from X] {$P^1_1$};
            \node (P21) [v:ghost,position=320:27mm from X] {$P^2_1$};
            \node (P31) [v:ghost,position=110:27mm from X] {$P^3_1$};

            \node (G3) [v:ghost,position=350:32mm from X] {$(G_3,\Omega_3)$};
        
        \end{pgfonlayer}{main}

        \begin{pgfonlayer}{foreground}
        \end{pgfonlayer}{foreground}

    \end{tikzpicture}
    \caption{A $(6,3)$-constriction of a society $(G,\Omega)$. Each of the three linkages $\mathcal{P}_1$ (depicted in \textcolor{CornflowerBlue}{blue}), $\mathcal{P}_2$ (depicted in \textcolor{BananaYellow}{yellow}), and $\mathcal{P}_3$ (depicted in \textcolor{AppleGreen}{green}) is a transaction on its corresponding restricted society that is flat on one side (indicated by the coloured areas).
    The dashed line indicates the society $(G_3,\Omega_3)$.
    Notice that $(G_3,\Omega_3)$ has a cylindrical rendition where the coloured parts are all disjoint from the unique vortex $c_3$.}
    \label{fig:Constriction}
\end{figure}
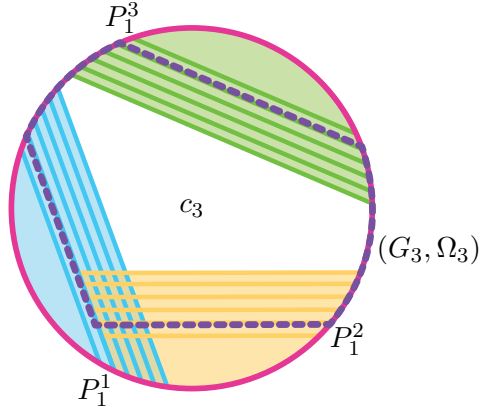

\begin{definition}[Constriction of a society]\label{def:constriction}
    Let $q,\ell$ be non-negative integers with $q \geq 3$, let $(G_0,\Omega_0)$ be a society, and let $\rho_0$ be the cylindrical rendition in a disk $\Delta_0$ with $\Delta_0$ its unique vortex and cell.
    We call $\mathfrak{N} = (G_0,\Omega_0,\mathfrak{P})$ a \emph{$(q,\ell)$-constriction of $(G_0,\Omega_0)$} if
    \begin{description}
        \item[~C1~] the set $\mathfrak{P} = \{ \mathcal{P}_1, \ldots , \mathcal{P}_\ell \}$ is a set of linkages in $G_0$ such that $\mathcal{P}_i = \{ P_1^i, \ldots , P_q^i \}$ for each $i \in [\ell]$,

        \item[~C2~] there exist societies $(G_1,\Omega_1), \ldots , (G_\ell,\Omega_\ell)$ such that for each $i \in [\ell]$, we have $G_i \subseteq G_{i-1} - V(P_1^i)$ and $\mathcal{P}_i$ is a planar, flat, isolated, separating transaction in $(G_{i-1},\Omega_{i-1})$ that is indexed naturally and flat on one side,\footnote{In the context of this definition the flat side will be the one that contains $P_1^i$.}

        \item[~C3~] for all $i \in [\ell]$, there exists a cylindrical rendition $\rho_i$ of $(G_i,\Omega_i)$ in a $\rho_{i-1}$-aligned disk $\Delta_i \subseteq \Delta_{i-1}$ around the vortex $c_i$, with $\{ P_3^i, \ldots , P_{q-1}^i \}$ being an unexposed transaction in $(G_i,\Omega_i)$ in $\rho_i$, and

        \item[~C4~] for all $i \in [\ell]$, the $\Delta_i$-society in $\rho_{i-1}$ is $(G_i, \Omega_i)$, we have $V(\Omega_i) \subseteq V(\Omega_{i-1}) \cup V(P_2^i)$, and $\widetilde{c}_i \subseteq \widetilde{c}_{i-1} \cup V(P_q^i)$.
    \end{description}
    See \cref{fig:Constriction} for an illustration.
    The society $(G_\ell, \Omega_\ell)$ is called the \emph{frontier of $\mathfrak{N}$}, with the \emph{frontier-rendition $\rho_\ell$ in the disk $\Delta_\ell$}.
    In particular, according to our definition there exists a cylindrical rendition $\rho$ of $(G_0,\Omega_0)$ in $\Delta_0$ around $c_i$ such that $\rho_\ell$ is a restriction of $\rho$ to $\Delta_\ell$.
    We call $\rho$ the \emph{total rendition for $\mathfrak{N}$}.
    Since we allow for $\ell = 0$, each society has a trivial $(q,0)$-constriction.
\end{definition}

We now show that we can always find a crooked transaction in the frontier, verify that the entire society has a cylindrical rendition of bounded depth, or we can extend a given constriction.

\begin{lemma}\label{lem:extendconstriction}
    Let $q,\ell$ be non-negative integers with $2 \leq q$.
    Let $(G,\Omega)$ be a society with a $(q+3,\ell)$-constriction $\mathfrak{N} = (G,\Omega,\mathfrak{P})$ with the frontier $(H,\Omega')$.
    Then either
    \begin{enumerate}
        \item there exists a crooked transaction of order at least $q$ in $(H,\Omega')$,

        \item there exists a cylindrical rendition of $(G,\Omega)$ of depth at most $12q+4$, or

        \item there exists a transaction $\mathcal{Q}$ in $(H,\Omega')$ such that $(G,\Omega,\mathfrak{P} \cup \{ \mathcal{Q} \})$ is a $(q+3,\ell+1)$-constriction.
    \end{enumerate}
    Moreover, there exists an algorithm that leads to one of these options in time $\mathbf{O}(q|E(G)| \log |V(G)|)$.
\end{lemma}
\begin{proof}
    Let $\rho = (\Gamma, \mathcal{D})$ be the total rendition for $\mathfrak{N}$ in a disk $\Delta$ and let $\rho_H$ be the frontier rendition of $(H,\Omega')$ in a disk $\Delta_H \subseteq \Delta$.
    If $(H,\Omega')$ has depth at most $12q+4$, then $(\Gamma, (\mathcal{D} \cup \{ \Delta_H \} ) \setminus \{ c \in \mathcal{D} ~\!\colon\!~ c \subseteq \Delta_H \})$ is a cylindrical rendition of $(G,\Omega)$ of depth at most $12q+4$.
    
    Thus we may suppose that $(H,\Omega')$ has depth at least $12q+4$.
    By applying \Cref{lem:planarandflatononesideorcrooked} to $(H,\Omega')$ we therefore either find a crooked transaction of order $q$ in $(H,\Omega')$, or we find a planar, flat, isolated separating transaction $\mathcal{P} = \{ P_1, \ldots , P_{q+4} \}$ of order $q+4$ in $(H,\Omega')$ that is flat on one side.
    Since the first option is one of our goals, we may suppose that we find $\mathcal{P}$, which we choose to be indexed naturally such that $P_1$ is contained on the side of $\mathcal{P}$ that has a vortex-free rendition.

    Let $I$ be the segment of $\Omega$ with the same endpoints as $P_{q+4}$ that contains all other endpoints of the paths in $\mathcal{P}$
    Further, let $H^*$ be the subgraph of $G$ consisting of the paths in $\mathcal{P}$ and the vertices in $I$.
    We let $G'$ be the union of $H^*$ and all $H^*$-bridges with at least one attachment in $H^* - V(P_{q+4})$ and let the cyclic ordering $\Omega''$ of $V(\Omega'') = I$ be defined naturally via the order of $\Omega$ on $I$.
    Since $\mathcal{P}$ is flat, isolated, separating, and flat on the side that contains $P_1$, the society $(G', \Omega'')$ has a vortex-free rendition $\rho'$ in a disk $\Delta' \subseteq \Delta_H$.

    We let $T_{q+3}$ be the trace of $P_{q+3}$ in $\rho'$ and let $\Delta_1$ be the closure of the disk in $\Delta' - T_{q+3}$ that contains the trace of $P_1$, with $G_1'$ being the crop of $G'$ to $\Delta_1$.
    Let $I'$ be the unique segment of $\Omega''$ that has the same endpoints as $P_{q+3}$ and contains both endpoints of $P_{q+4}$.
    This allows us to define a cylindrical rendition $\rho''$ of $(H, \Omega')$ as follows.
    We embed $G_1'$ as in $\rho'$ -- which does not involve any vortices, as $\rho'$ is vortex-free -- and introduce a single vortex $c'$ that contains the entirety of $H[(V(H) \setminus V(G_1)) \cup (N(\rho') \cap V(P_{q+3}))] - E(G_1')$ drawn onto $c'$ such that exactly the vertices in $N(\rho') \cap V(P_{q+3})$ and $I'$ appear on its boundary in a manner fitting $\Omega''$ and the order of vertices on $P_{q+3}$.
    Within $\rho''$, we let $T_2$ be the trace of $P_2$ and let $\Delta_{\ell+1}$ be the closure of the disk in $\Delta_H - T_2$ that contains the trace of $P_3$, with $(H',\Psi)$ being the $\Delta_{\ell+1}$-society in $\rho''$.

    The cylindrical rendition $\rho''$, the society $(H',\Psi)$, and the transaction $\mathcal{Q} = \{ P_1, \ldots , P_{q+3} \}$ in $(H,\Omega')$ together confirm that $(G,\Omega,\mathfrak{P} \cup \{ \mathcal{Q} \})$ is a $(q+3,\ell+1)$-constriction, which completes our proof.
\end{proof}

A simple induction using repeated applications of the above lemma yields the following corollary.
The worsening of the runtime as compared to \Cref{lem:extendconstriction} is due to our inability to bound the number of times we have to refine the constriction.

\begin{corollary}\label{cor:constrictionorrendition}
    Let $q,\ell$ be non-negative integers with $2 \leq q$ and let $(G,\Omega)$ be a society.
    Then
    \begin{enumerate}
        \item there exists a cylindrical rendition $\rho$ of $(G,\Omega)$ of depth at most $12q+4$, or
        
        \item there exists a $(q+3,\ell)$-constriction $\mathfrak{N}$ of $(G,\Omega)$ with the frontier $(H,\Omega')$ and a crooked transaction $\mathcal{Q}$ of order at least $q$ in $(H,\Omega')$.
    \end{enumerate}
    Moreover, there exists an algorithm that finds $\rho$ or $\mathfrak{N}$ and $\mathcal{Q}$ in time $\mathbf{O}(q|E(G)||V(G)| \log |V(G)|)$.
\end{corollary}

The first option of this corollary is of course something we would welcome in the context of \Cref{lem:crookedexistence}.
Thus our last task will be move the orthogonal transaction we find in the frontier of the constriction back onto the original society.
For this purpose Kawarabayashi et al.\ use a pair of so-called exchange lemmas for crooked transactions from \cite{RobertsonS1990Graph}.
We will ultimately only need one of the two, as the first -- Lemma 3.4 in \cite{KawarabayashiTW2021Quickly} and (3.2) in \cite{RobertsonS1990Graph} -- is only used to prove the second lemma.
The lemma that we will need is Lemma 3.5 in \cite{KawarabayashiTW2021Quickly} and (3.3) in \cite{RobertsonS1990Graph}.
As we do not need the techniques contained in its proof, we only give its statement.

\begin{proposition}[Robertson and Seymour \cite{RobertsonS1990Graph}]\label{lem:crookedexchange}
    Let $(G,\Omega)$ be a society and let $\mathcal{P}$ be a crooked transaction in $(G,\Omega)$.
    Let $Q_1$ and $Q_2$ be disjoint paths, each with one endpoint in $V(\Omega)$, the other endpoint in $V(\mathcal{P})$, and no other vertex in $V(\mathcal{P})$.
    There exists a path $P \in \mathcal{P}$ and subpaths $R_1, R_2$ of $P$ such that one of the following is a crooked transaction:
    \begin{enumerate}
        \item $(\mathcal{P} \setminus \{ P \}) \cup \{ R_1 \cup Q_1 \}$,

        \item $(\mathcal{P} \setminus \{ P \}) \cup \{ R_2 \cup Q_2 \}$, or

        \item $(\mathcal{P} \setminus \{ P \}) \cup \{ R \cup Q_1, R_2 \cup Q_2 \}$.
    \end{enumerate}
\end{proposition}

We will also need the following simple lemma about shrinking crooked transactions.

\begin{proposition}[Kawarabayashi, Thomas, and Wollan \cite{KawarabayashiTW2021Quickly}]\label{lem:crookedshrink}
    Let $(G,\Omega)$ be a society and let $\mathcal{P}$ be a crooked transaction of order at least five in $(G,\Omega)$.
    Then there exists a path $P \in \mathcal{P}$ such that $\mathcal{P} \setminus \{ P \}$ is a crooked transaction.
\end{proposition}

We are now ready to prove the last of the major steps mentioned earlier.
The arguments presented here are essentially taken directly from the proof of \Cref{lem:crookedexistence} in \cite{KawarabayashiTW2021Quickly}, though we added some observations on how the arguments affect the runtime of the process.

\begin{lemma}\label{lem:crookedrouteback}
    Let $q,\ell$ be integers with $q \geq 4$ and $\ell \geq 1$.
    Let $(G,\Omega)$ be a society with a $(q+3,\ell)$-constriction $\mathfrak{N} = (G,\Omega,\mathfrak{P} = \{ \mathcal{P}_1, \ldots , \mathcal{P}_\ell \})$ and a crooked transaction of order $q$ in the frontier $(H,\Psi)$ of $\mathfrak{N}$.
    
    Then there exists a crooked transaction $\mathcal{Q}'$ of order $q$ in the frontier $(H',\Psi')$ of the $(q+3,\ell-1)$-constriction $\mathfrak{N}' = (G,\Omega,\mathfrak{P} \setminus \{ \mathcal{P}_\ell \})$ such that the endpoints of $\mathcal{Q}'$ are found in $V(\Psi) \cap V(\Psi')$.
    
    In particular, there exists an algorithm that finds $\mathcal{Q}'$ in time $\mathbf{O}(q|E(G)||V(G)|)$.
\end{lemma}
\begin{proof}
    Let $\mathcal{Q}$ be the crooked transaction of order $q$ in $(H,\Psi)$.
    Furthermore, we let $\mathcal{P}_\ell = \{ P_1, \ldots , P_{q+3} \}$ be indexed naturally such that $P_1 \not\subseteq H$.
    We also let $\rho$ be the frontier rendition of $(H,\Psi)$ around the vortex $c$ in a disk $\Delta$.
    
    Recall that according to \Cref{def:constriction} the paths $P_3, P_4, \ldots , P_{q+2}$ are unexposed in $\rho$, vertices of $P_2$ form part of the set $V(\Psi)$, and vertices of $P_{q+3}$ form part of the boundary of $c$.
    Accordingly, we have $V(\Psi) \setminus V(\Psi') \subseteq V(P_2)$ and in particular, the endpoints of $P_2$ are both found in $V(\Psi)$.
    Let $I$ be the segment of $\Psi$ that has the same endpoints as $P_2$ and only contains vertices also found in $V(P_2)$.
    If all endpoints of $\mathcal{Q}$ are found outside of $I$, there is nothing to prove.
    Thus we may assume from this point forward, that some path in $\mathcal{Q}$ has an endpoint in $I$.

    Our goal from now on will be to iteratively minimise the size of $|E(\mathcal{Q}) \cup E(\{ P_3, P_4, \ldots , P_{q+3} \}|$, which is originally bounded by $|E(H)| \leq |E(G)|$, and show that any such linkage cannot have any endpoints in $I$.
    For later use, let $\mathcal{P}' = \{ P_3, P_4, \ldots , P_{q+3} \}$.
    
    Towards the goal of minimising the size of the set defined above, we will assume that the set is already of minimal size and point out how costly the construction of an improvement -- which will generally be achieved by finding a new crooked transaction that shares more edges with the paths in $\mathcal{P}'$ -- would be whenever we encounter a contradiction to the minimality of this quantity.
    
    Let us first confirm the following claim.

    \begin{claim}\label{claim:nounexposed}
        No path in $\mathcal{Q}$ is unexposed in $\rho$.
    \end{claim}
    \emph{Proof of \Cref{claim:nounexposed}:}
    As $\rho$ is a cylindrical rendition with the vortex $c$, any such path would in particular be grounded and thus have a trace $T$.
    This defines two disks within $\Delta - T$, one of which is disjoint from $c$.
    Let this disk be $d$.
    We can choose our unexposed path within $\mathcal{Q}$ to minimise $d$ (by containment within the other disks).
    This path must be a peripheral path of $\mathcal{Q}$, contradicting the fact that $\mathcal{Q}$ is crooked.
    \hfill$\blacksquare$

    This allows us to now easily deduce the following useful property.

    \begin{claim}\label{claim:meetallpaths}
        If a path $Q \in \mathcal{Q}$ has an endpoint in $I$, then it intersects each path within $\mathcal{P}'$.
        Furthermore, for each $i \in [3,q+3]$, at least two distinct paths in $\mathcal{Q}$ intersect $P_i$.
    \end{claim}
    \emph{Proof of \Cref{claim:meetallpaths}:}
    The first part of the claim is a consequence of \Cref{claim:nounexposed}, as the endpoint of $Q$ in $I$ and $\rho$ being a cylindrical rendition around $c$ otherwise forces $Q$ to be unexposed, leading to a contradiction.

    Towards the second part of the claim, choose some $i \in [3,q+3]$ such that $P_i$ is only intersected by a unique path $Q \in \mathcal{Q}$.
    This $P_i$ must exist due to the first part of our claim and our assumption that some path in $\mathcal{Q}$ has an endpoint $u$ in $I$.
    Let $R$ be a subpath of $P_i$ with one endpoint $v$ in $V(\Psi)$ and the other $w$ in $V(Q)$ such that $R$ is internally disjoint from $Q$.
    Since $P_i$ is disjoint from all paths in $\mathcal{Q} \setminus \{ Q \}$ and the path $uQwP_iv$ is unexposed, \Cref{claim:nounexposed} implies that no path in $\mathcal{Q}$ has an endpoint in the segment $I'$ of $\Psi$ that ends on $u$ and $v$ but does not contain the other endpoint $u'$ of $Q$.
    Therefore $(\mathcal{Q} \setminus \{ Q \}) \cup \{ u'QwP_iv \}$ is a crooked transaction in $(H,\Psi)$, which we could identify and construct in time $\mathbf{O}(q|E(G)|)$.
    Since at least one edge of $uQw$ is not contained in this new transaction, this contradicts the minimality condition on $\mathcal{Q}$, proving our claim.
    \hfill$\blacksquare$

    As a consequence of \Cref{claim:meetallpaths}, we now know that $V(\mathcal{Q}) \cap V(\mathcal{P}') \neq \emptyset$.

    For our next arguments, we call an endpoint $x$ of $P_i \in \mathcal{P}'$ \emph{available} if no path in $\mathcal{Q}$ has $x$ as an endpoint.
    Since $\mathcal{Q}$ has order $q$, there exist at least two available vertices.
    Let $x$ and $y$ be two such available vertices, where $x \in V(P_i)$ and $y \in V(P_j)$ for $i,j \in [3,q+3]$ (possibly with $i = j$).
    Further, let $R_x$ be a minimal subpath of $P_i$ ending in $x$ and a vertex of $V(\mathcal{Q})$ and let $R_y \subseteq P_j$ be defined analogously.

    Suppose $R_x$ and $R_y$ intersect.
    Then $i = j$ must hold and in particular $R_x$ and $R_y$ intersect in a single vertex.
    This however implies that a unique path from $\mathcal{Q}$ is intersecting $P_i$, which contradicts \Cref{claim:meetallpaths}.
    Thus $R_x$ and $R_y$ must be disjoint.

    We apply \Cref{lem:crookedexchange} to $\mathcal{Q}$, $R_x$, and $R_y$, which tells us that there exists a $Q'$ and two subpaths $Q_1, Q_2 \subseteq Q'$ such that one of the three outcomes of \Cref{lem:crookedexchange} holds for $Q'$, $Q_1$, $Q_2$, $R_x$, and $R_y$.
    Note that, given $\mathcal{Q}$, $R_x$, and $R_y$, we can determine which of the three outcomes of \Cref{lem:crookedexchange} holds in time $\mathbf{O}(q|E(G)|)$.

    Suppose that one of the first two outcomes holds and suppose further w.l.o.g.\ that $(\mathcal{Q} \cup \{ Q' \}) \cup \{ R_x \cup Q_1 \}$ is a crooked transaction.
    Due to the minimality assumptions on $\mathcal{Q}$, our new crooked transaction must contain the same edges as $\mathcal{Q}$.
    Thus, every edge of the path $Q' - Q_1$ is contained in $E(\bigcup_{P \in \mathcal{P}'} P)$.
    As a consequence $V(P_i) = V(R_x) \cup V(Q' - Q_1)$, which implies that $Q'$ is the unique path in $\mathcal{Q}$ that intersects $P_i$, contradicting \Cref{claim:meetallpaths}.

    We may therefore move on to the third option of \Cref{lem:crookedexchange}.
    Let $Q_3$ be the minimal subpath of $Q'$ with its endpoints in $V(R_x)$ and $V(R_y)$ such that $Q' = Q_1 \cup Q_2 \cup Q_3$.
    If $Q_3$ has at least one edge outside of $E(\bigcup_{P \in \mathcal{P}'} P)$ then our new crooked transaction contradicts the minimality of $\mathcal{Q}$, thanks to \Cref{lem:crookedshrink}.

    This means that $E(Q_3)$ is a proper subset of $E(\bigcup_{P \in \mathcal{P}'} P)$.
    Since $\mathcal{P}'$ is a linkage, $Q_3$ must therefore be contained in a single path found in $\mathcal{P}'$ and this can only be $P_i$ with $i = j$.
    This however again implies that $Q'$ is the only path intersecting $P_i$, contradicting \Cref{claim:meetallpaths} once more and completing our proof.
\end{proof}

We can employ induction together with this new lemma to derive a very useful corollary.

\begin{corollary}\label{cor:crookedrouteback}
    Let $q,\ell$ be positive integers with $4 \leq q$.
    Let $(G,\Omega)$ be a society with a $(q+3,\ell)$-constriction $\mathfrak{N}$ and a crooked transaction of order $q$ in the frontier of $\mathfrak{N}$.
    
    Then there exists a crooked transaction of order $q$ in $(G,\Omega)$.
    
    In particular, there exists an algorithm that finds this crooked transaction in time $\mathbf{O}(q|E(G)||V(G)|^2)$.
\end{corollary}

This allows us to easily derive the desired algorithmic version of \Cref{lem:crookedexistence}.
Note that using \Cref{lem:findtransactioninefficient} throughout our arguments in this section would immediately yield a version with a tighter bound on the depth and a slightly worse runtime.

\begin{lemma}\label{lem:crookedexistencealgo}
    Let $q$ be an integer with $q \geq 4$ and let $(G,\Omega)$ be a society.
    There exists a crooked transaction of order $q$ in $(G,\Omega)$ or $(G,\Omega)$ has a cylindrical rendition of depth at most $12q+4$.

    In particular, there exists an algorithm that returns either the crooked transaction or the cylindrical rendition in time $\mathbf{O}(q|E(G)||V(G)|^2)$.
\end{lemma}
\begin{proof}
    As noted in \Cref{def:constriction}, the society $(G,\Omega)$ itself trivially describes a $(q+3,0)$-constriction.
    This allows us to apply \Cref{cor:constrictionorrendition}.
    If this yields a cylindrical rendition, we are done.

    Otherwise, we find a crooked transaction of order $q$ in the frontier of a $(q+3,\ell)$-constriction of $(G,\Omega)$ for some non-negative integer $\ell$.
    If $\ell = 0$, we are done and otherwise \Cref{cor:crookedrouteback} yields the desired crooked transaction of order $q$ in $(G,\Omega)$.
\end{proof}

\section{Orthogonalising crooked transactions}\label{sec:orthogonalisecrooked}
As demonstrated in the last few sections, at the core of many of our proofs is the routing of large linkages through cycles, other linkages, and structures that more or less look like a grid.
This section is devoted to doing this in a more abstracted way than we previously had to.
In particular, we want to find a large orthogonal crooked transaction in societies with a cylindrical society and a large nest and often we also want to manipulate radial linkages associated with these objects.
Note that for these purposes \Cref{lem:crookedexistence} alone is not a strong enough statement, since it does not respect existing nests and renditions.

The results of this section are the technical underpinnings of the more complex structures and proofs we will present in the remaining sections.
The running times of the implicit algorithms we present here are central to the total running time of our algorithm for finding both the global decomposition and the witnesses for the local structure theorem.
This is a consequence of the fact that the lemmata and algorithms we present here will be used in the context of several iterative processes found in later arguments.

\subsection{Orthogonalising radial linkages}
The first type of object we will orthogonalise are radial linkages, which builds the basis for our later efforts in this section.
To illustrate why this is of interesting to us, we first present a useful intermediate lemma.
At several points towards establishing the Local Structure Theorem we will have to connect linkages that have been extracted from deeper inside the structure to their respective radial linkages while making sure that enough of this radial linkage still reaches all the way ``down'' to the vortices.
The following lemma is a consequence of \hyperref[prop:mengersthm]{Menger's Theorem} and fits exactly this purpose.

\begin{lemma}\label{lemma:connected_linkages}
Let $s,r,k,\ell$ be positive integers with $s \geq r+3$.
Let $(G,\Omega)$ be a society with a cylindrical rendition $\rho$ and a nest $\mathcal{C} = \{ C_1, \ldots , C_s \}$ around the vortex $c_0$.
Moreover, let $\mathcal{L}$ and $\mathcal{R}$ each be radial linkages of order $r$ in $(G,\Omega)$ such that both are orthogonal to $\mathcal{C}$ and let $I = [\ell,k] \subseteq [2,s]$ be an interval with $|I| = r+2$.

Then there exists a radial linkage $\mathcal{P}$ of order $r$ in $(G,\Omega)$ such that
\begin{enumerate}
    \item $\mathcal{P}$ is orthogonal to $\{ C_i ~\!\colon\!~ i \in [s] \setminus I \}$ with endpoints on $C_1$,
    
    \item $H_{\ell}\cap \bigcup\mathcal{P}$ is a subgraph of $H_{\ell}\cap \mathcal{L}$, where $H_{\ell}$ is the inner graph of $C_{\ell}$ in $\rho$.
    In particular, the endpoints of $\mathcal{P}$ on $V(C_1)$ coincide with the endpoints of $\mathcal{L}$ on $V(C_1)$, and
    
    \item $H_k \cap \bigcup \mathcal{P}$ is a subgraph of $H_k \cap \mathcal{R}$, where $H_k$ is the outer graph of $C_k$ in $\rho$.
    In particular, the endpoints of $\mathcal{P}$ on $V(\Omega)$ coincide with the endpoints of $\mathcal{R}$ on $V(\Omega)$.
\end{enumerate}
Moreover, there exists an algorithm that finds $\mathcal{P}$ in time $\mathbf{O}(r|E(G)|)$.
\end{lemma}
\begin{proof}
Let $(G',\Omega')$ be the $C_k$-society in $\rho$ and let $X$ be the set of all vertices $x$ such that $x$ is the first vertex in $V(G') \cap \pi(N(\rho))$ encountered when traversing along some $R\in\mathcal{R}$ starting in $V(\Omega)$ and let $\rho'$ be the restriction of $\rho$ to $(G',\Omega')$.
Notice that $X$ is a subset of the vertices of the outer graph of $C_{k-1}$.

Moreover, let $G''$ be the outer graph of $C_\ell$ in $\rho'$.
Finally, let $Y$ be the set of all vertices $y$ such that $y$ is the last vertex of $V(G'') \cap \pi(N(\rho))$ encountered when traversing along some path $L \in \mathcal{L}$ starting from $V(\Omega')$.
It follows that $Y$ is a subset of the vertices of the inner graph of $C_{\ell+1}$.

Notice that $G''$ has a vortex-free $\Delta$-decomposition, where $\Delta$ is a disk, such that
\begin{itemize}
    \item the vertices from $X$ are grounded and on a common face,
    \item the vertices from $Y$ are grounded and on a common face, and
    \item the cycles in $\{ C_{\ell+1}, \ldots , C_{k-1} \}$ are grounded and the trace of each $C_i$ with $i \in[\ell+1, k-1]$ separates $X$ from $Y$ in $\Delta$.
\end{itemize}
Now let $G'''\coloneqq G'\cap G''$.
Then, for each $x \in X$ there exists $R \in \mathcal{R}$ and a path $R' \subseteq R \cap G'''$ such that $R'$ connects $x$ and $V(C_\ell)$.
Similarly, for every $y \in Y$ there exists $L \in \mathcal{L}$ and a path $L' \subseteq L \cap G'''$ such that $L'$ connects $y$ and $V(C_k)$.

If there exists an $X$-$Y$-linkage of order $r$ in $G'''$, we may combine this linkage with the $V(\Omega)$-$X$-linkage consisting of subpaths of the paths in $\mathcal{R}$ contained in the outer graph of $C_k$ in $\rho$ and the $Y$-$V(C_1)$-linkage consisting of subpaths of $\mathcal{L}$ contained in the inner graph of $C_\ell$ in $\rho$.
This would then result in a $V(\Omega)$-$Y$-linkage as required.
Moreover, recall that via an algorithmic version of \Cref{prop:mengersthm} we are able to find such a linkage, if it exists, in time $\mathbf{O}(r|E(G)|)$.

So we may suppose, towards a contradiction, that there does not exist an $X$-$Y$-linkage of order $r$ in $G''$.
In this case, by \cref{prop:mengersthm}, there must exist a set $S \subseteq V(G')$ of size at most $r-1$ such that $G''-S$ does not contain an $X$-$Y$-path.
However, there must exist $L \in \mathcal{L}$, $R \in \mathcal{R}$, and $i \in [\ell+1,k-1]$ such that $H \coloneqq G'' \cap (L \cup R \cup C_i)$ does not contain a vertex from $S$.
Since both $\mathcal{L}$ and $\mathcal{R}$ are orthogonal to $\mathcal{C}$ it follows that $H$ is connected.
Since $V(H) \cap X \neq \emptyset$ and $V(H) \cap Y \neq \emptyset$ this is a contradiction to the choice of the set $S$.
\end{proof}

We can now move on to properly orthogonalising a radial linkage.
Before proving this lemma, we first define a property which we wish to preserve whenever we modify a radial linkage such that it becomes orthogonal to a nest.

Let $(G,\Omega)$ be a society, let $\rho$ be a cylindrical rendition of $(G,\Omega)$ in a disk $\Delta$, and let $\mathcal{P}$ and $\mathcal{Q}$ be two linkages of the same order.
If for each path $Q \in \mathcal{Q}$ there exists a path $P \in \mathcal{P}$ with the same endpoints as $Q$, we say that $\mathcal{P}$ and $\mathcal{Q}$ are \emph{end-identical}.

\begin{lemma}\label{lem:radialtoorthogonal}
    Let $s,r$ with $r \leq s$.
    Let $(G,\Omega)$ be a society with a cylindrical rendition $\rho$ around a vortex $c_0$ with a cozy nest $\mathcal{C}$ of order $s$ and a radial linkage $\mathcal{R}$ of order $r$ for $\mathcal{C}$.

    Then there exists a radial linkage $\mathcal{R}'$ of order $r$ for $\mathcal{C}$ that is orthogonal to $\mathcal{C}$ and end-identical to $\mathcal{R}$.

    Moreover, there exists an algorithm running in $\mathbf{O}(r|E(G)|)$-time that finds $\mathcal{R}'$.
\end{lemma}
\begin{proof}
    Note that since $\mathcal{R}$ is a radial linkage, there are no non-trivial $C_1$-paths in any of the paths in $\mathcal{R}$, as the only vertices used in $\mathcal{R}$ that are found in $V(C_1)$ are the endpoints of the paths that do not lie in $V(\Omega)$.
    Let the endpoints of the paths in $\mathcal{R}$ that lie in $V(\Omega)$ be $a_1, \ldots , a_r$ found in the given order in $\Omega$, such that $R_i$ has the endpoint $a_i$ for all $i \in [r]$.
    In the context of this proof, we will always traverse along the paths of our radial linkages starting from their unique endpoint on $V(\Omega)$.

    We now inductively construct radial linkages $\mathcal{R}_s, \dots , \mathcal{R}_1$ such that for each $i \in [s]$, the radial linkage $\mathcal{R}_i$ is end-identical to $\mathcal{R}$ and orthogonal to the cycles $\{ C_i, \dots , C_s \}$.

    Let us begin with the construction of $\mathcal{R}_s$.
    For each $i \in [r]$ let $x_s^i$ be the first vertex of $R_i \in \mathcal{R}$ on $C_s$ and let $y_s^i$ be its last.
    We claim that there exists a subpath $P_s^i$ of $C_s$ with endpoints $x_s^i$ and $y_s^i$ such that $P_s^i$ is disjoint from all paths $R \in \mathcal{R} \setminus \{ R_i \}$.
    
    To see this, first let $J_s$ be the trace of $C_s$ and $T_i$ be the trace of the $x_s^i$-$y_s^i$-subpath of $R_i$ in $\rho$.
    Since $\mathcal{C}$ is cozy and $R_i$ has exactly one vertex on $C_1$, the trace $T_i$ must be fully contained within the annulus $A_s$ whose boundary is $J_s$ and the trace of $C_1$.
    Moreover, if $\Delta_s^i$ is the $c_0$-disk of $J_s\cup T_i$, then, as $\mathcal{C}$ is cozy, every vertex or edge of any path $R\in\mathcal{R}\setminus\{ R_i\}$ that is drawn in $A_s$ must be drawn in $\Delta_s^i\setminus T_i$.
    Suppose this were not the case, then we note that since $\mathcal{C}$ is cozy, $A_s - \Delta_s^i$ defines a disk whose closure $D_s^i$ is bounded by $T_i$ and a segment of $J_s$.
    Should some path $R$ from $\mathcal{R} \setminus \{ R_i \}$ reach $D_s^i$ this immediately implies that $\mathcal{C}$ is not cozy, since $R$ must also reach $C_1$ and thus, there exists a $C_s$-path in the outer graph of $C_s$ in $R$.
    
    Now observe that there exists a segment $S_s^i$ of $J_s$ with endpoints $\pi^{-1}(x_s^i)$ and $\pi^{-1}(y_s^i)$ which is disjoint from $\Delta_s^i\setminus T_i$ and thus lies in the boundary of $D_s^i$.
    The subpath of $C_s$ whose trace is $S_s^i$ is exactly the path $P_s^i$ we desired.
    Let us now define $\mathcal{R}_s$ to be the linkage obtained from $\mathcal{R}$ by replacing, for each $i\in[r]$, the $x_s^i$-$y_s^i$-subpath of $R_i$ with the path $P_s^i$.
    Then $\mathcal{R}_s$ is a radial linkage of order $r$ for $\mathcal{C}$ which is end-identical to $\mathcal{R}$ and the intersection of any path in $\mathcal{R}_s$ with $C_s$ is a single, possibly trivial, path.
    Hence, $\mathcal{R}_s$ is orthogonal to $\{ C_s\}$ as required.

    Now assume that for some $j\in[3,s]$ we have already constructed the linkage $\mathcal{R}_j=\{ R_1^j,\dots,R_r^j\}$ such that $R_i^j$ has the same endpoints as $R_i$ for each $i\in[r]$.
    Let $J_{j-1}$ be the trace of $C_{j-1}$ and $A_{j-1}$ be the annulus whose boundary is $J_{j-1}$ and the trace of $C_1$.
    Moreover, for each $i\in[r]$ let $x_{j-1}^i$ be the first vertex of $R_i^j$ on $C_{j-1}$, $y_{j-1}^i$ be its last, and let $T_i$ be the trace of the $x_{j-1}^i$-$y_{j-1}^i$-subpath of $R_i^j$.
    It follows from the assumption that $\mathcal{R}_j$ is orthogonal to $\{ C_j,\dots,C_s\}$ and $\mathcal{C}$ is cozy, that $T_i$ is contained within $A_{j-1}$.
    Indeed, the same arguments as before now show that for each $i\in[r]$ there exists a $x_{j-1}^i$-$y_{j-1}^i$-subpath $P_{j-1}^i$ of $C_{j-1}$ which is disjoint from the paths in $\mathcal{R}_j\setminus\{ R_i^j\}$.
    Hence, we may define $R_i^{j-1}$ to be the path obtained from $R_i^j$ by replacing $x_{j-1}^iR_i^jy_{j-1}^i$ with $P_{j-1}^i$ for each $i\in[r]$ and then set $\mathcal{R}_{j-1}\coloneqq \{ R_{j-1}^i ~\!\colon\!~ i\in[r] \}$.
    
    The linkage $\mathcal{R}_1$ we ultimately arrive at is now the orthogonal, radial linkage for $\mathcal{C}$ we searched for and by construction we observe that $\mathcal{R}_1$ is end-identical to $\mathcal{R}$. 
    This completes the proof.
\end{proof}

\subsection{Identifying a crooked witness}
In our efforts to find an orthogonal crooked transaction we will of course want to use \hyperref[prop:mengersthm]{Menger's theorem}.
Using this liberally is however extremely dangerous, as Menger's theorem can completely alter the character of the transaction we will find.
For example, if we simply took the two end segments of a crooked transaction that we want to orthogonalise and ask \hyperref[prop:mengersthm]{Menger's theorem} to find us a large transaction between these end segments, the result may be a planar transaction.

Thus, we will now invest quite some effort into finding and preserving witnesses for the fact that the transaction we initially started with was crooked.
For now, we introduce two types of witnesses that certify that we are dealing with a crooked transaction, based on overpasses and doublecrosses introduced in \Cref{subsec:normalisingcrooked}.

\begin{definition}[Crooked witnesses]\label{def:crookedwitnessI}
    Let $s,p,j$ be positive integers with $2 \leq p$ and $j \leq s$.
    Let $(G, \Omega)$ be a society with a cylindrical rendition $\rho$ around a vortex $c_0$ in the disk $\Delta$ with a cozy nest $\mathcal{C} = \{ C_1, \ldots , C_s \}$ of order $s$.
    Further, let $\mathcal{P}$ be an exposed transaction of order $p$ in $(G,\Omega)$ and let $(G_j,\Omega_j)$ be the $C_j$-society in $\rho$.

    Let $X,Y$ be the two end segments of $\mathcal{P}$ in $\Omega$ and let $X'$, respectively $Y'$, be the set of vertices in $V(C_j) \cap \pi(N(\rho))$ that appear first when traversing the paths in $\mathcal{P}$ from $X$ to $Y$, respectively $Y$ to $X$.
    Note that $X'$ preserves the order of the paths induced by $X$ and the same is true for $Y'$ with respect to $Y$.
    Furthermore, we have $X' \cup Y' \subseteq V(\Omega_j)$.
    Let $I$ be the minimal segment of $\Omega_j$ that contains $X'$ and is disjoint from $Y'$
    and let $J$ be the minimal segment of $\Omega_j$ that contains $Y'$ and is disjoint from $X'$.
    Further, we let $A,B$ be the two distinct, disjoint segments of $\Omega_j$ that have one endpoint in $X'$ and the other in $Y'$ and are otherwise disjoint from $X' \cup Y'$.

    If $R$ is an $A$-$B$-path in $G_j - V(\bigcup \mathcal{P})$, we call $\{ R \}$ a \emph{snitch (for $\mathcal{P}$ on $C_j$)}.
    Two $A$-paths $P,P'$ and two $B$-paths $Q,Q'$ in $G_j - V(\bigcup \mathcal{P})$, such that the pair $P,P'$ form a cross in $(G_j,\Omega_j)$, the pair $Q,Q'$ form a cross in $(G_j,\Omega_j)$, and $P,P',Q,Q'$ are pairwise disjoint, cause us to call $\{ P,P',Q,Q' \}$ a \emph{doublecross (for $\mathcal{P}$ on $C_j$)} (recall \cref{fig:introWitnesses} for an illustration).
    We extend both of these notions in the natural way to snitches and doublecrosses for $\mathcal{P}$ on $(G,\Omega)$ which have all of their endpoints in $V(\Omega)$.
    Whenever there exists a set of paths $\mathcal{W}$ that constitute a snitch or a doublecross for $\mathcal{P}$ on $C_j$ (or $\Omega$), we say that there exists a \emph{crooked witness (for $\mathcal{P}$ on $C_j$ (or $(G,\Omega)$))}.
\end{definition}

Given a crooked transaction of modest size, it is very easy to find a transaction with a crooked witness in $(G,\Omega)$.

\begin{lemma}\label{lem:easycrookedwitness}
    Let $q$ be a positive integer, let $(G,\Omega)$ be a society, and let $\mathcal{Q}$ be a crooked transaction of order $5q-4$.
    
    Then $\mathcal{Q}$ contains a transaction $\mathcal{P}'$ of order $q$ and a crooked witness for $\mathcal{P}'$ in $(G,\Omega)$, which can be found in $\mathbf{O}(q|E(G)|)$-time.
\end{lemma}
\begin{proof}
    Note that $5q - 4 = 4((q+1) -2) + q$.
    Thus \Cref{lem:leapordoublecross} either yields a leap of order $q+1$ from $\mathcal{Q}$ or a doublecross of order $q+4$.
    Both options immediately yield the desired transaction $\mathcal{P}'$ of order $q$ together with a crooked witness.
\end{proof}

We will however instead require a crooked witness deep inside the nest, close to the vortex of the cylindrical rendition.
It is hard to give an intuitive justification for the effort we go through for this purpose at this point.
But such a witness is the key to our efforts of orthogonalising crooked transactions to (most of) a given nest.

\begin{lemma}\label{lem:crookedwitness}
    Let $s,p$ be positive integers with $2 \leq s \leq p$.
    Let $(G, \Omega)$ be a society with a cylindrical rendition $\rho = (\Gamma, \mathcal{D})$ around a vortex $c_0$ in the disk $\Delta$ with a cozy nest $\mathcal{C} = \{ C_1, \ldots , C_s \}$ of order $s$.
    Further, let $\mathcal{P}$ be a crooked transaction of order $8s-14+3p$ in $(G,\Omega)$.

    Then there exists an exposed transaction $\mathcal{P}' \subseteq \mathcal{P}$ of order $p$ that has a crooked witness on $C_2$.

    Furthermore, said transaction and its crooked witness can be found in time $\mathbf{O}( p |E(G)| )$.
\end{lemma}
\begin{proof}
    Let $(G',\Omega')$ be the $C_2$-society in $\rho$ and let $G_2$ be the outer graph of $C_2$ in $\rho$.
    Let $X,Y$ be the two end segments of $\mathcal{P}$ in $\Omega$, with $X'$, respectively $Y'$, being the set of vertices in $V(C_2)$ that appear first when traversing the paths in $\mathcal{P}$ from $X$ to $Y$, respectively $Y$ to $X$.
    As noted in \Cref{def:crookedwitnessI}, both $X'$ and $Y'$ preserve the order of the paths induced by $X$ and respectively $Y$ due to $\rho$ being a cylindrical rendition.
    We let $I$ be the minimal segment of $\Omega'$ that contains $X'$ and is disjoint from $Y'$
    and let $J$ be the minimal segment of $\Omega'$ that contains $Y'$ and is disjoint from $X'$.
    Let $u$ be one of the two endpoints of $I$ and let $P \in \mathcal{P}$ be the path with $u \in V(P)$.
    
    We start with an observation on subpaths of $P$ that immediately help us.
    Let $A,B$ be the two distinct, disjoint segments of $\Omega'$ that have one endpoint in $I$ and the other in $J$ and are otherwise disjoint from $I \cup J$, such that $u \in A$.

    \begin{claim}\label{claim:easyoverpasses}
        If $P$ contains one of the following types of paths, then we can find a snitch for a subtransaction of $\mathcal{P}$ of order $p$:
        \begin{itemize}
            \item an $A$-$B$-path within $G'$ that is disjoint from $I \cup J$,

            \item a $V(\Omega')$-path within $G'$ whose endpoints define two segments of $\Omega'$ such that one of them contains all of $I$, respectively all of $J$, and the other segment contains at least $p$ vertices of $Y$, respectively $X$, other than the endpoints of $P$, or

            \item an $I$-$J$-path within $G'$ such that for a segment $S$ of $\Omega'$ that has the same endpoints as this path, we have $||X' \cap S| - |Y' \cap S|| \geq p + 1$.
        \end{itemize}
    \end{claim}
    \emph{Proof of \Cref{claim:easyoverpasses}:}
        The fact that any of these paths yields a snitch is immediate from their specifications.
    \hfill$\blacksquare$

    Next, the following behaviour of $P$ with respect to vertices in $(I \cup J) \cap (X' \cup Y')$ is fairly basic but central to our later arguments.
    For any $Q \in \mathcal{P}$, let $I_Q$ be a smallest segment of $\Omega'$ that minimises the intersection of $I_Q$ and $B$, contains $A$, and contains all vertices in $V(Q) \cap (I \cup J)$.
    
    \begin{claim}\label{claim:properdivisionofpaths}
        For any path $Q \in \mathcal{P}$ with $X' \cup Y' \not\subseteq I_Q$, any $V(\Omega')$-path $Q'$ in $Q$ whose endpoints separate two vertices in $X' \cup Y'$ from each other in $\Omega'$ is found in $G'$.
    \end{claim}
    \emph{Proof of \Cref{claim:properdivisionofpaths}:}
        Let $Q'$ be such a path and suppose that $Q'$ is not found in $G'$, which implies that it is instead found in $G_2$.
        We let $a \in X' \cup Y'$ be the vertex within $I_Q$ that is separated from a vertex $b \in X' \cup Y'$ within $\Omega'$ by removing the endpoints of $Q'$.
        Furthermore, we can choose $b$ to be outside $I_Q$, since $X' \cup Y' \not\subseteq I_Q$.
        For each vertex in $X' \cup Y'$ there exists a unique path in $G_2 \cap \bigcup \mathcal{P}$ connecting it to $V(\Omega)$, since the vertices in $X' \cup Y'$ are chosen such that they are the first vertices encountered on a path from $\mathcal{P}$ when traversing it from either end.
        We note that $Q'$ must be internally disjoint from all such paths, since it is disjoint from all paths in $\mathcal{P} \setminus \{ Q \}$ and has both endpoints in $V(\Omega')$.
        However, since the endpoints of $Q'$ separate $a$ and $b$ within $\Omega'$ and $G_2$ is grounded in $\rho$, the path $Q'$ must intersect either the path linking $a$ to $V(\Omega)$, or the path linking $b$ to $V(\Omega)$, a contradiction.
    \hfill$\blacksquare$

    Now let us consider the situation in which all interactions between $I \cup J$ and $P$ occur fairly close to $A$.
    We claim that in this situation we can either find a snitch or set ourselves up for finding a doublecross later.
    Let us put the second option in more concrete terms and give a few more definitions before we state our claim.

    Let $\mathcal{L} \subseteq \mathcal{P}$ be a transaction of order $q$ and let $R,R'$ be a cross in $(G',\Omega')$.
    We say that $\mathcal{L}$ and $R,R'$ are a \emph{doublecross-setup of order $q$} if the two endpoints of $R,R'$ are all found in one of the two segments of $\Omega'$ that result from the deletion of the two minimal segments of $\Omega'$ that contain all the vertices in $X' \cap V(\bigcup \mathcal{L})$, respectively $Y' \cap V(\bigcup \mathcal{L})$, and none of the vertices in $Y'$, respectively $X'$.

    Our main focus for the remainder of the proof will be on subpaths of $P$ starting in $u$ and moving up until they have reached a certain distance in $I \cup J$ towards $B$.
    We start with a simple fact -- which derives from $\mathcal{P}$ being crooked -- that will prove to be valuable later on.
    
    Since $\mathcal{P}$ is crooked, if we let $u' \in Y'$ be the vertex of $V(P)$ found in $J$, we know that the segment of $\Omega'$ that has the endpoints $u,u'$ and is disjoint from $I \setminus \{ u \}$ must contain at least one vertex in $Y'$ other than $u'$.
    Let $u'' \in V(uPu') \cap (I \cup J)$, then we define $G_{u''}$ to be the union of $V(\Omega')$ and $uPu'' \cap G'$.

    \begin{claim}\label{claim:leftsideoverpassorcross}
        Let $q$ be an integer with $q \leq 2s - 3$.
        If $|I_P \cap (X' \cup Y')| \leq q$, we can either find an snitch for a subtransaction of $\mathcal{P}$ of order $p$, or we can find a doublecross-setup of order at least $(8s - 14 + 3p) - (2q + p - 1)$ in time $\mathbf{O}( p |E(G)| )$.
    \end{claim}
    \emph{Proof of \Cref{claim:leftsideoverpassorcross}:}
        We note that we have $|I_P \cap (X' \cup Y')| \geq 3$, since $u$ and $u'$ must be contained in $I_P$ and furthermore, as just noted prior to this claim, there exists a third vertex $v' \in Y' \cap I_P$ found in the segment of $\Omega'$ that contains $A$, $u$, and $u'$.
        We let $Q$ be the path that contains $v'$ and let $v \in X'$ be the other vertex of $Q$ found in $X' \cup Y'$.

        Suppose that $(G_{u'}, \Omega')$ contains a cross consisting of the paths $R,R'$, which we can easily check in linear time.
        Then the segment of $\Omega'$ that contains $A$ and all endpoints of $R,R'$ is a subset of $I_P$.
        Thus, by removing all paths in $\mathcal{P}$ that contain a vertex of $I_P \cap (X' \cup Y')$, we can find a transaction of order at least $(8s - 14 + 3p) - (q - 1)$ that together with $R,R'$ forms a doublecross-setup of order $(8s-14+3p) - (q - 1)$.

        Thus we may instead suppose that $(G_{u'},\Omega')$ in fact does not contain a cross.
        In particular, this implies that we may imagine $P$ as being drawn in a planar fashion within $\Gamma$, such that $P$ forms a curve $T$ that does not intersect itself, has its endpoints on the boundary of $\Delta$, and is otherwise disjoint from the boundary of $\Delta$.

        Consider the closures of the two disks within $\Delta - T$.
        Then exactly one of the two disks contains the drawing of $v'$.
        Let $\Delta_{v'}$ be the disk containing $v'$ and let $\Delta_v$ be the other disk.
        Other than the vertices in $V(P)$, each vertex of $G$ is drawn on exactly one of these two disks, since $\Gamma$ cannot draw vertices outside of $P$ onto the points representing the edges in $E(P)$.
        \Cref{claim:properdivisionofpaths} now implies that $\Delta_{v'}$ does not contain any of the vertices of $X' \setminus \{ u \}$, which means that $v$ is found in $\Delta_v$.
        (See \Cref{fig:witness} for an illustration of this configuration.)

        \begin{figure}[ht]
            \centering
            \scalebox{0.85}{
            \begin{tikzpicture}[scale=1.25]

            \pgfdeclarelayer{background}
		\pgfdeclarelayer{foreground}
			
		\pgfsetlayers{background,main,foreground}

            \begin{pgfonlayer}{background}
            \pgftext{\includegraphics[width=6cm]{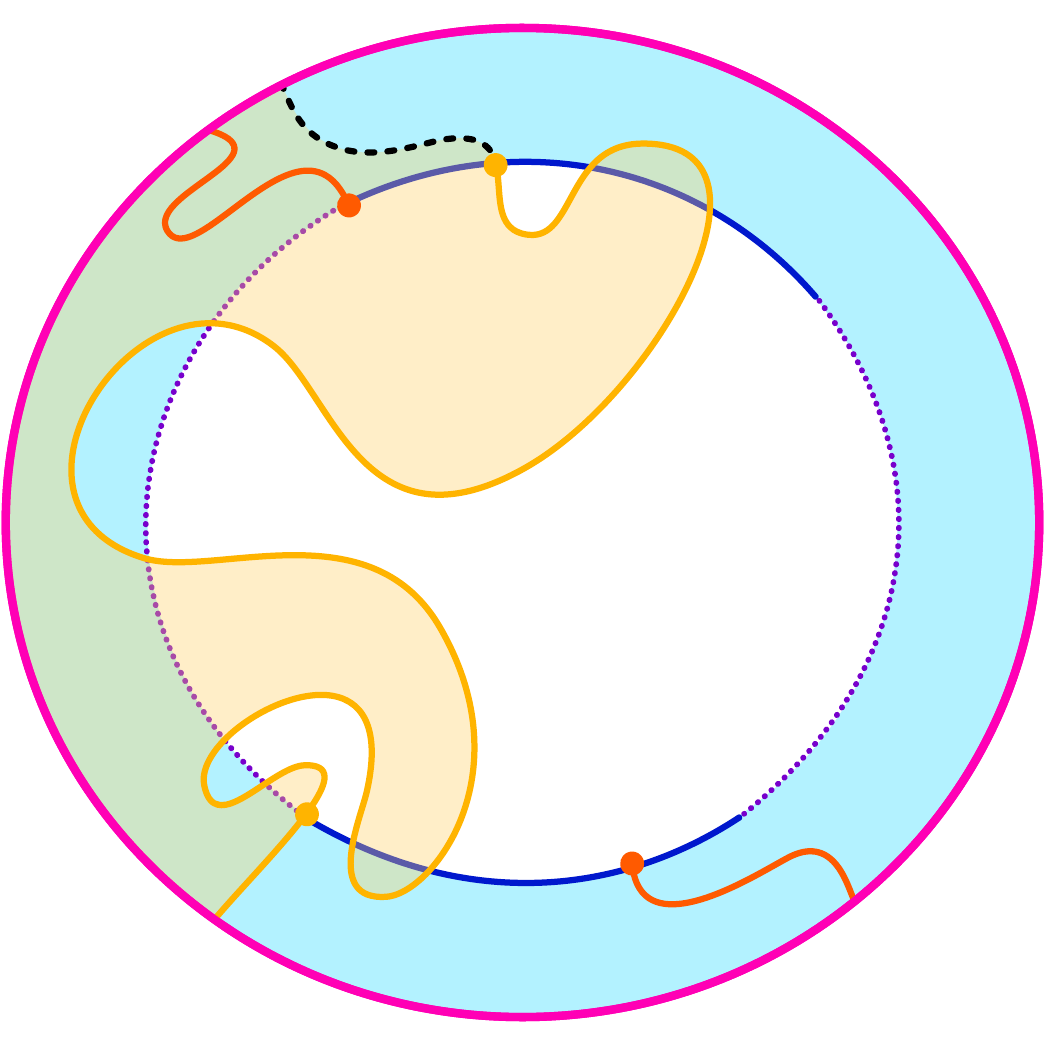}} at (C.center);
            \end{pgfonlayer}{background}
			
            \begin{pgfonlayer}{foreground}

            \node (THEPATH) at (-0.1,-1.05) [draw=none] {$P$};
            \node (DISKV) at (0.2,0.8) [draw=none] {$\Delta_v$};
            \node (DISKVPRIME) at (0.6,0.4) [draw=none] {$\Delta_{v'}$};
            \node (A) at (-2.3,0.5) [draw=none] {$A$};
            \node (B) at (2.35,0) [draw=none] {$B$};
            \node (V) at (-1,1.6) [draw=none] {$v$};
            \node (VPRIME) at (0.65,-1.7) [draw=none] {$v'$};
            \node (I) at (0,-2.25) [draw=none] {$I$};
            \node (J) at (1.55,1.7) [draw=none] {$J$};
            \node (U) at (-1.25,-1.9) [draw=none] {$u$};
            \node (UPRIME) at (0,2.3) [draw=none] {$u'$};
            
            \end{pgfonlayer}{foreground}
        
            \begin{pgfonlayer}{foreground}
            \end{pgfonlayer}{foreground}
        
            \end{tikzpicture}
            }
            \caption{An illustration of the situation we arrive in during the proof of \Cref{claim:leftsideoverpassorcross}. The faint yellow disk corresponds to $\Delta_v$, while the disk in the middle described by the dotted purple and thicker blue lines is the $C_2$-society. The dashed line in black indicates a line that will be followed onto $\Omega$ within the context of \Cref{claim:leftsideoverpassorcross}. From the perspective of \Cref{claim:farreachoverpassorcross} the reader can replace $u'$ with $u''$ with the dashed line in black indicating the line that $Q$ must pass through.}
            \label{fig:witness}
        \end{figure}
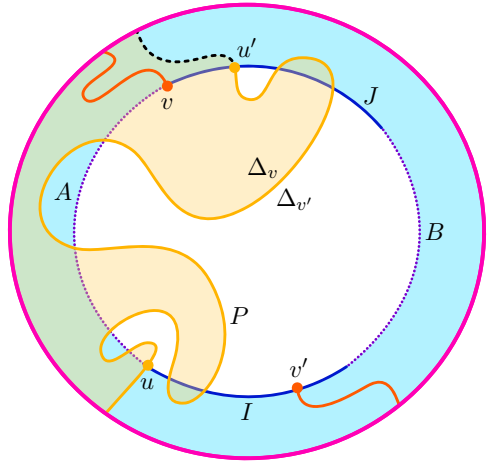

        As a consequence of this, the drawing of some $V(\Omega')$-path within $vQv'$ in $\Gamma$ must cross $T$.
        Note that $P$ contains two disjoint paths connecting both $u$ and $u'$ to $V(\Omega)$ within $G_2$.
        These two paths separate $v$ and $v'$ in $G_2$, since $P$ and $Q$ are disjoint.
        In particular, $u$ and $u'$ separate $v$ and $v'$ within $\Omega'$ thanks to our choice of $v'$ and the fact that $u$ is an endpoint of $I$.
        Thus there exists some pair of paths $P' \subseteq uPu'$ and $Q' \subseteq vQv'$ such that $P',Q'$ form a cross in $(G',\Omega')$.

        Note that $Q'$ must have at least one endpoint in $I_P$.
        Thanks to \Cref{claim:easyoverpasses}, we know that the minimal segment $I^*$ of $\Omega'$ that contains $I_P$ and all endpoints of $Q'$ either contains at most $2q + p$ vertices from $X' \cup Y'$ or $Q'$ is a snitch for a subtransaction of $\mathcal{P}$ of order $p$.
        In the second case, we are done.
        Otherwise, by removing the paths corresponding to the at most $2q + p$ endpoints in $I^*$ -- two of which belong to $P$ -- we can use $P',Q'$ to form a doublecross-setup of order $(8s-14+3p) - (2q+p-1)$.
        This proves our claim and following the arguments of our proof yields a procedure that finds the desired witnesses in time $\mathbf{O}( p |E(G)| )$.
    \hfill$\blacksquare$

    In general there is no guarantee that we are given a chance to apply \Cref{claim:leftsideoverpassorcross}.
    Thus we will have to argue that we can find a snitch or a doublecross-setup if $P$ wanders too far away from $A$.

    \begin{claim}\label{claim:farreachoverpassorcross}
        If $|I_P \cap (X' \cup Y')| \geq 2s-3$, we can either find a snitch for a subtransaction of $\mathcal{P}$ of order $p$, or we can find a doublecross-setup of order $(8s - 14 + 3p) - (2s - 3 + p)$.
    \end{claim}
    \emph{Proof of \Cref{claim:farreachoverpassorcross}:}
        Let $u'' \in V(P) \cap (I \cup J)$ be the vertex minimising the length of $uPu''$ such that the smallest segment $I_{u''}$ of $\Omega'$ that contains all vertices of $uPu''$ also contains at least $2s-3$ vertices from $X' \cup Y'$.
        Note that this may overshoot our goal of $2s-3$ vertices by at most $p$ vertices according to \Cref{claim:easyoverpasses}.
        If $(G_{u''},\Omega')$ contains a cross, we can therefore find a snitch or a doublecross-setup of the appropriate order using the same methods as in our proof of \Cref{claim:leftsideoverpassorcross}.
        Thus we may instead suppose that $(G_{u''},\Omega')$ does not contain a cross.
        
        Let $X'' = (I_{u''} \cap X') \setminus \{ u \}$ and let $Y'' = I_{u''} \cap Y'$.
        By pigeonhole principle, since $I_{u''}$ contains at least $2s-3$ vertices of $X' \cup Y'$, one of the sets $X''$ or $Y''$ contains at least $s-1$ vertices.
        W.l.o.g.\ we let this set be $X''$.
        Let $w$ be the endpoint of $I_{u''}$ that lies in $I$ and note that thus $uIw$ contains $X''$.
        This further implies that $uPw$ is a subpath of $P$ such that the smallest segment of $\Omega'$ that contains all vertices of $uPw$ also contains $X''$.

        Similar to the proof of \Cref{claim:leftsideoverpassorcross}, we note that we can assume that $Pw$ is drawn by $\Gamma$ as a curve $T$ in $\Delta$ that does not intersect itself and has one end on the boundary of $\Delta$, the other end on the boundary of $c_0$ and such that $T$ is disjoint from the boundary of $\Delta$ with the exception of one of its endpoints.
        For each $x \in X''$ let $x'$ be the vertex of $Y'$ that is found in the same path of $\mathcal{P}$ to which $x$ belongs and let $X^*$ be the collection of these vertices.
        Let $x \in X''$ be chosen arbitrarily, let $Q \in \mathcal{P}$ be the path with $x \in V(Q)$, and let $x' \in V(Q) \cap X^*$.
        
        Thanks to \Cref{claim:properdivisionofpaths} and as $x \in X''$, for $Q$ to reach $x' \in J$, the drawing of $Q$ in $\Gamma$ must either cross $T$ within $c_0$ and thus verify the existence of a cross in $(G',\Omega')$ consisting of two paths $P' \subseteq Pu_i$ and $Q' \subseteq Q$, or there exists a $V(\Omega')$-path $Q''$ in $Q$ that is found in $G_2$ such that there exists a segment of $\Omega'$ with the same endpoints as $Q''$ that contains $w$ and is disjoint from $A$, $B$, and $J$.

        For the second option, as $\mathcal{C}$ is cozy, we note that this path must intersect some cycle $C \in \mathcal{C} \setminus \{ C_1, C_2 \}$ such that $Q''$ does not contain any non-trivial $C$ paths in the outer graph of $C$ in $\rho$.
        Thus, there can exist at most $s-2$ vertices in $X''$ for which the second option holds.
        Since $|X''| > s-2$, there must therefore exist some vertex in $X''$ for which the first option is true.
        Hence, we may assume w.l.o.g.\ that the cross in $(G',\Omega')$ consisting of paths $P' \subseteq Pw$ and $Q' \subseteq Q$ exists.

        Via an entirely analogous argument to the one presented in our proof of \Cref{claim:leftsideoverpassorcross}, we can now see that either $Q'$ is a snitch for some subtransaction of $\mathcal{P}$ of order $p$, or there exists a doublecross-setup, with $P',Q'$ being the cross, for a transaction of order at least $8s-14+3p - (2s-3 + p)$.
        This proves our claim.
    \hfill$\blacksquare$

    We may therefore restrict ourselves to the case in which we have somehow found a doublecross-setup of order at least $(8s - 14 + 3p) - (4s - 7 + p)$ as a result of applying \Cref{claim:farreachoverpassorcross} or \Cref{claim:leftsideoverpassorcross}.
    Let $P',Q'$ be the cross within this doublecross-setup and let $\mathcal{L}$ be the transaction of order $4s-7+2p$ involved in it.

    We now want to reapply our arguments from \Cref{claim:farreachoverpassorcross} and \Cref{claim:leftsideoverpassorcross} to $\mathcal{L}$ from the perspective of $B$.
    For this purpose we use the endpoint $w \in B$ of $I$ that is distinct from $u$ for this purpose instead of $u$.
    Note that $w$ may be a vertex of a path $R \in \mathcal{P}$ with $P' \subseteq R$ or $Q' \subseteq R$ and thus in particular, it may be the case that $R \not\in \mathcal{L}$.
    However, we require the fact that $R$ cannot be a peripheral path of $\mathcal{P}$ to be able to apply \Cref{claim:leftsideoverpassorcross}.
    So we simply add $R$ back into the doublecross-setup if it was previously deleted.

    This clearly does not cause a problem if \Cref{claim:farreachoverpassorcross} or \Cref{claim:leftsideoverpassorcross} return a snitch.
    If they instead return a doublecross-setup, then we note that this doublecross-setup has order at least $p = (4s - 7 + 2p) - (4s - 7 + p)$.
    Let $P'',Q''$ be the second cross we found in this process.
    Note that the endpoints of $P'',Q''$ must be found in a segment of $\Omega'$ that contains $B$ and at most $4s - 7 + p$ vertices of $X' \cup Y'$.
    Similarly, the endpoints of $P',Q'$ are found in a segment of $\Omega'$ that contains $A$ and at most $4s - 7 + p$ vertices of $X' \cup Y'$.
    Thus, since $p$ is a positive integer and we started with a transaction of order $8s-14+3p$, the second cross $P'',Q''$ must be disjoint from $P',Q'$ and we have therefore actually found our desired doublecross.
\end{proof}

One can generalise the proof of this lemma in a straightforward way to find crooked witnesses on any given cycle $C_j$ of the nest, if we start with a crooked transaction of order $8s-(4j+2)+3p$.

\subsection{Orthogonalising transactions}
We now move on to presenting a surprisingly simple method that allows us to orthogonalise a given transaction to all but the innermost cycle of the nest.
The reason why this does not immediately solve all our problems is that, as mentioned earlier, this process does not preserve the type of a transaction.

\begin{lemma}\label{lemma:orthogonal_transaction}
Let $s,p$ be positive integers.
Let $(G, \Omega)$ be a society with a cylindrical rendition $\rho$ around a vortex $c_0$ in the disk $\Delta$ with a cozy nest $\mathcal{C} = \{ C_1, \ldots , C_s \}$ of order $s$.
Further, let $\mathcal{P}$ be an exposed transaction of order $s(p-1)+1$ in $(G,\Omega)$ with the end segments $X_1, X_2$.

Then there exists an exposed transaction $\mathcal{P}'$ of order $p$ such that
\begin{enumerate}
    \item $\mathcal{P}'$ is orthogonal to $\{ C_2, C_3, \dots , C_s \}$,
    
    \item $\mathcal{P}'$ connects vertices of $X_1 \cap V(\mathcal{P})$ to vertices of $X_2 \cap V(\mathcal{P})$, and

    \item the intersection of $\bigcup\mathcal{P}'$ with the union of the inner graph of $C_1$ in $\rho$ and $C_1$ is fully contained in $C_1\cup\bigcup\mathcal{P}$.
\end{enumerate}
Moreover, there exists an algorithm that finds $\mathcal{P}'$ in time $\mathbf{O}( p |E(G)| )$.
\end{lemma}
\begin{proof}
Let $G'$ be the union of $C_1$ and the intersection of the inner graph of $C_1$ in $\rho$ with $\bigcup \mathcal{P}$.
We define two radial linkages $\mathcal{L}_1$ and $\mathcal{L}_2$ as follows.
For $i \in [2]$ and every $P \in \mathcal{P}$, let $P_i$ be the unique $X_i$-$V(C_1)$-subpath of $P$.
We set $\mathcal{L}_i \coloneqq \{ P_i ~\!\colon\!~ P\in\mathcal{P} \}$.
Let $X_i'$ be the set of all endpoints of the paths in $\mathcal{L}_i$ on $C_1$.
Notice that $\mathcal{L}_1 \cup \mathcal{L}_2$ is also a radial linkage.

Now suppose there exists a linkage $\mathcal{Q}$ of order $p$ from $X_1'$ to $X_2'$ in $G'$.
By using an algorithmic version of \cref{prop:mengersthm} for finding $p$ paths, such a linkage can be found in time $\mathbf{O}( p |E(G)| )$.
Let $\mathcal{L}$ be the collection of all paths from $\mathcal{L}_1 \cup \mathcal{L}_2$ that share an endpoint with some path in $\mathcal{Q}$.
We now apply \cref{lem:radialtoorthogonal} to the nest $\mathcal{C}$ and the radial linkage $\mathcal{L}$ to obtain a radial linkage $\mathcal{L}'$ such that $\mathcal{L}$ and $\mathcal{L}'$ are end-identical, and $\mathcal{L}'$ is orthogonal to $\mathcal{C}$.
For each path $Q \in \mathcal{Q}$, let $Q'$ be the union of $Q$ and the two paths from $\mathcal{L}'$ that share an endpoint with $Q$.
Finally, let $\mathcal{P}'$ be the collection of all such paths $Q'$.
It follows that $\mathcal{P}'$ is an exposed transaction of order $p$ in $(G, \Omega)$ and $\mathcal{P}'$ is orthogonal to $\{ C_2, C_3, \dots , C_s \}$.
This is the outcome claimed in the assertion and thus, in this case, we are done. 

So now we may assume that there does not exist an $X_1'$-$X_2'$-linkage of order $p$ in $G'$.
By \cref{prop:mengersthm} there exists an $X_1'$-$X_2'$-separation $(A,B)$ of order at most $p-1$ in $G'$.
Let $S \coloneqq A \cap B$.

For every $P \in \mathcal{P}$, let $\widetilde{P}$ be the unique $X_1'$-$X_2'$-subpath of $P$.
Let $\widetilde{\mathcal{P}}$ be the collection of all paths $\widetilde{P}$ as defined above.

Let $F$ be a $V(C_1)$-path in $G$ that sticks out away from $c_0$ and let $F'$ be the unique subpath of $C_1$ such that $F'$ has the same endpoints as $F$ and $F \cup F'$ is the unique cycle in $C_1 \cup F$ whose trace does not separate $c_0$ and the boundary of $\Delta$.
We call $F'$ the \emph{gangplank} of $F$.
Let $\widetilde{P} \in \widetilde{\mathcal{P}}$.
A subpath $F'$ of $C_1$ is a \emph{gangplank} of $\widetilde{P}$ if there exists a subpath $F$ of $\widetilde{P}$ such that $F'$ is the gangplank of $F$.

We observe that $\widetilde{\mathcal{P}}$ is an $X_1'$-$X_2'$-linkage in $G$ of order $s(p-1)+1$.
Hence, there must exist a linkage $\widehat{\mathcal{P}} \subseteq \widetilde{\mathcal{P}}$ of order $(s-1)(p-1)+1$ that avoids the $p-1$ vertices in $S$.
Since $\mathcal{P}$ is exposed, every path in $\widehat{\mathcal{P}}$ must contain an edge which is drawn in the interior of $c_0$ and thus has at least one gangplank.

Let $S' \coloneqq S \cap V(C_1)$.
Suppose there exists a path $\widehat{P} \in \widehat{\mathcal{P}}$ such that all of its gangplanks avoid $S'$.
Let $J$ be the union of $G' \cap \widehat{P}$ and all gangplanks of $\widehat{P}$.
Then $J$ must be a connected subgraph of $G'$ that avoids $S$ but contains a vertex of $X_1'$ and a vertex of $X_2'$.
This is a contradiction to the choice of $S$.
Hence, every path $\widehat{P} \in \widehat{\mathcal{P}}$ has a gangplank that contains some vertex of $S'$.
With $|S'| \leq p-1$ and $|\widehat{\mathcal{P}}| \geq (s-1)(p-1)+1$ there must exist some $u \in S'$ such that $u$ is contained in a gangplank of $s$ paths $\widehat{P}_1, \dots , \widehat{P}_s \in \widehat{\mathcal{P}}$.
For each $i \in [s]$, let $F_i$ be a subpath of $\widehat{P}_i$ such that $u$ belongs to the gangplank $F'_i$ of $F_i$.

The trace of $F_i \cap F_i'$ bounds a disk $\Delta_i \subseteq \Delta$ that does not contain $c_0$ since $F_i$ sticks out away from $c_0$, for each $i \in [s]$.
Moreover, since the paths from $\widehat{\mathcal{P}}$ are pairwise vertex-disjoint, so are the paths $F_1, \ldots , F_s$.
It follows that for all distinct $i, j \in [s]$ we have either $\Delta_i \subsetneq \Delta_j$ or $\Delta_j \subsetneq \Delta_i$.
Thus, we may assume that the paths $F_1, \ldots , F_s$ are numbered such that $\Delta_1 \subsetneq \Delta_2 \subsetneq \dots \subsetneq \Delta_s$.

Since $\mathcal{C}$ is cozy, there cannot exist an $i \in [s]$ and a $j \in [s]$ such that $F_i$ contains a $C_j$-path that sticks out away from $c_0$ and is disjoint from the cycles in $\{ C_{j'} ~\!\colon\!~ j' \in [j+1,s] \}$.
It follows that for every $i \in [s]$ there exists a $j_{i} \in [2,s]$ such that $j_{i}$ is maximal with the property that $C_{j_{i}}$ contains a vertex of $F_{i}$ and every non-trivial $V(C_{j_i})$-subpath of $F_i$ sticks out towards $c_0$.
Hence, the trace of $C_{j_i}$ can intersect the disk $\Delta_i$ only in its boundary.
However, for any $i'<i$ the disk $\Delta_{i'}$ cannot contain a point from the boundary of $\Delta_i$.
Thus, $i'<i$ implies $j_{i'}<j_i$.
Since $s > s - 1 = |[2,s]|$, this is absurd.
Hence, this case cannot occur and our proof is complete.
\end{proof}

The first situation in which we will need the above lemma will require us to additionally protect a crooked witness.
This results in a slightly more technical statement which has a proof that is very close to that of \Cref{lemma:orthogonal_transaction}.
We include \Cref{lemma:orthogonal_transaction} and its proof here because we believe that it -- and in particular its proof -- may be of independent interest to readers and we will need a further variant of it later on.
Thus we only sketch the proof of the following lemma.

\begin{lemma}\label{lem:orthogonaltransactionwithwitness}
Let $s,p$ be positive integers with $s \geq 3$.
Let $(G, \Omega)$ be a society with a cylindrical rendition $\rho$ around a vortex $c_0$ in the disk $\Delta$ with a cozy nest $\mathcal{C} = \{ C_1, \ldots , C_s \}$ of order $s$.
Further, let $\mathcal{P}$ be an exposed transaction of order $(p+7)(s-2)+p$ in $(G,\Omega)$ with the end segments $X_1, X_2$ and let there exist a crooked witness $\mathcal{W}$ for $\mathcal{P}$ on $C_2$.

Then there exists an exposed transaction $\mathcal{P}'$ of order $p$ such that
\begin{enumerate}
    \item $\mathcal{W}$ is a crooked witness for $\mathcal{P}'$ on $C_2$,

    \item $\mathcal{P}'$ is orthogonal to $\{ C_3, C_4, \dots , C_s \}$, and
    
    \item $\mathcal{P}'$ connects vertices of $X_1 \cap V(\mathcal{P})$ to vertices of $X_2 \cap V(\mathcal{P})$.
\end{enumerate}
Moreover, there exists an algorithm that finds $\mathcal{P}'$ in time $\mathbf{O}( p |E(G)| )$.
\end{lemma}
\begin{proof}[Proof sketch]
    Suppose that $\mathcal{P}$ has a doublecross $\{ Q_1,Q_2,Q_3,Q_4 \}$, with $Q_1,Q_2$ and $Q_3,Q_4$ respectively forming crosses in the $C_2$-society $(G_2, \Omega_2)$.
    For each $i \in [4]$, let $x_i,x_i'$ be the two distinct endpoints of $Q_i$.
    Otherwise $\mathcal{P}$ has a snitch.
    In this case, we call the corresponding path $Q_1$, and let $x_1,x_1'$ be defined as above, with $x_i := x_1$, $x_i' := x_1'$, and $Q_i := Q_1$ for all $i \in [2,4]$.

    Let $G'' \coloneqq C_2 \cup G_2$, let $Z \coloneqq \{ x_i,x_i' ~\!\colon\!~ i \in [4] \}$, and note that $Z = V(C_2) \cap V( \bigcup_{i = 1}^4 Q_i )$.
    We let $G'$ be $G'' - ( V( \bigcup_{i = 1}^4 Q_i ) \setminus Z )$.
    
    We construct two radial linkages $\mathcal{L}_1, \mathcal{L}_2$ orthogonal to $\{ C_3, C_4, \dots , C_s \}$ from the paths in $\mathcal{P}$ as in the proof of \Cref{lemma:orthogonal_transaction} and let $X_1',X_2'$ be the endpoints of these linkages in $V(C_2)$.
    Let $\widetilde{P}$ be the $X_1'$-$X_2'$-path within $P$, for each $P \in \mathcal{P}$, and let $\widetilde{\mathcal{P}}$ be the collection of these paths.

    If there exists a linkage $\mathcal{Q}$ of order $p$ from $X_1'$ to $X_2'$ in $G'$, we proceed exactly as in the proof of \Cref{lemma:orthogonal_transaction} and thus find the desired transaction in time $\mathbf{O}( p |E(G)| )$.
    Therefore we need to argue that the initial size of $\mathcal{P}$ is sufficient to guarantee the existence of such a linkage $\mathcal{Q}$.

    Suppose $\mathcal{Q}$ does not exist, let $(A,B)$ be an $X_1'$-$X_2'$-separator of order at most $p-1$, which must exists according to \Cref{prop:mengersthm}, and let $S := (A \cap B) \cup Z$.
    Note that $|Z| \leq 8$.
    Thus we have $|S| \leq p+7$, since $(A,B)$ is a separator of order at most $p-1$.
    Hence there exists a linkage $\widehat{\mathcal{P}} \subseteq \widetilde{\mathcal{P}}$ of order at least $(p+7)(s-2)+1$ that avoids $S$, as $Z$ is disjoint from $\mathcal{P}$ to begin with.
    Since each of the paths in $\mathcal{P}$ is exposed, each path in $\widehat{\mathcal{P}}$ has at least one gangplank.
    Let $S' := S \cap V(C_2)$.
    
    We define \emph{gangplanks} for $V(C_2)$-paths that stick out away from $c_0$ and paths $\widetilde{P} \in \widetilde{\mathcal{P}}$ as in the proof of \Cref{lemma:orthogonal_transaction}.
    Once more, we conclude that there does not exist a path $\widehat{P} \in \widehat{\mathcal{P}}$ such that all of its gangplanks avoid $S'$.
    Hence, every path $\widehat{P} \in \widehat{\mathcal{P}}$ has a gangplank that contains some vertex of $S'$.
    As $|S'| \leq p+7$ and $|\widehat{\mathcal{P}}| \geq (p+7)(s-2)+1$, there exists some $u \in S'$ that is contained in a gangplank of $s-1$ paths $\widehat{P}_1, \ldots , \widehat{P}_{s-1} \in \widehat{\mathcal{P}}$.
    From here we can then reach a contradiction in just the same fashion as we did in the proof of \Cref{lemma:orthogonal_transaction}.
\end{proof}

Note that similar to \Cref{lem:crookedwitness}, this lemma could be generalised to allow for us to find a crooked witness on $C_j$ assuming we start with a transaction of order $(p+7)(s-j)+p$.

\subsection{Orthogonalising a crooked transaction}\label{sec:crookedmagic}
We are now almost ready to prove the main theorem of this section.
Roughly described, it will find a crooked transaction in a society with a cylindrical rendition and a nest which is orthogonal to all but a constant number of cycles in the nest.
Thus it is somewhat of an analogue to \Cref{lem:crookedexistencealgo}, although our theorem will also preserve the existing cylindrical rendition -- and thus the nest -- to a large degree.

Kawarabayashi et al.\ also need such a theorem and we have in fact already cited it earlier.
The full version of \Cref{prop:oldOrthogonality} (see Lemma 4.5 in \cite{KawarabayashiTW2021Quickly}) takes the same role as the theorem we are about to prove.
However, it causes the loss of a number of cycles from the nest that is \emph{linear} in the size of the crooked transaction we want.
As long as the size of the crooked transaction we require is a fixed value, this would be fine for our approach (see \Cref{lemma:reroute_crosses}).
However, the crooked transactions we want have a size which itself depends on the $K_t$-minor we exclude and thus using their lemma would cause the loss of large numbers of cycles each time we apply it.
Using the lemma from \cite{KawarabayashiTW2021Quickly} would therefore ultimately cause our functions to be exponential.

In reality we will need to generalise the main theorem of this section quite a bit more for our applications (see the start of \Cref{sec:multisocietycrooked}).
We choose to split this presentation because the proof of \Cref{thm:orthogonalexposedsorcery} adds another considerable layer of technical complexity comparable to the already very technical nature of the proof of \Cref{thm:orthogonalexposed}.
By first proving the result in a slightly simpler setting, we hope to aid the reader in understanding our arguments later on in the proof of \Cref{thm:orthogonalexposedsorcery}.

\begin{figure}[ht]
    \centering
    \begin{tikzpicture}

        \pgfdeclarelayer{background}
		\pgfdeclarelayer{foreground}
			
		\pgfsetlayers{background,main,foreground}

        \begin{pgfonlayer}{background}
        \node (C) [v:ghost] {{\includegraphics[width=12cm]{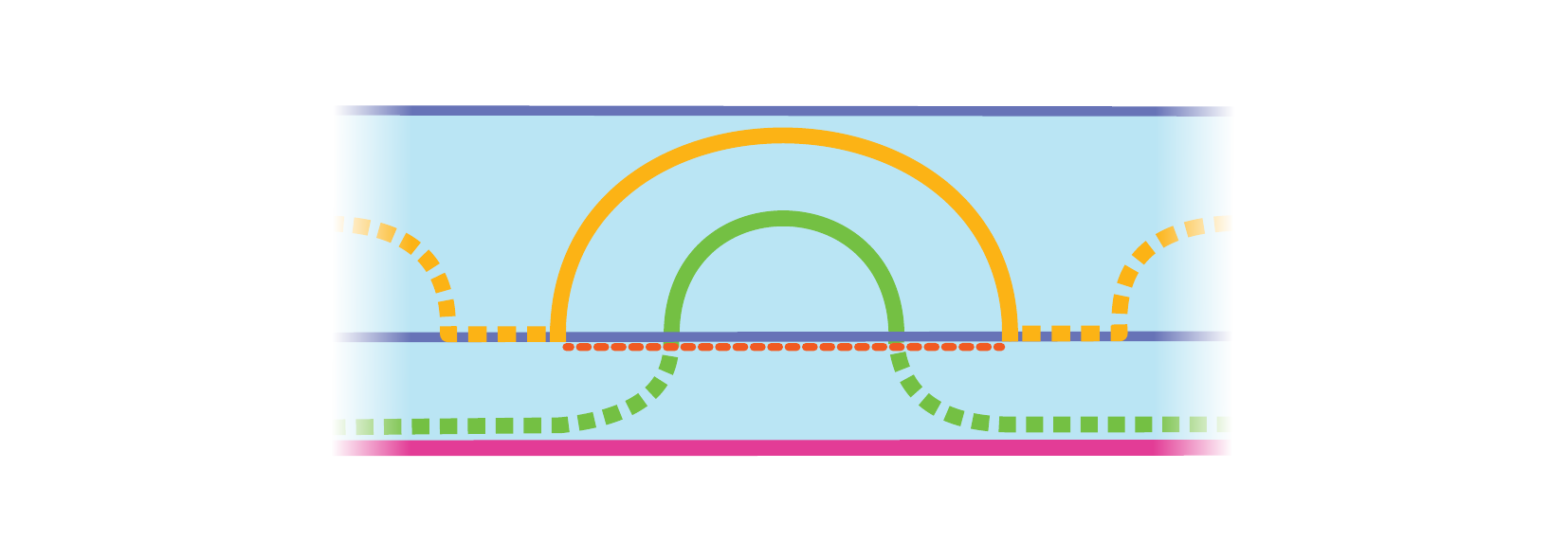}}};
        \end{pgfonlayer}{background}
			
        \begin{pgfonlayer}{main}

            \node (X) [v:ghost,position=0:0mm from C] {};

            \node (C1) [v:ghost,position=19:40mm from X] {$C_1$};
            \node (C2) [v:ghost,position=270:17mm from C1] {$C_2$};
            \node (Omega) [v:ghost,position=270:9mm from C2] {$\Omega$};

            \node (S) [v:ghost,position=270:7.5mm from X] {$S$};
            \node (Q) [v:ghost,position=0:11mm from X] {$Q'$};
            \node (R) [v:ghost,position=35:17mm from X] {$R'$};
        
        \end{pgfonlayer}{main}

        \begin{pgfonlayer}{foreground}
        \end{pgfonlayer}{foreground}

    \end{tikzpicture}
    \caption{Two paths $R'$ and $Q'$ (in \textcolor{BananaYellow}{yellow} and \textcolor{AppleGreen}{green} respectively) interacting with a nest $\mathcal{C}=\{ C_1,C_2\}$ in a way that constitutes a violation of $\mathcal{C}$ being cozy.
    To see this, notice how $Q'$ cannot leave the segment $S$ (depicted in \textcolor{CarrotOrange}{orange}) of $C_2$ without producing a $C_2$-path that sticks our (towards $\Omega$).
    This intuition is captured in \cref{lem:pathsbetweenc1andc2}.}
    \label{fig:ForcingBeaks}
\end{figure}
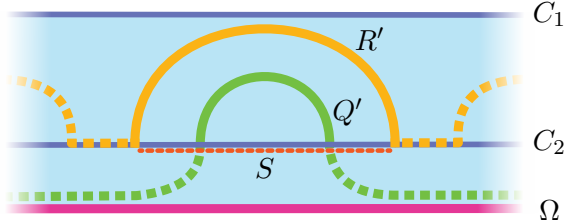

Let us now present a last technical prerequisite for our proof, before we get into the meat of the argument.
The following very technical lemma concerns the drawing of paths in cylindrical renditions with cozy nests.
It captures the intuition that in a cylindrical rendition with a cozy nest, we cannot find two paths that are nested and sit between two cycles of the nest (see \cref{fig:ForcingBeaks} for an illustration).
This is not entirely true, but does hold in most situations.

\begin{lemma}\label{lem:pathsbetweenc1andc2}
    Let $p$ be an integer with $p \geq 2$.
    Let $(G,\Omega)$ be a society with a cylindrical rendition $\rho$ around a vortex $c_0$ with a cozy nest $\mathcal{C} = \{ C_1, C_2 \}$ and let $\mathcal{P}$ be an exposed transaction of order $p$ in $(G,\Omega)$ with the distinct end segments $X,Y$.
    Further, let $X'$ be the set of vertices in $V(\bigcup \mathcal{P}) \cap V(C_2)$ that are seen first when traversing the paths in $\mathcal{P}$ from $X$ to $Y$ and let $Y'$ be defined analogously for traversals from $Y$ to $X$.
    Let $C_X \subseteq C_2$ and $C_Y \subseteq C_2$ be the two minimal, disjoint paths in $C_2$ such that $Z \subseteq V(C_Z)$ for both $Z \in \{ X', Y' \}$.

    Then for any pair of distinct paths $R,Q \in \mathcal{P}$ and any two $C_2$-subpaths $R' \subseteq R$ and $Q' \subseteq Q$ found in the union of $C_1$ and the outer graph of $C_1$ in $\rho$, if both $R'$ and $Q'$ each have at least one endpoint outside of $C_X \cup C_Y$, then the $c_0$-disk of $C_2 \cup R'$ contains the trace of $Q'$ in $\rho$ and vice versa.  
\end{lemma}
\begin{proof}
    Let $R,Q \in \mathcal{P}$ be distinct paths and let $R' \subseteq R$ and $Q' \subseteq Q$ be two $C_2$-subpaths found in the union of $C_1$ and the outer graph of $C_1$ in $\rho$ such that both $R'$ and $Q'$ each have at least one endpoint outside of $C_X \cup C_Y$.
    Let $v \in V(Q') \setminus V(C_X \cup C_Y)$ be one of these endpoints.

    Suppose the $c_0$-disk $d$ of $C_2 \cup R'$ does not contain the trace of $Q'$ in $\rho$.
    Let $C'$ be the cycle in $C_2 \cup R'$ whose $c_0$-disk in $\rho$ is $d$.
    Let $w_X \in V(Q) \cap X'$ and let $w_Y \in V(Q) \cap Y'$.
    
    Suppose that both $w_X$ and $w_Y$ are contained in the unique $C'$-path $L$ in $C_2$ of length at least one, which has the property that $C' \cup L = C_2 \cup R'$.
    As $Q$ is exposed, it must contain an edge of $\sigma(c_0)$.
    To reach $c_0$ from $w_X$ and $w_Y$, the path $Q$ must intersect $C'$, since $C'$ separates the internal vertices of $L$ from $\sigma(c_0)$.
    Because $Q \neq R$, this means that $Q$ must intersect $C' - V(R')$.
    However, this implies the existence of a $C_2$-path in $Q$ that sticks out away from $c_0$, contradicting the coziness of $\mathcal{C}$.

    We may therefore suppose that $w_X$ or $w_Y$ are not found in $L$.
    W.l.o.g.\ we assume that $w_Y$ lies outside of $L$.
    Since the endpoint $v \in V(Q')$ is not contained in $V(C_X \cup C_Y)$ the path $vQw_Y$ must somehow reach $C' - R'$.
    This once more implies the existence of a $C_2$-path in $Q$ that sticks out away from $c_0$, contradicting the coziness of $\mathcal{C}$.
    Thus the $c_0$-disk of $C_2 \cup R'$ must contain the trace of $Q'$ in $\rho$.
    
    The case in which $c_0$-disk of $C_2 \cup Q'$ does not contain the trace of $R'$ in $\rho$ can then be resolved analogously, completing our proof.
\end{proof}

We are now ready for the proof of our main theorem, which we split into three parts.
The first represents the interesting technical core of our argument, involving a quite fascinating discussion of how degenerately the strip of a planar, flat transaction can behave towards a cylindrical rendition if it is not orthogonal to it.
The second part will be the technical version of the statement, we later need in our application.
We mainly state the third part to bring our main theorem into a form resembling \Cref{lem:crookedexistencealgo}.
The proofs of both the second and third part will be comparatively brief and merciful.

\begin{lemma}\label{lem:orthogonalexposed}
    Let $p$ be integers with $p \geq 4$ and let $(G,\Omega)$ be a society with a cylindrical rendition $\rho$ with a cozy nest $\mathcal{C} = \{ C_1, C_2 \}$ around a vortex $c_0$ in the disk $\Delta$.
    Further, let there exist a planar, flat, isolated, separating transaction $\mathcal{P}$ of order $26p + 25$ with a crooked witness $\mathcal{W}$ on $C_2$, such that $V(\Omega)$ consists of the endpoints of $\mathcal{P}$ and $\mathcal{P}$ is flat on one side.
    
    Then there exists a crooked transaction $\mathcal{P}'$ of order $p$ in $(G,\Omega)$ with $\bigcup \mathcal{P}' \subseteq \bigcup \mathcal{P} \cup C_1 \cup C_2$.
    Moreover, there exists an algorithm that finds $\mathcal{P}'$ in time $\mathbf{O}( ps |E(G)| |V(G)|^2 )$.
\end{lemma}
\begin{proof}
    Let $p' \coloneqq 26p+25$ and let $\mathcal{P} = \{ P_1, \ldots , P_{p'} \}$ be indexed naturally.
    As a consequence of $\mathcal{P}$ being flat on one side, we may assume that the labels of the paths in $\mathcal{P}$ are chosen such that the following society $(H, \Psi)$ has a vortex-free rendition $\rho'$ into a disk $\Delta'$.
    Let $I$ be the unique segment of $\Omega$ that has the same endpoints as $P_{p'-1}$ and does not contain the endpoints of $P_{p'}$.
    Further, let $H''$ be the graph consisting of $\mathcal{P} \setminus \{ P_{p'} \}$ and the vertices in $I$, with $H$ being the union of $H''$ and all $H''$-bridges in $G$ with at least one attachment in $V(H'') \setminus V(P_{p'-1})$.
    As noted in \Cref{def:flatononeside}, $H$ must be disjoint from both $P_{p'}$ and the vertices in $V(\Omega) \setminus I$.
    The society $(H, \Psi)$ that we mentioned above is then defined by choosing $\Psi$ such that $V(\Psi) = I \cup (\pi(N(\rho)) \cap V(P_{p'-1}))$, with the order on $I$ being obtained naturally from $\Omega$ and the order of $\pi(N(\rho)) \cap V(P_{p'-1})$ being derived from a traversal of $P_{p'-1}$.

    For each $i \in [p']$ let $x_i,y_i$ be the two distinct endpoints of $P_i$, such that the two end segments $X,Y$ of $\mathcal{P}$ respectively contain all $x_i$s and all $y_i$s for the paths in $\mathcal{P}$.
    Let $x_i'$ be the first vertex in $V(C_2) \cap V(P_i)$ that is seen when traversing $P_i$ starting from $x_i$, with $y_i'$ being defined analogously for a traversal starting from $y_i$.
    Further, let $P_i^1$ be the maximal subpath of $P_i \cap C_2$ containing $x_i'$ and let $P_i^2$ be the maximal subpath of $P_i \cap C_2$ containing $y_i'$.

    \begin{figure}[ht]
    \centering
    \begin{tikzpicture}

        \pgfdeclarelayer{background}
		\pgfdeclarelayer{foreground}
			
		\pgfsetlayers{background,main,foreground}

        \begin{pgfonlayer}{background}
        \node (C) [v:ghost] {{\includegraphics[width=14cm]{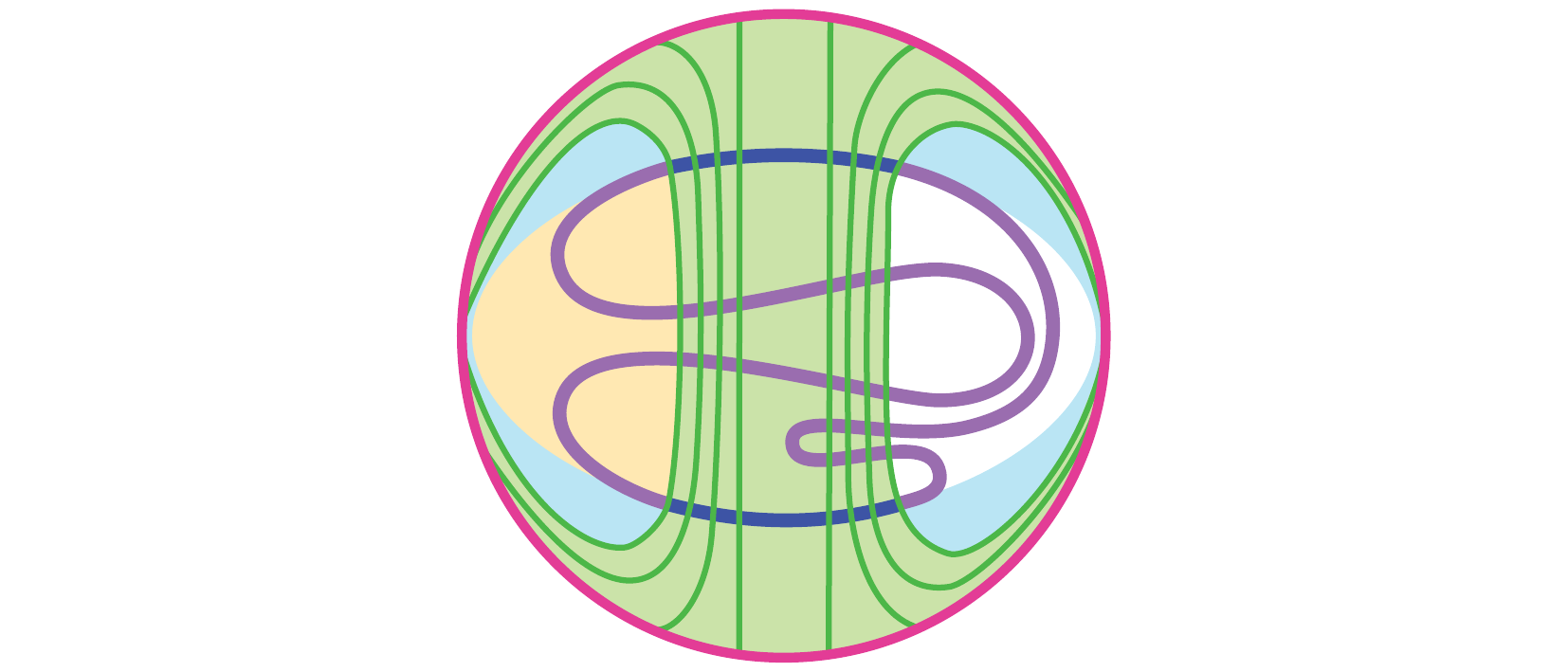}}};
        \end{pgfonlayer}{background}
			
        \begin{pgfonlayer}{main}

            \node (X) [v:ghost,position=0:0.3mm from C] {};

            \node (M) [v:ghost,position=0:0mm from X] {\textcolor{AO}{$\cdots$}};

            \node (N) [v:ghost,position=90:19mm from X] {$C_N$};
            \node (S) [v:ghost,position=270:19.5mm from X] {$C_S$};
            \node (W) [v:ghost,position=180:19mm from X] {$C_W$};
            \node (E) [v:ghost,position=328:21mm from X] {$C_E$};

            \node (Pp) [v:ghost,position=215:16mm from X] {$P_{p'}$};
            \node (P1) [v:ghost,position=0:11.5mm from X] {$P_1$};

        \end{pgfonlayer}{main}

        \begin{pgfonlayer}{foreground}
        \end{pgfonlayer}{foreground}

    \end{tikzpicture}
    \caption{The transaction $\mathcal{P}$ (the \textcolor{AppleGreen}{green} area) in the society $(G,\Omega)$ in the proof of \cref{lem:orthogonalexposed}. Here we depict a rendition $\psi$ of $(G,\Omega)$ such that the restriction of $\psi$ to the strip society of $\mathcal{P}$ is a vortex-free rendition of said society in the disk.
    The subpaths $C_N$ and $C_S$ of $C_2$ are depicted in \textcolor{MidnightBlue}{dark blue} and the subpaths $C_W$ and $C_E$ of $C_2$ are depicted in \textcolor{Amethyst}{purple}.
    The drawing indicates a possibility for the forced behaviour of $C_W$ and $C_E$ due to the presence of the crooked witness $\mathcal{W}$.}
    \label{fig:CardinalCyclePieces}
\end{figure}
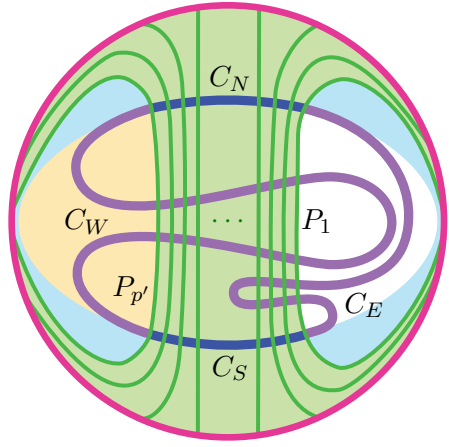
    
    We divide up most of $C_2$ into four sections as follows (see \cref{fig:CardinalCyclePieces} for an illustration):
    \begin{itemize}
        \item Let $C_N$ be the $P_1^1$-$P_{p'}^1$-path in $C_2$ with $P_j^1 \subseteq C_N$ for all $j \in [2,p'-1]$,

        \item let $C_E$ be the $P_{p'}^1$-$P_{p'}^2$-path in $C_2$ with $V(P_j^k) \cap V(C_E) = \emptyset$ for all $j \in [p'-1]$ and $k \in [2]$,

        \item let $C_S$ be the $P_{p'}^2$-$P_1^2$-path in $C_2$ with $P_j^2 \subseteq C_S$ for all $j \in [2,p'-1]$, and

        \item let $C_W$ be the $P_1^2$-$P_1^1$-path in $C_2$ with $V(P_j^k) \cap V(C_W)=\emptyset$ for all $j \in [2,p']$ and $k \in [2]$.
    \end{itemize}
    Note that by construction $H$ contains the two endpoints of $C_W$ but not those of $C_E$.
    The names of these pieces of $C_2$ are chosen to represent cardinal directions to help orient the reader (see also \cref{fig:CardinalCyclePieces}).
    We remark that if $C_2$ was orthogonal to $\mathcal{P}$, neither $C_E$ nor $C_W$ could intersect any paths in $\mathcal{P} \setminus \{ P_1, P_{p'} \}$.
    Furthermore, the endpoints of the paths in $\mathcal{W}$ must be found on $C_E$ and $C_W$.
    In light of these remarks, we claim the following.
    Recall that $p' = 26p + 25$.

    \begin{claim}\label{claim:onesidegoesthrough}
        At least one of the paths $C_E$ and $C_W$ intersects all paths in $\{ P_1, \ldots , P_{13p+12} \}$ or $\{ P_{13p+13}, \ldots , P_{p'-1} \}$.
    \end{claim}
    \emph{Proof of \Cref{claim:onesidegoesthrough}:}
        Depending on whether $\mathcal{W}$ is a snitch or a doublecross, we distinguish two cases.
        We start in the case in which $\mathcal{W}$ is a snitch and thus $\mathcal{W} = \{ P \}$ for some $C_2$-path $P$ found in the inner graph of $C_2$ in $\rho$.
        Recall that $P$ is disjoint from all paths in $\mathcal{P}$ and this also tells us that $P$ has one endpoint in $V(C_W)$ and the other in $V(C_E)$.
        Thus if $P \not\subseteq H$, then it must be disjoint from $H$, as $(H,\Psi)$ has a vortex-free rendition $\rho'$ in a disk $\Delta'$.
        In particular, this means that $C_W$ must intersect all paths in $\mathcal{P} \setminus \{ P_{p'} \}$, which satisfies our claim.

        We may therefore assume that $P \subseteq H$ and thus $P$ has a trace $T$ in $\rho'$.
        Since $P$ is disjoint from the paths in $\mathcal{P}$, the trace of each path in $\mathcal{P} \setminus \{ P_{p'-1}, P_{p'} \}$ in $\rho'$ splits $\Delta'$ into two disks, one of which contains $T$.
        Let $T_i$ be the trace of $P_i$ for $i \in [p' - 2]$ in $\rho'$.
        If $T$ is found in the disk of $\Delta' - T_1$ that does not contain $T_2$, then $C_E$ must intersect all paths in $\mathcal{P}$ to reach $T$ and we are done.
        Analogously, if $T$ is found in the disk of $\Delta' - T_{p'-3}$ that does not contain $T_1$, $C_W$ satisfies our requirements.
        Thus, there instead exists an $i \in [p' - 3]$ such that $T$ is found in the intersection of the disk in $\Delta' - T_i$ that contains $T_{i+1}$ and the disk in $\Delta' - T_{i+1}$ that contains $T_i$.
        As a consequence, to reach this area, either $C_W$ or $C_E$ must intersect at least half of $\mathcal{P} \setminus \{ P_{p'} \}$, as required.

        We may therefore move on to the case in which $\mathcal{W}$ is a doublecross.
        In this case, $\mathcal{W}$ contains two pairs of paths $P,P'$ and $Q,Q'$ such that $P,P'$ have all of their endpoints on $C_W$ and $Q,Q'$ have all of their endpoints on $C_W$, with both pairs forming crosses in the $C_2$-society in $\rho$.
        Both $P$ and $P'$ must intersect $C_1$ and both $C_1$ and $C_2$ intersect all paths in $\mathcal{P} \setminus \{ P_{p'} \}$.
        As a consequence, there exist four disjoint paths from $V(P \cup P')$ to $V(\Omega)$ in $H$.
        Thus there cannot exist a separation $(A,B)$ of order at most three in $H$ with $V(P \cup P') \subseteq A$ and $V(\Psi) \subseteq B$.
        Due to $\rho'$ being a vortex-free rendition, this implies that $V(P \cup P') \cap V(H) = \emptyset$.
        Therefore $C_W$ intersects all paths in $\mathcal{P} \setminus \{ P_{p'} \}$ to escape $H$ and reach $P$ and $P'$, which concludes the proof of our claim.
    \hfill$\blacksquare$
    
    Using \Cref{claim:onesidegoesthrough}, we let $\mathcal{P}' \subseteq \mathcal{P}$ be a transaction of order $q \coloneqq 2p+5$ such that $C_E$ or $C_W$ intersect all paths in $\mathcal{P}'$.
    Solely for the purpose of simplifying notation, we assume that $\mathcal{P}' = \{ P_1, \ldots , P_q \}$.
    We remark that $\bigcup \mathcal{P}' \subseteq H$.
    
    Let $\mathcal{P}_L = \{ P_1, \ldots , P_{12p+8} \}$ and let $\mathcal{P}_R = \{ P_{12p+9} , \ldots , P_{13p + 12} \}$.
    Since $\mathcal{P}$ is flat on one side and itself planar, flat, isolated, and separating, $(G,\Omega, \{ \mathcal{P}' \} )$ is a $(13p+12,1)$-constriction of $(G,\Omega)$.
    In particular, $\mathfrak{N} = (G,\Omega, \{ \mathcal{P}_R \} )$ is a $(p+4,1)$-constriction of $(G,\Omega)$ with the frontier $(G',\Omega')$ such that $\bigcup \mathcal{P}_L$ and $G'$ are disjoint.
    Further, let $G^*$ be the intersection of $\bigcup \mathcal{P}' \cup C_1 \cup C_2 \cup V(\Omega)$ and $G'$, and note that $V(\Omega') \subseteq G^*$, which makes $(G^*, \Omega')$ a society.

    \begin{figure}[ht]
    \centering
    \begin{tikzpicture}

        \pgfdeclarelayer{background}
		\pgfdeclarelayer{foreground}
			
		\pgfsetlayers{background,main,foreground}

        \begin{pgfonlayer}{background}
        \node (C) [v:ghost] {{\includegraphics[width=14cm]{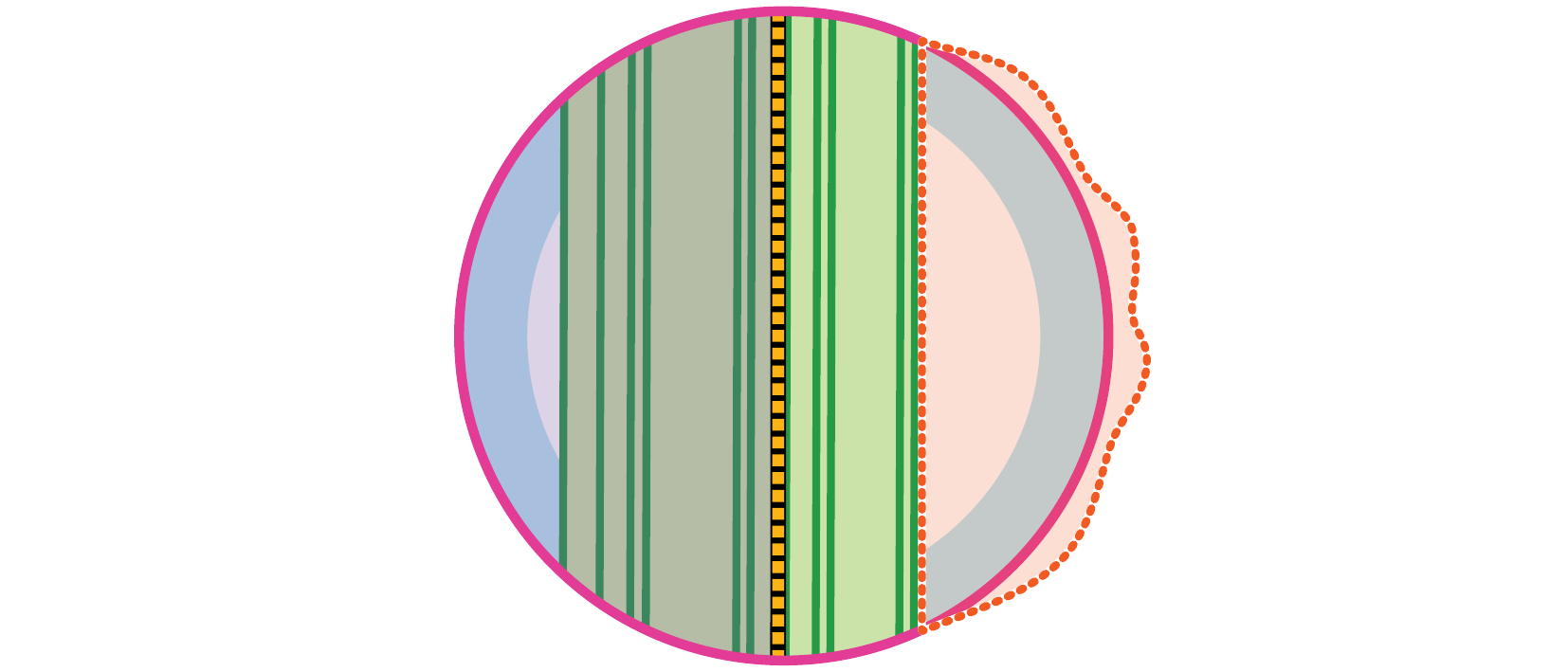}}};
        \end{pgfonlayer}{background}
			
        \begin{pgfonlayer}{main}

            \node (X) [v:ghost,position=0:0mm from C] {};

            \node (DR) [v:ghost,position=25:35mm from X] {$\Delta_R$};
            \node (DL) [v:ghost,position=155:35mm from X] {$\Delta_L$};

            \node (RhoL) [v:ghost,position=200:25.5mm from X] {$\rho_L$};
            \node (RhoR) [v:ghost,position=340:25.5mm from X] {$\rho_R$};

            \node (z) [v:ghost,position=90:31.5mm from X] {$z$};

            \node (c1) [v:ghost,position=0:18mm from X] {$c_1$};

            \node (PL) [v:ghost,position=129:12mm from X] {$\mathcal{P}_L$};
            \node (PR) [v:ghost,position=51:12mm from X] {$\mathcal{P}_R$};

            \node (Ldots) [v:ghost,position=180:8mm from X] {\textcolor{AO}{$\dots$}};
            \node (Rdots) [v:ghost,position=0:7.6mm from X] {\textcolor{AO}{$\dots$}};

            \node (P1) [v:ghost,position=231:32mm from X] {$P_1$};
            \node (P2) [v:ghost,position=240:32mm from X] {$P_2$};
            \node (P12) [v:ghost,position=270:32mm from X] {$P_{12p+10}$};
            \node (P12) [v:ghost,position=300:34mm from X] {$P_{13p+12}$};
        
        \end{pgfonlayer}{main}

        \begin{pgfonlayer}{foreground}
        \end{pgfonlayer}{foreground}

    \end{tikzpicture}
    \caption{Two disks $
    \Delta_L$ (on the left of $z$ and slightly \textcolor{BrilliantRose}{purple}) and $\Delta_R$ (on the right of $z$, the union of the \textcolor{AppleGreen}{green} area and the \textcolor{CarrotOrange}{orange} area).
    The society $(G,\Omega)$ is split into two pieces by the path $P_{12p+10}$ along its trace (the curve $z$).
    The left piece has a vortex-free rendition $\rho_L$ into $\Delta_L$ and the right piece can be interpreted as a constricted society with a cylindrical rendition $\rho_R$ and a unique vortex $c_0$ which contains part of the original infrastructure.}
    \label{fig:ConstrictingTheCrooked}
\end{figure}
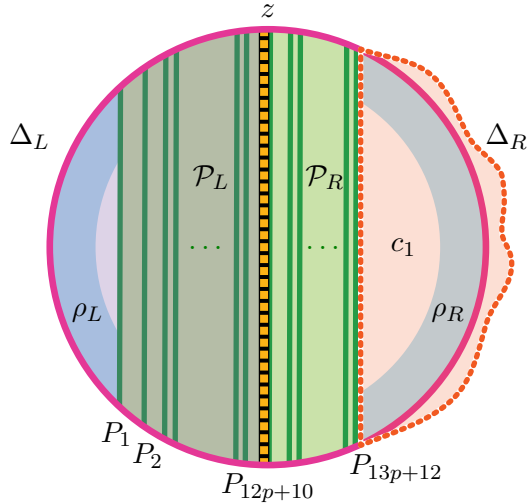

    Let $\rho_\mathfrak{N}$ be the total rendition for $\mathfrak{N}$ with the unique vortex $c_1$ in a disk $\Delta_\mathfrak{N}$.
    Note that, according to \Cref{def:constriction}, the graph $\sigma_{\rho_\mathfrak{N}}(c_1)$ is disjoint from $\bigcup ( \mathcal{P}' \setminus \{ P_{13p+12} \} )$.
    Following \Cref{def:constriction}, the trace of $P_{12p + 10}$ in $\rho_\mathfrak{N}$ splits $\Delta_\mathfrak{N}$ into two disks $\Delta_L$ and $\Delta_R$ such that the restriction of $\rho_\mathfrak{N}$ to $\Delta_R$ yields a cylindrical rendition $\rho_R$ of the society $(G^*,\Omega')$ in $\Delta_R$.
    We let $(H', \Psi')$ be the society resulting from the restriction $\rho_L$ of $\rho_\mathfrak{N}$ to $\Delta_L$ and let $Z = \pi(N(\rho_\mathfrak{N})) \cap V(P_{12p + 10})$.
    Note that $V(\Psi') \subseteq V(\Omega) \cup Z$.
    In fact, we have $V(\Psi') \cap V(\Omega') = Z$, $\bigcup \mathcal{P}_L \subseteq H'$, and $\bigcup (\mathcal{P}_R \setminus \{ P_{12p+9}, P_{12p+10} \} ) \subseteq G'$ (see \cref{fig:ConstrictingTheCrooked} for a visual guide to our setting).
    
    The intersection of $C_E$ or $C_W$ with $\mathcal{P}'$ guaranteed by \Cref{claim:onesidegoesthrough} leads us to a more structured object that describes a separating curve in $\rho_L$.

    \begin{claim}\label{claim:c2subpath}
        There exists a $Z$-path $J$ in $C_2 \cap H'$ with $V(J) \cap V(P_1) \neq \emptyset$ and there exists a $Y \in \{ N, S \}$ such that $J$ is disjoint from $C_Y$. 
    \end{claim}
    \emph{Proof of \Cref{claim:c2subpath}:}
        Let $C_X$ with $X \in \{ E, W \}$ be the subpath of $C_2$ that intersects all paths in $\mathcal{Q}$ and let $v \in V(P_1) \cap C_X$.
        We let $J$ be the $Z$-path in $C_2$ that contains $v$.
        
        Note first that since $Z$ separates the paths in $\mathcal{P}_L$ from the paths in $\mathcal{P}_R$ in $G$, we have $J \subseteq H'$.
        Since $C_X$ intersects all paths in $\mathcal{Q}$ and $Z \subseteq V(P_{12p+9})$, there exists a $v$-$Z$-path $J'$ in $J$ with $J' \subseteq C_X$.
        As both $C_N$ and $C_S$ contain vertices of $Z$, and $C_N$ and $C_S$ are separated by $C_E$ and $C_W$ in $C_2$, this implies that there exists a $Y \in \{ N, S \}$ such that $J$ is disjoint from $C_Y$, as desired.
    \hfill$\blacksquare$

    Using \Cref{claim:c2subpath}, we let $J$ be a $Z$-path in $C_2$ that intersects $P_1$ and assume w.l.o.g.\ that $C_S$ does not intersect $J$.
    The trace $T_J$ of $J$ within $\rho_L$ splits the disk $\Delta_L$ into two disks, of which we let $\Delta_J$ be the disk that contains none of the vertices in $V(\Psi')$ in the drawing associated with $\rho_L$.
    Let $(U,\chi)$ be the $\Delta_J$-society in $\rho_L$.
    For each $i \in [12p+8]$, let $v_i$ be the last vertex in $V(P_i) \cap V(J)$ when traversing $P_i$ starting from $x'_i$.
    Further, let $P_i'$ be $P_i \cap U$ for each $i \in [12p+8]$ and note that $P_i'$ may not necessarily be a connected subpath of $P_i$.
    Furthermore, due to $\rho_L$ being a vortex-free rendition, each $E(P_i') \setminus E(C_2)$ is non-empty for $i \in [2,12p+8]$.

    \begin{claim}\label{claim:c1intersection}
        Both $P_i'$ and $v_iP_iy_i'$ intersect $C_1$ for all $i \in [3,12p+7]$ and the $C_1$-path in $P_i$ that contains $v_i$ is internally disjoint from $C_2 - V(J)$.
    \end{claim}
    \emph{Proof of \Cref{claim:c1intersection}:}
        We have $\mathcal{P}' \subseteq \mathcal{P}$, which is an exposed transaction in $\rho$ since it has a crooked witness.
        Thus each path in $\mathcal{P}'$ contains an edge in $\sigma_{\rho}(c_0)$ and is thus exposed in $\rho$ as well.
        Recall that $\{ C_1, C_2 \}$ is a cozy nest in $\rho$.

        We first consider $P_2'$ and $P_{12p+8}'$, which both consist of components which are paths.
        For each $i \in [2,12p+8]$, let $Q_i'$ be the last component of $P_i'$, when traversing $P_i$ starting from $x'_i$, that contains at least one edge.
        Further, let $z_i$ be the endpoint of $Q_i'$ such that $x_2'P_iz_i$ contains $Q_i'$.
        Let $J'$ be the $z_2$-$z_{12p+8}$-path in $J$ and note that $\{ z_i ~\!\colon\!~ i \in [2,12p+8] \} \subseteq V(J)$, thanks to $\rho_L$ being a vortex-free rendition of $(H',\Psi')$.
        
        For all $i \in [3,12p+7]$, let $z_i'$ be the first vertex from $V(C_2)$ that is encountered when traversing $Q_i'$ starting from $z_i$.
        Note that it is possible that $z_i' \not\in V(J)$.
        We let $Q_i'' \coloneqq z_i'Q_i'z_i$ for all $i \in [2,12+8]$ and for the sake of simplicity, we let $R_1 \coloneqq Q_2''$ and $R_2 \coloneqq Q_{12p+8}''$.

        If both $R_1$ and $R_2$ intersect $C_1$, then let $Q''$ be the unique $C_1$-path in $R_1 \cup R_2 \cup C_2$ and note that for all $i \in [3,12p+7]$ the trace of $Q_i'$ must pass through the trace of $Q''$ to reach $z_i'$ if $Q_i'$ does not intersect $C_2$.
        This immediately implies the existence of a $C_2$-path sticking out away from $c_0$ in $\rho$, contradicting the coziness of $C_2$.

        Thus at least one path in $\{ R_1, R_2 \}$ does not intersect $C_1$.
        Suppose that for any $i \in [3,12p+7]$ the path $Q_i''$ does not intersect $C_1$.
        Then we note that $z_i$ and $z_i'$ separate $z_2$ and $z_{12p+8}$ on $C_2$.
        In particular, there exists a $j \in [2]$, such that the $c_0$-disk of $C_2 \cup Q_i''$ does not contain the trace of $R_i$.
        Furthermore, according to our assumptions on $C_S$ not intersecting $J$, we know that $J'$ is disjoint from both $C_S$ and $C_N$.
        Thus, via \Cref{lem:pathsbetweenc1andc2}, we arrive at a contradiction.

        Using analogous methods, we can show that the paths $v_iP_iy_i'$ for $i \in [3,12p+8]$ also behave as desired and allow us to confirm that the $C_1$-path in $P_i$ is internally disjoint from $C_2 - V(J)$.
    \hfill$\blacksquare$
    
    This claim now verifies that the following pair of vertices exists for each $i \in [3,12p+7]$.
    Let $u_i$ be the last vertex in $C_1 \cap P_i'$ encountered when traversing $P_i$ starting from $x_i'$ and let $w_i$ be the first vertex of $C_1 \cap P_i$ encountered when traversing $v_iP_iy_i'$ starting from $v_i$.
    To aid intuition, we note that, when traversing $P_i$ from $x_i'$ to $y_i'$, the vertices $u_i,v_i,w_i$ are encountered in the given order and in particular, $P_i^* \coloneqq u_iP_iw_i$ is both internally disjoint from $C_1$ and entirely disjoint from $C_2 - V(J)$.
    
    Let $A \coloneqq \{ u_i ~\!\colon\!~ i \in [3,12p+7] \}$ and $B \coloneqq \{ w_i ~\!\colon\!~ i \in [3,12p+7] \}$.
    We will now see that $A$ and $B$ -- which are both found in $V(C_1) \cap V(U)$ -- are in fact highly connected within $C_1$.

    \begin{claim}\label{claim:biglinkageinc1}
        There exists an $A$-$B$-linkage in $C_1$ of order $12p+5$.
    \end{claim}
    \emph{Proof of \Cref{claim:biglinkageinc1}:}
        We start by noting that $|A| = |B| = p-1$ and thus, if the claimed linkage exists, it connects every vertex in $A$ to a vertex in $B$ and vice versa.
        In fact, as we will show, for each $i \in [3,12p+7]$ there exists a $u_i$-$w_i$-path $Q_i \subseteq C_1$ that is disjoint from $(A \cup B) \setminus \{ u_i, w_i \}$.

        Let $i \in [3,12p+7]$ be arbitrary.
        The path $P_i^*$ has the endpoints $u_i$ and $w_i$, intersects $C_2$, and is internally disjoint from $C_1$.
        Thus $P_i^*$ must be found in the union $G_1$ of $C_1$ and the outer graph of $C_1$ in $\rho$.
        Furthermore, due to the coziness of $\mathcal{C}$ and the fact that $\mathcal{P}'$ is orthogonal to $C_3$, $P_i^*$ must also be found in the union of $C_2$ and the inner graph of $C_2$.
        We can split $C_1$ into two internally disjoint $u_i$-$w_i$-paths $R_i$ and $R_i'$ for each $i \in [3,p+1]$ and we choose to label these paths such that the trace of $R_i \cup P_i^*$ does not separate $N_{\rho}(\Omega)$ from $c_0$ in $\rho$, which must be true for exactly one of the two paths $R_i$ and $R_i'$.
        
        Suppose that there exist distinct $i,j \in [3,12p+7]$ such that $R_i$ and $R_j$ intersect.
        Note that $u_i,u_j,w_i,w_j$ are four distinct vertices.
        Thus $R_i$ must contain an endpoint of $R_j$ in its interior.
        W.l.o.g.\ we may thus assume that $u_j \in V(R_i) \setminus \{ u_i, w_i \}$.
        
        Should $w_j \not\in V(R_i)$, we note that $V(P_i^*)$ separates $u_j$ from $w_j$ in $G_1$, since $\rho$ is a cylindrical rendition.
        This contradicts the fact that $P_i^*$ and $P_j^*$ are disjoint.

        We may therefore further assume that $w_j \in V(R_i)$ in addition to $u_j \in V(R_i)$.
        However, in this case, $P_i^*$ separates $C_2 - V(P_i^*)$ from $u_j$ and $w_j$ in $G_1$.
        This in turn contradicts the fact that $v_j \in V(P_j^*) \cap V(C_2)$.
        Therefore, our desired $A$-$B$-linkage is in fact $\{ R_i ~\!\colon\!~ i \in [3,12p+7] \}$.
    \hfill$\blacksquare$

    Before we discuss the linkage this claim yields for us, we turn our attention back to $\rho_\mathfrak{N}$, $J$, $U$, and $G^*$, which were defined above.
    We apply \Cref{lem:crookedexistencealgo} to ask for a crooked transaction of order $p$ in $(G^*,\Omega')$.
    If this returns a crooked transaction, this lets us apply \Cref{lem:crookedrouteback} within $(G^*,\Omega')$ to find a crooked transaction $\mathcal{Q}'$ of order $p$ with all of its endpoints in $V(\Omega') \cap V(\Omega)$.
    Thus we have found a crooked transaction of order $p$ in $(G,\Omega)$
    In particular, due to the fact that we found the linkage constituting the crooked transaction in $G^*$, we know that $\bigcup \mathcal{Q}' \subseteq \bigcup \mathcal{P} \cup C_1 \cup C_2$, as desired.
    It turns out that our application of \Cref{lem:crookedrouteback} is in fact the most computationally costly aspect of our proof.

    Suppose that \Cref{lem:crookedrouteback} instead returns a cylindrical rendition $\rho^*$ of $(G^*,\Omega')$ of depth at most $12p+4$.
    Let $\mathcal{L}$ be the $A$-$B$-linkage of order $12p+5$, whose paths are all contained in $C_1$, that we found in \Cref{claim:biglinkageinc1}.
    The existence of $\mathcal{L}$ begs a simple question.
    How exactly is $\mathcal{L}$ drawn within $\rho_\mathfrak{N}$, whilst $\mathcal{L}$ and $C_2$ are disjoint?
    
    Let $A'$ be the last vertices of $V(U)$ that the paths in $\mathcal{L}$ see when starting from $A$ and note that $A' \subseteq Z$.
    Similarly, let $B'$ be the first vertices of $H' - V(U)$ that the paths in $\mathcal{L}$ see when traversing them from $A'$ to $B$.
    We again note that $B' \subseteq Z$.
    Using this perspective, it is easy to see that $\bigcup \mathcal{L} \cap G^*$ contains an $A'$-$B'$-linkage $\mathcal{L}'$ of order $12p+5$.
    Since $A' \cup B' \subseteq Z \subseteq V(\Omega')$, the linkage $\mathcal{L}'$ is also a transaction in $(G^*, \Omega')$.

    The definition of $(H',\Psi')$ and the associated vortex-free rendition $\rho_L$ allows us to combine $\rho^*$ with $\rho_L$ into a cylindrical rendition $\tau$ of $(H^*,\Omega)$, where $H^* = H' \cup G^* \subseteq G$, around a single vortex $c_1$ in $\Delta$.
    We note that $\bigcup \mathcal{L} \subseteq H^*$, $A \cup B \subseteq V(H^*)$, and $C_2 \subseteq H^*$.
    Suppose that $C_2$ is grounded in $\tau$, then the fact that $A$ and $B$ are separated by $V(J)$ in $H'$ in fact extends to $A$ and $B$ being separated by $V(C_2)$ in $H^*$, since $C_2$ is a cycle.
    This is a clear contradiction to the existence of $\mathcal{L}$, which is disjoint from $V(C_2)$.
    Thus $C_2$ has some edge that is drawn in $\sigma_{\tau}(c_1)$.

    As a consequence of this, $J$ is contained in a $\sigma_{\tau}(c_1)$-path $J'$ within $C_2$.
    Note that thanks to \Cref{claim:c1intersection} we know that for all $i \in [12p+5]$ the path $P_i^*$ is disjoint from $C_2 - V(J)$ and thus in particular from $J' - V(J)$.
    Since $\rho_L$ is a restriction of $\tau$, this implies that $J'$ separates $A$ from $B$ in $H^* - (V(\sigma_{\tau}(c_1)) \setminus \pi_{\tau}(\widetilde{c_1}))$.
    As $J'$ meets the boundary of $c_1$ in at least two nodes of $\tau$, we note that this fact implies that each path of $\mathcal{L}$ must also meet the boundary of $c_1$ in $\tau$.
    Hence $J'$ is grounded in $\tau$.
    We may thus let $T'$ be the trace of $J'$.
    Within $\Delta - (T' \cup c_1)$ there exists a unique disk $\Delta_{J'}$ whose closure contains $U$ and in particular $A$ in the drawing provided by $\tau$.
    We let $U'$ be the crop of $H^*$ to $\Delta_{J'}$ in $\tau$.
    
    This allows us to take $S = U' \cap \pi_{\tau}(\widetilde{c_1})$, which defines a segment of the vortex-society $(G_{c_1},\Omega_{c_1})$ of $c_1$ in $\tau$.
    Each path in $\mathcal{L}$ must also meet $V(\Omega_{c_1} \setminus S)$.
    Therefore there exists a transaction of order $12p+5$ in $(G_{c_1},\Omega_{c_1})$ with one of its end segments contained in $S$ and the other in $V(\Omega_{c_1}) \setminus S$.
    This clearly contradicts the depth of the cylindrical rendition $\rho^*$, which is a restriction of $\tau$, and thus we cannot find the rendition $\rho^*$ and must instead have found a crooked transaction when applying \Cref{lem:crookedexistencealgo} earlier.
\end{proof}

We can now apply this lemma to prove that the presence of a crooked witness for a transaction suffices to allow us to find an orthogonal crooked transaction, with some loss in the size of the desired transaction.

\begin{lemma}\label{lem:orthogonalexposedapply}
    Let the function $\mathsf{f}_{\ref{lem:orthogonalexposedapply}} \colon \mathbb{N}^2 \mapsto \mathbb{N}$ be defined such that $\mathsf{f}_{\ref{lem:orthogonalexposedapply}}(x,y) = 56 x y - 56 x + 52 y - 59$ for all positive integers $x,y$.

    Let $p,s$ be integers with $p \geq 4$ and $s \geq 3$ and let $(G,\Omega)$ be a society with a cylindrical rendition $\rho$ with a cozy nest $\mathcal{C} = \{ C_1, \ldots , C_s \}$ of order $s$ around a vortex $c_0$ in the disk $\Delta$.
    Further, let $\mathcal{Q}$ be a transaction of order $\mathsf{f}_{\ref{lem:orthogonalexposedapply}}(p,s)$ in $(G,\Omega)$ with a crooked witness $\mathcal{W}$ on $C_2$.
    
    Then there exists a crooked transaction $\mathcal{P}$ of order $p$ in $(G,\Omega)$ orthogonal to $\{ C_3, \ldots , C_s \}$.
    
    Moreover, there exists an algorithm that finds $\mathcal{P}$ in time $\mathbf{O}( ps |E(G)||V(G)|^2 )$.
\end{lemma}
\begin{proof}
    We refine $\mathcal{Q}$ by applying \Cref{lem:orthogonaltransactionwithwitness} to $\mathcal{Q}$ and $\mathcal{W}$.
    Since $56 p s - 56 p + 52 s - 59 = (56p + 52)(s-2) + (56p+45)$, this yields a transaction $\mathcal{P}_1$ of order $56p + 45$ such that $\mathcal{W}$ is a crooked witness for $\mathcal{P}_1$ on $C_2$ and $\mathcal{P}_1$ is orthogonal to $\mathcal{C}' \coloneqq \{ C_3, C_4, \ldots , C_s \}$.
    At this point an application of \Cref{lem:crookedorplanar} to $\mathcal{P}_1$ either yields a crooked transaction of order $p$ that is orthogonal to $\mathcal{C}'$ -- in which case we are done -- , or we find a planar transaction $\mathcal{P}_2 \subseteq \mathcal{P}_1$ of order $55p + 47$, since $\mathcal{P}_1$ has order $(55p + 47)+p-2$.
    Note that accordingly $\mathcal{P}_2$ is still orthogonal to $\mathcal{C}'$ and $\mathcal{W}$ remains a crooked witness for $\mathcal{P}_2$ on $C_2$.

    Now, we apply \Cref{lem:planartransactionflatting} to $\mathcal{P}_2$.
    Should this yield a crooked transaction $\mathcal{P}_3$ of order $p$, then this transaction is orthogonal to $\mathcal{C}'$, since $\mathcal{P}_3 \subseteq \mathcal{P}_1$, as desired in the second point of our statement.
    We may thus assume that we instead find a planar, flat, isolated, separating transaction $\mathcal{P} \subseteq \mathcal{P}_2$ of order $p' \coloneqq 26p+25$ that is flat on one side, since $\mathcal{P}_2$ has order $55p + 47 = 3p + 52p + 47$.

    Let $\mathcal{P} = \{ P_1, \ldots , P_{p'} \}$ be indexed naturally.
    As a consequence of $\mathcal{P}$ being flat on one side, we may assume that the labels of the paths in $\mathcal{P}$ are chosen such that the following society $(H, \Psi)$ has a vortex-free rendition $\rho'$ into a disk $\Delta'$.
    Let $I$ be the unique segment of $\Omega$ that has the same endpoints as $P_{p'-1}$ and does not contain the endpoints of $P_{p'}$.
    Further, let $H^*$ be the graph consisting of $\mathcal{P} \setminus \{ P_{p'} \}$ and the vertices in $I$, with $H$ being the union of $H^*$ and all $H^*$-bridges in $G$ with at least one attachment in $V(H^*) \setminus V(P_{p'-1})$.
    As noted in \Cref{def:flatononeside}, $H$ must be disjoint from both $P_{p'-1}$ and the vertices in $V(\Omega) \setminus I$.
    The society $(H, \Psi)$ that we mentioned above is then defined by choosing $\Psi$ such that $V(\Psi) = I \cup (\pi(N(\rho)) \cap V(P_{p'-1}))$, with the order on $I$ being obtained naturally from $\Omega$ and the order of $\pi(N(\rho)) \cap V(P_{p'-1})$ being derived from a traversal of $P_{p'-1}$.

    For each $i \in [p']$ let $x_i,y_i$ be the two distinct endpoints of $P_i$, let $x_i'$ be the last vertex in $\pi(N(\rho)) \cap V(C_3 \cap P_i)$ that is found when traversing $P_i$ from $x_i$ before encountering a vertex in $V(\sigma(c_0))$, and let $y_i'$ be defined analogously.
    Let $X' = \{ x_i' ~\!\colon\!~ i \in [p'] \}$, let $Y' = \{ y_i' ~\!\colon\!~ i \in [p'] \}$, let $(G_3,\Omega_3)$ be the $C_3$-society in $\rho$, and note that $X' \cup Y' \subseteq V(\Omega_3)$.
    We let $\Omega'$ be a cyclic ordering of $X' \cup Y'$ derived from the ordering of these vertices in $\Omega_3$.
    Thus both $X'$ and $Y'$ correspond to segments of $\Omega'$.
    Further, let 
    \[ G' \coloneq (G_3 - ( \bigcup_{i=1}^{q} (x_iP_ix_i' \cup y_iP_iy_i') - (X' \cup Y') )) . \]
    By definition $(G',\Omega')$ is a society.
    For each $i \in [p']$, we let $P_i' = x_i'P_iy_i'$ and we let $\mathcal{P}' = \{ P_1', \ldots , P_{p'}' \}$.
    Since $\mathcal{P}$ is orthogonal to $\mathcal{C}'$, we have $\bigcup \mathcal{P}' \subseteq G'$.
    The transaction $\mathcal{P}'$ in fact inherits many properties of $\mathcal{P}$. 

    \begin{claim}\label{claim:flatstrip}
        $\mathcal{P}'$ is planar, flat, isolated, separating, and flat on one side in $(G',\Omega')$.
    \end{claim}
    \emph{Proof of \Cref{claim:flatstrip}:}
        Let $(G_P, \Omega_P)$ be the $\mathcal{P}'$-strip in $(G',\Omega')$.
        The fact that $\mathcal{P}'$ is planar is immediate from its definition.
        We note that due to $\bigcup \mathcal{P}' \subseteq \bigcup \mathcal{P}$ and the fact that $(H,\Psi)$ has a vortex-free rendition, we know that $G_P \subseteq H$.
        Should $(G_P, \Omega_P)$ contain a cross consisting of two paths $P,P'$, these paths must therefore also be found in $H$.
        Due to our choice of $X \cup Y$, this pair of paths can be extended to a cross in $(H,\Psi)$, a contradiction.
        Thus, according to \Cref{prop:TwoPaths}, $(G_P,\Omega_P)$ has a vortex-free rendition, which means that $\mathcal{P}'$ is flat.
        The facts that $\mathcal{P}'$ is isolated, separating, and flat on one side are also consequences of $(H,\Psi)$ having a vortex-free rendition via analogous arguments.
    \hfill$\blacksquare$

    This now allows us to apply \Cref{lem:orthogonalexposed} to $\mathcal{P}'$ to find a crooked transaction $\mathcal{P}''$ of order $p$ whose endpoints lie in $X' \cup Y'$.
    Therefore, once we route back $\mathcal{P}''$ onto vertices in $V(\Omega)$ via the infrastructure in $\bigcup \mathcal{P} \setminus V(G')$, we have found a crooked transaction $\mathcal{Q}''$ of order $p$ on $(G,\Omega)$.
    By construction $\mathcal{Q}''$ is orthogonal to $\{ C_4, \ldots , C_s \}$.
    In addition, as the vertices in $V(C_3 \cap G^*)$ are all of degree 1 in $G^*$, if such a vertex is contained in $\bigcup \mathcal{Q}''$, they must be an endpoint of a path in $\mathcal{P}''$.
    As a consequence $\mathcal{Q}''$ is orthogonal to $\{ C_3, \ldots , C_s \}$.
    This completes our proof.
\end{proof}

We can now apply this lemma to prove the main theorem of this section, which can be seen as a strengthening of \Cref{lem:crookedexistencealgo}.
As mentioned earlier, we do not find an explicit application for this theorem in our later proofs and instead rely on \Cref{lem:orthogonalexposedapply}.

\begin{theorem}\label{thm:orthogonalexposed}
    Let $p,s$ be integers with $p \geq 4$ and $s \geq 3$ and let $(G,\Omega)$ be a society with a cylindrical rendition $\rho$ with a cozy nest $\mathcal{C} = \{ C_1, \ldots , C_s \}$ around a vortex $c_0$ in the disk $\Delta$.
    
    Then one of the following holds:
    \begin{enumerate}
        \item $(G,\Omega)$ has a cylindrical rendition $\rho'$ in $\Delta$ of depth at most $2016 p s - 2016 p + 1968 s - 2288$, or
        
        \item there exists a crooked transaction $\mathcal{P}$ of order $p$ in $(G,\Omega)$ orthogonal to $\{ C_3, \ldots , C_s \}$.
    \end{enumerate}
    Moreover, there exists an algorithm that finds $\rho'$ or $\mathcal{P}$ in time $\mathbf{O}( ps |E(G)||V(G)|^2 )$.
\end{theorem}
\begin{proof}
    We start by asking \Cref{lem:crookedexistencealgo} for a crooked transaction of order $168 p s - 168 p + 164 s - 191$ in $(G,\Omega)$.
    If this yields a cylindrical rendition of $(G,\Omega)$ with depth at most
    \[ \mathsf{f}_{\ref{thm:orthogonalexposed}}(p,s) = 2016 p s - 2016 p + 1968 s - 2288 =  12( 168 p s - 168 p + 164 s - 191 ) + 4 , \]
    we are done.
    Thus we can assume that \Cref{lem:crookedexistencealgo} yields a crooked transaction $\mathcal{P}_1$ of order
    \[ 168 p s - 168 p + 164 s - 191 = 8s - 14 + 3((56p + 52)(s-2) + (56p+45)) . \]
    By applying \Cref{lem:crookedwitness}, we find a transaction $\mathcal{P}_2 \subseteq \mathcal{P}_1$ of order $(56p + 52)(s-2) + (56p+45) = 56sp - 56p + 52s - 59$, together with a crooked witness $\mathcal{W}$ on $C_2$.
    This in turn allows us to apply \Cref{lem:orthogonalexposedapply}, to find the crooked transaction $\mathcal{P}$ we ask for in the statement.
\end{proof}

\section{Orthogonalisation with respect to several societies}\label{sec:multisocietycrooked}
As we will see in \Cref{sec:transactionmeshes}, we actually need \Cref{thm:orthogonalexposed} from \Cref{sec:orthogonalisecrooked} to be applicable in a setting in which we have multiple societies contained in one another.
Let us first define what we mean by this.

Let $s,k$ be a positive integers.
Let $(G_0, \Omega_0)$ be a society and let $\{ (G_i, \Omega_i) \}_{i \in [k]}$ be a family of $k$ societies, such that for all $i \in [k]$ the society $(G_i, \Omega_i)$ has a cylindrical rendition $\rho_i$ in a disk $\Delta_i$ with a cozy nest $\mathcal{C}_i$ of order $s$ around a vortex $c_i$ such that $V(\Omega_0) \cup V(\bigcup_{j \in [k] \setminus \{ j \}} G_j - (\sigma_{\rho_j}(c_j) - \widetilde{c_j}) \subseteq G_i - (\sigma_{\rho_i}(c_i) - \widetilde{c_i})$.
Furthermore, let it be true that $G_i \subseteq G_0$ for all $i \in [k]$.
If these conditions are met, we call $\{ ((G_i,\Omega_i),\rho_i,c_i,\Delta_i,\mathcal{C}_i) \}_{i\in[k]}$ a \emph{partial $(s,k)$-society configuration in $(G_0,\Omega_0)$}.

Should there further exist a cylindrical rendition $\rho_0$ of $(G_0,\Omega_0)$ in a disk $\Delta_0$ with a nest $\mathcal{C}_0$ of order $s$ around a vortex $c_0$, such that $V(\bigcup_{i \in [k]} G_i - (\sigma_{\rho_i}(c_i) - \widetilde{c_i}) \subseteq G_0 - (\sigma_{\rho_0}(c_0) - \widetilde{c_0})$, we call $\{ ((G_i,\Omega_i),\rho_i,c_i,\Delta_i,\mathcal{C}_i) \}_{i\in[0,k]}$ a \emph{full $(s,k)$-society configuration}.

We note that, letting $H_i'$ be the union of $C_1^i$ and the outer graph of $C_1^i$ in $\rho_i$, the graphs $H_i'$ and $H_j'$ are disjoint for all distinct $i,j \in [k]$, if we are dealing with a partial configuration, and for all distinct $i,j \in [0,k]$ in a full configuration.

This now allows us to define a new type of crooked witness.

\begin{definition}[Co-conspirators]\label{def:co-conspirators}
    Let $s,k,h$ be a positive integers with $s \geq h$.
    Let $\{ ((G_i,\Omega_i),\rho_i,c_i,\Delta_i,\mathcal{C}_i) \}_{i\in[k]}$ be a partial $(s,k)$-society configuration in the society $(G_0,\Omega_0)$.
    Further, let $\mathcal{P}$ be a transaction in $(G_0,\Omega_0)$.
    Let $X,Y$ be the end segments of $\mathcal{P}$, with $A,B$ being the two distinct segments of $\Omega_0$ that have one endpoint in $X$ and the other in $Y$.

    Let $Q$ be an $A$-$V(C_h^i)$-path with $i \in [k]$ in $G_0 - V(\bigcup \mathcal{P})$ that is internally disjoint from $V(\Omega_0)$.
    If there exists a pair of $B$-paths $P,P'$ in $G_0 - V(\bigcup \mathcal{P} \cup Q)$ that form a cross in $(G_0, \Omega_0)$, we call $\{ Q, P,P' \}$ \emph{turncoats (for $\mathcal{P}$ in $(G_0,\Omega_0)$ on level $h$ to $i$)}.
    Else, if there exists a $B$-$V(C_h^j)$-path $Q'$ with $j \in [k]$ in $G_0 - V(\bigcup \mathcal{P})$ that is internally disjoint from $V(\Omega_0) \cup V(Q)$, we call $\{ Q, Q' \}$ \emph{deserters (for $\mathcal{P}$ in $(G_0,\Omega_0)$ on level $h$ to $i$ and $j$)}.

    A set of paths $\mathcal{W}$ in $G_0$ are \emph{co-conspirators (for $\mathcal{P}$ in $(G_0,\Omega_0)$ on level $h$)} if $\mathcal{W}$ is a crooked witness for $\mathcal{P}$ in $(G_0,\Omega_0)$ on $(G_0,\Omega_0)$, or $\mathcal{W}$ are turncoats or deserters for $\mathcal{P}$ in $(G_0,\Omega_0)$ (on level $h$).
\end{definition}

The following is another surprisingly easy variant of \Cref{lemma:orthogonal_transaction} that deals with co-conspirators, similar to how \Cref{lem:orthogonaltransactionwithwitness} preserves the existence of a crooked witness.

\begin{lemma}\label{lemma:orthogonal_transaction_co-conspirators}
    Let $s,p,k$ be positive integers with $s \geq 3$.
    Let $\{ ((G_i,\Omega_i),\rho_i,c_i,\Delta_i,\mathcal{C}_i = \{ C_1^i, \ldots , C_s^i \}) \}_{i\in[k]}$ be a partial $(s,k)$-society configuration in the society $(G_0,\Omega_0)$.
    Let $\mathcal{P}$ be a transaction of order $(p+1)(s-2)+p$ in $(G_0,\Omega_0)$ with the end segments $X',Y'$, let $X$ be the set of endpoints of paths in $\mathcal{P}$ found in $X'$, and let $Y$ be defined analogously for $Y'$.
    Further, let $\mathcal{W}$ be co-conspirators for $\mathcal{P}$ in $(G_0,\Omega_0)$ (on level 2) such that the paths in $\mathcal{W}$ are internally disjoint from $V(C_2^i)$ for all $i \in [k]$.

    Then there exists a transaction $\mathcal{Q}$ in $(G_0,\Omega_0)$ that is an $X$-$Y$-linkage of order $p$ and $\bigcup \mathcal{Q}$ is disjoint from both $\bigcup \mathcal{W}$ and the union of $\bigcup_{i=1}^k (H_i^2 - V(C_2^i))$, where $H_i^2$ is the outer graph of $C_2^i$ in $\rho_i$.
    
    In particular, $\mathcal{Q}$ can be found in time $\mathbf{O}( p |E(G)| )$.
\end{lemma}
\begin{proof}[Proof sketch]
    Let $H_i = H_i^2 - V(C_2^i)$ for each $i \in [k]$ and let $G = G_0 - \bigcup_{i=1}^k V(H_i)$.
    We orthogonalise two radial linkages leading to $X$ and $Y$ respectively as in the proof of \Cref{lemma:orthogonal_transaction} and apply \hyperref[prop:mengersthm]{Menger's theorem} within $G$ to ask for an $X$-$Y$-linkage of order $p$.
    Of course if this yields a linkage, we are immediately done and have justified the runtime bound we claim in our statement.

    We may therefore instead assume that we find an $X$-$Y$-separator $S$ of order at most $p-1$ in $G$.
    First let $\widetilde{\mathcal{P}} \subseteq \mathcal{P}$ be a transaction of order at least $(p+1)(s-2)+1$ whose paths avoid $S$.
    Further, let $Z = V(\bigcup \mathcal{W}) \cap \bigcup_{i=1}^k V(C_2^i)$ be the set of endpoints the paths in $\mathcal{W}$ have in $\bigcup_{i=1}^k V(C_2^i)$.
    We note that $|Z| \leq 2$, since the paths in $\mathcal{W}$ are internally disjoint from $V(C_i^2)$ for all $i \in [0,k]$, with the worst-case being assumed by a pair of deserters.
    
    Let $i \in [k]$ be arbitrary.
    If $F$ is a $C_2^i$-path in $G_0$ that sticks out away from $c_i$ in $\rho_i$, we let $F'$ be the unique subpath of $C_1$ such that $F'$ has the same endpoints as $F$ and $F \cup F'$ is the unique cycle whose trace in $\rho_i$ does not separate $c_i$ and $\Omega_i$.
    As in the proof of \Cref{lemma:orthogonal_transaction}, we call $F'$ the \emph{gangplank} of $F$.
    Let $F''$ be the path in $C_2^i$ with the same endpoints as $F$ that is distinct from $F'$.
    We call $F''$ the \emph{gangboard} of $F$ that lies \emph{opposite to the gangplank $F'$}.
    A subpath $F'$ of $C_2^i$ that is a gangplank for some subpath of a path $\widetilde{P} \in \widetilde{\mathcal{P}}$ is also called a \emph{gangplank} of $\widetilde{P}$ and analogously, if $F'$ was a gangboard, it would be called a \emph{gangboard} of $\widetilde{P}$.

    Let $S' = (S \cap \bigcup_{i=1}^k V(C_2^i)) \cup Z$ and note that $|S'| \leq p+1$.
    Suppose there exists some path $\widetilde{P} \in \widetilde{\mathcal{P}}$ such that whenever a gangplank of $\widetilde{P}$ intersects $S'$, the gangboard opposite to it is disjoint from $S'$.
    Then the union of $\widetilde{P} \cap G$ and all of its gangplanks and gangboards that are disjoint from $S'$ forms a connected subgraph containing a vertex of $X$ and a vertex of $Y$.
    This contradicts the fact that $S$ is an $X$-$Y$-separator and thus no such path $\widetilde{P}$ may exist.

    Clearly, no path in $\widetilde{\mathcal{P}}$ can be found in $G - S$ and thus, every path in $\widetilde{\mathcal{P}}$ has a gangplank.
    In particular, thanks to our prior observation, we conclude that every path in $\widetilde{\mathcal{P}}$ has a gangplank that is intersected by $S'$ and for which the gangboard opposing it is also intersected by $S'$.
    As we have $|\widetilde{\mathcal{P}}| \geq (p+1)(s-2)+1$ and $|S'| \leq p+1$, this means that there exists some $u \in S'$ that is contained in the gangplank of $s-1$ paths $\widehat{P}_1, \dots , \widehat{P}_{s-1} \in \widehat{\mathcal{P}}$, with all of these gangplanks being contained in $C_2^i$ for some fixed $i \in [k]$.
    From this point we reach a contradiction via entirely analogous methods to those presented in the proof of \Cref{lemma:orthogonal_transaction} and thus we omit this part.
\end{proof}

We will also need a variant of \Cref{lem:pathsbetweenc1andc2}.
In our current context, the variant we need is actually easier to state than the original and the proof is analogous.
We thus present the following lemma without proof.

\begin{lemma}\label{lem:pathsbetweenc1andc2reprise}
    Let $p$ be an integer with $p \geq 2$.
    Let $(G,\Omega)$ be a society with a cylindrical rendition $\rho$ around a vortex $c_0$ with a cozy nest $\mathcal{C} = \{ C_1, C_2 \}$.
    Further, let $\mathcal{P}$ be a linkage in $G$ with its endpoints outside of the union $H_1$ of $C_1$ and the outer graph of $C_1$ in $\rho$.
    
    Then for any pair of distinct paths $R,Q \in \mathcal{P}$ and any two $C_2$-subpaths $R' \subseteq R$ and $Q' \subseteq Q$ found in $H_1$, the $c_0$-disk of $C_2 \cup R'$ contains the trace of $Q'$ in $\rho$ and vice versa.
\end{lemma}

We can now prove the main theorem, a variant of \Cref{thm:orthogonalexposed}, which we split into two parts.
The first of which, a lemma which we will be able to reuse in the next section, contains the main arguments, though it requires a rather specific setup to apply.
Our overall proof has many similarities to the proof of \Cref{thm:orthogonalexposed}, though it is also made more technical by the addition of another variable, which represents the number of societies we are dealing with.

The following lemma may seem extremely strange, since it all it seemingly does is replace co-conspirators on level 2 with co-conspirators on level 3 at some cost to the size of the transaction.
But this is in fact a massive difference for us.
As we will see in the proof, the fact that our transaction is disjoint from $C_3^i$ means that, if we have co-conspirators on level 3, we can actually say something meaningful about the location of the society in whose nests $C_3^i$ is found.
On the other hand the witness on level 2 tells us basically nothing about the location of this society.

\begin{lemma}\label{lem:orthogonalexposedsorcery}
    Let $p, s, k$ be non-negative integers with $p \geq 4$ and $s \geq 3$.
    Let $\{ ((G_i,\Omega_i),\rho_i,c_i,\Delta_i,\mathcal{C}_i = \{ C_1^i, \ldots , C_s^i \}) \}_{i\in[0,k]}$ be a full $(s,k)$-society configuration.
    Further, let $\mathcal{Q}$ be a transaction in $(G_0,\Omega_0)$ of order $130 k p + 3 p - 80 k + 9$ that is disjoint from the union of $C_3^i$ and the outer graph $G_3^i$ of $C_3^i$ in $\rho_i$ for all $i \in [k]$ and such that $\mathcal{Q}$ has co-conspirators $\mathcal{W}$ (on level 2).
    Finally, assume that $(G_0,\Omega_0)$ contains a crooked transaction of order 53.

    There exists a transaction $\mathcal{P}$ of order $p$, with co-conspirators $\mathcal{W}$ in $(G_0,\Omega_0)$ (on level 3), that is disjoint from the union of $C_3^i$ and the outer graphs $C_3^i$ in $\rho_i$ for all $i \in [k]$.
    
    Moreover, there exists an algorithm that finds $\mathcal{P}$ and $\mathcal{W}$ in time $\mathbf{O}( psk |E(G)||V(G)| )$.
\end{lemma}
\begin{proof}
    Let $q' \coloneqq 130 k p + 3 p - 80 k + 9$ and let $\mathcal{P}_4$ be the transaction in $(G_0,\Omega_0)$ of order $q'+5p-4$ that is disjoint from the union $G_3^i$ of $C_3^i$ and the outer graph $G_3^i$ of $C_3^i$ in $\rho_i$ for all $i \in [k]$ and such that $\mathcal{Q}$ has co-conspirators $\mathcal{W}$ (on level 2).
    First, we apply \Cref{lem:crookedorplanar} to $\mathcal{P}_4$.
    If this yields a crooked transaction of order $5p-4$, \Cref{lem:easycrookedwitness} finds us the desired objects and we are done.
    Thus we may assume that \Cref{lem:crookedorplanar} instead yields a planar transaction $\mathcal{P}_5 \subseteq \mathcal{P}_4$ of order $q'+2$, since $|\mathcal{P}_4| = q' + 5p - 4$.
    
    Next, we apply \Cref{lem:planartransactionflatting} to $\mathcal{P}_5$, noting that $q' + 2 = 130 k p + 3 p - 80 k + 11 = 3(p+4) + 2(k(65p - 40)+1) - 3$.
    Should this result in us finding a crooked transaction of order $p + 4$, then according to \Cref{lem:planartransactionflatting} we in particular find a crooked transaction $\mathcal{P}_5' \subseteq \mathcal{P}_5$ of order $q$ together with a crooked witness $\mathcal{W}'$ for $\mathcal{P}_5'$ on $(G_0,\Omega_0)$.
    If $\mathcal{W}' = \{ P \}$, then, should $P$ be disjoint from $C_3^i$ for all $i \in [k]$, we let $\mathcal{W} = \{ P \}$.
    Otherwise, let $x,y$ be the two endpoints of $P$, let $R_x$ be the unique $x$-$V(\bigcup_{i=1}^k C_3^i)$-path in $P$, let $R_y$ be the unique $y$-$V(\bigcup_{i=1}^k C_3^i)$-path in $P$, and let $\mathcal{W} = \{ R_x, R_y \}$.\footnote{Note that $R_x$ and $R_y$ may share an endpoint, if $P$ only intersects $\bigcup_{i = 1}^k C_2^i$ in a single vertex. This is an acceptable outcome for us.}
    On the other hand, if $\mathcal{W}'$ is a doublecross, then we let $P,P',Q,Q'$ be the four paths in $\mathcal{W}'$ that form a doublecross for $\mathcal{P}_3 \setminus \{ P,P',Q,Q' \}$ such that $P,P'$ and $Q,Q'$ each form a cross in $(G',\Omega')$.
    If $P \cup P'$ is disjoint from $C_3^i$ for all $i \in [k]$, we add $P,P'$ to $\mathcal{W}$ and otherwise, we let $R_P$ be a $(V(P \cup P') \cap V(\Omega'))$-$V(\bigcup_{i=1}^k C_3^i)$-path in $P \cup P'$ and add $R_P$ to $\mathcal{W}$.
    We proceed analogously with $Q \cup Q'$, which results in $\mathcal{W}$ being a doublecross, pair of turncoats, or pair of deserters for $\mathcal{P}_3$ on level 3, satisfying our goals.
    
    Thus, we can assume that \Cref{lem:planartransactionflatting} instead returns a planar, flat, isolated, separating transaction $\mathcal{P}_6 \subseteq \mathcal{P}_5$ in $(G_0,\Omega_0)$ of order $k(6p+3)+1$ that is flat on one side.
    Note that $\mathcal{W}$ remains a set of co-conspirators for $\mathcal{P}_6 \subseteq \mathcal{P}_4$.

    Let $p' \coloneqq  65p - 40$ and let $\mathcal{P}_6 = \{ P_1, \ldots , P_{k(p'-1)+1} \}$ be indexed naturally.
    Our next goal is to show that there exists a subtransaction of $\mathcal{P}_6$ whose strip is completely disjoint from the union $G_3^i$ of $C_3^i$ and the outer graph of $C_3^i$ in $\rho_i$ for all $i \in [k]$.

    \begin{claim}\label{claim:sortingoutforeignnests}
        There exists a transaction $\mathcal{P} \subseteq \mathcal{P}_6$ of order $p'$, with $(G_0^*,\Omega_0^*)$ being the $\mathcal{P}$-strip society in $(G_0,\Omega_0)$, such that $G^*_0$ is disjoint from $G_3^i$ for all $i \in [k]$.
    \end{claim}
    \emph{Proof of \Cref{claim:sortingoutforeignnests}:}
        The paths in $\mathcal{P}_6$ are all disjoint from $G_3^i$ for all $i \in [k]$.
        Thus this is a simple application of the pigeonhole principle, since $\mathcal{P}_6$ has order $k(p' - 1) + 1$.
    \hfill$\blacksquare$

    We let $\mathcal{P} \subseteq \mathcal{P}_6$ be the transaction described in \Cref{claim:sortingoutforeignnests} and assume that $\mathcal{P} = \{ P_1, \ldots , P_{p'} \}$ for the sole sake of easing notation.
    Let $p'' = 5p - 4$ and let $q = 12p'' + 8$, noting that $p' = q+p''+4$.
    As a consequence of $\mathcal{P}$ being planar, flat, isolated and flat on one side, the tuple $(G_0, \Omega_0, \mathcal{P})$ is a $(p',1)$-constriction of $(G_0,\Omega_0)$ and in particular, $\mathfrak{N} = (G_0,\Omega_0,\{ P_{q+1}, \ldots , P_{p'} \}$ is a $(p''+4,1)$-constriction of $(G_0,\Omega_0)$ with the frontier $(H^*, \Omega^*)$.
    Since $(G_0,\Omega_0)$ contains a crooked transaction of order greater than $52$, the society $(H^*,\Omega^*)$ cannot have a cylindrical rendition of depth at most $52$ and thus instead contains a crooked transaction of order 4.
    This in particular implies, using \Cref{lem:crookedrouteback}, that there exists a crooked transaction of order 4 and thus a cross in $(G_0,\Omega_0)$ that is entirely contained in $H^*$ and has its endpoints in $V(\Omega_0) \cap V(\Omega^*)$.

    Let $\rho_\mathfrak{N}$ be the total rendition in a disk $\Delta_\mathfrak{N}$ associated with $\mathfrak{N}$ and let $T_{q+2}$ be the trace of $P_{q+2}$ in $\rho_\mathfrak{N}$.
    There are two disks $\Delta_L$ and $\Delta_R$ in $\Delta_\mathfrak{N} - T_{q+2}$, the closure $\Delta_L$ of one of these two contains the trace of $P_{q+1}$ and the closure $\Delta_R$ of the other contains the trace of $P_{q+3}$.
    Note that the restriction of $\rho_\mathfrak{N}$ to $\Delta_R$ yields the frontier-rendition $\rho_R$ of $(H^*,\Omega^*)$.
    The restriction of $\rho_\mathfrak{N}$ to $\Delta_L$ in turn yields a vortex-free rendition $\rho_L$ of the society $(H,\Psi)$.
    We note that $\bigcup_{i=1}^{q} P_i \subseteq H$ and $\bigcup_{i=q+4}^{p'} P_i \subseteq H^*$.
    If we let $Z = \pi(N(\rho_\mathfrak{N})) \cap V(P_{q+2})$, we can further note that $V(\Psi) \cap V(\Omega^*) = Z$.

    At this point we would like more information on where exactly the graphs $G_3^i$ ended up for all $i \in [k]$.
    Suppose that $G_3^i \subseteq H$ for some $i \in [k]$.
    Then, if we let $(G'', \Omega'')$ be the $\mathcal{P}$-strip society in $(G_0,\Omega_0)$, we have $G_3^i \subseteq H - V(G'')$ thanks to \Cref{claim:sortingoutforeignnests}.
    Since $H$ is connected, there exists a $V(\Omega_0)$-$V(G_3^i)$-path $R$ in $H - (V(G'') - V(P_1))$, which lets us build a pair of turncoats for $\mathcal{P} \setminus \{ P_1, P_{p'-6}, P_{p'-5}, \ldots , P_{p'} \}$ using the cross in $(G',\Omega')$ described above.
    Since $p' - 8 \geq p$, this satisfies the second point of our statement and we are done.

    As a consequence, we must have $G_3^i \subseteq H^*$ for all $i \in [k]$.
    Recall that $\mathcal{W}$ are co-conspirators (on level 2) for $\mathcal{P}_6$ and since $\mathcal{P} \subseteq \mathcal{P}_6$, the set $\mathcal{W}$ also contains co-conspirators for $\mathcal{P}$.
    Since $\mathcal{P}$ is planar, flat, isolated, and separating, $\mathcal{W}$ cannot be a crooked witness for $\mathcal{P}$, since there cannot exist a cross in $(H,\Psi)$ and no overpass can jump over $\mathcal{P}$ to get to $H^*$.
    Thus $\mathcal{W}$ are turncoats or deserters.
    Suppose $\mathcal{W}$ are turncoats, then the existence of a vortex-free rendition of $(H,\Psi)$ still prohibits the existence of a cross in $(H,\Psi)$.
    We conclude that there must exist a fixed $i \in [k]$ such that a $V(\Omega_0)$-$V(\bigcup_{i=1}^k C_2^i)$-path $R \in \mathcal{W}$ exists that is found entirely in $H$.
    
    This is of course a little curious, since we have already ensured that $G_3^i$ lies in $H^*$.
    The consequence of this is that we now find that $C_2^i$ must be getting dragged through the entirety of $\mathcal{P}$ by $R$ and this inevitably results in the existence of a large linkage in $C_1^i$ that connects the two 'sides' of $C_2^i$.
    This then allows us to build a large crooked transaction with the desired properties in a manner that is entirely analogous to our arguments in the proof of \Cref{thm:orthogonalexposed}.
    We repeat this line of argument here for the convenience of the reader, starting with the claim about $C_2^i$ getting pulled through $\mathcal{P}$.

    \begin{claim}\label{claim:foreignc2subpath}
        There exists a $Z$-path $J$ in $C_2 \cap H$ with $V(J) \cap V(P_1) \neq \emptyset$. 
    \end{claim}
    \emph{Proof of \Cref{claim:foreignc2subpath}:}
        This is immediate from the combination of the following facts:
        The path $R$ has an endpoint in the segment of $\Psi$ that has both endpoints in $V(P_1)$.
        Furthermore $R$ is disjoint from $\bigcup \mathcal{P}$ and $R$ has its other endpoint in $V(C_2^i)$ for some $i \in [k]$. 
    \hfill$\blacksquare$

    Using \Cref{claim:foreignc2subpath}, we let $J$ be a $Z$-path in $C_2$ that intersects $P_1$.
    The trace $T_J$ of $J$ within $\rho_L$ splits the disk $\Delta_L$ into two disks, of which we let $\Delta_J$ be the disk that contains none of the vertices in $V(\Omega')$ in the drawing associated with $\rho_L$.
    Let $(U,\chi)$ be the $\Delta_J$-society in $\rho_L$.
    Let $X,Y$ be the two end segments of $\mathcal{P}$ and let $x_P \in X$ and $y_P \in Y$ be the two endpoints for each $P \in \mathcal{P}$.
    
    For each $j \in [12p''+8]$, let $v_j$ be the last vertex in $V(P_j) \cap V(J)$ when traversing $P_j$ starting from $x_{P_j}$.
    Further, let $P_j'$ be $P_j \cap U$ for each $i \in [12p''+8]$ and note that $P_j'$ may not necessarily be connected.
    Furthermore, due to $\rho_L$ being a vortex-free rendition, each $E(P_j') \setminus E(C_2^i)$ is non-empty for $j \in [2,12p''+8]$.

    \begin{claim}\label{claim:foreignc1intersection}
        Both $P_j'$ and $v_jP_jy_{P_j}$ intersect $C_1^i$ for all $j \in [3,12p''+8]$ and the $C_1^i$-path in $P_j$ that contains $v_j$ is internally disjoint from $C_2^i - V(J)$.
    \end{claim}
    \emph{Proof of \Cref{claim:foreignc1intersection}:}
        The proof of this claim is analogous to the proof of \Cref{claim:foreignc1intersection} if one substitutes the use of \Cref{lem:pathsbetweenc1andc2} with \Cref{lem:pathsbetweenc1andc2reprise}.
        Thus we omit it here.
    \hfill$\blacksquare$
    
    This claim now verifies that the following pair of vertices exists for each $j \in [3,12p''+8]$.
    Let $u_j$ be the last vertex in $C_1^i \cap P_j'$ encountered when traversing $P_j$ starting from $x_{P_j}$ and let $w_j$ be the first vertex of $C_1^i \cap P_j$ encountered when traversing $v_jP_jy_{P_j}$ starting from $v_j$.
    To aid intuition, we note that, when traversing $P_j$ from $x_{P_j}$ to $y_{P_j}$, the vertices $u_j,v_j,w_j$ are encountered in the given order and in particular, $P_j^* \coloneqq u_jP_jw_j$ is internally disjoint from $C_1^i$.
    
    Let $A \coloneqq \{ u_j ~\!\colon\!~ j \in [3,12p''+8] \}$ and $B \coloneqq \{ w_j ~\!\colon\!~ j \in [3,12p''+8] \}$.
    We will now see that $A$ and $B$ -- which are both found in $V(C_1^i) \cap V(U)$ -- are in fact highly connected within $C_1^i$.

    \begin{claim}\label{claim:biglinkageinc1reprise}
        There exists an $A$-$B$-linkage in $C_1^i$ of order $12p''+5$.
    \end{claim}
    \emph{Proof of \Cref{claim:biglinkageinc1reprise}:}
        The proof of this claim is completely analogous to the proof of \Cref{claim:biglinkageinc1}.
        Thus we omit it here.
    \hfill$\blacksquare$

    We now turn our attention back to $\rho_\mathfrak{N}$, $J$, and $U$, which were defined above and let $G^* \coloneqq H^* \cap (\bigcup \mathcal{P} \cup C_1^i \cup C_2^i)$.
    Next, we apply \Cref{lem:crookedexistencealgo} to ask for a crooked transaction of order $p''$ in $(G^*,\Omega^*)$.
    If this actually finds us such a crooked transaction, we are done after applying \Cref{lem:crookedrouteback} to it and afterwards applying \Cref{lem:easycrookedwitness}, which yields a transaction of order $p$ with a crooked witness, since $p'' = 5p-4$.

    Thus, we instead find a cylindrical rendition $\rho^*$ of $(G^*,\Omega^*)$ of depth at most $12p''+4$ in a disk.
    Let $\mathcal{L}$ be the $A$-$B$-linkage of order $12p''+5$, whose paths are all contained in $C_1^i$, that we found in \Cref{claim:biglinkageinc1reprise}.

    Let $A'$ be the last vertices of $V(U)$ that the paths in $\mathcal{L}$ see when starting from $A$ and note that $A' \subseteq Z$.
    Similarly, let $B'$ be the first vertices of $H' - V(U)$ that the paths in $\mathcal{L}$ see when traversing them from $A'$ to $B$.
    We again note that $B' \subseteq Z$.
    Using this perspective, it is easy to see that $\bigcup \mathcal{L} \cap G^*$ contains an $A'$-$B'$-linkage $\mathcal{L}'$ of order $12p''+5$.
    Since $A' \cup B' \subseteq Z \subseteq V(\Omega^*)$, the linkage $\mathcal{L}'$ is also a transaction in $(G^*, \Omega^*)$.

    The definition of $(H,\Psi)$ and the associated vortex-free rendition $\rho_L$ allows us to combine $\rho^*$ with $\rho_L$ into a cylindrical rendition $\tau$ of $H^{**} = H \cup G^* \subseteq G$ around a single vortex $c_1$ in $\Delta$.
    We note that $\bigcup \mathcal{L} \subseteq H^{**}$, $A \cup B \subseteq V(H^{**})$, and $C_2 \subseteq H^{**}$.
    Suppose that $C_2$ is grounded in $\tau$, then the fact that $A$ and $B$ are separated by $V(J)$ in $H$ extends to $A$ and $B$ being separated by $V(C_2)$ in $H^{**}$, since $C_2$ is a cycle.
    This is a clear contradiction to the existence of $\mathcal{L}$, which is disjoint from $C_2$.
    Thus $C_2$ has some edge drawn in $\sigma_{\tau}(c_1)$.

    As a consequence of this, $J$ is contained in a $\sigma_{\tau}(c_1)$-path $J'$ within $C_2$.
    Note that thanks to \Cref{claim:foreignc1intersection} we know that for all $i \in [12p''+5]$ the path $P_i^*$ is disjoint from $C_2 - V(J)$ and thus in particular from $J' - V(J)$.
    Since $\rho_L$ is a restriction of $\tau$, this implies that $J'$ separates $A$ from $B$ in $H^{**} - (V(\sigma_{\tau}(c_1)) \setminus \pi_{\tau}(\widetilde{c_1}))$.
    As $J'$ meets the boundary of $c_1$ in at least two nodes of $\tau$, we note that this fact implies that each path of $\mathcal{L}$ must also meet the boundary of $c_1$ in $\tau$.
    Hence $J'$ is grounded in $\tau$.
    We may thus let $T'$ be the trace of $J'$.
    Within $\Delta - (T' \cup c_1)$ there exists a unique disk $\Delta_{J'}$ whose closure contains $U$ and in particular $A$ in the drawing provided by $\tau$.
    We let $U'$ be the crop of $H^*$ to $\Delta_{J'}$ in $\tau$.
    
    This allows us to take $S = U' \cap \pi_{\tau}(\widetilde{c_1})$, which defines a segment of the vortex-society $(G_{c_1},\Omega_{c_1})$ of $c_1$ in $\tau$.
    Each path in $\mathcal{L}$ must also meet $V(\Omega_{c_1} \setminus S)$.
    Therefore there exists a transaction of order $12p''+5$ in $(G_{c_1},\Omega_{c_1})$ with one of its end segments contained in $S$ and the other in $V(\Omega_{c_1}) \setminus S$.
    This clearly contradicts the depth of the cylindrical rendition $\rho^*$, which is a restriction of $\tau$, and thus we cannot find the rendition $\rho^*$ and must instead have found a crooked transaction when applying \Cref{lem:crookedexistencealgo} earlier.
\end{proof}

We can now finally refine our results from \Cref{sec:orthogonalisecrooked} to respect the presence of several other nests and societies.
Note that with only a small amount of effort we could extend the below result to a full analogue of \Cref{thm:orthogonalexposed}.
We chose not to do this since the result is cumbersome to state -- like the other results of this section -- and is not needed for our arguments.

\begin{theorem}\label{thm:orthogonalexposedsorcery}
    Let the function $\mathsf{f}_{\ref{thm:orthogonalexposedsorcery}} \colon \mathbb{N}^3 \mapsto \mathbb{N}$ be defined such that $\mathsf{f}_{\ref{thm:orthogonalexposedsorcery}}(x,y,z) = 650 x y z - 650 x z - 400 y z + 400 z + 15 x y - 15 x + 55 y - 32$ for all positive integers $x,y,z$.

    Let $p,s,k$ be non-negative integers with $p \geq 4$ and $s \geq 3$.
    Let $\{ ((G_i,\Omega_i),\rho_i,c_i,\Delta_i,\mathcal{C}_i = \{ C_1^i, \ldots , C_s^i \}) \}_{i\in[0,k]}$ be a full $(s,k)$-society configuration.
    Further, let $\mathcal{P}_1$ be a crooked transaction in $(G_0,\Omega_0)$ of order $\mathsf{f}_{\ref{thm:orthogonalexposedsorcery}}(p,s,k)$ orthogonal to $\mathcal{C}' \coloneqq \{ C_3^0, \ldots , C_s^0 \}$.
    
    Then there exists a transaction $\mathcal{P}$ of order $p$ with co-conspirators $\mathcal{W}$ in $(G_0,\Omega_0)$ (on level 3) that is orthogonal to $\mathcal{C}'$, and disjoint from the union of $C_3^i$ and the outer graphs $C_3^i$ in $\rho_i$ for all $i \in [k]$.
    
    Moreover, there exists an algorithm that finds $\mathcal{P}$ and $\mathcal{W}$ in time $\mathbf{O}( psk |E(G)||V(G)|^2 )$.
\end{theorem}
\begin{proof}
    Let $q \coloneqq 130 k p + 3 p - 80 k + 9$.
    We can use \Cref{lem:leapordoublecross} to find a crooked transaction $\mathcal{P}_2 \subseteq \mathcal{P}_1$ that is a leap of order $(q+2)(s-2) + q + 9$ or a doublecross of order $(q+2)(s-2) + q + 9$, which is possible since $\mathcal{P}_1$ has order $\mathsf{f}_{\ref{thm:orthogonalexposedsorcery}}(p,s,k) = 5 q s - 5 q + 10 s + 13 = 4((q+2)(s-2) + q + 9 - 2) + (q+2)(s-2) + q + 9 - 4$.
    Since $\mathcal{P}_2 \subseteq \mathcal{P}_1$, it is in particular orthogonal to $\mathcal{C}'$.

    Let $X,Y$ be the two end segments of $\mathcal{P}_2$ in $\Omega_0$ and let $X',Y'$ be the set of endpoints contained in $X$ and $Y$ respectively.
    For each $P \in \mathcal{P}_2$, we let $x_P' \in X'$ and $y_P' \in Y'$ be its endpoints.
    Further, we let $x_P$ be the last vertex of $\pi(N(\rho)) \cap V(C_3^0)$ that is found on $P \in \mathcal{P}_2$ when traversing it from $x_P$ until a vertex of $V(\sigma(c_0))$ is encountered and let $y_P$ be defined analogously.
    We let $X'' \coloneqq \{ x_P ~\!\colon\!~ P \in \mathcal{P}_2 \}$ and $Y'' \coloneqq \{ y_P ~\!\colon\!~ P \in \mathcal{P}_2 \}$.
    Let $(G_3, \Omega_3)$ be the $C_3$-society in $\rho_0$ and note that $X'' \cup Y'' \subseteq V(\Omega_3)$.
    We define the ordering $\Omega'$ with $V(\Omega') = X'' \cup Y''$ by taking the cyclic ordering induced by $\Omega_3$ on $V(\Omega')$, according to which $X''$ and $Y''$ each correspond to segments of $\Omega'$.
    Further, let
    \[ G' \coloneqq G_3 - (\bigcup_{P \in \mathcal{P}_2} (x_P'Px_P \cup y_P'Py_P) - (X'' \cup Y'')) . \]
    By definition $(G',\Omega')$ is a society.
    We let $\mathcal{P}_3 \coloneqq \{ x_PPy_P ~\!\colon\!~ P \in \mathcal{P}_2 \}$ and note that $\mathcal{P}_3$ is a crooked transaction on $(G',\Omega')$ of the same order as $\mathcal{P}_3$, since $\mathcal{P}_3$ is orthogonal to $\mathcal{C}'$.
    In particular, if $\mathcal{P}_2$ is a leap then $\mathcal{P}_3$ is a leap and similarly, if $\mathcal{P}_2$ is a doublecross then $\mathcal{P}_3$ is also a doublecross.
    We note that this certifies that $(G',\Omega')$ contains a crooked transaction of order at least $53$, since $p \geq 4$ and $s \geq 3$.
    This fact will be required later for us to apply \Cref{lem:orthogonalexposedsorcery}.

    If $\mathcal{P}_3$ is a leap, let $P \in \mathcal{P}_3$ be the overpass in $\mathcal{P}_3$ (see \Cref{subsec:normalisingcrooked}).
    Should $P$ be disjoint from the union of $C_2^i$ for all $i \in [k]$, we let $\mathcal{W} = \{ P \}$.
    Otherwise, let $R_x$ be the unique $x_P$-$V(\bigcup_{i=1}^k C_2^i)$-path in $P$, let $R_y$ be the unique $y_P$-$V(\bigcup_{i=1}^k C_2^i)$-path in $P$, and let $\mathcal{W} = \{ R_x, R_y \}$.
    On the other hand, if $\mathcal{P}_3$ is a doublecross, then we let $P,P',Q,Q'$ be the four paths in $\mathcal{P}_3$ that form a doublecross for $\mathcal{P}_3 \setminus \{ P,P',Q,Q' \}$ such that $P,P'$ and $Q,Q'$ each form a cross in $(G',\Omega')$.
    If $P \cup P'$ is disjoint from $C_2^i$ for all $i \in [k]$, we add $P,P'$ to $\mathcal{W}$ and otherwise, we let $R_P$ be a $(V(P \cup P') \cap V(\Omega'))$-$V(\bigcup_{i=1}^k C_2^i)$-path in $P \cup P'$ and add $R_P$ to $\mathcal{W}$.
    We proceed analogously with $Q \cup Q'$, which results in $\mathcal{W}$ being a doublecross, pair of turncoats, or pair of deserters for $\mathcal{P}_3$.
    
    Independent of whether $\mathcal{P}_3$ is a leap or a doublecross, we can thus find a transaction $\mathcal{P}_4' \subseteq \mathcal{P}_3$ of order $(q+2)(s-2) + q$ with co-conspirators $\mathcal{W}$ (on level 2) such that the paths in $\mathcal{W}$ are internally disjoint from $V(\bigcup_{i=1}^k C_2^i)$.
    This allows us to apply \Cref{lemma:orthogonal_transaction_co-conspirators} to $\mathcal{P}_4'$ and $\mathcal{W}$, resulting in a transaction $\mathcal{P}_4$ in $(G',\Omega')$ of order $q$ that is an $X''$-$Y''$-linkage.
    Furthermore, $\bigcup \mathcal{P}_4$ is disjoint from both $\bigcup \mathcal{W}$ and the union of $\bigcup_{i=1}^k (H_i^2 - V(C_2^i))$, where $H_i^2$ is the outer graph of $C_2^i$ in $\rho_i$.

    Thus, since $\mathcal{P}_4$ is a transaction of order $q = 130 k p + 3 p - 80 k + 11$, we may now apply \Cref{lem:orthogonalexposedsorcery} to find a transaction $\mathcal{P}_5$ of order $p$, with co-conspirators $\mathcal{W}'$ in $(G',\Omega')$ (on level 3), that is disjoint from the union of $C_3^i$ and the outer graphs $C_3^i$ in $\rho_i$ for all $i \in [k]$.
    Since $V(\Omega') = X'' \cup Y''$, both $\mathcal{P}_5$ and $\mathcal{W}'$ can now be routed back onto $(G_0,\Omega_0)$ whilst ensuring that the resulting transaction $\mathcal{P}_6$, for which there now exist co-conspirators $\mathcal{W}''$ in $(G_0,\Omega_0)$ (on level 3), is orthogonal to $\{ C_3^0, \ldots , C_s^0 \}$.
    Furthermore, $\mathcal{P}_6$ is disjoint from the union of $C_3^i$ and the outer graph of $C_3^i$ in $\rho_i$ for all $i \in [k]$.
    This completes our proof.
\end{proof}

\section{Transaction meshes}\label{sec:transactionmeshes}
The following definition is inspired by the notion of ``overflow'' that was recently introduced by Gorsky, Kawarabayashi, Kreutzer, and Wiederrecht in \cite{GorskyKKW2024Packing,Gorsky2024Structure} in the context of a GMST-like result for directed graphs that exclude a large fractional packing of directed cycles of even length.
They are able to iteratively refine a society whilst keeping the number of apices under a fixed bound regardless of the number of iterations this process goes through.
We were not able to adapt this approach for our proof whilst keeping the bounds on the functions involved polynomial.
The idea of iteratively refining a decomposition with large strips of ``flat'' transactions to gradually embed more and more of the society in each iteration is nonetheless important to our approach.

Roughly speaking, a transaction mesh in a society \((G, \Omega)\) represents a well-structured vortex-free rendition of a society \((H, \Omega)\) with \(H \subseteq G\) in a disk with holes such that every vertex of \(H\) incident with an edge of \(G\) not in \(H\) is a node on the boundary of some hole. One can also think of a transaction mesh as a variant of a rendition in a disk, where the ``vortices'' \(v_i\) are not associated with distinct graphs \(\sigma(v_i)\), but rather with one graph which is attached to nodes in the sets \(\tilde{v_i}\). Before we dive into the (quite technical) definition, let us briefly describe the infrastructure of a transaction mesh.

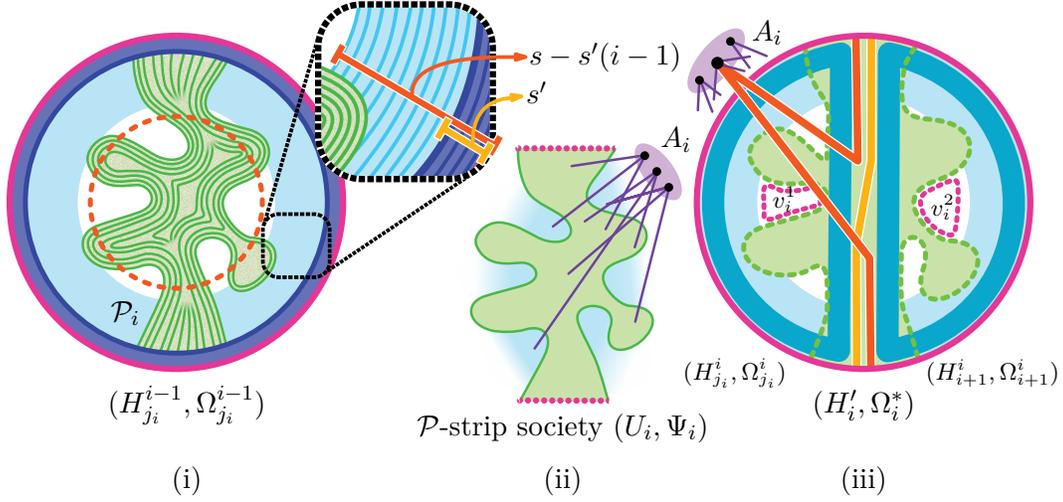
\begin{figure}[ht]
    \centering
    \begin{tikzpicture}

        \pgfdeclarelayer{background}
		\pgfdeclarelayer{foreground}
			
		\pgfsetlayers{background,main,foreground}

        \begin{pgfonlayer}{background}
        \node (C) [v:ghost] {{\includegraphics[width=14cm]{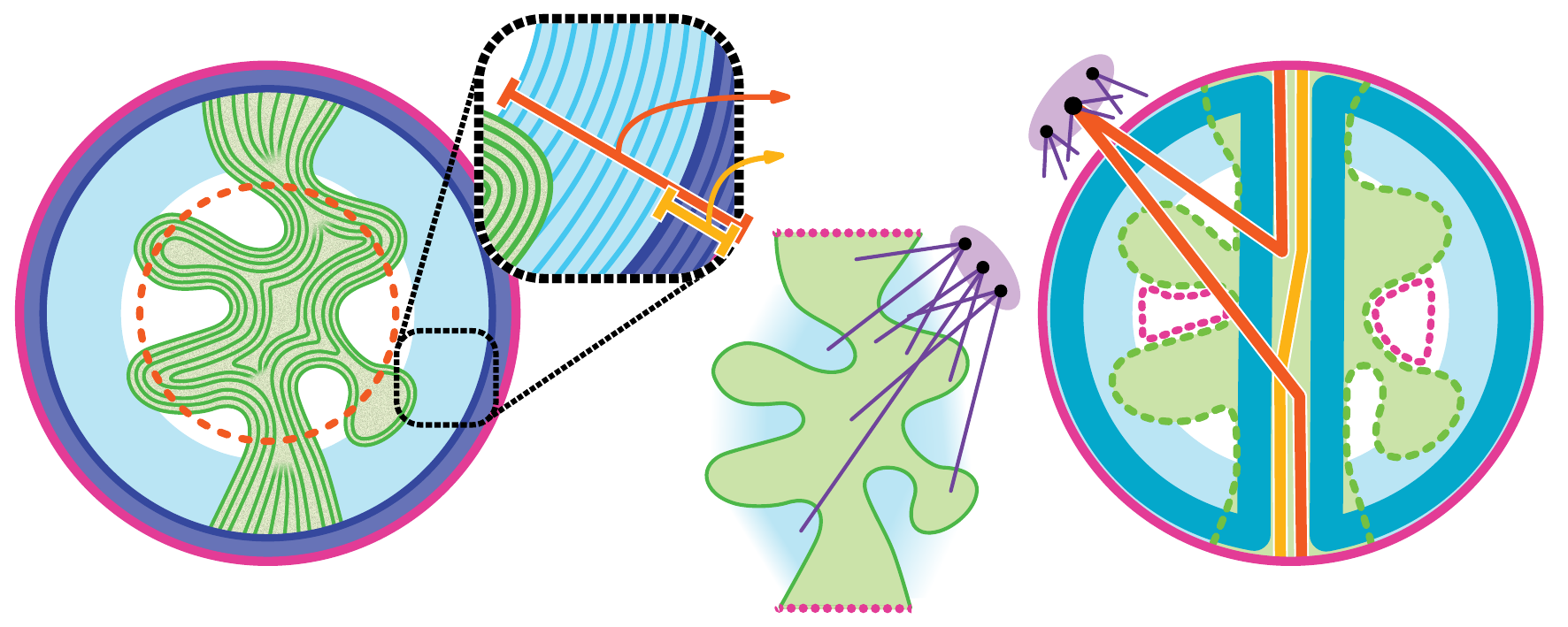}}};
        \end{pgfonlayer}{background}
			
        \begin{pgfonlayer}{main}

            \node (X) [v:ghost,position=0:0.3mm from C] {};

            \node (ss) [v:ghost,position=84:14.7mm from X] {$s'$};
            \node (s) [v:ghost,position=89:4.7mm from ss] {};
            \node (s2) [v:ghost,position=0:9mm from s] {$s-s'(i-1)$};
            
            \node (Bottom) [v:ghost,position=270:33mm from X] {};
            \node (Bottom2) [v:ghost,position=270:4mm from Bottom] {};

            \node (i) [v:ghost,position=180:45mm from Bottom] {};
            \node (ii) [v:ghost,position=0:5mm from Bottom] {};
            \node (iii) [v:ghost,position=0:45mm from Bottom] {};

            \node (iR) [v:ghost,position=180:45mm from Bottom2] {(i)};
            \node (iiR) [v:ghost,position=0:5mm from Bottom2] {(ii)};
            \node (iiiR) [v:ghost,position=0:45mm from Bottom2] {(iii)};

            \node (Pi) [v:ghost,position=114:20mm from i] {$\mathcal{P}_i$};

            \node (Hjii-1) [v:ghost,position=90:6mm from i] {$(H_{j_i}^{i-1},\Omega_{j_i}^{i-1})$};

            \node (Ai2) [v:ghost,position=70:44mm from ii] {$A_i$};

            \node (strip) [v:ghost,position=90:3mm from ii] {$\mathcal{P}$-strip society $(U_i,\Psi_i)$};

            \node (HiP) [v:ghost,position=90:6mm from iii] {$(H'_{i},\Omega_{i}^*)$};
            \node (Ai2) [v:ghost,position=103:57mm from iii] {$A_i$};

            \node (Hiji) [v:ghost,position=149:20mm from iii] {\footnotesize $(H_{j_i}^{i},\Omega_{j_i}^{i})$};
            \node (Hii+1) [v:ghost,position=31:20mm from iii] {\footnotesize $(H_{i+1}^{i},\Omega_{i+1}^{i})$};

            \node (vi1) [v:ghost,position=107.5:35mm from iii] {\footnotesize $v_i^1$};
            \node (vi2) [v:ghost,position=72:34mm from iii] {\footnotesize $v^2_i$};
        
        \end{pgfonlayer}{main}

        \begin{pgfonlayer}{foreground}
        \end{pgfonlayer}{foreground}

    \end{tikzpicture}
    \caption{The construction of the societies \((H_{j_i}^{i}, \Omega_{j_i}^{i})\) and \((H_{i+1}^{i}, \Omega_{i+1}^{i})\) from \((H_{j_i}^{i-1}, \Omega_{j_i}^{i-1})\). (i) We take the \(C_{s-s'i}^{j_i,i-1}\)-society \((H_i^*, \Omega_i^*)\) in \((H_{j_i}^{i-1}, \Omega_{j_i}^{i-1})\), and a planar transaction \(\mathcal{P}_i\) in \((H_i^*, \Omega_i^*)\). As a result the initial order of the nest \(s-s'(i-1)\) decreases by \(s'\). (ii) The \(\mathcal{P}_i\)-society becomes flat after deleting a small set \(A_i\). (iii) We take a rendition \(\tau_i\) of the society \((H_i', \Omega_i^*)\) which captures the nest of \((H_i^*, \Omega_i^*)\) and the (flat) \(\mathcal{P}_i\)-society.
    The rendition has two vortices \(v_i^1\) and \(v_i^2\), which are separated by the traces of the paths in \(\mathcal{P}_i\). If \(A_i \neq \emptyset\), then reintroducing its vertices to the strip creates a cross in the middle of the strip society, which separates the societies
    \((H_{j_i}^{i}, \Omega_{j_i}^i)\) and \((H_{i+1}^{i}, \Omega_{i+1}^i)\).}
    \label{fig:InsideTransactionMesh}
\end{figure}

A transaction mesh parametrised with \(\ell=0\) is simply a cylindrical rendition \(\rho\) of \((G, \Omega)\) around a vortex \(c_0\) with a large nest, and in such a transaction mesh, we denote by \((H_1^0, \Omega_1^0)\) the society where \(H_1^0\) is the part of \(G\) drawn outside the vortex, and \(\Omega_1^0 = \Omega\). In general, a transaction mesh
has societies \((H_1^{\ell}, \Omega_1^{\ell}), \ldots, (H_{\ell+1}^{\ell}, \Omega_{\ell+1}^{\ell})\), and a collection of \(\ell\) transactions.
When the parameter increments from \(\ell-1\) to \(\ell\), we add a transaction \(\mathcal{P}_\ell\) to the mesh. Roughly speaking, for some \(j_\ell \in [\ell]\), the transaction is found in the society obtained from \((H_{j_\ell}^{\ell-1}, \Omega_{j_\ell}^{\ell-1})\) by adding the part of \(G\) which is outside the transaction mesh. The transaction \(\mathcal{P}_\ell\) is planar, exposed, and flat (after deleting a small set \(A_\ell \subseteq V(G)\)).
Upon adding the transaction to the mesh, we split the society \((H_{j_\ell}^{\ell-1}, \Omega_{j_\ell}^{\ell-1})\) into two societies \((H_{j_{\ell}}^{\ell}, \Omega_{\ell}^{\ell})\) and \((H_{\ell+1}^{\ell}, \Omega_{\ell+1}^{\ell-1})\) (more precisely, we find these societies with slightly smaller nests in the union of \((H_{j_\ell}^{\ell-1}, \Omega_{j_\ell}^{\ell-1})\) and the \(\mathcal{P}_j\)-strip society), see \cref{fig:InsideTransactionMesh}. The remaining societies \((H_j^{\ell}, \Omega_j^{\ell})\) with \(j \in [\ell+1] \setminus \{j_\ell, \ell+1\}\) are obtained by slightly shrinking the societies \((H_j^{\ell-1}, \Omega_j^{\ell-1})\).

One further notable feature of the definition is the treatment of the set $A_\ell$ we mentioned earlier.
For technical reasons, we will distinguish between this set being empty -- in which case we somehow got our flat transaction $\mathcal{P}_\ell$ ``for free'' -- and this set containing vertices.
In the second case, we demand that this set in fact is part of a cross that attaches ``in the middle'' of the strip of $\mathcal{P}_\ell$.
Of course, since $\mathcal{P}_\ell$ is flat, this cross disappears as soon as we remove $A_\ell$.
The location of this cross allows us to protect it from the ``vortices'' of the transaction mesh.
This technicality is a herald of things to come, as \Cref{def:Mtree} introduces three more types of objects we track in case $A_\ell$ is empty.
Our reason for doing this is so that we can prove later on that, if we exclude the existence of a $K_t$-minor in $G$, the transaction mesh cannot grow for an unbounded number of iterations.

\begin{definition}[Transaction mesh]\label{def:transactionmesh}
    Let $s,s',a,p,\ell$ be non-negative integers with $s'\ell \leq s$, $a\ell \leq s$, as well as $2s + 2 \leq p$, and let $(G,\Omega)$ be a society.
    A tuple $\mathfrak{M} = (H, \Omega, \mathcal{C}, \mathfrak{P})$ is an \emph{$(s,s',a,p,\ell)$-transaction mesh} in $(G,\Omega)$
    as follows:
    
    The society $(G,\Omega = \Omega_1^0)$ has a cylindrical rendition $\rho$ around the vortex $c_0$ and $\mathcal{C}_1^0 \coloneqq \mathcal{C} = \{ C_1^{1,0}, \ldots , C_s^{1,0} \}$ is a cozy nest of order $s$ in $\rho$, we let $H_1^0$ be the union of $C_1^{1,0}$ and the outer graph of $C_1^{1,0}$ in $\rho$, the set $\mathfrak{P} = \{ \mathcal{P}_1, \ldots , \mathcal{P}_\ell \}$ is a collection of linkages, each of order $p$ and $\mathcal{C} \cup \bigcup_{\mathcal{P} \in \mathfrak{P}} \mathcal{P} \subseteq H$.

    Furthermore, for each $i \in [\ell]$, we require the following:
    \begin{description}
        \item[M1~] There exist societies $(H_1^{i-1}, \Omega_1^{i-1}), \ldots , (H_i^{i-1}, \Omega_i^{i-1})$, each with a cylindrical rendition $\rho_j^{i-1} $ in a disk\footnote{Technically all of these societies have vortex-free renditions, but it helps to have a reference point that defines the inside of the cycle $C_1^{j,i-1}$, which we can easily enforce by simply inventing an otherwise empty vortex that we paste into the middle of our rendition.} around a vortex $c_j^{i-1}$ and cozy nests $\mathcal{C}_1^{i-1}, \ldots , \mathcal{C}_i^{i-1}$ of order $s - s'(i-1)$, and two indices $J_i = \{ j_i, i+1 \} \subseteq [ i+1 ]$ such that for every $j \in [i] \setminus \{ j_i \}$ we have
        $\mathcal{C}_j^i = \{ C_1^{j,i}, \ldots , C_{s - (s'i + 1)}^{j,i} \}$, where $C_k^{j,i} = C_k^{j,i-1}$ for each $k \in [s - (s'i + 1)]$, and $(H^i_j, \Omega^i_j)$ is the $C_{s - s'i}^{j,i-1}$-society in $(H^{i-1}_j, \Omega^{i-1}_j)$.

        \item[M2~] The linkage $\mathcal{P}_i$ is a planar, exposed transaction in $(H^*_i \cup \bigcup \mathcal{P}_i, \Omega^*_i)$, where $(H^*_i, \Omega^*_i)$ is the $C_{s - s'i}^{j_i,i-1}$-society in $\rho_{j_i}^{i-1}$ of $(H^{i-1}_{j_i}, \Omega^{i-1}_{j_i})$, such that $\rho_i^*$ is the restriction of $\rho_{j_i}^{i-1}$ to $H^*_i$.
        
        \item[M3~] There exists a set $A_i \subseteq V(G)$ with $|A_i| \leq a$, we let $G''$ be the component of $G - (\bigcup_{j = 1}^i V(\Omega_j^{i-1}) \cup \bigcup_{j = 1}^i A_j \cup V(\Omega^*_i))$ that contains $C_1^{j_i,i-1}$, let $G' \coloneqq G[G'' \cup \bigcup_{j = 1}^i V(\Omega_j^{i-1}) \cup V(\Omega^*_i)]$, and let $(U_i,\Psi_i)$ be the $\mathcal{P}_i$-strip society of $(G',\Omega^*_i)$.

        There exists a vortex-free rendition of $(U_i,\Psi_i)$ in a disk and $U_i$ is disjoint from $\bigcup_{j \in [i] \setminus \{ j_i \}} \mathcal{C}_j^i$.
        The society $(H_i' \coloneqq (H^*_i \cup U_i) - \bigcup_{j = 1}^i A_j, \Omega^*_i)$ has a rendition $\tau_i$ in a disk $\Delta_i^*$ with exactly two vortices $v^i_1,v^i_2$, such that $\tau_i$ and $\rho_{j_i}^{i-1}$ agree on the rendition of $H^*_i - \bigcup_{j = 1}^i A_j$, we have $v^i_1,v^i_2 \subseteq c_{j_i}^{i-1}$, and for each $P \in \mathcal{P}_i$ the trace of $P$ in $\tau_i$ separates $v^i_1$ and $v^i_2$ in $\Delta_i^*$. 

        \item[M4~] The linkage $\mathcal{P}_i = \{ P_1^i , \ldots , P_p^i \}$ is indexed naturally and taut in $(U_i,\Psi_i)$.
        
        \item[M5~] Let $Q^i_L \coloneqq  P^i_{\lfloor \nicefrac{p}{2} \rfloor}$, and $Q^i_R \coloneqq  P^i_{\lceil \nicefrac{(p+1)}{2} \rceil}$.
        If $A_i \neq \emptyset$, $(H_i'[V(U_i') \cup A_i], \Psi_i')$ has a cross, where $(U_i', \Psi_i')$ is the $\{ Q^i_L, Q^i_R \}$-society in $(U_i, \Psi_i)$,\footnote{Note in particular that, due to \textbf{\textsf{M6}}, this cross must contain some vertex of $A_i$.} and no vertex in $A_i$ has a neighbour in $V(H^*_i)$.
        
        \item[M6~] For the societies $(H^i_{j_i}, \Omega^i_{j_i})$ and $(H^i_{i+1}, \Omega^i_{i+1})$ defined as follows, we have:
        \begin{enumerate}
            \item The graph $H'_{i,j_i}$ is the crop of $H_i'$ to the closure $\Delta_{i,j_i}$ of the disk in $\Delta^*_i - T_L^i$ that contains the trace of $P_1^i$, where $T_L^i$ is the trace of $Q_L^i$ in $\tau_i$.
            
            \item Analogously, $H'_{i,i-1}$ is the crop of $H_i'$ to the closure $\Delta_{i,i+1}$ of the disk in $\Delta^*_i - T_R^i$ that contains the trace of $P_p^i$, where $T_R^i$ is the trace of $Q_R^i$ in $\tau_i$.

            \item For both $j \in J_i$, we let $(H^i_j, \Omega^i_j)$ be the $\Delta_{i,j}$-society in $\tau_i$ and let $\rho_j^i$ be the restriction of $\tau_i$ to $H^i_j$, where the unique vortex $c_j^i$ is one of $v_1^i,v_2^i$.\footnote{Note that this means that $V(\Omega^i_{j_i}) \subseteq V(\Omega^{i-1}_{j_i}) \cup V(Q_L^i)$ and $V(\Omega^i_{i+1}) \subseteq V(\Omega^{i-1}_{j_i}) \cup V(Q_R^i)$.}

            There exists a nest $\mathcal{C}_j = \{ C_1^j, \ldots , C_{s - (s'i + 1)}^j \} \subseteq ( \bigcup ( \{ C_1^{j_i,i-1}, \ldots , C_{s - (s'i + 1)}^{j_i,i-1} \} \cup \mathcal{P}_i ) \cap H'_{i,j}$ of order $s - (s'i + 1)$ in $\rho_j^i$ and another nest $\mathcal{C}^i_j$ of the same order as $\mathcal{C}_j$ in $\rho_j^i$, such that $\bigcup \mathcal{C}^i_j$ is found in the union of $C_1^j$ and the outer graph of $C_1^j$ in $\rho_j^i$.
        \end{enumerate}
    \end{description}
    Finally, $\bigcup_{k=0}^{\ell} ( U_k \cup \bigcup_{j=1}^{k+1} H_j^k ) \subseteq H \subseteq G - \bigcup_{i=1}^\ell A_i$, the society $(H,\Omega)$ has a vortex-free rendition in a disk.
    Furthermore, the society $(H - \bigcup_{i \in [\ell+1]} (H_i^\ell - V(\Omega_i^\ell)), \Omega)$ has a vortex-free rendition $\rho^\dagger$ in a disk in which $\bigcup_{i \in [\ell+1]} V(\Omega_i^\ell) \subseteq N(\rho^\dagger)$, such that for each $i \in [\ell+1]$, $\rho^\dagger$ is a vortex-free rendition of $(H - \bigcup_{j \in [\ell+1] \setminus \{ i \}} (H_j^\ell - V(\Omega_j^\ell)), \Omega_i)$ in a disk.

    For each $i \in [\ell + 1]$ the society $(H_i^\ell,\Omega_i^\ell)$ is called an \emph{area of $\mathfrak{M}$} with the \emph{area-nest} $\mathcal{C}_i^\ell$, and for each $j \in [\ell]$, let $A_j$ be the \emph{apices of $\mathcal{P}_j$ in $\mathfrak{M}$}.
    For all objects we drop the mention of $\mathfrak{M}$ if it is clear from the context.
    We explicitly allow for $\ell = 0$, in which case any society with a cylindrical rendition and a nest of size at least $s$ has an $(s,s',a,p,\ell)$-transaction mesh.
\end{definition}

We remark that for any $i \in [\ell]$ and any two distinct $j,k \in [\ell+1]$ the two graphs $H^i_j$ and $H^i_k$ are disjoint by construction.
This is useful for us, since we ultimately want to completely disconnect the parts of $G$ that the areas describe.
To properly state this, we define the part of $G$ that we consider to be attached to a given area.

\begin{definition}[Districts]\label{def:district}
    Let $s,s',a,p,\ell$ be non-negative integers and let $(G,\Omega)$ be a society, with an $(s,s',a,p,\ell)$-transaction mesh $\mathfrak{M} = (H, \Omega, \mathcal{C}, \mathfrak{P})$.
    Further, let $\{ (H_i, \Omega_i) \}_{i \in [\ell+1]}$ be the areas of $\mathfrak{M}$, with the area-nests $\mathcal{C}_i$ and apices $A_i$.

    For all $i \in [\ell + 1]$, we let $G_i$ be the union of $H_i$ and all $\bigcup \mathcal{C}_i$-bridges in $G - (\bigcup_{j=1}^\ell ( V(\Omega_j) \cup A_j ) \cup V(\Omega_{\ell+1}))$.
    We call $(G_i, \Omega_i)$ the \emph{district corresponding to $(H_i, \Omega_i)$} in $\mathfrak{M}$ and call $\mathcal{C}_i$ the \emph{district-nest of $(G_i, \Omega_i)$}.
    We refer to the collection of all districts corresponding to areas of $\mathfrak{M}$ as the \emph{districts of $\mathfrak{M}$}.
    If it is clear that we are referring to a nest of an area or district, we do not further specify that it is an area- or district-nest.
\end{definition}

Every district $(G_i, \Omega_i)$ of $\mathfrak{M}$ as above has a cylindrical rendition $\rho_i$ around a unique vortex $c_i$ such that the district-nest $\mathcal{C}_i$ is a cozy nest of order $s - s'\ell$ around $c_i$ in $\rho_i$.

One observation that is fundamental to this purpose of transaction meshes is that $(G,\Omega)$ can now be described entirely via $H$ and the districts of $\mathfrak{M}$.

\begin{observation}\label{obs:transactionmeshisgraph}
    Let $s,s',a,p,\ell$ be non-negative integers and let $(G,\Omega)$ be a society, with an $(s,s',a,p,\ell)$-transaction mesh $\mathfrak{M} = (H, \Omega, \mathcal{C}, \mathfrak{P})$.
    Further, let $\{ (G_i, \Omega_i) \}_{i \in [\ell+1]}$ be the districts of $\mathfrak{M}$.

    We have $G = H \cup \bigcup_{i \in [\ell+1]} G_i$.
\end{observation}

\subsection{Padded transaction meshes}\label{sec:paddedtransactionmeshes}
In some parts of our proof, we will need more infrastructure than a nest can provide.
We can reasonably expect this additional infrastructure to exist, since we use transactions to build the transaction mesh that are much larger than what we eventually require as the order of the district-nests.
For this purpose, we introduce the notion of a padding of a nest.

\begin{definition}[Padded societies]\label{def:paddedsocieties}
    Let $s,p,\ell$ be non-negative integers with $1 \leq s \leq p$.
    Let $(G,\Omega)$ be a society with a cylindrical rendition $\rho$ in a disk around a vortex $c_0$ with a nest $\mathcal{C}$ of order $s$.
    Furthermore, let $\mathfrak{P} = \{ \mathcal{P}_1, \ldots , \mathcal{P}_\ell \}$ be a set of unexposed transactions in $(G,\Omega)$ of order $p$ in $\rho$.

    Then we say that $(\mathcal{C}, \mathfrak{P})$ is an \emph{$(s,p,\ell)$-padding for $(G,\Omega)$}, if $\mathcal{C}$ and each $\mathcal{P} \in \mathfrak{P}$ is cozy and for all $i \in [\ell]$, the transaction $\mathcal{P}_i = \{ P_1^i, \ldots , P_p^i \}$ can be indexed naturally such that the $c_0$-disk of $P_1^i$ in $\rho$ is contained in the $c_0$-disk of $P_j^i$ for all $j \in [2,p]$, and for all $j \in [s]$ the path $P^i_{p+j-s}$ intersects $C_j$.
    
    Due to all parts of $(\mathcal{C},\mathfrak{P})$ being grounded, the graph $\bigcup \mathcal{C} \cup \bigcup_{\mathcal{P} \in \mathfrak{P}} \mathcal{P}$ contains a unique $c_0$-disk whose trace is described by a cycle $C$ which we call the \emph{innermost cycle of $(\mathcal{C},\mathfrak{P})$}.
    Note that any given society with a nest of order $s$ has an $(s,q,0)$-padding for any positive integer $q$.
\end{definition}

For our later task of building a transaction mesh, we require versions of two results which we have seen earlier.
First, we want to extend the notion of crooked witnesses to paddings and then present a version of \Cref{lem:crookedwitness}.
The definition needed for this leans heavily on the idea that we can find nests in paddings that sensibly follow the structure of the infrastructure of the padding.

\begin{remark}\label{rem:innernest}
    Let $s,p,\ell$ be non-negative integers with $1 \leq s \leq p$.
    Let $(G,\Omega)$ be a society with a cylindrical rendition $\rho$ in a disk around a vortex $c_0$ with an $(s,p,\ell)$-padding $(\mathcal{C} = \{ C_1, \ldots , C_s\}, \mathfrak{P} = \{ \mathcal{P}_1, \ldots , \mathcal{P}_\ell \})$.
    Further, let $\mathcal{P}_i = \{ P_1^i, \ldots , P_p^i \}$ be indexed naturally such that the $c_0$-disk of $P_p^i$ contains the trace of all other paths in $\mathcal{P}_i$.

    For each $i \in [s]$, we let $C_i'$ be the cycle whose trace bounds the closure of the $c_0$-disk in $C_i \cup \bigcup_{j \in [\ell]} P_i^j$.
    Then $\mathcal{C}' = \{ C_1', \ldots , C_s' \}$ is a nest in $(G,\Omega)$.
\end{remark}

We call $\mathcal{C}'$ from the above remark the \emph{inner nest} of $(\mathcal{C},\mathfrak{P})$, which in particular contains the innermost cycle of $(\mathcal{C},\mathfrak{P})$.
Note that $\mathcal{C}'$ is unlikely to be cozy if $\ell \neq 0$ and especially if $s < p$.
Though in the case in which $\ell = 0$, $\mathcal{C}'$ is simply $\mathcal{C}$ itself.
There is however a society in which most of $\mathcal{C}'$ is cozy and we will work in this society later.
We therefore take a moment to prove this fact.

\begin{lemma}\label{lem:innernestcozyininnersociety}    
    Let $s,p,\ell$ be non-negative integers with $2 \leq s \leq p$.
    Let $(G,\Omega)$ be a society with a cylindrical rendition $\rho$ in a disk around a vortex $c_0$ with an $(s,p,\ell)$-padding $(\mathcal{C}, \mathfrak{P})$, with $\mathcal{C}' = \{ C_1', \ldots , C_s' \}$ being its inner nest.
    Further, let $(G',\Omega')$ be the $C_s$-society in $\rho$ and let $\rho'$ be the restriction of $\rho$ to $G'$.

    Then $\mathcal{C}' \setminus \{ C_s' \}$ is a cozy nest in $\rho'$ around $c_0$.
\end{lemma}
\begin{proof}
    Let $(\mathcal{C} = \{ C_1, \ldots , C_s \}, \mathfrak{P} = \{ \mathcal{P}_1, \ldots , \mathcal{P}_\ell \})$.
    Further, let $\mathcal{P}_i = \{ P_1^i, \ldots , P_p^i \}$ be indexed naturally such that the $c_0$-disk of $P_p^i$ contains the trace of all other paths in $\mathcal{P}_i$.
    
    Suppose that the outer graph of some $C_i'$ in $\rho'$ with $i \in [s-1]$ contains a $C_i'$-path $P$ that sticks out away from $c_0$ and does not intersect any cycle in $(\mathcal{C}' \setminus \{ C_s' \}) \setminus \{ C_i' \}$.
    Then we know by definition that $C_i' \subseteq C_i \cup \bigcup_{j \in [\ell]} P_i^j$ and $C_{i+1}' \subseteq C_{i+1} \cup \bigcup_{j \in [\ell]} P_{i+1}^j$, with the trace of either of these cycles bounding the $c_0$-disk of the subgraph of $G$ they are respectively found within.
    Due to this, we know that $P$ does not intersect $C_{i+1} \cup \bigcup_{j \in [\ell]} P_{i+1}^j$.
    We further note that each cycle and path in $(\mathcal{C}, \mathfrak{P})$ and $\mathcal{C}'$ is grounded.
    Thus $P$ contains a $C$-path $P'$ that sticks out away from $c_0$ for some $C \in \mathcal{C}$ such that $P'$ is disjoint from all cycles in $\mathcal{C} \setminus \{ C \}$, or a $Q$-path $P''$ that sticks out away from $c_0$ for some $P'' \in \mathcal{P} \in \mathfrak{P}$ such that $P''$ is disjoint from all paths in $\mathcal{P} \setminus \{ Q \}$.
    Either option contradicts the coziness of the objects in $(\mathcal{C}, \mathfrak{P})$ and thus $\mathcal{C}'$ must be cozy itself.
\end{proof}

Using the inner nest $\mathcal{C}' = \{ C_1', \ldots , C_s' \}$ of some padding, we can now say that a set $\mathcal{W}$ of paths is a crooked witness for some transaction $\mathcal{P}$ in $(G,\Omega)$ on $C_j'$ for some $j \in [s]$.
This allows us to state our analogue of \Cref{lem:crookedwitness} for this setting.
The proof of this lemma is entirely analogous to the proof of \Cref{lem:crookedwitness}, so we omit it.

\begin{lemma}\label{lem:crookedwitnesspadded}
    Let $s,p,q,\ell$ be non-negative integers with $2 \leq s \leq p$.
    Let $(G,\Omega)$ be a society with a cylindrical rendition $\rho$ in a disk around a vortex $c_0$ with an $(s,p,\ell)$-padding $(\mathcal{C},\mathfrak{P})$, with the inner nest $\mathcal{C}' = \{ C_1', \ldots , C_s' \}$, and let $\mathcal{P}$ be a crooked transaction of order $8p - 14 + 3q$ in $(G,\Omega)$.

    Then there exists an exposed transaction $\mathcal{P}' \subseteq \mathcal{P}$ of order $q$ that has a crooked witness on $C_2'$.

    Moreover, said transaction and its crooked witness can be found in time $\mathbf{O}(q|E(G)|)$.
\end{lemma}

We mildly extend the notion of a radial linkage to a society $(G,\Omega)$ with an $(s,p,\ell)$-padding $(\mathcal{C},\mathfrak{P})$, with $C$ being its innermost cycle, by calling any $V(\Omega)$-$V(C)$-linkage $\mathcal{R}$ a \emph{radial linkage for $(\mathcal{C},\mathfrak{P})$ (in $(G,\Omega)$)}.

The second lemma we need is yet another variant of \Cref{lemma:orthogonal_transaction}.
Here we are more modest than in previous lemmas and only want to ensure that the transaction we find behaves in a very controlled way once it enters into the inner graph of $C_1$.
Note that here we do not orthogonalise the radial linkage we use to construct the linkage in this statement.

\begin{lemma}\label{lemma:orthogonal_transaction_padded}
    Let $s,p,q,k$ be non-negative integers with $2 \leq s \leq p$.
    Let $\{ ((G_i,\Omega_i),\rho_i,c_i,\Delta_i,\mathcal{C}_i = \{ C_1^i, \ldots , C_s^i \}) \}_{i\in[0,k]}$ be a full $(s,k)$-society configuration, such that for each $i \in [0,k]$ the tuple $(\mathcal{C}_i,\mathfrak{P}_i)$ is an $(s,p,\ell_i)$-padding for $(G_i,\Omega_i)$ with $\ell_i \leq \ell$ and the inner nest $\mathcal{C}_i' = \{ C_1^i, \ldots , C_s^i \}$.
    Let $\mathcal{P}$ be a transaction of order $(q+7)(p-2)+q$ in $(G_0,\Omega_0)$ with the end segments $X',Y'$, let $X$ be the set of vertices in $V(C_2^0)$ that are seen first when traversing the paths in $\mathcal{P}$ starting from $X'$, and let $Y$ be defined analogously for $Y'$.
    Further, let $\mathcal{W}$ be a crooked witness for $\mathcal{P}$ on $C_2^0$ in $(G_0,\Omega_0)$.

    Then there exists an $X$-$Y$-linkage $\mathcal{Q}$ of order $q$ in the union of $C_2^0$ and the inner graph of $C_2^0$ in $\rho_i$, such that $\bigcup \mathcal{Q}$ is disjoint from both $\bigcup \mathcal{W}$ and the union of $\bigcup_{i=1}^k (H_i^2 - V(C_2^i))$, where $H_i^2$ is the outer graph of $C_2^i$ in $\rho_i$.
    
    In particular, $\mathcal{Q}$ can be found in time $\mathbf{O}( q |E(G)| )$.
\end{lemma}
\begin{proof}[Proof sketch]
    Let $H_i = H_i^2 - V(C_2^i)$ for each $i \in [0,k]$ and let $G = G_0 - \bigcup_{i=0}^k V(H_i)$.
    We apply \hyperref[prop:mengersthm]{Menger's theorem} within $G$ to ask for an $X$-$Y$-linkage of order $p$.
    Of course if this yields a linkage, we are immediately done and have justified the runtime bound we claim in our statement.

    We may therefore instead assume that we find an $X$-$Y$-separator $S$ of order at most $q-1$ in $G$.
    First let $\widetilde{\mathcal{P}} \subseteq \mathcal{P}$ be a transaction of order at least $(q+7)(p-2)+1$ whose paths avoid $S$.
    Further, let $Z = V(\bigcup \mathcal{W}) \cap \bigcup_{i=0}^k V(C_2^i)$ be the set of endpoints the paths in $\mathcal{W}$ have in $\bigcup_{i=0}^k V(C_2^i)$.
    We note that $|Z| \leq 8$, with the worst-case correspond to a doublecross.
    
    \emph{Gangplanks and gangboards} for any $i \in [0,k]$ and any $C_2^i$-path in $G_0$ that sticks out away from $c_i$ in $\rho_i$ are then defined as in the proof sketch of \Cref{lemma:orthogonal_transaction_co-conspirators}.
    Let $S' = (S \cap \bigcup_{i=0}^k V(C_2^i)) \cup Z$ and note that $|S'| \leq q+7$.
    Suppose there exists some path $\widetilde{P} \in \widetilde{\mathcal{P}}$ such that whenever a gangplank of $\widetilde{P}$ intersects $S'$, the gangboard opposite to it is disjoint from $S'$.
    Then the union of $\widetilde{P} \cap G$ and all of its gangplanks and gangboards that are disjoint from $S'$ forms a connected subgraph containing a vertex of $X$ and a vertex of $Y$.
    This contradicts the fact that $S$ is an $X$-$Y$-separator and thus no such path $\widetilde{P}$ may exist.

    Clearly, no path in $\widetilde{\mathcal{P}}$ can be found in $G - S$ and thus, every path in $\widetilde{\mathcal{P}}$ has a gangplank.
    In particular, thanks to our prior observation, we conclude that every path in $\widetilde{\mathcal{P}}$ has a gangplank that is intersected by $S'$ and for which the gangboard opposing it is also intersected by $S'$.
    As we have $|\widetilde{\mathcal{P}}| \geq (q+7)(p-2)+1$ and $|S'| \leq q+7$, this means that there exists some $u \in S'$ that is contained in the gangplank of $p-1$ paths $\widehat{P}_1, \dots , \widehat{P}_{p-1} \in \widehat{\mathcal{P}}$, with all of these gangplanks being contained in $C_2^i$ for some fixed $i \in [0,k]$.
    From this point we reach a contradiction via entirely analogous methods to those presented in the proof of \Cref{lemma:orthogonal_transaction}, using the fact that all transactions in the padding $(\mathcal{C}_i,\mathfrak{P}_i)$ are cozy, and thus we omit this part.
\end{proof}

Before we start discussing padded societies in the context of transaction meshes, we discuss how to find a slightly less thick padding derived by moving into the $C_{s-k}$-society via the following lemma.

\begin{lemma}\label{lem:narrowing}
    Let $s,p,\ell,k$ be non-negative integers with $1 \leq k < s \leq p$.
    Let $(G,\Omega)$ be a society with an $(s,p,\ell)$-padding $(\mathcal{C} = \{ C_1, \ldots , C_s \},\mathfrak{P})$ in a cylindrical rendition $\rho$ in a disk $\Delta$ around a vortex $c_0$. Then the $C_{s-k}$-society $(G',\Omega')$ in $\rho$ has an $(s-k-1,p-k-1,\ell)$-padding.
\end{lemma}
\begin{proof}
    Let $\mathfrak{P} = \{ \mathcal{P}_1, \ldots , \mathcal{P}_\ell \}$ and $\mathcal{P}_i = \{ P_1^i, \ldots , P_p^i \}$ for all $i \in [\ell]$.
    To show that our lemma holds, we mainly need to show that we preserve the parameter $\ell$ when moving into the $C_{s-k}$-society $(G',\Omega')$ in $\rho$ and that we can preserve the intersection properties of the cycles in $\mathcal{C}$ and the outermost $s - k - 1$ paths in each transaction.
    
    For this purpose it suffices to prove that for all $i \in [\ell]$ and $j \in [p - k - 1]$, there exists a unique $C_{s-k}$-path in $P_j^i$ that sticks out towards $c_0$ in $\rho$.
    We show that this holds by induction over $k$, with $k = 0$ being trivial.
    Suppose the statement does not hold and $P_i^j$ does indeed have two distinct $C_{s-k}$-paths $P,P'$ that stick out towards $c_0$ in $\rho$.
    Consider the $P$-path and the $P'$-path in $C_{s-k}$.
    Due to the coziness of $\mathcal{P}_i$ both of these paths must be intersected by $P_{j-1}^i$, which implies that $j < p$.
    Let $u$ be a vertex of $P_{j-1}^i$ in the $P$-path in $C_{s-k}$ and let $u'$ be a vertex of $P_{j-1}^i$ in the $P'$-path $C_{s-k}$.
    We let $Q$ be the $u$-$u'$-path in $P_{j-1}^i$.
    There must exist an $C_{s-k}$-path that sticks out away from $c_0$ in $Q$.
    For this not to contradict the coziness of $\mathcal{C}$, $Q$ must intersect $C_{s-k-1}$, which implies that $k > 0$.
    But this implies that $P_{j-1}^i$ has more than one $C_{s-k-1}$-path sticking out towards $c_0$, contradicting our induction hypothesis.
\end{proof}

The $(s-k-1,p-k-1,\ell)$-padding is called the \emph{$k$-narrowed padding derived from $(\mathcal{C},\mathfrak{P})$}.

Next, we want to discuss how exactly we hope to find paddings in our areas and districts and what we expect of these objects.
For this sake, we introduce a notion that helps us orient ourselves inside of a transaction mesh.
If a transaction $\mathcal{P}_i$ intersects an area, then at most half of that transaction intersects that area by construction.
We establish what half this is for each transaction.

\begin{definition}[Perspective of an area or district]\label{def:perspective}
    Let $s,s',a,p,\ell$ be non-negative integers and let $(G,\Omega)$ be a society, with an $(s,s',a,p,\ell)$-transaction mesh $\mathfrak{M} = (H, \Omega, \mathcal{C}, \mathfrak{P} = \{ \mathcal{P}_1, \ldots , \mathcal{P}_\ell \})$.
    Further, for all $i \in [\ell + 1]$, let $(H_i, \Omega_i)$ be the areas of $\mathfrak{M}$.

    For each $i \in [\ell + 1]$ and $j \in [\ell]$, we let $h_{i,j}$ be 0 if $E(H_i) \cap E(\bigcup \mathcal{P}_j) = \emptyset$.
    Otherwise, we let $h_{i,j} = 1$ if $E(H_i) \cap E(\bigcup P_1^j) \neq \emptyset$ and we let $h_{i,j} = p$ if $E(H_i) \cap E(P_p^j) \neq \emptyset$.
    By construction exactly one of these two options must hold if $E(H_i) \cap E(\bigcup \mathcal{P}_j) \neq \emptyset$, as we have $V(\Omega_i) \subseteq V(\bigcup \mathcal{C} \cup \bigcup_{k=1}^\ell (Q_L^k \cup Q_R^k))$ and $V(\Omega_i)$ is a separator in $H$, since $(H,\Omega)$ has a vortex-free rendition in a disk.

    We call $(h_{i,j})_{j \in [\ell]}$ the \emph{perspective of $(H_i, \Omega_i)$ in $\mathfrak{M}$} and call $(h_{i,j})_{i \in [\ell +1], j \in [\ell]}$ the \emph{perspectives in $\mathfrak{M}$}.
    The \emph{perspective of a district in $\mathfrak{M}$} is simply the perspective of the area in $\mathfrak{M}$ that it corresponds to.
\end{definition}

Using the notion of perspectives, we observe that each district in a transaction mesh with $\ell \geq 1$ has a large radial linkage derived from the transaction with the highest index that intersects it.

\begin{observation}\label{obs:radial_linkages_for_districts}
    Let $s,s',a,p,q$ be non-negative integers with $s \leq q \leq \lfloor \nicefrac{p}{4} \rfloor$, and let $\ell$ be a positive integer.
    Let $(G,\Omega)$ be a society, with an $(s,s',a,p,\ell)$-transaction mesh $\mathfrak{M} = (H, \Omega, \mathcal{C}, \mathfrak{P} = \{ \mathcal{P}_1, \ldots , \mathcal{P}_\ell \})$.
    Further, we let $(H_i, \Omega_i)$ be an area of $\mathfrak{M}$, let $\mathcal{C}_i$ be the area-nest of order $s$ for $(H_i, \Omega_i)$, and let $(h_{i,j})_{j \in [\ell]}$ be the perspective of $(H_i, \Omega_i)$ in $\mathfrak{M}$.

    Let $j \in [\ell]$ be the maximal such that $h_{i,j} \neq 0$, then there exists an radial linkage $\mathcal{R}_i$ of order $q$ for $\mathcal{C}_i$ such that $\bigcup \mathcal{R}_i \subseteq \bigcup \{ P_1^j, \ldots , P_q^j \}$, if $h_{i,j} = 1$, or $\bigcup \mathcal{R}_i \subseteq \bigcup \{ P_{p + 1 - q}^j, \ldots , P_p^j \}$, if $h_{i,j} = p$.
\end{observation}

For $\ell \geq 1$, we call the radial linkage $\mathcal{R}_i$ the \emph{radial area-linkage of $(H_i, \Omega_i)$ in $\mathfrak{M}$} and transfer this object to the district corresponding to $(H_i, \Omega_i)$ in $\mathfrak{M}$, where we call it the \emph{radial district-linkage}.

In addition to this radial linkage, we also want to say that each district and area has their own padding.
In particular, we want to have something close to an $(s,\lfloor \nicefrac{p}{2} \rfloor,\ell_i)$-padding for each district.
We therefore define a more restricted version of transaction meshes as follows.

\begin{definition}[Padded transaction meshes]\label{def:paddedtransactionmeshes}
    Let $s,s',a,p,\ell$ be non-negative integers with $\max(a,s')\ell \leq s$ and $2s+2-(s'\ell) \leq p$.
    Let $\mathfrak{M} = (H, \Omega, \mathcal{C}, \mathfrak{P} = \{ \mathcal{P}_1, \ldots , \mathcal{P}_\ell \})$ be an $(s,s',a,p,\ell)$-transaction mesh in the society $(G,\Omega)$ with the districts $\{ (G_i,\Omega_i) \}_{i \in [\ell+1]}$ each with a cylindrical rendition $\rho_i$ around a vortex $c_i$ with the district nest $\mathcal{C}_i$, and let $(h_{i,j})_{i \in [\ell +1], j \in [\ell]}$ be the perspectives of $\mathfrak{M}$.
    For all $i \in [\ell]$, let $q_1^i = \lfloor \nicefrac{p}{2} - s'(\ell - i) \rfloor$, let $q_p^i = \lceil \nicefrac{(p+1)}{2} + s'(\ell - i) \rceil$, and for all $i \in [\ell+1]$, let 
    \[ I_i \coloneqq \{ j \in [\ell] ~\!\colon\!~ V(P^j_{q_{h_{i,j}}}) \cap V(\Omega_i) \neq \emptyset \} . \]

    We say that $\mathfrak{M}$ is \emph{padded} if for all $i \in [\ell+1]$, there exists an $(s - (s'\ell + 1),\lfloor \nicefrac{p}{2} - (s'\ell + 1) \rfloor,|I_i|)$-padding $(\mathcal{C}_i,\mathfrak{P}_i = \{ \mathcal{P}_j^i ~\!\colon\!~ j \in I_i \})$ for $(G_i,\Omega_i)$ such that $\mathcal{C}_i$ is the outer nest of $(\mathcal{C}_i,\mathfrak{P}_i)$, and for all $v \in V(\Omega_i) \cap V(\bigcup_{j=1}^\ell P^j_{q_{h_{i,j}}})$ there exists a transaction $\mathcal{P}_j^i \in \mathfrak{P}_i$ such that $v$ lies outside of the $c_0$-disk of all paths $P \in \mathcal{P}_j^i$ in $\rho_i$.
\end{definition}

We note that it is possible to derive from \Cref{def:transactionmesh} that every transaction mesh is in fact already padded by following the construction we have basically already laid out in the above definition.
Later on, we will prove this explicitly.
Note that $\mathfrak{P}_i$ will be non-empty for all districts $(G_i,\Omega_i)$ as along as $\ell \neq 0$ due to the way in which we construct new areas in \textbf{\textsf{M6}} of \Cref{def:transactionmesh}.

\subsection{The tree-like structure of transaction meshes}\label{sec:Mtree}
For the proof of the core theorem of \Cref{sec:societyclassification}, we want to introduce a meta-structure that helps us analyse the development of a transaction mesh.
Note that here we do not need to care about paddings.
To facilitate this definition, we will have to recall some objects from \Cref{def:transactionmesh} that we will use to inductively define the object we are interested in, but will not be explicitly accessed during our later proofs.
Our goal here is to establish four parameters that describe properties of a transaction mesh.
Over the course of \Cref{sec:societyclassification}, we will show that in any situation we encounter we can increase one of these parameters and establish that this will ultimately terminate in one of a selection of good outcomes for us.

Before we introduce the definition, we first want to discuss the role the contents of a particular set which we will call $\mathcal{Z}$ will play.
The objects we choose here will seem odd to the reader, especially because there is no bound on their size despite the fact that we are deleting them.
This is due to the fact that we can easily separate these objects from the rest of the graph, because they are bounded by a separator of small size that is well-connected into the remainder of the graph.
To find these, we will later use the following extension of Menger's theorem that is presented in \cite{KawarabayashiTW2021Quickly} as Lemma 6.8.
The runtime-estimate we add to their statement of the theorem is derived from their proof, which essentially consists of applying an algorithmic version of \Cref{prop:mengersthm} $\mathbf{O}(k)$ times.

\begin{proposition}[Kawarabayashi, Thomas, and Wollan \cite{KawarabayashiTW2021Quickly}]\label{prop:findblob}
    Let $k$ be a non-negative integer, let $G$ be a graph, and let $X,Y$ be two non-empty subsets of $V(G)$.
    Assume that there do not exist $k + 1$ vertex-disjoint $X$-$Y$-paths in $G$.
    Then there exists a separation $(A, B)$ of $G$ with $X \subseteq A$ and $Y \subseteq B$ such that for every $v \in A \cap B$, there exists a path from $v$ to $X$ in $G[A]$ and moreover, exactly one of the following holds:
    \begin{enumerate}
        \item $|A \cap B| = k$ and there exists an $(A \cap B)$-$Y$-linkage of order $k$ in $G[B]$, or

        \item $|A \cap B| < k$ and $A \cap B \subseteq Y$.
    \end{enumerate}
    Moreover, there exists an algorithm that, given $G$, $X$, and $Y$ as input, finds a separation $(A, B)$ corresponding to one of these two options in time $\mathbf{O}(k^2|E(G)|)$.
\end{proposition}

Our use of the set $\mathcal{Z}$ here is inspired by a construction Kawarabayashi et al.\ use in steps 15 to 19 of their proof of Theorem 9.1 in \cite{KawarabayashiTW2021Quickly}.
Unless the reader wishes to read the entirety of \cite{KawarabayashiTW2021Quickly}, we do not recommend consulting this proof for a better understanding of our definition.

If $T$ is a tree, we let $L(T)$ be the set of its leaves and let $I(T) = V(T) \setminus L(T)$ be its \emph{internal vertices}.
Let $w$ be designated as the \emph{root} of $T$, then the \emph{children} of a vertex $p \in V(T)$ are the neighbours of $p$ in $T$ that are separated from $w$ by $p$.
We also call $p$ the \emph{parent} of its children.
The \emph{depth} of a vertex $v$ in a tree $T$ with the root $w$ is the length of the unique $w$-$v$-path in $T$.
In particular, the root itself has depth 0 in $T$.

Let $S$ be a non-empty set.
We say that a set $\mathcal{F} \subseteq 2^S$ \emph{partitions} $S$ if the sets in $\mathcal{F}$ are pairwise disjoint and non-empty, and $\bigcup_{F \in \mathcal{F}} F = S$.

\begin{definition}[$\mathfrak{M}$-tree]\label{def:Mtree}
    Let $s,s',a,p,\ell,k$ be non-negative integers with $\max(a,s')\ell \leq s$ and $2s+2 \leq p$.
    Let $\mathfrak{M} = (H, \Omega, \mathcal{C}, \mathfrak{P})$ be an $(s,s',a,p,\ell)$-transaction mesh in the society $(G,\Omega)$ with the areas $(H_i, \Omega_i)$ and districts $(G_i,\Omega_i)$ for $i \in [\ell+1]$.
    For each $i \in [0,\ell]$, we denote the set of vortices in \(\mathfrak{M}\) as $V_i = \{ v_j^i ~\!\colon\!~ j \in [i+1] \}$
    , set $V^i = \bigcup_{j=0}^i V_j$, and thus have $v_1^0 \in V^i$ for all $i \in [0,\ell]$.

    We will define the contents of the tuple $\mathfrak{T} = ( \mathfrak{M}, T, \phi, \mathcal{A}, \mathcal{Z}, \mathcal{F}, \mathfrak{L} )$, which we call the \emph{$\mathfrak{M}$-tree in $G'$} with \emph{precision $k$}, by induction over $\ell$ as follows.
    In all cases, the graph $T$ is a tree with $V(T) = V^\ell$ rooted in $v_1^0$, the function $\phi$ maps $V(T)$ to a collection of societies that includes the districts of $\mathfrak{M}$, the sets $\mathcal{A}= \{ A_1, \ldots , A_\ell \}$ and $\mathcal{Z} = \{ Z_1, \ldots , Z_\ell \}$ contain subsets of $V(G)$, the set $\mathcal{F}$ partitions $[\ell+1]$, and $\mathfrak{L} = \{ \mathcal{L}_1, \ldots , \mathcal{L}_\ell \}$ is a set of linkages of order $k+1$ in $G$.
    
    For $\ell = 0$, we define $\mathfrak{T}$ such that $T$ is a tree with $V(T) = \{ v_1^0 \}$, the function $\phi$ maps $1$ to the sole district of $\mathfrak{N}$, the sets $\mathcal{A}, \mathcal{Z}, \mathfrak{L}$ are empty, and $\mathcal{F} = \{ \{ 1 \} \}$.

    Let $\mathfrak{M}' = (H, \Omega, \mathcal{C}, \mathfrak{P} \setminus \{ \mathcal{P}_\ell \})$ be the $(s,s',a,p,\ell - 1)$-transaction mesh in $(G,\Omega)$ be derived from $\mathfrak{M}$ with $\ell \geq 1$ such that $(H_i^{\ell - 1}, \Omega_i^{\ell - 1})$ is the area with district $(G_i^{\ell - 1}, \Omega_i^{\ell - 1})$ in $\mathfrak{M}'$ and $H_i \subseteq H_i^{\ell-1}$ for all $i \in [\ell]$.
    We let $\mathcal{C}_i^{\ell-1} = \{ C_1^{i,\ell-1}, \ldots , C_{s - (s'(\ell-1) + 1)}^{i,\ell-1} \}$ be the area-nest of $(H_i^{\ell - 1}, \Omega_i^{\ell - 1})$ for all $i \in [\ell]$.
    Furthermore, let $\mathfrak{T}' = ( \mathfrak{M}', T', \psi, \mathcal{A}', \mathcal{Z}', \mathcal{F}', \mathfrak{L}' )$ be the $\mathfrak{M}'$-tree with precision $k$.
    
    We define $\mathfrak{T}$ such that
    \begin{description}
        \item[~T1~] let $j_\ell \in [\ell]$ be chosen such that $\mathcal{P}_\ell$ is a transaction in the $C_{s - s'\ell}^{j_\ell,\ell-1}$-society of $(G_{j_\ell}^{\ell-1}, \Omega_{j_\ell}^{\ell-1})$, then
            \begin{enumerate}
                \item for each $i \in [\ell] \setminus \{ j_\ell \}$, we let $v_i^\ell$ be the unique child of $v_i^{\ell-1}$ and set $\phi(v_i^\ell) = (G_i, \Omega_i)$, and

                \item we let $v_{j_\ell}^\ell$ and $v_{\ell+1}^\ell$ be the two children of $v_{j_\ell}^{\ell-1}$ and set $\phi(v_{j_\ell}^\ell) = ( G_{j_\ell}, \Omega_{j_\ell} )$, as well as $\phi(v_{\ell+1}^\ell) = ( G_{\ell+1} , \Omega_{\ell+1} )$,\footnote{Accordingly, the vertices in $V_i$ have depth exactly $i$ in $T$ for each $i \in [0,\ell]$.}
            \end{enumerate}
            
        \item[~T2~] we let $\mathcal{A} = \mathcal{A}' \cup \{ A_\ell \}$, with $A_\ell$ being the apices of $\mathcal{P}_\ell$ in $\mathfrak{M}$,

        \item[~T3~] if $A_\ell \neq \emptyset$, we let $\mathcal{Z} = \mathcal{Z}' \cup \{ Z_\ell \coloneqq  \emptyset \}$ and $\mathfrak{L} = \mathfrak{L}' \cup \{ \mathcal{L}_\ell \coloneqq  \emptyset \}$,

        \item[~T4~] we let $\mathsf{z} \in [\ell - 1]$ be maximal such that $A_i = \emptyset$ for all $i \in [\ell - \mathsf{z}, \ell-1]$,

        \item[~T5~] if $A_\ell = \emptyset$ and there exists a minimal, non-empty $V(\Omega_{j_\ell})$-$V(\Omega_{i+1})$-separator $Z \subseteq V(G'_\ell) \setminus V(H)$ of order at most $k$ in $G'_\ell \coloneqq \bigcup_{i \in F'} (G_i - ( \bigcup_{i \in [\ell - \mathsf{z}, \ell -1]} Z_i ) )$,\footnote{\Cref{def:district} already ensures that each district is disjoint from $\bigcup_{S \in \mathcal{A}} S$.} then
        \begin{enumerate}
            \item let $S \subseteq V(G'_\ell) \setminus V(H)$ be one side of the separation $(S,T)$ of order $2a+2$ in $G'_\ell - (Z \setminus \{ w \})$ with $w \in S$ for some $w \in Z$, such that $S$ induces a connected subgraph in $G_\ell' - (Z \setminus \{ w \})$, the set $S \cap T$ can be partitioned into $S_{j_\ell}, S_{i+1}$ with $|S_{j_\ell} \cap T| = |S_{i+1} \cap T| = a + 1$, there exists an $S_x$-$\bigcup_{h=1}^{\ell+1} V(\Omega_h)$-linkage $\mathcal{L}_x'$ in $G'_\ell[T]$ of order $a + 1$ for both $x \in \{ j_\ell, i+1 \}$, $\mathcal{L}_{j_\ell}' \cup \mathcal{L}_{i+1}'$ is an $S$-$\bigcup_{h=1}^{\ell+1} V(\Omega_h)$-linkage of order $2a + 2$, and $w$ is an $S_{j_\ell}$-$S_{i+1}$-separator in $G_\ell'[S]$, and

            \item we let $\mathcal{Z} = \mathcal{Z}' \cup \{ Z_\ell \coloneqq S \cup Z \}$ and $\mathfrak{L} = \mathfrak{L}' \cup \{ \mathcal{L}_\ell \coloneqq \emptyset \}$,
        \end{enumerate}

        \item[~T6~] if $A_\ell = \emptyset$ and there exists no $V(\Omega_{j_\ell})$-$V(\Omega_{i+1})$-path in $G'_\ell$, we let $\mathcal{Z} = \mathcal{Z}' \cup \{ Z_\ell \coloneqq  \emptyset \}$ and $\mathfrak{L} = \mathfrak{L}' \cup \{ \mathcal{L}_\ell \coloneqq  \emptyset \}$,

        \item[~T7~] if $A_\ell = \emptyset$ and there exists a $V(\Omega_{j_\ell})$-$V(\Omega_{i+1})$-linkage $\mathcal{L}$ of order $k+1$ in $G'_\ell$, we let $\mathcal{Z} = \mathcal{Z}' \cup \{ Z_\ell \coloneqq  \emptyset \}$ and $\mathfrak{L} = \mathfrak{L}' \cup \{ \mathcal{L}_\ell \coloneqq  \mathcal{L} \}$,
        
        \item[~T8~] let $Z' = \bigcup_{i \in [\ell - \mathsf{z}, \ell - 1]} Z_i \cup Z_\ell$, and let $\mathcal{F}$ be a partitioning of $[\ell+1]$ into as few sets as possible such that for all $F \in \mathcal{F}$ with $|F| \geq 2$ there exists a $V(\Omega_i)$-$V(\Omega_j)$-path in $G_\ell'' \coloneqq \bigcup_{i=1}^{\ell+1} G_i - Z'$ for all distinct $i,j \in F$.
    \end{description}
    We define the \emph{$k$-signature of $\mathfrak{M}$} to be a 4-tuple $(a',z,f,h)$ of non-negative integers, such that
    \begin{enumerate}
        \item $a'$ is the number of sets in $\mathcal{A}$ that are non-empty,

        \item $z$ is 0 if $A_\ell \neq \emptyset$ and otherwise $z$ is $z' \coloneqq |\{ Z' \in \mathcal{Z}' \cup \{ Z_\ell \} \colon Z' \neq \emptyset \}|$,

        \item $f \coloneqq | \{ F \in \mathcal{F} ~\!\colon\!~ |F| \geq 2, \text{ or } ( F = \{ i \} \text{ and } (G_i - Z',\Omega_i) \text{ has a cross)} \} |$, and

        \item $h \coloneqq |\{ i \in [\ell] ~\!\colon\!~ \mathcal{L}_i \neq \emptyset \text{ and there exists no } j \in [\ell] \text{ with } i < j \text{ and } A_j \cup Z_j \neq \emptyset . \}|$.
    \end{enumerate}
    Given two signatures $\Bar{x} = (x_1,x_2,x_3,x_4)$ and $\Bar{y} = (y_1,y_2,y_3,y_4)$, we say that $\Bar{x} > \Bar{y}$ if there exists an $i \in [2]$ such that for all $j \in [i-1]$, we have $x_j = y_j$ and we have $x_i > y_i$, or we have $x_i = y_i$ for both $i \in [2]$ and $x_3+x_4 > y_3+y_4$.
    If we have $x_i = y_i$ for all $i \in [4]$, we say that $\Bar{x} = \Bar{y}$ and we declare that $\Bar{x} \geq \Bar{y}$ if we have $\Bar{x} > \Bar{y}$ or $\Bar{x} = \Bar{y}$.
\end{definition}

The structure of the $\mathfrak{M}$-tree will generally play a supporting role in our analysis, giving us access to the apices of the transactions of $\mathfrak{M}$ and information on how the districts of $\mathfrak{M}$ connect to each other.

Note that the set $Z_\ell$ we find in \textbf{\textsf{T5}} of \Cref{def:Mtree} has a size we cannot bound in $k$.
It will however be convenient for us to delete all of it.
Nonetheless for later counting arguments, we establish that for a given set $Z' = S \cup Z \in \mathcal{Z}$ of the $\mathfrak{M}$-tree $\mathfrak{T} = ( \mathfrak{M}, T, \phi, \mathcal{A}, \mathcal{Z}, \mathcal{F}, \mathfrak{L} )$ with precision $k$, where $(S,T)$ and $Z$ are defined as in \textbf{\textsf{T5}} of \Cref{def:Mtree}, we let $B_\mathfrak{T}(Z) \coloneqq Z \cup (S \cap T)$ and note that $|B_\mathfrak{T}(Z)| \leq 3k+2$.

\section{Perilous districts}\label{sec:ApathsArgument}
The purpose of this section is to establish an argument that allows us to say, given a transaction mesh with large enough infrastructure, that the total number of districts in the transaction mesh with interesting behaviour is bounded or there exists a $K_t$-minor in our graph.

This manifests in the following statement.

\begin{theorem}\label{thm:ApathsArgument}
    Let $t,s,s',a,p,\ell$ be positive integers with $s \geq (2t+24)(8t^3+20)+s'\ell$ and $p \geq (64t^3-128)(2t+24)$.
    Let $(G,\Omega)$ be a society with a cylindrical rendition $\rho$ in a disk $\Delta$, a nest $\mathcal{C}$ of order $s$ in $\rho$ around the vortex $c_0$, and an $(s,s',a,p,\ell)$-transaction mesh $\mathfrak{M}$.
    Further, for each $i \in [\ell +1]$, let $(G_i,\Omega_i)$ be a district of $\mathfrak{M}$, with $G' \coloneqq  \bigcup_{i=1}^{\ell+1} G_i$.

    Then one of the following holds
    \begin{enumerate}
        \item the graph $G$ contains a $K_t$-minor model controlled by a mesh whose horizontal paths are subpaths of distinct cycles from $\mathcal{C}$, or

        \item there exist sets $S \subseteq V(G')$ and $I = I_1 \cup I_2 \subseteq [\ell+1]$ with $|S| \leq 8t^3-16$, $|I_1| \leq 64t^3-128$, and $|I_2| < \nicefrac{1}{2}(t-3)(t-4)$ such that each $\{ V(\Omega_i) ~\!\colon\!~ i \in [\ell+1] \}$-path in $G' - S$ is a $\{ V(\Omega_i) ~\!\colon\!~ i \in I_1 \}$-path and each district $(G_j - S,\Omega_j)$ with $j \in [\ell+1] \setminus (I_1\cup I_2)$ has a vortex-free rendition in a disk.
    \end{enumerate}
    Furthermore, there exists an algorithm that, given $\mathfrak{M}$ and $\rho$ as input, finds one of the outcomes above in time $\mathbf{O}(t^4|E(G)|)$.
\end{theorem}

To analyse the structure we are dealing with, let us first introduce a lemma that allows us to ``capture'' a bounded number of districts.
We show that we can root a bounded number of districts, together with their nests, on the society of our transaction mesh.
Later, we will reuse the construction from \cref{lemma:cliqes_in_extended_Dyckwalls} to obtain a large mesh where each of these districts belongs to its private face on the middle row.
These constructions will allow us to directly apply \cref{KtFromJumps} and force a $K_t$-minor if too many of the districts we selected were equipped with crosses.

\begin{figure}[ht]
    \centering
    \begin{tikzpicture}

        \pgfdeclarelayer{background}
		\pgfdeclarelayer{foreground}
			
		\pgfsetlayers{background,main,foreground}

        \begin{pgfonlayer}{background}
        \node (C) [v:ghost] {{\includegraphics[width=12cm]{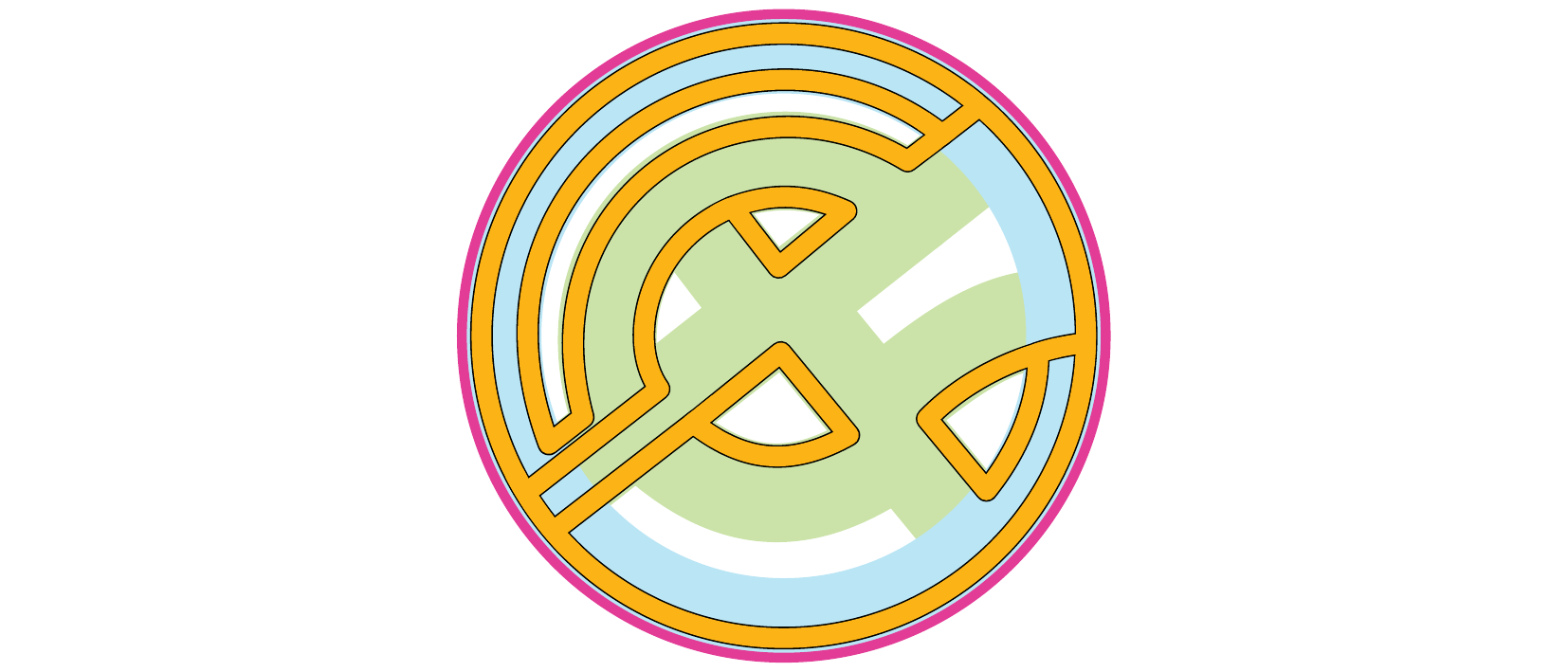}}};
        \end{pgfonlayer}{background}
			
        \begin{pgfonlayer}{main}

            \node (X) [v:ghost,position=0:0mm from C] {};
        
        \end{pgfonlayer}{main}

        \begin{pgfonlayer}{foreground}
        \end{pgfonlayer}{foreground}

    \end{tikzpicture}
    \caption{An illustration of the outcome of \cref{lemma:rooting_areas}, catching some of the districts in a transaction mesh and connecting their nests down to a part of the outer-most nest.}
    \label{fig:CatchingDistricts}
\end{figure}

\begin{lemma}\label{lemma:rooting_areas}
Let $k,q,s,s',a,p,\ell$ be non-negative integers with $1 \leq k \leq \ell$, $s \geq (4k+5)q + s'\ell$ and $p \geq 16kq$.
Let $(G,\Omega)$ be a society with an $(s,s',a,p,\ell)$-transaction mesh $\mathfrak{M} = (H, \Omega, \mathcal{C} = \{ C_1, \ldots , C_s \},\mathfrak{P})$ and let $\mathcal{A} = \{ (G_i,\Omega_i) ~\!\colon\!~ i \in [k] \}$ be a selection\footnote{Notice that in the definition of transaction meshes the areas are ordered. We slightly abuse notation here -- and many times later on -- by assigning new names $(G_i,\Omega_i)$ to the $k$ ``special'' areas we want to deal with.
This is safe because we do not care about the other areas in this process.} of $k$ pairwise distinct areas of $\mathfrak{M}$.

Then there exists an extended $q$-surface-wall $D$ with signature $(0,0,k)$ in $H$ such that the base cycles of $D$ are the cycles $C_{s-s'-q-1}, \dots , C_{s-s'} \in \mathcal{C}$ and there exists a bijection $\phi$ between the set $[k]$ and the vortex segments $S_v$ of $D$ such that for each $i \in [k]$, the nest of the vortex segment $\phi(i)$ is a subset of the district-nest of $(G_i,\Omega_i)$ (See \cref{fig:CatchingDistricts} for an illustration.).

Moreover, there exists an algorithm that finds $D$ in time $\mathbf{O}(qk|E(G)|)$.
\end{lemma}
\begin{proof}
By \cref{obs:radial_linkages_for_districts} each area $(G_i,\Omega_i)$ has a radial linkage $\mathcal{R}_i$ of order $\nicefrac{1}{4}p \geq 4kq$ for the area-nest $\mathcal{C}_i = \{ C^i_1, \dots , C^i_{s-s'\ell-1} \}$ of $(G_i,\Omega_i)$.
Moreover, $\mathcal{R}_i$ is made up of subpaths of one of the transactions among $\mathfrak{P}$.

For each $i \in [k]$, let $X_i'$ be the set of first vertices of $\mathcal{R}_i$ we encounter on the cycle $C^i_{5q+1} \in \mathcal{C}_i$ when traversing the paths of $\mathcal{R}_i$ starting from $C_1^i$.
Then let $X_i \subseteq X_i'$ be an arbitrary selection of $4q$ of these vertices.
Moreover, let $X \coloneqq \bigcup_{i \in [k]}X_i$.

Now let $\mathcal{P} \in \mathfrak{P}$ be a transaction of $\mathfrak{M}$ such that there exists a subpath $Y'$ of $C_{s - s'}$ that contains both endpoints of every member of $\mathcal{P}$.
Notice that such a transaction $\mathcal{P}$ must exist by the definition of transaction meshes and our choice of $\ell \geq 1$.
Let $Y$ be the set of all endpoints of $\mathcal{P}$ on $C_{s-s'}$.

It follows from \hyperref[prop:mengersthm]{Menger's Theorem} that there exists a linkage $\mathcal{L}$ of size $|X| = 4kq$ in $H$ between $Y$ and $X$, because there are at least $4kq$ cycles in $\mathcal{C}$, as well as $4kq$ cycles in each $\mathcal{C}_i$, at least $4kq$ paths from $\mathcal{P}$, and at least $4kq$ paths in each $\mathcal{R}_i$.
At least one of these families would need to be fully separated from the rest by a separator of size at most $4kq-1$, which is impossible.

Notice that $\mathcal{L}$ is fully contained in the inner graph of $C_{s-s'}$.

Next, we need to turn $\mathcal{L}$ into a linkage that is orthogonal to the innermost $q$ cycles of each area-nest as well as the cycles $C_{s-s'-(q-1)},\dots,C_{s-s'}$. 
For each $i\in[k]$, we apply \cref{lem:radialtoorthogonal} to $\mathcal{R}_i$ to obtain a radial linkage $\mathcal{R}_i'$ orthogonal to the cycles $C^i_1,\dots,C^i_{5q+1}$ and which is end-identical with $\mathcal{R}_i$.
Let $X_i^{\star}$ be the set of first vertices encountered on the cycle $C^i_{q+1}$ when following along the paths of $\mathcal{R}_i'$ starting from $C^i_1$ and let $X^\star = \bigcup_{i\in [k]}X_i^\star$.
Then, we follow along $\mathcal{L}$, starting from its endpoints on $C_{s-s'}$ until we meet the cycle $C_{s-s'-q-(4kq-1)}$.
Let $\mathcal{F}$ be the resulting linkage between $C_{s-s'}$ and $C_{s-s'-q-(4kq-1)}$.
Another application of \cref{lem:radialtoorthogonal} yields a linkage $\mathcal{F}'$ of the same size which is orthogonal to $\{ C_{s-s'-q-(4kq-1)},\dots,C_{s-s'} \}$ and end-identical with $\mathcal{F}$.
Finally, let $Y^{\star}$ be the set of all vertices on $C_{s-s'-q}$ we encounter when traversing along the paths of $\mathcal{F}'$ starting from $C_{s-s'}$.

Now consider the graph $H^{\star}$ obtained by taking the intersection of the inner graph of $C_{s-s'-q}$ and the outer graphs of the cycles $C^i_{q+1}$ for all $i\in[k]$ -- taken in the vortex-free rendition of $(H,\Omega)$ that exists according to \Cref{def:transactionmesh} -- by adding the cycles $C_{s-s'-q}$ and $C^i_{q+1}$ back in for all $i\in[k]$.
Observe that $H^\star$ still contains the cycles $C_{s-s'-q},\dots,C_{s-s'-q-(4kq-1)}$ as well as, for each $i\in[k]$, the cycles $C^i_{q+1},\dots,C^i_{5q+1}$.
It follows, than an application of \cref{prop:mengersthm} yields an $X^{\star}$-$Y^{\star}$-linkage $\mathcal{F}^{\star}$ of size $4qk$ in $H^{\star}$.
We may now combine $\mathcal{F}^{\star}$ with the orthogonal linkages $\mathcal{R}_i'$, $i\in[k]$, and the orthogonal linkage $\mathcal{F}'$ to obtain a linkage $\mathcal{L}^{\star}$ whose intersection with any of the cycles $C^i_1,\dots,C^i_q$ for all $i\in[k]$, and the cycles $C_{s-s'-(q-1)},\dots,C_{s-s'}$ is a single, non-empty path.
Moreover, $\mathcal{L}^{\star}$ has $4q$ endpoints on $C^i_1$ for each $i\in[k]$.

It is now easy to see that $\{ C_{s-s'},\dots,C_{s-s'-(q-1)}\}$ together with $\mathcal{L}^{\star}$ and the cycles $\{ C^i_1,\dots,C^i_q\}$ for each $i\in[k]$ form our desired extended $q$-surface-wall.
\end{proof}

Combining \cref{lemma:rooting_areas} with the construction of the graph $H'_{w,2t}$ in extended surface-walls with vortices we have used in the proof of \cref{lemma:cliqes_in_extended_Dyckwalls} yields a $(w\times 2t)$ mesh $M$ such that each of our ``special'' areas is highly connected to its private face in the middle row of $M$.
This allows us to call upon \cref{KtFromJumps} to obtain the following corollary.

\begin{corollary}\label{cor:many-to-many}
Let $t,s,s',a,p,\ell,k$ be non-negative integers with $k = 2t^3-4$, $1 \leq t \leq \ell$, $s \geq (4k+5)(2t+24)+s'\ell$, and $p \geq 16k(2t+24)$, and let $(G,\Omega)$ be a society with an $(s,s',a,p,\ell)$-transaction mesh $\mathfrak{M} = (H, \Omega, \mathcal{C}, \mathfrak{P})$.
Moreover, let $\mathcal{A} = \{ (G_i,\Omega_i) ~\!\colon\!~ i \in [k] \}$ be a selection of $k$ pairwise distinct areas of $\mathfrak{M}$.

If there exists a collection $\mathcal{Q}$ of $t^3-2$ pairwise disjoint paths $Q_j$ in $G$ such that $\mathcal{Q}$ is a family of $\{ V(G_i) ~\!\colon\!~ i \in [k] \}$-paths and each $Q_j$ is internally disjoint from $H$, then there exists a $K_t$-minor model in $G$ that is controlled by a mesh whose horizontal paths are subpaths of distinct cycles from $\mathcal{C}$.

Moreover, there exists an algorithm that finds this minor model in time $\mathbf{O}(t^4|E(G)|)$.
\end{corollary}

We are now ready for the proof of \cref{thm:ApathsArgument}.

\begin{proof}[Proof of \Cref{thm:ApathsArgument}]
For each $i \in [\ell+1]$, let $\mathcal{C}_i = \{ C_1^i, \dots , C^i_{s-s'\ell-1} \}$ be the district nest of $(G_i,\Omega_i)$.
We set $\mathcal{X} \coloneqq \{ C_1^i, C_2^i  ~\!\colon\!~ i \in [\ell+1] \}$, $\mathcal{X}' \coloneqq \{ C_1^i  ~\!\colon\!~ i \in [\ell+1] \}$, and $\mathcal{G} \coloneqq ( C_1^1, C_2^1, C_1^2, \ldots , C_2^{\ell+1} )$.
Notice that each member of $\mathcal{X}$ is a cycle and thus the maximum degree of any member of $\mathcal{X}$ is $2$.

Suppose $P$ is a $\mathcal{X}$-path in $G'$.
Then, if there exists some $i \in [\ell+1]$ such that $P$ has an endpoint in $V(C_2^i)$, the other endpoint of $P$ must lie in $V(C_1^i)$, since $C_1^i$ separates $C_2^i$ from all other members of $\mathcal{X}$, since each of our areas has a cylindrical rendition with $\mathcal{C}_i$ being its nest of grounded cycles.
Thus any $\mathcal{G}$-jump has to be an $\mathcal{X}'$-path.

We now ask \cref{PPjumps} for a subsequence $\mathcal{G}'$ of $\mathcal{G}$ such that there exist $t^3-2$ pairwise vertex-disjoint $\mathcal{G'}$-jumps in $G'$ with endpoints in $2t^3-4$ pairwise $\mathcal{G}'$-independent graphs from $\mathcal{X}$.

If these $\mathcal{G}'$-jumps exist, let $\mathcal{L}$ be the collection of these paths.
According to our observation above, $\mathcal{L}$ is also a collection of $\mathcal{X}'$-paths. 
Let $J \subseteq [\ell+1]$ be the selection of the $2t^3-4$ districts that have endpoints of paths from $\mathcal{L}$.
By \cref{cor:many-to-many} this means we can find, in time $\mathbf{O}(t^4|E(G)|)$, a $K_t$-minor model controlled by a mesh whose horizontal paths are subpaths of distinct members of $\mathcal{C}$ as desired.

Hence, we may assume that such a family $\mathcal{L}$ does not exist.
Instead, \cref{PPjumps} returns a set $S \subseteq V(G')$ with $|S| \leq 8t^3-16$ and a set $\mathcal{Y} \subseteq \mathcal{X}$ with $|\mathcal{Y}| \leq 64t^3-128$ such that every $\mathcal{G}$-jump in $G' - S$ must have both endpoints in $\bigcup \mathcal{Y}$.

For each $i \in [\ell+1]$, the predecessor of $C_1^i$ in $\mathcal{G}$, if it exists, is $C_2^{i-1}$ and the successor of $C_1^i$ in $\mathcal{G}$ is $C_2^i$.
Thus any $\mathcal{X}'$-path in $G' - S$ is a $\mathcal{G}$-jump and thus must have both endpoints in $\bigcup \mathcal{Y}$.
Furthermore, for every $i \in [\ell+1]$ the cycle $C^i_1$ separates $V(\Omega_i)$ from the district vortex of $(G_i,\Omega_i)$.
Let $I' \coloneqq \{ i \in [\ell+1] ~\!\colon\!~ C^i_1 \in \mathcal{Y} \}$.
It follows that in $G' - S$ every $\{ V(\Omega_i) ~\!\colon\!~ i \in [\ell+1] \}$-path must be a $\{ V(\Omega_i) ~\!\colon\!~ i \in I' \}$-path.

This implies that no component of $G'' \coloneqq \bigcup_{i \in [\ell+1] \setminus I'} G_i - S$ contains vertices of $V(\Omega_i)$ and $V(\Omega_j)$ for any pair of distinct $i,j \in [\ell+1] \setminus I'$.
Let $J \subseteq [\ell+1] \setminus I'$ be the collection of all integers $j \in [\ell+1] \setminus I'$ such that $(G_i - S,\Omega_i)$ has a cross.
If $|J| \geq \nicefrac{1}{2}(t-3)(t-4)$ it follows, via an application of \cref{lemma:rooting_areas} to $\mathfrak{M}$ and the selection $J$, followed by an application of \cref{lemma:cliqes_in_extended_Dyckwalls}, that $G$ contains a $K_t$-minor model which is controlled by a mesh whose horizontal paths are subpaths of distinct cycles from $\mathcal{C}$.
Note that the order of the extended surface wall we require is in $\mathbf{O}(t^4)$, which further justifies our claimed runtime.
Thus, we may assume $|J| < \nicefrac{1}{2}(t-3)(t-4)$.
By setting $I \coloneqq I' \cup J$ we are done.
\end{proof}

\section{The society classification theorem}\label{sec:societyclassification}
In this section we prove the society classification theorem, which is the key tool to proving the local structure theorem that will itself appear in \Cref{sec:localstructure}.
Before we state it, let us briefly give an imprecise intuition behind this theorem.
It says that given a cylindrical rendition of a society with a sufficiently large nest with a radial linkage, we can guarantee one of four outcomes:
\begin{enumerate}
    \item We find a $K_t$-minor that attaches conveniently to a mesh we find in our nest,

    \item we find a large crosscap transaction (under a small set of apices),
    
    \item we find a large handle transaction (under a small set of apices), or

    \item we can completely decompose the entire vortex into a few vortices of bounded depth (using a larger, but still relatively small set of apices) and these vortices are themselves surrounded by nests that attach neatly to an extended surface-wall.
\end{enumerate}
The name of the theorem derives from the idea that it classifies societies into one of four options, each of which is in some way beneficial to us.
We will give estimates on the order of the functions involved in this theorem in its proof in \Cref{sec:classifyingsocieties}.

\begin{theorem}[Society classification theorem]\label{thm:societyclassification}
    There exist polynomial functions $\mathsf{apex}^\mathsf{genus}_{\ref{thm:societyclassification}}, \mathsf{loss}_{\ref{thm:societyclassification}} \colon \mathbb{N} \rightarrow \mathbb{N}$, $\mathsf{nest}_{\ref{thm:societyclassification}}, \mathsf{cost}_{\ref{thm:societyclassification}} \colon \mathbb{N}^2 \rightarrow \mathbb{N}$, and $ \mathsf{apex}^\mathsf{fin}_{\ref{thm:societyclassification}}$,  $\mathsf{depth}_{\ref{thm:societyclassification}} \colon \mathbb{N}^3\to\mathbb{N}$, such that the following holds.
    
    Let $t,s,p,k$ be positive integers with $s \geq \mathsf{nest}_{\ref{thm:societyclassification}}(t,k)$.
    Let $(G,\Omega)$ be a society with a cylindrical rendition $\rho$ in a disk $\Delta$, a cozy nest $\mathcal{C} = \{ C_1, \ldots , C_s \}$ of order $s$ in $\rho$ around the vortex $c_0$.
    Further, let $(G', \Omega')$ be the $C_{s - \mathsf{cost}_{\ref{thm:societyclassification}}(t,k)}$-society in $\rho$.

    Then $G'$ contains a set $A \subseteq V(G')$ such that one of the following exists:
    \begin{enumerate}
        \item A $K_t$-minor model in $G$ controlled by a mesh whose horizontal paths are subpaths of distinct cycles from $\mathcal{C}$.

        \item A flat crosscap transaction $\mathcal{P}$ of order $p$ in $(G'-A,\Omega')$, with $|A| \leq \mathsf{apex}^\mathsf{genus}_{\ref{thm:societyclassification}}(t)$, and a nest $\mathcal{C}'$ in $\rho$ of order $s - (\mathsf{loss}_{\ref{thm:societyclassification}}(t)+\mathsf{cost}_{\ref{thm:societyclassification}}(t,k))$ around $c_0$ to which $\mathcal{P}$ is orthogonal.

        \item A flat handle transaction $\mathcal{P}$ of order $p$ in $(G'-A,\Omega')$, with $|A| \leq \mathsf{apex}^\mathsf{genus}_{\ref{thm:societyclassification}}(t)$, and a nest $\mathcal{C}'$ in $\rho$ of order $s - (\mathsf{loss}_{\ref{thm:societyclassification}}(t)+\mathsf{cost}_{\ref{thm:societyclassification}}(t,k))$ around $c_0$ to which $\mathcal{P}$ is orthogonal.

        \item A rendition $\rho'$ of $(G - A, \Omega)$ in $\Delta$ with breadth $b \in [\nicefrac{1}{2}(t-3)(t-4)-1]$ and depth at most $\mathsf{depth}_{\ref{thm:societyclassification}}(t,k,p)$, $|A| \leq \mathsf{apex}^\mathsf{fin}_{\ref{thm:societyclassification}}(t,k,p)$, and an extended $k$-surface-wall $D$ with signature $(0,0,b)$, such that $D$ is grounded in $\rho'$, the base cycles of $D$ are the cycles $C_{s-\mathsf{cost}_{\ref{thm:societyclassification}}(t,k)-1-k},\dots,C_{s-\mathsf{cost}_{\ref{thm:societyclassification}}(t,k)-1}$, and there exists a bijection between the vortices $v$ of $\rho'$ and the vortex segments $S_v$ of $D$, where $v$ is the unique vortex contained in the disk $\Delta_{S_v}$ defined by the trace of the inner cycle of the nest of $S_v$ where $\Delta_{S_v}$ is chosen to avoid the trace of the simple cycle of $D$.
    \end{enumerate}
    
    In particular, the set $A$, the $K_t$-minor model, the transaction $\mathcal{P}$, the rendition $\rho'$, and the extended surface-wall $D$ can each be found in time $\mathbf{poly}( t + s + p + k ) |E(G)||V(G)|^2$.
\end{theorem}

Before we delve into the quite extensive procedure to prove this theorem, we give a rough guide to our approach.
In the broadest of terms, we will simply try to build a transaction mesh in $(G,\Omega)$ and along the way, if we fail to keep building the transaction mesh, we will find the four outcomes stated in the theorem.

The actual proof is then split into three parts:
\emph{(i)} We have to show that we can do something with a transaction mesh that we know cannot be extended anymore, meaning that we can no longer find transaction large enough to increase the transaction mesh for one more iteration.
This will actually be the exclusive source for the last two options of the statement.
(See \Cref{sec:boundeddepth} for this part.)

\emph{(ii)} If we do end up finding a big transaction that looks promising, we have to show that we can actually build a new transaction mesh using it.
This in particular means we have to delve into the internals of \Cref{def:transactionmesh} and also \Cref{def:Mtree}, as the $\mathfrak{M}$-tree will be crucial to the third part of the proof.
(See \Cref{sec:buildingtransactionmesh} for this part.)

\emph{(iii)} Should we never run into the first part and instead keep finding infrastructure to increase the order of our transaction mesh, we also run into trouble.
Since we lose a little bit of the transactions and in particular, the nest of the society we started with, at some point this process will destroy all that we have built.
Thus, in this last part, we show that there is actually an explicit bound on the order of a transaction mesh we can build with our methods without encountering a $K_t$-minor model.
(See \Cref{sec:signatures} for this part).

The proof that puts all of the intermediate results we find in these subsections together and ultimately yields \Cref{thm:societyclassification} is then found in \Cref{sec:classifyingsocieties}.

\subsection{What to do when every district has bounded depth}\label{sec:boundeddepth}
Though the fourth option of \Cref{thm:societyclassification} may seem to be the most involved at first glance, we will actually discuss it first as it is the option that has the lowest requirement for analysis of the internals of \Cref{def:transactionmesh} and \Cref{def:Mtree}.
Mainly, we are interested in how to use linkages that jump between districts of our transaction mesh to find handle or crosscap transactions and if we cannot find such a linkage, we separate all districts from each other and employ the tools from \Cref{sec:ApathsArgument} to verify the details of the fourth item of \Cref{thm:societyclassification}.

This explanation is of course an oversimplification.
As we will see later, one of the principle technical hurdles in this subsection will be that we need to simultaneously orthogonalise two flat transactions within the same society.
This is used in the one spot of our proof in which handle transactions can appear.
Surprisingly, the machinery we have developed earlier is not quite enough to solve this problem.
Instead we will have to argue ad hoc that this is possible, though the methods we use then again yield quite modest bounds on the amount of paths from the transactions and cycles from the nest this costs us.

We will start by slightly generalising \Cref{lemma:orthogonal_transaction} to linkages that connect two societies, which requires us to first define such linkages.

\begin{definition}[Connecting linkages]\label{def:connectinglinkage}
    Let $G$ be a graph and let $(G_1,\Omega_1),(G_2,\Omega_2)$ be two societies with $G_1, G_2 \subseteq G$ and $V(\Omega_1) \cap V(\Omega_2) = \emptyset$, such that for both $i \in [2]$ there exists a cylindrical rendition of $(G_i,\Omega_i)$ in a disk with a nest $\mathcal{C}_i$ around a vortex $c_i$ and we have $G_i \subseteq \sigma_{\rho_{3-i}}(c_{3-i}) - \widetilde{c}_{3-i}$.
    
    A $V(\Omega_1)$-$V(\Omega_2)$-linkage $\mathcal{P}$ in $G$ is called a \emph{connecting linkage (between $(G_1,\Omega_1)$ and $(G_2,\Omega_2)$)}.
    We call $\mathcal{P}$ \emph{orthogonal to $\mathcal{C}_1$ and $\mathcal{C}_2$} if the intersection of each path of $\mathcal{P}$ with each cycle of $\mathcal{C}_i$, for both $i \in [2]$, consists of a single, non-empty path.
    The \emph{end segments of $\mathcal{P}$} are two minimal segments of $\Omega_1$ and $\Omega_2$ that contain the endpoints of $\mathcal{P}$, which are unique as long as $\mathcal{P}$ has order at least 3.
\end{definition}

The proof for the next lemma again follows the proof of \Cref{lemma:orthogonal_transaction} closely.
Thus we only sketch the differences, as in the proofs of similar lemmas throughout this article.

\begin{lemma}\label{lemma:orthogonal_connecting_linkage}
Let $s,p,q$ be positive integers and let $k$ be a non-negative integer.
Let $\{ ((G_i,\Omega_i),\rho_i,c_i,\Delta_i,\mathcal{C}_i) \}_{i \in [0,k]}$ be a full $(s,k)$-society configuration, such that for each $i \in [0,k]$ the tuple $(\mathcal{C}_i,\mathfrak{P}_i)$ is an $(s,p,\ell_i)$-padding for $(G_i,\Omega_i)$ with $\ell_i \leq \ell$ and the inner nest $\mathcal{C}_i' = \{ C_1^i, \ldots , C_s^i \}$.
Further, let $\mathcal{P}$ be a connecting linkage in $G$ of order $p(q-1)+1$ between $(G_0,\Omega_0)$ and $(G_1,\Omega_1)$ with the end segments $X_0, X_1$ such that $X_i \subseteq V(\Omega_i)$ for both $i \in \{ 0,1 \}$.

Then there exists a connecting linkage $\mathcal{P}'$ in $G$ of order $q$ between $(G_0,\Omega_0)$ and $(G_1,\Omega_1)$ such that
\begin{enumerate}
    \item $\mathcal{P}'$ is orthogonal to $\{ C_2^i, \dots , C_s^i \}$ for both $i \in \{ 0,1 \}$,
    
    \item $\bigcup \mathcal{P}'$ is disjoint from the union of $C_2^i$ and the outer graph of $C_2^i$ in $\rho_i$ for all $i \in [2,k]$, and
    
    \item $\mathcal{P}'$ connects vertices of $X_0 \cap V(\mathcal{P})$ to vertices of $X_1 \cap V(\mathcal{P})$.
\end{enumerate}
Moreover, there exists an algorithm that finds $\mathcal{P}'$ in time $\mathbf{O}( (p+k) |E(G)| )$.
\end{lemma}
\begin{proof}[Proof sketch]
    For all $i \in [0,k]$, let $H_i$ be the graph derived from outer graph of $C_1^i$ in $\rho_i$ by deleting $V(C_1^i)$.
    We let $G'$ be $G - V(\bigcup_{i=0}^k H_i)$.
    For both $i \in \{ 0, 1 \}$ and each $P \in \mathcal{P}$ let $P_i$ be the unique $X_i$-$C_1^i$-path in $P$.
    We set $\mathcal{L}_i = \{ P_i ~\!\colon\!~ P \in \mathcal{P} \}$, for both $i \in \{ 0,1 \}$, and let $X_i'$ be the set of all endpoints of $\mathcal{L}_i$ in $V(C_1^i)$.
    By applying \Cref{lem:radialtoorthogonal} to both $\mathcal{L}_0$ and $\mathcal{L}_1$, we find two new radial linkages $\mathcal{L}_0'$ and $\mathcal{L}_1'$, such that $\mathcal{L}_i'$ is an $X_i$-$X_i'$-linkage and $\mathcal{L}_i'$ is orthogonal to $\mathcal{C}_i$ for both $i \in \{ 0,1 \}$.

    If we now apply \hyperref[prop:mengersthm]{Menger's Theorem} in $G'$ to find an $X_0'$-$X_1'$-linkage $\mathcal{P}'$ of order $q$, then we are done if we actually find $\mathcal{P}'$, via arguments analogous to those presented in the proof of \Cref{lemma:orthogonal_transaction}.
    Thus we must instead find an $X_0'$-$X_1'$-separation $(A,B)$ of order at most $q-1$ in $G'$.
    We let $S \coloneq A \cap B$ and for every $P \in \mathcal{P}$, we let $\widetilde{P}$ be the unique $X_0'$-$X_1'$-subpath of $P$, with $\widetilde{\mathcal{P}}$ being the collection of all paths $\widetilde{P}$.
    By definition $\widetilde{\mathcal{P}}$ is an $X_0'$-$X_1'$-linkage of order $(p-1)(q-1)+1$ in $G$.

    We may now define gangplanks and gangboards as in the proof of \Cref{lemma:orthogonal_transaction_co-conspirators}, in reference to $c_i$ and $\rho_i$ if the path $\widetilde{P} \in \widetilde{\mathcal{P}}$ intersects $C_1^i$ for $i \in [0,k]$.
    If there existed a path $\widetilde{P} \in \widetilde{\mathcal{P}}$ such that the path and all of its gangplanks avoided $S$, we would reach a contradiction to $S$ being a separator.
    Thus, just as in the proof of \Cref{lemma:orthogonal_transaction} and \Cref{lemma:orthogonal_transaction_padded}, we note that there must exist some vertex $u \in S$ which is contained in the gangplanks of at least $p$ distinct paths from $\widetilde{\mathcal{P}}$.
    This circumstances now leads us towards a contradiction in exactly the same fashion as it did in our proof of \Cref{lemma:orthogonal_transaction_padded}.
\end{proof}

The principal agenda of this subsection is to establish a small set of tools to handle certain routing problems of radial linkages within a padded transaction mesh $\mathfrak{M}$.
The purpose of these tools is to allow for the projection of a large crosscap transaction within a district of $\mathfrak{M}$ back to the original society of the transaction mesh.
Moreover, sometimes we will find large linkages between two distinct districts of $\mathfrak{M}$ and we will also have to be able to route such linkages back onto the original society while, at the same time, ensuring the existence of a big transaction within the infrastructure of $\mathfrak{M}$ itself that separates the two districts while being fully disjoint from the newly routed linkages.
This second part is crucial for the construction of a handle transaction.
The main difficulty here also reveals the sole reason for our introduction of paddings:
In later steps of our proofs we will need these linkages to be much larger than the original nest to avoid an exponential blow-up of our bounds.
To ensure this, it is not enough to rely on the infrastructure provided by the regular district nests, instead we require the full power of the padding to provide enough connectivity for the routing of such large linkages.

Let $s,s',a,p,\ell,k$ be non-negative integers with $\max(a,s')\leq s$ and $2s+2-(s'\ell)\leq p$.
Let $\mathfrak{M}=(H,\Omega,\mathcal{C},\mathfrak{P}=\{ \mathcal{P}_1,\dots,\mathcal{P}_{\ell}\})$ be a padded $(s,s',a,p,\ell)$-transaction mesh in the society $(G,\Omega)$.
Finally, let $k\in[\lfloor\nicefrac{p}{2}\rfloor-s'\ell]$.
The \emph{$k$-splitstream} of $\mathfrak{M}$ is the pair $(\mathfrak{S}_k,\mathfrak{T}_k)$ where
\begin{align*}
\mathfrak{S}_k & \coloneq \big\{ \mathcal{S}_i = \{ P^i_j \colon j\in[k]\cup[p-k+1,p] \} \colon i\in[\ell] \big\}\text{, and}\\
\mathfrak{T}_k & \coloneqq \big\{ \mathcal{S}_i = \{ P^i_j \colon j\in[k+1,\lfloor\nicefrac{p}{2}\rfloor-s'\ell]\cup[\lceil \nicefrac{(p+1)}{2}\rceil-(s'\ell+k),p-k]  \} \colon i\in[\ell] \big\}.
\end{align*}
We call $\mathcal{S}_k$ the \emph{surface stream} and $\mathcal{T}_k$ the \emph{deep stream} of the $k$-splitstream of $\mathfrak{M}$.

Notice, that for $\mathfrak{X}\in\{\mathfrak{S},\mathfrak{T} \}$, the graph $X = \bigcup \mathfrak{X}_k$ is not necessarily connected.
However, from any vertex $v\in V(X)$ there is a path in $X$ from $v$ to $V(\Omega)$.

We begin by showing that a large radial linkage in a district of $\mathcal{M}$ can -- essentially -- be projected onto the original society by only using an appropriately sized surface stream of $\mathfrak{M}$.

\begin{lemma}\label{lemma:RadialLinkageOnTheSurface}
Let $t,s,s',a,p,\ell,k$ be non-negative integers with $\max(a,s')\ell\leq s$, $2s+2-(s'\ell)\leq p$, and $k\in[\lfloor\nicefrac{p}{2}\rfloor-s'\ell]$.
Let $\mathfrak{M}=(H,\Omega,\mathcal{C},\mathfrak{P}=\{ \mathcal{P}_1,\dots,\mathcal{P}_{\ell}\})$ be a padded $(s,s',a,p,\ell)$-transaction mesh in the society $(G,\Omega)$ with districts $\{ (G_i,\Omega_i) \}_{i\in[\ell+1]}$ and cylindrical rendition $\rho$ with nest $\mathcal{C}$.
Moreover, let $i\in[\ell+1]$ and $(\mathcal{C}_i,\mathfrak{P}_i)$ be the padding of $(G_i,\Omega_i)$ and $\mathcal{R}$ be a radial linkage of order at least $k$ in $(G_i,\Omega_i)$.
Finally, let $X$ be the set of endpoints of $\mathcal{R}$ on the innermost cycle of $(\mathcal{C}_i,\mathfrak{P}_i)$ and let $(G',\Omega')$ be the $C_{s-s'\ell}$-society in $\rho$.

Then there exists an $X$-$V(\Omega')$-linkage $\mathcal{L}$ of order $k$ in $\bigcup\mathcal{C} \cup \bigcup \mathfrak{S}_k \cup \bigcup\mathcal{R}$ where $\mathfrak{S}_k$ is the surface stream of the $k$-splitstream of $\mathfrak{M}$.
Moreover, the linkage $\mathcal{L}$ can be found in time $\mathbf{O}(k|E(G)|)$.
\end{lemma}

\begin{proof}
Let $H'\coloneqq \bigcup\mathcal{C} \cup \bigcup \mathfrak{S}_k \cup \bigcup\mathcal{R}$.
Suppose there does not exist an $X$-$V(G)$-linkage $\mathcal{L}$ of order $k$ in $G'$.
Then, by \cref{prop:mengersthm}, there exists a set $Y$ of size at most $k-1$ such that there is no path from $X$ to $V(\Omega)$ in $H'-Y$.
It follows that there exists a path $R\in\mathcal{R}$ that is disjoint from $Y$.
By choice of $Y$ we know that $R$ cannot have an endpoint on $V(\Omega')$.
It follows that there exists $\iota\in[\ell]$ such that the endpoint of $R$ that is not in $X$ belongs to a path of $\mathcal{P}_{\iota}$.
Hence, $\mathcal{R}$ intersects $k$ paths of the sublinkage $\mathcal{S}_{\iota}\in\mathfrak{S}_k$ of $\mathcal{P}_{\iota}$.

Now notice that there exists a sequence $\sigma=\langle \mathcal{Z}_{i_j} \rangle_{j\in[\ell']}$, $\ell'\leq \ell$ and $i_j\in[\ell]$ for all $j\in[\ell']$ such that
\begin{enumerate}
    \item for each $i\in[\ell']$, $\mathcal{Z}_{i_j}\subseteq \mathcal{S}_{i_j}$ is a linkage of order $k$,
    \item $i_1 = \iota$,
    \item each path in $\mathcal{Z}_{i_{\ell'}}$ contains a vertex of $V(\Omega')$, and
    \item for all $j\in[\ell'-1]$, every path of $\mathcal{Z}_{i_j}$ intersects every path of $\mathcal{Z}_{i_{j+1}}$.
\end{enumerate}
The existence of $\sigma$ follows directly from the definition of transaction meshes and $k$-splitstreams.
Now, however, as $|Y|\leq k-1$, each $\mathcal{Z}_{i_j}$ contains a path $Z_{i_j}$ disjoint from $Y$.
It follows that $R\cup \bigcup_{j\in[\ell']}Z_{i_j}$ is connected, disjoint from $Y$ and contains a vertex of both $X$ and $V(\Omega')$.
Hence, there exists an $X$-$V(\Omega')$-path in $H'-Y$ which is absurd.

Therefore, we must be able to find our desired linkage $\mathcal{L}$ by applying \cref{prop:mengersthm}.
\end{proof}

Next, we will show that, given a transaction mesh $\mathfrak{M}$, a $k$-splitstream $(\mathfrak{S},\mathfrak{T})$ of $\mathfrak{M}$, and two ``special'' districts $A$ and $B$ of $\mathfrak{M}$ we can always find a large transaction on $(G,\Omega)$ which separates $A$ and $B$ in $\mathfrak{M}$ and is entirely contained in the deep stream of $(\mathfrak{S},\mathfrak{T})$.
This transaction will be the base for the construction of handle transactions later on.

\begin{lemma}\label{lemma:SplitTwoDistricts}
Let $t,s,s',a,p,\ell,k$ be non-negative integers with $\max(a,s')\ell\leq s$, $2s+2-(s'\ell)\leq p$, and $k\in[\lfloor\nicefrac{p}{2}\rfloor-s'\ell]$.
Let $\mathfrak{M}=(H,\Omega,\mathcal{C},\mathfrak{P}=\{ \mathcal{P}_1,\dots,\mathcal{P}_{\ell}\})$ be a padded $(s,s',a,p,\ell)$-transaction mesh in the society $(G,\Omega)$ with districts $\{ (G_i,\Omega_i) \}_{i\in[\ell+1]}$ and cylindrical rendition $\rho$ with nest $\mathcal{C}$.

Let $i_1,i_2\in[\ell+1]$ be distinct and let $(H',\Omega')$ the $C_{s-s'\ell}$-society in $\rho$ with respect to $H$, i.e.\@ $H'\subseteq H$.
Then there exists a transaction $\mathcal{Q}$ of order $h\coloneqq \lfloor\nicefrac{p}{2}\rfloor-s'\ell-k$ in $(H',\Omega')$ that separates the areas $H_{i_1}$ and $H_{i_2}$ in $H$ where $(H_{i_1},\Omega_{i_1})$ and $(H_{i_2},\Omega_{i_2})$ are distinct areas of $\mathfrak{M}$ such that $\mathcal{Q}$ is completely contained in the deep stream of the $k$-splitstream of $\mathfrak{M}$.
Moreover, the transaction $\mathcal{Q}$ can be found in time $\mathbf{O}(p|E(G)|)$.
\end{lemma}

\begin{proof}
We will iteratively construct a sequence of constrictions $(B_i,\Lambda_i)$ of the society $(H',\Omega')$ by considering the linkages $\mathcal{P}_i$ one at a time.
Each time, we add $\mathcal{P}_{i+1}$ to $(B_i,\Lambda_i,\mathfrak{L}_i)$ we ask if $\mathcal{P}_{i+1}$ separates $H_{i+1}$ and $H_{i+2}$ from one another in $B_i$.
If this is the case we make use of the linkages constricting $(B_i,\Lambda_i)$ to extend $\mathcal{P}_{i+1}$ back to $V(\Lambda_i)$.
Otherwise, we know that one of the two ``sides'' of $\mathcal{P_{i+1}}$ in $B_i$ is fully disjoint from both of our special areas.
In this case, we can further restrict the society.
Since $H_{i_1}$ and $H_{i_2}$ define distinct areas in $\mathfrak{M}$, this process must eventually find a large transaction $\mathcal{Q}'$ separating both.
At this point it suffices to restrict the attention to the linkages found in the deep stream which correspond to the linkages used in the construction of $\mathcal{Q}'$ in order to obtain the desired transaction $\mathcal{Q}$.

Notice that, by definition of transaction meshes, there exists a vortex-free rendition $\rho'$ of $(H',\Omega')$ in the disk.
We set $(B_0,\Lambda_0,\emptyset)\coloneqq (H',\Omega',\emptyset)$.

Now, suppose for some $j\in[0,\ell-1]$ we have already constructed the $(\lfloor\nicefrac{p}{2}\rfloor-s'\ell,z_j)$-constriction $(B_j,\Lambda_j,\mathfrak{L}_j)$ of $(B_0,\Lambda_0)$ with $z_j\in[0,j]$ satisfying the following properties:
\begin{enumerate}
    \item we may index the set $\mathfrak{L}_j=\{ \mathcal{L}_{i} \colon i \in I_j\subseteq [j]\}$ such that every path of $\mathcal{L}_{i}$ is a subpath of some path in $\mathcal{P}_{i}$ where $|I_j|=z_j$,
    \item for each $i\in[j]\setminus I_j$, $\mathcal{P}_i$ is disjoint from $B_j$, 
    \item for all $i\in[j+1,\ell]$, $\mathcal{P}_i$ is either fully contained in $B_j$ or disjoint from $B_j$, and
    \item $H_{i_x}\subseteq B_j$ for both $x\in[2]$.
\end{enumerate}
Since both $H_{i_1}$ and $H_{i_2}$ are contained in $B_j$ there must be some $\iota\in[j+1,\ell]$ such that $\mathcal{P}_{\iota}$ is contained in $B_j$.
Let us select $\iota$ to be the smallest such index, then $\mathcal{P}_{\iota}$ is a transaction on $(B_j,\Lambda_j)$.
We consider two cases.
\smallskip

\textbf{Case 1:} $\mathcal{P}_{\iota}$ separates $H_{i_1}$ and $H_{i_2}$ in $H$.

In this case let $\mathfrak{T}_k$ be the deep stream of the $k$-splitstream of $\mathfrak{M}$.
Notice that $\mathcal{T}_{\iota}\subseteq \mathcal{P}_{\iota}$, i.\@e.\@ the subtransaction of $\mathcal{P}_{\iota}$ that makes up the deep stream, still separates $H_{i_1}$ and $H_{i_2}$ in $B_j$.
Moreover, it is now straight forward to see that there exists a transaction $\mathcal{Q}\subseteq \bigcup\mathfrak{T}_k$ of order $\lfloor\nicefrac{p}{2}\rfloor-s'\ell-k=\nicefrac{1}{2}|\mathcal{T}_{\iota}|$ on $(B_0,\Lambda_0)$ that contains half of the paths of $\mathcal{S}_{\iota}$.

We make two observations on $\mathcal{Q}$.

First, $\mathcal{Q}$ can be found inductively by first finding a corresponding transaction $\mathcal{Q}_j$ in $(B_j,\Lambda_j)$ and then extending $\mathcal{Q}_j$ towards $V(\Lambda_0)$ by retracing the sequence of constrictions we constructed so far.

Second, $\mathcal{Q}$ separates $H_{i_1}$ and $H_{i_2}$ in $B_0=H'$.

Hence, in this case we have found the desired transaction and may conclude.
\smallskip

\textbf{Case 2:} $\mathcal{P}_{\iota}$ does not separate $H_{i_1}$ and $H_{i_2}$ in $H$.

In this case there exists $y\in\{ 1,p\}$ such that every $V(H_{i_1})$-$V(P^{\iota}_{p+1-y})$-path in $H$ intersects $P^{\iota}_y$.
Without loss of generality we may assume $y=1$.
We now define $\mathcal{L}_{\iota}'\coloneqq \{P^{\iota}_1,\dots,P^{\iota}_{\lfloor\nicefrac{p}{2}\rfloor-s'\ell} \}$.
Consider the path $R = P^{\iota}_{\lfloor\nicefrac{p}{2}\rfloor-s'\ell+1}$.
Still, $R$ separates both $H_{i_x}$, $x\in[2]$, from $P^{\iota}_p$.
Indeed, if we consider the restriction $\rho_j$ of $\rho'$ to the society $(B_j,\Lambda_j)$, then $\rho_j$ is a vortex-free rendition of $(B_j,\Lambda_j)$ in a disk $\Delta_j$ and $R$ is grounded in $\rho_j$.
Let $\gamma$ be the trace of $R$.
Then $\gamma$ separates $\Delta_j$ into a disk $\Delta_{\iota}$ and a disk $\Delta'$ where $\Delta_{\iota}$ contains both $H_{i_1}$ and $H_{i_2}$.
Let $B_{\iota}$ be the subgraph of $B_j$ drawn on $\Delta_{\iota}$ and let $\Lambda_{\iota}$ be obtained from the segment of $\Lambda_j$ whose vertices are drawn on $\Delta_{\iota}$ by extending the order by the natural order of the ground vertices of $R$.
Finally, for each $i\in I_j$ let $\mathcal{L}'_j$ be the restriction of $\mathcal{L}_j$ to $\Delta_{\iota}$ and let $I_{\iota}\coloneqq \{ \iota\}\cup \{ i\in I_j \colon \mathcal{L}'_i\neq\emptyset \}$.
It follows that $(B_{\iota},\Lambda_{\iota},\{\mathcal{L}'_i \colon i\in I_{\iota} \})$ is the desired constriction of $(B_0,\Lambda_0)$.
\end{proof}

Having now found a large transaction that separates two districts $A$ and $B$ in the graph $H$ from the transaction mesh, we want to simultaneously route both $A$ and $B$ back onto the society with linkages that are disjoint from the transaction that separates them. 

\begin{lemma}\label{lemma:UrinHitsTheFan}
Let $t,s,s',a,p,\ell,k$ be non-negative integers with $\max(a,s')\ell\leq s$, $2s+2-s'\ell\leq p$, and $k\in[\lfloor\nicefrac{p}{2}\rfloor-s'\ell]$ such that $\lfloor\nicefrac{p}{2}\rfloor-s'\ell\geq 3k+1$.
Let $\mathfrak{M}=(H,\Omega,\mathcal{C},\mathfrak{P}=\{ \mathcal{P}_1,\dots,\mathcal{P}_{\ell}\})$ be a padded $(s,s',a,p,\ell)$-transaction mesh in the society $(G,\Omega)$ with districts $\{ (G_i,\Omega_i) \}_{i\in[\ell+1]}$ and a cylindrical rendition $\rho$ with nest $\mathcal{C}$.

Let $i_1,i_2\in[\ell+1]$ be distinct and for each $j\in[2]$, let $\mathcal{R}_j$ be a radial linkage in $(G_{i_j},\Omega_{i_j})$ of order $k$.
Let $(H_{i_j},\Omega_{i_j})$ be the area corresponding to the district $(G_{i_j},\Omega_{i_j})$ for both $j\in[2]$.

Finally, let $(H',\Omega')$ the $C_{s-s'\ell}$-society in $\rho$ with respect to $H$ and let $G'$ be the inner graph of $C_{s-s'\ell}$.
Then there exist the following three objects:
\begin{enumerate}
    \item a transaction $\mathcal{Q}$ of order $\lfloor\nicefrac{p}{2}\rfloor-s'\ell-3k$ on $(H',\Omega')$ that separates $H_{i_1}$ and $H_{i_2}$ in $H$, and
    \item two linkages $\mathcal{L}_i$ in $H$, $i\in[2]$, each of order $k$ such that for each $i\in[2]$, $\mathcal{L}_i$ is a linkage between $V(\Omega')$ and the endpoints of $\mathcal{R}_i$ on the inner cycle of $(\mathcal{C}_{i_j},\mathfrak{P}_{i_j})$.
\end{enumerate}
Moreover, $V(\mathcal{Q})$, $V(\mathcal{L}_1)$, and $V(\mathcal{L}_2)$ are pairwise disjoint and the linkages $\mathcal{Q}$, $\mathcal{L}_1$, and $\mathcal{L}_2$ can be found in time $\mathbf{O}(p|E(G)|)$.
\end{lemma}

\begin{proof}
We begin with two applications of \cref{lemma:RadialLinkageOnTheSurface} to the radial linkages $\mathcal{R}_1$ and $\mathcal{R}_2$.
As a result, we obtain two linkages $\mathcal{L}_1'$ and $\mathcal{L}_2'$ such that for each $j\in[2]$, $|\mathcal{R}_j|=|\mathcal{L}_{j}'|$,
the endpoints of $\mathcal{L}_j'$ not on $C{s-s'\ell}$ are precisely the endpoints of $\mathcal{R}_j$ on the inner cycle of $(G_{i_j},\Omega_{i_j})$, and $\mathcal{L}_j'$ is completely contained within the surface stream of the $k$-splitstream of $\mathfrak{M}$.

So far, we have two good candidates for the final linkages we are looking for.
However, we are missing the transaction on $(G',\Omega')$ that separates our two special districts, and at this stage we have no guarantee that $\mathcal{L}_1'$ and $\mathcal{L}_2'$ are disjoint.

Next, we apply \cref{lemma:SplitTwoDistricts} to $\mathfrak{M}$ and the two areas $(H_{i_1},\Omega_{i_1})$ and $(H_{i_2},\Omega_{i_2})$ using the same $k$ to split the stream.
The result is a transaction $\mathcal{Q}'=\{ Q_1,\dots, Q_{ \lfloor\nicefrac{p}{2}\rfloor-s'\ell-k} \}$, indexed naturally such that $Q_1$ separates all other paths in $\mathcal{Q}'$ from $H_{i_1}$ in $H$, which is fully contained in the deep stream of the $k$-splitstream of $\mathfrak{M}$ and which separates $(H_{i_1},\Omega_{i_1})$ and $(H_{i_2},\Omega_{i_2})$ in $H$.

What is left is to adjust $\mathcal{L}_1'$ and $\mathcal{L}_2'$.
We define three subtransactions of $\mathcal{Q}'$ as follows.
\begin{align*}
    \mathcal{Q}_1 & \coloneqq \{ Q_1,\dots,Q_k \}\\
    \mathcal{Q} & \coloneqq \{ Q_{k+1},\dots,Q_{ \lfloor\nicefrac{p}{2}\rfloor-s'\ell-2k} \}\text{, and}\\
    \mathcal{Q}_2 & \coloneqq \{ Q_{ \lfloor\nicefrac{p}{2}\rfloor-s'\ell-2k+1},\dots,Q_{ \lfloor\nicefrac{p}{2}\rfloor-s'\ell-k} \}
\end{align*}
Then $\mathcal{Q}$ is the desired transaction of order $\lfloor\nicefrac{p}{2}\rfloor-s'\ell-3k\geq 1$ as required by the statement of our lemma.

For each $j\in[2]$ let $H'_j$ be the component of $H-\bigcup\mathcal{Q}$ that contains $H_{i_j}$.
Now for each $j\in[2]$ consider the graph
\begin{align*}
F_j\coloneqq \big(\bigcup\mathcal{L}_j' \cap H_j\big) \cup \bigcup\mathcal{Q}_j
\end{align*}
and let $X_j$ be the endpoints of $\mathcal{L}_j'$ on the inner cycle of $(G_{i_j},\Omega_{i_,j})$, while $Y_j$ is the set of remaining endpoints of $\mathcal{L}_j'$ on $V(\Omega')\cap V(H_j)$ together with the set of endpoints of $\mathcal{Q}_j$ (which all belong to $V(\Omega')\cap H_j$).

We claim that for both $j\in[2]$ there is an $X_j$-$Y_j$-linkage of order $k$ in $F_j$.
To see this, suppose there is no such linkage.
Without loss of generality, we may assume that $j=1$ since the arguments are fully symmetric.
Then \cref{prop:mengersthm} yields the existence of a set $S_1$ of order at most $k-1$ such that $F_1$ does not have any $X_1$-$Y_1$-path.
Since $\mathcal{L}_1'$ has order $k$ there is a path $L\in \mathcal{L}_1'$ which is disjoint from $S_1$.
If $L\subseteq H_1$ then we would have a contradiction to our choice of $S_1$, hence $L$ must contain a vertex that does not belong to $H_1$.
It follows, that there exists an $X_1$-$V(Q_{k+1})$-subpath $L'$ of $L$.
Since $(H,\Omega)$ has a vortex-free rendition in the disk we must have that $L'$ intersects every path in $\mathcal{Q}_1$.
With $|\mathcal{Q}_1|=k$ it holds, however, that there is some $P\in\mathcal{Q}_1$ which is disjoint from $S_1$.
Since $P$ contains a vertex of $V(\Omega')\cap V(H_1)$ and $L'$ intersects $P$, we have now found an $X_1$-$Y_1$-path in $F_1-S_1$ which is absurd.
Hence, for each $j\in[2]$, an application of \cref{prop:mengersthm} yields the desired linkage $\mathcal{L}_j$ and our proof is complete.
\end{proof}

The next statement is quite substantial in size, but it intuitively breaks down to the following:
If all of the districts of the transaction mesh already have bounded depth and we attempt to separate them from each other, we will either find a handle or crosscap transaction by utilising \cref{lemma:UrinHitsTheFan}, or we are able to construct the rendition that we seek in the last point of \Cref{thm:societyclassification}.
For the latter option most of the technical aspects are handled by \Cref{lemma:rooting_areas}.

This lemma features all items from \Cref{thm:societyclassification} and will actually be the only way in which we reach the last two points of the theorem in our proof.
The case which produces a handle will give us a lot of grief, as noted at the beginning of this section.

\begin{lemma}\label{lem:everydistrictboundeddepth}
    There exist functions $\mathsf{loss}_{\ref{lem:everydistrictboundeddepth}}$, $\mathsf{nest}_{\ref{lem:everydistrictboundeddepth}} \colon \mathbb{N} \rightarrow \mathbb{N}$,  $\mathsf{link}_{\ref{lem:everydistrictboundeddepth}}\colon\mathbb{N}^3\to\mathbb{N}$, and $\mathsf{apex}_{\ref{lem:everydistrictboundeddepth}} \colon \mathbb{N}^4\to\mathbb{N}$ with $\mathsf{loss}_{\ref{lem:everydistrictboundeddepth}}(x)$, $\mathsf{nest}_{\ref{lem:everydistrictboundeddepth}}(x)$, $\mathsf{link}_{\ref{lem:everydistrictboundeddepth}}(x,y,z)$. $\mathsf{apex}_{\ref{lem:everydistrictboundeddepth}}(w,x,y,z) \in \mathbf{poly}(w+x+y+z)$ for all integers $x,y,z$, such that the following holds.

    Let $t,s,k,a,a',p,d,q$ be non-negative and $\ell$ be positive integers with $p \geq 2s+2-k(\ell+1) + \mathsf{link}_{\ref{lem:everydistrictboundeddepth}}(t,s,q)$, and $s \geq \mathsf{nest}_{\ref{lem:everydistrictboundeddepth}}(t) + a + k(\ell+1) +1$.
    Let $(G,\Omega)$ be a society with a cylindrical rendition $\rho$ in a disk $\Delta$, a cozy nest $\mathcal{C} = \{ C_1, \ldots , C_s \}$ of order $s$ in $\rho$ around the vortex $c_0$, and an $(s,k,a,p,\ell)$-transaction mesh $\mathfrak{M}$, such that the union of the apices for the transactions in $\mathfrak{P}$ have order at most $a'$.
    Further, let $\{(G_i,\Omega_i)\}_{i \in [\ell+1]}$ be the districts of $\mathfrak{M}$, with cylindrical renditions $\rho_i$ and nests $\mathcal{C}_i = \{ C_1^i, \ldots , C_{s - (k\ell + 1)}^i \}$.
    Further, let $(G'', \Omega'')$ be the $C_{s - k(\ell+1)}$-society in $\rho$.

    If for all $i \in [\ell+1]$ the $C_{s - k(\ell+1)}^i$-society $(G_i',\Omega_i')$ of $(G_i,\Omega_i)$ in $\rho_i$ has depth at most $d$, then there exists a set $A \subseteq V(G)$ and  one of the following holds:
    \begin{enumerate}
        \item a $K_t$-minor model in $G$ controlled by a mesh whose horizontal paths are subpaths of distinct cycles from $\mathcal{C}$ (in this case $A=\emptyset$), 
    
        \item we have $|A|\leq 8t^3 + a'$ and there exists a flat crosscap transaction $\mathcal{P}$ of order $q$ in $(G''-A,\Omega'')$, and a nest $\mathcal{C}'$ in $\rho$ of order $s - (\mathsf{loss}_{\ref{lem:everydistrictboundeddepth}}(t)+k(\ell+1)+1)$ around $c_0$ to which $\mathcal{P}$ is orthogonal,

        \item we have $|A|\leq 8t^3 + a'$ and there exists a flat handle transaction $\mathcal{P}$ of order $q$ in $(G''-A,\Omega'')$, and a nest $\mathcal{C}'$ in $\rho$ of order $s - (\mathsf{loss}_{\ref{lem:everydistrictboundeddepth}}(t)+k(\ell+1)+1)$ around $c_0$ to which $\mathcal{P}$ is orthogonal, or

        \item we have $|A| \leq \mathsf{apex}_{\ref{lem:everydistrictboundeddepth}}(t,s,p,q)$ and there exists a rendition $\rho'$ of $(G - A, \Omega)$ in $\Delta$ with breadth $b \in [\nicefrac{1}{2}(t-3)(t-4)-1]$ and depth at most $d$, and an extended $k$-surface-wall $D$ with signature $(0,0,b)$, such that $D$ is grounded in $\rho'$, the base cycles of $D$ are the cycles $C_{s-k\ell-1-k},\dots,C_{s-k\ell-1}$, and there exists a bijection between the vortices $v$ of $\rho'$ and the vortex segments $S_v$ of $D$, where $v$ is the unique vortex contained in the disk $\Delta_{S_v}$ defined by the trace of the inner cycle of the nest of $S_v$ where $\Delta_{S_v}$ is chosen to avoid the trace of the simple cycle of $D$.
    \end{enumerate}

    In particular, the set $A$, the $K_t$-minor model, the transaction $\mathcal{P}$, the rendition $\rho'$, and the extended surface-wall $D$ can each be found in time $\mathbf{poly}(t+k+p+q+s)|E(G)|$.
\end{lemma}

\begin{proof}
We give the following estimate on the functions from the statement of our claims including one function we use for auxiliary purposes.
    \begin{align*}
        \mathsf{link}_{\ref{lem:everydistrictboundeddepth}}(t,s,q) & \coloneqq 2s(q+1) + 16t^3 + 6\big(s\cdot ( (97t^3(q+40))^2+2 )+ 8s +1-1\big)+1\\
        \mathsf{aux\text{-}link}(t,s,p,q) & \coloneqq  p\cdot \big( s\cdot ( (97t^3(q+40))^2+2 )+8s \big)+1\\
        \mathsf{apex}_{\ref{lem:everydistrictboundeddepth}}(t,s,p,q) & \coloneqq 8t^3 + (64t^3 -128)^2\cdot \mathsf{aux\text{-}link}(t,p,q) -16 \in \mathbf{O}(t^{12}psq^2)\\
        \mathsf{nest}_{\ref{lem:everydistrictboundeddepth}}(t) & \coloneqq 16t + 8t^3 + 10\\
        \mathsf{loss}_{\ref{lem:everydistrictboundeddepth}}(t) & \coloneq 8t^3+10
    \end{align*}

    We start by applying \Cref{thm:ApathsArgument}.
    If this yields a $K_t$-minor model then it also has the extra properties we seek and thus we are done.
    As a consequence we may suppose that we instead find sets $S \subseteq V(G')$ and $I = I_1 \cup I_2 \subseteq [\ell+1]$ with $|S| \leq 8t^3-16$, $|I_1| \leq 64t^3-128$, and $|I_2| < \nicefrac{1}{2}(t-3)(t-4)$ such that each $\{ V(\Omega_i) ~\!\colon\!~ i \in [\ell+1] \}$-path in $G' - S$ is a $\{ V(\Omega_i) ~\!\colon\!~ i \in I_1 \}$-path and each district $(G_j - S,\Omega_j)$ with $j \in [\ell+1] \setminus (I_1\cup I_2)$ has a vortex-free rendition in a disk.
    
    Let $G' \coloneqq \bigcup_{i \in [\ell+1]} G_i'$.    
    For all $i\neq j \in I_1$, we apply \hyperref[prop:mengersthm]{Menger's theorem} to find a $V(\Omega_i')$-$V(\Omega_j')$-linkage of order at least $\mathsf{aux\text{-}link}(t,s,p,q)$ in $G'$.

    Let us first discuss the case where none of these applications of \cref{prop:mengersthm} yield the desired linkage.
    
    In this case, we may assume that after at most $|I_1|^2$ iterations, we have found a $V(\Omega_i')$-$V(\Omega_j')$-separator $S_{i,j}$ of size a most $\mathsf{aux\text{-}link}(t,s,p,q)$ for any two distinct $i,j \in I_1$.
    Let $S' \coloneq S \cup \bigcup_{i \in I_1} \bigcup_{j \in I_1 \setminus \{ i \}} S_{i,j}$ and note that we have $|S| \leq |S| + \mathsf{aux\text{-}link}(t,s,p,q)|I_1|^2 \leq \mathsf{apex}_{\ref{lem:everydistrictboundeddepth}}(t,s,p,q)$.
    Moreover, it follows that in $G' - S'$ there do not exists any $V(\Omega_i')$-$V(\Omega_j')$-paths for $i,j \in [\ell+1]$ with $i \neq j$.

    We now update our graphs as follows.
    For each $i\in[\ell+1]$ let $\Omega_i''$ be the restriction of $\Omega_i'$ to $V(\Omega_i')\setminus S'$ and let $G_i''$ be the union of all components of $G_i'-S'$ that contain a vertex of $V(\Omega_i'')$.
    Notice that for each $i \in [\ell+1]$, both $V(\Omega_i'')$ separates $V(\Omega)$ from $V(G_i'')$ in $G'-S'$.

    Since each district $(G_i - S,\Omega_i)$ with $i \in [\ell+1] \setminus I_2$ has a vortex-free rendition in a disk, this remains true for $(G_i'',\Omega_i'')$.

    Let $\{(H_i,\Omega_i)\}_{i \in [\ell+1]}$ be the areas of $\mathfrak{M}$.
    Further, let $\{(H_i',\Omega_i')\}_{i \in [\ell+1]}$ be the family of $C_{s - s'(\ell+1)}^i$-societies in the restrictions of $\rho_i$ to $H_i$ for each $i \in [\ell+1]$.
    From \Cref{def:transactionmesh} we know that the society $(H - \bigcup_{i \in [\ell+1]} (H_i - V(\Omega_i)), \Omega)$ has a vortex-free rendition $\rho^*$ in a disk $\Delta$ in which $\bigcup_{i \in [\ell+1]} V(\Omega_i) \subseteq N(\rho^*)$, such that $\rho^*$ is a vortex-free rendition of $(H - \bigcup_{i \in [\ell+1]} (H_i - V(\Omega_i)), \Omega_i)$ in a disk for each $i \in [\ell+1]$.
    Using the renditions $\rho_i$ for all $i \in [\ell+1]$, we can extend $\rho^*$ to a rendition $\rho^{**}$ of $(H - \bigcup_{i \in [\ell+1]} (H_i' - V(\Omega_i')), \Omega)$ in $\Delta$ in which $\bigcup_{i \in [\ell+1]} V(\Omega_i') \subseteq N(\rho^{**})$, such that $\rho^{**}$ is a vortex-free rendition of $(H - \bigcup_{i \in [\ell+1]} (H_i' - V(\Omega_i')), \Omega_i')$ in a disk for each $i \in [\ell+1]$.

    Let $H^* \coloneqq ( H_i' - V(\Omega_i'')) \cup \bigcup_{i \in [\ell+1] \setminus I_2} G_i'' ) - S'$ and set $X\coloneqq \bigcup_{i\in[\ell+1]}(V(G_i)\setminus (V(G_i'')\cup S'))$.
    Notice that every neighbour of $X$ in $G$ must belong to $S'$ by construction.
    The rendition $\rho^{**}$ can then be further augmented to a vortex-free rendition $\rho^{***}$ of $(H^*, \Omega)$ in a $\Delta$ in which $\bigcup_{i \in I_2} V(\Omega_i'') \subseteq N(\rho^{***})$, such that $\rho^{***}$ is a vortex-free rendition of $(H^*, \Omega_i')$ in a disk for each $i \in I_2$.
    Finally, to construct the rendition $\rho'$ of $(G - S', \Omega)$ in a $\Delta$ from the fourth item of our statement, we can now add each society $(G_i'', \Omega_i'')$ with $i \in I_2$ into $\rho^{***}$ as a vortex and draw the graph $G[X]$ in the interior of an arbitrary cell.
    Since $|I_2| < \nicefrac{1}{2}(t-3)(t-4)$, the rendition $\rho'$ has breadth at most $\nicefrac{1}{2}(t-3)(t-4) -1$ and depth at most $d$.
    The extended $k$-surface wall we seek is then recovered by applying \Cref{lemma:rooting_areas} to the collection $\{ (H_i,\Omega_i) \}_{i \in I_2}$.

    We continue with the case where we find a $V(\Omega'_{i_1})$-$V(\Omega'_{i_2})$-linkage $\mathcal{P}_1$ of order $\mathsf{aux\text{-}link}(t,s,p,q)$ in $G''$ for some $i_1\neq i_2\in I_1$.
    Set $J\coloneqq [\ell+1]\setminus \{ i_1,i_2\}$.

    In this case we return the set $S$ and all previously found sets $S_{i,j}$ to the graph.
    We may now apply \cref{lemma:orthogonal_connecting_linkage} to $\mathcal{P}_1$ and the societies $(G_{i_1}',\Omega_{i_1}')$, $(G_{i_2}',\Omega_{i_2}')$, and $\{ (G_j',\Omega_j') \}_{j\in J}$ together with their corresponding paddings as provided by $\mathfrak{M}$.
    As a result, we obtain a connecting linkage $\mathcal{P}_2$ of order $(97t^3(p(q+1)+41))^2+8s+1$ such that $\mathcal{P}_2$ is orthogonal to the nests $\{ C_2^{i_x},\dots,C_{s-k(\ell+1)-1}^{i_x}\}$ for both $x\in[2]$ and disjoint from the cycles $C^j_2$ for all $j\in J$.

    Next, for each $x\in[2]$ let $\mathcal{L}_x$ be the radial linkage in $(G'_{i_x},\Omega'_{i_x})$ obtained by following along the paths of $\mathcal{P}_2$ starting on $V(\Omega_{i_x}')$ and ending when $C_2^{i_x}$ is encountered for the first time.
    Let $X_x$ be the set of endpoints of the linkage $\mathcal{L}_x$ on $C_2^{i_x}$ for each $x\in[2]$ and let $\mathcal{M}$ be the set of all $X_1$-$X_2$-paths found in the paths of $\mathcal{P}_2$.
    We then apply \cref{lemma:UrinHitsTheFan} to the two districts $(G_{i_1}',\Omega_{i_1}')$ and $(G_{i_2}',\Omega_{i_2}')$.
    This results in two linkages $\mathcal{R}_1$ and $\mathcal{R}_2$ of order $(97t^3(p(q+1)+41))^2+8s+1$, such that $\mathcal{R}_x$ connects $V(\Omega')$ to $X_x$ for both $x\in [2]$.
    Moreover, \cref{lemma:UrinHitsTheFan} returns a planar transaction $\mathcal{Q}_1$ in $(G'',\Omega'')$ of order
    \begin{align*}
        & \left\lfloor \frac{p}{2} \right\rfloor - k(\ell+1) - 3\cdot\big( s\cdot ( (97t^3(q+40))^2+2 )+ 8s +1\big)\\
        & \geq s(q+1)+1 + 8t^3
    \end{align*}
    such that $\mathcal{Q}_1$ separates $V(\Omega_{i_1}')$ and $V(\Omega_{i_2}')$ in $H \cap G''$.
    Additionally, we have that $V(\mathcal{Q}_1)$, $V(\mathcal{R}_1)$, and $V(\mathcal{R}_2)$ are pairwise disjoint.
    
    Notice that $\mathcal{M}$ shares at most its endpoints with the graph $H$ and is therefore, in particular, vertex-disjoint from $\mathcal{Q}_1$.
    Indeed, we may now reunite $\mathcal{M}$ with the linkages $\mathcal{R}_1$ and $\mathcal{R}_2$ to obtain a transaction $\mathcal{P}_3$ of order $s\cdot ( (97t^3(q+40))^2+2 )+ 8s +1$ on $(G'',\Omega'')$ that is completely disjoint from the transaction $\mathcal{Q}_1$.
    Moreover, every path from $\mathcal{P}_3$ forms a cross on $(G'',\Omega'')$ with any path from $\mathcal{Q}_1$.

    At this point we have almost reached our goal.
    The transaction $\mathcal{Q}_1$ is already fully contained in the vortex-free rendition of $(H,\Omega)$ and therefore its strip society is flat and isolated in $(G'',\Omega'')$.
    Our next goal is to replace both $\mathcal{Q}_1$ and $\mathcal{P}_3$ with transactions that are orthogonal to most of the nest of $(G'',\Omega'')$.
    This requires some amount of set up.

    From here on we consider $\rho''$ to be the vortex-free rendition of $(H''\coloneqq H\cap G'',\Omega'')$ into the disk $\Delta$ as induced by $\rho$.
    
    Let $A'_1$ and $A'_2$ be the two minimal segments of $\Omega''$ such that $\mathcal{Q}_1$ is an $A_1'$-$A_2'$-linkage.
    Similarly, let $B'_1$ and $B_2'$ be the two minimal segments of $\Omega''$ such that $\mathcal{P}_3$ is a $B_1'$-$B_2'$-linkage.
    Then, let $A_1$ and $A_2$ be the sets of the first vertices encountered on $C_2$ when following along the paths in $\mathcal{Q}_1$ starting in $A_1'$ and $A_2'$ respectively.
    Similarly, let $B_1$ and $B_2$ be the sets of the first vertices encountered on $C_2$ when following along the paths in $\mathcal{P}_3$ starting in $B_1'$ and $B_2'$ respectively.
    For each $i\in[2]$ let $I_i'$ be the shortest subpath of $C_2$ that contains all vertices of $A'_i$ and is disjoint from $A_{3-i}'$.
    
    Finally, let $I_i$ be the shortest subpath of $I_i'$ that contains all vertices of $\mathcal{Q}_1$ that belong to $I_i'$.

    \begin{claim}\label{claim:entryAndExit}
    There exist $\mathcal{S}_{i,j}\subseteq \mathcal{P}_3$ such that $I_i$ does not contain any vertices from $V(\mathcal{R}_{j})\setminus V(\mathcal{S}_{i,j})$.
    Moreover, we have that $|\mathcal{S}_{i,j}|\leq s$.
    \end{claim}
    
  \emph{Proof of \Cref{claim:entryAndExit}:}
    Let us first discuss the set $\mathcal{S}_{i,j}$.
    For each $x\in[2]$ let $b_x$ be the endpoint of $I_i$ that is closest to $B_x$ along $I_i'$.

    We greedily construct a sequence $R_1,\dots,R_y$ of paths from $\mathcal{R}_j$ such that for every $z\in[y]$, $R_1$ contains the vertex of $I_i$ that is closest to $b_{3-j}$ along $I_i$ among all paths in $\mathcal{R}_j\setminus \{ R_x  \colon x\in[z-1]\}$.
    Let $y$ be as large as possible.
    Notice that we may assume $y\geq 1$ as otherwise $\mathcal{R}_j$ does not intersect $I_i$.

    For each $z\in[y]$ we construct a $V(\Omega'')$-path $T_z$ and a subdisk $\Delta_z$ of $\Delta$ such that $\Delta_{z-1}\subseteq \Delta_z$ (here $\Delta_0 = \emptyset$) as follows.

    Start with $z=1$.
    Let $v_1$ be the vertex of $I_i$ that belongs to $\mathcal{R}_j$ which is closest to $b_{3-j}$ along $I_i$ and let $R_1\in\mathcal{R}_j$ be the path containing $v_1$.
    We follows along $R_1$ starting in $B_j'$ until we meet $v_1$, let $R'_1$ be the resulting subpath of $R_1$.
    Now we follow along $I_i$ until one of the following happens:
    Let us no fix $u$ to be one of the following:
    either, $u\in A_i$ is an endpoint of some $V(\Omega'')$-$A_i$-subpath $L$ of a path in $\mathcal{Q}_1$, or $u=b_{3-j}$ and is the endpoint of some path $B_{3-j}$-$B_{3-j}'$-subpath $L$ of some path in $\mathcal{R}_{3-j}$.
    Notice that we call the path we found $L$ in either case as we only use this path to define the disks $\Delta_z$ and we will always use the same $L$.
    Let $M_1$ be the $v_1$-$u$-subpath of $I_i'$.
    Then $T_1\coloneqq R'_1\cup M_1\cup L$ is a grounded $V(\Omega'')$-path.
    Moreover, $T_1$ separates $\Delta$ into two disks, one of which contains the vertices of $A_{3-j}$.
    Let $\Delta_1$ be the disk disjoint from $A_{3-j}$.

    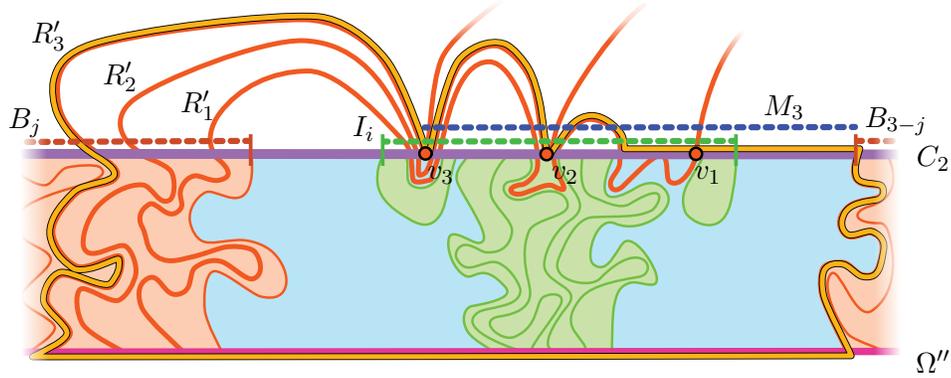
\begin{figure}[ht]
    \centering
    \begin{tikzpicture}

        \pgfdeclarelayer{background}
		\pgfdeclarelayer{foreground}
			
		\pgfsetlayers{background,main,foreground}

        \begin{pgfonlayer}{background}
        \node (C) [v:ghost] {{\includegraphics[width=13cm]{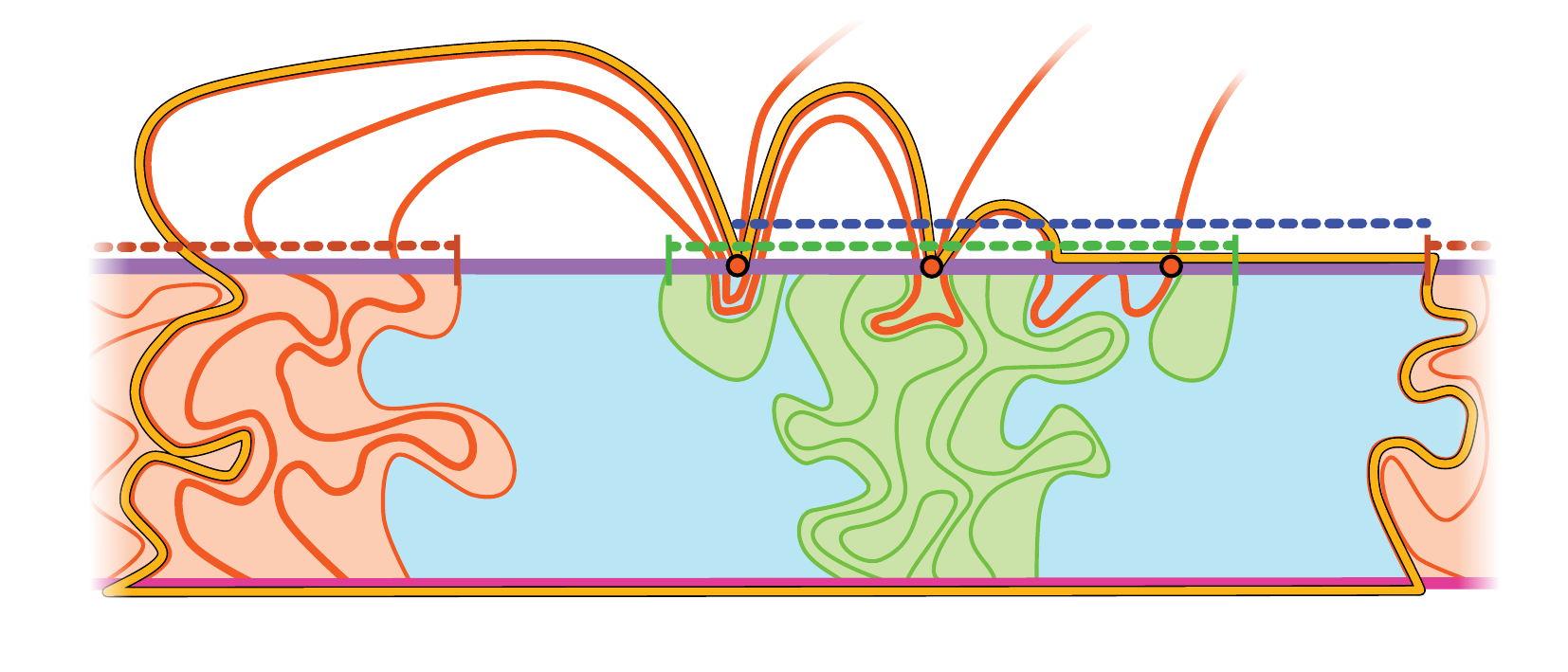}}};
        \end{pgfonlayer}{background}
			
        \begin{pgfonlayer}{main}

            \node (X) [v:ghost,position=0:0.3mm from C] {};

            \node (Omega) [v:ghost,position=341:67mm from X] {$\Omega''$};

            \node (C2) [v:ghost,position=90:27mm from Omega] {$C_2$};

            \node (B3j) [v:ghost,position=135:7mm from C2] {$B_{3-j}$};

            \node (M3) [v:ghost,position=173:15mm from B3j] {$M_3$};

            \node (v1) [v:ghost,position=196:26mm from B3j] {$v_1$};
            \node (v2) [v:ghost,position=180:19mm from v1] {$v_2$};
            \node (v3) [v:ghost,position=180:16.5mm from v2] {$v_3$};

            \node (Ii) [v:ghost,position=149:12mm from v3] {$I_i$};

            \node (Bj) [v:ghost,position=179:45mm from Ii] {$B_j$};

            \node (R3) [v:ghost,position=75:12mm from Bj] {$R_3'$};
            \node (R2) [v:ghost,position=330:11mm from R3] {$R_2'$};
            \node (R1) [v:ghost,position=340:11mm from R2] {$R_1'$};
        
        \end{pgfonlayer}{main}

        \begin{pgfonlayer}{foreground}
        \end{pgfonlayer}{foreground}

    \end{tikzpicture}
    \caption{A diagram of the iterative construction in the proof of \cref{claim:entryAndExit} within the proof of \cref{lem:everydistrictboundeddepth}.
    The picture shows the situation after the vertex $v_3$ together with the paths $R_3'$ and $M_3$ (indicated by the \textcolor{MidnightBlue}{dark blue} dashed line) have been constructed.
    The bold \textcolor{BananaYellow}{yellow} line indicates the boundary of the disk $\Delta_3$.
    Notice how, below the position of $v_3$ all previously constructed paths are forced to enter the nest and leave it again.
    This accumulation of $V(C_2)$-paths is the cause for our bound.}
    \label{fig:orthogonaliseHandle}
\end{figure}

    Now suppose we have already constructed the paths $R'_1,\dots,R'_z$, $M_1,\dots,M_z$, $T_1,\dots,T_z$, and the disks $\Delta_1,\dots,\Delta_z$ for some $z\in [1,y]$.

    Suppose for the moment that $z\neq y$.
    Then let Let $v_{z+1}$ be the vertex of $I_i$ that belongs to $\mathcal{R}_j\setminus\{ R_x \colon x\in [z]\}$ which is closest to $b_{3-j}$ along $I_i$ and let $R_{z+1}\in\mathcal{R}_j\setminus\{ R_x \colon x\in [z]\}$ be the path containing $v_{z+1}$.
    We follows along $R_{z+1}$ starting in $B_j'$ until we meet $v_{z+1}$, let $R'_{z+1}$ be the resulting subpath of $R_{z+1}$.
    From here our construction diverges from the case $z=1$.
    
    If $R'_{z+1}$ is contained in $\Delta_z$ we set $\Delta_{z+1}\coloneqq \Delta_z$, $M_{z+1}\coloneqq M_z$, and $T_{z+1}\coloneqq T_{z}$.
    In this case we are done with the construction step for $z+1$.
    
    So we may assume that $R'_{z+1}$ is not contained in $\Delta_{z}$. 
    We now let $M_{z+1}$ be the $v_{z+1}$-$b_{3-j}$-subpath of $I_i$.
    Second, let us follow along $I_i$ towards $b_{3-j}$ until we encounter one of the previously constructed $T_x$.
    Since $\rho''$ is a vortex-free rendition of $(H'',\Omega'')$ in $\Delta$ we must have that $x=z$ by our construction and because $R'_{z+1}$ is not contained in $\Delta_z$.
    We now follow along $T_z$ until we meet the endpoint of $L$ on $\Omega''$.
    The resulting $V(\Omega'')$-path is set to be $T_{z+1}$ and we define $\Delta_{z+1}$ to be the subdisk of $\Delta$ which contains the trace of $T_{z+1}$ in its boundary and fully contains the disk $\Delta_z$.
    This completes the construction.
    See \cref{fig:orthogonaliseHandle} for an illustration.

    So now let us discuss the case where $z=y$.
    We claim that $y\leq s$.
    By construction, every path $R'_x$ with $x\in [y]$ has its endpoint on $\Omega''$ in $\Delta_y$.
    Since $\mathcal{R}_j$ is contained in the surface stream of the splitstream used by our original application of \cref{lemma:UrinHitsTheFan}, while $\mathcal{Q}'$ must be contained in the corresponding deep stream, it follows that all vertices of $I_i$ belong to said deep stream and thus none of the paths $R'_{x}$ can have its other endpoints in $\Delta_{y}$.
    However, this means that each $R_x$, $x\in[y-1]$, must leave $\Delta_{y}$ and so it must eventually pass through the path $M_{y}$ for the last time.
    Let $w$ be the endpoint of $M_{y}$ that is closest to $b_j$ on $I_i$.
    Since $I_i$ is disjoint from $B_j$ we know that, with exception for the smallest index $z'$ index where $M_y=M_{z'}$, $R'_x$, $x\in[y]\setminus \{z'\}$, contains a grounded $V(C_2)$-path $U_x$ that is contained in the outer graph of $C_2$ and whose endpoints are separated by $w$ on $I_i'$.
    Now, with $\mathcal{C}$ being cozy, we have that $y-1\leq s-1$ and thus, $y\leq s$ as desired.

    We set $\mathcal{S}_{i,j}=\{ P\in\mathcal{P}_3 \colon \text{there is }z\in[y]\text{ with }R_z\subseteq P\}$ and observe that $I_i$ does not contain any vertices from $V(\mathcal{R}_{j})\setminus V(\mathcal{S}_{i,j})$ by the maximality of $y$.
    This completes the proof of \Cref{claim:entryAndExit}.
    \hfill$\blacksquare$
    \medskip
    
    Now let, for each $j\in[2]$, $J_j$ be the $V(I_1')$-$V(I_2')$-subpath of $C_2$ that contains $B_j$.
    Notice that the only place where $\mathcal{P}$ can intersect $C_2$ outside of $I_1$ and $I_2$ are $J_1$ and $J_2$.
    We are about to change that.

    As before, let us fix $i,j\in[2]$ for the following arguments.

    \begin{claim}\label{claim:entryAndExit2}
    There exist $\mathcal{F}_{i,j}\subseteq \mathcal{P}_3\setminus \bigcup_{i',j'\in[2]}\mathcal{S}_{i',j'}$ such that the minimal subpath $J_j'\subseteq J_j$ that contains all vertices of $V(\mathcal{R}_j)\setminus \big(V(\mathcal{F}_{i,j})\cup \bigcup_{i',j'\in[2]}V(\mathcal{S}_{i',j'})\big)$ does not contain any vertices of $\mathcal{Q}_1$ that can be reached by a subpath of some path in $\mathcal{Q}_1$ starting on $I_i$ without intersecting $I_{3-i}$.
    Moreover, we have that $|\mathcal{F}_{i,j}|\leq s$.
    \end{claim}
    
    \emph{Proof of \Cref{claim:entryAndExit2}:}
    For $x\in [2]$ let $a_x$ be the endpoint of $J_j$ that belongs to $A_x$.
    Now let $v$ be the vertex of $J_j$ that is closest to $a_{3-i}$ and contained in some path $Q\in \mathcal{Q}_1$ such that there exists an $I_i$-$v$-subpath $Q''$ of $Q$ that is disjoint from $I_{3-i}$.
    From the endpoint of $Q''$ traverse along $I_i$ towards $A_i$ and from there follow along a subpath of some path in $\mathcal{Q}_1$ to $V(\Omega'')$, let $Q'$ be the resulting $v$-$V(\Omega'')$-path.
    Notice that $Q'$ is disjoint from $I_{3-i}$.
    Then traverse along $J_j$, starting in $v$ and towards $a_i$, until we encounter the first vertex $b\in B_j$ and let $M$ be the subpath of $J_j$ traversed this way.
    Let $R$ be the corresponding $b$-$V(\Omega'')$-path contained in some path of $\mathcal{R}_j$.
    Now notice that $Q'\cup M\cup R$ is a grounded $V(\Omega'')$-path whose trace splits $\Delta$ into two disks, one of which, say $\Delta_j$, is disjoint from $B_{3-j}$.
    
    Let $\mathcal{F}_{i,j}$ be the collection of all paths from $\mathcal{P}_3$ that contain a vertex drawn on $\Delta_j$.

    Notice, that every vertex of $B_j$ that belongs to the $v$-$a_i$-subpath of $J_j$ belongs to $\Delta_j$ since $\mathcal{R}_j$ and $\mathcal{Q}_1$ are vertex-disjoint.
    Moreover, every path $X$ in $\mathcal{F}_{i,j}$ contains a subpath $Y$ with one end in $\Delta_j$ and the other in $V(C_2^{i_j})$.
    Since $\Delta_j$ contains a vertex of $\mathcal{Q}'$ which must be contained in the deep stream of the splitstream we envoked when we called \cref{lemma:UrinHitsTheFan} earlier, this other endpoint of a path $Y$ as above does not belong to $\Delta_j$.
    It follows that every path $X\in\mathcal{F}_{i,j}$ contains a $V(C_2)$ subpath $U$ which is contained in the outer graph of $C_2$ in $\rho''$ and whose endpoints are connected by a subpath of $J_j$ that contains $v$.
    Therefore, as $\mathcal{C}$ is cozy, we obtain $|\mathcal{F}_{i,j}|\leq s$ as desired.

    Moreover, by choice of $v$ we have that every vertex of $V(J_j)\cap V(\mathcal{Q}_1)$ that was reachable from $I_i$ with a subpath of some path in $\mathcal{Q}_1$ without intersecting $I_{3-i}$ must be drawn on $\Delta_j$.
    Therefore, the proof of our claim is complete.
    \hfill$\blacksquare$
    \medskip
    
    Let $\mathcal{P}_4 \coloneqq \mathcal{P}_3 \setminus \bigcup_{i,j\in[2]}\big( \mathcal{F}_{i,j}\cup\mathcal{S}_{i,j}\big)$.
    Notice that, by \cref{claim:entryAndExit} and \cref{claim:entryAndExit2} we have that
    \begin{align*}
        |\mathcal{P}_4| \geq |\mathcal{P}_3| - 8s \geq s\cdot ( (97t^3(q+40))^2+2 )+1.
    \end{align*}

    Next, for each $j\in[2]$ let $B_j''$ be the subset of $B_j$ containing all vertices that belong to $\mathcal{P}_4$.
    Moreover, for each $i\in[2]$ let now $N'_i$ be the $B_1''$-$B_2''$-subpath containing $A_i$ and let $N_i$ be the minimal subpath of $N_i'$ that contains all vertices of $N_i'$ that belong to $\mathcal{Q}'$.
    Notice that, as a consequence of \cref{claim:entryAndExit} and \cref{claim:entryAndExit2}, every vertex of $\mathcal{Q}_1$ in $C_2$ must now be contained in $N_1 \cup N_2$.

    Consider the nest $\mathcal{N}_1 \coloneqq \{ C_2,\dots,C_{s-k(\ell+1)-1}\}$ and recall that $\mathcal{P}_4$ must be an exposed transaction in the restriction $\rho'''$ of $\rho$ to $(G'',\Omega'')$.
    We have that still, $c_0$ is the unique vortex of $\rho'''$ and now $\mathcal{N}_1$ is a cozy nest around $c_0$ in $\rho'''$.

    We apply \cref{lemma:orthogonal_transaction} to $\mathcal{P}_4$, the nest $\mathcal{N}_1$ and the rendition $\rho'''$ of the society $(G''\Omega'')$.
    Let $\mathcal{P}_5$ be the resulting transaction of order at least $(97t^3(q+40))^2+3$ which is orthogonal to $\mathcal{N}\setminus \{ C_2\}$.

    Since, by iii) of \cref{lemma:orthogonal_transaction} we know that $\mathcal{P}_5$ intersects the union of the inner graph of $C_2$ in $\rho'''$ and $C_2$ only in $C_2\cup\bigcup\mathcal{P}_4$, we obtain from the discussion above that any vertex that $\mathcal{P}_5$ may share with $\mathcal{Q}_1$ must belong to $N_1\cup N_2$.
    However, since both $N_1$ and $N_2$ are disjoint from $\mathcal{P_4}$, all vertices in $V(N_1\cup N_2)$ had degree $2$ in $C_2\cup\bigcup\mathcal{P}_4$.
    It follows that if there is a path $P\in\mathcal{P}_5$ that contains a vertex of $N_i$ for $i\in[2]$, then $N_i\subseteq P$.
    Therefore, we may now obtain a transaction $\mathcal{P}_6$ of order at least $(97t^3(q+40))^2+1$ by discarding at most two paths from $\mathcal{P}_5$ such that $\mathcal{P}_6$ is vertex-disjoint from $\mathcal{Q}_1$.

    Next, we call upon \cref{lem:monotonetransaction} to find a monotone transaction $\mathcal{P}_7\subseteq \mathcal{P}_6$ of order at least $97t^3(q+40)+1$.
    With $\mathcal{P}_7$ being orthogonal to $\mathcal{N}_1\setminus \{ C_2\}$ we may now evoke \cref{lemma:findflattransaction}.
    This results either in a $K_t$-minor controlled by the mesh defined by $\mathcal{P}_7$ and $\mathcal{N}_1\setminus\{ C_2\}$ -- in which case we are done -- or in a set $A\subseteq V(G'')$ such that the strip society of a transaction $\mathcal{P}_8\subseteq \mathcal{P}_7$ of order at least $q$ is flat and isolated in $(G''-A,\Omega'')$.
    We may assume to find the set $A$ and the transaction $\mathcal{P}_8$.
    
    We discard at most $8t^3$ many paths from $\mathcal{Q}_1$ to obtain a new transaction $\mathcal{Q}_2$ of order at least $s(q+1)+1$.

    Suppose the $\mathcal{P}_8$-strip in $(G''-A,\Omega'')$ contains a vertex of some path from $\mathcal{Q}_2$.
    Since $\mathcal{Q}_2$ is disjoint from $A$ and a transaction in $(G''-A,\Omega)$ that is also disjoint from $\mathcal{P}_8$ this would imply that the $\mathcal{P}_8$-strip in $(G'',\Omega'')$ contains a vertex of $A_1$.
    As $A_1$ is necessarily disjoint from the end segments of $\mathcal{P}_8$ on $\Omega''$ this is a contradiction to the strip society of $\mathcal{P}_8$ being isolated in $(G''-A,\Omega'')$.

    Next, let $N_3$ and $N_4$ be the two subpaths of $C_{8t^3+2}$ that are contained in the $\mathcal{P}_8$-strip in $(G''-A,\Omega'')$.
    By the discussion above we have that $N_3$ and $N_4$ are disjoint from the paths in $\mathcal{Q}_2$.
    We now apply \cref{lemma:orthogonal_transaction} to $\mathcal{Q}_2$ and obtain a planar transaction $\mathcal{Q}_3$ of order $q+2$ which is orthogonal to $\mathcal{N}_2 \coloneqq \mathcal{N}_1\setminus \{ C_2,\dots,C_{8t^3+2}\}$.
    Similar to before, the only way how $\mathcal{Q}_3$ can intersect $\mathcal{P}_8$ is via the paths $N_3$ and $N_4$.
    Moreover, if some path of $\mathcal{Q}_3$ contains a vertex of $N_i$ for some $i\in[2]$ then it must contain the entire path $N_i$.
    Hence, by removing at most two paths from $\mathcal{Q}_3$ we obtain a planar transaction $\mathcal{Q}_4$ of order at least $q$ which is orthogonal to $\mathcal{N}_2$.
    Notice that the strip society of $\mathcal{Q}_4$ in $(G''-A,\Omega'')$ is flat and isolated.

    We then apply \cref{lem:reconciliation} to the nest $\mathcal{N}_2$ and the transaction $\mathcal{P}_8$ in $\rho'''$ together with the two boundary paths of $\mathcal{Q}_4$.
    In case $\mathcal{P}_8$ is a crosscap transaction, the proof, at this point, is complete.
    Otherwise, $\mathcal{P}_8$ is planar.
    In this case notice that we have now found a vortex-free rendition $\rho^{\star}$ of the union of the outer graph of $C_{8t^3+9}$ together with the $\mathcal{P}_8$-strip in $(G''-A,\Omega'')$ and the $\mathcal{Q}_4$-strip in $(G''-A,\Omega)$.
    Hence, we have found an orthogonal and flat handle transaction of the required order and our proof is complete.
\end{proof}

\subsection{Building a transaction mesh and updating the \texorpdfstring{$\mathfrak{M}$}{M}-tree}\label{sec:buildingtransactionmesh}
For our proof of \Cref{thm:societyclassification} we will need a tool that lets us extend the structure of a transaction mesh for one more iteration.
These tools must be applicable whenever we cannot apply our results from \Cref{sec:boundeddepth}, which means that for this section we can generally assume that our transaction meshes have a district that does not have bounded depth.
To allow us to integrate the transactions we can find there, we will need to do a lot of clean-up work.
This is one of the most technical parts of our proof as we will need to interact with the details in \Cref{def:transactionmesh} and \Cref{def:Mtree}.
In addition to this, we have finally reached the point at which we need \Cref{thm:orthogonalexposedsorcery}.

We start with a lemma that is simple in spirit, but has a very complicated statement.
Essentially we want to ensure that, if we have some district with sufficient depth, we can either find a $K_t$-minor, a flat crosscap transaction on the original society, or we can find a planar, flat transaction which we can integrate into our transaction mesh.
We delay the integration into the transaction mesh to the lemma following this one.
The third item can best be understood by directly comparing it to the items \textbf{\textsf{M2}} to \textbf{\textsf{M5}} of \Cref{def:transactionmesh}.
The proof of this lemma brings together most of the lemmas and theorems we have proven up to this point.

\begin{lemma}\label{lem:freeflattransactionorprotectedcross}
    There exist polynomial functions $\mathsf{loss}_{\ref{lem:freeflattransactionorprotectedcross}} \colon \mathbb{N} \rightarrow \mathbb{N}$, $\mathsf{apex}^\mathsf{link}_{\ref{lem:freeflattransactionorprotectedcross}}, \mathsf{apex}^\mathsf{genus}_{\ref{lem:freeflattransactionorprotectedcross}} \colon \mathbb{N}^2 \rightarrow \mathbb{N}$, $\mathsf{nest}_{\ref{lem:freeflattransactionorprotectedcross}} \colon \mathbb{N}^5 \rightarrow \mathbb{N}$, and $\mathsf{depth}_{\ref{lem:freeflattransactionorprotectedcross}} \colon \mathbb{N}^6 \rightarrow \mathbb{N}$, such that the following holds.

    Let $t,s,s',a,a_1,a_2,p,p',q,\ell$ be non-negative integers with $\mathsf{nest}_{\ref{lem:freeflattransactionorprotectedcross}}(t,a_1,a_2,s',\ell) \leq s$ and $2s+2 - s'\ell \leq q \leq s(q+1)+1 \leq p$.
    
    Let $(G,\Omega)$ be a society with a cylindrical rendition $\rho$ in a disk $\Delta$, a cozy nest $\mathcal{C} = \{ C_1, \ldots , C_s \}$ of order $s$ in $\rho$ around the vortex $c_0$, and a padded $(s,s',a,p,\ell)$-transaction mesh $\mathfrak{M} = (H,\Omega,\mathcal{C},\mathfrak{P})$.
    Further, let $\{(G_i,\Omega_i)\}_{i \in [\ell+1]}$ be the districts of $\mathfrak{M}$ each with a cylindrical rendition $\rho_i$ in a disk around the vortex $c_i$ and an $(s - (s'\ell + 1),\lfloor \nicefrac{p}{2} \rfloor - (s'\ell + 1),\ell_i)$-padding $(\mathcal{C}_i = \{ C_1^i, \ldots , C_{s - (s'\ell + 1)}^i \}, \mathfrak{P}_i)$ with the inner nest $\mathcal{C}_i' = \{ C_1^{i*}, \ldots , C_{s - (s'\ell + 1)}^{i*} \}$.
    Let $Z \subseteq V(G) - V(H)$ with $|Z| \leq a_1$ and let the union of the apices for the transactions in $\mathfrak{P}$ have order at most $a_2$.
    Let $i \in [\ell+1]$ be chosen such that $(G_i'-Z,\Omega_i')$ has depth at least $\mathsf{depth}_{\ref{lem:freeflattransactionorprotectedcross}}(t,s,q,p',\ell,a_1)$, where $(G_i',\Omega_i')$ is the $C_{s - (s'(\ell+1))}^{i*}$-society of $(G_i,\Omega_i)$.
    Finally, let $G' \coloneqq \bigcup_{i \in [\ell+1]} G_i$ and let $(G'', \Omega'')$ be the $C_{s - s'\ell}$-society in $\rho$.

    There there exists a set $A \subseteq V(G'')$, and
    \begin{enumerate}
        \item a $K_t$-minor model in $G$ controlled by a mesh whose horizontal paths are subpaths of distinct cycles from $\mathcal{C}$,
    
        \item a flat crosscap transaction $\mathcal{P}$ of order $q$ in $(G'' - (A \cup Z),\Omega'')$, whilst $|A| \leq \mathsf{apex}^\mathsf{genus}_{\ref{lem:freeflattransactionorprotectedcross}}(t,a_2)$, and a nest $\mathcal{C}'$ in $\rho$ of order $s - \mathsf{loss}_{\ref{lem:freeflattransactionorprotectedcross}}(t)$ around $c_0$ to which $\mathcal{P}$ is orthogonal, or

        \item there exists a planar, isolated, flat transaction $\mathcal{P}$ of order $p'$ in $(G_i' - (A \cup Z),\Omega_i')$, with $|A \cup Z| \leq \mathsf{apex}^\mathsf{link}_{\ref{lem:freeflattransactionorprotectedcross}}(t,a_1)$ such that:
            \begin{enumerate}
                \item There exist co-conspirators $\mathcal{W}$ for $\mathcal{P}$ in $(G_i',\Omega_i')$ (on level 3).

                \item The $\mathcal{P}$-society $(U_i,\Psi_i)$ in $(G_i' - (A \cup Z),\Omega_i')$ has a vortex-free rendition $\rho_i^S$ in a disk.

                \item If $H_8^{i*}$ is the union of $C_8^{i*}$ and the outer graph of $C_8^{i*}$ in the restriction $\rho_i'$ of $\rho_i$ to $G_i' - (A \cup Z)$, then $(U_i \cup H_8^{i*}, \Omega_i')$ has a rendition $\rho_i^*$ in a disk $\Delta$ with exactly two vortices $v_1,v_2$ such that the trace of each path in $\mathcal{P}$ in $\rho_i^*$ separates $v_1$ and $v_2$ in $\Delta$.
                Furthermore, the $\rho_i'$ and $\rho_i^*$ agree on the rendition of $H_8^{i*}$ and $\rho_i^S$ and $\rho_i'$ agree on the rendition of $U_i$.

                \item For all $j \in [\ell+1] \setminus \{ i \}$, the cycle $C_3^{j*}$ is disjoint from $U_i$.

                \item The transaction $\mathcal{P}$ is taut in $(U_i,\Psi_i)$.
                
                \item For all $j \in [\ell+1] \setminus \{ i \}$, no vertex in $A$ has a neighbour in the union of $C_{|A|+2}^j$ and the outer graph of $C_{|A|+2}^j$, and either
                \begin{enumerate}
                    \item the set $A$ is empty and $V(U_i')$ is disjoint from $Z$, where $(U_i',\Psi_i)$ is the $\mathcal{P}$-strip society in $(G_i',\Omega_i')$, or

                    \item letting $Q^i_L \coloneqq  P^i_{\lfloor \nicefrac{p}{2} \rfloor}$, $Q^i_R \coloneqq  P^i_{\lceil \nicefrac{(p+1)}{2} \rceil}$, and $(U_i', \Psi_i')$ be the $\{ Q^i_L, Q^i_R \}$-society in $(U_i, \Psi_i)$, the society $(G_i'[V(U_i') \cup A \cup Z], \Psi_i')$ has a cross and $A \subseteq V(G')$.
                \end{enumerate}
            \end{enumerate} 
    \end{enumerate}
    In particular, the $K_t$-minor model, the set $A$, and $\mathcal{P}$ can each be found in time $\mathbf{poly}(t+a_1+a_2+p+s'+\ell+p') |E(G)||V(G)|^2$.
\end{lemma}
\begin{proof}
    We give the following estimates for the functions in the statement, including functions used for the purpose of keeping the statement of the functions modestly sized.
    Let $s'' \coloneqq s - (s'(\ell + 1) + 1)$.
    \begin{align*}
        \mathsf{loss}_{\ref{lem:freeflattransactionorprotectedcross}}(t) \coloneqq                      &\ 8t^3 + 7 \\
        \mathsf{nest}_{\ref{lem:freeflattransactionorprotectedcross}}(t,a_1,a_2,s',\ell) \coloneqq      &\ s'\ell + 8t^3 + a_1 + a_2 + 8\\
        \mathsf{apex}^\mathsf{link}_{\ref{lem:freeflattransactionorprotectedcross}}(t,a_1) \coloneqq    &\ 8t^3 + a_1 \\
        \mathsf{apex}^\mathsf{genus}_{\ref{lem:freeflattransactionorprotectedcross}}(t,a_2) \coloneqq   &\ 8t^3 + a_2 \\
        \mathsf{link}^c_{\ref{lem:freeflattransactionorprotectedcross}}(s,q) \coloneqq                  &\ 97t^3(s(q+1)+41)+1 \\
        \mathsf{link}^p_{\ref{lem:freeflattransactionorprotectedcross}}(t,p',\ell,a_1) \coloneqq        &\ 97t^3(\ell(2(8t^3+a_1+1)p'+120)+40)+1 \\
        \mathsf{crook}^1_{\ref{lem:freeflattransactionorprotectedcross}}(t,s,q,p',\ell,a_1) \coloneqq   &\ \mathsf{f}_{\ref{thm:orthogonalexposedsorcery}}((\mathsf{link}^p_{\ref{lem:freeflattransactionorprotectedcross}}(t,a_1,p',\ell) - 1)(\mathsf{link}^c_{\ref{lem:freeflattransactionorprotectedcross}}(s,q) - 1)+1,s'' - 1,\ell) \\
        \mathsf{crook}^2_{\ref{lem:freeflattransactionorprotectedcross}}(t,s,q,p',\ell,a_1) \coloneqq   &\ \mathsf{f}_{\ref{lem:orthogonalexposedapply}}(\mathsf{crook}^1_{\ref{lem:freeflattransactionorprotectedcross}}(t,s,q,p',\ell,a_1),s''-1) \\
        \mathsf{depth}_{\ref{lem:freeflattransactionorprotectedcross}}(t,s,q,p',\ell,a_1) \coloneqq     &\ 12(8 \lfloor \nicefrac{p}{2} \rfloor - 14 +3(( \mathsf{crook}^2_{\ref{lem:freeflattransactionorprotectedcross}}(t,s,q,p',\ell,a_1) + 7)(\lfloor \nicefrac{p}{2} \rfloor - 2) + \\
                                                                                                        &\ \mathsf{crook}^2_{\ref{lem:freeflattransactionorprotectedcross}}(t,s,q,p',\ell,a_1)))+4
    \end{align*}

    Since the stack of functions building up to $\mathsf{depth}_{\ref{lem:freeflattransactionorprotectedcross}}(t,s,q,p',\ell,a_1)$ are quite complex, let us briefly give an estimate of the order of this function.
    It is easy to see that $\mathsf{link}^c_{\ref{lem:freeflattransactionorprotectedcross}}(s,q) \in \mathbf{O}(t^3sq)$ and $\mathsf{link}^p_{\ref{lem:freeflattransactionorprotectedcross}}(t,p',\ell,a_1) \in \mathbf{O}(t^6p'\ell + t^3a_1p'\ell)$.
    As a result we have $\mathsf{crook}^1_{\ref{lem:freeflattransactionorprotectedcross}}(t,s,q,p',\ell,a_1) \in \mathbf{O}( t^9p'qs^2\ell^2 + t^6a_1p'qs^2\ell^2 )$, using the fact that $s'' \leq s$.
    This is then means that $\mathsf{crook}^2_{\ref{lem:freeflattransactionorprotectedcross}}(t,s,q,p',\ell,a_1) \in \mathbf{O}(t^9p'qs^3\ell^2 + t^6a_1p'qs^3\ell^2)$, which yields $\mathsf{depth}_{\ref{lem:freeflattransactionorprotectedcross}}(t,s,q,p',\ell,a_1) \in \mathbf{O}(t^9p'qs^3\ell^2 + t^6a_1p'qs^3\ell^2 + p)$.
    
    Let $i \in [\ell+1]$ be the integer for which the $C_{s - (s'(\ell+1))}^{i*}$-society $(G_i' - Z,\Omega_i')$ of $(G_i,\Omega_i)$ has depth at least $\mathsf{depth}_{\ref{lem:freeflattransactionorprotectedcross}}(t,s,q,p',\ell,a_1)$, and let $s'' \coloneqq s - (s'(\ell + 1) + 1)$.
    Further, for all $i \in [\ell]$, let $\mathfrak{P}_i = \{ \mathcal{P}_1^i, \ldots , \mathcal{P}_{\ell_i}^i \}$ and for all $j \in [\ell_i]$, let $\mathcal{P}_j^i = \{ P_1^{j,i}, \ldots , P_{\lfloor \nicefrac{p}{2} - (s'\ell + 1) \rfloor}^{j,i} \}$.

    For now, let $c \coloneqq \mathsf{crook}^2_{\ref{lem:freeflattransactionorprotectedcross}}(t,s,q,p',\ell,a_1)$ for ease of reading.
    We start by applying \Cref{lem:crookedexistencealgo} to $(G_i' - Z,\Omega_i')$, which result in us finding a crooked transaction $\mathcal{Q}_1$ in $(G_i' - Z,\Omega_i')$, since $(G_i' - Z,\Omega_i')$ has depth at least $8 \lfloor \nicefrac{p}{2} \rfloor - 14 +3((c + 7)(\lfloor \nicefrac{p}{2} \rfloor - 2) + c)$.\footnote{We slightly overestimate what we need here. We could replace the occurrences of $\lfloor \nicefrac{p}{2} \rfloor$ with $\lfloor \nicefrac{p}{2} \rfloor - (s'\ell + 1)$ throughout this paragraph, since the transactions in our padding have order $\lfloor \nicefrac{p}{2} \rfloor - (s'\ell + 1)$. But this difference disappears in our estimates for the order of our functions. So we choose to make our expressions simpler instead.}
    For each $j \in [\ell]$, let $\mathcal{C}_j' = \{ C_1^{j*}, \ldots , C_{s''}^{j*} \}$ be the inner nest of the $s'$-narrowing of $(\mathcal{C}_j, \mathfrak{P}_j)$ (see \Cref{lem:narrowing}).
    We may now use \Cref{lem:crookedwitnesspadded} to find a transaction $\mathcal{Q}_2 \subseteq \mathcal{Q}_1$ of order $(c + 7)(\lfloor \nicefrac{p}{2} \rfloor - 2) + c$ that has a crooked witness $\mathcal{W}$ on $C_2^{i*}$.

    Let $X,Y$ be the two end segments of $\mathcal{Q}_2$ in $(G_i' - Z,\Omega_i')$ and let $\mathcal{R}_Z$ be the $Z$-$V(C_2^{i*})$-linkage of order $|\mathcal{Q}_2|$ for both $Z \in \{ X,Y \}$ and let $Z'$ be the set of endpoints of the paths of $\mathcal{R}_Z$ in $V(C_2^{i*})$.
    Using \Cref{lemma:orthogonal_transaction_padded}, there exists an $X'$-$Y'$-linkage $\mathcal{Q}_3'$ of order $c = \mathsf{crook}^2_{\ref{lem:freeflattransactionorprotectedcross}}(t,s,q,p',\ell,a_1)$ in the union of $C_2^{i*}$ and the inner graph of $C_2^{i*}$ in $\rho_i$ such that $\bigcup \mathcal{Q}_3'$ is disjoint from both $\bigcup \mathcal{W}$ and $\bigcup_{j \in [\ell+1] \setminus \{ i \}} H_j^3$, where for each $j \in [\ell+1]$ the graph $H_j^3$ is the union of $C_3^{j*}$ and the outer graphs of $C_3^{j*}$ in $\rho_j$.

    Now, for each $j \in [\ell+1]$, let $(G_j'',\Omega_j'')$ be the $C_{s''}^{j*}$-society in $\rho_j$ and note that this society retains a cozy nest of order $s'' - 1$ around $c_0$, namely $\mathcal{C}_j' \setminus \{ C_s''^{j*} \}$ according to \Cref{lem:innernestcozyininnersociety}.
    We can construct a transaction $\mathcal{Q}_3$ in $(G_i'',\Omega_i'')$ of order $|\mathcal{Q}_3'|$ from $\mathcal{Q}_3'$ and $\mathcal{R}_X \cup \mathcal{R}_Y$ on $(G_i'',\Omega_i'')$ by following the paths in $\mathcal{R}_X \cup \mathcal{R}_Y$ from their endpoints in $X' \cup Y'$ down to $V(\Omega_i'')$.
    We let $X''$ be the set of endpoints of the paths in $\mathcal{Q}_3$ that is found in paths from $\mathcal{R}_X$ and let $Y''$ be analogously defined via $\mathcal{R}_Y$.
    The transaction $\mathcal{Q}_3$ is disjoint from $\bigcup \mathcal{W}$, which is in turn a crooked witness for it, since $\mathcal{R}_X \cup \mathcal{R}_Y$ are found entirely inside the outer graph of $C_2^{i*}$, aside from the vertices in $X' \cup Y'$.
    This also confirms that $\mathcal{Q}_3$ is disjoint from $\bigcup_{j \in [\ell+1] \setminus \{ i \}} H_j^3$.

    We let $\Omega_i^*$ be defined such that $V(\Omega_i^*) = V(\Omega_i'') \cap V(\bigcup \mathcal{Q}_3)$ and the order itself is derived by restricting $\Omega_i''$ to $V(\Omega_i^*)$, and let $G_i^* \coloneq G_i'' - Z$.
    Note that $V(\Omega_i^*) \subseteq X'' \cup Y''$.
    We have now ensured that all necessary conditions are met to allow us to apply \Cref{lem:orthogonalexposedapply} to $\{ (G_j'' - Z, \Omega_j'') \}_{j \in [\ell+1] \setminus \{ i \}}$ and $(G_i^*,\Omega_i^*)$.
    As a result we find a crooked transaction $\mathcal{Q}_4$ in $(G_i^*,\Omega_i^*)$ of order $\mathsf{crook}^1_{\ref{lem:freeflattransactionorprotectedcross}}(t,s,q,p',\ell,a_1)$ that is orthogonal to $\mathcal{C}_i'' \coloneqq \mathcal{C}_i' \setminus \{ C_1^{i*}, C_2^{i*}, C_{s''}^{i*} \}$.
    In turn, this allows us to apply \Cref{thm:orthogonalexposedsorcery} to find a transaction $\mathcal{Q}_5$ in $(G_i^*,\Omega_i^*)$ of order $(\mathsf{link}^p_{\ref{lem:freeflattransactionorprotectedcross}}(t,a_1,p',\ell) - 1)(\mathsf{link}^c_{\ref{lem:freeflattransactionorprotectedcross}}(s,q) - 1)+1$ with co-conspirators $\mathcal{W}'$ (on level 3), such that $\mathcal{Q}_5$ is orthogonal to $\mathcal{C}_i''$ and disjoint from the union of $C_3^{j*}$ and the outer graph of $C_3^{j*}$ in $\rho_j$ for all $j \in [\ell+1] \setminus \{ i \}$.

    Note that due to the presence of $\mathcal{W}'$, we know that $\mathcal{Q}_5$ is exposed in the restriction $\rho_i^*$ of $\rho_i$ to $G_i^*$.
    Next, we apply \Cref{lem:monotonetransaction} to find a monotone transaction $\mathcal{Q}_6 \subseteq \mathcal{Q}_5$ which is either a crosscap transaction of order $\mathsf{link}^c_{\ref{lem:freeflattransactionorprotectedcross}}(s,q)$ or a planar transaction of order $\mathsf{link}^p_{\ref{lem:freeflattransactionorprotectedcross}}(t,a_1,p',\ell)$.
    We treat these two cases separately.
    
    \textbf{$\mathcal{Q}_6$ is a crosscap transaction.}
    We note that thanks to $\mathcal{Q}_6 \subseteq \mathcal{Q}_5$, the transaction $\mathcal{Q}_6$ is orthogonal to $\mathcal{C}_i''$.
    This allows us to apply \Cref{lemma:findflattransaction} to either find a $K_t$-minor model with the properties we desire, in which case we are done, or find a transaction $\mathcal{Q}_7 \subseteq \mathcal{Q}_6$ -- which is thus still a crosscap transaction -- of order $s(q+1)+1$ and a set $A' \subseteq V(\sigma_{\rho_i}(c_i))$ with $|A'| \leq 8t^3$, such that $\mathcal{Q}_7$ is isolated and flat in $(G^*_i - A',\Omega^*_i)$.
    We can extend $\mathcal{Q}_7$ to a isolated, flat crosscap transaction in $(G_i',\Omega_i')$ via the paths in $\mathcal{R}_X \cup \mathcal{R}_Y$, though this transaction may not be orthogonal to any part of $\mathcal{C}_i$.
    Let $\mathcal{R}^\star$ be the collection of all $V(\Omega_i')$-$V(C_3^{i*})$-paths in $\bigcup \mathcal{Q}_7$ and note that $\mathcal{R}^\star$ has order $2(s(q+1)+1)$.
    We let $X^\star$ be the endpoints of $\mathcal{R}^\star$ on $C_3^{i*}$. 

    We now apply \Cref{lemma:RadialLinkageOnTheSurface} to $\mathcal{R}^\star$ with $k = \lfloor\nicefrac{p}{2}\rfloor-s'\ell$ to find an $X^\star$-$V(\Omega'')$-linkage in $H$,
    such that $\mathcal{R}^\star$ can be combined with $\mathcal{Q}_7$ into a crosscap transaction $\mathcal{Q}_8$ in $(G'',\Omega'')$.
    Let $A_1, \ldots , A_\ell$ be the apices of the transactions in $\mathfrak{P}$.
    According to \Cref{def:transactionmesh}, there exists a vortex-free rendition $\rho^\dagger$ of the society $(H - \bigcup_{i \in [\ell+1]} V(H_i^\ell - V(\Omega_i^\ell)), \Omega)$ in a disk in which $\bigcup_{i \in [\ell+1]} V(\Omega_i^\ell) \subseteq N(\rho^\dagger)$, such that for each $i \in [\ell+1]$, $\rho^\dagger$ is a vortex-free rendition of $(H - \bigcup_{j \in [\ell+1] \setminus \{ i \}} (H_j^\ell - V(\Omega_j^\ell)), \Omega_i)$ in a disk.
    The existence of $\rho^\dagger$, combined with the fact that the strip society of $\mathcal{Q}_7$ in $(G^*_i,\Omega_i')$ has a vortex-free rendition, due to $\mathcal{Q}_7$ being flat, implies that the strip society $(U^\star,\Psi^\star)$ in $(G''' \coloneqq G'' - (\bigcup_{j=1}^\ell A_j \cup Z \cup A'), \Omega'')$ has a vortex-free rendition.
    Thus $\mathcal{Q}_8$ is flat and isolated in $(G''',\Omega'')$.
    We note that $A^\star \coloneqq \bigcup_{j=1}^\ell A_j \cup Z \cup A'$ has order at most $8t^3 + a_1 + a_2$.

    We can now apply \Cref{lem:reconciliation} to $\mathcal{Q}_8$ and the restriction $\rho''$ of $\rho$ to $(G''',\Omega'')$ to find a vortex-free rendition $\rho_8$ of $(G^8,\Omega'')$ into a surface $S$ that is isomorphic to the projective plane with a boundary, where $G^8$ is the union of $U^\star$, $C_{8t^3 + a_1 + a_2 + 7}$, and the outer graph of $C_{8t^3 + a_1 + a_2+7}$ in $\rho''$.
    In particular, the trace $T$ of $C_{8t^3 + a_1 + a_2+7}$ in $\rho''$ splits $S$ into an annulus $A$, whose closure contains $C_{8t^3 + a_1 + a_2+7}$ and the outer graph of $C_{8t^3 + a_1 + a_2+7}$ and another copy $S'$ of the projective plane with a boundary.
    We may now modify the surface in which $\rho''$ is drawn by first restricting it to $A$ and then placing a disk $d$ with the boundary $T$ against $A$.
    We can then draw the contents of $S'$ into $d$ as we wish.
    This allows us to apply \Cref{lemma:orthogonal_transaction} to $\mathcal{Q}_8$ within $(G^8,\Omega'')$ to find a transaction $\mathcal{Q}_9$ in $G^8$ that is orthogonal to $\{ C_{8t^3 + a_1 + a_2 + 8}, \ldots , C_{s-s'-1} \}$.
    Since $\mathcal{Q}_8$ has order $s(q+1)+1$, we can guarantee that $\mathcal{Q}_9$ has order at least $q+2$.
    Note that, with the exception of at most two paths that may be using $C_{8t^3 + a_1 + a_2+7}$ to travel through $d$, all of the paths in $\mathcal{Q}_9$ must travel entirely through $U^\star$ within $d$ and thus, there exists a transaction $\mathcal{Q}_{10} \subseteq \mathcal{Q}_9$, with $|\mathcal{Q}_{10}| \geq |\mathcal{Q}_9| - 2$, such that $\mathcal{Q}_{10}$ is a flat, isolated crosscap transaction in $(G'',\Omega'')$ and $\mathcal{Q}_{10}$ is orthogonal to $\{ C_{8t^3 + a_1 + a_2 + 8}, \ldots , C_{s-s'-1} \}$.
    As $\mathcal{Q}_9$ has order $q+2$ and $\mathcal{Q}_{10}$ accordingly has order at least $q$, this satisfies the second point of our statement, allowing us to move on to the other case.

    \textbf{$\mathcal{Q}_6$ is a planar transaction.}
    In this case, we note that $\mathcal{Q}_6$ has order $97t^3(\ell(2(8t^3+a_1+1)p'+120)+40)+1$ and apply \Cref{lemma:findflattransaction}.
    If this returns a $K_t$-minor model with the properties we desire, we are done.
    Otherwise, we find a transaction $\mathcal{Q}_6' \subseteq \mathcal{Q}_6$ of order $\ell(2(8t^3+a_1+1)p'+120)$ and a set $A \subseteq \sigma_{\rho_i'}(c_i)$ with $|A| \leq 8t^3$, such that $\mathcal{Q}_6'$ is flat and isolated in $(G_i^* - A,\Omega_i^*)$.
    Let $(U,\Psi)$ be the $\mathcal{Q}_6$-strip in $(G_i^* - A, \Omega_i^*)$.

    If we have $V(H_j^3) \cap V(U) \neq \emptyset$ for any $j \in [\ell+1] \setminus \{ i \}$, we must have $H_j^3 \subseteq U$, as $\mathcal{Q}_6' \subseteq \mathcal{Q}_6$ and $\bigcup \mathcal{Q}_6$ is disjoint from all such graphs.
    Thus, since $\mathcal{Q}_6'$ has order $\ell( 2(8t^3+a_1+1)p'+120 )$, there exists a transaction $\mathcal{Q}_7 \subseteq \mathcal{Q}_6'$ of order $2(8t^3+a_1+1)p'+120$ that is planar, isolated, flat, and whose strip society $(U',\Psi')$ in $(G_i^* - A,\Omega_i^*)$ does not intersect any $H_j^3$ with $j \in [\ell+1] \setminus \{ i \}$, meaning $V(H_j^3) \cap V(U') = \emptyset$ for all of these graphs.
    Thus $\mathcal{Q}_7$ is a planar, isolated, flat transaction and whose strip $(H_7,\Psi_7)$ in $(G_i^* - A, \Omega_i^*)$ also does not intersect any $H_j^3$ with $j \in [\ell+1] \setminus \{ i \}$.

    We start with a few observations in relation to the third option of our statement.
    First, we note that $\mathcal{Q}_7$ has co-conspirators, since $\mathcal{Q}_5$ has the co-conspirators $\mathcal{W}'$ in $(G_i^*, \Omega_i^*)$, which we can inherit.
    Furthermore, $C_3^{j*}$ is disjoint from $U'$, since in particular $H_j^3$ is disjoint from $U'$, as argued above.
    Thus we can guarantee the properties outlined in points a), b), and d) for $\mathcal{Q}_7$ at this point, though the transaction we are talking about is currently on the wrong society.
    We will remedy this last point later on whilst preserving the properties we have just touched on.

    Moving further, we assume that $A$ is chosen to be minimal such that $\mathcal{Q}_7$ is flat and isolated in $(G_i^* - A, \Omega_i^*)$, which can be easily enforced in time $\mathbf{O}(a|V(G)|)$ by returning individual vertices to $G_i^* - A$ and checking whether the properties still hold via \Cref{prop:TwoPaths} and \Cref{prop:mengersthm}.
    We want to show that $A$ is disjoint from the union $H_j'$ of $C_{|A|+1}^{j*}$ and the outer graph of $C_{|A|+1}^{j*}$ in $\rho_j$ for all $j \in [\ell+1] \setminus \{ i \}$.
    For this purpose, we may assume that $\ell \geq 1$, as this otherwise trivially holds.

    Suppose that there exists a vertex $u \in A$ such that $u \in V(H_j')$ for some $j \in [\ell+1] \setminus \{ i \}$.
    Then according to the minimality of $A$, the transaction $\mathcal{Q}_7$ is not flat or not isolated in $(G_i^u \coloneqq G_i^* - (A \setminus \{ u \}), \Omega_i^*)$.
    
    First, suppose that $\mathcal{Q}_7$ is not flat and thus there exists a pair of paths $P,Q$ forming a cross in $(G_i^u, \Omega_i^*)$.
    We must have $u \in V(P \cup Q)$, since $\mathcal{Q}_7$ was previously flat.
    In particular, there exists a $C_1^{j*}$-path $P' \subseteq P \cup Q$ in the outer graph of $C_1^{j*}$ in $\rho_j$ that contains $u$.
    Furthermore, for $P'$ to be in $H_j'$ it has to be in some bridge $B$ for the graph $\bigcup \mathcal{Q}_7 \cup V(\Omega_i^*)$ that has an attachment $v$ outside of the boundary paths of $\mathcal{Q}_7$.
    We note that $v$ lies outside of the outer graph of $C_1^{j*}$ in $\rho_j$.
    Since $u$ lies in $H_j'$, there exists some cycle in $C \in \mathcal{C}_j'$ that is disjoint from $A$ with $(V(C) \cap V(P')) \setminus u \neq \emptyset$ and $C \subseteq B$.
    Thus, in $G_i^* - A$, there exists a bridge $B' \subseteq B$ with $v \in V(B')$ and there exists a $C_1^{j*}$-$C_{|A|+1}^{j*}$-path $P'' \subseteq P'$ with $P'' \subseteq B'$.
    This contradicts the fact that we ensured that $H_7$ is disjoint from $H_j^3$ earlier.

    The argument showing that $\mathcal{Q}_7$ must be isolated as well proceeds analogously.
    Hence $A$ must be disjoint from $H_j'$ for all $j \in [\ell+1] \setminus \{ i \}$.
    The fact that $A$ has no neighbours in the union of $C_{|A|+2}^{j*}$ and the outer graph of $C_{|A|+2}^{j*}$ in $\rho_j$ then immediately follows from this.
    This confirms the first part of point f).
    
    Recall that we let $(H_7,\Psi_7)$ be the strip of $\mathcal{Q}_7$ in $(G_i^* - A,\Omega_i^*)$.
    Let $G_8$ be the union of $C_8^{i*}$ and the inner graph of $C_8^{i*}$ in $\rho_i$.
    \Cref{lem:reconciliation} tells us that we can find a rendition $\rho'$ of $(H_8 \coloneqq H_7 \cup G_8, \Omega_i')$ into a disk $\Delta$ such that the restriction of $\rho'$ to $G_8$ is the restriction of $\rho_i$ to $G_8$, the entire graph $H_8$ is grounded in $\rho'$, and $\rho'$ has exactly two vortices $v_1,v_2$.
    Furthermore, we note that for each path $Q \in \mathcal{Q}_7$ there exists a pair of paths $R_X \in \mathcal{R}_X$ and $R_Y \in \mathcal{R}_Y$, such that $Q$ contains the endpoint of $R_X$ in $X'$ and $Q$ contains the endpoint of $R_Y$ in $Y'$.
    Using these, we note that the rendition $\rho'$ even has the property that for each $Q \in \mathcal{Q}_7$ there exists a path in $H_8$ with one endpoint in $X$ and the other in $Y$ whose trace separates $v_1$ and $v_2$ in the disk $\Delta$.
    This construction will ultimately be what tells us that we can guarantee point b) of our statement.

    It is now time to move $\mathcal{Q}_7$ back onto $V(\Omega_i'')$.
    Note that $V(\Omega_i^*) \cap V(\Omega_i'')$ must be empty and each vertex in $V(\Omega_i^*)$ is found inside of some path in $\mathcal{R}_X \cup \mathcal{R}_Y$.
    We now ask for an $X$-$Y$-linkage $\mathcal{Q}_8$ of order $\nicefrac{1}{2}((8t^3+a_1+1)p'+60)$ in $\bigcup (\mathcal{R}_X \cup \mathcal{R}_Y \cup \mathcal{Q}_7)$ via \Cref{prop:mengersthm}.
    Suppose w.l.o.g.\ that this yields a separator $S$ of order at most $\nicefrac{1}{2}((8t^3+a_1+1)p'+60) - 1$ within $\bigcup (\mathcal{R}_X \cup \mathcal{R}_Y \cup \mathcal{Q}_7)$.
    Note that all elements of $\mathcal{Q}_7$ are pairwise disjoint and the same is true for the elements of $\mathcal{R}_X \cup \mathcal{R}_Y$.
    However, the fact that $|S| < \nicefrac{1}{2} ((8t^3+a_1+1)p'+60)$ implies that there exists a vertex $u \in S$ that is contained in three elements of $\mathcal{Q}_8 \cup \mathcal{R}_X \cup \mathcal{R}_Y$.
    As just pointed out, this leads to a contradiction.

    The linkage $\mathcal{Q}_8$ is in particular a transaction in $(G_i' - (A \cup Z), \Omega_i')$ of order $(8t^3+a_1+1)p'+60$, each path of which is both contained and grounded in $H_8$.
    Thus $\mathcal{Q}_8$ is planar, isolated, and flat in $(G_i' - (A \cup Z), \Omega_i')$ as well.
    Furthermore, the trace of each path in $\mathcal{Q}_8$ separates the vortices $v_1,v_2$ of $\rho'$ in $\Delta$.
    We let $\mathcal{Q}_8 = \{ Q_1', \ldots , Q_{(8t^3+a_1+1)p'+60}' \}$ be indexed naturally such that the trace of $Q_1'$ separates the trace of $v_1$ from the traces of all other paths in $\mathcal{Q}_8$.
    Further, let $\mathcal{Q}_L = \{ Q_1', \ldots , Q_{31}' \}$, let $\mathcal{Q}_R = \{ Q_{(8t^3+a_1+1)p'+30}', \ldots , Q_{(8t^3+a_1+1)p'+60}' \}$, and let $\mathcal{Q}_M = \mathcal{Q}_8 \setminus (\mathcal{Q}_L \cup \mathcal{Q}_R)$.

    Recall that $\mathcal{W}'$ are co-conspirators (on level 3) for $\mathcal{Q}_7$ in $(G_i^*,\Omega_i^*)$.
    Let $I,I'$ be the two end segments of $\mathcal{Q}_7$ and let $A,B$ be the two segments of $\Omega_i^*$ with one endpoint in $I$ and the other in $I'$.
    Suppose that $\mathcal{W}'$ is a snitch, turncoat, or pair of deserters, then there exists a path $P \in \mathcal{W}'$ such that it is either the only path in $\mathcal{W}'$ with endpoints in $A$ or the only path in $\mathcal{W}'$ with endpoints in $B$.
    By construction of $\mathcal{Q}_8$ and the fact that $\mathcal{Q}_8$ and $C_{s''}^{i*}$ are both grounded in $\rho'$, we note that the only two paths in $\mathcal{Q}_8$ which can reach these endpoints on $C_{s''}^{i*}$ without intersecting another path from $\mathcal{Q}_8$ are $Q_1'$ and $Q_{(8t^3+a_1+1)p'+60}'$.
    Independent of whether we find a path with these endpoints on $A$ or $B$, we can simply follow along $C_{s''}^{i*}$ until we meet $Q_1'$ and $Q_{(8t^3+a_1+1)p'+60}'$ and then continue traversing this path until we meet $V(\Omega_i')$.
    Thus, if $\mathcal{W}'$ is a snitch or a pair of deserters, we can actually find co-conspirators (on level 3) for $\mathcal{Q}_8$.
    Should $\mathcal{W}'$ be a turncoat, we can at least route one half of the co-conspirators back onto $\Omega_i'$.

    Therefore we may suppose that $\mathcal{W}'$ contains a pair of paths $P,P'$ which forms a cross in $(G_i^*, \Omega_i^*)$ and w.l.o.g.\ we may assume that $P,P'$ has all of their endpoints in $A$.
    We further assume that we have chosen $A$ such that the endpoints of $P,P'$ are separated within $\Delta$ by the trace of $Q_1'$ in $\rho'$ from the traces of all other paths in $\rho'$ and in particular, from the vortex $v_2$ of $\rho'$.
    Let $T_{32}$ be the trace of $Q_{32}'$ in $\rho'$, let $\Delta_{31}$ be the closure of the disk within $\Delta - T-{32}$ that contains $v_1$ and let $\Delta_1$ be the $v_1$-disk of $C_{s''}^{j*}$ in $\rho'$.\footnote{This is a mild abuse of notation, as $\rho'$ has two vortices. We can easily enforce this to be interpreted correctly by deleting $v_2$ and $\sigma_{\rho'}(v_2)$ from the rendition and the graph before trying to determine the $v_1$-disk of $C_{s''}^{j*}$.}
    Let $\Delta_L$ be the disk in the intersection of $\Delta_{32}$ and $\Delta_1$ that contains $v_1$.
    We note that $\Delta_L$ is in fact $\rho'$-aligned due to its boundary being defined via the union of two grounded objects within $\rho'$.
    
    We let $(H_L,\Omega_L)$ be the $\Delta_L$-society in $\rho'$ and let $\rho_L$ be the restriction of $\rho'$ to $H_L$.
    For each $i \in [21]$, we let $C_i^*$ be the cycle whose trace defines the $v_1$-cycle of $C_{s'' - (22 - i)}^{i*} \cup Q_{10+i}$.
    The set $\{ C_1^*, \ldots , C_{21}^* \}$ is a set of pairwise disjoint cycles forming a nest in $(H_L,\Omega_L)$, which is a simple consequence of the properties the elements of $\mathcal{C}_i'$ and $\mathcal{Q}_8$ have towards each other in $\rho'$.
    We can use this set to apply \Cref{lemma:makenestcozy} to find a cozy nest $\mathcal{C}' = \{ C_1', \ldots , C_{21}' \}$ in $\rho_L$ such that $\bigcup_{j=1}^{21} C_i'$ is contained in the union of $C_1^*$ and the outer graph of $C_1^*$ in $\rho_L$.
    From this last property we can easily derive the existence of a radial linkage $\mathcal{R}$ for $\mathcal{C}'$ of order $10$ found entirely in $\bigcup_{j=1}^{10} Q_j'$.
    \Cref{lem:radialtoorthogonal} allows us to turn $\mathcal{R}$ into a radial linkage $\mathcal{R}'$ for $\mathcal{C}'$ that is orthogonal to $\mathcal{C}'$ and end-identical to $\mathcal{R}'$.
    We let $\mathcal{R}'' \subseteq \mathcal{R}'$ be the four paths that use the endpoints on $C_1'$ that are found in $Q_7 \cup Q_8 \cup Q_9 \cup Q_{10}$.

    Note that $P,P'$ is disjoint from the graph $H_7$ of the strip of $\mathcal{Q}_7$, since it is disjoint from $\mathcal{Q}_7$ and $\mathcal{Q}_7$ is isolated.
    Thus, we know that $(H_L' \coloneqq H_L \cup ((P \cup P') - (V(H_8) \setminus V(H_L))), \Omega_L)$ contains a cross $Q,Q'$ with $Q \subseteq P$ and $Q' \subseteq P'$.
    Let $\rho_L'$ be a rendition of $(H_L',\Omega_L)$ derived from $\rho_L$ by drawing all vertices in $V(H_L') \setminus V(H_L)$ and edges in $E(H_L') \setminus E(H_L)$ into $v_1$.
    Further, let $(G_L,\Omega_L^*)$ be the $C_{17}'$-society in $\rho_L'$ and let $x_1,x_2,x_3,x_4 \in V(C_{17}')$ be the four vertices that are first seen when traversing each of the paths in $\mathcal{R}''$ starting from their endpoints in $V(\Omega_L)$.
    This allows us to apply \Cref{lemma:reroute_crosses} using $(G_L,\Omega_L^*)$, $\rho_L'$, $\{ C_1', \ldots , C_{16}' \}$, $\mathcal{R}'$, $\mathcal{R}''$, and the existence of $Q,Q'$ to find a cross consisting of two paths $P'',Q''$ in $(G_L,\Omega_L^*)$ on the vertices $x_1,x_2,x_3,x_4$.

    Let $\mathcal{R}'''$ be the $\{ x_j ~\!\colon\!~ j \in [4] \}$-$V(C_{21}')$-linkage of order 4 in $\bigcup \mathcal{R}''$.
    Additionally, let $\mathcal{R}^*$ be $V(C_{17}')$-$V(\Omega_i')$-linkage of order 4 in $Q_7 \cup Q_8 \cup Q_9 \cup Q_{10}$.
    We claim that $\bigcup ( \mathcal{R}''' \cup \mathcal{R}^* \cup \{ C_{18}', C_{19}', C_{20}', C_{21}' \} )$ contains a $\{ x_j ~\!\colon\!~ j \in [4] \}$-$V(\Omega_i')$-linkage of order 4, which would allow us to find a cross on $(G_i', \Omega_i')$ containing $P'',Q''$.
    To find this, we apply \Cref{prop:mengersthm}.
    Suppose we find a $\{ x_j ~\!\colon\!~ j \in [4] \}$-$V(\Omega_i')$-separator $S$ of order 3 in this graph instead of the linkage we wish for.
    Then note that size the elements of the three sets we used to construct this graph are pairwise disjoint, one of each of these elements is not intersected by $S$.
    Since all paths in $\mathcal{R}'''$ intersect all cycles in $\{ C_{18}', C_{19}', C_{20}', C_{21}' \}$ and these in turn intersect all paths in $\mathcal{R}^*$, this leads to a contradiction.

    Analogous arguments of course also recover a cross on $B$ towards $(G_i', \Omega_i')$ via the paths in $\mathcal{Q}_R$.
    Thus we know that we can confirm point a) of the third part of our statement and can therefore guarantee the existence of co-conspirators for $\mathcal{Q}_8$.

    Let $\mathcal{Q}_M = \{ Q_1'', \ldots , Q_{(8t^3+a_1+1)p'-2}'' \} \subseteq \mathcal{Q}_8$ be indexed naturally and let $\mathcal{Q}_8' = \{ Q_{\lfloor \nicefrac{p'}{2} \rfloor}'', \ldots , Q_{(8t^3+a_1+1)p'-2 - (\lceil \nicefrac{(p'+1)}{2} \rceil - 1)}'' \}$, which has order $(8t^3+a_1)p'$.
    Further, let $(U'', \Psi'')$ be the $\mathcal{Q}_8'$-strip in $(G_i' - A, \Omega_i')$.
    
    If the society $(U^* \coloneqq G_i'[V(U'') \cup A \cup Z], \Psi'')$ contains a cross, the transaction $\mathcal{Q}_M \setminus \{ Q_{\lfloor \nicefrac{p'}{2} \rfloor + 1}'', \ldots , Q_{(8t^3+a_1+1)p'-2 - \lceil \nicefrac{(p'+1)}{2} \rceil}'' \}$ has order at least $p'$ and satisfies point f) in the third part of our statement.
    
    Otherwise, \Cref{prop:TwoPaths} tells us that $(U^*, \Psi'')$ has a vortex-free rendition in a disk.
    Note that $A \cup Z$ is disjoint from $\bigcup \mathcal{Q}_8'$.
    Thus, since $|A \cup Z| \leq 8t^3+a_1$ and $\mathcal{Q}_8'$ has order $(8t^3+a_1)p'$, there exists a transaction $\mathcal{Q}_9 \subseteq \mathcal{Q}_8'$ of order $p'$, whose strip society $(H_9,\Psi_9)$ in $(G_i', \Omega_i')$ is disjoint from $A \cup Z$, meaning $(A \cup Z) \cap V(H_9) = \emptyset$.
    This in particular means that $A \cup Z$ does not have any neighbours in $H_9$ that are not vertices of the boundary paths of $\mathcal{Q}_9$.
    Thus $(H_9,\Psi_9)$ has a vortex-free rendition, which means that we can choose $A' = \emptyset$ and now also satisfy point f) of our statement, with $\mathcal{Q}_9$ and $A'$.

    We note that we can easily extend $\rho'$ to also cover the part of the graph found in $G_i' - V(G_i^* \cup A \cup Z)$, due to our definition of $(G_i^*,\Omega_i^*)$ via the $C_{s''}^{i*}$-society in $\rho_i$.
    Thus we can guarantee that point b) holds independent of whether we satisfy the subpoints of f) via $\mathcal{Q}_M \setminus \{ Q_{\lfloor \nicefrac{p'}{2} \rfloor + 1}'', \ldots , Q_{(8t^3+a_1+1)p'-2 - \lceil \nicefrac{(p'+1)}{2} \rceil}'' \}$ or $\mathcal{Q}_9$.
    All that remains is to apply \Cref{corollary:maketransactiontaut} to find a taut transaction, which allows us to satisfy e) without violating the other properties we fought for.
    This completes our proof.
\end{proof}

We now show how to use the transaction we find in the third option of \Cref{lem:freeflattransactionorprotectedcross} to extend our existing transaction mesh by one more transaction.
However, simply building this new transaction mesh is not enough, we will prove here that we can sensibly extend the $\mathfrak{M}$-tree of the previous transaction mesh.
This should finally shed some light into the origin and purpose of the various options found in \Cref{def:Mtree} and in particular, give some intuition as to what the signature of a transaction mesh is tracking.
We will go into much more detail regarding $\mathfrak{M}$-trees and their signatures in the next subsection.

\begin{lemma}\label{lem:buildtransactionmesh}
    There exists polynomial functions $\mathsf{loss}_{\ref{lem:buildtransactionmesh}} \colon \mathbb{N}^2 \rightarrow \mathbb{N}$, $\mathsf{apex}^\mathsf{link}_{\ref{lem:buildtransactionmesh}} \colon \mathbb{N}^3 \rightarrow \mathbb{N}$, $\mathsf{nest}_{\ref{lem:buildtransactionmesh}} \colon \mathbb{N}^6 \rightarrow \mathbb{N}$, and $\mathsf{depth}_{\ref{lem:buildtransactionmesh}} \colon \mathbb{N}^8 \rightarrow \mathbb{N}$, such that the following holds.
    
    Let $t,s,s',a,a',p,q,k,z',\ell$ be non-negative integers $\mathsf{nest}_{\ref{lem:buildtransactionmesh}}(t,z',a',s',\ell,k) \leq s$, $8k+2a'+8+s(q-1) \leq \nicefrac{p}{4}$, $z' \leq \ell$, and $\mathsf{apex}^\mathsf{link}_{\ref{lem:freeflattransactionorprotectedcross}}(t,z',k) \leq a$.
    
    Let $(G,\Omega)$ be a society with a cylindrical rendition $\rho$ in a disk $\Delta$, a cozy nest $\mathcal{C} = \{ C_1, \ldots , C_s \}$ of order $s$ in $\rho$ around the vortex $c_0$, and a padded $(s,s',a,p,\ell)$-transaction mesh $\mathfrak{M}$ with the $k$-signature $\Bar{x} = (x_1,x_2,x_3,x_4)$, such that $x_2 \leq z'$ and the union of the apices for the transactions in $\mathfrak{P}$ have order at most $a'$.
    Further, let $\{(G_i,\Omega_i)\}_{i \in [\ell+1]}$ be the districts of $\mathfrak{M}$ with cylindrical renditions $\rho_i$ in disks and an $(s - (s'\ell + 1),\lfloor \nicefrac{p}{2} \rfloor - (s'\ell + 1),\ell_i)$-padding $(\mathcal{C}_i = \{ C_1^i, \ldots , C_{s - (s'\ell + 1)}^i \}, \mathfrak{P}_i)$ with the inner nest $\mathcal{C}_i' = \{ C_1^{i*}, \ldots , C_{s - (s'\ell + 1)}^{i*} \}$.
    Let $i \in [\ell+1]$ be chosen such that $(G_i'-Z,\Omega_i')$ has depth at least $\mathsf{depth}_{\ref{lem:buildtransactionmesh}}(t,s,q,p,a',\ell,z',k)$, where $(G_i',\Omega_i')$ is the $C_{s - (s'(\ell+1))}^{i*}$-society of $(G_i,\Omega_i)$.
    Finally, let $G' \coloneqq \bigcup_{i \in [\ell+1]} G_i$ and let $(G'', \Omega'')$ be the $C_{s - s'\ell}$-society in $\rho$.

    There there exists a set $A \subseteq V(G'')$, and
    \begin{enumerate}
        \item a $K_t$-minor model in $G$ controlled by a mesh whose horizontal paths are subpaths of distinct cycles from $\mathcal{C}$,
    
        \item a flat crosscap transaction $\mathcal{P}$ of order $q$ in $(G''-A,\Omega'')$, with $|A| \leq \mathsf{apex}^\mathsf{genus}_{\ref{lem:freeflattransactionorprotectedcross}}(t,a')$ and a nest $\mathcal{C}'$ in $\rho$ of order $s - \mathsf{loss}_{\ref{lem:freeflattransactionorprotectedcross}}(t)$ around $c_0$ to which $\mathcal{P}$ is orthogonal, or

        \item we have $A \subseteq G'$ with $|A| \leq \mathsf{apex}^\mathsf{link}_{\ref{lem:freeflattransactionorprotectedcross}}(t,z',k)$, there exists a transaction $\mathcal{P}$ of order $p - 2\mathsf{loss}_{\ref{lem:buildtransactionmesh}}(t,k)$ in $(G_i', \Omega_i')$ with co-conspirators $\mathcal{W}$ (on level 3), and a padded $(s - (\mathsf{loss}_{\ref{lem:buildtransactionmesh}}(t,k) + s'),s',a,p - 2\mathsf{loss}_{\ref{lem:buildtransactionmesh}}(t,k),\ell+1)$-transaction mesh $\mathfrak{M}' = (H',\Omega, \mathcal{C}', \mathfrak{P} \cup \{ \mathcal{P}_{\ell+1} \coloneqq \mathcal{P} \})$ with the $k$-signature $\Bar{y}$, such that $\Bar{y} \geq \Bar{x}$.
    \end{enumerate}
    
    In particular, the $K_t$-minor model, the set $A$, the transaction $\mathcal{P}$, the co-conspirators $\mathcal{W}$, and the transaction mesh $\mathfrak{M}'$ can each be found in time $\mathbf{poly}( t+p+s'+\ell+p' ) |E(G)||V(G)|^2$.
\end{lemma}
\begin{proof}
    We give the following estimates for the functions in the statement.
    Let $p' \coloneqq p - 2\mathsf{loss}_{\ref{lem:buildtransactionmesh}}(t,k)$.
    \begin{align*}
        \mathsf{loss}_{\ref{lem:buildtransactionmesh}}(t,k) \coloneqq                       &\ \max(8t^3 + 10, 2k+12) \\
        \mathsf{apex}^\mathsf{link}_{\ref{lem:buildtransactionmesh}}(t,z',k) \coloneqq      &\ \mathsf{apex}^\mathsf{link}_{\ref{lem:freeflattransactionorprotectedcross}}(t,z'(3k+2)) \\
        \mathsf{nest}_{\ref{lem:buildtransactionmesh}}(t,z',a',s',\ell,k) \coloneqq           &\ \mathsf{loss}_{\ref{lem:buildtransactionmesh}}(t,k) + \mathsf{nest}_{\ref{lem:buildtransactionmesh}}(t,z'(3k+2),a',s',\ell)\\
        \mathsf{depth}_{\ref{lem:buildtransactionmesh}}(t,s,q,p,a',\ell,z',k) \coloneqq       &\ \mathsf{depth}_{\ref{lem:freeflattransactionorprotectedcross}}(t,s,q,p',\ell,z'(3k+2))
    \end{align*}

    We again provide some estimates for the order of some of these functions, using the estimates we derived in the proof of \Cref{lem:freeflattransactionorprotectedcross}.
    Two simpler estimates are for $\mathsf{nest}_{\ref{lem:buildtransactionmesh}}(t,z',a',s',\ell,k) \in \mathbf{O}( \max(t^3,k) + s'\ell + z'k + a' )$ and for $\mathsf{apex}^\mathsf{link}_{\ref{lem:buildtransactionmesh}}(t,z',k) \in \mathbf{O}( t^3 + z'k )$.
    The most complex estimate yields $\mathsf{depth}_{\ref{lem:buildtransactionmesh}}(t,s,q,p',\ell,z',k) \in \mathbf{O}( t^9pqs^2\ell^2 + t^6z'kpqs^2\ell^2 )$, using the fact that $p' \leq p$.

    Let $\mathfrak{T} = ( T, \mathfrak{M}, \phi, \mathcal{A}, \mathcal{Z} = \{ Z_1, \ldots , Z_\ell \}, \mathcal{F}, \mathfrak{L} )$ be the $\mathfrak{M}$-tree of precision $k$ that witnesses the signature $\Bar{x}$ of $\mathfrak{M}$ and let $\mathfrak{P} = \{ \mathcal{P}_1, \ldots , \mathcal{P}_\ell \}$.
    Furthermore, let $\mathsf{z} \in [\ell]$ be the maximal integer such that $A_i = \emptyset$ for all $i \in [(\ell + 1) - \mathsf{z},\ell]$.
    We let $Z \coloneqq \bigcup_{j = (\ell + 1) - \mathsf{z}}^\ell B_\mathfrak{T}(Z_i)$.
    Note that according to \textbf{\textsf{T5}} we have $Z \subseteq V(G') \setminus V(H)$ and furthermore, we have $|Z| \leq (3k+2)z'$, due to $x_3 \leq z'$.
    
    Let $(G_j',\Omega_j')$ be the $C_{s - (s'(\ell+1))}^j$-society of the district $(G_j,\Omega_j)$ with $j \in [\ell+1]$ which has depth at least $\mathsf{depth}_{\ref{lem:buildtransactionmesh}}(t,s,q,p',\ell,z',k)$.
    Note that we use the index $j$ here in reference to the index $j_i$ that first appears in \textbf{\textsf{M1}} of \Cref{def:transactionmesh}.
    In the context of this proof our $j$ will correspond to $j_{\ell+1}$ in the construction of our new transaction mesh.
    We now apply \Cref{lem:freeflattransactionorprotectedcross} at $(G_j,\Omega_j)$ using our $Z$.
    
    Of course the first two options of that statement correspond directly to the first two options of the present statement and thus, if we reach either of them, we are done.
    Otherwise, we find the planar, isolated, flat transaction $\mathcal{P}$ and the set $A \subseteq V(G')$ with $|A| \leq \mathsf{apex}^\mathsf{link}_{\ref{lem:buildtransactionmesh}}(t,z',k)$ which together have the properties described in the third option of the statement of \Cref{lem:freeflattransactionorprotectedcross}.
    Let $\mathcal{P} = \{ P_1, \ldots , P_{p'} \}$ and let $\mathcal{C}= \{ C_1, \ldots , C_s \}$.
    
    Our goal will be to add $\mathcal{P}_{\ell+1} \coloneq \mathcal{P}$ to $\mathfrak{M}$, make some modifications to $\mathfrak{M}$ and then verify the conditions set out in \Cref{def:transactionmesh}, to show that this yields a new transaction mesh $\mathfrak{M}'$.
    We will also have to confirm that $\mathfrak{M}'$ remains padded.
    Following this, we need to update the signature and analyse the $\mathfrak{M}'$-tree to confirm the third point of our statement.
    
    Our analysis will be split into two cases, depending on whether or not $A$ is empty.
    Nonetheless, the construction of $\mathfrak{M}'$ and the verification of its properties is largely the same.
    Thus we begin with discussing the properties of $\mathfrak{M}'$ we can guarantee in both cases.
    
    We let $a'' \in \{ 10, |A| + 10 \}$ and for now let $H'$ be a graph whose definition we specify later on.
    Both of these objects will be made concrete within the two cases we later consider.
    Let $\mathfrak{P}'$ be derived from $\mathfrak{P}$ by iteratively removing the boundary paths of each $\mathcal{P} \in \mathfrak{P}$ for $a''$ iterations.
    Our first goal is to show that in most cases $\mathfrak{M}' = (H',\Omega, \mathcal{C} \setminus \{ C_1, \ldots , C_{a''} \}, \mathfrak{P}' \cup \{ \mathcal{P}_{\ell+1} \coloneqq \mathcal{P} \})$ is the $(s-(a''+s'),s',a,p-2a'',\ell+1)$-transaction mesh we are after.

    To verify this, we need to check the properties demanded of $\mathfrak{M}'$ in \Cref{def:transactionmesh} and \Cref{def:paddedtransactionmeshes}.
    Let $\{(H_i,\Omega_i)\}_{i \in [\ell+1]}$ be the areas of $\mathfrak{M}$ with the $(s - (s'\ell + 1),\lfloor \nicefrac{p}{2} - (s'\ell + 1) \rfloor,\ell_i)$-padding $(\mathcal{C}_i = \{ C_1^i, \ldots , C_{s - (s'\ell + 1)}^i \}, \mathfrak{P}_i = \{ \mathcal{P}_1^i, \ldots , \mathcal{P}_{\ell_i}^i \})$ and for all $j \in [\ell_i]$ let $\mathcal{P}_j^i = \{ P_1^{j,i}, \ldots , P_{\lfloor \nicefrac{p}{2} - (s'\ell + 1) \rfloor}^{j,i} \}$.
    For all $i \in [\ell + 1] \setminus \{ j \}$, we can construct the new area $(H_i^{\ell+1}, \Omega^{\ell+1})$ and set $\mathcal{C}_i^{\ell+1} = \{ C_{a''+1}^i, \ldots , C_{s - (s'(\ell+1) +1)}^i \}$, largely following the instructions laid out in \textbf{\textsf{M1}} of \Cref{def:transactionmesh}.
    This allows us to let the new padding be defined by taking the $s'$-narrowing of $(\mathcal{C}_i^{\ell+1}, \{ \{ P_{a''+1}^{1,i}, \ldots , P_{s - (s'(\ell+1) +1)}^{1,i} \}, \ldots , \{ P_{a''+1}^{\ell_i,i}, \ldots , P_{s - (s'(\ell+1) +1)}^{\ell_i,i} \} \} )$ (see \Cref{lem:narrowing}).
    Thus for $i \in [\ell + 1] \setminus \{ j \}$ we have confirmed that the area $(H_i^{\ell+1}, \Omega^{\ell+1})$ is padded in accordance with \Cref{def:paddedtransactionmeshes}.

    Let $(H_j',\Omega_j')$ be the $C_{s - (s'(\ell+1))}^j$-society of $(H_j,\Omega_j)$ in the restriction of $\rho_j$ to $H_j$.
    The transaction $\mathcal{P}$ is planar and exposed in $(H_j' \cup \bigcup \mathcal{P}, \Omega_j')$, as well as being disjoint from $\mathcal{C}_i^{\ell+1}$ for all $i \in [\ell + 1] \setminus \{ j \}$, according to point iii)d) of \Cref{lem:freeflattransactionorprotectedcross}.
    This confirms \textbf{\textsf{M2}} of \Cref{def:transactionmesh}.

    For \textbf{\textsf{M3}}, we can place $A$ into the role of $A_{\ell+1}$, note that the graph $G'$ defined there corresponds to the component of $G'$ defined in our statement that contains $H_j'$.
    We know that the strip-society $(U_j,\Psi_j)$ of $\mathcal{P}$ in $(G',\Omega_j')$ has a vortex-free rendition $\rho_j^S$ thanks to point iii)b) of \Cref{lem:freeflattransactionorprotectedcross}.
    The rendition $\tau_j$ of $(H_j' \cup U_j,\Omega_j')$ in a disk $\Delta_j$ demanded in the remainder of \textbf{\textsf{M3}} is confirmed the exist by point iii)c) of \Cref{lem:freeflattransactionorprotectedcross}.
    We restate the properties of $\tau_j$ in the context of the objects we defined here for the convenience of the reader.
    The rendition $\tau_j$ agrees with $\rho_j^S$ on $U_j$ and agrees with $\rho_j$ on the union of $C_8^{i*}$ and the outer graph of $C_8^{i*}$ in $\rho_j$.
    Furthermore, the trace in $\tau_j$ of each path in $\mathcal{P}$ separates the unique pair of vortices $v_1,v_2$ of $\tau_j$ in $\Delta_j$.

    The items \textbf{\textsf{M4}} and \textbf{\textsf{M5}} of \Cref{def:transactionmesh} follow directly from the properties of $\mathcal{P}$ and $A$ laid out in \Cref{lem:freeflattransactionorprotectedcross}.
    Thus we are left with having to confirm \textbf{\textsf{M6}} of \Cref{def:transactionmesh} and having to construct a nice padding for the resulting pair of areas.
    
    We can define the new areas $(H_j^{\ell+1}, \Omega_j^{\ell+1}), (H_{\ell+2}^{\ell+1}, \Omega_{\ell+2}^{\ell+1})$ in a straightforward way by using the properties of $\tau_j$ and following the instructions laid out in i), ii), and iii) of \textbf{\textsf{M6}}, where we let $Q_L \coloneqq P_{\lfloor \nicefrac{p'}{2} \rfloor}$, and $Q_R \coloneqq  P_{\lceil \nicefrac{(p'+1)}{2} \rceil}$.
    Let $\rho_L$ be the restriction of $\tau_j$ to $H_j^{\ell+1}$ and let $\rho_R$ be the restriction of $H_{\ell+2}^{\ell+1}$.
    We assume that $v_1,v_2$ are named such that $\rho_L$ is a cylindrical rendition of $(H_j^{\ell+1}, \Omega_j^{\ell+1})$ in a disk around $v_1$ and $\rho_R$ is a cylindrical rendition of $(H_{\ell+2}^{\ell+1}, \Omega_{\ell+2}^{\ell+1})$ in a disk around $v_2$.
    To construct the new nest for $(H_j^{\ell+1}, \Omega_j^{\ell+1})$, for each $i \in [a''+1, s - (s'(\ell+1) +1)]$, we let $C_i^{j,\ell+1}$ be cycle whose $v_1$-disk in $\rho_L$ is inclusion-wise maximal with respect to the set of cycles in $C_i^j \cup P_{(\lfloor \nicefrac{p'}{2} \rfloor + i) - s - (s'(\ell+1) +2)}$.
    We let $\mathcal{C}_j^{\ell+1*} = \{ C_i^{j,\ell+1} ~\!\colon\!~ i \in [a''+1, s - (s'(\ell+1) +1)] \}$ and define $\mathcal{C}_{\ell+2}^{\ell+1*}$ analogous for $(H_{\ell+2}^{\ell+1}, \Omega_{\ell+2}^{\ell+1})$ and $\rho_R$.
    This allows us to simply apply \Cref{lemma:makenestcozy} to make the nests $\mathcal{C}_j^{\ell+1*}, \mathcal{C}_{\ell+2}^{\ell+1*}$ cozy, yielding two nests $\mathcal{C}_j^{\ell+1}, \mathcal{C}_{\ell+2}^{\ell+1}$ and confirming the remaining parts of point iii) in \textbf{\textsf{M6}}.
    
    To construct the appropriate padding for these two areas we first take the $s'$-narrowing $(\mathcal{C}_j^{\ell+1\star}, \{ \{ P_{a''+1}^{1,j}, \ldots , P_{s - (s'(\ell+1) +1)}^{1,j} \}, \ldots , \{ P_{a''+1}^{\ell_j,j}, \ldots , P_{s - (s'(\ell+1) +1)}^{\ell_j,j} \} \} )$ of the $(s - (s'\ell + 1),\lfloor \nicefrac{p}{2} - (s'\ell + 1) \rfloor,\ell_i)$-padding of $(H_j^{\ell+1}, \Omega_j^{\ell+1})$.
    The $(s - (s'(\ell+1) + 1)-a'',\lfloor \nicefrac{p}{2} - (s'(\ell+1) + 1-a'') \rfloor,\ell_j^*)$-padding $(\mathcal{C}_j^{\ell+1}, \mathfrak{P}_j^{\ell+1})$ for $(H_j^{\ell+1}, \Omega_j^{\ell+1})$ is then formed by letting, for each maximal segment $S$ of $\Omega_j^{\ell+1}$ with $S \subseteq V(P_{s - s'(\ell+1)}^{i,j}$ for some $i \in [\ell_j]$, $\mathcal{P}_i^{j,\ell+1}$ be a transaction on $V(\Omega_j^{\ell+1}$ in $\bigcup_{h=a''+1}^{s - s'(\ell+1) - 1} P_h^{i,j}$ and adding these transactions to $\mathfrak{P}_j^{\ell+1}$.
    These exists thanks to the fact that the paths involved in the construction are grounded in $\rho_L$.
    We note in particular that using the same arguments as in the proof of \Cref{lem:narrowing}, one can see that there exist at most $\ell_j$ such segments $S$ and thus $\ell_j^* \leq \ell_j + 1$.
    Analogous constructions yield an $(s - (s'(\ell+1) + 1)-a'',\lfloor \nicefrac{p}{2} - (s'(\ell+1) + 1-a'') \rfloor,\ell_{\ell+2}^*)$-padding for $(H_{\ell+2}^{\ell+1}, \Omega_{\ell+2}^{\ell+1})$.
    This confirms that our new transaction mesh is padded.

    For the parts that need divergent analysis, we start with the easier of the two cases.

    \textbf{Non-empty apex set:}
    Suppose that $A \neq \emptyset$.
    We let $a'' = |A| + 10$ in this case.
    Let $(U_j,\Psi_j)$ be the $\mathcal{P}$-society in $(G_j' - (A \cup Z), \Omega_j')$.
    Then according to \Cref{lem:freeflattransactionorprotectedcross}iii)f)2), if we let $(U_j', \Psi_j')$ be the $\{ Q_L, Q_R \}$-society in $(U_j, \Psi_j)$, the society $(G_j'[V(U_j') \cup A \cup Z], \Psi_j')$ has a cross and $A \subseteq V(G')$.
    
    We let $H'$ be the union of $U_j$ and $H - A$.
    The society $(H',\Omega)$ then has a vortex-free rendition since $(H,\Omega)$ has a vortex-free rendition according to \Cref{def:transactionmesh} and $(H_j' \cup U_j,\Omega_j')$ has a vortex-free rendition, as we just noted.
    Furthermore, using $\tau_j$ and the special rendition $\rho^\dagger$ from \Cref{def:transactionmesh} for $\mathfrak{M}$, it is easy to derive the vortex-free rendition $\rho^{\dagger\dagger}$ of the society $(H' - \bigcup_{i \in [\ell+2]} (H_i^{\ell+1} - V(\Omega_i^{\ell+1})), \Omega)$ in a disk in which $\bigcup_{i \in [\ell+1]} V(\Omega_i^{\ell+1}) \subseteq N(\rho^{\dagger\dagger})$, such that for each $i \in [\ell+1]$, $\rho^{\dagger\dagger}$ is a vortex-free rendition of $(H' - \bigcup_{j \in [\ell+1] \setminus \{ i \}} (H_j^{\ell+1} - V(\Omega_j^{\ell+1})), \Omega_i^{\ell+1})$ in a disk.
    We then let $A_{\ell+1} \coloneqq A \cup Z$, which are the apices associated with $\mathcal{P}$.
    As a result we have confirmed that $\mathfrak{M}'$ is a padded $(s-(a''+s'),s',a,p-2a'',\ell+1)$-transaction mesh, as desired.

    We now turn our attention towards the the $\mathfrak{M}$-tree $\mathfrak{T} = ( T, \mathfrak{M}, \phi, \mathcal{A}, \mathcal{Z}, \mathcal{F}, \mathfrak{L} )$ of precision $k$.
    Here our construction for the new $\mathfrak{M}'$-tree $\mathfrak{T}'$ is fairly straightforward.
    In particular, we let $\mathfrak{T}' = ( T', \mathfrak{M}', \phi', \mathcal{A} \cup \{ A_{\ell+1} \}, \mathcal{Z} \cup \{ Z_{\ell+1} \coloneqq \emptyset \}, \mathcal{F}', \mathfrak{L} \cup \{ \mathcal{L}_{\ell+1} \coloneqq \emptyset \} )$, where $T'$ and $\phi'$ are derived from $T$ and $\phi$ in the way described by \textbf{\textsf{T1}} of \Cref{def:Mtree} and $\mathcal{F}'$ is similarly derived directly via the instructions laid out in \textbf{\textsf{T8}} of \Cref{def:Mtree}.
    Note that $\mathfrak{T}'$ also has precision $k$.
    From the definition of the $k$-signature $\Bar{y}$ of $\mathfrak{M}'$ given in \Cref{def:Mtree}, we can now easily derive the fact that $\Bar{y} \geq \Bar{x}$, since $A_{\ell+1} \neq \emptyset$.
    This concludes the analysis of this case.

    \textbf{Empty apex set:}
    Suppose instead that $A = \emptyset$.
    We let $a'' = 10$.
    Note that \textbf{\textsf{M8}} of \Cref{def:transactionmesh} holds trivially, since $A = \emptyset$.
    Thus, using $A_{\ell+1} = \emptyset$, and $H'$ as the union of $H$ and $U_j$, where $(U_j,\Psi_j)$ is the $\mathcal{P}$-society in $(G_j', \Omega_j')$, we confirm that $\mathfrak{M}'$ is a $(s-(a''+s'),s',a,p-2a'',\ell+1)$-transaction mesh using analogous arguments to the previous case.
    We will have to replace $\mathfrak{M}'$ in the last subcase we consider to actually build a proper $\mathfrak{M}'$-tree.
    But for the other subcases $\mathfrak{M}'$ will suffice.

    We now again have to consider how the new $\mathfrak{M}'$-tree $\mathfrak{T}'$ can be constructed.
    For this purpose, we can once more extend $T$ and $\phi$ of the $\mathfrak{M}$-tree $\mathfrak{T} = ( T, \mathfrak{M}, \phi, \mathcal{A}, \mathcal{Z}, \mathcal{F}, \mathfrak{L} )$ of precision $k$ in a straightforward fashion to $T'$ and $\phi'$ in the way described by \textbf{\textsf{T1}} of \Cref{def:Mtree}.
    The remainder of our arguments now focuses on how the two new districts $(G_j^{\ell+1}, \Omega_j^{\ell+1}), (G_{\ell+2}^{\ell+1}, \Omega_{\ell+2}^{\ell+1})$, respectively corresponding to $(H_j^{\ell+1}, \Omega_j^{\ell+1}), (H_{\ell+2}^{\ell+1}, \Omega_{\ell+2}^{\ell+1})$, in $\mathfrak{M}'$ are connected to each other in $G^* \coloneqq \bigcup G_i^* - Z$, where $\{(G_i^*,\Omega_i^*)\}_{i \in [\ell+2]}$ are the districts of $\mathfrak{M}'$ with $G_j^* = G_j^{\ell+1}$ and $G_{\ell+2}^* = G_{\ell+2}^{\ell+1}$.

    First, suppose that there exists a $V(\Omega_j^*)$-$V(\Omega_{\ell+2}^*)$-linkage $\mathcal{L}_{\ell+1}$ of order $k+1$ in $G^*$.
    In this case we can let $\mathfrak{T}' = ( T', \mathfrak{M}', \phi', \mathcal{A} \cup \{ A_{\ell+1} \}, \mathcal{Z} \cup \{ Z_{\ell+1} \coloneqq \emptyset \}, \mathcal{F}', \mathfrak{L} \cup \{ \mathcal{L}_{\ell+1} \} )$, where $\mathcal{F}'$ is constructed via the instructions laid out in \textbf{\textsf{T8}} of \Cref{def:Mtree}.
    We note that the order of $\mathcal{F}'$ cannot have changed when compared to $\mathcal{F}$, as we only moved $s'$ cycles inwards, leaving each area connected, since $s \geq \mathsf{nest}_{\ref{lem:buildtransactionmesh}}(t,z',a'',s',\ell)$ in each area and we did not delete any vertices otherwise.
    Thus, if $\Bar{y}$ is the $k$-signature of $\mathfrak{T}$, then we have $\Bar{y} > \Bar{x}$, as the non-empty linkage $\mathcal{L}_{\ell+1}$ was added to $\mathfrak{L}$.

    Thus, we may instead suppose that there exists a $V(\Omega_j^*)$-$V(\Omega_{\ell+2}^*)$-separator $Z'$ in $G^*$.
    We note that $V(G^*) \cap V(H) = \emptyset$ and we can choose $Z'$ to be minimal.
    The case in which $Z' = \emptyset$ has to be considered separately.
    Here, we let $\mathfrak{T}' = ( T', \mathfrak{M}', \phi', \mathcal{A} \cup \{ A_{\ell+1} \}, \mathcal{Z} \cup \{ Z_{\ell+1} \coloneqq \emptyset \}, \mathcal{F}', \mathfrak{L} \cup \{ \mathcal{L}_{\ell+1} \coloneqq \emptyset \} )$, where $\mathcal{F}'$ is constructed via the instructions laid out in \textbf{\textsf{T8}} of \Cref{def:Mtree}.
    The number of sets in $\mathcal{F}'$ cannot decrease since the strip of $\mathcal{P}$ may split a family into two, but it does so within a single district and thus produces a new district for the family.
    In fact, if the number of families stays the same, then, due to the presence of the co-conspirators for $\mathcal{P}$, this results in us finding more districts that host crosses.
    Thus, if $\Bar{y}$ is the $k$-signature of $\mathfrak{T}$, we have we have $\Bar{y} > \Bar{x}$.
    
    We may therefore move on to the case in which $Z' \neq \emptyset$.
    This means we mainly need to confirm that \textbf{\textsf{T5}} of \Cref{def:Mtree} holds.
    Consider a partitioning $\mathcal{F}^*$ of $[\ell+2]$ into as few sets as possible such that for all $F \in \mathcal{F}^*$ with $|F| \geq 2$ there exists a $V(\Omega_h^{\ell+1})$-$V(\Omega_{h'}^{\ell+1})$-path in $G^* - Z'$ for all distinct $h,h' \in F$.
    We then let $F_j,F_{\ell+2} \in \mathcal{F}^*$ be the sets with $j \in F_j$ and $\ell+2 \in F_{\ell+2}$, which must be distinct as $Z'$ is a $V(\Omega_h^{\ell+1})$-$V(\Omega_{h'}^{\ell+1})$-separator in $G^*$.
    Let $J^*$ be the union of $\bigcup_{i \in F_j \cup F_{\ell+2}} V(\Omega_i^{\ell+1})$ and all $\bigcup_{i \in F_j \cup F_{\ell+2}} V(\Omega_i^{\ell+1})$-bridges in $G^* - Z'$.
    Further, let $w \in Z'$ be arbitrary and let $J_w \coloneqq G^*[V(J^*) \cup \{ w \}]$.
    
    Note that due to the minimality of $Z'$, there exists a $V(\Omega_j^*)$-$V(\Omega_{\ell+2}^*)$-path in $J_w$.
    In particular, every $V(\Omega_j^*)$-$V(\Omega_{\ell+2}^*)$-path in $J_w$ must use $w$.
    Accordingly, there exists a separation $(A_j,A_{\ell+2})$ in $J_w$ with $V(\Omega_j^*) \subseteq A_j$, $V(\Omega_{\ell+2}^*) \subseteq A_{\ell+2}$, and $A_j \cap A_{\ell+2} = \{ w \}$.
    We let $W_i = \{ V(\Omega_h^{\ell+1} ~\!\colon\!~ h \in F_i \}$ for both $i \in \{ j,\ell+2 \}$ and note that $W_i \subseteq V(J_w[A_i])$.
    Since $\mathfrak{M}'$ is a $(s-a'',s',a,p-2a'',\ell+1)$-transaction mesh, there exist radial linkages $\mathcal{R}_i$ of order $\nicefrac{p}{4} \geq 2k+2$ for $\mathcal{C}_i^{\ell+1}$ for all $i \in [\ell+2]$, according to \Cref{obs:radial_linkages_for_districts}.
    
    For both $i \in \{ j, \ell+2 \}$, we apply \Cref{prop:findblob} in $J_w[A_i]$ to try to separate $\{ w \}$ and $W_i$, with $\{ w \}$ in the role of $X$ in the statement.
    This finds a separation $(S_i,T_i)$ such that $W_i \subseteq T_i$ and $\{ w \} \subseteq S_i$.
    Let $S'_i$ be $S_i \cap T_i$.
    We claim that $S'$ is disjoint from the union of $C^h_{2k+a''+3}$ and the outer graph of $C^h_{2k+a''+3}$ in $\rho_h'$ for all $h \in [\ell+2]$.
    Suppose this is not the case and thus, due to $\rho_h'$ being a cylindrical rendition, there exists some $v \in V(C_{2k+a''+3}^h) \cap S'$.
    However, this also implies that for each $h' \in [a''+1, \ldots , 2k+a''+3]$ the set $V(C_{h'}^h) \cap S'$ is non-empty, since we could otherwise get to $W_i$ via $\mathcal{R}_h$ which clearly contradicts the fact that $|S'| \leq k$.
    This in particular implies that $S' \not\subseteq W_i$. 

    Thus \Cref{prop:findblob} tells us that $|S_i \cap T_i| = k$ and there exists an $(S_i \cap T_i)$-$W_i$-linkage $\mathcal{L}_i'$ of order $k+1$ in $J_w[T_i]$.
    Thus, $(S,T)$ with $S = S_j \cup S_{\ell+2}$ and $T = T_j \cup T_{\ell+2}$, is a separation with $W_i \subseteq T$ for both $i \in \{ j , \ell +2 \}$ and $|S \cap T| \leq 2k+2$, such that $\mathcal{L}_j' \cup \mathcal{L}_{\ell+2}'$ is an $S$-$T$-linkage of order $2k+2$ in $J_w[T]$.
    We let $Z_{\ell+1}' = Z' \cup (S \cap T)$ and note that $|Z'| \leq 3k+2$.
    In particular, our earlier argument for $S'_i$ also tells us that both $S$ and $Z$ must be disjoint from the union of $C^i_{2k+a''+3}$ and the outer graph of $C^i_{2k+a''+3}$ in $\rho_i'$ for all $i \in [\ell+2]$.

    Thus, we can now define a new transaction mesh $\mathfrak{M}^* = (H^*,\Omega, \mathcal{C} \setminus \{ C_1, \ldots , C_{2k+a''+2} \}, \mathfrak{P}' \cup \{ \mathcal{P}_{\ell+1} \coloneqq \mathcal{P} \})$, with $H^*$ being $H - Z_{\ell+1}'$, which inherits all necessary properties from $\mathfrak{M}^*$ to ensure that it is a $(s-(2k+12+s'),s',a,p-2a'',\ell+1)$-transaction mesh.
    Finally, we can let $\mathfrak{T}' = ( T', \mathfrak{M}^*, \phi', \mathcal{A} \cup \{ A_{\ell+1} \}, \mathcal{Z} \cup \{ Z_{\ell+1}' \}, \mathcal{F}^*, \mathfrak{L} \cup \{ \mathcal{L}_{\ell+1} \coloneqq \emptyset \} )$, where $\mathcal{F}^*$ is a partitioning of $[\ell+1]$ into as few sets as possible such that for all $F \in \mathcal{F}$ with $|F| \geq 2$ there exists a $V(\Omega_i^*)$-$V(\Omega_j^*)$-path in $\bigcup_{i=1}^{\ell+1} G_i^* - (Z \cup Z_{\ell+1}')$ for all distinct $i,j \in F$, fitting the conditions laid out in \textbf{\textsf{T8}} of \Cref{def:Mtree}.
    Since we added a non-empty set to $\mathcal{Z}$, we conclude that $\Bar{x}$ is smaller than the $k$-signature $\Bar{y}$ of $\mathfrak{M}^*$, finishing our proof.
\end{proof}

\subsection{Analysing signatures}\label{sec:signatures}
One of the main problems we need to solve in relation to transaction meshes is that we cannot keep building them for too many iterations.
As a transaction mesh grows, we are slowly gathering apices that may destroy parts of the infrastructure of the transaction mesh and there is even an in-built loss of infrastructure that keeps gnawing away at our nest.
The goal of this subsection is therefore to show that for any given integer $t$ there is a concrete number of iterations polynomial in $t$ after which we must find a $K_t$-minor.
This allows us to later set up initial nest to be big enough to endure that many iterations.

We will be heavily relying on \Cref{thm:ApathsArgument} within this section and in particular, we will rely on the bounds for the sets presented in the theorem's second option.
Thus we simply codify these as a value depending on $t$.

\begin{definition}\label{def:Apathsfunction}
    Let $t$ be a positive integer.
    We let
    \[ \mathsf{d}_t \coloneq 72t^3 + \nicefrac{1}{2}(t-3)(t-4) - 144 . \]
\end{definition}

Note that $\mathsf{d}_t$ is exactly the maximum number of elements contained in the union of the objects $S$ and $I$ defined in the second option of \Cref{thm:ApathsArgument}.

Our first -- and by far the easiest lemma in this subsection -- demonstrates nicely how exactly we plan to use \Cref{thm:ApathsArgument}.
\begin{lemma}\label{lem:toomanyfamilies}
    Let $t,s,s',a,p,\ell$ be non-negative integers with $s \geq (2t+24)(8t^3+20)+s'\ell$ and $p \geq (64t^3-128)(2t+24)$.
    Let $(G,\Omega)$ be a society with an $(s,s',a,p,\ell)$-transaction mesh $\mathfrak{M} = (H, \Omega, \mathcal{C}, \mathfrak{P})$ and let $(a',z,f,h)$ be the $\mathsf{d}_t$-signature of $\mathfrak{M}$.

    If $f > \mathsf{d}_t$ the graph $G$ contains a $K_t$-minor model controlled by a mesh whose horizontal paths are subpaths of distinct cycles from $\mathcal{C}$.
    
    Furthermore, there exists an algorithm that finds this $K_t$-minor model in time $\mathbf{O}(t^4|E(G)|)$.
\end{lemma}
\begin{proof}
    Let $\mathfrak{T} = ( T, \mathfrak{M}, \phi, \mathcal{A}, \mathcal{Z}, \mathcal{F}, \mathfrak{L} )$ be the $\mathfrak{M}$-tree of precision $\mathsf{d}_t$ that witnesses the signature of $\mathfrak{M}$.
    Applying \Cref{thm:ApathsArgument} either yields a $K_t$-minor model with our desired properties, or two sets $S \subseteq V(G')$ and $I \subseteq [\ell+1]$ with $|S| \leq 8t^3-16$ and $|I| \leq 64t^3+\nicefrac{1}{2}(t-3)(t-4)-128$ with the properties laid out in the statement of \Cref{thm:ApathsArgument}.
    Note that the set $S$ can only separate at most $8t^3-16$ families in $\mathcal{F}$ into societies with vortex-free renditions and only $64t^3+\nicefrac{1}{2}(t-3)(t-4)-128$ of the remaining families can be covered by having a district included in $I$.
    As $f > \mathsf{d}_t$, this leads to a contradiction to the existence of these sets.
\end{proof}

We next want to show that the parameter $h$ in the signature of a transaction mesh cannot grow indefinitely.
For this purpose we show an intermediate lemma that argues that two sets $S$ and $I$ defined as in the second option of \Cref{thm:ApathsArgument} must be at least as large as $h$.

\begin{lemma}\label{lem:theDgrows}
    Let $t,s,s',a,p,\ell$ be non-negative integers.
    Let $(G,\Omega)$ be a society with an $(s,s',a,p,\ell)$-transaction mesh $\mathfrak{M} = (H, \Omega, \mathcal{C}, \mathfrak{P})$, let $(G_1,\Omega_1), \ldots , (G_{\ell+1}, \Omega_{\ell+1})$ be the districts of $\mathfrak{M}$, and let $G' \coloneqq  \bigcup_{i=1}^{\ell+1} G_i$.
    Further, let $(a',z,f,h)$ be the $\mathsf{d}_t$-signature of $\mathfrak{M}$.

    Then for any pair of sets $S \subseteq V(G')$ and $I \subseteq [\ell+1]$ such that each $\{ V(\Omega_i) ~\!\colon\!~ i \in [\ell+1] \}$-path in $G' - S$ is a $\{ V(\Omega_i) ~\!\colon\!~ i \in I \}$-path and each district $(G_j,\Omega_j)$ with $j \in \{ V(\Omega_i) ~\!\colon\!~ i \in [\ell+1] \setminus I \}$ has a vortex-free rendition in a disk, we have $|S|+|I| \geq h$.
\end{lemma}
\begin{proof}
    Let $\mathfrak{T} = ( T, \mathfrak{M}, \phi, \mathcal{A}, \mathcal{Z}, \mathcal{F}, \mathfrak{L} = \{ \mathcal{L}_1, \ldots , \mathcal{L}_\ell \} )$ be the $\mathfrak{M}$-tree of precision $\mathsf{d}_t$ that witnesses the signature of $\mathfrak{M}$ and let $\mathfrak{P} = \{ \mathcal{P}_1, \ldots , \mathcal{P}_\ell \}$.

    We will prove this statement by induction over $h$, noting that it is trivially true for $h = 0$.
    Let $\mathfrak{M}' = (H, \Omega, \mathcal{C}, \mathfrak{P} \setminus \{ \mathcal{P}_\ell \})$ be the $(s,s',a,p,\ell - 1)$-transaction mesh in $(G,\Omega)$, which can be derived from $\mathfrak{M}$ since $\ell \geq h \geq 1$, such that $(G_1', \Omega_1'), \ldots , (G_\ell', \Omega_\ell')$ are the districts of $\mathfrak{M}'$.
    For each $i \in [\ell]$, we let $\mathcal{C} = \{ C_1^{i,\ell}, \ldots , C_{s - (s'i +1)}^{i,\ell} \}$ be the district nest of $(G_i', \Omega_i')$.
    Since $h \geq 1$, we know that the $\mathsf{d}_t$-signature of $\mathfrak{M}'$ is $(a',z,f,h')$ with $h \geq h' \geq h-1$.
    If $h' = h$, then our statement holds by induction.
    Thus we may suppose that $h' = h-1$.
    
    Suppose that our lemma is false and somehow there exists a pair of sets $S \subseteq V(G')$ and $I \subseteq [\ell+1]$ as in the statement with $|S|+|I| < h$.
    As $h \geq 1$, we know that $\mathcal{L}_\ell$ is a linkage of order $\mathsf{d}_t+1$ that joins $V(\Omega_{i+1})$ to $V(\Omega_{j_i})$ for some $j_i \in [\ell]$.
    In particular, since $\mathsf{d}_t+1 > \mathsf{d}_t$, there exists a path in $\mathcal{L}_\ell$ that is not hit by any vertex in $S$ and thus, we have $i+1, j_i \in I$.

    By definition of $\mathfrak{M}$ and $\mathfrak{M}'$, the linkage $\mathcal{P}_\ell$ is found in the district $(G_{j_i}', \Omega_{j_i}')$ of $\mathfrak{M}'$.
    Let $J = I \setminus \{ i+1 \}$ and $G'' = \bigcup_{i=1}^{\ell} G_i'$.

    \begin{claim}\label{claim:stillAset}
        Each $\{ V(\Omega_i') ~\!\colon\!~ i \in [\ell] \}$-path in $G'' - S$ is a $\{ V(\Omega_i') ~\!\colon\!~ i \in J \}$-path and each district $(G_j',\Omega_j')$ with $j \in \{ V(\Omega_i') ~\!\colon\!~ i \in [\ell] \setminus J \}$ has a vortex-free rendition in a disk.
    \end{claim}
    \emph{Proof of \Cref{claim:stillAset}:}
    First, let us consider an arbitrary $j \in \{ V(\Omega_i') ~\!\colon\!~ i \in [\ell] \setminus J \}$ and note that $j \in \{ V(\Omega_i) ~\!\colon\!~ i \in [\ell+1] \setminus I \}$.
    Thus there already exists a vortex-free rendition of the $C_{s - s'(i + 1)}^{i,\ell}$-society of $(G_i',\Omega_i')$, which can easily be extended to also encompass the outer graph of $C_{s - s'(i + 1)}^{i,\ell}$.

    Thus we may suppose that there somehow exists a $\{ V(\Omega_i') ~\!\colon\!~ i \in [\ell] \}$-path in $G'' - S$ that has at least one endpoint outside of $\bigcup_{i \in J} V(\Omega_i')$.
    Suppose this endpoint is found in $j' \in [\ell] \setminus J$.
    Since $j_i \in J$ and the areas of $\mathfrak{M}$ and $\mathfrak{M}'$ have vortex-free renditions, this path contains a subpath that is a $\{ V(\Omega_i) ~\!\colon\!~ i \in [\ell] \}$-path and in particular a $\bigcup_{i \in I} V(\Omega_i)$-$V(\Omega_{j'})$-path.
    This contradicts each $\{ V(\Omega_i) ~\!\colon\!~ i \in [\ell] \}$-path in $G' - S$ being a $\{ V(\Omega_i) ~\!\colon\!~ i \in I \}$-path, as $j' \not\in I$.
    \hfill$\blacksquare$

    We note that $|J| = |I| -1$ and thus we have $|S|+|J| < h-1$, since $|S|+|I| < h$.
    However, since by \Cref{claim:stillAset}, the sets $S$ and $J$ perform as required in our statement, we should in fact have $|S|+|J| \geq h-1$ according to our induction hypothesis.
    This is absurd and completes our proof.
\end{proof}

Using this lemma it becomes easy to conclude that $h$ cannot grow indefinitely without producing a large $K_t$-minor model for us.

\begin{corollary}\label{cor:theDcantbetoolarge}
    Let $t,s,s',a,p,\ell$ be non-negative integers with $s \geq (2t+24)(8t^3+20)+s'\ell$ and $p \geq (64t^3-128)(2t+24)$.
    Let $(G,\Omega)$ be a society with an $(s,s',a,p,\ell)$-transaction mesh $\mathfrak{M} = (H, \Omega, \mathcal{C}, \mathfrak{P})$ and let $(a',z,f,h)$ be the $\mathsf{d}_t$-signature of $\mathfrak{M}$.
    
    If $h > \mathsf{d}_t$ the graph $G$ contains a $K_t$-minor model controlled by a mesh whose horizontal paths are subpaths of distinct cycles from $\mathcal{C}$.
    
    Furthermore, there exists an algorithm that finds this $K_t$-minor model in time $\mathbf{O}(t^4|E(G)|)$.
\end{corollary}
\begin{proof}
    Let $(G_1,\Omega_1), \ldots , (G_{\ell+1}, \Omega_{\ell+1})$ be the districts of $\mathfrak{M}$ and let $G' \coloneqq  \bigcup_{i=1}^{\ell+1} G_i$.
    According to \Cref{lem:theDgrows}, if $h \geq \mathsf{d}_t$, we know in particular that for any pair of sets $S \subseteq V(G')$ and $I \subseteq [\ell+1]$ such that each $\{ V(\Omega_i) ~\!\colon\!~ i \in [\ell+1] \}$-path in $G' - S$ is a $\{ V(\Omega_i) ~\!\colon\!~ i \in I \}$-path and each district $(G_j,\Omega_j)$ with $j \in \{ V(\Omega_i) ~\!\colon\!~ i \in [\ell+1] \setminus I \}$ has a vortex-free rendition in a disk, we have $|S|+|I| \geq h$.
    However, as $h > \mathsf{d}_t$, this implies that we can find a $K_t$-minor model via \Cref{thm:ApathsArgument} as desired.
\end{proof}

Next, we take on the parameter $z$ of the signature of a transaction mesh.
To see why this parameter should also be limited in its growth, we point to the fact that, similar to $h$, $z$ guarantees the existence of a lot of connectivity between the new districts.
As the transaction mesh grows and whilst we are not encountering a cross that lets $a'$ grow, the parts of $G$ that provide this connectivity are disjoint from each other by construction.
As a consequence, two sets $S$ and $I$ as in the second option of \Cref{thm:ApathsArgument} can again not cover all of the connecting infrastructure between the districts of our transaction mesh.
Whilst this is simple to claim, actually proving this is another matter.
This is without a doubt the most complex part of this subsection.

\begin{lemma}\label{lem:toomanyblobs}
    Let $t,s,s',a,p,\ell$ be non-negative integers with $s \geq (2t+24)(8t^3+20)+s'\ell$ and $p \geq (64t^3-128)(2t+24)$.
    Let $(G,\Omega)$ be a society with an $(s,s',a,p,\ell)$-transaction mesh $\mathfrak{M} = (H, \Omega, \mathcal{C}, \mathfrak{P})$ and let $(a',z,f,h)$ be the $\mathsf{d}_t$-signature of $\mathfrak{M}$.
    
    If $z \geq \mathsf{d}_t$ the graph $G$ contains a $K_t$-minor model controlled by a mesh whose horizontal paths are subpaths of distinct cycles from $\mathcal{C}$.
    Furthermore, there exists an algorithm that finds this $K_t$-minor model in time $\mathbf{O}(t^4|E(G)|)$.
\end{lemma}
\begin{proof}
    Let $\mathfrak{T} = ( T, \mathfrak{M}, \phi, \mathcal{A}, \mathcal{Z} = \{ Z_1, \ldots , Z_\ell \}, \mathcal{F}, \mathfrak{L} )$ be the $\mathfrak{M}$-tree of precision $\mathsf{d}_t$ that witnesses the signature of $\mathfrak{M}$ and let $\mathfrak{P} = \{ \mathcal{P}_1, \ldots , \mathcal{P}_\ell \}$.
    Furthermore, let $\mathsf{z} \in [\ell]$ be the maximal integer such that $A_i = \emptyset$ for all $i \in [(\ell +1 ) - \mathsf{z},\ell]$ and note that since $z \geq \mathsf{d}_t$, we have $|[(\ell +1) - \mathsf{z}, \ell]| \geq \mathsf{d}_t$.

    By \Cref{def:Mtree}, we know that for any two distinct $i,j \in [(\ell +1) - \mathsf{z}, \ell]$ the sets $Z_i$ and $Z_j$ are mutually disjoint and in particular, all $Z_i$ are disjoint from $V(H)$.
    Again following \Cref{def:Mtree}, each $Z_i$ can be divided into two sets $Z_i'$ and $X_i$ such that $Z_i' \cap X_i = \{ z_i \}$, we have $|Z_i'| \leq \mathsf{d}_t$, and $X_i$ is one side of a separation $(X_i,Y_i)$ as described in \textbf{\textsf{T5}} of \Cref{def:Mtree}, with $X_i$ taking the role of $S$ in \textbf{\textsf{T5}}\,i) and $X^1_i, X^2_i$ being the partitioning of $X_i \cap Y_i$ into the two sets of size $\mathsf{d}_t+1$.

    We first prove the following claim that shows that these $X_i$ still connect different districts, with many disjoint paths emanating from $X^1_i$ and $X^2_i$ into disjoint collections of districts.
    Let $(G_1,\Omega_1), \ldots , (G_{\ell+1},\Omega_{\ell+1})$ be the districts of $\mathfrak{M}$ and let $G' := \bigcup_{j \in [\ell+1]} G_j - \bigcup_{j \in [\ell+1]} A_j$.
    We note that, due to our choice of $\mathsf{z}$, we have $G' = \bigcup_{j \in [\ell+1]} G_j - \bigcup_{j \in [\ell - \mathsf{z}]} A_j$.

    \begin{claim}\label{claim:blobpreservation}
        For each $i \in [(\ell +1) - \mathsf{z}, \ell]$ with $Z_i \neq \emptyset$, there exist two disjoint, non-empty sets of integers $I_1^i,I_2^i \subseteq [\ell+1]$ such that for both $j \in [2]$, there exists an $X^j_i$-$( \bigcup_{j' \in I_j^i}V(\Omega_{j'}) )$-linkage of order $\mathsf{d}_t+1$ in $G' - \bigcup_{i' \in [(\ell +1) - \mathsf{z} , i-1] } Z_{i'}$.
    \end{claim}
    \emph{Proof of \Cref{claim:blobpreservation}:}
    Let $\mathfrak{M}' = (H,\Omega,\mathcal{C},\mathfrak{P} \setminus \{ \mathcal{P}_j \colon j \in [i+1, \ell] \})$ be the $(s,s',a,p,i)$-transaction mesh derived from $\mathfrak{M}$ with the districts $(G_1^i,\Omega_1^i), \ldots , (G_{i+1}^i,\Omega_{i+1}^i)$.
    Furthermore, let $(H, \Omega, \mathcal{C}, \mathfrak{P} \setminus \{ \mathcal{P}_j \colon j \in [i, \ell] \})$ be the $(s,s',a,p,i-1)$-transaction mesh derived from $\mathfrak{M}$ with the districts $(G_1^{i-1},\Omega_1^{i-1}), \ldots , (G_i^{i-1},\Omega_i^{i-1})$ and let $j_i \in [i]$ be the integer such that $\mathcal{P}_i$ is found in the $C_{s - s'i}^{j_1,i-1}$-society of $(G_{j_i}^{i-1}, \Omega_{j_i}^{i-1})$, where $\{ C_1^{j_i,i-1}, \ldots , C_{s - (s'(i-1) -1)}^{j_i,i-1} \}$ is the district nest of $(G_{j_i}^{i-1}, \Omega_{j_i}^{i-1})$.

    By our choice of $Z_i$, according to \textbf{\textsf{T5}}\,i) of \Cref{def:Mtree}, there exists an $X^1_i$-$\bigcup_{h=1}^{i+1} V(\Omega_h^i)$-linkage $\mathcal{L}^1$ of order $\mathsf{d}_t+1$ in $G'' = \bigcup_{j \in [i+1]} G_j^i - ( \bigcup_{j \in [\ell - \mathsf{z}]} A_j \cup \bigcup_{j \in [(\ell +1) - \mathsf{z}, i]} Z_i )$ and an $X^2_i$-$\bigcup_{h=1}^{i+1} V(\Omega_h^i)$-linkage $\mathcal{L}^2$ of order $\mathsf{d}_t+1$ in $G''$, such that $\mathcal{L}^1 \cup \mathcal{L}^2$ is itself a linkage that is disjoint from $X_i$, except for one endpoint of each path.
    In fact, the existence of the separation in \textbf{\textsf{T5}}\,i) implies that we can partition $[i+1]$ into two sets $W_{j_i}$ and $W_{i+1}$, such that $h \in W_h$ for both $h \in \{ j_i,i+1 \}$ and $\mathcal{L}^1$ is an $X_1^i$-$\bigcup_{h \in W_{j_i}} V(\Omega_h^i)$-linkage and $\mathcal{L}^2$ is an $X_2^i$-$\bigcup_{h \in W_{i+1}} V(\Omega_h^i)$-linkage.
    Note that $X_i$ is entirely disjoint from $H$ as well.
    
    Thus to find the desired paths, we can follow each of the paths in $\mathcal{L}^1 \cup \mathcal{L}^2$ starting from $X_i$ until they meet a vertex of $\bigcup_{j' \in [\ell+1]} V(\Omega_{j'})$.
    This is guaranteed to happen, since these paths must intersect $H$, as we have $\bigcup_{h=1}^{i+1} V(\Omega_h^i) \subseteq V(H)$ and $H$ has a vortex-free rendition.
    We collect the resulting paths originating from $\mathcal{L}^1$ in $\mathcal{R}^1$ and those originating from $\mathcal{L}^2$ in $\mathcal{R}^2$ and let $I_j^i \coloneqq \{ j \in [\ell+1] \colon R \in \mathcal{R}^i \text{ and } R \text{ has an endpoint in } X_j . \}$ for both $j \in [2]$.

    Suppose there exists some $j \in I_1^i \cap I_2^i$.
    Then let $R \in \mathcal{R}^1 \cup \mathcal{R}^2$ be one of the paths that has an endpoint $v$ in $V(\Omega_j)$.
    In particular, we must have $v \in V(G_j \cap H)$ and thus $R$ must intersect all paths of the district-nest of $(G_j, \Omega_j)$.
    Thus there exists an $X_1^i$-$X_2^i$-path in $G''$ that avoids $Z_i \setminus (X_1^i \cup X_2^i)$, contradicting the fact that $V(\Omega_{j_i}^i)$ and $V(\Omega_{i+1}^i)$ lie in different connected components of $G''$.
    \hfill$\blacksquare$

    This pair of sets $I_1^i,I_2^i$ describes districts which were all interconnected in $G' - \bigcup_{i' \in [(\ell +1) - \mathsf{z} , i-1] } Z_{i'}$, but after the deletion of $Z_i$ the districts corresponding to $I_1^i$ and those corresponding to $I_2^i$ have been disconnected from each other.

    We now apply \Cref{thm:ApathsArgument} to $G'$.
    If this yields a $K_t$-minor, we are done.
    Thus we may suppose that this instead results in us finding sets $S \subseteq V(G')$ and $J \subseteq [\ell+1]$ with $|S| \leq 8t^3-16$ and $|J| \leq 64t^3+\nicefrac{1}{2}(t-3)(t-4)-128$ such that each $\{ V(\Omega_i) ~\!\colon\!~ i \in [\ell+1] \}$-path in $G' - S$ is a $\{ V(\Omega_i) ~\!\colon\!~ i \in J \}$-path and each district $(G_j,\Omega_j)$ with $j \in [\ell+1] \setminus J$ has a vortex-free rendition in a disk.
    Note that by definition $|S| + |J| \leq \mathsf{d}_t \leq |I|$.

    We observe that according to \textbf{\textsf{T5}}\,i) of \Cref{def:Mtree}, $Z_i$ and $Z_j$ are disjoint for any two distinct $i,j \in I$.
    Thus, if we let $I' \subseteq I$ be the set of integers $i \in I$ such that $Z_i \cap S = \emptyset$, we have $|J| \leq |I'|$.

    With the help of \Cref{claim:blobpreservation} and the fact that $X_i$ induces a connected subgraph of $G$, we observe that for all $i \in I'$ there exists a $\{ V(\Omega_j) ~\!\colon\!~ j \in I_1^i \}$-$\{ V(\Omega_j) ~\!\colon\!~ j \in I_2^i \}$-path in $G' - S$, since $|S| < \mathsf{d}_t + 1$.
    Thus for each $i \in I'$, the set $J$ must include two distinct integers $i_1^i \in I_1^i$ and $i_2^i \in I_2^i$.

    Let $I = \{ i_1, \ldots , i_z \}$ such that for any two $h,j \in [z]$ we have $i_h < i_j$ if and only if $h < j$.
    We inductively define a forest $T_i$, for $i \in [0,z]$, that describes the structure of the sets in $\mathcal{F}$.
    For this purpose, let $\mathcal{F}^i$ be a partitioning of $[\ell+1]$ into as few sets as possible such that for each $F \in \mathcal{F}^i$ with $|F| \geq 2$ there exists a $V(\Omega_i)$-$V(\Omega_j)$-path in $G'_i \coloneqq G' - \bigcup_{i' \in [(\ell +1) - \mathsf{z} , i-1] } Z_{i'}$ for all distinct $i,j \in F$.
    Thus we have $G'_0 = G'$.

    We initialise our construction by letting $T_0$ be the graph in which the sets in $\mathcal{F}^0$ correspond to the vertices $V(T_0) = \{ v^0_F ~\!\colon\!~ F \in \mathcal{F}^0 \}$ and there are no edges in $T_0$.
    For each $i \in [z]$, we let $T_i$ be defined by adding the vertices $\{ v^i_F ~\!\colon\!~ F \in \mathcal{F}^i \}$ to $T_{i-1}$ and adding an edge between $v^i_F$ and $v^{i-1}_{F'}$ if $F \subseteq F'$.
    We claim that this results in each vertex that is added to construct $T_i$ having degree 1.
    This is a simple consequence of the fact that $G'_i = G'_i - Z_{i_{i-1}}$, as both $\mathcal{F}^i$ and $\mathcal{F}^{i-1}$ reflect the structure of the connected components containing the districts of these graphs.
    In particular, $\mathcal{F}^i$ is a refinement of $\mathcal{F}^{i-1}$.
    We let $W \coloneq T_z$ and we further consider all trees in $W$ to be rooted at their unique vertex in $\{ v^0_F ~\!\colon\!~ F \in \mathcal{F}^0 \}$.

    For each $i_j \in I$, let $F_j \in \mathcal{F}^{j-1}$ be the set with $I_1^j \cup I_2^j \subseteq F_j$, which must exist since $Z_j$ is non-empty.
    Note that $v^{j-1}_{F_j}$ has at least two children as $Z_j$ separates two sets $V(\Omega_i)$ and $V(\Omega_h)$ with $i \in I_2^j$ and $h \in I_2^j$, according to \Cref{claim:blobpreservation}.
    Using this fact, we let $W_1^j,W_2^j$ be two components of $W - v^{j-1}_{F_j}$, for each $i_j \in I$, such that for both $h \in [2]$, the tree $W_h^j$ contains a unique child of $v^j_{F_j}$ and $W_h^j$ contains a leaf $v^z_F$ with $F \cap I_h^j \neq \emptyset$.
    We call this leaf $v^j_h$.
    Thanks to the choices we made during our constructions, we can now verify the following claim, which leads us to the end of our proof.

    \begin{claim}\label{claim:toomanyblobs}
        The set $\{ v^j_1,v^j_2 ~\!\colon\!~ j \in [z] \}$ contains at least $z + 1$ vertices.
    \end{claim}
    \emph{Proof of \Cref{claim:toomanyblobs}:}
    We will prove this by proving a stronger statement by induction on $z$.
    In particular, we want to show that the set $\mathcal{S}_h \coloneqq \{ W^j_1,W^j_2 ~\!\colon\!~ j \in [h] \}$ contains $h+1$ pairwise disjoint trees for all $h \in [z]$.
    The case $h=1$ is immediately settled by \Cref{claim:blobpreservation}.
    Suppose therefore that $h > 1$ and $\mathcal{S}_{h-1} \subseteq \mathcal{S}_h$ contains $h$ pairwise disjoint trees $W_1, \ldots , W_h$.
    If either $W^h_1$ or $W^h_2$ are disjoint from the trees $W_1, \ldots, W_h$, we are done.

    Suppose w.l.o.g.\ that $W_1$ intersects $W^h_1 \cup W^h_2$.
    As the root of $W_1$ has depth at most $h-1$, this must mean that $W_1$ contains the common parent of the roots of $W^h_1$ and $W^h_2$.
    Thus, by definition of the elements of $\mathcal{S}_h$, we must have $W^h_1 \cup W^h_2 \subseteq W_1$ and $W^h_1 \cup W^h_2$ is therefore disjoint from the trees $W_2, \ldots , W_h$.
    We conclude that $W^h_1, W^h_2, W_2, \ldots , W_h$ are $h+1$ pairwise disjoint trees in $\mathcal{S}_h$.
    \hfill$\blacksquare$

    However, as $z \geq \mathsf{d}_t$, this contradicts the properties of the set $J$ we got from our application of \Cref{thm:ApathsArgument}.
    Thus we will find a $K_t$-minor with the desired properties.
\end{proof}

Our last lemma will actually not require us to use \Cref{thm:ApathsArgument}.
The argument we present here is fairly technical, since we will have to delve into the internals of \Cref{def:transactionmesh}.
But the intuition of our approach is simple.
Similar to \Cref{lemma:rooting_areas}, we will find an extended surface-wall that captures each of the areas inside of the transactions that host crosses which intersect the apices of $\mathfrak{M}$ (see \textbf{\textsf{M5}} of \Cref{def:transactionmesh}) and since $a' \geq k$, this lets us find a $K_t$-minor with the desired properties thanks to \Cref{lemma:cliqes_in_extended_Dyckwalls}.

\begin{lemma}\label{lem:toomanyprotectedcrosses}
    Let $t,s,s',a,p,\ell$ be non-negative integers with $k \coloneqq \nicefrac{1}{2}(t-3)(t-4)$, $s' \geq 5t + 74$, $s \geq s'\ell+(8t+96)k+2t+25$, and $p \geq 2s$.
    Let $(G,\Omega)$ be a society with an $(s,s',a,p,\ell)$-transaction mesh $\mathfrak{M} = (H, \Omega, \mathcal{C}, \mathfrak{P})$ and let $(a',z,f,h)$ be the $\mathsf{d}_t$-signature of $\mathfrak{M}$.

    If $a' \geq k$ the graph $G$ contains a $K_t$-minor model controlled by a mesh whose horizontal paths are subpaths of distinct cycles from $\mathcal{C}$.
    Furthermore, there exists an algorithm that finds this $K_t$-minor model in time $\mathbf{O}(t^3|E(G)|)$.
\end{lemma}
\begin{proof}
    Let $\mathfrak{T} = ( T, \mathfrak{M}, \phi, \mathcal{A} = \{ A_1, \ldots , A_\ell \}, \mathcal{Z}, \mathcal{F}, \mathfrak{L} )$ be the $\mathfrak{M}$-tree of precision $\mathsf{d}_t$ that witnesses the signature of $\mathfrak{M}$ and let $\mathfrak{P} = \{ \mathcal{P}_1, \ldots , \mathcal{P}_\ell \}$.
    We let $(H_j^i, \Omega_j^i)$ for all $i \in [0,\ell]$ and $j \in [i+1]$ be defined as in \Cref{def:transactionmesh}, with $\mathcal{C}_i^j = \{ C_1^{i,j}, \ldots , C_{s - (s'i + 1)}^{i,j} \}$ being the corresponding nests.
    Let $k \coloneqq \nicefrac{1}{2}(t-3)(t-4)$.
    
    Let $i \in [\ell]$ be chosen such that $A_i \neq \emptyset$.
    As in \Cref{def:transactionmesh}, let $j_i \in [i]$ be the integer such that $\mathcal{P}_i = \{ P_1^i, \ldots , P_p^i \}$ is a planar transaction that is naturally indexed in $(H_i \cup \bigcup \mathcal{P}, \Omega_i)$, where $(H_i, \Omega_i)$ is the $C_{s - s'i}^{j_i,i-1}$-society in the cylindrical rendition of $(H_{j_i}^{i-1}, \Omega_{j_i}^{i-1})$.
    Let $G''$ be the component of $G - (\bigcup_{j = 1}^i V(\Omega_j^{i-1}) \cup \bigcup_{j = 1}^i A_j \cup V(\Omega^*_i))$ that contains $C_1^{j_i,i-1}$, let $G' \coloneqq G[G'' \cup \bigcup_{j = 1}^i V(\Omega_j^{i-1}) \cup V(\Omega^*_i)]$, let $Q^i_L \coloneqq  P^i_{\lfloor \nicefrac{p}{2} \rfloor}$, and $Q^i_R \coloneqq  P^i_{\lceil \nicefrac{(p+1)}{2} \rceil}$.
    Further, let $(U_i, \Psi_i)$ be the $\mathcal{P}_i$-strip society in $(G', \Omega_i)$ and let $(U_i', \Psi_i')$ be the $\{ Q^i_L, Q^i_R \}$-society in $(U_i, \Psi_i)$.
    According to the point \textbf{\textsf{M3}} in \Cref{def:transactionmesh}, $U_i$ is disjoint from $\bigcup_{j \in [i] \setminus \{ j_i \}} \mathcal{C}_j^i$, and since $A_i \neq \emptyset$, the society $(G_i[V(U_i') \cup A_i], \Psi_i')$ has a cross.
    We let $G_i \coloneqq G'[V(U_i') \cup A_i]$.

    Recall that $(H,\Omega)$ has a vortex-free rendition $\rho_H$ in a disk $\Delta$ according to \Cref{def:transactionmesh}.
    We are interested in two more societies in relation to the $i \in [\ell]$ we have chosen above.
    First, consider the paths $S^i_L \coloneqq P^i_{\lfloor \nicefrac{p}{2} \rfloor - (s'-1)}$ and $S^i_R \coloneqq  P^i_{\lceil \nicefrac{(p+1)}{2} \rceil + s' - 1}$.
    Recall that $V(\Omega_i) \subseteq V(C_{s - s'i}^{j_i,i-1})$ and note that $\{ S^i_L, S^i_R \}$ is a transaction in $(H_i \cup S^i_L \cup S^i_R, \Omega_i)$
    Thus $S^i_L \cup S^i_R \cup C_{s - s'i}^{j_i,i-1}$ contain a cycle $C'_i$ such that $S^i_L \cup S^i_R \subseteq C'_i$ and the trace of $C'_i$ in $\rho_H$ separates the traces of $Q^i_L$ and $Q^i_R$ from the boundary of $\Delta$.
    As a consequence, the disk $\Delta_i$ within $\Delta$ defined by the trace of $C'_i$ that contains the traces of $Q^i_L$ and $Q^i_R$ defines the $\Delta_i$-society $(H_i', \Omega_i')$.

    We observe that $V(\Psi_i') \subseteq V(\Omega_i')$ and thus $(H_i' \cup G_i, \Omega_i')$ contains a cross.
    Let $T_L^i,T_R^i$ be the traces of $Q^i_L,Q^i_R$ in $\rho_H$, let $e \in E(Q^i_L) \setminus E(H_i)$ be an edge of $Q^i_L$ found outside of $E(H_i)$ -- which must exist since $Q^i_L$ comes from an exposed transaction -- and let $e_T$ be the curve in $T_L^i$ representing $e$.
    We use $e_T$ to define the inclusion-wise minimal subdisk $\Delta_i' \subseteq \Delta_i$ defined by the trace in $\rho_H$ of a cycle in $S^i_L \cup S^i_R \cup C_{2t + 25}^{j,i-1} \subseteq H$ that contains $e_T$.
    The cycle whose trace in $\rho_H$ bounds $\Delta_i'$ must contain edges of both $S^i_L$ and $S^i_R$.
    We let $(H_i'', \Omega_i'')$ be the $\Delta_i'$-society in $\rho_H$ and let $G_i' \coloneqq H_i'' \cup G_i$.
    
    Consider \textbf{\textsf{M5}} of \Cref{def:transactionmesh} and in particular, the fact that the vertices of $A_i$ have no neighbours in $V(H_i)$.
    We conclude that, since $(H_i' \cup G_i, \Omega_i')$ contains a cross, $(G_i', \Omega_i')$ contains a cross.
    
    Note that $H_i''$ contains subpaths of $P_j^i$ for each $j \in [\lfloor \nicefrac{p}{2} \rfloor - (s'-2), \lceil \nicefrac{(p+1)}{2} \rceil + s' - 2]$ and thus $H_i''$ contains subpaths of $2s' - 2 = 8t+48$ paths from $\mathcal{P}_i$.
    Furthermore, $H_i''$ contains subpaths of the cycles $C_1^{j,i-1}, \ldots , C_{2t + 24}^{j,i-1}$.
    For all $j \in [2t+24]$, this allows us to define a cycle $C_j^i$ via the boundary of the inclusion-wise minimal subdisk of $\Delta_i$ that contains $e_T$ such that
    \[ C_j^{i\star} \subseteq C^{j,i-1}_j \cup P^i_{\lfloor \nicefrac{p}{2} \rfloor - (4t + 48 + j)} \cup P^i_{\lceil \nicefrac{(p+1)}{2} \rceil + 4t + 48 + j} . \]
    We note that the cycles $C_1^{i\star}, \ldots , C_{2t+24}^{i\star}$ are pairwise disjoint.
    Since $G_i'$ is derived from the $\{ Q^i_L, Q^i_R \}$-society in $(U_i, \Psi_i)$ and the vertices $A_i$ have no neighbours in $V(H_i)$, the society $(G_i', \Omega_i')$ has a cylindrical rendition in $\Delta_i'$ in which $\{ C_1^{i\star}, \ldots , C_{2t+24}^{i\star} \}$ is a nest of order $2t+24$.
    This allows us to use \Cref{lemma:makenestcozy} to find a cozy nest $\mathcal{C}_i' = \{ C_1^i, \ldots , C_{2t+24}^i \}$ of order $2t+24$ in $(G_i',\Omega_i')$.
    Furthermore, using the paths $P_j^i$ with $j \in [\lfloor \nicefrac{p}{2} \rfloor - (4t + 48), \lceil \nicefrac{(p+1)}{2} \rceil + 4t + 48]$, we can find a radial linkage $\mathcal{R}_i$ for $\mathcal{C}_i'$ of order $8t+96$ in $G_i'$.

    Let $I \coloneq \{ i \in [\ell] ~\!\colon\!~ A_i \neq \emptyset \}$ and let $I' \subseteq I$ be a subset of $I$ of size exactly $k$, which exists thanks to $a' \geq k$.
    We have thus, for every $i \in I$, constructed a society $(G_i', \Omega_i')$ with a cylindrical rendition in a disk $\Delta_i'$ in which a nest $\mathcal{C}'_i$ of order $10t+121 = 5(2t + 24) + 1$ and a radial linkage $\mathcal{R}_i$ for $\mathcal{C}'_i$ of order $8t+96 = 4(2t + 24)$ exist.
    Furthermore, $(G_i', \Omega_i')$ contains a cross.
    By construction, for any two distinct $i,j \in I$, the graphs $G_i'$ and $G_j'$ are disjoint.
    We let $X_i$ be the $8t+96$ endpoints of $\mathcal{R}_i$ in $V(\Omega_i')$, where we note that we can choose $\mathcal{R}_i$ such that $X_i$ are found in a segment of $\Omega_i'$ that is disjoint from any path $P_j^i \in \mathcal{P}_i$ with $j \not\in [\lfloor \nicefrac{p}{2} \rfloor - (4t + 48), \lceil \nicefrac{(p+1)}{2} \rceil + 4t + 48]$.
    We set $X \coloneqq \bigcup_{i \in I'} X_i$.

    Let $\mathcal{P} \in \mathfrak{P}$ be a transaction of $\mathfrak{M}$ such that there exists a subpath $Y'$ of $C_{s-s'}$ that contains both endpoints of every member of $\mathcal{P}$.
    Such a transaction must exist, since $\ell \geq 1$, according to \Cref{def:transactionmesh}.
    Let $Y$ be the set of all vertices that are the first vertex encountered on $C_{s-s'}$ when traversing along some member of $\mathcal{P}$ starting from $Y'$ in either direction of traversal.

    Using \hyperref[prop:mengersthm]{Menger's Theorem}, we conclude that there exists a linkage $\mathcal{L}$ of order $(8t+96)k$ between $X$ and $Y$, as there are at least that many cycles of $\mathcal{C}$ that have an index lower than $s - s'$, as well as at least that many cycles of each $\mathcal{C}'_i$, paths from $\mathcal{P}$, and paths in each $\mathcal{R}_i$, each themselves contained in $\mathcal{P}_i \in \mathfrak{P}$.
    At least one of these families would need to be fully separated from the others by a separator of size at most $(8t+96)k - 1$, which is impossible.
    Moreover, $\mathcal{L}$ is contained in the inner graph of $C_{s - s'}$ and disjoint from $G_i'' - V(\Omega_i'')$.

    Next, we need to turn $\mathcal{L}$ into a linkage that is orthogonal to the $2t+24$ innermost cycles of each nest $\mathcal{C}'_i$, as well as the cycles $C_{s-(s'+2t+23)}, \ldots , C_{s-s'}$.
    For each $i \in I$, we apply \Cref{lem:radialtoorthogonal} to $\mathcal{R}_i$ and $\mathcal{C}'_i$ to get an orthogonal, radial linkage $\mathcal{R}_i'$ for $\mathcal{C}'_i$ that is end-identical to $\mathcal{R}_i$.
    Let $X_i^\star$ be the set of first vertices encounter on the cycle $C_{2t+25}^i$ when following along the paths in $\mathcal{R}_i'$ starting from $C_1^i$ and let $X^\star = \bigcup_{i \in I} X_i^\star$.
    Then, we follow along $\mathcal{L}$, starting from its endpoints in $C_{s-s'}$ until we meet the cycle $C_{s-(s'+2t+23+(8t+96)k)}$.
    Let $\mathcal{J}$ be the resulting $V(C_{s-s'})$-$V(C_{s-(s'+2t+23+(8t+96)k))}$-linkage between.
    Another application of \Cref{lem:radialtoorthogonal} yields a linkage $\mathcal{J}'$ of the same size which is orthogonal to $\{ C_{s-(s'+2t+23+(8t+96)k)}, \ldots , C_{s-s'} \}$ and end-identical to $\mathcal{J}$.
    Finally, let $Y^\star$ be the set of all vertices on $C_{s-(s'+2t+23)}$ we encounter when traversing along the paths of $\mathcal{J}'$ starting from $C_{s-s'}$.

    Consider the graph $H^\star$ obtained by taking the intersection of the inner graph of $C_{s-(s'+2t+23)}$ and the outer graphs of the cycles $C^i_{2t+25}$ for all $i \in I$ -- taken in the vortex-free rendition of $(H,\Omega)$ in a disk that exists according to \Cref{def:transactionmesh} -- by adding the cycles $C_{s-(s'+2t+23)}$ and $C^i_{2t+25}$ back in for all $i \in I$.
    The graph $H^\star$ still contains the cycles $C_{s-(s'+2t+23+(8t+96)k)}, \ldots , C_{s-s'}$ as well as, for each $i \in I$, the cycles $C^i_{2t+25}, \ldots , C^i_{10t+121}$.
    An application of \Cref{prop:mengersthm} yields an $X^\star$-$Y^\star$-linkage $\mathcal{J}^\star$ of order $(8t+96)k$ in $H^\star$.

    We can now combine $\mathcal{J}^\star$ with the orthogonal linkages $\mathcal{R}_i'$ for all $i \in I$ and the orthogonal linkage $\mathcal{J}'$ to obtain a $V(C_{s-s'})$-$V(\bigcup_{i \in I} C_1^i)$-linkage $\mathcal{L}^\star$ whose intersection with any of the cycles $C_1^i, \ldots , C_{2t+24}^i$, for all $i \in I$, and the cycles $C_{s-(s' + 2t + 23)}$ is a single, non-empty path.
    Moreover $\mathcal{L}^\star$ has $8t+96$ endpoints on each $C_1^i$ for each $i \in I$.

    As a consequence $\{ C_{s - s'}, \ldots , C_{s - (s'+2t+23)}\}$ together with $\mathcal{L}^\star$ and the cycles $\{ C^i_1, \ldots , C^i_{2t+24} \}$ for each $i \in I$ form an extended $(2t+24)$-surface-wall.
    This then yields the desired $K_t$-minor via \Cref{cor:cliquesinlandscapes}.
\end{proof}

\subsection{Classifying societies}\label{sec:classifyingsocieties}
Now that we have gathered numerous tools, we are finally ready to stick them all together to prove the main theorem of this section.
This mainly consists of us applying \Cref{lem:buildtransactionmesh} until we either run into the conditions in which \Cref{lem:everydistrictboundeddepth} tells us that we are done, or until we have accumulated enough iterations that we can use the results from \Cref{sec:signatures} to show that we can find a $K_t$-minor model with the properties we desire.

\begin{proof}[Proof of \Cref{thm:societyclassification}]
    We give the following estimates for the functions in the statement, including some additional functions useful in the proof.
    Recall the definition of $\mathsf{d}_t$ from \Cref{def:Apathsfunction} and note that $\mathsf{d}_t \in \mathbf{O}(t^3)$.
    \begin{align*}
        \mathsf{iter}_{\ref{thm:societyclassification}}(t) \coloneqq                    &\ (t-3)(t-4)\mathsf{d}_t^2 \in \mathbf{O}(t^8) \\
        \mathsf{apex}_{\ref{thm:societyclassification}}(t) \coloneqq                    &\ 8t^3 + 3\mathsf{d}_t^2 + 2\mathsf{d}_t \in \mathbf{O}(t^6) \\
        \mathsf{apex}^\mathsf{genus}_{\ref{thm:societyclassification}}(t) \coloneqq     &\ 8t^3 + \nicefrac{1}{2}(t-3)(t-4) \mathsf{apex}_{\ref{thm:societyclassification}}(t) \in \mathbf{O}(t^8) \\
        \mathsf{loss}^{\mathsf{iter}}_{\ref{thm:societyclassification}}(t) \coloneqq    &\ \mathsf{loss}_{\ref{lem:buildtransactionmesh}}(t,\mathsf{d}_t) \in \mathbf{O}(t^3) \\
        \mathsf{loss}_{\ref{thm:societyclassification}}(t) \coloneqq                    &\ \mathsf{loss}_{\ref{lem:everydistrictboundeddepth}} \in \mathbf{O}(t^3) \\
        \mathsf{cost}_{\ref{thm:societyclassification}}(t,k) \coloneqq                  &\  \max(k, 5t+74) \in \mathbf{O}(t+k) \\
        \mathsf{nest}_{\ref{thm:societyclassification}}(t,k) \coloneqq                  &\  (8t+96)(\nicefrac{1}{2}(t-3)(t-4)) + (2t +24)(8t^3 + 20) \ + \\
                                                                                        &\  \mathsf{cost}_{\ref{thm:societyclassification}}(t,k)\mathsf{iter}_{\ref{thm:societyclassification}}(t) + \mathsf{nest}_{\ref{lem:everydistrictboundeddepth}}(t) \ + \\
                                                                                        &\ \mathsf{nest}_{\ref{lem:buildtransactionmesh}}(t,\mathsf{d}_t,\mathsf{apex}^\mathsf{genus}_{\ref{thm:societyclassification}}(t),\mathsf{cost}_{\ref{thm:societyclassification}}(t,k),\mathsf{iter}_{\ref{thm:societyclassification}}(t),\mathsf{d}_t) \in \mathbf{O}(t^9+t^8k)\\
        \mathsf{link}_{\ref{thm:societyclassification}}(t,k,p) \coloneqq                &\ (64t^3 - 128)(2t+24) + 32\mathsf{d}_t + 8\mathsf{apex}^\mathsf{genus}_{\ref{thm:societyclassification}}(t) + \\
                                                                                        &\ 4\mathsf{nest}_{\ref{thm:societyclassification}}(t,k)(p-1) +  \mathsf{link}_{\ref{lem:everydistrictboundeddepth}}(t,\mathsf{nest}_{\ref{thm:societyclassification}}(t,k),p) \in \mathbf{O}(t^{14} + t^9p + t^8kp) \\
        \mathsf{apex}^\mathsf{fin}_{\ref{thm:societyclassification}}(t,k,p) \coloneqq   &\ \mathsf{apex}_{\ref{lem:everydistrictboundeddepth}}(t,\mathsf{nest}_{\ref{thm:societyclassification}}(t,k),\mathsf{link}_{\ref{thm:societyclassification}}(t,k,p),p)\\
        \mathsf{depth}_{\ref{thm:societyclassification}}(t,k,p) \coloneqq               &\ \mathsf{depth}_{\ref{lem:buildtransactionmesh}}(t,\mathsf{nest}_{\ref{thm:societyclassification}}(t,k),p,\mathsf{link}_{\ref{thm:societyclassification}}(t,k,p),\mathsf{apex}^\mathsf{genus}_{\ref{thm:societyclassification}}(t),\mathsf{iter}_{\ref{thm:societyclassification}}(t),\mathsf{d}_t,\mathsf{d}_t)
    \end{align*}

    We note that $\mathsf{apex}^\mathsf{fin}_{\ref{thm:societyclassification}}(t,k,p) \in \mathbf{O}(t^{35}p^2 + t^{30}p^3 + t^{29}kp^3 + t^{34}kp^2 + t^{28}k^2p^3)$ and furthermore, $\mathsf{depth}_{\ref{thm:societyclassification}}(t,k,p) \in \mathbf{O}(t^{50}p + t^{45}p^2 + t^{44}kp^2 + t^{49}kp + t^{44}kp^2 + t^{43}k^2p^2)$.
    These estimates are a little cumbersome to state but do help us with slightly optimising our analysis for the bounds of the local structure theorem later on.

    Let $p' \coloneqq \mathsf{link}_{\ref{lem:everydistrictboundeddepth}}(t,k,p)$, let $s' \coloneqq \mathsf{cost}_{\ref{thm:societyclassification}}(t,k)$, let $a \coloneqq \mathsf{apex}_{\ref{thm:societyclassification}}(t) $, and let $a' \coloneqq \mathsf{loss}^{\mathsf{iter}}_{\ref{thm:societyclassification}}(t)$.
    
    We start by recalling the observation from \Cref{def:transactionmesh} according to which any society $(G,\Omega)$ with a cylindrical rendition and a nest of order $s$ defines an $(s,s',a,p',0)$-transaction mesh.
    For our $(G,\Omega)$, we let this transaction mesh be $\mathfrak{M}_0$.
    We will now iteratively apply \Cref{lem:buildtransactionmesh} until we find one of our desired outcomes.
    The main task of our proof now is to use the tools we proved in \Cref{sec:signatures} to show that this happens after at most $\mathsf{iter}_{\ref{thm:societyclassification}}(t)$ iterations.

    Let $i \in [0,\mathsf{iter}_{\ref{thm:societyclassification}} - 1]$ and suppose we have already built a sequence of padded transaction meshes $\mathfrak{M}_0, \mathfrak{M}_1, \ldots , \mathfrak{M}_j$ such that for all $j \in [0,i]$, $\mathfrak{M}_j$ is an $(s - j(a'+s'),s',a,p'-2a'j,j)$-transaction mesh with the $\mathsf{d}_t$-signature $\Bar{x}_j$.
    In particular, we suppose that we also have a sequence of $\mathsf{d}_t$-signatures $\Bar{x}_0, \Bar{x}_1, \ldots , \Bar{x}_j$, such that $\Bar{x}_{h-1} \leq \Bar{x}_h$ for all $h \in [i]$.
    For each $j \in [0,i]$, we let $\Bar{x}_j = (a''_j, z_j,f_j,h_j)$.

    If $a''_i \geq \nicefrac{1}{2}(t-3)(t_4)$, then according to \Cref{lem:toomanyprotectedcrosses} we can find a $K_t$-minor model fitting the conditions we laid out in our statement.
    Thus we have $a''_j < \nicefrac{1}{2}(t-3)(t_4)$ for all $j \in [i]$.
    If there exists an interval $[j,l] \subseteq [i]$ such that $a''_h = a''_{h'}$ for all $h,h' \in [j,l]$ and $z_l \geq \mathsf{d}_t$, then, according to \Cref{lem:toomanyblobs}, this again yields the $K_t$-minor model we desire and thus our proof is complete.

    This allows us to conclude two things:
    Firstly, that the union of the apices of the transactions in $\mathfrak{P}_i$ has size at most $a\nicefrac{1}{2}(t-3)(t_4)$.
    Let $\mathfrak{T}_i = ( T_i, \mathfrak{M}_i, \phi_i, \mathcal{A}_i, \mathcal{Z}_i = \{ Z_1^i, \ldots , Z_\ell^i \}, \mathcal{F}_i, \mathfrak{L}_i )$ be the $\mathfrak{M}_i$-tree and let $\mathsf{z} \in [i]$ be the maximal integer such that $A_j = \emptyset$ for all $j \in [(i + 1) - \mathsf{z},i]$.
    Then secondly, $Z_i \coloneqq \bigcup_{j = (\ell + 1) - \mathsf{z}}^\ell B_\mathfrak{T}(Z_j^i)$ has order at most $3(\mathsf{d}_t^2)+2\mathsf{d}_t$.
 
    We now let $\{(G_j,\Omega_j)\}_{j \in [i+1]}$ be the districts of $\mathfrak{M}_i$ each with a cylindrical rendition $\rho_j$, and let $\mathcal{C}_j = \{ C_1^j, \ldots , C_{s - (is' + 1)}^j \}$ be the district-nest of $(G_j,\Omega_j)$ for all $j \in [i+1]$, with $(G_j',\Omega_j')$ being $C_{s - (i+1)s'}^j$-society in $\rho_j$.

    Suppose that for all $j \in [i+1]$ the society $(G_j',\Omega_j')$ has depth less than $\mathsf{depth}_{\ref{thm:societyclassification}}(t,k,p)$.
    We first discuss the case in which $i=0$.
    In this case, we can simply ask for a $V(\Omega)$-$\sigma(c_0)$-linkage of order $4k$ in $G$ via \Cref{prop:mengersthm}.
    Should this not exist, we receive a separator $S$ with $|S| < 4k$, which allows us to find a vortex-free rendition of $(G,\Omega)$ derived from $\rho$ by chucking the remains of $\sigma(c_0)$ into some arbitrary cell.
    This clearly satisfies the fourth option in our statement.
    Thus we may assume that there exist a radial linkage of order $4k$ for $\mathcal{C}$ in $G$, which we can orthogonalise via \Cref{lem:radialtoorthogonal}, again allowing us to satisfy the fourth option in our statement.
    Should $i \neq 0$, we can accomplish the same thing by applying \Cref{lem:everydistrictboundeddepth}.
    This yields one of the four options that we desire in our statement, thanks to definition of $\mathsf{nest}_{\ref{thm:societyclassification}}$ and our choice of $p$.

    Hence we can assume that $(G_j',\Omega_j')$ has depth at least $\mathsf{depth}_{\ref{thm:societyclassification}}(t,k,p)$ for some $j \in [i+1]$, which allows us to apply \Cref{lem:buildtransactionmesh}, again due to our definition for $\mathsf{nest}_{\ref{thm:societyclassification}}$, as well as the constants $p$ and $a$.
    This either yields a $K_t$-minor that satisfies the first point of our theorem, a crosscap transaction of order $p$ that satisfies the second point, or we can find a transaction $\mathcal{P}_{i+1}$ in $(G_j',\Omega_j')$ with a set of co-conspirators which builds a padded $(s - (i+1)(a'+s'),s',a,p' -2a'(i+1),i+1)$-transaction mesh $\mathfrak{M}_{i+1}$ in $(G,\Omega)$ with the $\mathsf{d}_t$-signature $\Bar{x}_{i+1}$.

    As a consequence, we are either done after at most $\mathsf{iter}_{\ref{thm:societyclassification}}(t)-1$ iterations, or we can find transaction meshes $\mathfrak{M}_0, \mathfrak{M}_1, \ldots , \mathfrak{M}_{\mathsf{iter}_{\ref{thm:societyclassification}}(t)}$ with $\mathsf{d}_t$-signatures defined as above.
    Recall that we have already determined that $a''_j < \nicefrac{1}{2}(t-3)(t_4)$ and additionally know that there does not exist an interval $[j,l] \subseteq [\mathsf{iter}_{\ref{thm:societyclassification}}(t)]$ such that $a''_h = a''_{h'}$ for all $h,h' \in [j,l]$ and $z_l \geq \mathsf{d}_t$.

    Suppose there exists an interval $[j,l] \subseteq [\mathsf{iter}_{\ref{thm:societyclassification}}(t)]$ with $|[j,l]| \geq 2\mathsf{d}_t$ such that $a''_i = a''_{i'}$ and $z_i = z_{i'}$ for all $i,i' \in [j,l]$.
    As a consequence, according to the definition of $\mathsf{d}_t$-signatures in \Cref{def:Mtree}, we must have $h_i \leq h_{i+1}$ for all $i \in [j,l-1]$.
    If $h_l > \mathsf{d}_t$, then \Cref{cor:theDcantbetoolarge} tells us that we can find a $K_t$-minor with the desired properties.
    Thus, we can also assume that $h_i \leq \mathsf{d}_t$ for all $i \in [j,l]$.
    We claim that under this restriction, $\mathfrak{M}_l$ allows us to find a $K_t$-minor model with the desired properties.

    \begin{claim}\label{claim:auseforcrookedtransactions}
        For each $i \in [j,l-1]$, if $h_i = h_{i+1}$ then $f_i < f_{i+1}$.
    \end{claim}
    \emph{Proof of \Cref{claim:auseforcrookedtransactions}:}
        According to the definition of $\mathsf{d}_t$-signatures in \Cref{def:Mtree}, for each $i \in [j+1,l]$ for which we have $h_{i-1} = h_i$, the two new districts in $\mathfrak{M}_i$ constructed via $\mathcal{P}_i$ were separated without the deletion of any vertices.
        In particular, no new vertices are deleted starting from the construction of $\mathfrak{M}_j$ all the way until we have constructed $\mathfrak{M}_l$.
        Recall that \Cref{lem:buildtransactionmesh} provides us with a set of co-conspirators $\mathcal{W}_i$ for each $\mathcal{P}_i$.

        Now for each $i \in [j,l]$, we let $(G_{j_i},\Omega_{j_i})$ be the district of $\mathfrak{M}_{i-1}$ with the cylindrical rendition $\rho_i$ such that $\mathcal{P}_i$ is a transaction in the $C^i_{s - (i+1)s')}$-society $(G_{j_i}',\Omega_{j_i}')$ of $\rho_i$.
        Further, let $(G_1^i, \Omega_1^i),(G_2^i,\Omega_2^i)$ be the two unique, distinct districts of $\mathfrak{M}_i$ such that $V(\Omega_{i'}^i) \cap V(\Omega_{j_i}') \neq \emptyset$ for both $i' \in [2]$.
        Note that therefore, among the co-conspirators $\mathcal{W}_i$ there exist both paths with endpoints in $V(\Omega_{i'}^i)$ for both $i' \in [2]$.
        Due to $V(\Omega_1^i)$ and $V(\Omega_2^i)$ being separated in the union $U_i$ of the districts of $\mathfrak{M}_i$, we know that $\mathcal{W}_i$ cannot be a snitch.

        Let $F \in \mathcal{F}_{i-1}$ be the set with $j_i \in F$ from $\mathfrak{T}_{i-1}$, which we defined earlier.
        Due to $V(\Omega_1^i)$ and $V(\Omega_2^i)$ being separated in $U_i$ and the fact that we did not remove any vertices inside of the vortices of the districts to construct $\mathfrak{M}_i$ from $\mathfrak{M}_{i-1}$, we know that $|\mathcal{F}_i| = |\mathcal{F}_{i-1}| + 1$.
        In particular, if we let $F_1, F_2 \in \mathcal{F}_i$ be the two distinct sets such that $F_{i'}$ contains the integer corresponding to the district $(G_{i'}^i,\Omega_{i'}^i)$ for both $i' \in [2]$, we have $\mathcal{F}_{i-1} \setminus \{ F \} = \mathcal{F}_i \setminus \{ F_1, F_2 \}$.

        Now, if for both $i' \in [2]$ we either have that $(G_{i'}^i,\Omega_{i'}^i)$ has a cross or $|F_{i'}| \geq 2$, our claim is prove.
        Thus we may suppose this is not the case.
        This in particular means that $\mathcal{W}_i$ cannot contain any crosses and must thus be a pair of deserters.
        However, these deserters now themselves verify that these two districts are connected to other districts, which means that we must have $|F_{i'}| \geq 2$ for both $i' \in [2]$.
        Thus our claim holds.
    \hfill$\blacksquare$

    As a consequence of this claim and the fact that $h_l \leq \mathsf{d}_t$, we know that $f_l > \mathsf{d}_t$, since $|[j,l]| \geq 2\mathsf{d}_t$.
    Thanks to \Cref{lem:toomanyfamilies}, we know that this means that we can find a $K_t$-minor model as desired.
    Thus since $\mathsf{d}_t > 72t^3 + \nicefrac{1}{2}(t-3)(t-4) - 145 = 8t^3-16 + 64t^3-128 + \nicefrac{1}{2}(t-3)(t-4) - 1$, we must instead find the $K_t$-minor model.
    As a consequence no interval $[j,l] \subseteq [\mathsf{iter}_{\ref{thm:societyclassification}}(t)]$ with $|[j,l]| \geq 2\mathsf{d}_t$ such that $a''_i = a''_{i'}$ and $z_i = z_{i'}$ for all $i,i' \in [j,l]$ can exist.
    However, our choice for $\mathsf{iter}_{\ref{thm:societyclassification}}(t)$ makes it impossible for all of the restrictions we have placed on the signatures of our transactions meshes to hold at the same time.
    We must therefore be able to find a $K_t$-minor model in the way we desire and thus our proof is complete.
\end{proof}


\section{The local structure theorem}\label{sec:localstructure}

We finally have everything in place to start discussing the proof for the Local Structure Theorem.
On a high level, our proof will start with a large mesh.
We then call the Flat Wall Theorem, that is \cref{thm:flatmesh}.
If this does not yield a clique-minor, we have found a large flat mesh $M$.
A large portion of $M$ will now be reserved to make sure our tangle agrees with the -- eventually -- embedded part of the graph.
The remaining part of $M$ is used to define a first society $(G,\Omega)$ together with a large nest and a radial linkage such that $(G,\Omega)$ has a cylindrical rendition and everything that does not belong to the compass of our flat mesh $M$ is confined into this vortex.

From here on the proof follows the route already illustrated in \cref{fig:intro_flowchart}.
We will maintain a $\Sigma$-landscape $\Lambda$ with a single vortex where $\Sigma$ starts out as the sphere.
We iteratively apply the Society Classification Theorem (that is \cref{thm:societyclassification}) to $(G,\Omega)$.
This yields one of four outcomes;
\begin{enumerate}
    \item A clique-minor which would allow us to terminate the process,
    \item a small apex set and a flat crosscap transaction,
    \item a small apex set and a flat handle-transaction, or
    \item a small apex set and a dissolution of our unique vortex into a small number of bounded depth vortices, each with a nest and a linkage reaching back into the nest of $(G,\Omega)$.
\end{enumerate}
In cases (ii) and (iii) we will call upon one of two technical lemmas -- one for the orientable and one for the non-orientable case -- to either attach a handle or a crosscap to $\Sigma$, thereby increasing its Euler-genus, to route a small part of the transaction back onto our original mesh, and to install a new society $(G'\Omega')$, again with a large nest and a unique vortex.
By doing this we maintain our landscape and are able to iterate until either the Euler-genus of our surface is so large that \cref{cor:universal-surface-walls} provides us with the forbidden minor, or \cref{thm:societyclassification} returns a rendition of $(G',\Omega')$ in the disk with bounded breadth and depth.

Recall the definitions of landscapes and layouts from \cref{subsec:landscapes}.
Our version of the Local Structure Theorem reads as follows.

\begin{theorem}\label{thm:strongest_localstructure}
There exist functions $\mathsf{apex}_{\ref{thm:strongest_localstructure}},\mathsf{depth}_{\ref{thm:strongest_localstructure}}\colon\mathbb{N}^2\to\mathbb{N}$ and $\mathsf{mesh}_{\ref{thm:strongest_localstructure}}\colon\mathbb{N}^3\to\mathbb{N}$ such that for all non-negative integers $k$, $t\geq 5$, and $r$, every graph $H$ on $t$ vertices, every graph $G$ and every $\mathsf{mesh}_{\ref{thm:strongest_localstructure}}(t,r,k)$-mesh $M\subseteq G$ one of the following holds.
\begin{enumerate}
    \item \textsl{$G$ has an $H$-minor controlled by $M$}, or
    \item $G$ has a $k$-$(\mathsf{apex}_{\ref{thm:strongest_localstructure}}(t,k),\nicefrac{1}{2}(t-3)(t-4),\mathsf{depth}_{\ref{thm:strongest_localstructure}}(t,k),r)$-$\Sigma$-layout $\Lambda$ centred at $M$ where $\Sigma$ is a surface \textsl{where $H$ does not embed} and $\delta$ is $(M-A)$-central where $A$ is the apex set of $\Lambda$.
\end{enumerate}
Moreover, it holds that

{\centering
  $ \displaystyle
    \begin{aligned}
        \mathsf{apex}_{\ref{thm:strongest_localstructure}}(t,k),~ \mathsf{depth}_{\ref{thm:strongest_localstructure}}(t,k) \in \mathbf{O}\big((t+k)^{112}\big), \text{ and } \mathsf{mesh}_{\ref{thm:strongest_localstructure}}(t,r,k) \in \mathbf{O}\big((t+k)^{115}r \big) .
    \end{aligned}
  $
\par}

There also exists an algorithm that, given $t,k,r$, a graph $H$, a graph $G$ and a mesh $M$ as above as input finds one of these outcomes in time $\mathbf{poly}(t+k+r)|E(G)||V(G)|^2$.
\end{theorem}

In the following we first present a proof of a slightly weaker statement, namely \cref{thm:localstructure}.
The proof can easily be adapted to a proof of \cref{thm:strongest_localstructure} as explained in \cref{subsec:wollasnsconjectures}.

\subsection{Extending the surface}\label{subsec:surfaceextension}

In order to be able to extend our surface whenever we find a crosscap or handle transaction, we require a way to synchronise this crosscap or handle with the surface we are currently on.
So far, all of our tools have been tailored to make our graph similar to the surface as possible.
This time around the surface has to follow the lead of our graph.

We import two crucial lemmas from Kawarabayashi et al.\ \cite{KawarabayashiTW2021Quickly} which are strong enough to not require any additional refinement.
Indeed, their lemmas are already exactly in the format we need them to be.
Since both lemmas are long and rather technical to state we dedicate this section purely to their statement and some explanation.

\begin{figure}[ht]
    \centering
    \begin{tikzpicture}

        \pgfdeclarelayer{background}
		\pgfdeclarelayer{foreground}
			
		\pgfsetlayers{background,main,foreground}

        \begin{pgfonlayer}{background}
            \pgftext{\includegraphics[width=8cm]{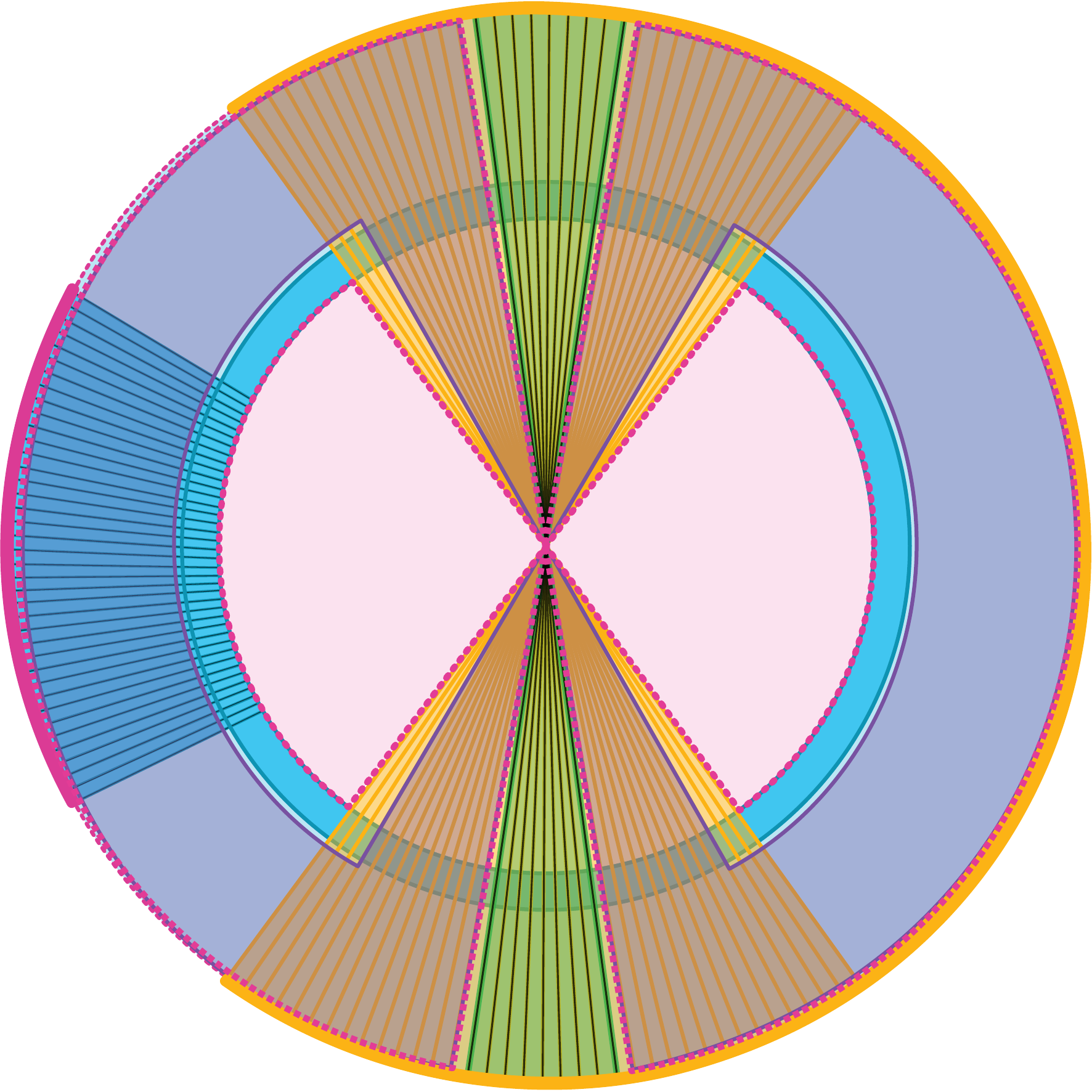}};
        \end{pgfonlayer}{background}
			
        \begin{pgfonlayer}{main}
        \node (C) [v:ghost] {};

        \node (Omega) [v:ghost,position=135:45mm from C] {$(G,\Omega)$};

        \node (X1) [v:ghost,position=180:43mm from C] {$X_1$};
        \node (X2) [v:ghost,position=0:43mm from C] {$X_2$};

        \node (c_11) [v:ghost,position=0:15mm from C] {$c_1'$};
        \node (c_1) [v:ghost,position=225:31mm from C] {$c_1$};

        \node (PP) [v:ghost,position=112.5:43mm from C] {$\mathcal{P}$};
        \node (PPP) [v:ghost,position=270:43mm from C] {$\mathcal{P}'$};

        \node (RR) [v:ghost,position=200:43mm from C] {$\mathcal{R}$};

        \node (CC) [v:ghost,position=40:25.7mm from C] {$\mathcal{C}$};
        \node (CCC) [v:ghost,position=335:33.5mm from C] {$\mathcal{C}'$};
        
        \end{pgfonlayer}{main}
        
        \begin{pgfonlayer}{foreground}
        \end{pgfonlayer}{foreground}

    \end{tikzpicture}
    \caption{The outcome of \cref{lemma:integrate-crosscap} when applied to the society $(G,\Omega)$ with nest $\mathcal{C}$ and radial linkage $\mathcal{R}$. The purple overlay indicates how the nest $\mathcal{C}'$ surrounds the vortex $c_1'$ in the cylindrical rendition of the vortex society of $c_1$.}
    \label{fig:crosscapflip}
\end{figure}
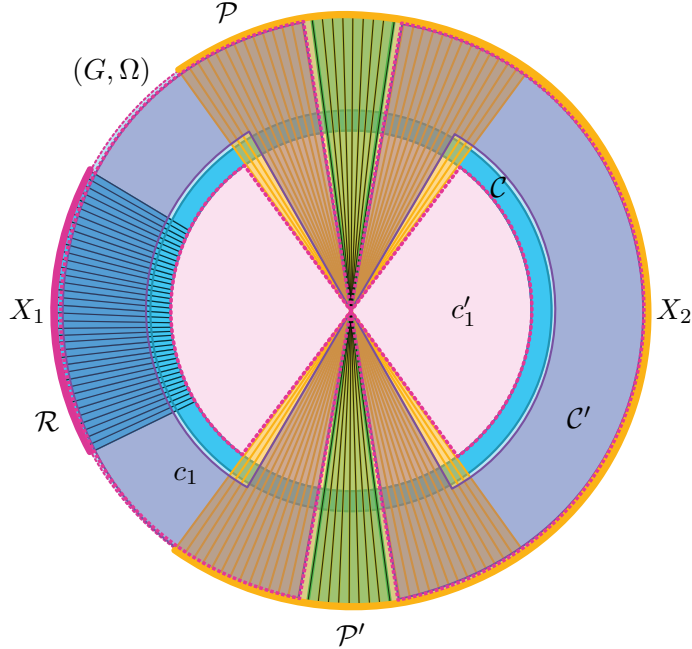

\begin{proposition}[Kawarabayashi, Thomas, and Wollan \cite{KawarabayashiTW2021Quickly}]\label{lemma:integrate-crosscap}
Let $s$ and $p$ be non-negative integers.
Let $(G,\Omega)$ be a society with a cylindrical rendition $\rho_0$ in the disk $\Delta$ with a nest $\mathcal{C}=\{ C_1,\dots,C_{s+9}\}$ around the vortex $c_0$.
Let $X_1,X_2$ be disjoint segments of $\Omega$ such that there exist
\begin{itemize}
    \item a radial linkage $\mathcal{R}$ orthogonal to $\mathcal{C}$ starting in $X_1$, and
    \item a flat crosscap transaction $\mathcal{P}$ of order at least $p+2s+7$ with all endpoints in $X_2$ and disjoint from $\mathcal{R}$.
\end{itemize}
Let $\Sigma^*$ be a surface, homeomorphic to the projective plane minus an open disk, which is obtained from $\Delta$ by adding a crosscap to the interior of $c_0$.
\smallskip
Then there exists a crosscap transaction $\mathcal{P}'\subseteq \mathcal{P}$ of order $p$, consisting of the middle $p$ paths of $\mathcal{P}$, and a rendition $\rho_1$ of $(G,\Omega)$ in $\Sigma^*$ with a unique vortex $c_1$ and the following hold:
\begin{enumerate}
    \item $\mathcal{P}'$ is disjoint from $\sigma(c_1)$,

    \item the vortex society of $c_1$ in $\rho_1$ has a cylindrical rendition $\rho_1'$ with a nest $\mathcal{C}'=\{ C_1',\dots,C_s'\}$ around the unique vortex $c_1'$,

    \item every element of $\mathcal{P}$ has an endpoint in $V(\sigma_{\rho_1'}(c_1'))$,

    \item $\mathcal{R}$ is orthogonal to $\mathcal{C}'$ and for every $i\in[s]$ and every $R\in\mathcal{R}$, $C_i'\cap R= C_{i+8}\cap R$. Moreover,

    \item let $\mathcal{R}=\{ R_1,\dots,R_{\ell}\}$. For each $i\in[\ell]$ let $x_i$ be the endpoint of $R_i$ in $X_1$, and let $y_i$ be the last vertex of $R_i$ on $c_1$ when traversing along $R_i$ starting from $x_i$; then if $x_1,x_2,\dots,x_{\ell}$ appear in $\Omega$ in the order listed, then $y_1,y_2,\dots,y_{\ell}$ appear on $\widetilde{c_1}$ in the order listed.

    \item Finally, let $\Delta'$ be the open disk bounded by the trace of $C_{s+8}$ in $\rho_0$.
    Then $\rho_0$ restricted to $\Delta\setminus\Delta'$ is equal to $\rho_1$ restricted to $\Delta\setminus \Delta'$.
\end{enumerate}
Moreover, there exists an algorithm that computes this outcome in time $\mathbf{poly}(sp|\mathcal{R}|)|V(G)|$.
\end{proposition}

To summarise, \cref{lemma:integrate-crosscap} allows to isolate some part of a flat crosscap transaction and to pull it ''inside'' the already embedded part of the graph.
By doing so, the lemma cuts out a disk from the initial disk $\Delta$, then glues in a Moebius-strip, and finally cuts out another disk to allow for a new vortex together with its nest.
This newly created society which contains the new nest and vortex is guaranteed to be fully disjoint from the ``protected'' part of the crosscap transaction that was pulled in.
Moreover, the lemma allows to propagate a radial linkage through the newly added Moebius-strip.
These are the key features that will allow for the construction of our surface-wall later on.
See \cref{fig:crosscapflip} for an illustration.

Please note that the original statement of \cref{lemma:integrate-crosscap} in \cite{KawarabayashiTW2021Quickly} requires $\mathcal{R}$ to end on the boundary of the vortex $c_0$ but only needs $s+8$ cycles in the initial nest instead of $s+9$.
By requiring one additional cycle we do not need $\mathcal{R}$ to extend all the way to $c_0$ since we can declare the $c_0$-disk of the innermost cycle to be the vortex.
This, however, costs one cycle which is taken account for by the slight increase of the number of required cycles in the nest.
This also explains why we shift, in iv), by $+8$ instead of $+7$ as in the original statement.
\medskip

Next we introduce the analogue of \cref{lemma:integrate-crosscap} for handle transactions.
While the statement is somewhat different, the idea remains the same.

\begin{figure}[ht]
    \centering
    \begin{tikzpicture}

        \pgfdeclarelayer{background}
		\pgfdeclarelayer{foreground}
			
		\pgfsetlayers{background,main,foreground}

        \begin{pgfonlayer}{background}
            \pgftext{\includegraphics[width=8cm]{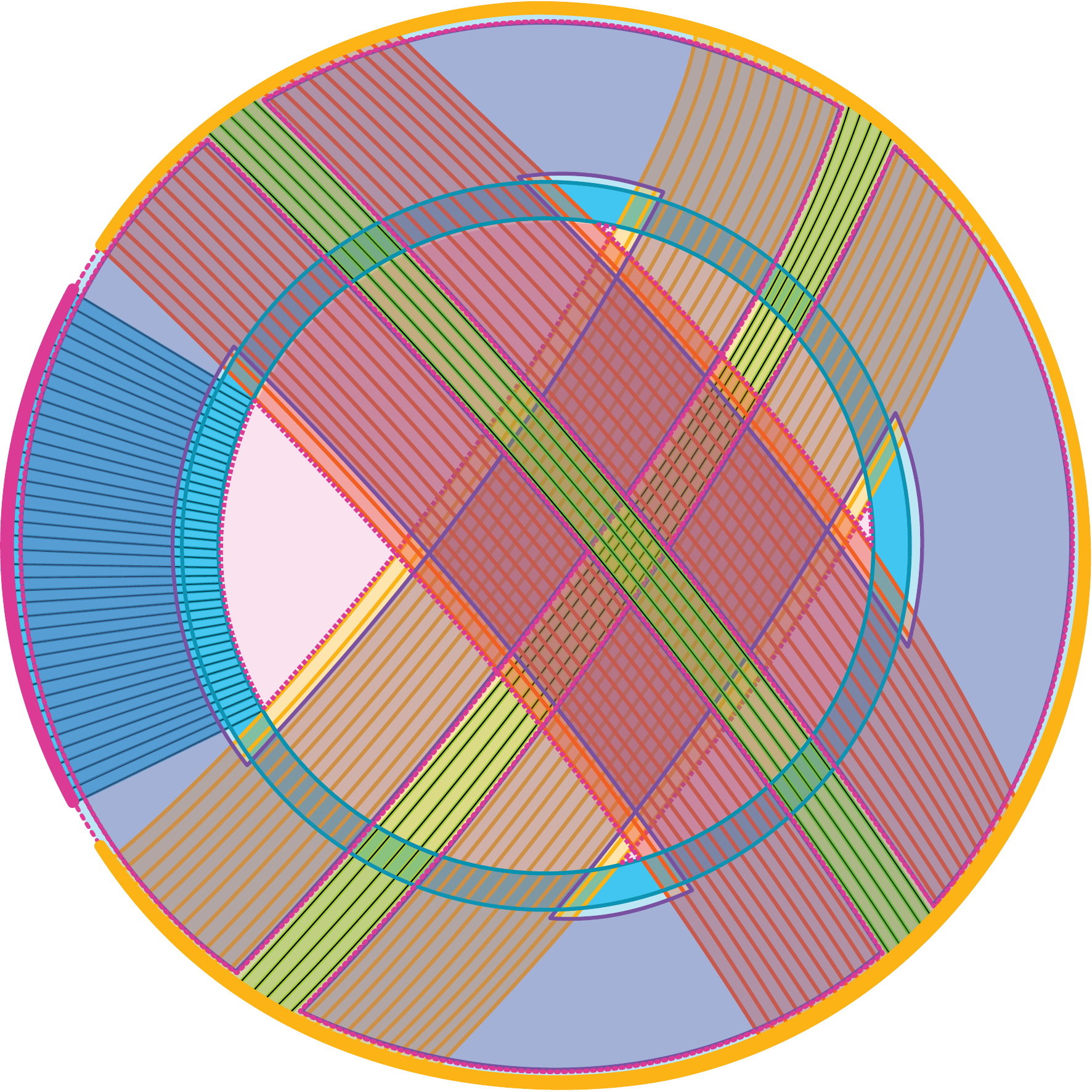}};
        \end{pgfonlayer}{background}
			
        \begin{pgfonlayer}{main}
        \node (C) [v:ghost] {};

        \node (Omega) [v:ghost,position=90:44mm from C] {$(G,\Omega)$};

        \node (X1) [v:ghost,position=180:43mm from C] {$X_1$};
        \node (X2) [v:ghost,position=280:43mm from C] {$X_2$};

        \node (c_11) [v:ghost,position=180:17mm from C] {$c_1'$};
        \node (c_1) [v:ghost,position=225:31mm from C] {$c_1$};

        \node (PP1) [v:ghost,position=135:43mm from C] {$\mathcal{P}_1$};
        \node (PPP1) [v:ghost,position=313.5:43.5mm from C] {$\mathcal{P}_1'$};
        
        \node (PP2) [v:ghost,position=70:43mm from C] {$\mathcal{P}_2$};
        \node (PPP2) [v:ghost,position=240:43mm from C] {$\mathcal{P}_2'$};

        \node (RR) [v:ghost,position=200:43mm from C] {$\mathcal{R}$};

        \node (CC) [v:ghost,position=170:25.7mm from C] {$\mathcal{C}$};
        \node (CCC) [v:ghost,position=0:33.5mm from C] {$\mathcal{C}'$};
        
        \end{pgfonlayer}{main}
        
        \begin{pgfonlayer}{foreground}
        \end{pgfonlayer}{foreground}

    \end{tikzpicture}
    \caption{The outcome of \cref{lemma:integrate-handle} when applied to the society $(G,\Omega)$ with nest $\mathcal{C}$ and radial linkage $\mathcal{R}$. The purple overlay indicates how the nest $\mathcal{C}'$ surrounds the vortex $c_1'$ in the cylindrical rendition of the vortex society of $c_1$.}
    \label{fig:handleflip}
\end{figure}
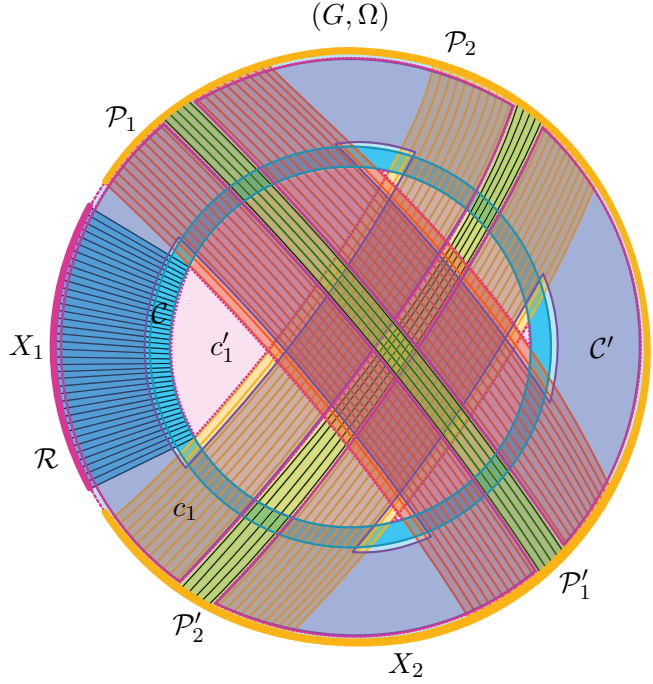

\begin{proposition}[Kawarabayashi, Thomas, and Wollan \cite{KawarabayashiTW2021Quickly}]\label{lemma:integrate-handle}
Let $s$ and $p$ be non-negative integers.
Let $(G,\Omega)$ be a society with a cylindrical rendition $\rho_0$ in the disk $\Delta$ with a nest $\mathcal{C}=\{ C_1,\dots,C_{s+9}\}$ around the vortex $c_0$.
Let $X_1,X_2$ be disjoint segments of $\Omega$ such that there exist
\begin{itemize}
    \item a radial linkage $\mathcal{R}$ orthogonal to $\mathcal{C}$ starting in $X_1$, and
    \item a flat handle transaction $\mathcal{P}$ of order at least $2p+4s+12$ with all endpoints in $X_2$ and disjoint from $\mathcal{R}$.
\end{itemize}
Let $\mathcal{P}_1$ and $\mathcal{P}_2$ be the two planar transactions such that $\mathcal{P}=\mathcal{P}_1\cup\mathcal{P}_2$.
Let $\Sigma^+$ be a surface, homeomorphic to the torus minus an open disk, which is obtained from $\Delta$ by adding a handle to the interior of $c_0$.
\smallskip
Then there exist transactions $\mathcal{P}_1'\subseteq \mathcal{P}_1$ and $\mathcal{P}_2'\subseteq\mathcal{P}_2$ of order $p$, each $\mathcal{P}_i'$ consisting of the middle $p$ paths of $\mathcal{P}_i$, $i\in[2]$, such that $\mathcal{P}'=\mathcal{P}_1'\cup\mathcal{P}_2'$ is a handle transaction, and a rendition $\rho_1$ of $(G,\Omega)$ in $\Sigma^+$ with a unique vortex $c_1$ and the following hold:
\begin{enumerate}
    \item $\mathcal{P}'$ is disjoint from $\sigma(c_1)$,

    \item the vortex society of $c_1$ in $\rho_1$ has a cylindrical rendition $\rho_1'$ with a nest $\mathcal{C}'=\{ C_1',\dots,C_s'\}$ around the unique vortex $c_1'$,

    \item every element of $\mathcal{P}$ has an endpoint in $V(\sigma_{\rho_1'}(c_1'))$,

    \item $\mathcal{R}$ is orthogonal to $\mathcal{C}'$ and for every $i\in[s]$ and every $R\in\mathcal{R}$, $C_i'\cap R= C_{i+8}\cap R$. Moreover,

    \item let $\mathcal{R}=\{ R_1,\dots,R_{\ell}\}$. For each $i\in[\ell]$ let $x_i$ be the endpoint of $R_i$ in $X_1$, and let $y_i$ be the last vertex of $R_i$ on $c_1$ when traversing along $R_i$ starting from $x_i$; then if $x_1,x_2,\dots,x_{\ell}$ appear in $\Omega$ in the order listed, then $y_1,y_2,\dots,y_{\ell}$ appear on $\widetilde{c_1}$ in the order listed.

    \item Finally, let $\Delta'$ be the open disk bounded by the trace of $C_{s+8}$ in $\rho_0$.
    Then $\rho_0$ restricted to $\Delta\setminus\Delta'$ is equal to $\rho_1$ restricted to $\Delta\setminus \Delta'$.
\end{enumerate}
Moreover, there exists an algorithm that computes this outcome in time $\mathbf{poly}(sp|\mathcal{R}|)|V(G)|$.
\end{proposition}

For an illustration of the setting of \cref{lemma:integrate-handle} see \cref{fig:handleflip}.
As before, our statement of \cref{lemma:integrate-handle} is slightly changed when compared to its original form as stated by Kawarabayashi et al.\ \cite{KawarabayashiTW2021Quickly}.
Again, we increase the size of the initial nest from $s+8$ by one to $s+9$ to not have to maintain that our radial linkage reaches all the way to the vortex $c_0$.

\subsection{The proof of the Local Structure Theorem}\label{subsec:localstructureproof}

We proceed with a proof of a first version of the Local Structure Theorem with polynomial bounds.

\begin{theorem}\label{thm:localstructure}
There exist functions $\mathsf{apex}_{\ref{thm:localstructure}},\mathsf{depth}_{\ref{thm:localstructure}}\colon\mathbb{N}^2\to\mathbb{N}$ and $\mathsf{mesh}_{\ref{thm:localstructure}}\colon\mathbb{N}^3\to\mathbb{N}$ such that for all non-negative integers $k$, $t\geq 5$, and $r$, every graph $G$ and every $\mathsf{mesh}_{\ref{thm:localstructure}}(t,r,k)$ mesh $M\subseteq G$ one of the following holds.
\begin{enumerate}
    \item $G$ has a $K_t$-minor controlled by $M$, or
    \item $G$ has a $k$-$(\mathsf{apex}_{\ref{thm:localstructure}}(t),\nicefrac{1}{2}(t-3)(t-4),\mathsf{depth}_{\ref{thm:localstructure}}(t,k),r)$-$\Sigma$-layout centred at $M$ where $\Sigma$ is a surface of Euler-genus at most $t^2$.
\end{enumerate}
Moreover, it holds that

{\centering
  $ \displaystyle
    \begin{aligned}
        \mathsf{apex}_{\ref{thm:localstructure}}(t,k),~ \mathsf{depth}_{\ref{thm:localstructure}}(t,k) \in \mathbf{O}\big((t+k)^{112}\big), \text{ and } \mathsf{mesh}_{\ref{thm:localstructure}}(t,r,k) \in \mathbf{O}\big((t+k)^{115}+t^3r \big) .
    \end{aligned}
  $
\par}

There also exists an algorithm that, given $t,k,r$, a graph $G$ and a mesh $M$ as above as input finds one of these outcomes in time $\mathbf{poly}(t+k+r)|E(G)||V(G)|^2$.
\end{theorem}

Note that, while \cref{thm:localstructure} is already strong enough to prove the conjecture of Kawarabayashi, Thomas, and Wollan for polynomial bounds for the GMST, the bound on the Euler-genus it guarantees is not optimal.
Moreover, it is only concerned with clique-minors rather than arbitrary graphs $H$.
A more refined and even stronger version of the Local Structure Theorem, namely the main theorem of this section, can be found in \cref{thm:strongest_localstructure}.
In \cref{subsec:wollasnsconjectures} we explain how to adapt the proof of \cref{thm:localstructure} to derive \cref{thm:strongest_localstructure}.
Here, we also implement tight bounds on the Euler-genus.

Let us begin by providing estimates for the functions mentioned in \cref{thm:localstructure}.
For this, we first need to fix four intermediate values as follows.
Let $\mathsf{c}$ be the constant from \cref{cor:universal-surface-walls} such that every $x$-surface-wall $W$ with $h$ handles and $c$ crosscaps where $2h+c=t^2$ and $x=\mathsf{c}t^{12}$ contains a $K_t$-minor controlled by the cylindrical wall of $W$.
The following functions determine the order of the transactions and the order of the nests we will need to find throughout the proof.
Both functions depend on a parameter $g$ we use to keep track on the number of ``genus-increasing'' steps we have performed in order to facilitate an inductive proof.
\begin{align*}
    \mathsf{radial}(g,t,k)\coloneqq~& (g+2)(8k+8\mathsf{c}t^{12}+4+ \nicefrac{4k}{2}(t-3)(t-4)+1)\\
    \mathsf{nest}(g,t,k)\coloneqq~ &(g+2)\big(8k+8\mathsf{c}t^{12} + \mathsf{radial}(g-1,t,k) +14+ \mathsf{cost}_{\ref{thm:societyclassification}}(t,4(k+\mathsf{c}t^{12}+1))\\
    &+\mathsf{loss}_{\ref{thm:societyclassification}}(t) \big) + \mathsf{nest}_{\ref{thm:societyclassification}}(t,4(k+\mathsf{c}t^{12}+1)) + \nicefrac{4k}{2}(t-3)(t-4)+1\\
    \mathsf{transaction}(g,t,k)\coloneqq~ & 4k+4\mathsf{c}t^{12}+24 + 4\mathsf{nest}(g-1,t,k)+2\mathsf{radial}(g-1,t,k)
\end{align*}
With these values fixed, the functions of \cref{thm:localstructure} are as follows.
\begin{align*}
    \mathsf{apex}_{\ref{thm:localstructure}}(t,k) \coloneqq~& t^2\cdot \mathsf{apex}^\mathsf{genus}_{\ref{thm:societyclassification}}(t) + \mathsf{apex}^\mathsf{fin}_{\ref{thm:societyclassification}}(t,4(k+\mathsf{c}t^{12}+1),\mathsf{transaction}(t^2,t,k)) + 16t^3\\
    \mathsf{depth}_{\ref{thm:localstructure}}(t,k) \coloneqq~& \mathsf{depth}_{\ref{thm:societyclassification}}(t,4(k+\mathsf{c}t^{12}+1),\mathsf{transaction}(t^2,t,k))\\
    \mathsf{mesh}_{\ref{thm:localstructure}}(t,r,k) \coloneqq ~& 100t^3\Big(\mathsf{apex}_{\ref{thm:localstructure}}(t,k)+2\mathsf{nest}(t^2,t,k)+2t+17+r\\&
    +\nicefrac{1}{2}(t-3)(t-4)(2\mathsf{depth}_{\ref{thm:societyclassification}}(t,4(k+\mathsf{c}t^{12}+1),\mathsf{transaction}(t^2,t,k))+1)\Big)
\end{align*}
Using the bounds provided in the proof of \Cref{thm:societyclassification}, we can give the following estimates for the functions above:
We have $\mathsf{radial}(g,t,k) \in \mathbf{O}(gk + gt^{12})$ and $\mathsf{nest}(g,t,k), \mathsf{transaction}(g,t,k) \in \mathbf{O}(t^{20} + t^{12}g^2 + t^8k + g^2k)$. 
This in turn allows us to determine that $\mathsf{apex}_{\ref{thm:localstructure}}(t,k) \in \mathbf{O}((t+k)^{112})$, $\mathsf{depth}_{\ref{thm:localstructure}}(t,k) \in \mathbf{O}(t+k)^{107}$, and $\mathsf{mesh}_{\ref{thm:localstructure}}(t,r,k) \in \mathbf{O}( (t+k)^{115} + t^3r )$.

\begin{proof}[Proof of \cref{thm:localstructure}]
Let $M$ be a $\mathsf{mesh}_{\ref{thm:localstructure}}(t,k,r)$-mesh in $G$.
We begin by applying \cref{thm:flatmesh} to $M$.
This yields, in time $\mathbf{poly}(t+k+r)|V(G)|$, either a $K_t$-minor controlled by $M$ or we find a set $A_0\subseteq V(G)$ of size at most $16t^3$ and an $x$-submesh $M'$ of $M$ with
\begin{align*}
 x \geq ~& \frac{1}{100t^3}\mathsf{mesh}_{\ref{thm:localstructure}}(t,k,r)-2t-4\\
 \geq ~& \mathsf{apex}_{\ref{thm:localstructure}}(t,k) +9+ r + 2\mathsf{nest}(t^2,t,k)\\
 ~&+\nicefrac{1}{2}(t-3)(t-4)(2\mathsf{depth}_{\ref{thm:societyclassification}}(t,4(k+\mathsf{c}t^{12}+1),\mathsf{transaction}(t^2,t,k))+1)
\end{align*}
Let $\mathcal{C}_0'=\{ {C^0}'_1,\dots,{C^0}'_{\mathsf{nest}(t^2,t,k)}\}$ be the set of the $\mathsf{nest}(t^2,t,k)$ vertex-disjoint cycles obtained by iteratively peeling of the perimeters ${C^0}'_i$ of the $x-2(i-1)$-submeshes of $M'$ that are obtained by this procedure, for $i\in[\mathsf{nest}(t^2,t,k)]$.
Moreover, let $M_0$ be the $\big(\mathsf{apex}_{\ref{thm:localstructure}}(t,k)+\nicefrac{1}{2}(t-3)(t-4)(2\mathsf{depth}_{\ref{thm:societyclassification}}(t,4(k+\mathsf{c}t^{12}+1),\mathsf{transaction}(t^2,t,k))+1)+11+ r\big)$-submesh of $M'$ that is left over after removing the cycles in $\mathcal{C}'_0$ from $M'$ and let $U_0$ be the perimeter of $M_0$.

\begin{observation}\label{obs:localstrucobs1}
Let $I\subseteq [\mathsf{nest}(t^2,t,k)]$ with $|I|=k$ and let $O$ be a cylindrical $k$-mesh with base cycles ${C^0}'_i$, $i\in I$, and rails taken from subpaths of the paths of $M$.
Then the tangle $\mathcal{T}_O$ is a truncation of the tangle $\mathcal{T}_{M_0}$.
\end{observation}

Let $\Sigma_0$ be the sphere.
As $M'$ is flat in $G$, there exists a $\Sigma_0$-decomposition $\delta_0$ for $G$ with one vortex such that $M'$ is grounded in $\delta$ and the trace of the perimeter of $M'$ bounds a closed disk $\Delta_0'\subseteq\Sigma_0$ which is disjoint from the vortex of $\delta$ and every vertex $x\in V(M')$ that is a node in $\delta$ is drawn in $\Delta_0'$.

Let $X_0$ be the set of all vertices of the perimeter $U_0$ of $M_0$ that are grounded and let $\Delta_0''$ be the \textsl{open} subdisk of $\Delta_0$ which is bounded by the trace of $U_0$ and disjoint from the vortex of $\delta_0$.
Moreover, let $G_0$ be the subgraph of $G$ induced by all vertices of $G$ drawn on $\Delta_0'\setminus \Delta_0''$.
Finally, let $\Omega$ be the cyclic ordering of $X_0$ obtained by traversing along $U_0$ in clockwise direction.
It follows that the society $(G_0,\Omega_0)$ has a cylindrical rendition $\rho_0$ in the disk $\Delta_0$ with the nest $\mathcal{C}_0'$ around the unique vortex $c_0$.
Indeed, we may select $\mathsf{radial}(t^2,t,k)$ many paths from the mesh $M'$ to obtain a radial linkage $\mathcal{R}_0'$ which links vertices from distinct rows on $U_0$ to the cycle ${C^0}'_1$

By using \cref{lemma:makenestcozy} we may find a \textsl{cozy} nest $\mathcal{C}_0=\{ C_1^0,\dots,C^0_{\mathsf{nest}(t^2,t,k)}\}$ around $c_0$ in $\rho_0$ in time $\mathbf{O}(\mathsf{nest}(t^2,t,k)|E(G)|^2)$.
By using \cref{lem:radialtoorthogonal} we can also find a radial linkage $\mathcal{R}_0$ of order $\mathsf{radial}(t^2,t,k)$ which is end-identical with $\mathcal{R}_0'$ and orthogonal to $\mathcal{C}_0$ in time
$\mathbf{O}(\mathsf{radial}(t^2,t,k)\mathsf{nest}(t^2,t,k)|E(G)|)$.
Moreover, if $\overline{G_0}$ is the graph $G-(V(G_0)\setminus V(\Omega)$, then $(\overline{G_0},\Omega_0)$ has a vortex-free rendition in the surface obtained from $\Sigma_0$ by removing an open disk.
Finally, let $\mathcal{C}^*\coloneqq \{ C^0_{\mathsf{nest}(t^2,t,k)},\dots,C^0_{\mathsf{nest}(t^2,t,k)-k-\mathsf{c}t^{12}+1}\}$ be the outermost $4(k+\mathsf{c}t^{12}+1)$ cycles of $\mathcal{C}^0$.
Then $(G_0,\Omega_0,\mathcal{C}^*,\mathcal{R}_0,\emptyset)$ is a surface configuration of $(G_0,\Omega_0)$ for the surface $\Sigma_0$ of strength $(4(k+\mathsf{c}t^{12}+1),\mathsf{radial}(t^2,t,k))$.

\begin{observation}\label{obs:localstrucobs2}
Let $I\subseteq [\mathsf{nest}(t^2,t,k)]$ with $|I|=k$ and let $O$ be a cylindrical $k$-mesh with base cycles ${C^0}_i$, $i\in I$, and rails taken from subpaths of the paths of $M$.
Then the tangle $\mathcal{T}_O$ is a truncation of the tangle $\mathcal{T}_{M_0}$.
\end{observation}

\paragraph{Inductively refining an near embedding: The set-up.}
Observe that the tuple $\Lambda_0=(A_0,M_0,D_0,\delta_0)$ is, in fact, a $\Sigma_0$-landscape of detail $4(k+\mathsf{c}t^{12}+1)$. where $D_0$ is a $(4(k+\mathsf{c}t^{12}+1))$-surface-wall without handles or crosscaps whose base cycles coincide with the cycles in $\mathcal{C}^*$, this wall exists as a consequence of \cref{obs:surface-configs-to-walls}.

We now prove that that given the following collection of objects for some $i\in[0,t^2-1]$, we can apply \cref{thm:societyclassification} to the society $(G_i,\Omega_i)$ to either find a $K_t$-minor, create all of the objects below for $i+1$, or, through the last outcome of  \cref{thm:societyclassification}, create a $\Sigma_i$-layout centred at $M$.
Since we have already shown how to reach the point $i=0$, we will eventually be able to assume that this last outcome of \cref{thm:societyclassification} comes true.

For now, let us introduce the objects we are looking for.
Let $i\in[0,t^2-1]$ be given together with
\begin{itemize}
    \item an apex set $A_i\subseteq V(G)$ with $A_{i-1}\subseteq A_i$, $A_{-1}\coloneqq\emptyset$, and $|A_i|\leq i\cdot \mathsf{apex}^\mathsf{genus}_{\ref{thm:societyclassification}}(t)+16t^3$,
    \item a surface $\Sigma_i$ obtained from the sphere by adding a total of $i$ handles and crosscaps in some combination,
    \item a $\Sigma_i$-decomposition $\delta_i$ with a unique vortex $c_i$,
    \item a $(4(k+\mathsf{c}t^{12}+1))$-surface wall $D_i$ with the amount of handles and crosscaps used to obtain $\Sigma_i$ from the sphere such that the base cycles of $D_i$ coincide with the cycles of $\mathcal{C}^*$,
    \item a society $(G_i,\Omega_i)$ such that there exists a $\delta_i$-aligned disk $\Delta_i$ whose boundary intersects $\delta_i$ exactly in $V(\Omega_i)$ and the restriction $\rho_i$ of $\delta_i$ to $\Delta_i$ is a cylindrical rendition of $(G_i,\Omega_i)$ with $c_i$ being its unique vortex,
    \item a cozy nest $\mathcal{C}_i=\{C^i_1,\dots,C^i_{\mathsf{nest}(t^2-i,t,k)} \}$ of order $\mathsf{nest}(t^2-i,t,k)$ in $\rho_i$ around $c_i$,
    \item a family of transactions $\mathfrak{P}_i=\{ \mathcal{P}_1,\dots,\mathcal{P}_i\}$ on $(G_0,\Omega_0)$ such that $\mathcal{P}_j$ is a handle or crosscap transaction,
    \item a radial linkage $\mathcal{R}_i$ of order $\mathsf{radial}(t^2-i,t,k)$ orthogonal to $\mathcal{C}_i\cup\mathcal{C}^*$ whose endpoints on $(G_0,\Omega_0)$ coincide with some of the endpoints of $\mathcal{R}_0$ and which is disjoint from the paths in $\bigcup\mathfrak{P}_i$, and
    \item all objects above are chosen such that $(G_0-A_i,\Omega_0,\mathcal{C}^*,\mathcal{R}_i,\mathfrak{P}_i)$ is a $\Sigma_i$-configuration of strength $(4(k+\mathsf{c}t^{12}+1),\mathsf{radial}(t^2-i,t,k),p_1,\dots,p_i)$ with $p_j=4k+4\mathsf{c}t^{12}$.
\end{itemize}

By \cref{obs:surface-configs-to-walls}, the last point ensures the existence of the $(4(k+\mathsf{c}t^{12}+1))$-surface-wall $D_i$.
Moreover, if we ever reach the point $i=t^2$, \cref{cor:universal-surface-walls} guarantees that we have found a $K_t$-minor controlled by $D_i$.

\paragraph{Inductively refining an near embedding: Construction}
So suppose for some $i\in[0,t^2-1]$ the list of objects as above has already been constructed.

We begin by applying \cref{thm:societyclassification} to $(G_i,\Omega_i)$.
What follows is a case distinction on the different possible outcomes.
\medskip

\textbf{Case 1:} The outcome is \textsl{i)} of \cref{thm:societyclassification} which yields a $K_t$-minor controlled, and therefore controlled, by a mesh whose horizontal paths are subpaths of distinct cycles of $\mathcal{C}_i$.
It is easy to see that this implies that this $K_t$-minor is controlled by $M_0$ as desired.
\medskip

\textbf{Case 2:}
Let $(G',\Omega')$ be the $C^i_{\mathsf{nest}(t^2-i,t,k)-\mathsf{cost}_{\ref{thm:societyclassification}}(t,4(k+\mathsf{c}t^{12}+1))-1}$-society in $\rho_i$.
Case 2 contains two subcases that can be treated in almost identical ways.
Those are outcome \textsl{ii)} and outcome \textsl{iii)} of \cref{thm:societyclassification} where the both outcomes yield a set $A\subseteq V(G_i)$ of size at most $\mathsf{apex}^\mathsf{genus}_{\ref{thm:societyclassification}}(t)$.
\smallskip

Moreover, in the first outcome $\mathcal{Q}$ is a flat crosscap transaction of order $\mathsf{transaction}(t^2-i,t,k)$ in $(G'-A,\Omega')$ together with a nest $\mathcal{C}'$ in $\rho_i$ of order $\mathsf{nest}(t^2-i,t,k)-\mathsf{cost}_{\ref{thm:societyclassification}}(t,4(k+\mathsf{c}t^{12}+1))-\mathsf{loss}_{\ref{thm:societyclassification}}(t)$ around $c_i$ to which $\mathcal{Q}$ is orthogonal.
\smallskip

In the second outcome, $\mathcal{Q}$ is a flat handle transaction of order $\mathsf{transaction}(t^2-i,t,k)$ in $(G'-A,\Omega')$ together with a nest $\mathcal{C}'$ in $\rho_i$ of order $\mathsf{nest}(t^2-i,t,k)-\mathsf{cost}_{\ref{thm:societyclassification}}(t,4(k+\mathsf{c}t^{12}+1))-\mathsf{loss}_{\ref{thm:societyclassification}}(t)$ around $c_i$ to which $\mathcal{Q}_1$ and $\mathcal{Q}_2$ are orthogonal.
\smallskip

In both cases we can reduce the order of $\mathcal{Q}$ by $2\mathsf{radial}(t^2-i-1,t,k)$ by shedding off, in case one the last $2\mathsf{radial}(t^2-i+1,t,k)$ paths of $\mathcal{Q}$ and in case two, from each $\mathcal{Q}_j$, $j\in[2]$, the last $\mathsf{radial}(t^2-i-1,t,k)$.
Let $\mathcal{Q}'$, $\mathcal{Q}'_1$, and $\mathcal{Q}'_2$, be the resulting transactions.
Moreover, notice that we may select a radial linkage $\mathcal{L}$ of order $\mathsf{radial}(t^2-i-1,t,k)$ in $(G'-A,\Omega')$ which is orthogonal to $\mathcal{C}'$.
Indeed, we may select $\mathcal{Q}'$ and $\mathcal{L}$ such that there are segments $I_1$, $I_2$ of $\Omega'$ where $\mathcal{Q}'$ has all endpoints in $V(I_1)$, $\mathcal{L}$ has its endpoints in $V(\Omega')$ in $V(I_2)$, and $I_1$ and $I_2$ are disjoint. 
It follows that $|\mathcal{Q}'|=4k+4\mathsf{c}t^{12}+14 + 4\mathsf{nest}(t^2-i-1,t,k)$ and $|\mathcal{Q}_1'|=|\mathcal{Q}_2'|=2k+2\mathsf{c}t^{12}+7 + 2\mathsf{nest}(t^2-i-1,t,k)$ in case $\mathcal{Q}'$ is a handle transaction.

In the following we discuss only how to proceed in the case where $\mathcal{Q}'$ is a crosscap transaction.
The second case, namely the one where $\mathcal{Q}'$ is a handle transactions, can be handles analogously with the exception that, instead of \cref{lemma:integrate-crosscap} one needs to use \cref{lemma:integrate-handle} and some of the choices we make below need to be made for the transactions $\mathcal{Q}'_i$ and the objects derived from them individually.
Upon closer inspection, the reader will see that the size of $\mathcal{Q}'$ (and in particular of $\mathcal{Q}_1'$ and $\mathcal{Q}_2'$) was chosen to allow the application of any of these two lemmas.

Let us assume $\mathcal{Q}'=\{ Q_1,\dots,Q_{4k+4\mathsf{c}t^{12}+14 + 4\mathsf{nest}(t^2-i-1,t,k)}\}$ to be numbered such that $Q_j$ separates $Q_{j-1}$ from $Q_{j+1}$ within the strip of $\mathcal{Q}'$ for all $j\in[2,4k+4\mathsf{c}t^{12}+14 + 4\mathsf{nest}(t^2-i-1,t,k)-1]$.
Let $a\coloneqq 2\mathsf{nest}(t^2-i-1,t,k)+7$ and $b\coloneqq 2\mathsf{nest}(t^2-i-1,t,k)+4k+4\mathsf{c}t^{12}+8$.
Moreover, let $\mathcal{A}\coloneqq \{ Q_1,\dots,A_{a-1}\}$, $\mathcal{B}=\{ Q_{b+1},\dots,Q_{4k+4\mathsf{c}t^{12}+14 + 4\mathsf{nest}(t^2-i-1,t,k)}\}$, and $\mathcal{I}\coloneqq \{ Q_a,\dots,Q_b\}$.

We assume the cycles in $\mathcal{C}'$ to be numbered $C_1',\dots, C_{\mathsf{nest}(t^2-i,t,k)-\mathsf{cost}_{\ref{thm:societyclassification}}(t,4(k+\mathsf{c}t^{12}+1))-\mathsf{loss}_{\ref{thm:societyclassification}}(t)}'$ from innermost to outermost.
Let $I=[\mathsf{nest}(t^2-i,t,k)-\mathsf{cost}_{\ref{thm:societyclassification}}(t,4(k+\mathsf{c}t^{12}+1))-\mathsf{loss}_{\ref{thm:societyclassification}}(t)-3,\mathsf{nest}(t^2-i,t,k)-\mathsf{cost}_{\ref{thm:societyclassification}}(t,4(k+\mathsf{c}t^{12}+1))-\mathsf{loss}_{\ref{thm:societyclassification}}(t)-8k+8\mathsf{c}t^{12}-\mathsf{radial}(t^2-i-1,t,k)-9]\subseteq [\mathsf{nest}(t^2-i,t,k)-\mathsf{cost}_{\ref{thm:societyclassification}}(t,4(k+\mathsf{c}t^{12}+1))-\mathsf{loss}_{\ref{thm:societyclassification}}(t)]$.
Then $|I|=8k+8\mathsf{c}t^{12}+\mathsf{radial}(t^2-i-1,t,k)+6=2|\mathcal{I}|+|\mathcal{L}|+2$.
Observe that this means, in particular,
\begin{align*}
    2|\mathcal{I}|+|\mathcal{L}| =~& 8k+8\mathsf{c}t^{12} +4 + \mathsf{radial}(t^2-i-1,t,k) \leq \mathsf{radial}(t^2-i,t,k) = |\mathcal{R}_i|.
\end{align*}
Let us now select $\mathcal{I}'$ to be all $V(\Omega')$-$V(C_1')$-subpaths of the paths in $\mathcal{I}$.
It follows that $\mathcal{I}'$ is a radial linkage of order $2|\mathcal{I}|$.
Moreover, $\mathcal{I}'\cup\mathcal{L}$ is a radial linkage of order $8k+8\mathsf{c}t^{12}+\mathsf{radial}(t^2-i-1,t,k)+4=|I|-2$ which is orthogonal to $\mathcal{C}'$.
This allows us to call \cref{lemma:connected_linkages} for $\mathcal{R}_i$, $\mathcal{I}'\cup\mathcal{L}$, $\mathcal{C}'$, and $I$.
As a result we obtain, in time $\mathbf{poly}(t+k)|E(G)|$, a radial linkage $\mathcal{R}'$ of order $|I|-2$ which shares its endpoints on $V(\Omega')$ with the endpoints of $\mathcal{R}_i$ and its endpoints on $C_1'$ with $\mathcal{I}'\cup\mathcal{L}$.
Moreover, $\mathcal{R}'$ is orthogonal to $\mathcal{C}'\setminus \{ C_i ~\colon~ i\in I\}$ and within the inner graph, let us call it $G''$, of $C_{\mathsf{nest}(t^2-i,t,k)-\mathsf{cost}_{\ref{thm:societyclassification}}(t,4(k+\mathsf{c}t^{12}+1))-\mathsf{loss}_{\ref{thm:societyclassification}}(t)-8k+8\mathsf{c}t^{12}-\mathsf{radial}(t^2-i-1,t,k)-9}'$ it is disjoint from $\mathcal{Q}'\setminus\mathcal{I}$.
Let $\mathcal{J}\subseteq \mathcal{R}'$ be all those paths that do not meet $\mathcal{I}$ within $G''$.
Then $|\mathcal{J}|=|\mathcal{R}'|-8k+8\mathsf{c}t^{12}-4=\mathsf{radial}(t^2-i-1,t,k)$.

Next we inspect our remaining nest a bit closer.
Let $\mathcal{C}''\subseteq \mathcal{C}$ be the set of the innermost $\mathsf{nest}(t^2-i,t,k)-\mathsf{cost}_{\ref{thm:societyclassification}}(t,4(k+\mathsf{c}t^{12}+1))-\mathsf{loss}_{\ref{thm:societyclassification}}(t)-8k+8\mathsf{c}t^{12}-\mathsf{radial}(t^2-i-1,t,k)-10$ cycles of $\mathcal{C}'$.
Let $\Omega''$ be the cyclic ordering of the ground vertices of $C'_{\mathsf{nest}(t^2-i,t,k)-\mathsf{cost}_{\ref{thm:societyclassification}}(t,4(k+\mathsf{c}t^{12}+1))-\mathsf{loss}_{\ref{thm:societyclassification}}(t)-8k+8\mathsf{c}t^{12}-\mathsf{radial}(t^2-i-1,t,k)-9}$ and let $(G'',\Omega'')$ the the resulting society with cylindrical rendition $\rho''$, which is the restriction of $\rho'$ to the $c_i$-disk $\Delta'$ in $\Delta_i$ bounded by the trace of $C'_{\mathsf{nest}(t^2-i,t,k)-\mathsf{cost}_{\ref{thm:societyclassification}}(t,4(k+\mathsf{c}t^{12}+1))-\mathsf{loss}_{\ref{thm:societyclassification}}(t)-8k+8\mathsf{c}t^{12}-\mathsf{radial}(t^2-i-1,t,k)-9}$.
Moreover, let $\mathcal{J}'$ and $\mathcal{Q}''$ be the restrictions of $\mathcal{J}$ and $\mathcal{Q}'$ to $G''$.
\begin{align*}
|\mathcal{C}''|=~&\mathsf{nest}(t^2-i,t,k)-\mathsf{cost}_{\ref{thm:societyclassification}}(t,4(k+\mathsf{c}t^{12}+1))-\mathsf{loss}_{\ref{thm:societyclassification}}(t)-8k+8\mathsf{c}t^{12}\\
&-\mathsf{radial}(t^2-i-1,t,k)-10\\
=~& (t^2-i-1)\big(8k+8\mathsf{c}t^{12}+\mathsf{radial}(t^2-i,t,k)+24+\mathsf{cost}_{\ref{thm:societyclassification}}(t,4(k+\mathsf{c}t^{12}+1))+\mathsf{loss}_{\ref{thm:societyclassification}}(t) \big)\\
&+ \mathsf{nest}_{\ref{thm:societyclassification}}(t,4(k+\mathsf{c}t^{12}+1)) +\nicefrac{4k}{2}(t-3)(t-4)+15\\
\geq~& (t^2-i-1)\big(8k+8\mathsf{c}t^{12}+\mathsf{radial}(t^2-i-1,t,k)+24+\mathsf{cost}_{\ref{thm:societyclassification}}(t,4(k+\mathsf{c}t^{12}+1))+\mathsf{loss}_{\ref{thm:societyclassification}}(t) \big)\\
&+ \mathsf{nest}_{\ref{thm:societyclassification}}(t,4(k+\mathsf{c}t^{12}+1))+\nicefrac{4k}{2}(t-3)(t-4)+15\\
=~& \mathsf{nest}(t^2-i-1)+14
\end{align*}
and
\begin{align*}
|\mathcal{Q}''| =~& 4k+4\mathsf{c}t^{12}+14+4\mathsf{nest}(t^2-i-1,t,k)\\
=~& 2 (2k+2\mathsf{c}t^{12})+4\mathsf{nest}(t^2-i-1,t,k)+14
\end{align*}
Notice that the calculations above suffice for the application of \cref{lemma:integrate-crosscap}.
Indeed, we exceed the necessary numbers by more than $2\mathsf{nest}(t^2-i-1,t,k)$.
If instead we were dealing with a handle transaction we would next want to apply \cref{lemma:integrate-handle}.
Here we would have constructed two transactions, namely $\mathcal{Q}''_1$ and $\mathcal{Q}''_2$, each of order at least $2k+2\mathsf{c}t^{12}+2\mathsf{nest}(t^2-i-1,t,k)+7$ and thus, together, they form a transaction $\mathcal{Q}''$ of the order above.
This explains out choice for the size of $\mathcal{Q}''$.

So we are now ready to apply \cref{lemma:integrate-crosscap} to $(G'',\Omega'')$, $\mathcal{C}''$, $\mathcal{J}''$, and $\mathcal{Q}''$ with $p=4k+4\mathsf{c}t^{12}+2$, rendition $\rho'$ and the disk $\Delta'$.
As a result we obtain the following list of objects.
Let $\Sigma^*$ be a surface homeomorphic to the projective plane minus an open disk which is obtained from $\Delta'$ by adding a crosscap to the interior of $c_i$.
Then there exists $\mathcal{I}'\subseteq \mathcal{Q}''$ which is exactly the restriction of $\mathcal{I}$ to $G''$, and a rendition $\rho'''$ of $(G'',\Omega'')$ in $\Sigma^*$ with a unique vortex $c_{i+1}'$ and the following properties hold:
\begin{itemize}
    \item $\mathcal{I}'$ is disjoint from $\sigma(c_{i+1}')$,
    \item the vortex society $(G_{i+1},\Omega_{i+1})$ of $c_{i+1}'$ in $\rho'''$ has a cylindrical rendition $\rho_{i+1}$ with nest $\mathcal{C}_{i+1}=\{ C_1^{i+1},\dots,C^{i+1}_{\mathsf{nest}(t^2-i-1)} \}$ around the unique vortex $c_{i+1}$,
    \item every element of $\mathcal{Q}''$ has an endpoint in $V(\sigma_{\rho_{i+1}(c_{i+1})})$, and
    \item $\mathcal{J}''$ is orthogonal to $\mathcal{C}_{i+1}$. Moreover,
    \item let $\mathcal{J}''=\{J_1,\dots, J_{\mathsf{radial}(t^2-i-1,t,k)}\}$.
    For each $j\in[\mathsf{radial}(t^2-i-1,t,k)]$ let $x_j$ be the endpoint of $J_j$ in $V(\Omega'')$ and let $y_j$ be the endpoint of $J_j$ on $c_{i+1}'$; then if $x_1,x_2,\dots,x_{\mathsf{radial}(t^2-i-1,t,k)}$ appear on $\Omega''$ in the order listed, then $y_1,y_2,\dots,y_{\mathsf{radial}(t^2-i-1,t,k)}$ appear on $\widetilde{c}_{i+1}'$ in the order listed.
    \item Finally, let $\Delta''$ be the open disk bounded by the trace of of the outermost cycle of $\mathcal{C}''$ in $\rho'$.
    Then $\rho'$ restricted to $\Delta'\setminus \Delta''$ is equal to $\rho'''$ restricted to $\Delta'\setminus\Delta''$.
\end{itemize}
Now let $\delta_{i+1}$ be obtained by first unifying the renditions $\rho'''$ and $\rho_{i+1}$ along the vortex society of $c_{i+1}'$, then combining the resulting rendition of $(G'',\Omega'')$ in the disk $\Delta'$ with the rendition $\rho'$ along to boundary of $\Delta''$, and then reintegrating $\rho'$ into $\delta_i$.
Moreover, let $\Sigma_{i+1}$ be obtained from $\Sigma_i$ be first removing the interior of $\Delta'$ and replacing it with $\Sigma^*$.
Notice that this means that $\Sigma_{i+1}$ is obtained from $\Sigma_i$ by adding a single crosscap.
We also set $A_{i+1}\coloneqq A_i\cup A$ and obtain 
\begin{align*}
|A_{i+1}|\leq~& i\cdot\mathsf{apex}^\mathsf{genus}_{\ref{thm:societyclassification}}(t)+16t^3+\mathsf{apex}^\mathsf{genus}_{\ref{thm:societyclassification}}(t)\\
\leq ~& (i+1)\mathsf{apex}^\mathsf{genus}_{\ref{thm:societyclassification}}(t)+16t^3.
\end{align*}
So the first three points of our invariant are maintained.
\smallskip

Notice that, by construction, we have that the radial linkage $\mathcal{R}'$ can be extended onto the restriction of $\mathcal{R}_i$ to the outer graph of the closed curve obtained by following along the vertices on $\Omega'$.
This allows us to, firstly, extend the crosscap transaction $\mathcal{I}$ of order $4k+4\mathsf{c}t^{12}+2$ along $\mathcal{R}_i$ to obtain a crosscap transaction $\mathcal{P}_{i+1}$ on $(G_0,\Omega_0)$ which is disjoint from $\mathcal{P}_j$, $j\in[i]$,
secondly, to extend $\mathcal{R}'$ along $\mathcal{R}_i$ to form the radial linkage $\mathcal{R}_{i+1}$ which is orthogonal to both $\mathcal{C}^*$ and $\mathcal{C}_{i+1}$ and satisfies
\begin{align*}
 |\mathcal{R}_{i+1}| = \mathsf{radial}(t^2-i-1,t,k),
\end{align*}
its endpoints coincide with some of the endpoints of $\mathcal{R}_0$, and it is disjoint from the paths in $\mathfrak{P}_{i+1}=\mathfrak{P}_i \cup \{ \mathcal{P}_{i+1} \}$.
Indeed, it is straight forward to see that $(G_0-A_{i+1},\Omega_0,\mathcal{C}^*,\,\mathcal{R}_{i+1},\mathfrak{P}_{i+1})$ is a $\Sigma_{i+1}$-configuration of strength $(4(k+\mathsf{c}t^{12}+1),\mathsf{radial}(t^2-i,t,k),p_1,\dots,p_{i+1})$ with $p_j=4k+4\mathsf{c}t^{12}$.
So the last three points of our invariant are also satisfied.
\smallskip

By \cref{obs:surface-configs-to-walls} the $\Sigma_{i+1}$-configuration found in the previous paragraph yields the existence of a $(4(k+\mathsf{c}t^{12}+1))$-surface wall $D_{i+1}$ with the amount of crosscaps and handles used to obtain $\Sigma_{i+1}$ from the sphere such that the base cycles of $D_{i+1}$ coincide with the cycles of $\mathcal{C}^*$.
This, together with the society $(G_{i+1},\Omega_{i+1})$, the cylindrical rendition $\rho_{i+1}$ and the nest $\mathcal{C}_{i+1}$ around the unique vortex $c_{i+1}$ of both $\rho_{i+1}$ and $\delta_{i+1}$, satisfies the remaining three points of our invariant and thus this step is complete.
\smallskip

Observe that, in the case where $i=t^2-1$ we have reached a situation where $K_t$ embeds in $\Sigma_{i+1}$.
This means that \cref{cor:universal-surface-walls} implies the existence of a $K_t$-model controlled by $D_{i+1}$.
So in this case we would be done immediately.
Hence, we may assume that $i+1\in[t^2-1]$.
\bigskip

\textbf{Case 3:}
The remaining case to be discussed is outcome \textsl{iv)} of \cref{thm:societyclassification}.
Here we are given a set $A \subseteq V(G_i)$ with $|A| \leq \mathsf{apex}_{\ref{thm:societyclassification}}(t,4(k+\mathsf{c}t^{12}+1),\mathsf{transaction}(t^2,t,k))$, a rendition $\rho'$ of $(G_i - A, \Omega_i)$ in $\Delta_i$ with breadth $b \in [\nicefrac{1}{2}(t-3)(t-4)-1]$ and depth at most $\mathsf{depth}_{\ref{thm:societyclassification}}(t,4(k+\mathsf{c}t^{12}+1),\mathsf{transaction}(t^2-i,t,k))$, and an extended $(4(k+\mathsf{c}t^{12}+1))$-Dyck-wall $D$ with signature $(0,0,b)$, such that $D$ is grounded in $\rho'$, the base cycles of $D$ are the cycles $C_{\mathsf{nest}(t^2-i,t,k)-\mathsf{cost}_{\ref{thm:societyclassification}}(t,4(k+\mathsf{c}t^{12}+1))-1-k-\mathsf{c}t^{12}}^i,\dots,C_{\mathsf{nest}(t^2-i,t,k)-\mathsf{cost}_{\ref{thm:societyclassification}}(t,4(k+\mathsf{c}t^{12}+1))-1}^i$, and there exists a bijection between the vortices $v$ of $\rho'$ and the vortex segments $S_v$ of $D$, where $v$ is the unique vortex contained in the disk $\Delta_{S_v}$ defined by the trace of the inner cycle of the nest of $S_v$ where $\Delta_{S_v}$ is chosen to avoid the trace of the simple cycle of $D$.

Notice that the union of all rails of $D$ form a radial linkage $\mathcal{R}''$ of order $16b(k+\mathsf{c}t^{12})$ such that every vortex segment provides exactly $16(k+\mathsf{c}t^{12})$ of these paths.
Moreover, there are at precisely $4b(k+\mathsf{c}t^{12})+4$ cycles among $C_{\mathsf{nest}(t^2-i,t,k)-\mathsf{cost}_{\ref{thm:societyclassification}}(t,4(k+\mathsf{c}t^{12}+1))-1-k-\mathsf{c}t^{12}}^i,\dots,C_{\mathsf{nest}(t^2-i,t,k)-\mathsf{cost}_{\ref{thm:societyclassification}}(t,4(k+\mathsf{c}t^{12}+1))-1}^i$.
Now let $\mathcal{R}'$ be a radial linkage of order $4b(k+\mathsf{c}t^{12})$ formed by selecting $4(k+\mathsf{c}t^{12}+1)$ rails from each vortex segment of $D$.

Notice that, by definition of $D$, for every $j\in[\mathsf{nest}(t^2-i,t,k)-\mathsf{cost}_{\ref{thm:societyclassification}}(t,4(k+\mathsf{c}t^{12}+1))-1-k-\mathsf{c}t^{12},\mathsf{nest}(t^2-i,t,k)-\mathsf{cost}_{\ref{thm:societyclassification}}(t,4(k+\mathsf{c}t^{12}+1))-1]$ and every vortex $v$ of $\rho$, both the nest of $S_v$ and $v$ itself are disjoint from $C^i_j$.

Let $I=[\mathsf{nest}(t^2-i,t,k)-\mathsf{cost}_{\ref{thm:societyclassification}}(t,4(k+\mathsf{c}t^{12}+1))-1-k-\mathsf{c}t^{12}+2,\mathsf{nest}(t^2-i,t,k)-\mathsf{cost}_{\ref{thm:societyclassification}}(t,4(k+\mathsf{c}t^{12}+1))-1]]$.
Then we have $|I|=4b(k+\mathsf{c}t^{12})+2$.
Therefore, we may call upon \cref{lemma:connected_linkages} for $\mathcal{C}'$, $I$, $\mathcal{R}'$ and $\mathcal{R}_i$ to obtain a radial linkage $\mathcal{R}_{i+1}$ of order $4(k+\mathsf{c}t^{12})$ whose endpoints on the outermost cycles of the nests of the $S_v$ coincide with the vertices of the rails of $S_v$ for every vortex segment $S_v$ of $\rho'$ and whose other endpoints coincide with the endpoints of $\mathcal{R}_i$ on $V(\Omega_0)$.
Moreover, the endpoints of the paths among $\mathcal{R}_i$ that lead to $S_v$ appear consecutively on $\Omega_0$.
Hence, we may obtain an extended $k$-surface wall $D_{i+1}$ from $D_i$ by discarding some of the cycles and paths in each of the handle and crosscap segments and integrating the nests of the $S_v$ along the paths in $\mathcal{R}_{i+1}$.
We may also obtain a $\Sigma_i$-decomposition $\delta_{i+1}$ of breadth at most $b$ and depth at most
\begin{align*}
    \mathsf{depth}_{\ref{thm:societyclassification}}(t,4(k+\mathsf{c}t^{12}+1),\mathsf{transaction}(t^2-i,t,k))\leq \mathsf{depth}_{\ref{thm:societyclassification}}(t,4(k+\mathsf{c}t^{12}+1),\mathsf{transaction}(t^2,t,k))
\end{align*}
for $G-A_{i+1}$ where $A_{i+1\coloneqq A_i\cup A}$ with
\begin{align*}
    |A_{i+1}| \leq~& |A_i| + \mathsf{apex}_{\ref{thm:societyclassification}}(t,4(k+\mathsf{c}t^{12}+1),\mathsf{transaction}(t^2,t,k))\\
    \leq~& i\cdot \mathsf{apex}^\mathsf{genus}_{\ref{thm:societyclassification}}(t) + \mathsf{apex}_{\ref{thm:societyclassification}}(t,4(k+\mathsf{c}t^{12}+1),\mathsf{transaction}(t^2,t,k)) + 16t^3.
\end{align*}
It follows, in particular from the choice of the size of $M_0$, that $\Lambda=(A_{i+1},M_0,D_{i+1},\delta_i)$ is a $k$-$(\mathsf{apex}_{\ref{thm:localstructure}}(t^2,t,k),\nicefrac{1}{2}(t-3)(t-4),\mathsf{depth}_{\ref{thm:localstructure}}(t,k),r)$-$\Sigma_i$-layout centred at $M$ as desired.
With this, our proof is complete.
\end{proof}

\subsection{Structure with respect to a mesh}\label{subsec:localstructurewrtmesh}

The $\Sigma$-decomposition, or better the $\Sigma$-layout of $G$ we obtain from \cref{thm:localstructure} is ``\textsl{centred at $M$}'' where $M$ is the mesh we start with.
This property is, at least at first glance, different from the property Robertson and Seymour \cite{RobertsonS2003Graph}, and later Kawarabayashi et al.\ \cite{KawarabayashiTW2021Quickly}, established in their versions of the Local Structure Theorem.
They were looking for a property called ``\textsl{$W$-central}'' where $W$ is a large wall.
In this subsection we establish that the $\Sigma$-layouts provided by \cref{thm:localstructure} also imply the property of being $M$-central, where we lift the definition of $W$-central of Robertson and Seymour to arbitrary meshes instead of restricting ourselves to walls.

For this, recall that, given some $r$-mesh $M$ and a separation $(A,B)$ of order at most $r-1$ we said that $X\in \{ A,B\}$ is the \textsl{majority side} of $(A,B)$ if $X\setminus Y$, where $Y\in\{ A,B\}\setminus\{ X\}$, contains the vertices of both a vertical and a horizontal path of $M$.
Moreover, we observed that the orientation $\mathcal{T}_M$ of $\mathcal{S}_r$ is a tangle.

Let $\Sigma$ be a surface and let $\delta$ be a $\Sigma$-decomposition of a graph $G$ containing an $r$-mesh $M$.
We say that $\delta$ is \emph{$M$-central} if there is no cell $c\in C(\delta)$ such that $V(\sigma(c))$ contains the majority side of a separation from $\mathcal{T}_M$.
Similarly, let $A\subseteq V(G)$, $|A|\leq r-1$, let $\Sigma'$ be a surface and $\delta'$ be a $\Sigma'$-decomposition of $G-A$.
Then we say that $\delta'$ is \emph{$(M-A)$-central} for $G$ if no cell of $\delta'$ contains the majority side of a separation from $\mathcal{T}_M\cap\mathcal{S}_{r-|A|}$.

The property of being $M$-central for a mesh can be stated more abstractly, fully with respect to a tangle and not a mesh.
That is, let $\mathcal{T}$ be a tangle of order $r$ in a graph $G$,let $\Sigma$ be a surface, and $\delta$ be a $\Sigma$-decomposition for $G$.
We say that $\delta$ is \emph{$\mathcal{T}$-central} if no cell of $\delta$ contains the big side of a separation from $\mathcal{T}$.
Given a set $A\subseteq V(G)$, $|A|\leq r-1$, a surface $\Sigma'$ and a $\Sigma'$-decomposition $\delta'$ for $G-A$, we say that $\delta'$ is \emph{$(\mathcal{T}-A)$-central} if no cell of $\delta'$ contains the big side of a separation from $\mathcal{T}\cap\mathcal{S}_{r-|A|}$.

In \cite{KawarabayashiTW2021Quickly} the authors prove the following useful result.
The version of their lemma, namely Lemma 15.3, we state here is slightly enhanced as it also contains the additional information proven in Claim 1 of the proof of Lemma 15.3 in \cite{KawarabayashiTW2021Quickly}.

\begin{proposition}[Kawarabayashi, Thomas, and Wollan \cite{KawarabayashiTW2021Quickly}]\label{lemma:tanglecentrality}
Let $b,d,r\geq 0$ be integers with $r\geq b(2d+1)+6$.
Let $G$ be a graph, $\Sigma$ be a surface and $\delta$ be a $\Sigma$-decomposition for $G$ of breadth at most $b$ and depth at most $d$.
Moreover, let $M$ be an $r$-mesh in $G$ which is flat in $\delta$.
Finally, let $c_1,\dots, c_b\in C(\delta)$ be the vortices of $\delta$ and, for each $i\in[b]$, let $(G_i,\Omega_i)$ be the vortex society of $c_i$.

Then there exists a set $A\subseteq \bigcup_{i\in[b]}V(G_i)$ of size at most $b(2d+1)$ with $|A\cap V(G_i)|\leq 2d+1$ such that, for every tangle $\mathcal{T}$ with $\mathcal{T}_M\subseteq \mathcal{T}$, the restriction of $\delta$ to $G-A$ is $(\mathcal{T}-A)$-central.

Moreover, there exists an algorithm that finds such a set $A$, given $G$, $\delta$, and $M$, in time $\mathbf{poly}(r)|E(G)|$.
\end{proposition}

Now \cref{lemma:tanglecentrality} allows us to strengthen \cref{thm:localstructure} by ensuring that the $\Sigma$-decomposition found is also $M$-central where $M$ is the initial mesh.

\begin{theorem}\label{thm:stronger_localstructure}
There exist functions $\mathsf{apex}_{\ref{thm:stronger_localstructure}},\mathsf{depth}_{\ref{thm:stronger_localstructure}}\colon\mathbb{N}^2\to\mathbb{N}$ and $\mathsf{mesh}_{\ref{thm:stronger_localstructure}}\colon\mathbb{N}^3\to\mathbb{N}$ such that for all non-negative integers $k$, $t\geq 5$, and $r$, every graph $G$ and every $\mathsf{mesh}_{\ref{thm:stronger_localstructure}}(t,r,k)$ mesh $M\subseteq G$ one of the following holds.
\begin{enumerate}
    \item $G$ has a $K_t$-minor controlled by $M$, or
    \item $G$ has a $k$-$(\mathsf{apex}_{\ref{thm:stronger_localstructure}}(t),\nicefrac{1}{2}(t-3)(t-4),\mathsf{depth}_{\ref{thm:stronger_localstructure}}(t,k),r)$-$\Sigma$-layout $\Lambda$ centred at $M$ where $\Sigma$ is a surface of Euler-genus at most $t^2$ and \textsl{$\delta$ is $(M-A)$-central where $A$ is the apex set of $\Lambda$}.
\end{enumerate}
Moreover, it holds that

{\centering
  $ \displaystyle
    \begin{aligned}
        \mathsf{apex}_{\ref{thm:stronger_localstructure}}(t,k),~ \mathsf{depth}_{\ref{thm:stronger_localstructure}}(t,k) \in \mathbf{O}\big((t+k)^{112}\big), \text{ and } \mathsf{mesh}_{\ref{thm:stronger_localstructure}}(t,r,k) \in \mathbf{O}\big((t+k)^{115}r \big) .
    \end{aligned}
  $
\par}

There also exists an algorithm that, given $t,k,r$, a graph $G$ and a mesh $M$ as above as input finds one of these outcomes in time $\mathbf{poly}(t+k+r)|E(G)||V(G)|^2$.
\end{theorem}

\begin{proof}
This strengthening follows immediately from \cref{thm:localstructure} by applying \cref{lemma:tanglecentrality} after making the following observations.

First, if the outcome of \cref{thm:localstructure} is the clique-minor, there is nothing to show.

Now let $\Lambda$ be the $k$-$(\mathsf{apex}_{\ref{thm:localstructure}}(t),\nicefrac{1}{2}(t-3)(t-4),\mathsf{depth}_{\ref{thm:localstructure}}(t,k),r)$-$\Sigma$-layout with mesh $M_0$ and $\Sigma$-decomposition $\delta$.
By definition of $\Sigma$-layouts, we know that $M_0$ is a $\big(\mathsf{apex}_{\ref{thm:localstructure}}(t,k)+\nicefrac{1}{2}(t-3)(t-4)(2\mathsf{depth}_{\ref{thm:societyclassification}}(t,k+\mathsf{c}t^{12},\mathsf{transaction}(t^2,t,k))+1)+11+ r\big)$-submesh of $M$ and thus $\mathcal{T}_{M_0}\subseteq \mathcal{T}_M$.

Finally, let $A_1$ be the apex set of $\Lambda$ and notice that $V(M_0)\cap A_1=\emptyset$.
This implies that the restriction of $\mathcal{T}_{M_0}$ to $G-A_1$ is still a tangle of order $\nicefrac{1}{2}(t-3)(t-4)(2\mathsf{depth}_{\ref{thm:societyclassification}}(t,k+\mathsf{c}t^{12},\mathsf{transaction}(t^2,t,k))+1)+11+r$.

If we now let $A_2$ be the set produced by \cref{lemma:tanglecentrality} we may simply set $A\coloneqq A_1\cup A_2$ and the proof is complete.
\end{proof}

\subsection{Two conjectures of Wollan}\label{subsec:wollasnsconjectures}

In his survey on explicit bounds on the GMST, Wollan states two conjectures regarding the Local Structure Theorem \cite{Wollan2022Explicit}.
The first one asks if a Local Structure Theorem with polynomial bounds is possible while also having a tight bound on the Euler-genus of the surfaces that may occur.
This conjecture reads as follows.

Any large mesh naturally gives rise to a model $\eta$ of a grid.
Let $H$ and $G$ be graphs with $t = |V(H)|$.
We say that a mesh $M$ \emph{grasps} a model $\mu$ of $H$ if for every branch set $\mu(v)$, $v\in V(H)$, there exist distinct indices $r_1,\dots,r_t$ of the rows of $M$ and distinct indices $c_1,\dots,c_t$ of the columns of $M$ such that $\mu(v)$ contains the intersection of the $r_i$th row and $c_i$th column of $M$ of for all $i\in[t]$.

Notice that any minor model that is grasped by a large mesh is also controlled by it.
The reverse, however, is not immediately true.
Kawarabayshi, Thomas, and Wollen proved their variant of the Local Structure Theorem in terms of $K_t$-minors being \textsl{grasped} rather than controlled because they sought a way to avoid the use of tangles in their proofs \cite{Wollan2025Personal}.
Since their theorems were stated using this notion, it is natural that Wollan stated his conjectures also for the notion of a minor being grasped.

\begin{conjecture}[Wollan \cite{Wollan2022Explicit}]\label{con:wollan1}
There exist polynomials $p,q$ which satisfy the following.
Let $r,t\geq 0$ be integers and let $H$ be a graph with $t$ vertices.
Let $R=p(t,r)+q(t)$.
Let $G$ be a graph and $W$ be an $R$-wall in $G$.

Then either $G$ has a model of an $H$-minor grasped by $W$, or there exists a subset $A$ of at most $q(t)$ vertex of $G$, a surface $\Sigma$ in which $H$ does not embed and a $W$-central $\Sigma$-decomposition $\delta$ of $G-A$ of depth and breadth at most $q(t)$ and an $r$-subwall $W'$ of $W$ such that $V(W')\cap A=\emptyset$ and $W'$ is flat in $\delta$.
\end{conjecture}

His second conjecture is a special case of \cref{con:wollan1} for graphs that embed in bounded genus surfaces.
He claims that \cref{con:wollan1} would imply the following.

Let $G$ be a graph drawn without crossing edges in a surface $\Sigma$.
A closed curve $\gamma$ in a surface $\Sigma$ is \emph{genus-reducing} if the surface (or the two surfaces) obtained from $\Sigma$ by cutting along $\gamma$ and capping the (up to two) resulting holes has smaller Euler-genus than $\Sigma$.

The \emph{representativity} of the drawing is the minimum number $k$ such that every genus-reducing curve that intersects the drawing in vertices only must intersect at least $k$ vertices.

\begin{conjecture}[Wollan \cite{Wollan2022Explicit}]\label{con:wollan2}
There exists a polynomial $p$ which satisfies the following.
Let $H$ be a graph on $t\geq 1$ vertices.
Let $\Sigma$ be a surface in which $H$ can be embedded.
If $G$ is a graph drawn in $\Sigma$ with representativity at least $p(t)$, then $G$ contains $H$ as a minor.
\end{conjecture}

In this section we prove a slightly weaker variant of \cref{con:wollan1} where we replace ``grasped'' by ``controlled''.
It turns out that this notion is still strong enough to imply \cref{con:wollan2}.

\begin{proof}[Proof of \Cref{thm:strongest_localstructure}]
\Cref{thm:strongest_localstructure} is essentially a further refinement of \cref{thm:stronger_localstructure} by using \cref{cor:universal-surface-walls} more carefully.

For the sake of slightly better bounds, we look back to the proof of \cref{thm:localstructure}.
Indeed, recall from the proof of \cref{thm:localstructure} that, during the construction phase of the surface, we maintain a $(k+\mathsf{c}t^{12})$-surface wall $D_i$ that has exactly the number of handles and crosscaps has $\Sigma_i$ which is the surface we have found in the $i$th step.
In the proof of \cref{thm:localstructure} we discussed that we can never reach $i=t^2$ because then \cref{cor:universal-surface-walls} would imply that we have found our $K_t$-minor.
However, \cref{cor:universal-surface-walls} actually implies that we can never reach a surface $\Sigma_i$ where $H$ embeds since otherwise we would find an $H$-minor controlled by $D_i$.
This observation is enough the prove our claim.
\end{proof}

Notice that \cref{thm:strongest_localstructure} directly implies our second main theorem, namely \cref{thm:intro_local_structure}.
With this, one major milestone of this paper is achieved.

\paragraph{The proof of \cref{con:wollan2}.}
The following proof of \cref{con:wollan2} was presented to us by Wollan in a discussion of the results of this paper.
He kindly agreed for us to include it here \cite{Wollan2025Personal}.

We will prove the following theorem, where the value of $\mathsf{r}(g,t)$ is defined as further below.

\begin{theorem}\label{thm:wollan2}
For all positive integers $g$ and $t$ the following holds:
Every graph that can be drawn with representativity at least $\mathsf{r}(g,t)\in\mathsf{poly}(t+g)$ in a surface $\Sigma$ of Euler-genus at most $g$ contains, as a minor, every graph $H$ on at most $t$ vertices that can be embedded in $\Sigma$.
\end{theorem}

Let us give a short outline before we introduce some required definitions and results.
The idea is as follows:
Fix $H$ to be a graph on $t$ vertices and let $\mathsf{p}(t)$ (roughly) be the sum of $\mathsf{mesh}_{\ref{thm:strongest_localstructure}}(t,3,t^2)$ $\mathsf{apex}_{\ref{thm:strongest_localstructure}}(t,t^2)$ and $y(t)\cdot t^2\mathsf{depth}_{\ref{thm:strongest_localstructure}}(t,t^2)$, where the value of $y(t)$ will become clear later on.

Let $G$ be a graph that embeds on a surface $\Sigma$ where $H$ minimally embeds such that any embedding of $G$ in $\Sigma$ has representativity at least $\mathsf{p}(t)$.
Due to the high representatitivty, $G$ must have large treewidth and thus, contain a huge mesh which allows for the application of \cref{thm:strongest_localstructure}.
Let $\Sigma'$ be the surface for which \cref{thm:strongest_localstructure} returns a $\Sigma'$-decomposition $\delta$.

Notice that the Euler-genus of $\Sigma'$ must be smaller than the Euler-genus of $\Sigma$ as otherwise we would have that $H$ is a minor of $G$.

Next we may assume that all the non-vortex cells of $\delta$ are empty as otherwise we would find a subdivision of $K_{3,3}$ or $K_5$ that can be removed from $G$ by deleting a small set of vertices.
This is impossible as this would drop the Euler-genus of $G$ by at least $1$.

Finally, we consider any vortex $c_0$ of $\delta$.
Using a known trick for bounded depth vortices allows us to either find many disjoint crosses arranged in series around the boundary of $c_0$, or a small set of vertices (in size bounded by $y(t)$) whose deleting removes all such crosses and allows to fully embed the vortex into a disk.
In the first case, one can find a $K_{3,t^2}$ as a minor of $G$, certifying that the Euler-genus of $G$ is larger than the one of $\Sigma$, a contradiction.
So the second case must hold.
This, however, means that $G$ embeds into $\Sigma'$ after deleting a small set of vertices which is a contradiction to the representativity of $G$ on $\Sigma$.
\smallskip

First, we introduce some necessary properties of embeddings of high representativity.

\begin{proposition}[Demaine, Fomin, Hajiaghayi, and Thilikos \cite{DemaineFHT2005Subexponential}]\label{prop:meshesOnSurfaces}
Let $k$ be a positive integer.
Let $G$ be a graph embedded\footnote{We mean here explicitly \textsl{$2$-cell embedded}.} in a surface $\Sigma$ not homeomorphic to a sphere, with representativity at least $4k$.
Then $G$ contains a $k$-mesh.
\end{proposition}

We shall need a relation between the Euler-genus of a graph $G$ and the Euler-genus of a surface where $G$ can be drawn with high representativity.

\begin{proposition}[Seymour, Thomas \cite{SeymourT1996Uniqueness}]\label{prop:RepresentativeGenus}
Let $k$ be a positive integer.
Let $G$ be a $3$-connected graph with a cross-free drawing $\Gamma$ of representativity $k$ on a surface $\Sigma$ and let $\Gamma'$ be a cross-free drawing of $G$ of representativity at least $3$ in another surface $\Sigma'$.
Suppose $k \geq \{ 320,5\log \mathsf{eg}(\Sigma)\}$, then either there is a homeomorphism from $\Sigma$ to $\Sigma'$ taking $\Gamma$ to $\Gamma'$, or $\mathsf{eg}(\Sigma') \geq \mathsf{eg}(\Sigma) + 10^{-4}k^2$.
\end{proposition}

Notice that \cref{prop:RepresentativeGenus} requires $3$-connectivity.
To deal with this, we import another result, this time due to Robertson and Vitray.

Let $G$ be a $2$-connected graph and let $G_1,G_2\subseteq G$ be two edge-disjoint subgraphs of $G$ such that $G=G_1\cup G_2$ which intersect in exactly two vertices, namely $x$ and $y$.
For both $i\in[2]$ let $G_i'$ the the subgraph of $G$ obtained from $G_i$ by adding an $x$-$y$-path $P^i_{xy}\subseteq G_{3-i}$.
We call the operation defined above a \emph{$2$-split} of $G$.
A subgraph $H$ of $G$ is a \emph{$3$-connected block}\footnote{There are several ways to define $3$-blocks and the like. We use this definition as it is the one used in \cref{prop:representAndConnect1}.} if it can be obtained from $G$ by a sequence of $2$-splits.
Any $3$-connected block of a $2$-connected graph $G$ is either a subdivision of a $3$-connected graph, or a graph consisting of three parallel paths sharing their endpoints.
The later are called the \emph{trivial} $3$-connected blocks of $G$.

Let $G$ be a graph.
A \emph{$3$-connected block} of $G$ is a $3$-connected block of a block of $G$.

\begin{proposition}[Robertson, Vitray \cite{RobertsonV1990Representativity}]\label{prop:representAndConnect1}
Let $G$ be a graph with a cross-free drawing $\Gamma$ of representativity $k\geq 3$ in a surface $\Sigma$.
Then there exists a unique non-trivial $3$-connected block $H$ of $G$ such that the subdrawing $\Gamma'\subseteq \Gamma$ of $H$ is a cross-free drawing of representativity $k$ of $H$ in $\Sigma$.
\end{proposition}

\begin{observation}\label{obs:genusApex}
Let $k$ be a positive integer, $\Sigma$ be a surface that is not the sphere, and $G$ be a graph with a cross-free drawing $\Gamma$ of representativity $k$ in a $\Sigma$.
Then, for every set $S\subseteq V(G)$ of at most $k-1$ vertices, the drawing $\Gamma'\subseteq \Gamma$ of $G-S$ is of representativity at least $k-|A|\geq 1$.
\end{observation}

Combining \cref{prop:RepresentativeGenus}, \cref{prop:representAndConnect1}, and \cref{obs:genusApex} yields the following corollary about deleting small vertex sets from highly representative graphs.

\begin{corollary}\label{cor:preserve:RepAndGenus}
Let $G$ be a graph and $\Sigma$ be a surface of Euler-genus $g\geq 1$.
Let $r \geq 320+5\lceil\log g\rceil + 100k+t$ and let $\Gamma$ be a  cross-free drawing of $G$ on $\Sigma$ with representatitivity $r$.
Then for every set $A\subseteq V(G)$ with $|A|\leq t$, and every cross-free drawing $\Gamma'$ of $G-A$ with representativity at least $3$ on some surface $\Sigma'$ it holds that either
\begin{enumerate}
    \item $\Sigma$ and $\Sigma'$ are homeomorphic, or
    \item $\mathsf{eg}(\Sigma') \geq \mathsf{eg}(\Sigma) + k^2$.
\end{enumerate}
\end{corollary}

Next, we need to know small vertex separators in highly representative drawings are well-behaved with respect to the drawing.

\begin{proposition}[Mohar \cite{Mohar1995Uniqueness}]\label{prop:HighRepSmallSep}
Let $k$ be a positive integer.
Let $G$ be a $2$-connected graph with a cross-free drawing $\Gamma$ of representativity at least $k$ in a surface $\Sigma$.
Let $X\subseteq V(G)$ be a set of at most $k-1$ vertices such that $G-X$ is disconnected but $G-Y$ is connected for all $Y\subsetneq X$.
Then there exist edge-disjoint connected subgraphs $G_1,G_2\subseteq G$ where $G=G_1\cup G_2$ and $V(G_1)\cap V(G_2)=X$ such that
\begin{enumerate}
    \item the subdrawing $\Gamma_1\subseteq \Gamma$ of $G_1$ is a drawing of $G_1$ on $\Sigma$, and
    \item there exists a disk $\Delta\subseteq \Sigma$ whose boundary intersects $\Gamma$ exactly in the vertices of $X$ such that the subdrawing $\Gamma_2\subseteq \Gamma$ of $G_2$ is a drawing of $G_2$ in $\Delta$.
\end{enumerate}
\end{proposition}

\begin{figure}[ht]
    \centering
    \begin{tikzpicture}

        \pgfdeclarelayer{background}
		\pgfdeclarelayer{foreground}
			
		\pgfsetlayers{background,main,foreground}

        \begin{pgfonlayer}{background}
        \node (C) [v:ghost] {{\includegraphics[width=16cm]{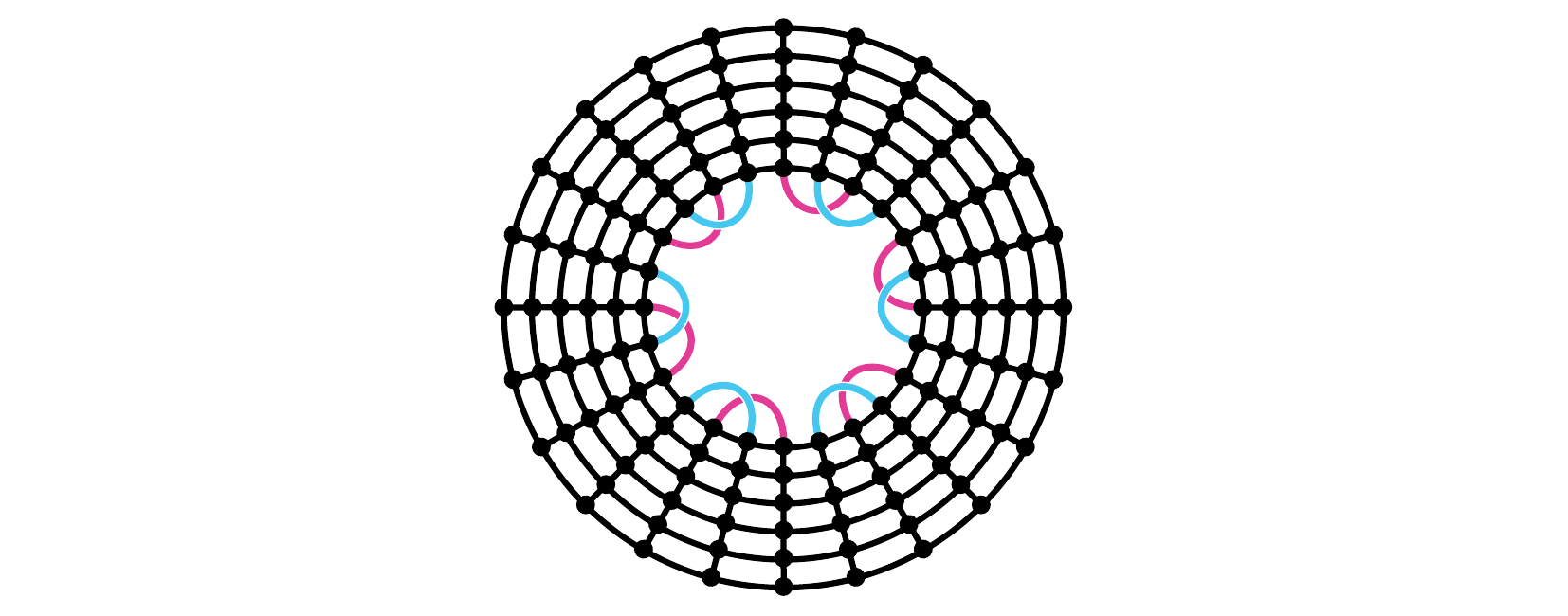}}};
        \end{pgfonlayer}{background}
			
        \begin{pgfonlayer}{main}
        \end{pgfonlayer}{main}

        \begin{pgfonlayer}{foreground}
        \end{pgfonlayer}{foreground}

      \end{tikzpicture}
    \caption{The shallow vortex grid of order $6$.}
    \label{fig:shallowVortex}
\end{figure}

Finally, we need a way to deal with the vortices.
We have already seen a grid-like representation of a vortex in \textsl{(iii)} of \cref{fig:three_ingredients}.
The \emph{shallow vortex grid} of order $k$ is obtained from the $(k,4k)$-cylindrical grid by adding pairwise crossing edges to the inner-most concentric cycle.
See \cref{fig:shallowVortex} for an example.
The following theorem is a combination of a repeated application of Lemma 3.4 and then one application of Lemma 3.7 from \cite{ThilikosW2024Killing}.

\begin{proposition}[Thilikos, Wiederrecht \cite{ThilikosW2024Killing}]\label{prop:KillTheV}
Let $k\leq t$ be positive integers.
There exists a universal constant $c_{\ref{prop:KillTheV}}$ such that, if $(G,\Omega)$ is a society with a cylindrical rendition of depth at most $t$ and a nest of order at least $12k^2+c_{\ref{prop:KillTheV}}$ around the vortex $c_0$, then one of the following holds.
\begin{enumerate}
    \item There exists a separation $(A,B)$ of order at most $12t(k-1)$ with $V(\Omega)\cap B\subseteq A\cap B$, such that if $\Omega'$ is the restriction of $\Omega$ to $A\setminus B$, then $(G[A\setminus B],\Omega')$ has a vortex-free rendition in a disk, or
    \item $G$ contains the shallow vortex grid of order $k$ as a minor.
\end{enumerate}
\end{proposition}

\begin{figure}[ht]
    \centering
    \begin{tikzpicture}

        \pgfdeclarelayer{background}
		\pgfdeclarelayer{foreground}
			
		\pgfsetlayers{background,main,foreground}

        \begin{pgfonlayer}{background}
        \node (C) [v:ghost] {{\includegraphics[width=13cm]{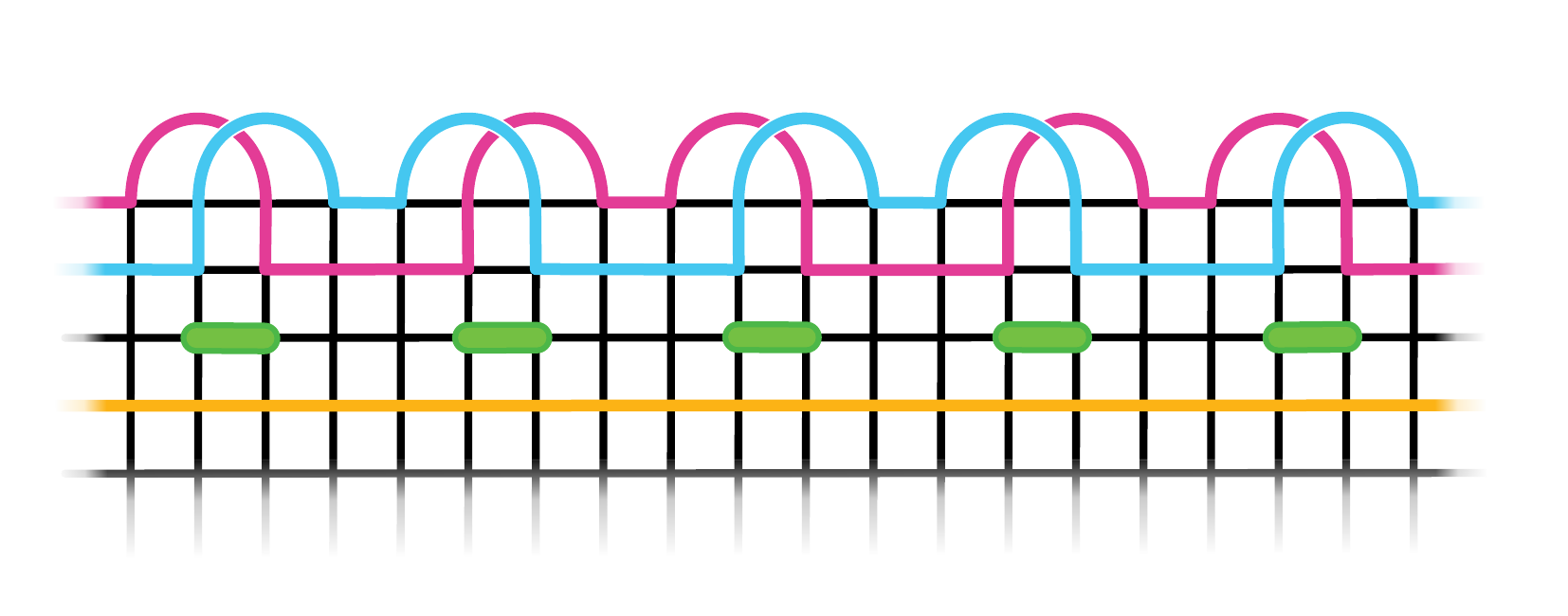}}};
        \end{pgfonlayer}{background}
			
        \begin{pgfonlayer}{main}
        \end{pgfonlayer}{main}

        \begin{pgfonlayer}{foreground}
        \end{pgfonlayer}{foreground}

    \end{tikzpicture}
    \caption{A $K_{3,5}$ as a minor in a shallow vortex grid. The three vertices of high degree are obtained by contracting the \textcolor{HotMagenta}{magenta}, the \textcolor{CornflowerBlue}{blue}, and the \textcolor{BananaYellow}{yellow} path respectively.
    The vertices of the large independent set are obtained by contracting each of the \textcolor{AppleGreen}{green} edges.}
    \label{fig:shallowVortexK3t}
\end{figure}

It is relatively easy to see that large shallow vortex grids contain $K_{3,t}$-minors.
See \cref{fig:shallowVortexK3t} for an illustration.

\begin{observation}\label{obs:K3tinVortex}
Let $k\geq 1$ be an integer.
The shallow vortex grid of order $k$ contains the graph $K_{3,k}$ as a minor.
\end{observation}

The following lemma is folklore and an immediate consequence of Euler's Formula \cite{Euler1758Elementa}.

\begin{proposition}\label{prop:K3tLargeGenus}
For every positive integer $g$, the Euler-genus of $K_{3,2g+3}$ is at least $g+1$.
\end{proposition}

With this, we can formulate another corollary for later use.

\begin{corollary}\label{cor:noVortex}
Let $g$ be a non-negative integer.
Let $G$ be a graph embedded on a surface of Euler genus $g$.
Then $G$ does not contain the shallow vortex grid of order $2g+3$ as a minor.
\end{corollary}

We are now ready to define the function we need for our proof of \cref{con:wollan2}.
Recall the constant $c_{\ref{prop:KillTheV}}$ from \cref{prop:KillTheV}.
\begin{align*}
    \mathsf{x}(g) & \coloneqq 12(2g+3)^2+c_{\ref{prop:KillTheV}}+3\\
    \mathsf{r}(g,t) & \coloneqq  4\cdot \mathsf{mesh}_{\ref{thm:strongest_localstructure}}(t,3,\mathsf{x}(g)) + \mathsf{apex}_{\ref{thm:strongest_localstructure}}(t,\mathsf{x}(g))\\
    & + 6(t-3)(t-4)(2g+2)\cdot\mathsf{depth}_{\ref{thm:strongest_localstructure}}(t,\mathsf{x}(g)) + 100t + 5\lceil\log g\rceil + 320
\end{align*}

\begin{proof}[Proof of \cref{thm:wollan2}]
Fix $g$ and $t$ and let $G$, $\Sigma$ be as in the statement and let $\Gamma$ be the cross-free drawing of $G$.
We may assume that $t\geq 3$.

By \cref{prop:meshesOnSurfaces} and the definition of $\mathsf{r}$ we know that $G$ contains a $\mathsf{mesh}_{\ref{thm:strongest_localstructure}}(t,3,\mathsf{x}(g))$-mesh $M$.

We now apply \cref{thm:strongest_localstructure} to $G$ and $M$.
This either yields a $K_t$-minor certifying that $G$ indeed contains \textsl{all} graphs on $t$ vertices as a minor, or we gain a set $A_1\subseteq V(G)$ of size at most $\mathsf{apex}(t,\mathsf{x}(g))$, a surface $\Sigma'$ and a $\mathsf{x}(g)$-$(\mathsf{apex}_{\ref{thm:strongest_localstructure}}(t,\mathsf{x}(g)),\nicefrac{1}{2}(t-3)(t-4),\mathsf{depth}_{\ref{thm:strongest_localstructure}}(t,\mathsf{x}(g),3)$-$\Sigma'$-layout $\Lambda$ for $G$ such that $K_t$ does not embed in $\Sigma'$.

By definition of layouts, we know that $G-A_1$ contains a $\mathsf{x}(g)$-surface wall representing $\Sigma'$ and therefore, the embedded part of $G-A_1$ in $\Sigma'$ has representativity at least $\mathsf{x}(g) \geq 3$.
Moreover, every vortex of $\Lambda$ has a nest of order $\mathsf{x}(g)\geq 12(2g+3)^2+c_{\ref{prop:KillTheV}}+1$.

By \cref{cor:noVortex} we know that $G$ cannot contain the shallow vortex grid of order $2g+3$ as a minor since $G$ has an embedding in $\Sigma$.
Let $\ell\in[0,\nicefrac{1}{2}(t-3)(t-4)]$ be the number of vortices of $\Lambda$.
So, by our observation on the nest size of each vortex, we may now extract the societies $(G_i,\Omega_i)$ defined by the inner graphs of the outermost cycle of the nests of our vortices and apply \cref{prop:KillTheV} to each of them.
As a result, for each $i\in[\ell]$ we obtain a set $Z_i\subseteq V(G_i)$ of size at most $12(2g+2)\mathsf{depth}_{\ref{thm:strongest_localstructure}}$ such that $(G_i-Z_i,\Omega_i')$ has a cylindrical rendition in the disk.

We may now observe three things.
First, all of the $Z_i$ are fully disjoint from the base cylinder and the handle and crosscaps of the $\mathsf{x}(g)$-surface wall representing $\Sigma'$ in $G-A_1$.
Hence, the embedded part of $G-A_1-\bigcup_{i\in[\ell]}Z_i$ still has representatitivty at least $\mathsf{x}(g)$.
Second, $A_2\coloneqq A_1\cup \bigcup_{i\in[\ell]}Z_i$ has size at most $\mathsf{apex}_{\ref{thm:strongest_localstructure}}(t,\mathsf{x}(g)) + 6(t-3)(t-4)(2g+2)\mathsf{depth}_{\ref{thm:strongest_localstructure}}(t,\mathsf{x}(g)$.
Third, the $\Sigma'$-decomposition of $G-A_1$ in $\Lambda$ can be adjusted to be a vortex-free $\Sigma'$-decomposition $\delta$ of $G-A_2$.

In the next step we prove, that we may assume that no cell of $\delta$ contains a vertex in its interior.
In other words, $\delta$ is an embedding of $G-A_2$ on $\Sigma'$ with representativity at least $\mathsf{x}(g)$.
To see this suppose there is a cell $c$ of $\delta$ such that $|V(\sigma(c))\setminus \widetilde{c}|\geq 1$.
Notice that, due to the previous step, we know that $c$ is not a vortex cell.
Hence, $|\widetilde{c}|\leq 3$.
Now there exists a minimal set $X\subseteq A_2$ such that $(V(\sigma(c))\cup X,V(G)\setminus (V(\sigma(c))\setminus\widetilde{c}))$ is a separation.
Let $Y\coloneqq X\cup\widetilde{c}$.
As $|Y|+3<\mathsf{r}(g,t)$, \cref{prop:HighRepSmallSep} implies that there exists a disk $\Delta\subseteq \Sigma$ such that $G[V(\sigma(c))\cup X]$ is drawn by $\Gamma$ on $\Delta$ with exactly the vertices in $Y$ being drawn on the boundary of $\Delta$.
It follows, that also $\sigma(c)$ has such a drawing and thus, we may replace $c$ with a vortex-free rendition of $\sigma(c)$ in $\Delta$ such that no cell of this rendition has vertices in its interior.
By repeating this step for each cell $c$ of $\delta$ with $V(\sigma(c))\setminus\widetilde{c}\neq\emptyset$ we can ensure that indeed no cell of the resulting $\Sigma'$-decomposition has vertices in its interior.
Hence, $G-A_2$ has a cross-free drawing $\Gamma'$ in $\Sigma'$ which is of representativity at least $\mathsf{x}(g) \geq 3$.

By \cref{cor:preserve:RepAndGenus} and since $\mathsf{r}(g,t)\geq 320 +5\lceil\log g\rceil + 100t + |A_2|$ it follows now that either $\Sigma$ and $\Sigma'$ are homeomorphic, a case we have already excluded, or $\mathsf{eg}(\Sigma') \geq g +t^2$.
In this second case, however, it follows that $K_t$ can be embedded in $\Sigma'$ which is also a contradiction and completes our proof.
\end{proof}

\section{From local to global}\label{sec:localtoglobal}

We are now ready to prove \cref{thm:mainthm_simplest}.
Indeed, as we did for the Local Structure Theorem, we prove a much more powerful result that will directly imply \cref{thm:mainthm_simplest}.

Our proof follows along a well established route for the globalisation of local structure theorems based on walls or tangles.
In fact, we directly take the outline of this proof from Kawarabayashi et al.\ \cite{KawarabayashiTW2021Quickly} who took it from Diestel, Kawarabayashi, Muller, and Wollan \cite{DiestelKMW2012Excluded}.
A precursor of this proof can already be found in Graph Minors X \cite{RobertsonS1991Graph} where Robertson and Seymour proved a similar statement in a hypergraph setting.

Let $k\geq 4$, $w$, $b$, and $a$.
We say that a graph $G$ has an \emph{$(a,b,w)$-near embedding} in a surface $\Sigma$ \emph{of detail $k$} if there exist integers $r$ and $d\leq \nicefrac{1}{2}(w-1)$ such that $G$ has a $k$-$(a,b,d,r)$-$\Sigma$-layout $(A,M,D,\delta)$ where
\begin{enumerate}
    \item all vertices of $G-A$ that do not belong to the interior of a vortex cell of $\delta$ are grounded, and
    \item every vortex of $\delta$ has a linear decomposition of \textsl{width} at most $w$.
\end{enumerate}
Notice that it is possible that $V(G)\subseteq A$ where $A$ is the apex set of the $k$-$(a,b,d,r)$-$\Sigma$-layout as above.
In such a case we will assume $\Sigma$ to be the \emph{empty surface}.
Notably, only a graph without vertices embeds in the empty surface.

Let $G$ be a graph and let $(T, \beta)$ be a tree-decomposition of $G$.
For any given $t \in V(T)$, the \emph{torso} $G_t$ of $G$ at the node $t$ is constructed by taking $G_t' \coloneqq G[\beta(t)]$ and adding, for each pair $u,v\in\beta(t)$, the edge $uv$ if it does not already exist and if there exists $t'\in N_T(t)$ with $u,v\in\beta(t')$.

\begin{theorem}\label{thm:GMST}
There exist functions $\mathsf{adhesion}_{\ref{thm:GMST}},\mathsf{apex}_{\ref{thm:GMST}},\mathsf{vortex}_{\ref{thm:GMST}}\colon\mathbb{N}^2\to\mathbb{N}$ such that for every graph $H$ on $t\geq 1$ vertices, every $k\geq 3$ and every graph $G$ one of the following holds:
\begin{enumerate}
    \item $G$ contains $H$ as a minor, or
    \item there exists a tree-decomposition $(T,\beta)$ for $G$ of adhesion at most $\mathsf{adhesion}_{\ref{thm:GMST}}(t,k)$ such that for every $x\in V(T)$ the torso $G_x$ of $G$ at $x$ has an $(\mathsf{apex}_{\ref{thm:GMST}}(t,k),\nicefrac{1}{2}(t-3)(t-4),\mathsf{vortex}_{\ref{thm:GMST}}(t,k))$-near embedding of detail $k$ into a surface where $H$ does not embed.
\end{enumerate}
Moreover, it holds that

{\centering
  $ \displaystyle
    \begin{aligned}
        \mathsf{adhesion}_{\ref{thm:GMST}}(t,k),\mathsf{apex}_{\ref{thm:GMST}}(t,k),~ \mathsf{vortex}_{\ref{thm:GMST}}(t,k) \in \mathbf{O}\big((t+k)^{2300}\big).
    \end{aligned}
  $
\par}
There also exists an algorithm that, given $k$, $H$, and $G$ as input, finds either an $H$-minor model in $G$ or a tree-decomposition $(T,\beta)$ as above in time $2^{\mathbf{poly}(t+k)}|V(G)|^{3}|E(G)|\log |V(G)|$.
\end{theorem}

In order to prove \cref{thm:GMST} we need one additional tool.
That is, we need a way to algorithmically decide if a given set $X\subseteq V(G)$ of bounded size has a $\nicefrac{2}{3}$-balanced separator of a given size.

\begin{proposition}[Reed \cite{Reed1992Finding}]\label{prop:balancedseps}
There exists an algorithm that takes as input an integer $k$, a graph $G$, and a set $X\subseteq V(G)$ of size at most $3k+1$ and finds, in time $2^{\mathbf{O}(k)}|E(G)|$, either a $\nicefrac{2}{3}$-balanced separator of size at most $k$ for $X$ or correctly determines that $X$ is $(k,\nicefrac{2}{3})$-linked in $G$.
\end{proposition}

Our proof of \cref{thm:GMST} will be an inductive proof of a slightly stronger statement as follows.

Let $G$ be a graph.
A tuple $(T,r,\beta)$ is a \emph{rooted} tree-decomposition for $G$ if $(T,\beta)$ is a tree-decomposition for $G$ and $r\in V(T)$.
Notions like width and adhesion naturally transfer from tree-decompositions to rooted tree-decompositions.

\begin{theorem}\label{thm:GMST_induction}
There exist functions $\mathsf{link}_{\ref{thm:GMST_induction}},\mathsf{apex}_{\ref{thm:GMST_induction}},\mathsf{vortex}_{\ref{thm:GMST_induction}}\colon\mathbb{N}^2\to\mathbb{N}$ such that for every graph $H$ on $t\geq 1$ vertices, every $k\geq 3$, every graph $G$, and every vertex set $X\subseteq V(G)$ with $|X|\leq 3\mathsf{link}_{\ref{thm:GMST_induction}}(t,k)+1$ one of the following holds:
\begin{enumerate}
    \item $G$ contains $H$ as a minor, or
    \item there exists a rooted tree-decomposition $(T,r,\beta)$ for $G$ of adhesion at most 
    \begin{align*}
    3\mathsf{link}_{\ref{thm:GMST_induction}}(t,k)+\mathsf{apex}_{\ref{thm:GMST_induction}}(t,k)+\mathsf{vortex}_{\ref{thm:GMST_induction}}(t,k)+3
    \end{align*}
    such that 
    \begin{enumerate}
        \item for every $x\in V(T)$ the torso $G_x$ of $G$ at $x$ has a $(4\mathsf{link}_{\ref{thm:GMST_induction}}(t,k)+\mathsf{apex}_{\ref{thm:GMST_induction}}(t,k),\nicefrac{1}{2}(t-3)(t-4),2\mathsf{link}_{\ref{thm:GMST_induction}}(t,k)+\mathsf{vortex}_{\ref{thm:GMST_induction}}(t,k))$-near embedding of detail $k$ into a surface where $H$ does not embed, and
        \item let $A_r$ be the apex set for the torso $G_r$ of $G$ at $r$, then $X\subseteq A_r$.
    \end{enumerate}
\end{enumerate}
Moreover, it holds that

{\centering
  $ \displaystyle
    \begin{aligned}
        \mathsf{apex}_{\ref{thm:GMST_induction}}(t,k),~ \mathsf{vortex}_{\ref{thm:GMST_induction}}(t,k) \in \mathbf{O}\big((t+k)^{115}\big)\text{ and }\mathsf{link}_{\ref{thm:GMST_induction}}(t,k)\in\mathbf{O}\big( (t+k)^{2300}\big).
    \end{aligned}
  $
\par}
There also exists an algorithm that, given $k$, $H$, and $G$ as input, finds either an $H$-minor model in $G$ or a tree-decomposition $(T,\beta)$ as above in time \(2^{\mathbf{poly}(t+k)}|V(G)|^3|E(G)|\log|V(G)|\).
\end{theorem}

\begin{proof}
Let us start by introducing the functions involved.
\begin{align*}
    \mathsf{apex}_{\ref{thm:GMST_induction}}(t,k) \coloneqq~& \mathsf{apex}_{\ref{thm:strongest_localstructure}}(t,k)\\
    \mathsf{vortex}_{\ref{thm:GMST_induction}}(t,k) \coloneqq~& 2\mathsf{depth}_{\ref{thm:strongest_localstructure}}(t,k)+1\\
    \mathsf{link}_{\ref{thm:GMST_induction}}(t, k)\coloneqq~& \mathsf{c}_1\mathsf{mesh}_{\ref{thm:strongest_localstructure}}(t,k,k)^{20}
\end{align*}
where $\mathsf{c}_1$ is the first constant from \cref{thm:algogrid}. Note that \(\mathsf{mesh}_{\ref{thm:strongest_localstructure}}(t,k,k) \in \mathbf{O}((t+k)^{115}+t^3r)\), so \(\mathsf{link}_{\ref{thm:GMST_induction}}(t, k) \in \mathbf{O}((t+k)^{2300})\) as required.

We proceed by induction on $|V(G)\setminus X|$.

First observe that in case $|V(G)|\leq 3\mathsf{link}_{\ref{thm:GMST_induction}}(t,k)+1$ we may simply define $T$ to be the tree on one vertex, say $r$, and set $\beta(r)\coloneqq V(G)$.
This satisfies the second outcome of the assertion.
Hence, we may assume $|V(G)|\geq 3\mathsf{link}_{\ref{thm:GMST_induction}}(t,k)+2$.
\smallskip

Next, suppose $|X|\leq 3\mathsf{link}_{\ref{thm:GMST_induction}}(t,k)$.
In this case there exists some vertex $v\in V(G)\setminus X$ and by setting $X'\coloneqq X\cup \{ v\}$ we obtain that $|V(G)\setminus X'| < |V(G)\setminus X|$.
If we apply the induction hypothesis to $G$ and $X'$, we either find $H$ as a minor in $G$ or the resulting rooted tree-decomposition also satisfies the second outcome of the assertion for $X$.
Hence, we may assume that $|X|=3\mathsf{link}_{\ref{thm:GMST_induction}}(t,k)+1$.
\smallskip

Let us apply \cref{prop:balancedseps} to $X$ in $G$. Note that it takes \(2^{\mathbf{poly}(t+k)}|E(G)|\) time to compute.
Then one of two outcomes occur.
\begin{enumerate}
    \item there exists a $\nicefrac{2}{3}$-balanced separator $S$ for $X$ in $G$ with $|S|\leq \mathsf{link}_{\ref{thm:GMST_induction}}(t,k)$, or
    \item $X$ is $\big(\mathsf{link}_{\ref{thm:GMST_induction}}(t,k),\nicefrac{2}{3}\big)$-linked in $G$.
\end{enumerate}
\medskip

\textbf{Case 1:} There exists a $\nicefrac{2}{3}$-balanced separator $S$ for $X$ in $G$ with $|S|\leq \mathsf{link}_{\ref{thm:GMST_induction}}(t,k)$.
\smallskip

Let $G_1',\dots, G_{\ell}'$ be the components of $G-S$ and for each $i\in[\ell]$ let $G_i\coloneqq G[V(G_i')\cup S]$.
Moreover, for each $i\in[\ell]$ let $X_i'\coloneqq (V(G_i)\cap X)\cup S$.

It follows that
\begin{align*}
    |X_i'|\leq~& \lfloor\nicefrac{2}{3}(3\mathsf{link}_{\ref{thm:GMST_induction}}(t,k)+1)\rfloor + \mathsf{link}_{\ref{thm:GMST_induction}}(t,k)\\
    \leq~& 3\mathsf{link}_{\ref{thm:GMST_induction}}(t,k).
\end{align*}
We now describe how to define the rooted tree-decomposition $(T,r,\beta)$ in case we do not find an $H$-minor in $G$.
The definition of both $T$ and $\beta$ are done iteratively.
First we simply introduce the node $r$ and set $\beta(r)\coloneqq X\cup S$.
Notice that we may include all of $X\cup S$ into the apex set for the near embedding for the torso of $G$ at the node $r$.
This is true since
\begin{align*}
    |X\cup S| \leq~& |X| + |S| \leq 3\mathsf{link}_{\ref{thm:GMST_induction}}(t,k) + 1 + \mathsf{link}_{\ref{thm:GMST_induction}}(t,k)\\
    \leq~ & 4\mathsf{link}_{\ref{thm:GMST_induction}}(t,k)+1.
\end{align*}

For each $i\in[\ell]$ where $V(G_i)=X_i'$ we introduce a node $x_i$, make it adjacent to $r$ and set $\beta(x_i)\coloneqq X_i'$.
What remains are those $i\in[\ell]$ where there exists $v\in V(G_i)\setminus X_i'$.
For each such $i\in[\ell]$ select $v_i\in V(G_i)\setminus X_i'$ arbitrarily and set $X_i\coloneqq X_i'\cup\{ v_i\}$.
We now have that $|V(G_i)\setminus X_i|<|V(G)\setminus X|$ and thus, by the induction hypothesis, $G_i$ either contains $H$ as a minor and we satisfy the first outcome of the assertion, or there exists a rooted tree-decomposition $(T_i,r_i,\beta_i)$ for $G_i$ that satisfies the second outcome of our assertion.
If the first case occurs for any $i\in[\ell]$ we are done, so we may assume that the second case holds.
For each $i\in[\ell]$ where $v_i$ was defined we add the tree $T_i$ to our tree $T$, make $r_i$ adjacent to $r$, and set $\beta(x)\coloneqq \beta_i(x)$ for each $x\in V(T_i)$.

It is straight forward that $(T,r,\beta)$ constructed this way is indeed a rooted tree-decomposition satisfying the requirements of the assertion.
\medskip

\textbf{Case 2:} $X$ is $\big(\mathsf{link}_{\ref{thm:GMST_induction}}(t,k),\nicefrac{2}{3}\big)$-linked in $G$.
\smallskip

In this case we may call upon \cref{thm:algogrid} to find a $\mathsf{mesh}_{\ref{thm:strongest_localstructure}}(t,k,k)$-mesh $M$ in $G$ such that the tangle $\mathcal{T}_M$ is a truncation of the tangle $\mathcal{T}_X$. Note that this takes \(2^{\mathbf{poly}(t+k)}|V(G)|^2|E(G)|\log |V(G)|\) time.
We then apply \cref{thm:strongest_localstructure} to $M$. This takes \(2^{\mathbf{poly}(t+k)}|V(G)|^2|E(G)|\) time.
This either yields an $H$-minor model in $G$ controlled by $M$, and therefore by $\mathcal{T}_X$, or there exists some surface $\Sigma$ where $H$ does not embed and we find a $k$-$(\mathsf{apex}_{\ref{thm:strongest_localstructure}}(t),\nicefrac{1}{2}(t-3)(t-4),\mathsf{depth}_{\ref{thm:strongest_localstructure}}(t, k),k)$-$\Sigma$-layout with landscape $\Lambda=(A,M',D,\delta)$ centred at $M$ where $\delta$ is $(M-A)$-central.

We define $A_r\coloneqq A\cup X$ and observe that
\begin{align*}
    |A_r| \leq |A|+|X| \leq 3\mathsf{link}_{\ref{thm:GMST_induction}}(t,k) + \mathsf{apex}_{\ref{thm:strongest_localstructure}}(t,k) +1.
\end{align*}

Let $u_1,\dots,u_{q}$, $q\leq \nicefrac{1}{2}(t-3)(t-4)$, be the vortices of $\delta$.
By \cref{prop:depth_to_lin_decomp} we may find, for each $i\in[q]$, a linear decomposition $(Y^i_1,Y^i_2,\dots,Y^i_{n_i})$ of adhesion at most $\mathsf{depth}_{\ref{thm:strongest_localstructure}}(t,k)$. This takes \(\mathbf{poly}(t+k)|E(G)|\log |V(G)|\) time.
Since $\delta$ is $(M-A)$-central, it is also $(\mathcal{T}_M-A)$-central.
Moreover, since $\mathcal{T}_M$ is a truncation of $\mathcal{T}_X$ it follows that $\delta$ is also $(\mathcal{T}_X-A)$-central.
Hence,
\begin{align*}
|Y^i_j\cap X| \leq \lfloor\nicefrac{2}{3}|X|\rfloor \leq 2\mathsf{link}_{\ref{thm:GMST_induction}}(t,k) 
\end{align*}
for all $i\in[q]$ and all $j\in[n_i]$.

Let $U$ be the set of all ground vertices of $\delta$ together with all vertices that appear in sets $Y^i_j\cap Y^i_{j+1}$ for some $i\in[q]$ and $j\in[n_i-1]$.

As in the first case, we will construct the rooted tree-decomposition $(T,r,\beta)$ by iteratively defining the tree $T$ and the mapping $\beta$.
We begin by introducing the vertex $r$ together with $\beta(r)\coloneqq U\cup A_r$.
With the discussion above, it follows that
\begin{align*}
    |Y^i_j \cap (U\cup A_r)| \leq~& 2\mathsf{link}_{\ref{thm:GMST_induction}}(t,k) + 2\mathsf{depth}_{\ref{thm:strongest_localstructure}}(t,k)+1\\
    \leq~& 2\mathsf{link}_{\ref{thm:GMST_induction}}(t,k)+\mathsf{vortex}_{\ref{thm:GMST_induction}}(t,k)).
\end{align*}

For each $i\in[q]$ and $j\in[n_i]$ let $F_{i,j}\coloneqq G[A\cup Y^i_j]$.
Then
\begin{align*}
    |V(F_{i,j})\cap (U\cup A_r)| \leq~& 2\mathsf{link}_{\ref{thm:GMST_induction}}(t,k)+\mathsf{vortex}_{\ref{thm:GMST_induction}}(t,k))+\mathsf{apex}_{\ref{thm:GMST_induction}}(t,k))\\
    \leq~& 3\mathsf{link}_{\ref{thm:GMST_induction}}(t,k)
\end{align*}
since $\mathsf{mesh}_{\ref{thm:strongest_localstructure}}(t,k)\geq \mathsf{vortex}_{\ref{thm:GMST_induction}}(t,k))+\mathsf{apex}_{\ref{thm:GMST_induction}}(t,k))$.
So we may define $X_{i,j}'\coloneqq V(F_{i,j})\cap (U\cup A_r)$.
As before in Case 1 there are two cases to be considered.
If $V(F_{i,j})=X_{i,j}'$ then we introduce a vertex $r_{i,j}$ adjacent to $r$ and set $\beta(r_{i,j})\coloneqq V(F_{i,j})$.
Otherwise there exists some $v_{i,j}\in V(F_{i,j})\setminus X_{i,j}'$ and we set $X_{i,j}\coloneqq X_{i,j}'\cup\{ v_{i,j}\}$.
Now we have
\begin{align*}
    |V(F_{i,j})\setminus X_{i,j}| < |V(G)\setminus X|
\end{align*}
for all $i\in[q]$ and $j\in[n_i]$.
So the induction hypothesis produces either an $H$-minor, in which case we are done, or a rooted tree-decomposition $(T_{i,j},r_{i,j},\beta_{i,j})$ satisfying the second outcome of the assertion.
So we may include each of the trees $T_{i,j}$ into our tree $T$, make $r_{i,j}$ adjacent to $r$, and set $\beta(x)\coloneqq \beta_{i,j}(x)$ for each $x\in V(T_{i,j})$.
This takes care of the vortices.

What is left is to discuss the graphs $J_c'\coloneqq \sigma(c)$ for the non-vortex cells of $\delta$.
Recall that $\delta$ is $(\mathcal{T}_X-A)$-central.
So $J_c'$ cannot contain more than $2\mathsf{link}_{\ref{thm:GMST_induction}}(t,k)$ vertices of $X$ for any cell $c\in C(\delta)$ that is not a vortex.
Hence, we may set $J_c\coloneqq G[A\cup V(J_c')]$ and $X'_c\coloneqq V(J_c)\cap A_r$.
It follows that 
\begin{align*}
|X_c'|\leq~& 2\mathsf{link}_{\ref{thm:GMST_induction}}(t,k) + \mathsf{apex}_{\ref{thm:GMST_induction}}(t,k)\\
\leq~& 3\mathsf{link}_{\ref{thm:GMST_induction}}(t,k). 
\end{align*}
For each non-vortex cell $c$ where $V(J_c)=X_c'$ we introduce a vertex $r_c$ adjacent to $r$ and set $\beta(r_c)\coloneqq X_c'$.
So we are left with the non-vortex cells $c$ where we can find a vertex $v_c\in V(J_c)\setminus X_c'$.
Here we set $X_c\coloneqq X_c'\cup \{ v_c\}$ and obtain
\begin{align*}
|V(J_c)\setminus X_c| < |V(G)\setminus X|
\end{align*}
similar to the cases before.
From the induction hypothesis we now wither obtain an $H$-minor of $G$ or a rooted tree-decomposition $(T_c,r_c,\beta_c)$ for each $J_c$ that satisfies the requirements of the second outcome of the assertion.
As before we introduce the trees $T_c$ to our tree $T$, make $r_c$ adjacent to $r$, and set $\beta(x)\coloneqq \beta_c(x)$ for each $x\in V(T_c)$.
With this, our construction is finally complete.
Moreover, it follows from the construction that the final $(T,r,\beta)$ is indeed a rooted tree-decomposition with the properties as required.

The inductive proof gives rise to a recursive algorithm. In total, it takes \(O(|V(G)|)\) recursive calls to find the tree-decomposition of the entire graph \(G\), and each recursive call takes \(2^{\mathbf{poly}(t+k)}|V(G)|^2|E(G)|\log |V(G)|\) time, so the total running time is \(2^{\mathbf{poly}(t+k)}|V(G)|^3|E(G)|\log|V(G)|\).
\end{proof}

\section{Computing the Graph Minor Structure Decomposition in polynomial time}\label{sec:polytime}

The running time of the algorithm given in \cref{thm:GMST} is polynomial in $|V(G)|$ but exponential in $t$ and $k$.
In this section we discuss how this algorithm can be adjusted to have a total running time of $(t+k+|V(G)|)^{\mathbf{O}(1)}$ while still guaranteeing polynomial (in $t$ and $k$) bounds on all parameters.

The overall strategy to achieve this has two pieces.
First we observe that every single step after finding a large wall can already be realised in time $(t+k+|V(G)|)^{\mathbf{O}(1)}$.
Second, we employ techniques from different places within the algorithmic theory of graph minors to find a polynomial-time approximation algorithm that, given a set $X$, either finds a $\nicefrac{2}{3}$-balanced separator of size $\mathbf{O}(k\log k)$ for $X$ or determines that $X$ is $(k,\nicefrac{2}{3})$-linked.
In the later case, this algorithm then proceed with finding a large wall whose tangle agrees with the tangle of $X$.

Towards finding this algorithm from the second piece of our strategy, we need the following approximation result of Amir \cite{Amir2010Approximation}.

\begin{proposition}[Amir \cite{Amir2010Approximation}]\label{prop:approximateBalancedSet}
There exist a positive constant $c_{\ref{prop:approximateBalancedSet}}$ and an algorithm that takes as input a graph $G$, a set $X\subseteq V(G)$, and an integer $k\geq 1$, and either finds a $\nicefrac{2}{3}$-balanced separator $S\subseteq V(G)$ for $X$ with $|S|\leq c_{\ref{prop:approximateBalancedSet}}k\log k$ or correctly determines that $X$ is $(k,\nicefrac{2}{3})$-linked in time $\mathbf{O}(k^2\log k\cdot |V(G)|^2|E(G)|)$.
\end{proposition}

We will also need a procedure that, given a graph $G$ and an integer $k$, either finds a $k$-wall in $G$ or correctly determines that the treewidth of $G$ is smaller than a given bound polynomial in $k$.

\begin{proposition}[Chekuri, Chuzhoy \cite{ChekuriC2016Polynomial}]\label{prop:approxGrid}
There exist positive constants $a_{\ref{prop:approxGrid}}$ and $b_{\ref{prop:approxGrid}}$ and a randomised algorithm that takes as input a graph $G$ and a positive integer $k$ and either outputs a $k$-wall in $G$ or determines that $\mathsf{tw}(G) \leq a_{\ref{prop:approxGrid}}k^{98}(\log k)^{^{b_{\ref{prop:approxGrid}}}}$ with high probability in time $(|V(G)|+k)^{\mathbf{O}(1)}$.
\end{proposition}

So far, we know how to -- approximately -- check for a given set $S$ if it is highly linked or not.
Moreover, if $S$ is highly linked we know that $G$ must have large treewidth and thus, \cref{prop:approxGrid} finds a large wall in $G$.
However, what we actually need is a version of \cref{thm:algogrid} that will find a large wall whose tangle agrees with the tangle of $S$.

A closer look at the proof of \cref{prop:approxGrid} reveals the following.
The construction of the wall found by Chekuri and Chuzhoy operates in several phases.

\paragraph{Well-linked sets.}
First, they employ a strategy that either finds that the input graph $G$ has small (in $k$) treewidth, or it finds a highly (in $k$) linked set.
This set is then passed through several refinement steps (see subsections 2.1 to 2.8 in \cite{ChekuriC2016Polynomial}) to ultimately obtain a large (in $k$) set that is ``node-well-linked'' (see Theorem 2.21 in \cite{ChekuriC2016Polynomial}).
This set is obtained from a highly-linked set $S$ in a way that preserves its tangle.

The definition of node-well-linkedness of Chekuri and Chuzhoy is the same as the definition of \textsl{well-linkedness} used by Kawarabayashi et al.\ from \cite{KawarabayashiTW2021Quickly}.
They define such sets as follows.
Let $G$ be a graph and $S\subseteq V(G)$ be a set of vertices.
We say that $S$ is \emph{well-linked} if for every partition $S_1,S_2$ of $S$ into two sets there does not exist a separation $(A_1,A_2)$ of $G$ such that $A_i\subseteq A_i$ for both $i\in[2]$ and $|A_1\cap A_2|<\min\{|S_1|,|S_2|\}$.

Let us briefly show that the notions of highly linked sets and well-linked sets do really translate into each other while preserving their tangles.

Let $\alpha\in[\nicefrac{2}{3},1)$ be fixed and $k$ be a positive integer.
Furthermore, let $S$ be a $(k,\alpha)$-linked set $S$ in a graph $G$.
A set $F\subseteq V(G)$ is said to be \emph{$S$-free} in $G$ if for every separation $(A_1,A_2)$ of order at most $|F|$ it holds that, if $F\subseteq A_i$ for some $i\in[2]$, then $|A_i\cap S|>\alpha|S|$.

\begin{proposition}[Thilikos and Wiederrecht \cite{ThilikosW2024Excluding}]\label{prop:freesets}
Let $k$ be a positive integer, $\alpha\in[\nicefrac{2}{3},1)$, $G$ be a graph, and $S\subseteq V(G)$ be a $(k,\alpha)$-linked set.
Then there exists an $S$-free set $F$ of size $k-1$ in $G$.
Moreover, there exists an algorithm that, given $G$ and $S$, computes $F$ in time $\mathbf{O}(k|V(G)|^2)$.
\end{proposition}

So, given a highly linked set $S$, we are able to find an $S$-free set $F$ of large size in polynomial time.
Next, we need a specific property of $S$-free sets.

We may now use the fact that every highly linked set $X$ induces a tangle $\mathcal{T}_X$ to rephrase Lemma 14.2 from \cite{KawarabayashiTW2021Quickly} in terms of highly linked sets instead of tangles.

\begin{proposition}[Kawarabayashi, Thomas, and Wollan \cite{KawarabayashiTW2021Quickly}]\label{prop:freeIsWell}
Let $\alpha\in[\nicefrac{2}{3},1)$, $k\geq 1$ be an integer, $G$ be a graph and $S\subseteq V(G)$ be a $(k,\alpha)$-linked set.
If $F\subseteq V(G)$ is $S$-free with $|F|<k$ then $F$ is well-linked.
\end{proposition}

Next, we need to assert that a $S$-free set $F$ of large enough size are enough to define tangles themselves.
This step is necessary because in what follows we will start to edit the graph $G$ by iteratively deleting edges.
While this process will be able to maintain the well-linkedness of $F$ but not necessarily the fact that $S$ is highly linked.

\begin{lemma}\label{lemma:welllinkedIsHighlyLinked}
Let $k$ be a positive integer and $F$ be a well-linked set of size $3k$ in a graph $G$.
Then $\mathcal{T}_F\coloneqq \{ (A,B)\in\mathcal{S}_k \colon |B\cap F|>2k \}$ is a tangle of order $k$ in $G$.

Moreover, if there exists $\alpha\in [\nicefrac{2}{3},1)$ and a $(3k+1,\alpha)$-linked set $S$ such that $F$ is $S$-free, then $\mathcal{T}_F\subseteq \mathcal{T}_S$.
\end{lemma}

\begin{proof}
We begin by showing that for every separation $(A,B)\in\mathcal{S}_k$ one of $(A,B)$ or $(B,A)$ belongs to $\mathcal{T}_F$.
This follows directly from the fact that $|A\cap B|<k$ and $|F|=3k$ and the definition of well-linked sets as follows:
We have that $\min\{ |A\cap F|,|B\cap F| \}<k$ as otherwise $|A\cap B|\geq k$ or we would have a contradiction to $F$ being well-linked.
Moreover, with $|F|=3k$ we have that $\max\{|A\cap F|,|B\cap F|\}>2k$.

So every separation $(A,B)\in\mathcal{S}_k$ is oriented by $\mathcal{T}_F$.
It remains to show that three small sides cannot cover all of $G$.
Let $(A_1,B_1),(A_2,B_2),(A_3,B_3)\in\mathcal{T}_F$.
Then, for each $i\in[3]$, we have $|A_i\cap F|\leq k-1$ and thus, $|(A_1\cup A_2\cup A_3)\cap F|\leq 3k-3<3k$ and $\mathcal{T}_F$ is ineed a tangle.

At last we need to show that, if there exists $\alpha\in [\nicefrac{2}{3},1)$ and a $(3k+1,\alpha)$-linked set $S$ such that $F$ is $S$-free, then $\mathcal{T}_F\subseteq \mathcal{T}_S$.
Suppose there exists $(A,B)\in\mathcal{T}_F$ such that $(B,A)\in\mathcal{T}_S$.
From the definition of $\mathcal{T}_F$ we have that $|A\cap B|<k$ and $|B\cap F|>2k$ and therefore $|A\cap F|\leq k-1$.
It follows that $(A,B\cup F)$ is a separation of order less than $2k$.
As $F\subseteq B\cup F$, the definition of $S$-free sets implies that $|(B\cup F)\cap S|>\alpha|S|$ which, in turn, implies that $(A,B\cup F)\in\mathcal{T}_S$.
However, as $(B,A)\in\mathcal{T}_S$ by assumption, we now have that $G[B]\cup B[A]=G$ which contradicts the fact that $\mathcal{T}_S$ is a tangle.
Hence, $\mathcal{T}_F\subseteq \mathcal{T}_S$ as desired.
\end{proof}

\paragraph{Degree reduction.}
While the refinements mentioned above are possible in general graphs, the actual methods employed by Chekuri and Chuzhoy require the graph to have bounded maximum degree.
This is achieved by passing the graph $G$ and its highly linked set through a degree-reducing procedure that preserves the well-linkedness property of a large subset of the vertices.
This procedure was introduced by Chekuri and Ene in \cite{ChekuriE2013Polylogarithmic} and is based on the so-called \textsl{Cut-Matching Game} of Khandekar, Rao, and Vazirani \cite{KhandekarRV2009Graph} using bounds by Orecchia, Schulman, Vazirani, and Vishnoi \cite{OrecchiaSVV2008Partitioning}.
This game starts by removing all edges of the graph and then iteratively reintroduces edges.
It can be shown (see \cite{ChekuriE2013Polylogarithmic} and \cite{ChekuriC2016Polynomial}) that this game can be used to efficiently (via a randomised algorithm) reduce the maximum degree of an input graph to $\mathbf{O}(\log^3 k)$ while dropping the treewidth by only a poly-logarithmic factor.
Indeed, if one starts out with a well-linked set $Z$ of $k$ vertices in a graph $G$, one ends up with a well linked set $Z'\subseteq Z$ of at least $\Omega(\nicefrac{k}{\operatorname{poly log}k})$ vertices in a graph $G'\subseteq G$ of maximum degree in $\mathbf{O}(\log^3 k)$ such that $Z'$ is well-linked in $Z$.

\paragraph{Path-of-sets system.}
From here the strategy of Chekuri and Chuzhoy is to make use of their well-linked, let us call it $S$, set and grow a structure they call a ``\textsl{path-of-sets system}''.
They show that such a path-of-sets system of large enough \textsl{order} always contains a large wall (see Corollary 3.4 in \cite{ChekuriC2016Polynomial}).
Such a path-of-sets system consists of a sequence $X_1,\dots,X_n$ of well-linked sets, all obtained from $S$ such that there are huge linkages between $X_i$ and $X_{i+1}$ for all $i\in[n-1]$.
As a consequence, the wall $W$ they find passes through $S$ in a way such that there is at least one row $P$ of $W$ that contains a vertex of $S$ in each intersection of $P$ with any column of $W$.

\begin{observation}\label{obs:wellLinkedMesh}
Let $k\geq 3$ be an integer and let $M$ be a $k$-mesh with horizontal paths $P_1,\dots,P_k$ and vertical paths $Q_1,\dots,Q_k$ in a graph $G$.
Then any set $X$ of $k$ vertices such that there there exists $i\in[k]$ and $X$ contains a vertex of $P_i\cap Q_j$ for all $j\in[k]$ is well-linked in $G$.
\end{observation}

It follows that, if $W$ is the wall found within a path-of-sets system grown from a well-linked set $S$, then $\mathcal{T}_W \subseteq \mathcal{T}_S$.
Chekuri and Chuzhoy provide an efficient randomised algorithm that takes as input a graph of bounded maximum degree and a well-linked set $S$ and then finds a path-of-sets system as above (see Theorem 6.1 in \cite{ChekuriKS2005Multicommodity}).
\medskip

From the discussion above, we may deduce the following strengthening of the polynomial Grid Theorem of Chekuri and Chuzhoy.

\begin{proposition}[Chekuri, Chuzhoy \cite{ChekuriC2016Polynomial}]\label{prop:approxAlignedGrid}
There exist positive constants $a_{\ref{prop:approxAlignedGrid}}$ and $b_{\ref{prop:approxAlignedGrid}}$ and a randomised algorithm that takes as input a graph $G$, a positive integer $k$, and a $(a_{\ref{prop:approxAlignedGrid}}k^{98}(\log k)^{^{b_{\ref{prop:approxAlignedGrid}}}},\nicefrac{2}{3})$-linked set and outputs a $k$-wall $W$ in $G$ such that $\mathcal{T}_W\subseteq \mathcal{T}_S$ with high probability in time $(|V(G)|+k)^{\mathbf{O}(1)}$.
\end{proposition}

With this, we are now able to provide a variant of \cref{thm:GMST} with an algorithm that runs in time polynomial in $t+k+|V(G)|$ and computes either a $K_t$-minor or the desired tree-decomposition where all involved parameters (i.e. adhesion, Euler-genus, number of apices, and number and depth of vortices) are polynomials only depending on $t$ (and $k$ in the case of the breadth of the vortices).

\begin{theorem}\label{thm:GMSTpolytime}
There exist functions $\mathsf{adhesion}_{\ref{thm:GMSTpolytime}},\mathsf{apex}_{\ref{thm:GMSTpolytime}},\mathsf{vortex}_{\ref{thm:GMSTpolytime}}\colon\mathbb{N}^2\to\mathbb{N}$ such that for every graph $H$ on $t\geq 1$ vertices, every $k\geq 3$ and every graph $G$ one of the following holds:
\begin{enumerate}
    \item $G$ contains $H$ as a minor, or
    \item there exists a tree-decomposition $(T,\beta)$ for $G$ of adhesion at most $\mathsf{adhesion}_{\ref{thm:GMSTpolytime}}(t,k)$ such that for every $x\in V(T)$ the torso $G_x$ of $G$ at $x$ has an $(\mathsf{apex}_{\ref{thm:GMSTpolytime}}(t,k),\nicefrac{1}{2}(t-3)(t-4),\mathsf{vortex}_{\ref{thm:GMSTpolytime}}(t,k))$-near embedding of detail $k$ into a surface where $H$ does not embed.
\end{enumerate}
Moreover, it holds that

{\centering
  $ \displaystyle
    \begin{aligned}
        \mathsf{adhesion}_{\ref{thm:GMST}}(t,k),\mathsf{apex}_{\ref{thm:GMST}}(t,k),~ \mathsf{vortex}_{\ref{thm:GMST}}(t,k) \in \mathbf{O}\big((t+k)^{11500}\operatorname{poly log}(t+k)\big).
    \end{aligned}
  $
\par}
There also exists a randomised algorithm that, given $k$, $H$, and $G$ as input, finds either an $H$-minor model in $G$ or a tree-decomposition $(T,\beta)$ as above with high probability in time $(t+k+|V(G)|)^{\mathbf{O}(1)}$.
\end{theorem}

The proof of \cref{thm:GMSTpolytime} is almost the same as the proof of \cref{thm:GMST} which is why we only provide a rough sketch here.
The core idea is to prove a strengthening such as \cref{thm:GMST_induction} but replace the algorithm that finds a balanced separator for a bounded size set $X$ of concludes that $X$ is highly linked with \cref{prop:approximateBalancedSet} and replace \cref{thm:algogrid} with \cref{prop:approxAlignedGrid}.
To achieve the desired outcome, the functions of \cref{thm:GMSTpolytime} on the required linkedness of the set $X$ needs to be increased to fit the requirement of roughly $\mathbf{O}(k^{100})$ by \cref{prop:approxAlignedGrid} over the roughly $\mathbf{O}(k^{20})$ required by \cref{thm:algogrid}.
Moreover, due to the approximative nature of \cref{prop:approximateBalancedSet}, the size of the root set $X$, as introduced in \cref{thm:GMST_induction}, needs to increase from $3\cdot x+1$, where $x$ is the required linkedness, to $3\cdot c_{\ref{prop:approximateBalancedSet}}\cdot x\log x+1$ in order to ensure that the inductive procedure stays within the required bounds.

\section{Conclusion}\label{sec:conclusion}

We have proved that the Graph Minor Structure Theorem and the Local Structure Theorem hold with polynomial bounds. This is a tremendous breakthrough in the area of graph minors, which improves many results relying on these theorems. Nevertheless, there are still many related challenging problems, which remain open, and we conclude the paper by discussing some of them.

\paragraph{The Grid Theorem.} Some of the biggest challenges in the area of graph minors concern the Grid Theorem. Recent years have brought many interesting developments on the qualitative side, including refinements of the theorem, which characterise the structure of graphs
excluding a planar graph of a specific type, such as an apex-forest~\cite{DujmovicEJPD2020Minor,HodorLMR2024quickly}, a rectangular \((t \times k)\)-grid for a fixed \(t\)~\cite{HuynhJMSW2022Excluding,Rambaud2025Excluding}, or a planar graph of bounded treedepth~\cite{DujmovicHHJLMMRW2023Gridminor}. This open-ended line of research is quite active, and it would be exciting to see even more structural generalisations of the Grid Theorem.

On the quantitive side, the big question is a tight(er) bound for the Grid Theorem. The gap between the \(\Omega(k^2 \log k)\) lower bound and the \(\mathbf{O}(k^9 \operatorname{polylog} k)\) upper bound remains very wide.
Robertson and Seymour~\cite{RobertsonST1994Quickly} conjectured that the lower bound \(\Omega(k^2 \log k)\) is best possible, and Demaine et al.~\cite{DemaineHK2005Algorithmic} speculate that the optimal bound is \(\Theta(k^3)\).

\begin{problem}
   What are the asymptotics of the function \(f(k)\) describing the smallest treewidth forcing a \((k \times k)\)-grid-minor?
\end{problem}

We are particularly interested in the algorithmic variant of this question.
In the simplest setting, we are interested in algorithms, which given a graph of treewidth at least \(f(k)\) find a model of the \((k \times k)\)-grid.
However, as we have explained in \Cref{sec:polytime}, it is important to be able to not only find an arbitrary model of a grid in a graph, but one which is related to a given tangle or highly linked set.
Kawarabayashi et al.~\cite{KawarabayashiTW2021Quickly} derived an algorithm that finds a model related to a given tangle of order roughly $f(k^2)$, assuming that they are given an algorithm that finds a $k \times k$-grid minor model in a graph of treewidth $f(k)$.
It took us a considerable effort in \Cref{sec:polytime} to derive an algorithm finding a wall related to a given tangle of order \(\Omega(k^{98}\operatorname{poly log} k)\), but for that we had to dive into the details of the algorithm by \cite{ChekuriC2016Polynomial}.
For that reason, we appeal to the brave researchers working on the Grid Theorem to develop more ``localised'' algorithms which find a grid-minor related to a given tangle, or a similar object, such as a highly linked set.

Two advantages of highly linked sets over tangles is that they have an efficient representation in computer memory, and we can efficiently verify if a set is highly linked (see \cref{prop:approximateBalancedSet}).
Recall that every \((k, \alpha)\)-linked set \(X\) determines a tangle of order \(\Omega(k)\), and any tangle of order \(k\) determines a well-linked set of order \(\Omega(k)\), which is \((\Omega(k), \nicefrac{2}{3})\)-linked. In this sense, tangles and highly linked sets are ``approximately equivalent''.

\begin{definition}[\textsc{\(f(k)\)-Grid-Minor}]
    For any given function $f \colon \mathbb{N} \rightarrow \mathbb{N}$, we define the \textsc{\(f(k)\)-Grid-Minor} problem, in which you are given a graph \(G\), an integer $k$, and an \((f(k), \nicefrac{2}{3})\)-linked set \(X \subseteq V(G)\), and the task is to find a model \(\eta\) of the \((k \times k)\)-grid
    such that there exists pairwise vertex-disjoint paths \(P_1, \ldots, P_k\) in \(G\)
    with each \(P_i\) being an \(X\)-\(\eta((i, 1))\)-path.
\end{definition}
Note that we are choosing to state this problem for linked sets instead of tangles, as we believe this should be more concrete than formulating the problem for a tangle.
However, a solution for this problem will also result in a grid minor model which is controlled, in the sense we define in this paper, by the tangle of the linked set that is given as an input.

We concentrate on the setting where \(f(k)\) is a polynomial function, and the running time of an algorithm is \(\mathbf{O}(k^c\cdot |V(G)|^d)\) for some constants \(c, d \ge 0\), but we note that the setting with a \(\mathbf{O}_k(|V(G)|^c) = \mathbf{O}(g(k)\cdot |V(G)|^c)\) running time with an arbitrary function \(g(k)\) is also interesting.

Using the terminology introduced above, \cref{prop:approxAlignedGrid} asserts that there is a polynomial-time randomised algorithm for the \textsc{\(f(k)\)-Grid-Minor} with \(f(k) = \mathbf{O}(k^{98}\operatorname{poly log} k)\).
The first challenge is improving the exponent \(98\).

\begin{problem}
    What is the least constant \(c\) such that \textsc{\(f(k)\)-Grid-Minor}
    admits a polynomial-time algorithm for some function \(f(k) = k^{c+\mathbf{o}(1)}\)? 
\end{problem}

Another problem concerns the explicit bounds on the time-complexity of the \textsc{\(f(k)\)-Grid-Minor} problem. The algorithm by Chekuri and Chuzhoy~\cite{ChekuriC2016Polynomial} is polynomial-time, but an explicit bound on the running time is not known. This is the only part of our algorithm for the GMST which prevents us from providing an explicit bound on the (polynomial) running time.

\begin{problem}
    What is the least constant \(d\) such that \textsc{\(f(k)\)-Grid-Minor}
    admits a \(\mathbf{O}(k^c|V(G)|^{d+\mathbf{o}(1)})\)-time algorithm for some polynomial function \(f(k)\) and some constant \(c \ge 0\)? 
\end{problem}

Finally, it would be interesting to see a non-randomised algorithm for \textsc{\(f(k)\)-Grid-Minor}. As before, an algorithm with an explicit bound on the running time would be desirable.

\begin{problem}
    Does \textsc{\(f(k)\)-Grid-Minor} admit a fully deterministic polynomial-time algorithm for some polynomial function \(f(k)\)? 
\end{problem}

\paragraph{Algorithmic Challenges.}
There have been reports indicating that it might be possible to compute the structural decomposition guaranteed by the GMST in almost linear time in the size of the graph, allowing for an arbitrary runtime in the size of the excluded graph \cite{KorhonenPS2024Minor}.
Such a result would certainly be very interesting and bears the question of how and whether this can be improved.
We leave the formulation of the details of these questions to the potential authors of such a result.

In the context of our proof, we want to highlight one particular bottleneck related to finding certain transactions.
The runtime of our algorithm for the LST depends chiefly on the runtime of the Society Classification Theorem, which in turn depends mainly on the runtime of our variants of Lemma 3.6 from \cite{KawarabayashiTW2021Quickly}.
Our analysis in \Cref{sec:findcrooked} of the algorithm associated with the proof of this lemma resulted in a cubic runtime in the size of the graph.
According to our analysis of our proofs, literally every part of our proof up to \Cref{sec:localtoglobal} that does not use \Cref{lem:crookedrouteback} or \Cref{lem:crookedexistencealgo} runs in linear time in the size of the graph.\footnote{There is also a small exception in our uses of \Cref{lemma:makenestcozy} and \Cref{corollary:maketransactiontaut}. We believe that both of these can easily be pushed to run in linear time with fairly low effort. The same cannot be said for \Cref{lem:crookedrouteback} and \Cref{lem:crookedexistencealgo}.}
This brings us to the following conjecture.

\begin{conjecture}\label{con:findcrooked}
    There exists a polynomial function $f \colon \mathbb{N} \rightarrow \mathbb{N}$, such that the following holds.

    Let $p$ be a positive integer and let $(G,\Omega)$ be a society.
    Then there exists an algorithm running in time $\mathbf{poly}(p)|E(G)|$ that either finds a crooked transaction of order $p$ in $(G,\Omega)$ or a cylindrical rendition $\rho$ of $(G,\Omega)$ in a disk with depth at most $f(p)$.
\end{conjecture}

As evidence towards this conjecture, we point towards our earlier observation that crooked transactions generalise the concept of a cross in a society.
Thus the main result from Graph Minors IX \cite{RobertsonS1990Grapha} is a generalisation of the \hyperref[prop:TwoPaths]{Two Paths Theorem}.
For this theorem there exists several algorithms with nearly linear running times (see \cite{Tholey2006Solving,Tholey2009Improved}) and one algorithm with linear running time \cite{KawarabayashiLR2015Connectivity}.
Thus we are hopeful that this may extend to crooked transactions.

We believe a resolution to this conjecture will ultimately have to find algorithms running in $\mathbf{poly}(p)|E(G)|$-time for \Cref{lem:crookedrouteback} and \Cref{lem:findtransactionefficient}.
This first of these two is somewhat cumbersome to state, but as the second is a fairly natural problem even outside of the realm of graph minors, we state two associated conjectures here.

\begin{conjecture}\label{con:findtransaction}
    There exists a polynomial function $f \colon \mathbb{N} \rightarrow \mathbb{N}$, such that the following holds.

    Let $p$ be a positive integer and let $(G,\Omega)$ be a society with depth at least $f(p)$.
    Then there exists an algorithm running in time $\mathbf{poly}(p)|E(G)|$ that finds a transaction of order $p$ in $(G,\Omega)$.
\end{conjecture}

\begin{conjecture}\label{con:findtransactionorlineardecomp}
    There exists a polynomial function $f \colon \mathbb{N} \rightarrow \mathbb{N}$, such that the following holds.

    Let $p$ be a positive integer and let $(G,\Omega)$ be a society.
    Then there exists an algorithm running in time $\mathbf{poly}(p)|E(G)|$ that either finds a transaction of order $p$ in $(G,\Omega)$ or a linear decomposition of $(G,\Omega)$ with adhesion less than $f(p)$.
\end{conjecture}

Note that both \Cref{con:findcrooked} and \Cref{con:findtransactionorlineardecomp} imply \Cref{con:findtransaction}, whilst neither \Cref{con:findtransaction}, nor \Cref{con:findtransactionorlineardecomp} imply \Cref{con:findcrooked}.
\Cref{con:findtransaction} is likely to be both the easiest to resolve and most accessible of our conjectures.
We further conjecture that the function $f$ in all three of these statements should actually be linear.

\paragraph{Finding a simpler proof.}
One notable feature of the KTW-proof of the GMST is that they endeavoured to keep their proof elementary.
Elementary in this context means that all theorems not proven in their paper(s) can be found in standard graph theory textbooks.
Thus their paper should in theory be accessible to graduate students or even particularly determined undergraduate students.
Note that in this context \cite{KawarabayashiTW2018New} and \cite{KawarabayashiTW2021Quickly} are considered to be a single continuous text containing the proof of the GMST.

We did not attempt to put this restriction upon ourselves and we are under no delusions that our proofs will be seen as simple or elementary.
Within our proof there are however several ideas that we believe simplify the approach of KTW if applied well and, while their proof was certainly a massive leap forward, we do not believe that their approach can be considered to be easy to understand either.
Reiterating Lov{\'a}sz's statement from \cite{Lovasz2005Graph}, we therefore still believe it would be quite important to have a simpler proof of the GMST with explicit bounds.
In particular, a proof that is accessible to graduate students -- maybe even to undergraduate students -- and could be taught within one university semester would be a very valuable addition to the literature.
We believe that it should be possible to find such a proof whilst still delivering modest bounds for all the involved functions, meaning at worst exponential but possibly even polynomial.
As in both our work and that of KTW, this proof may even start from the assumption that the Grid Theorem is true.

Here we want to also explicitly express a demand for a simpler proof of the well-quasi ordering theorem for graph minors \cite{RobertsonS2004Graph}.
This theorem is an important entry point for a set of conjectures that expand on the well-quasi ordering properties of graphs with respect to minors and other relations (see \cite{PaulPT2023Graph,PaulPT2023Universal}).
A more accessible proof would help to facilitate more research in this important area of (meta) graph theory.

\paragraph{Tight bounds.} The definition of \(\Sigma\)-decomposition involves four parameters which should be bounded in the LST: \textsl{(i)} the Euler genus of the surface, \textsl{(ii)} the number of apex vertices, \textsl{(iii)} the number of vortices (called breadth in the case of \(\Sigma\)-decomposition), and \textsl{(iv)} the depth of the vortices. All known proofs of the LST are very complex, so aiming for tight bounds for these parameters may seem to be a hopeless endeavour. Yet, surprisingly, the nature of the proof makes it possible to attain a tight bound for some of these parameters (though not for \(\mathsf{mesh}(t, r)\)).

Already the original proof by Robertson and Seymour, gives an optimal bound for (i). More precisely, we know that the excluded minor \(H\) does not embed in the  surface \(\Sigma\) returned by the LST. This is optimal in the sense that for every surface \(\Sigma\) in which \(H\) does not embed, there exist \(H\)-minor-free graphs with meshes for which the LST requires either the surface \(\Sigma\) or a more complex one (obtained by adding handles and/or crosscaps to \(\Sigma\)). Furthermore, the surface-wall constructed by the proof ensures that
for any graph \(H'\) embeddable in \(\Sigma\), if the parameter \(k\) is sufficiently large, then \(H'\) is a minor of \(G\).

As for the apex vertices (ii), the situation is more complicated.
The Flat Wall Theorem is known to hold with less than \(t-4\) apex vertices~\cite{Chuzhoy2015Improved,GiannopoulouT2013Optimizing,KawarabayashiTW2018New}, which is an optimal bound (a grid with \(t-5\) apex vertices adjacent to each other and to the vertices of the grid contains no \(K_t\)-minor).
In general, the currently known techniques give no hope for getting a bound on the number of apex vertices which is anywhere close to optimal, but we can bound the number of ``bad'' apex vertices. Let us call an apex vertex \emph{major} if it is adjacent to a vertex which does not belong to any vortex. Dvořák and Thomas~\cite{DvorakT2016Listcoloring} showed that if the forbidden minor \(H\) becomes planar after deleting \(a\) vertices, then the LST is satisfied with less than \(a\) major apex vertices (if \(H = K_t\), then \(a=t-4\)). This variant has some interesting applications, including the recent proof of the Clustered Hadwiger Conjecture~\cite{DujmovicEMW2023Proof}.
Wood~\cite{Wood2024Personal} asked if the LST holds with polynomial bounds and less than \(a\) major apex vertices. A positive answer to this question would imply that \(K_t\)-minor-free graphs are \((t-1)\)-colourable with clusters of size polynomial  in \(t\). We suspect that our proof can be modified so that it has the property desired by Wood, although in such case the breadth of the \(\Sigma\)-decomposition would slightly increase. 
\begin{conjecture}
    There exists a version of the LST with polynomial bounds and a number of major apices smaller than the least integer $a$ such that $H$ becomes planar after deleting $a$ apex vertices.
\end{conjecture}

Regarding (iii), the proof by Kawarabayashi et al.~\cite{KawarabayashiTW2021Quickly} gives a good explicit \(2t^2\) bound on the number of vortices. In this paper, we get a slightly better bound of \(\nicefrac{1}{2}(t-3)(t-4)\) by showing that \(K_t\) is a minor of a graph which can be drawn in the plane with \(\nicefrac{1}{2}(t-3)(t-4)\) pairs of edges crossing. Our proof suggests that the optimal bound is the minimum crossing number of a graph containing \(H\) as a minor. Here, the \emph{crossing number} \(\operatorname{cr}(H)\) of a graph \(H\) is the lowest number of edge crossings in a drawing of \(H\) on the plane. Let \(\operatorname{cr}'(H)\) denote the lowest value of \(\operatorname{cr}(H')\) where \(H'\) ranges over all graphs containing \(H\) as a minor.
Note that in general \(\operatorname{cr}'(H)\) can be much smaller than \(\operatorname{cr}(H)\). For example,  \(\operatorname{cr}(K_t) = \Omega(t^4)\)~\cite{DeKlerkMPRS2006improved}, whereas \(\operatorname{cr}'(K_t) \le \nicefrac{1}{2}(t-3)(t-4)\), as we have shown.
We conjecture the following.

\begin{conjecture}
    There exists a version of the LST with polynomial bounds where the breadth of $\rho$ is smaller than $\operatorname{cr}'(H)$.
\end{conjecture}

Note that in order to prove the conjecture, it suffices to show that \cref{lemma:cliqes_in_extended_Dyckwalls}
holds with \(2t+8\) replaced by any polynomial function, \(\nicefrac{1}{2}(t-3)(t-4)\) replaced by \(\operatorname{cr}'(H)\), and \(K_t\) replaced by \(H\).

\begin{figure}
    \centering
    \includegraphics{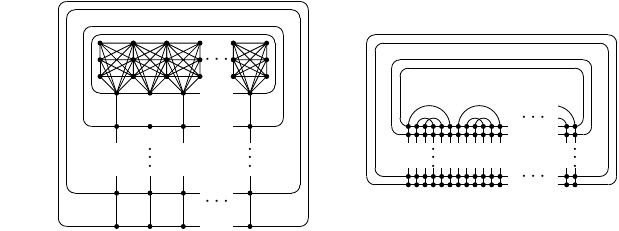}
    \caption{Two constructions of vortex grids of depth \(3\), which are obstructions to the existence of a \(\Sigma\)-decomposition of width \(2\) where $\Sigma$ is any fixed surface. Graphs from the left construction contain all graphs from the right construction, but we were not able to prove the converse.}
    \label{fig:vortex-grids}
\end{figure}

Finally, the depth of the vortices. Thilikos and Wiederrecht~\cite{ThilikosW2024Killing} showed that the only obstruction to a vortex-free \(\Sigma\)-decomposition with bounded parameters is the presence of a ``shallow vortex-grid'' as a minor. There are several ways in which we can construct obstructions to a \(\Sigma\)-decomposition of small depth (see \cref{fig:vortex-grids}), and it is not clear whether there exists a nice construction of graphs characterising graphs admitting a \(\Sigma\)-decomposition of depth, say, at most \(3\).

Our proof yields vortices of depth depending on the parameter \(k\). Before we ask about a tight bound on the depth, we should ask if the depth can be bounded in terms of \(t\) alone. We believe that this is possible.
\begin{conjecture}
    The LST can be proven with polynomial bounds and with depth of the \(\Sigma\)-decomposition bounded in terms of \(t\) alone.
\end{conjecture}

While we can bound the parameters of a \(\Sigma\)-decomposition separately, probably in an optimal \(\Sigma\)-decomposition they never need to reach the upper bounds simultaneously. We would like to conclude the paper with a somewhat provocative conjecture, which generalises all of the conjectures above.
It provides more fine-grained bounds on the parameters of the \(\Sigma\)-decomposition returned by the LST, with bounds depending on possible near embeddings of the forbidden minor \(H\).

\begin{conjecture}
    The LST can be proven with polynomial bounds so that when we apply it to \((H, G, M)\) and it returns \((\Sigma, A, M', \rho)\), then there does not exist an \((a, b, w)\)-near embedding of \(H\) in \(\Sigma\) such that \(a\) is at most the number of major apex vertices in \(A\), \(b\) is at most the breadth of \(\rho\), and \(w\in\mathbf{o}(d)\) where \(d\)
    is the depth of \(\rho\).
\end{conjecture}

\bigskip

\textbf{Acknowledgements.}
We thank Reinhard Diestel, Gwena\"el Joret, Stefan Kober, Stephan Kreutzer, Sang-il Oum, Dimitrios Thilikos, and Yelena Yuditsky for helpful and inspiring discussions.
Additionally, we thank Paul Wollan for several pleasant conversations which provided background information on \cite{KawarabayashiTW2021Quickly} and for communicating the ideas of the proof of \cref{thm:wollan2} to us.
The second author thanks Grzegorz Guśpiel.

We would also like to thank Jim Geelen.
While he was not directly involved in this project, the first grasp any of us got on the deeper levels of Robertson and Seymour's theory of graph minors was through his excellent series of lectures.
We are grateful to him for opening this world to us and for generously sharing his notes and recordings.

\bibliographystyle{alphaurl}
\bibliography{arxiv_literature}

\newcommand{\etalchar}[1]{$^{#1}$}
\begin{thebibliography}{dKMP{\etalchar{+}}06}

\bibitem[AFJ{\etalchar{+}}25]{AprileFJKSWY2025Integer}
Manuel Aprile, Samuel Fiorini, Gwena{\"e}l Joret, Stefan Kober, Micha{\l}~T. Seweryn, Stefan Weltge, and Yelena Yuditsky.
\newblock Integer programs with nearly totally unimodular matrices: The cographic case.
\newblock In {\em Proceedings of the 2025 {{Annual ACM-SIAM Symposium}} on {{Discrete Algorithms}} ({{SODA}})}, pages 2301--2312, Philadelphia, PA, January 2025. {Society for Industrial and Applied Mathematics}.
\newblock \href {https://doi.org/10.1137/1.9781611978322} {\path{doi:10.1137/1.9781611978322}}.

\bibitem[AK10]{AdlerK2010Lower}
Isolde Adler and Philipp~Klaus Krause.
\newblock A lower bound for the tree-width of planar graphs with vital linkages, November 2010.
\newblock \href {https://arxiv.org/abs/1011.2136} {\path{arXiv:1011.2136}}.

\bibitem[AKK{\etalchar{+}}11]{AdlerKKLST2011Tight}
Isolde Adler, Stavros~G. Kolliopoulos, Philipp~Klaus Krause, Daniel Lokshtanov, Saket Saurabh, and Dimitrios Thilikos.
\newblock Tight {{Bounds}} for {{Linkages}} in {{Planar Graphs}}.
\newblock In {\em Automata, {{Languages}} and {{Programming}}}, volume 6755, pages 110--121. Springer Berlin Heidelberg, Berlin, Heidelberg, 2011.
\newblock \href {https://doi.org/10.1007/978-3-642-22006-7_10} {\path{doi:10.1007/978-3-642-22006-7_10}}.

\bibitem[Ami10]{Amir2010Approximation}
Eyal Amir.
\newblock Approximation {{Algorithms}} for {{Treewidth}}.
\newblock {\em Algorithmica}, 56(4):448--479, April 2010.
\newblock \href {https://doi.org/10.1007/s00453-008-9180-4} {\path{doi:10.1007/s00453-008-9180-4}}.

\bibitem[Arn85]{Arnborg1985Efficient}
Stefan Arnborg.
\newblock Efficient algorithms for combinatorial problems on graphs with bounded decomposability --- {{A}} survey.
\newblock {\em BIT}, 25(1):1--23, March 1985.
\newblock \href {https://doi.org/10.1007/BF01934985} {\path{doi:10.1007/BF01934985}}.

\bibitem[BB72]{BerteleB1972Nonserial}
Umberto Bertel{\`e} and Francesco Brioschi.
\newblock {\em Non-Serial Dynamic Programming}.
\newblock Academic Press, New York, June 1972.

\bibitem[BB73]{BerteleB1973Nonserial}
Umberto Bertel{\`e} and Francesco Brioschi.
\newblock On non-serial dynamic programming.
\newblock {\em Journal of Combinatorial Theory, Series A}, 14(2):137--148, March 1973.
\newblock \href {https://doi.org/10.1016/0097-3165(73)90016-2} {\path{doi:10.1016/0097-3165(73)90016-2}}.

\bibitem[BHJ18]{BruhnHJ2018Frames}
Henning Bruhn, Matthias Heinlein, and Felix Joos.
\newblock Frames, ${A}$-{P}aths, and the {E}rdős-{P}ósa {P}roperty.
\newblock {\em SIAM Journal on Discrete Mathematics}, 32(2):1246--1260, January 2018.
\newblock \href {https://doi.org/10.1137/17M1148542} {\path{doi:10.1137/17M1148542}}.

\bibitem[CC16]{ChekuriC2016Polynomial}
Chandra Chekuri and Julia Chuzhoy.
\newblock Polynomial {{Bounds}} for the {{Grid-Minor Theorem}}.
\newblock {\em Journal of the ACM}, 63(5):1--65, December 2016.
\newblock \href {https://doi.org/10.1145/2820609} {\path{doi:10.1145/2820609}}.

\bibitem[CE13]{ChekuriE2013Polylogarithmic}
Chandra Chekuri and Alina Ene.
\newblock Poly-logarithmic {{Approximation}} for {{Maximum Node Disjoint Paths}} with {{Constant Congestion}}.
\newblock In {\em Proceedings of the {{Twenty-Fourth Annual ACM-SIAM Symposium}} on {{Discrete Algorithms}}}, pages 326--341. {Society for Industrial and Applied Mathematics}, January 2013.
\newblock \href {https://doi.org/10.1137/1.9781611973105.24} {\path{doi:10.1137/1.9781611973105.24}}.

\bibitem[Chu15]{Chuzhoy2015Improved}
Julia Chuzhoy.
\newblock Improved {{Bounds}} for the {{Flat Wall Theorem}}.
\newblock In {\em Proceedings of the {{Twenty-Sixth Annual ACM-SIAM Symposium}} on {{Discrete Algorithms}}}, pages 256--275. {Society for Industrial and Applied Mathematics}, October 2015.
\newblock \href {https://doi.org/10.1137/1.9781611973730.20} {\path{doi:10.1137/1.9781611973730.20}}.

\bibitem[CKS05]{ChekuriKS2005Multicommodity}
Chandra Chekuri, Sanjeev Khanna, and F.~Bruce Shepherd.
\newblock Multicommodity flow, well-linked terminals, and routing problems.
\newblock In {\em Proceedings of the Thirty-Seventh Annual {{ACM}} Symposium on {{Theory}} of Computing}, pages 183--192, Baltimore MD USA, May 2005. ACM.
\newblock \href {https://doi.org/10.1145/1060590.1060618} {\path{doi:10.1145/1060590.1060618}}.

\bibitem[CLPP23]{Cohen-AddadLPP2023Planar}
Vincent {Cohen-Addad}, Hung Le, Marcin Pilipczuk, and Micha{\l} Pilipczuk.
\newblock Planar and {{Minor-Free Metrics Embed}} into {{Metrics}} of {{Polylogarithmic Treewidth}} with {{Expected Multiplicative Distortion Arbitrarily Close}} to 1*.
\newblock In {\em 2023 {{IEEE}} 64th {{Annual Symposium}} on {{Foundations}} of {{Computer Science}} ({{FOCS}})}, pages 2262--2277, Santa Cruz, CA, USA, November 2023. IEEE.
\newblock \href {https://doi.org/10.1109/FOCS57990.2023.00140} {\path{doi:10.1109/FOCS57990.2023.00140}}.

\bibitem[Cou90]{Courcelle1990Monadic}
Bruno Courcelle.
\newblock The monadic second-order logic of graphs. {{I}}. {{Recognizable}} sets of finite graphs.
\newblock {\em Information and Computation}, 85(1):12--75, March 1990.
\newblock \href {https://doi.org/10.1016/0890-5401(90)90043-H} {\path{doi:10.1016/0890-5401(90)90043-H}}.

\bibitem[CT21]{ChuzhoyT2021Tighter}
Julia Chuzhoy and Zihan Tan.
\newblock Towards tight(er) bounds for the {{Excluded Grid Theorem}}.
\newblock {\em Journal of Combinatorial Theory, Series B}, 146:219--265, January 2021.
\newblock \href {https://doi.org/10.1016/j.jctb.2020.09.010} {\path{doi:10.1016/j.jctb.2020.09.010}}.

\bibitem[CX22]{CurticapeanX2022Parameterizing}
Radu Curticapean and Mingji Xia.
\newblock Parameterizing the {{Permanent}}: {{Hardness}} for fixed excluded minors.
\newblock In {\em Proceedings of the 2022 {{Symposium}} on {{Simplicity}} in {{Algorithms}} ({{SOSA}})}, pages 297--307, Philadelphia, PA, January 2022. {Society for Industrial and Applied Mathematics}.
\newblock \href {https://doi.org/10.1137/1.9781611977066} {\path{doi:10.1137/1.9781611977066}}.

\bibitem[DDO{\etalchar{+}}04]{DeVosDOSRSV2004Excluding}
Matt DeVos, Guoli Ding, Bogdan Oporowski, Daniel~P. Sanders, Bruce Reed, Paul Seymour, and Dirk Vertigan.
\newblock Excluding any graph as a minor allows a low tree-width 2-coloring.
\newblock {\em Journal of Combinatorial Theory, Series B}, 91(1):25--41, May 2004.
\newblock \href {https://doi.org/10.1016/j.jctb.2003.09.001} {\path{doi:10.1016/j.jctb.2003.09.001}}.

\bibitem[DEJ{\etalchar{+}}20]{DujmovicEJPD2020Minor}
Vida Dujmovi\'{c}, David Eppstein, Gwena\"{e}l Joret, Pat Morin, and David~R. Wood.
\newblock Minor-closed graph classes with bounded layered pathwidth.
\newblock {\em SIAM Journal on Discrete Mathematics}, 34(3):1693--1709, 2020.
\newblock \href {https://doi.org/10.1137/18M122162X} {\path{doi:10.1137/18M122162X}}.

\bibitem[DEMW23]{DujmovicEMW2023Proof}
Vida Dujmovi{\'c}, Louis Esperet, Pat Morin, and David~R. Wood.
\newblock Proof of the {{Clustered Hadwiger Conjecture}}.
\newblock In {\em 2023 {{IEEE}} 64th {{Annual Symposium}} on {{Foundations}} of {{Computer Science}} ({{FOCS}})}, pages 1921--1930, Santa Cruz, CA, USA, November 2023. IEEE.
\newblock \href {https://doi.org/10.1109/FOCS57990.2023.00116} {\path{doi:10.1109/FOCS57990.2023.00116}}.

\bibitem[DF95]{DowneyF1995Parameterized}
Rodney~G. Downey and Michael~R. Fellows.
\newblock Parameterized {{Computational Feasibility}}.
\newblock In {\em Feasible {{Mathematics II}}}, pages 219--244. Birkh{\"a}user Boston, Boston, MA, 1995.
\newblock \href {https://doi.org/10.1007/978-1-4612-2566-9_7} {\path{doi:10.1007/978-1-4612-2566-9_7}}.

\bibitem[DFHT05]{DemaineFHT2005Subexponential}
Erik~D. Demaine, Fedor~V. Fomin, MohammadTaghi Hajiaghayi, and Dimitrios~M. Thilikos.
\newblock Subexponential {P}arameterized {A}lgorithms on {B}ounded-{G}enus {G}raphs and ${H}$-{M}inor-{F}ree {G}raphs.
\newblock {\em Journal of the ACM}, 52(6):866--893, November 2005.
\newblock \href {https://doi.org/10.1145/1101821.1101823} {\path{doi:10.1145/1101821.1101823}}.

\bibitem[DGK07]{DawarGK2007Locally}
Anuj Dawar, Martin Grohe, and Stephan Kreutzer.
\newblock Locally {{Excluding}} a {{Minor}}.
\newblock In {\em 22nd {{Annual IEEE Symposium}} on {{Logic}} in {{Computer Science}} ({{LICS}} 2007)}, pages 270--279, Wroclaw, Poland, 2007. IEEE.
\newblock \href {https://doi.org/10.1109/LICS.2007.31} {\path{doi:10.1109/LICS.2007.31}}.

\bibitem[DHH{\etalchar{+}}23]{DujmovicHHJLMMRW2023Gridminor}
Vida Dujmovi{\'c}, Robert Hickingbotham, J{\k e}drzej Hodor, Gwena{\"e}l Joret, Hoang La, Piotr Micek, Pat Morin, Cl{\'e}ment Rambaud, and David~R. Wood.
\newblock The grid-minor theorem revisited, July 2023.
\newblock \href {https://arxiv.org/abs/2307.02816} {\path{arXiv:2307.02816}}.

\bibitem[DHK05]{DemaineHK2005Algorithmic}
Erik~D. Demaine, MohammadTaghi Hajiaghayi, and Ken-ichi Kawarabayashi.
\newblock Algorithmic {{Graph Minor Theory}}: {{Decomposition}}, {{Approximation}}, and {{Coloring}}.
\newblock In {\em 46th {{Annual IEEE Symposium}} on {{Foundations}} of {{Computer Science}} ({{FOCS}}'05)}, pages 637--646, Pittsburgh, PA, USA, 2005. IEEE.
\newblock \href {https://doi.org/10.1109/SFCS.2005.14} {\path{doi:10.1109/SFCS.2005.14}}.

\bibitem[DHK10]{DemaineHK2010Decomposition}
Erik~D. Demaine, MohammadTaghi Hajiaghayi, and Ken-ichi Kawarabayashi.
\newblock Decomposition, {{Approximation}}, and {{Coloring}} of {{Odd-Minor-Free Graphs}}.
\newblock In {\em Proceedings of the {{Twenty-First Annual ACM-SIAM Symposium}} on {{Discrete Algorithms}}}, pages 329--344. {Society for Industrial and Applied Mathematics}, January 2010.
\newblock \href {https://doi.org/10.1137/1.9781611973075.28} {\path{doi:10.1137/1.9781611973075.28}}.

\bibitem[DHK11]{DemaineHK2011Contraction}
Erik~D. Demaine, MohammadTaghi Hajiaghayi, and Ken-ichi Kawarabayashi.
\newblock Contraction {D}ecomposition in ${H}$-{M}inor-{F}ree {G}raphs and {A}lgorithmic {A}pplications.
\newblock In {\em Proceedings of the forty-third annual ACM symposium on Theory of computing}, pages 441--450, San Jose California USA, June 2011. ACM.
\newblock \href {https://doi.org/10.1145/1993636.1993696} {\path{doi:10.1145/1993636.1993696}}.

\bibitem[Die24]{Diestel2024Tangles}
Reinhard Diestel.
\newblock {\em Tangles: {{A Structural Approach}} to {{Artificial Intelligence}} in the {{Empirical Sciences}}}.
\newblock Cambridge University Press, 1 edition, May 2024.
\newblock \href {https://doi.org/10.1017/9781009473323} {\path{doi:10.1017/9781009473323}}.

\bibitem[dKMP{\etalchar{+}}06]{DeKlerkMPRS2006improved}
E.~de~Klerk, J.~Maharry, D.~V. Pasechnik, R.~B. Richter, and G.~Salazar.
\newblock Improved bounds for the crossing numbers of ${K}_{m,n}$ and ${K}_n$.
\newblock {\em SIAM Journal on Discrete Mathematics}, 20(1):189--202, 2006.
\newblock \href {https://doi.org/10.1137/S0895480104442741} {\path{doi:10.1137/S0895480104442741}}.

\bibitem[DKMW12]{DiestelKMW2012Excluded}
Reinhard Diestel, Ken-ichi Kawarabayashi, Theodor M{\"u}ller, and Paul Wollan.
\newblock On the excluded minor structure theorem for graphs of large tree-width.
\newblock {\em Journal of Combinatorial Theory, Series B}, 102(6):1189--1210, November 2012.
\newblock \href {https://doi.org/10.1016/j.jctb.2012.07.001} {\path{doi:10.1016/j.jctb.2012.07.001}}.

\bibitem[DT16]{DvorakT2016Listcoloring}
Zden{\v e}k Dvo{\v r}{\'a}k and Robin Thomas.
\newblock List-coloring apex-minor-free graphs, December 2016.
\newblock \href {https://arxiv.org/abs/1401.1399} {\path{arXiv:1401.1399}}.

\bibitem[Dyc88]{Dyck1888Beitraege}
Walther Dyck.
\newblock {Beitr{\"a}ge zur Analysis situs: I. Aufsatz. Ein- und zweidimensionale Mannigfaltigkeiten}.
\newblock {\em Mathematische Annalen}, 32(4):457--512, December 1888.
\newblock \href {https://doi.org/10.1007/BF01443580} {\path{doi:10.1007/BF01443580}}.

\bibitem[ES35]{ErdosS1935Combinatorial}
Paul Erd{\H o}s and George Szekeres.
\newblock A combinatorial problem in geometry.
\newblock {\em Compositio Mathematica}, 2:463--470, 1935.

\bibitem[Eul58]{Euler1758Elementa}
Leonhard Euler.
\newblock Elementa doctrinae solidorum.
\newblock {\em Novi Commentarii academiae scientiarum Petropolitanae}, 4:109--150, 1758.

\bibitem[FF56]{FordF1956Maximal}
L.~R. Ford and D.~R. Fulkerson.
\newblock Maximal {{Flow Through}} a {{Network}}.
\newblock {\em Canadian Journal of Mathematics}, 8:399--404, 1956.
\newblock \href {https://doi.org/10.4153/CJM-1956-045-5} {\path{doi:10.4153/CJM-1956-045-5}}.

\bibitem[FHDM24]{FuladiHD2024Short}
Niloufar Fuladi, Alfredo Hubard, and Arnaud De~Mesmay.
\newblock Short {{Topological Decompositions}} of {{Non-orientable Surfaces}}.
\newblock {\em Discrete \& Computational Geometry}, 72(2):783--830, September 2024.
\newblock \href {https://doi.org/10.1007/s00454-023-00580-3} {\path{doi:10.1007/s00454-023-00580-3}}.

\bibitem[FJWY22]{FioriniJWY2022Integer}
Samuel Fiorini, Gwena{\"e}l Joret, Stefan Weltge, and Yelena Yuditsky.
\newblock Integer programs with bounded subdeterminants and two nonzeros per row.
\newblock In {\em 2021 {{IEEE}} 62nd {{Annual Symposium}} on {{Foundations}} of {{Computer Science}} ({{FOCS}})}, pages 13--24, Denver, CO, USA, February 2022. IEEE.
\newblock \href {https://doi.org/10.1109/FOCS52979.2021.00011} {\path{doi:10.1109/FOCS52979.2021.00011}}.

\bibitem[Ful23]{Fuladi2023Embedded}
Niloufar Fuladi.
\newblock {\em Embedded {{Graphs}}: {{Crossings}} and {{Decompositions}}}.
\newblock Doctoral {{Thesis}}, Universit{\'e} Gustave Eiffel, 2023.

\bibitem[FW99]{FrancisW1999Conways}
George~K. Francis and Jeffrey~R. Weeks.
\newblock Conway's {{ZIP Proof}}.
\newblock {\em The American Mathematical Monthly}, 106(5):393--399, May 1999.
\newblock \href {https://doi.org/10.1080/00029890.1999.12005061} {\path{doi:10.1080/00029890.1999.12005061}}.

\bibitem[Gal64]{Gallai1964MaximumMinimum}
T.~Gallai.
\newblock {Maximum-Minimum S{\"a}tze und verallgemeinerte Faktoren von Graphen}.
\newblock {\em Acta Mathematica Academiae Scientiarum Hungaricae}, 12(1-2):131--173, March 1964.
\newblock \href {https://doi.org/10.1007/BF02066678} {\path{doi:10.1007/BF02066678}}.

\bibitem[GGW15]{GeelenGW2015Highly}
Jim Geelen, Bert Gerards, and Geoff Whittle.
\newblock The {{Highly Connected Matroids}} in {{Minor-Closed Classes}}.
\newblock {\em Annals of Combinatorics}, 19(1):107--123, March 2015.
\newblock \href {https://doi.org/10.1007/s00026-015-0251-3} {\path{doi:10.1007/s00026-015-0251-3}}.

\bibitem[GH23]{GavoilleH2023MinorUniversal}
Cyril Gavoille and Claire Hilaire.
\newblock Minor-{{Universal Graph}} for {{Graphs}} on {{Surfaces}}, May 2023.
\newblock \href {https://arxiv.org/abs/2305.06673} {\path{arXiv:2305.06673}}.

\bibitem[GHR18]{GeelenHR2018Explicit}
Jim Geelen, Tony Huynh, and R.~Bruce Richter.
\newblock Explicit bounds for graph minors.
\newblock {\em Journal of Combinatorial Theory, Series B}, 132:80--106, September 2018.
\newblock \href {https://doi.org/10.1016/j.jctb.2018.03.004} {\path{doi:10.1016/j.jctb.2018.03.004}}.

\bibitem[GKKW24]{GorskyKKW2024Packing}
Maximilian Gorsky, Ken-ichi Kawarabayashi, Stephan Kreutzer, and Sebastian Wiederrecht.
\newblock Packing {{Even Directed Circuits Quarter-Integrally}}.
\newblock In {\em Proceedings of the 56th {{Annual ACM Symposium}} on {{Theory}} of {{Computing}}}, pages 692--703, Vancouver BC Canada, June 2024. ACM.
\newblock \href {https://doi.org/10.1145/3618260.3649682} {\path{doi:10.1145/3618260.3649682}}.

\bibitem[GKR13]{GroheKR2013Simple}
Martin Grohe, Ken-ichi Kawarabayashi, and Bruce~A. Reed.
\newblock A {{Simple Algorithm}} for the {{Graph Minor Decomposition}} - {{Logic}} meets {{Structural Graph Theory}}--.
\newblock In {\em Proceedings of the {{Twenty-Fourth Annual ACM-SIAM Symposium}} on {{Discrete Algorithms}}}, pages 414--431. {Society for Industrial and Applied Mathematics}, January 2013.
\newblock \href {https://doi.org/10.1137/1.9781611973105.30} {\path{doi:10.1137/1.9781611973105.30}}.

\bibitem[GM12]{GroheM2012Structure}
Martin Grohe and D{\'a}niel Marx.
\newblock Structure theorem and isomorphism test for graphs with excluded topological subgraphs.
\newblock In {\em Proceedings of the Forty-Fourth Annual {{ACM}} Symposium on {{Theory}} of Computing}, pages 173--192, New York New York USA, May 2012. ACM.
\newblock \href {https://doi.org/10.1145/2213977.2213996} {\path{doi:10.1145/2213977.2213996}}.

\bibitem[GNW21]{GeelenNW2021Excluding}
Jim Geelen, Peter Nelson, and Zach Walsh.
\newblock Excluding a line from $\mathbb{C}$-representable matroids, January 2021.
\newblock \href {https://arxiv.org/abs/2101.12000} {\path{arXiv:2101.12000}}.

\bibitem[Gor24]{Gorsky2024Structure}
Maximilian Gorsky.
\newblock {\em The Structure of (Even) Directed Cycles}.
\newblock PhD thesis, Technische Universit{\"a}t Berlin, September 2024.
\newblock \href {https://doi.org/10.14279/depositonce-21276} {\path{doi:10.14279/depositonce-21276}}.

\bibitem[Gro03]{Grohe2003Local}
Martin Grohe.
\newblock Local {{Tree-Width}}, {{Excluded Minors}}, and {{Approximation Algorithms}}.
\newblock {\em Combinatorica}, 23(4):613--632, December 2003.
\newblock \href {https://doi.org/10.1007/s00493-003-0037-9} {\path{doi:10.1007/s00493-003-0037-9}}.

\bibitem[GST23]{GolovachST2023ModelChecking}
Petr~A. Golovach, Giannos Stamoulis, and Dimitrios~M. Thilikos.
\newblock Model-{{Checking}} for {{First-Order Logic}} with {{Disjoint Paths Predicates}} in {{Proper Minor-Closed Graph Classes}}.
\newblock In {\em Proceedings of the 2023 {{Annual ACM-SIAM Symposium}} on {{Discrete Algorithms}} ({{SODA}})}, pages 3684--3699, Philadelphia, PA, January 2023. {Society for Industrial and Applied Mathematics}.
\newblock \href {https://doi.org/10.1137/1.9781611977554} {\path{doi:10.1137/1.9781611977554}}.

\bibitem[GT13]{GiannopoulouT2013Optimizing}
Archontia~C. Giannopoulou and Dimitrios~M. Thilikos.
\newblock Optimizing the {{Graph Minors Weak Structure Theorem}}.
\newblock {\em SIAM Journal on Discrete Mathematics}, 27(3):1209--1227, January 2013.
\newblock \href {https://doi.org/10.1137/110857027} {\path{doi:10.1137/110857027}}.

\bibitem[Hal76]{Halin1976$S$functions}
Rudolf Halin.
\newblock ${S}$-functions for graphs.
\newblock {\em Journal of Geometry}, 8(1-2):171--186, March 1976.
\newblock \href {https://doi.org/10.1007/BF01917434} {\path{doi:10.1007/BF01917434}}.

\bibitem[Hea90]{Heawood1890MapColour}
Percy~J Heawood.
\newblock Map-{{Colour Theorem}}.
\newblock {\em Quarterly Journal of Pure and Applied Mathematics}, 24:332--338, 1890.

\bibitem[HJM{\etalchar{+}}22]{HuynhJMSW2022Excluding}
Tony Huynh, Gwena{\"e}l Joret, Piotr Micek, Micha{\l}~T. Seweryn, and Paul Wollan.
\newblock Excluding a {{Ladder}}.
\newblock {\em Combinatorica}, 42(3):405--432, June 2022.
\newblock \href {https://doi.org/10.1007/s00493-021-4592-8} {\path{doi:10.1007/s00493-021-4592-8}}.

\bibitem[HLMR24]{HodorLMR2024quickly}
J{\k{e}}drzej Hodor, Hoang La, Piotr Micek, and Cl{\'e}ment Rambaud.
\newblock Quickly excluding an apex-forest, April 2024.
\newblock \href {https://arxiv.org/abs/2404.17306} {\path{arXiv:2404.17306}}.

\bibitem[Jun70]{Jung1970Verallgemeinerung}
Heinz~A. Jung.
\newblock Eine verallgemeinerung des $n$-fachen zusammenhangs für graphen.
\newblock {\em Mathematische Annalen}, 187(2):95--103, June 1970.
\newblock \href {https://doi.org/10.1007/BF01350174} {\path{doi:10.1007/BF01350174}}.

\bibitem[JW13]{JoretW2013Complete}
Gwena{\"e}l Joret and David~R. Wood.
\newblock Complete graph minors and the graph minor structure theorem.
\newblock {\em Journal of Combinatorial Theory, Series B}, 103(1):61--74, January 2013.
\newblock \href {https://doi.org/10.1016/j.jctb.2012.09.001} {\path{doi:10.1016/j.jctb.2012.09.001}}.

\bibitem[Kau11]{Kauffman2011Seven}
Louis~H. Kauffman.
\newblock Seven knots and knots in the seven-color map.
\newblock In {\em Homage to a Pied Puzzler}, pages 75--83. A K Peters, Ltd., Wellesley, Massachusetts, 2011.

\bibitem[Kaw07]{Kawarabayashi2007Half}
Ken-ichi Kawarabayashi.
\newblock Half integral packing, {{Erd{\H o}s-Pos{\'a}-property}} and graph minors.
\newblock In {\em {{SODA}} '07: {{Proceedings}} of the Eighteenth Annual {{ACM-SIAM}} Symposium on {{Discrete}} Algorithms}, SODA '07, page 1187–1196, USA, 2007. Society for Industrial and Applied Mathematics.

\bibitem[KK20]{KawarabayashiK2020Linear}
Ken-ichi Kawarabayashi and Yusuke Kobayashi.
\newblock Linear min-max relation between the treewidth of an {{$H$-minor-free}} graph and its largest grid minor.
\newblock {\em Journal of Combinatorial Theory, Series B}, 141:165--180, March 2020.
\newblock \href {https://doi.org/10.1016/j.jctb.2019.07.007} {\path{doi:10.1016/j.jctb.2019.07.007}}.

\bibitem[KLR15]{KawarabayashiLR2015Connectivity}
Ken-ichi Kawarabayashi, Zhentao Li, and Bruce~A. Reed.
\newblock Connectivity {{Preserving Iterative Compaction}} and {{Finding}} 2 {{Disjoint Rooted Paths}} in {{Linear Time}}, September 2015.
\newblock \href {https://arxiv.org/abs/1509.07680} {\path{arXiv:1509.07680}}.

\bibitem[KM07]{KawarabayashiM2007Recent}
Ken-ichi Kawarabayashi and Bojan Mohar.
\newblock Some {{Recent Progress}} and {{Applications}} in {{Graph Minor Theory}}.
\newblock {\em Graphs and Combinatorics}, 23(1):1--46, February 2007.
\newblock \href {https://doi.org/10.1007/s00373-006-0684-x} {\path{doi:10.1007/s00373-006-0684-x}}.

\bibitem[KPS24a]{KorhonenPS2024Minora}
Tuukka Korhonen, Micha{\l} Pilipczuk, and Giannos Stamoulis.
\newblock Minor {{Containment}} and {{Disjoint Paths}} in {{Almost-Linear Time}}.
\newblock In {\em 2024 {{IEEE}} 65th {{Annual Symposium}} on {{Foundations}} of {{Computer Science}} ({{FOCS}})}, pages 53--61, Chicago, IL, USA, October 2024. IEEE.
\newblock \href {https://doi.org/10.1109/FOCS61266.2024.00014} {\path{doi:10.1109/FOCS61266.2024.00014}}.

\bibitem[KPS24b]{KorhonenPS2024Minor}
Tuukka Korhonen, Micha{\l} Pilipczuk, and Giannos Stamoulis.
\newblock Minor {{Containment}} and {{Disjoint Paths}} in almost-linear time, April 2024.
\newblock \href {https://arxiv.org/abs/2404.03958} {\path{arXiv:2404.03958}}.

\bibitem[KRV09]{KhandekarRV2009Graph}
Rohit Khandekar, Satish Rao, and Umesh Vazirani.
\newblock Graph partitioning using single commodity flows.
\newblock {\em Journal of the ACM}, 56(4):1--15, June 2009.
\newblock \href {https://doi.org/10.1145/1538902.1538903} {\path{doi:10.1145/1538902.1538903}}.

\bibitem[KTW18]{KawarabayashiTW2018New}
Ken-ichi Kawarabayashi, Robin Thomas, and Paul Wollan.
\newblock A new proof of the flat wall theorem.
\newblock {\em Journal of Combinatorial Theory, Series B}, 129:204--238, March 2018.
\newblock \href {https://doi.org/10.1016/j.jctb.2017.09.006} {\path{doi:10.1016/j.jctb.2017.09.006}}.

\bibitem[KTW21]{KawarabayashiTW2021Quickly}
Ken-ichi Kawarabayashi, Robin Thomas, and Paul Wollan.
\newblock Quickly excluding a non-planar graph, January 2021.
\newblock \href {https://arxiv.org/abs/2010.12397} {\path{arXiv:2010.12397}}.

\bibitem[KW10]{KawarabayashiW2010Shorter}
Ken-ichi Kawarabayashi and Paul Wollan.
\newblock A shorter proof of the graph minor algorithm: The unique linkage theorem.
\newblock In {\em Proceedings of the Forty-Second {{ACM}} Symposium on {{Theory}} of Computing}, pages 687--694, Cambridge Massachusetts USA, June 2010. ACM.
\newblock \href {https://doi.org/10.1145/1806689.1806784} {\path{doi:10.1145/1806689.1806784}}.

\bibitem[KW11]{KawarabayashiW2011Simpler}
Ken-ichi Kawarabayashi and Paul Wollan.
\newblock A simpler algorithm and shorter proof for the graph minor decomposition.
\newblock In {\em Proceedings of the Forty-Third Annual {{ACM}} Symposium on {{Theory}} of Computing}, pages 451--458, San Jose California USA, June 2011. ACM.
\newblock \href {https://doi.org/10.1145/1993636.1993697} {\path{doi:10.1145/1993636.1993697}}.

\bibitem[Liu22]{Liu2022Packing}
Chun-Hung Liu.
\newblock Packing topological minors half-integrally.
\newblock {\em Journal of the London Mathematical Society}, 106(3):2193--2267, October 2022.
\newblock \href {https://doi.org/10.1112/jlms.12633} {\path{doi:10.1112/jlms.12633}}.

\bibitem[Lov05]{Lovasz2005Graph}
L{\'a}szl{\'o} Lov{\'a}sz.
\newblock Graph minor theory.
\newblock {\em Bulletin of the American Mathematical Society}, 43(1):75--86, October 2005.
\newblock \href {https://doi.org/10.1090/S0273-0979-05-01088-8} {\path{doi:10.1090/S0273-0979-05-01088-8}}.

\bibitem[LPVV01]{LazarusPVV2001Computing}
Francis Lazarus, Michel Pocchiola, Gert Vegter, and Anne Verroust.
\newblock Computing a canonical polygonal schema of an orientable triangulated surface.
\newblock In {\em Proceedings of the Seventeenth Annual Symposium on {{Computational}} Geometry}, pages 80--89, Medford Massachusetts USA, June 2001. ACM.
\newblock \href {https://doi.org/10.1145/378583.378630} {\path{doi:10.1145/378583.378630}}.

\bibitem[LS15]{LeafS2015Treewidth}
Alexander Leaf and Paul Seymour.
\newblock Tree-width and planar minors.
\newblock {\em Journal of Combinatorial Theory, Series B}, 111:38--53, March 2015.
\newblock \href {https://doi.org/10.1016/j.jctb.2014.09.003} {\path{doi:10.1016/j.jctb.2014.09.003}}.

\bibitem[LSZ20]{LokshtanovSZ2020Efficient}
Daniel Lokshtanov, Saket Saurabh, and Meirav Zehavi.
\newblock Efficient {{Graph Minors Theory}} and {{Parameterized Algorithms}} for ({{Planar}}) {{Disjoint Paths}}.
\newblock In {\em Treewidth, {{Kernels}}, and {{Algorithms}}}, volume 12160, pages 112--128. Springer International Publishing, Cham, 2020.
\newblock \href {https://doi.org/10.1007/978-3-030-42071-0_9} {\path{doi:10.1007/978-3-030-42071-0_9}}.

\bibitem[Men27]{Menger1927Zur}
Karl Menger.
\newblock Zur allgemeinen {{Kurventheorie}}.
\newblock {\em Fundamenta Mathematicae}, 10(1):96--115, 1927.
\newblock \href {https://doi.org/10.4064/fm-10-1-96-115} {\path{doi:10.4064/fm-10-1-96-115}}.

\bibitem[Moh95]{Mohar1995Uniqueness}
Bojan Mohar.
\newblock Uniqueness and minimality of large face-width embeddings of graphs.
\newblock {\em Combinatorica}, 15(4):541--556, December 1995.
\newblock \href {https://doi.org/10.1007/BF01192526} {\path{doi:10.1007/BF01192526}}.

\bibitem[Moh01]{Mohar2001Graph}
Bojan Mohar.
\newblock Graph minors and graphs on surfaces.
\newblock In {\em Surveys in {{Combinatorics}}}, number 288 in London {{Mathematical Society Lecture Note Series}}, pages 145--163. Cambridge University Press, 2001.
\newblock \href {https://doi.org/10.1017/CBO9780511721328.008} {\path{doi:10.1017/CBO9780511721328.008}}.

\bibitem[MT01]{MoharT2001Graphs}
Bojan Mohar and Carsten Thomassen.
\newblock {\em Graphs on Surfaces}.
\newblock Johns {{Hopkins}} Studies in the Mathematical Sciences. Johns Hopkins University Press, Baltimore, 2001.

\bibitem[Nor15]{Norin2015New}
Sergey Norin.
\newblock New tools and results in graph minor structure theory.
\newblock In {\em Surveys in {{Combinatorics}} 2015}, pages 221--260. Cambridge University Press, 1 edition, July 2015.
\newblock \href {https://doi.org/10.1017/CBO9781316106853.008} {\path{doi:10.1017/CBO9781316106853.008}}.

\bibitem[OSVV08]{OrecchiaSVV2008Partitioning}
Lorenzo Orecchia, Leonard~J. Schulman, Umesh~V. Vazirani, and Nisheeth~K. Vishnoi.
\newblock On partitioning graphs via single commodity flows.
\newblock In {\em Proceedings of the Fortieth Annual {{ACM}} Symposium on {{Theory}} of Computing}, pages 461--470, Victoria British Columbia Canada, May 2008. ACM.
\newblock \href {https://doi.org/10.1145/1374376.1374442} {\path{doi:10.1145/1374376.1374442}}.

\bibitem[PPT23a]{PaulPT2023Graph}
Christophe Paul, Evangelos Protopapas, and Dimitrios~M. Thilikos.
\newblock Graph {{Parameters}}, {{Universal Obstructions}}, and {{WQO}}, April 2023.
\newblock \href {https://arxiv.org/abs/2304.03688} {\path{arXiv:2304.03688}}.

\bibitem[PPT23b]{PaulPT2023Universal}
Christophe Paul, Evangelos Protopapas, and Dimitrios~M. Thilikos.
\newblock An {{Overview}} of {{Universal}} {{Obstructions}} for {{Graph}} {{Parameters}}, April 2023.
\newblock \href {https://arxiv.org/abs/2304.14121} {\path{arXiv:2304.14121}}.

\bibitem[PPTW24]{PaulPTW2024Obstructionsa}
Christophe Paul, Evangelos Protopapas, Dimitrios~M. Thilikos, and Sebastian Wiederrecht.
\newblock Obstructions to {{Erd{\"o}s-P{\'o}sa Dualities}} for {{Minors}}.
\newblock In {\em 2024 {{IEEE}} 65th {{Annual Symposium}} on {{Foundations}} of {{Computer Science}} ({{FOCS}})}, pages 31--52, Chicago, IL, USA, October 2024. IEEE.
\newblock \href {https://doi.org/10.1109/FOCS61266.2024.00013} {\path{doi:10.1109/FOCS61266.2024.00013}}.

\bibitem[Ram25]{Rambaud2025Excluding}
Cl{\'e}ment Rambaud.
\newblock Excluding a rectangular grid, January 2025.
\newblock \href {https://arxiv.org/abs/2501.11617} {\path{arXiv:2501.11617}}, \href {https://doi.org/10.48550/arXiv.2501.11617} {\path{doi:10.48550/arXiv.2501.11617}}.

\bibitem[Ree92]{Reed1992Finding}
Bruce~A. Reed.
\newblock Finding approximate separators and computing tree width quickly.
\newblock In {\em Proceedings of the Twenty-Fourth Annual {{ACM}} Symposium on {{Theory}} of Computing - {{STOC}} '92}, pages 221--228, Victoria, British Columbia, Canada, 1992. ACM Press.
\newblock \href {https://doi.org/10.1145/129712.129734} {\path{doi:10.1145/129712.129734}}.

\bibitem[RS83]{RobertsonS1983Graph}
Neil Robertson and Paul~D. Seymour.
\newblock Graph minors. {{I}}. {{Excluding}} a forest.
\newblock {\em Journal of Combinatorial Theory, Series B}, 35(1):39--61, August 1983.
\newblock \href {https://doi.org/10.1016/0095-8956(83)90079-5} {\path{doi:10.1016/0095-8956(83)90079-5}}.

\bibitem[RS84]{RobertsonS1984Grapha}
Neil Robertson and Paul Seymour.
\newblock Graph {{Width}} and {{Well-Quasi-Ordering}}: A {{Survey}}.
\newblock In {\em Progress in {{Graph Theory}}}, volume~2, pages 399--406. Academic Press, Toronto, Orlando, 1984.

\bibitem[RS85]{RobertsonS1985Graph}
Neil Robertson and Paul~D. Seymour.
\newblock Graph minors -- a survey.
\newblock In {\em Surveys in {{Combinatorics}} 1985}, pages 153--171. Cambridge University Press, 1 edition, July 1985.
\newblock \href {https://doi.org/10.1017/CBO9781107325678.009} {\path{doi:10.1017/CBO9781107325678.009}}.

\bibitem[RS86]{RobertsonS1986Grapha}
Neil Robertson and Paul~D. Seymour.
\newblock Graph minors. {{V}}. {{Excluding}} a planar graph.
\newblock {\em Journal of Combinatorial Theory, Series B}, 41(1):92--114, August 1986.
\newblock \href {https://doi.org/10.1016/0095-8956(86)90030-4} {\path{doi:10.1016/0095-8956(86)90030-4}}.

\bibitem[RS90a]{RobertsonS1990Grapha}
Neil Robertson and Paul~D. Seymour.
\newblock Graph minors. {{IV}}. {{Tree-width}} and well-quasi-ordering.
\newblock {\em Journal of Combinatorial Theory, Series B}, 48(2):227--254, April 1990.
\newblock \href {https://doi.org/10.1016/0095-8956(90)90120-O} {\path{doi:10.1016/0095-8956(90)90120-O}}.

\bibitem[RS90b]{RobertsonS1990Graph}
Neil Robertson and Paul~D. Seymour.
\newblock Graph minors. {{IX}}. {{Disjoint}} crossed paths.
\newblock {\em Journal of Combinatorial Theory, Series B}, 49(1):40--77, June 1990.
\newblock \href {https://doi.org/10.1016/0095-8956(90)90063-6} {\path{doi:10.1016/0095-8956(90)90063-6}}.

\bibitem[RS91]{RobertsonS1991Graph}
Neil Robertson and Paul~D. Seymour.
\newblock Graph minors. {{X}}. {{Obstructions}} to tree-decomposition.
\newblock {\em Journal of Combinatorial Theory, Series B}, 52(2):153--190, July 1991.
\newblock \href {https://doi.org/10.1016/0095-8956(91)90061-N} {\path{doi:10.1016/0095-8956(91)90061-N}}.

\bibitem[RS95a]{RobertsonS1995Graphb}
Neil Robertson and Paul~D. Seymour.
\newblock Graph {{Minors}}. {{XII}}. {{Distance}} on a {{Surface}}.
\newblock {\em Journal of Combinatorial Theory, Series B}, 64(2):240--272, July 1995.
\newblock \href {https://doi.org/10.1006/jctb.1995.1034} {\path{doi:10.1006/jctb.1995.1034}}.

\bibitem[RS95b]{RobertsonS1995Graph}
Neil Robertson and Paul~D. Seymour.
\newblock Graph {{Minors}}. {{XIII}}. {{The Disjoint Paths Problem}}.
\newblock {\em Journal of Combinatorial Theory, Series B}, 63(1):65--110, January 1995.
\newblock \href {https://doi.org/10.1006/jctb.1995.1006} {\path{doi:10.1006/jctb.1995.1006}}.

\bibitem[RS03a]{RobertsonS2003Grapha}
Neil Robertson and Paul~D. Seymour.
\newblock Graph {{Minors}}. {{XVI}}. {{Excluding}} a non-planar graph.
\newblock {\em Journal of Combinatorial Theory, Series B}, 89(1):43--76, September 2003.
\newblock \href {https://doi.org/10.1016/S0095-8956(03)00042-X} {\path{doi:10.1016/S0095-8956(03)00042-X}}.

\bibitem[RS03b]{RobertsonS2003Graph}
Neil Robertson and Paul~D. Seymour.
\newblock Graph {{Minors}}. {{XVIII}}. {{Tree-decompositions}} and well-quasi-ordering.
\newblock {\em Journal of Combinatorial Theory, Series B}, 89(1):77--108, September 2003.
\newblock \href {https://doi.org/10.1016/S0095-8956(03)00067-4} {\path{doi:10.1016/S0095-8956(03)00067-4}}.

\bibitem[RS04]{RobertsonS2004Graph}
Neil Robertson and Paul~D. Seymour.
\newblock Graph {{Minors}}. {{XX}}. {{Wagner}}'s conjecture.
\newblock {\em Journal of Combinatorial Theory, Series B}, 92(2):325--357, November 2004.
\newblock \href {https://doi.org/10.1016/j.jctb.2004.08.001} {\path{doi:10.1016/j.jctb.2004.08.001}}.

\bibitem[RS09]{RobertsonS2009Graph}
Neil Robertson and Paul~D. Seymour.
\newblock Graph minors. {{XXI}}. {{Graphs}} with unique linkages.
\newblock {\em Journal of Combinatorial Theory, Series B}, 99(3):583--616, May 2009.
\newblock \href {https://doi.org/10.1016/j.jctb.2008.08.003} {\path{doi:10.1016/j.jctb.2008.08.003}}.

\bibitem[RS12]{RobertsonS2012Graph}
Neil Robertson and Paul~D. Seymour.
\newblock Graph {{Minors}}. {{XXII}}. {{Irrelevant}} vertices in linkage problems.
\newblock {\em Journal of Combinatorial Theory, Series B}, 102(2):530--563, March 2012.
\newblock \href {https://doi.org/10.1016/j.jctb.2007.12.007} {\path{doi:10.1016/j.jctb.2007.12.007}}.

\bibitem[RST94]{RobertsonST1994Quickly}
Neil Robertson, Paul~D. Seymour, and Robin Thomas.
\newblock Quickly {{Excluding}} a {{Planar Graph}}.
\newblock {\em Journal of Combinatorial Theory, Series B}, 62(2):323--348, November 1994.
\newblock \href {https://doi.org/10.1006/jctb.1994.1073} {\path{doi:10.1006/jctb.1994.1073}}.

\bibitem[RV90]{RobertsonV1990Representativity}
Neil Robertson and Richardson Vitray.
\newblock Representativity of surface embeddings.
\newblock In {\em Paths, Flows, and {{VLSI-layout}}}, volume~9 of {\em Algorithms {{Combin}}.} Springer Berlin, Bonn, 293-328 edition, 1990.

\bibitem[Sey80]{Seymour1980Disjoint}
Paul~D. Seymour.
\newblock Disjoint paths in graphs.
\newblock {\em Discrete Mathematics}, 29(3):293--309, 1980.
\newblock \href {https://doi.org/10.1016/0012-365X(80)90158-2} {\path{doi:10.1016/0012-365X(80)90158-2}}.

\bibitem[Shi80]{Shiloach1980Polynomial}
Yossi Shiloach.
\newblock A {{Polynomial Solution}} to the {{Undirected Two Paths Problem}}.
\newblock {\em Journal of the ACM}, 27(3):445--456, July 1980.
\newblock \href {https://doi.org/10.1145/322203.322207} {\path{doi:10.1145/322203.322207}}.

\bibitem[SST24]{SauST2024More}
Ignasi Sau, Giannos Stamoulis, and Dimitrios~M. Thilikos.
\newblock A more accurate view of the {{Flat Wall Theorem}}.
\newblock {\em Journal of Graph Theory}, 107(2):263--297, October 2024.
\newblock \href {https://doi.org/10.1002/jgt.23121} {\path{doi:10.1002/jgt.23121}}.

\bibitem[ST96]{SeymourT1996Uniqueness}
P.~D. Seymour and Robin Thomas.
\newblock Uniqueness of highly representative surface embeddings.
\newblock {\em Journal of Graph Theory}, 23(4):337--349, December 1996.
\newblock \href {https://doi.org/10.1002/(SICI)1097-0118(199612)23:4<337::AID-JGT2>3.0.CO;2-S} {\path{doi:10.1002/(SICI)1097-0118(199612)23:4<337::AID-JGT2>3.0.CO;2-S}}.

\bibitem[SW89]{SeeseW1989Grids}
D.G Seese and W~Wessel.
\newblock Grids and their minors.
\newblock {\em Journal of Combinatorial Theory, Series B}, 47(3):349--360, December 1989.
\newblock \href {https://doi.org/10.1016/0095-8956(89)90033-6} {\path{doi:10.1016/0095-8956(89)90033-6}}.

\bibitem[Tho80]{Thomassen19802Linked}
Carsten Thomassen.
\newblock 2-{{Linked Graphs}}.
\newblock {\em European Journal of Combinatorics}, 1(4):371--378, December 1980.
\newblock \href {https://doi.org/10.1016/S0195-6698(80)80039-4} {\path{doi:10.1016/S0195-6698(80)80039-4}}.

\bibitem[Tho06]{Tholey2006Solving}
Torsten Tholey.
\newblock Solving the 2-{{Disjoint Paths Problem}} in {{Nearly Linear Time}}.
\newblock {\em Theory of Computing Systems}, 39(1):51--78, February 2006.
\newblock \href {https://doi.org/10.1007/s00224-005-1256-9} {\path{doi:10.1007/s00224-005-1256-9}}.

\bibitem[Tho09]{Tholey2009Improved}
Torsten Tholey.
\newblock Improved {{Algorithms}} for the 2-{{Vertex Disjoint Paths Problem}}.
\newblock In {\em {{SOFSEM}} 2009: {{Theory}} and {{Practice}} of {{Computer Science}}}, volume 5404, pages 546--557. Springer Berlin Heidelberg, Berlin, Heidelberg, 2009.
\newblock \href {https://doi.org/10.1007/978-3-540-95891-8_49} {\path{doi:10.1007/978-3-540-95891-8_49}}.

\bibitem[TW24a]{ThilikosW2024Excluding}
Dimitrios~M. Thilikos and Sebastian Wiederrecht.
\newblock Excluding {{Surfaces}} as {{Minors}} in {{Graphs}}, February 2024.
\newblock \href {https://arxiv.org/abs/2306.01724v4} {\path{arXiv:2306.01724v4}}.

\bibitem[TW24b]{ThilikosW2024Killing}
Dimitrios~M. Thilikos and Sebastian Wiederrecht.
\newblock Killing a {{Vortex}}.
\newblock {\em Journal of the ACM}, 71(4):1--56, August 2024.
\newblock \href {https://doi.org/10.1145/3664648} {\path{doi:10.1145/3664648}}.

\bibitem[Val79]{Valiant1979Complexity}
Leslie~G. Valiant.
\newblock The {{Complexity}} of {{Enumeration}} and {{Reliability Problems}}.
\newblock {\em SIAM Journal on Computing}, 8(3):410--421, August 1979.
\newblock \href {https://doi.org/10.1137/0208032} {\path{doi:10.1137/0208032}}.

\bibitem[Wag37]{Wagner1937Ueber}
Klaus Wagner.
\newblock {{\"U}ber eine Eigenschaft der ebenen Komplexe}.
\newblock {\em Mathematische Annalen}, 114(1):570--590, December 1937.
\newblock \href {https://doi.org/10.1007/BF01594196} {\path{doi:10.1007/BF01594196}}.

\bibitem[Wol22]{Wollan2022Explicit}
Paul Wollan.
\newblock Explicit {{Bounds}} for {{Graph Minors}}.
\newblock In {\em Surveys in {{Combinatorics}} 2022}, pages 215--236. Cambridge University Press, May 2022.
\newblock \href {https://doi.org/10.1017/9781009093927.008} {\path{doi:10.1017/9781009093927.008}}.

\bibitem[Wol24]{Wollan2024Personal}
Paul Wollan.
\newblock Personal {{Communication}}, September 2024.

\bibitem[Wol25]{Wollan2025Personal}
Paul Wollan.
\newblock Personal {{Communication}}, March 2025.

\bibitem[Woo24]{Wood2024Personal}
David Wood.
\newblock Personal {{Communication}}, October 2024.

\end{thebibliography}

\end{document}